\documentclass{amsbook}

\usepackage{sty/preamble}

\usepackage{geometry}
\setamsbookgeometry 
\savegeometry{abg} 

\title[Homotopy Theory of Enriched Mackey Functors]{Homotopy Theory of Enriched Mackey Functors}
\newcommand{\subtitle}{Closed Multicategories, Permutative Enrichments, and Algebraic Foundations for Spectral Mackey Functors}

\authorinfoNJDY

\hypersetup{pdfauthor=\authors}

\date{08 December 2022}


\makeatletter
\def\@titleimg{%
  \begin{center}

\vfill
\begin{tikzpicture}
  \def\d{5}
  \def\w{1.58}
  \begin{scope}
    \clip (2*\w+\d,\d-3.75*\w/2) rectangle ++(12*\w,6);
    \foreach \x in {0,...,15}{
      \draw
      (\x*\w,0) -- ++(\d,\d) node (a) {} -- ++(\d,-\d);
      \draw[line width=2pt]
      (\x*\w,0) -- ++(\d,\d) node (a) {} -- ++(\d,-\d);
      \draw (a) [fill=black] circle(2.5pt);
    }
  \end{scope} 
\end{tikzpicture}
\vspace{-10pc}

  \end{center}
}
\makeatother

\begin{document}
\pdfbookmark{Title}{Title}
\frontmatter

\begin{abstract}
  Mackey functors provide the coefficient systems for equivariant cohomology theories. More generally, enriched presheaf categories provide a classification and organization for many stable model categories of interest. Changing enrichments along $K$-theory multifunctors provides an important tool for constructing spectral Mackey functors from Mackey functors enriched in algebraic structures such as permutative categories.

This work gives a detailed development of diagrams, presheaves, and Mackey functors enriched over closed multicategories. Change of enrichment, including the relevant compositionality, is treated with care. This framework is applied to the homotopy theory of enriched diagram and Mackey functor categories, including equivalences of homotopy theories induced by $K$-theory multifunctors. Particular applications of interest include diagrams and Mackey functors enriched in pointed multicategories, permutative categories, and symmetric spectra.

\end{abstract}

\maintitlepage 

\cleardoublepage
\thispagestyle{empty}
\vspace*{13.5pc}
\begin{center}
  \begin{itemize}
  \item[] The first author dedicates this book to Nemili, Linus, and Kavya.
  \item[] The second author dedicates this book to Jacqueline.
  \end{itemize}
\end{center}
\cleardoublepage

\pdfbookmark{\contentsname}{Contents}
\tableofcontents

\newcommand{\sect}[1]{\section*{#1}}
\newcommand{\subsect}[1]{\bigskip\noindent\textbf{#1}.}

\newcommand{\prefacepartNumName}[1]{\Cref{#1}: \nameref{#1}}
\newcommand{\prefacepartSubsect}[1]{\subsect{\Cref{#1}}\ \nameref{#1}\\}
\newcommand{\prefacechapNumName}[1]{\medskip\begin{center}\Cref{#1}: \nameref{#1} \end{center}}
\newcommand{\prefaceappNumName}[1]{\Cref{#1}.\ \nameref{#1}}

\chapter*{Preface}

\chapquote{``Just as the social insects build marvellously designed intricate structures by apparently carrying materials around at random so have the mathematicians built a marvellously articulated body of abstract concepts by following their individual instincts with an eye to what their colleagues are doing."}{George W.\ Mackey}{What Do Mathematicians Do? Paris, 1982}

This work develops techniques and basic results concerning the homotopy theory of enriched diagrams and enriched Mackey functors.
Presentation of a category of interest as a diagram category has become a standard and powerful technique in a range of applications.
Diagrams that carry enriched structures provide deeper and more robust applications.
With an eye to such applications, we provide further development of both the categorical algebra of enriched diagrams, and the homotopy theoretic applications in $K$-theory spectra.

The title of this work refers to certain enriched presheaves, known as Mackey functors, whose homotopy theory classifies that of equivariant spectra.
More generally, certain stable model categories are classified as modules---in the form of enriched presheaves---over categories of generating objects.
We provide further review of this motivating context in \cref{ch:motivations} below.

The main body of this work provides a detailed study of enriched diagrams, including enriched presheaves and enriched Mackey functors, and their homotopy theory.
\cref{part:background} provides background on the homotopy-theoretic context, including $K$-theory functors and the homotopy theory of multicategories.
\cref{part:multicat} extends this material to the homotopy theory of pointed multicategories, providing for the later applications a setting that is both conceptually and technically more natural.

The categorical algebra of enriched diagrams is the subject of \cref{part:enrpresheave} and may be of independent interest.
It extends the theory of enrichment over a symmetric monoidal category in two ways.
First, it gives a careful exposition of the theory of enrichment over a multicategory, including fundamental definitions of closed multicategory and self-enrichment.
Second, it carefully explains change of enrichment along a multifunctor and the resulting diagram change of enrichment.
Compositionality of the latter is more subtle, and treated in detail.

Applications to the homotopy theory of enriched diagrams, Mackey functors, and change of enrichment are the focus of \cref{part:homotopy-mackey}.
These arise from enrichments over permutative categories and pointed multicategories, with change of enrichment along $K$-theory multifunctors.

\sect{Audience}

This work is aimed at graduate students and researchers with an interest in category theory, algebraic $K$-theory, and homotopy theory. 
Our highly detailed exposition is designed to make this work accessible to a wide audience.

\sect{Part and Chapter Summaries}

This work consists of the following.
\begin{itemize}
\item \prefacepartNumName{ch:motivations}
\item \prefacepartNumName{part:background}
\item \prefacepartNumName{part:multicat}
\item \prefacepartNumName{part:enrpresheave}
\item \prefacepartNumName{part:homotopy-mackey}
\item \Cref{ch:prelim,ch:prelim_enriched,ch:prelim_multicat} on \nameref{ch:prelim}, \nameref{ch:prelim_enriched}, and \nameref{ch:prelim_multicat}
\item \prefacepartNumName{ch:questions}
\end{itemize}

Below is a brief summary of each part and chapter.
Following these summaries, we outline the interdependence of their content.
In the main text, each chapter starts with an introduction that describes more thoroughly its content and connections with other chapters.

\prefacechapNumName{ch:motivations}

This chapter reviews the use of enriched diagrams and enriched Mackey functors in equivariant homotopy theory and the theory of stable model categories.
The results outlined here are not used directly in the main body of the work below, but they provide an important motivating context.
The goal of this chapter is, therefore, to outline some of the main ideas and provide numerous references to the literature for further treatment.

\prefacepartSubsect{part:background}
This part provides background that is essential for this text, and is more specialized than that of the appendices.
The main inputs for modern $K$-theory spectra are multicategories and permutative categories reviewed in \cref{ch:ptmulticat}.
The construction of $K$-theory functors---also known as infinite loop space machinery---is reviewed in \cref{ch:Kspectra}.
More modern applications to the homotopy theory of multicategories are summarized in \cref{ch:multperm}.

\prefacechapNumName{ch:ptmulticat}

This chapter reviews the 2-category of small multicategories, including three important special cases.
These are pointed multicategories, left $\Mone$-modules, and permutative categories with multilinear functors.
These variants are related by various free, forgetful, and endomorphism functors that will be used throughout the rest of this work.

\prefacechapNumName{ch:Kspectra}

This chapter reviews the $K$-theory functors $\Kse$ and $\Kem$, due to Segal and Elmendorf-Mandell, respectively.
These are also called infinite loop space machines because they produce connective spectra from permutative categories and multicategories.
Each is constructed as a composite of other functors, via certain diagram categories, that we describe.

\prefacechapNumName{ch:multperm}

This chapter reviews equivalences of homotopy theories between $\Multicat$, the category of small multicategories and multifunctors, $\permcatst$, the category of small permutative categories and strict monoidal functors, and $\permcatsu$, the category of small permutative categories and strictly unital symmetric monoidal functors.
These equivalences are given by a free left adjoint to the endomorphism functor.
This material provides important foundation for that of \cref{part:multicat}. 

\prefacepartSubsect{part:multicat}
The purpose of these two chapters is to extend the equivalences of homotopy theories from \cref{ch:multperm} to the context of pointed multicategories.
\Cref{ch:ptmulticat-sp} develops the essential extensions to the pointed case.
\Cref{ch:ptmulticat-alg} develops the multifunctoriality results of the pointed free functor, along with multinaturality of the adjunction unit and counit.

\prefacechapNumName{ch:ptmulticat-sp}

This chapter extends the material of \cref{ch:multperm} to a pointed free construction, $\Fst$, from pointed multicategories to permutative categories.
This is not a restriction, along the inclusion of pointed multicategories among all multicategories, but an extension, along the functor that adjoins a disjoint basepoint.
Essential results, such as the adjunction with the endomorphism construction and compatibility with stable equivalences, are likewise extended from  \cref{ch:multperm}.

\prefacechapNumName{ch:ptmulticat-alg}

This chapter shows that the pointed free construction from \cref{ch:ptmulticat-sp}, $\Fst$, is a non-symmetric multifunctor.
Furthermore, $\Fst$ provides equivalences of homotopy theories between categories of non-symmetric algebras in pointed multicategories and permutative categories.
This is the basis for applications to enriched diagrams in \cref{ch:mackey_eq}.

\prefacepartSubsect{part:enrpresheave}
This part covers the categorical algebra of enrichment over multicategories (\cref{ch:menriched}), change of enrichment along a multifunctor (\cref{ch:change_enr}), closed multicategories (\cref{ch:gspectra}), and self-enrichment thereof (\cref{ch:std_enrich}).
These are combined in \cref{ch:gspectra_Kem} to develop the diagram change of enrichment for non-symmetric multifunctors (\cref{gspectra-thm-v}) and the presheaf change of enrichment for multifunctors (\cref{gspectra-thm-v-cor}).
These results are the foundation for applications to homotopy theory in \cref{part:homotopy-mackey}.

\prefacechapNumName{ch:menriched}

This chapter gives the basic definitions and results for enrichment in a non-symmetric multicategory $\M$.
\Cref{EndV-enriched} shows that this material agrees, in the case that $\M$ is the endomorphism multicategory of a monoidal category, with classical enriched category theory over $\V$.
The main application takes $\M$ to be the multicategory of permutative categories and strictly unital symmetric monoidal functors, which is not an endomorphism multicategory.
\Cref{sec:enriched-in-perm,sec:perm-self-enr,sec:eval-perm} treat this case in detail.

\prefacechapNumName{ch:change_enr}

This chapter develops the first collection of results around change of enrichment along a (possibly non-symmetric) multifunctor.
As in \cref{ch:menriched}, it is shown that this theory extends the classical theory for enrichment over (possibly symmetric) monoidal categories.
Compositionality and 2-functoriality for the change-of-enrichment constructions are treated in \cref{sec:functoriality-change-enr,sec:twofunc-change-enr}, respectively.

\prefacechapNumName{ch:gspectra}

This chapter extends the general multicategorical enrichment theory from \cref{ch:menriched} to enrichment over \emph{closed} multicategories.
Because enrichment over permutative categories is both illustrative of the general theory and essential for the further applications, this chapter focuses on that case in detail.

\prefacechapNumName{ch:std_enrich}

This chapter describes the theory of self-enrichment for closed multicategories, and of standard enrichment for multifunctors between closed multicategories.
The self-enrichment of the multicategory of permutative categories, from \cref{ch:gspectra}, is a special case.
Compositionality of standard enrichment is discussed in \cref{sec:fun-std-enr-multi}, and applied to the factorization of Elmendorf-Mandell $K$-theory in \cref{sec:factor-Kemse}.

\prefacechapNumName{ch:gspectra_Kem}

This chapter provides the main results of \cref{part:enrpresheave}.
These make use of the preceding material on enrichment over (closed) multicategories, and apply it to categories of enriched diagrams and enriched Mackey functors.
A key detail, both here and in the homotopical applications of \cref{part:homotopy-mackey} is that non-symmetric multifunctors provide a diagram change of enrichment, but not necessarily a change of enrichment for enriched Mackey functors (presheaves).
The essential reason is that symmetry of a multifunctor is required for commuting the $(-)^\op$ in the domain of enriched presheaves with change of enrichment.
\Cref{sec:presheaf-K,sec:mult-mackey-spectra} give applications to Elmendorf-Mandell $K$-theory, with attention to the relevant symmetry conditions among other details.

\prefacepartSubsect{part:homotopy-mackey}
The two final chapters of this work apply the preceding categorical algebra and homotopical constructions.
\Cref{ch:mackey} develops applications to the homotopy theory of enriched diagrams and Mackey functors in general.
\Cref{ch:mackey_eq} gives further detailed applications to enriched diagrams and Mackey functors of pointed multicategories and permutative categories.

\prefacechapNumName{ch:mackey}

This chapter establishes the general theory for a pair of non-symmetric multifunctors $(E,F)$ to provide inverse equivalences of homotopy theories between enriched diagram categories.
The main result is \cref{mackey-gen-xiv} and does not require $E$ or $F$ to satisfy the symmetry condition of a multifunctor.
A similar result for enriched Mackey functor categories, in \cref{mackey-xiv-cor}, requires that $E$, but not necessarily $F$, is a multifunctor.
This is important for the applications, \cref{mackey-xiv-pmulticat,mackey-xiv-mone} below.
There, $E$ is an endomorphism multifunctor and $F$ is a corresponding free non-symmetric multifunctor.

\prefacechapNumName{ch:mackey_eq} 

This chapter applies the general theory from \cref{ch:mackey} to change of enrichment along the inverse equivalences of homotopy theories developed in \cref{part:multicat}.
The main results, \cref{mackey-xiv-pmulticat,mackey-xiv-mone,mackey-pmulti-mone},
establish equivalences of homotopy theories for enriched diagram categories and Mackey functor categories over
pointed multicategories,
permutative categories, and
$\Mone$-modules.

\bigskip
\begin{center}
  \nameref{part:appendices}
\end{center} 
This work includes the following four appendices of supplemental background material and further open questions.
\begin{description}
\item[\prefaceappNumName{ch:prelim}]
  This appendix reviews basic concepts related to monoidal categories and 2-categories.
\item[\prefaceappNumName{ch:prelim_enriched}]
  Here we review the classical theory of categories enriched over monoidal categories.
\item[\prefaceappNumName{ch:prelim_multicat}]
  This appendix gives background on multicategories, including enriched multicategories, endomorphism multicategories, and pointed multicategories.
\item[\prefaceappNumName{ch:questions}]
  In this appendix we discuss a number of open questions related to the topics of this work.
  They provide further motivation for the main text.
\end{description}

\sect{Chapter Interdependence}

The material in \cref{ch:motivations} is not prerequisite for the main text, but is part of the broader context in which this work is situated.
\Cref{ch:prelim,,ch:prelim_enriched,,ch:prelim_multicat} contain background material that will be used throughout.

As noted in the summaries above, part of this work involves abstract categorical algebra of multicategorical enrichments, which may be of independent interest.
In the following table, we separate chapters and sections into columns according to whether they involve only (multi-)categorical algebra, or additional homotopy-theoretic concepts.
\begin{center}
{\renewcommand{\arraystretch}{1.1}%
\setlength{\tabcolsep}{3ex}%
\begin{tabular}{l|ll}
  & \textbf{(Multi-) Categorical Algebra}
  & \textbf{Homotopy Theory}
  \\ \hline
  \\[-.5pc]
  & \cref{ch:prelim,,ch:prelim_enriched,,ch:prelim_multicat}
  &
  \\
  \\
  \cref{part:background}
  & \cref{ch:ptmulticat}
  & \cref{ch:Kspectra,ch:multperm}
  \\
  \\
  \cref{part:multicat}
  &
  & \cref{ch:ptmulticat-sp,ch:ptmulticat-alg}
  \\
  \\
  \cref{part:enrpresheave}
  & \cref{ch:menriched,,ch:change_enr,,ch:gspectra}
  &
  \\
  & \cref{sec:selfenr-clmulti,sec:std-enr-multifunctor,sec:fun-std-enr-multi}
  & \cref{sec:factor-Kemse}
  \\
  & \cref{sec:enr-diag-psh,sec:enr-diag-change-enr,sec:diag-psh-change-enr-functors,sec:fun-change-enr-diag}
  & \cref{sec:presheaf-K,sec:mult-mackey-spectra}
  \\
  \\
  \cref{part:homotopy-mackey}
  & \cref{sec:comparing-enr-diag-cat,sec:unidg,sec:coudg}
  & \cref{sec:heq-enrdiag-cat}
  \\
  &
  & \cref{ch:mackey_eq} 
\end{tabular}
}
\end{center}

Each entry in the above table depends on those to its left and above.
Thus, the material in the left column may be read independently of that in the right.
The introduction of each chapter contains subsections titled \emph{Connection with Other Chapters} and \emph{Background} that give more detailed discussions of the respective dependencies.
\sect{Related  Literature}

Here we list a selection of references for background or further reading.
\begin{description}
\item[\textnormal{2-Dimensional Categories}] 
  \cite{johnson-yau}
\item[\textnormal{Monoidal Categories and Enriched Multicategories}]
  \cite{cerberusIII,yau-operad}
\item[\textnormal{Stable Homotopy Theory}]
  \cite{barnes_roitzheim,mayconcise}
\item[\textnormal{Equivariant Homotopy Theory}]
  \cite{tomdieck,lewis-may-steinberger,alaska-notes,hhrbook} 
\item[\textnormal{Algebraic $K$-Theory Spectra}]
  \cite{cerberusIII,johnson-yau-multiK,bousfield_friedlander,may-permutative,mandell_inverseK,quillenKI,segal,thomason,waldhausen}
\item[\textnormal{Multifunctorial $K$-Theory}]
  \cite{elmendorf-mandell,elmendorf-mandell-perm,cerberusIII,johnson-yau-Fmulti,johnson-yau-invK,johnson-yau-multiK,johnson-yau-permmult}
\item[\textnormal{Spectral Mackey Functors}]
  \cite{schwede-shipley_stable,bohmann_osorno-mackey,barwick,merling-malkiewich-Athy,barwick-glasman-shah,merling-malkiewich,guillou_may,gmmo}
\end{description}

\setcounter{chapter}{-1} 
\chapter{Motivations from Equivariant Topology}
\label{ch:motivations}

In this chapter we describe context from equivariant topology and the theory of stable model categories that motivates our further study of multicategorically enriched categories, enriched diagrams, enriched Mackey functors, and change of enrichment.

\begin{convention}\label{convention:G-finite}
  Assume throughout this chapter that $G$\label{not:G} is a finite group.
  See \cref{remark:finiteness,remark:self-duality} for further comments on this convention.
\end{convention}

\subsection*{Connection with Main Content}
The purpose of this chapter is to indicate the role that categorical diagrams---particularly Mackey functors---play in equivariant homotopy theory.
None of the mathematics in this present work depends on the content of this chapter, but the attendant applications are a key motivation.

For example, the Burnside 2-category $\GE$ (\cref{definition:Burnside2}) is enriched in the multicategory of permutative categories, $\permcatsu$ (\cref{sec:multpermcat}).
We give a treatment of
\begin{itemize}
\item categories enriched in closed multicategories, in \cref{ch:menriched},
\item change of enrichment, in \cref{ch:change_enr}, 
\item the closed multicategory structure of $\permcatsu$ in \cref{ch:gspectra}, and
\item self-enrichment for closed multicategories in \cref{ch:std_enrich}.
\end{itemize}

In the Guillou-May \cref{theorem:GM}, the domain of spectral Mackey functors, $(\GE_{\bbK})^\op$, is given by a change of enrichment $(-)_{\bbK}$ and requires a distinction between enriched diagrams, with domain $\GE_{\bbK}$, and enriched Mackey functors, with domain $(\GE_{\bbK})^\op$.
We describe the relevant subtleties further in \cref{remark:GE-non-self-dual,rk:BO7.5}.

Similarly, but in a more abstract context, the spectral presheaves in the Schwede-Shipley Characterization \cref{theorem:schwede-shipley} have domain $\EP^\op$.
The input $\EP$ is the spectral endomorphism category of a set of compact generators $P$ for a simplicial, cofibrantly generated, proper, and stable model category $\M$.

We give a general treatment of enriched diagrams and enriched Mackey functors, including interactions with change of enrichment, in \cref{ch:gspectra_Kem}.
We develop techniques and applications for the corresponding homotopy theory in \cref{ch:mackey,ch:mackey_eq}.

\subsection*{Chapter Summary}
A substantive treatment of equivariant homotopy theory is well beyond our current scope.
At the end of this introduction we give a list of key references.
The remaining content in this chapter is restricted to those definitions and results that provide motivating context for our work below.

\cref{sec:equivar-top} concerns equivariant spaces.
\begin{itemize}
\item The orbit category of $G$ is denoted $\OrbG$; see \cref{definition:orbit-cat}.

\item Elmendorf's \cref{theorem:elmendorf} shows that the homotopy theory of $G$-spaces is equivalent to that of topological presheaves on $\OrbG$.
\end{itemize}

\cref{sec:Burnside1} concerns Abelian Mackey functors.
\begin{itemize}
\item The Burnside ring of $G$ is denoted $\GA$.  Its elements are isomorphism classes of finite $G$-sets with disjoint union and Cartesian product; see \cref{definition:Burnside1}.
  
\item The Burnside category of $G$ is denoted $\GB$.  Its morphisms are isomorphism classes of spans between finite $G$-sets.
  Disjoint union provides an enrichment over Abelian groups; see \cref{definition:Burnside1}.
\item Abelian Mackey functors are enriched presheaves on the Burnside category; see
  \cref{definition:Ab-G-Mackey}.
\end{itemize}

\cref{sec:Burnside2,sec:stmodcat} concern spectral Mackey functors.
\begin{itemize}
\item The Burnside 2-category of $G$ is denoted $\GE$.  Its 1- and 2-cells are categories of spans between finite $G$-sets.
  Disjoint union, together with a choice of pullbacks and whiskering by a strict unit, provides an enrichment over permutative categories; see
  \cref{definition:Burnside2}.
  
\item Spectral Mackey functors are enriched presheaves on a spectral enrichment of the Burnside 2-category; see
  \cref{definition:Sp-G-Mackey}.
  
\item The Guillou-May \cref{theorem:GM} shows that the homotopy theory of $G$-spectra is equivalent to that of spectral Mackey functors.

\item The Schwede-Shipley Characterization \cref{theorem:schwede-shipley} shows that the homotopy theory of a simplicial, cofibrantly generated, proper, and stable model category is equivalent to that of spectral presheaves on an endomorphism category of generating objects.
\end{itemize}

\subsection*{References}
Main references for equivariant homotopy theory include, at least, the following.
We include further specialized references at relevant points in the discussion below.
\begin{itemize}
\item The text by tom~Dieck \cite{tomdieck} lays the foundations for equivariant homotopy theory of spaces, including equivariant (co)homology theories known as Bredon cohomology.
\item The monograph \cite{lewis-may-steinberger}, by Lewis, May, and Steinberger, gives the foundational treatment of equivariant stable homotopy theory, particularly equivariant spectra.
\item The CBMS Alaska conference proceedings \cite{alaska-notes} refines and significantly extends the preceding theory, including more development of the closed monoidal structure for equivariant spectra.
\item The recent textbook account by Hill-Hopkins-Ravenel \cite{hhrbook} provides a more modern perspective, with thorough treatment of norm operations and the slice filtration that are essential in their solution of the Kervaire invariant problem \cite{hhr-kervaire}.
\end{itemize}

\section{Equivariant Spaces and Presheaves on the Orbit Category}\label{sec:equivar-top}

Recall \cref{convention:G-finite} that $G$ is assumed to be a finite group.

\begin{remark}\label{remark:finiteness}
  Many, but not all, of the concepts below extend to more general cases of interest, such as $G$ being a compact Lie group or a general topological group.
  The most important exception is our definition of the Burnside category in \cref{definition:Burnside2,definition:Burnside1}, which depends on finiteness of $G$.
  See \cref{remark:self-duality} for further comments and references regarding that point.
\end{remark}

\begin{definition}\label{definition:diagram-presheaf}
  Suppose $\C$ and $\M$ are categories, with $\C$ small.
  A \emph{diagram of shape $\C$}\index{diagram}, or \emph{$\C$-diagram} in $\M$ is a functor
  \[
    \C \to \M.
  \]
  A \emph{presheaf on $\C$}\index{presheaf} or \emph{$\C$-presheaf} in $\M$ is a diagram of shape $\C^\op$ in $\M$, where $\C^\op$ is the opposite category of $\C$.
  That is, a presheaf on $\C$ is a functor
  \[
    \C^\op \to \M.
  \]
  The phrase ``in $\M$'' is often omitted when $\M$ is clear from context.
  Morphisms between diagrams and presheaves are natural transformations, and so
  \[
    \Cat( \C , \M ) \andspace
    \Cat( \C^\op , \M )
  \]
  are the respective categories of diagrams and presheaves on $\C$.
  If $\M$ is a symmetric monoidal closed category (\cref{def:closedcat}) or, more generally, a closed multicategory (\cref{def:closed-multicat}), then there are corresponding enriched variants described in \cref{def:enr-diag-cat}.
\end{definition}

\begin{definition}\label{definition:orbit-cat}
  The \emph{orbit category}\index{orbit!category} of a group $G$, denoted $\OrbG$,\label{not:OrbG} consists of the following.
  Its objects are the $G$-orbits $G/H$, where $H$ is a subgroup of $G$, and its morphisms are the $G$-equivariant morphisms.
\end{definition}

\begin{remark}
  Note that each $G$-equivariant map
  \[
    f\cn G/H \to G/K \inspace \OrbG
  \]
  determines and is determined by an element $g \in G$, where $f(eH) = gK$, such that $g^\inv Hg \subset K$.
  Thus, the morphisms in $\OrbG$ are given by subconjugacy relations.
\end{remark}

\begin{definition}[$G$-Spaces]\label{definition:Gtop}
  A \emph{$G$-space}\index{G-space@
    $G$-space}\index{equivariant space} is a topological space on which $G$ acts continuously.
  Morphisms of $G$-spaces are continuous functions that commute with the $G$-action.
  The category of $G$-spaces and their morphisms is denoted $\TopG$\label{not:TopG}.
\end{definition}

\begin{definition}[Fixed Points]\label{definition:fixed-points}
  For each $G$-space $X$, and for each subgroup $H$ in $G$, the \emph{$H$-fixed point space}\index{fixed point!space}, denoted $X^H$,\label{not:XH} consists of the subspace of points on which $H$ acts trivially.
  As a $G$-space, $X^H$ can be defined equivalently as the space of $G$-equivariant morphisms
  \[
    \TopG(G/H,X),
  \]
  where $G/H$ has the discrete topology.
  The assignment
  \[
    G/H \mapsto X^H
  \]
  determines a presheaf of spaces on the orbit category $\OrbG$,
  \begin{equation}\label{eq:PhiX}
    \Phi X \cn \OrbG^\op \to \TopG
  \end{equation}
  called the \emph{fixed point functor}\index{fixed point!functor}.
\end{definition}

Discussion of equivariant homotopy and (co)homology is beyond our current scope, but the following gives an indication of the role that the orbit category plays in equivariant topology.
\begin{explanation}
  The coefficient systems\index{coefficient system} for Bredon cohomology\index{Bredon cohomology}\index{cohomology!Bredon} of $G$-spaces are given by presheaves 
  \[
    A\cn\OrbG^\op \to \Ab,\label{not:Ab}
  \]
  where $\Ab$ is the category of Abelian groups and group homomorphisms.
  In particular, for a $G$-space $X$, the composite with $\pi_n$ for $n \ge 2$ yields a coefficient system
  \[
    \OrbG^\op \fto{\Phi X} \Top \fto{\pi_n} \Ab.\dqed
  \]
\end{explanation}

The following result due to Elmendorf \cite{elmendorf-systems} gives a different indication of the importance of presheaves on the orbit category.
\begin{theorem}[{\cite{elmendorf-systems,GMRenriched}}]\label{theorem:elmendorf}
  The fixed points functor, $\Phi$, induces a Quillen equivalence\label{not:htyQ}
  \[
    \Phi\cn \TopG \fto{\hty_Q} (\TopG\mh\Cat)(\OrbG^\op,\Top)
  \]
  between the category of $G$-equivariant topological spaces and the category of topological presheaves on $\OrbG$.
\end{theorem}
As we outline below, presheaves on the Burnside (2-)category, which are known as Mackey functors, fill an analogous role in the generalization to stable equivariant homotopy.

\section{The Burnside Category and Abelian \texorpdfstring{$G$}{G}-Mackey Functors}\label{sec:Burnside1}

The Burnside category (\cref{definition:Burnside1} below) extends the orbit category of $G$ using spans of finite $G$-sets.
The key motivation for this expansion of $\OrbG$ is to account for the restriction, induction, and transfer morphisms on finite $G$-sets.
Further explanation and examples of this perspective can be found in \cite{webb-mackey} and \cite[Sections~8.1 and~8.2]{hhrbook}.

\begin{definition}[Finite $G$-Sets]\label{definition:FinG}
  Let $\FinG$\label{not:FinG}\index{finite $G$-sets}\index{G-sets@$G$-sets!finite} denote the following skeleton of the category of finite $G$-sets.
  The objects of $\FinG$ are pairs $(\ufs{n},\al)$, where $n$ is a natural number, $\ufs{n} = \{ 1, \ldots, n \}$, and 
  \[
    \al \cn G \to \Si_n
  \]
  is a group homomorphism.
  We regard $X = (\ufs{n},\al)$ as a $G$-set with the action
  \[
    g\cdot i = \al(g)(i)
  \]
  for $g \in G$ and $i \in \ufs{n}$.
  The morphisms $f \cn (\ufs{n},\al) \to (\ufs{p},\be)$ in $\FinG$ are $G$-equivariant morphisms.
  That is, $f$ is a map of sets $\ufs{n} \to \ufs{p}$ such that
  \[
    \be(g)(f(i)) = f\big( \al(g)(i) \big) \forspace g \in G \andspace i \in \ufs{n}.
  \]

  We call $n$ the \emph{cardinality} of $X = (\ufs{n},\al)$ and write $|X| = n$.
  Additionally, we define the following.
  \begin{enumerate}
  \item\label{it:FinG-disj} The disjoint union of finite $G$-sets, $\bincoprod$, defines a permutative structure with unit given by the empty $G$-set.
    We write $(\ufs{0},!)$ for the empty finite set and the unique action homomorphism $G \to \Si_0$.
    
  \item\label{it:FinG-cart} The Cartesian product, together with the lexicographic ordering\index{lexicographic ordering}\index{ordering!lexicographic}
    \begin{equation}\label{eq:lexord}
      \ufs{n}\times\ufs{p} \iso \ufs{np} \quad\text{via}\quad
      (i,j)\mapsto p(i-1) + j,
    \end{equation}
    defines a second permutative structure on $\FinG$.
    Its unit is the terminal $G$-set $(\ufs{1},!)$, consisting of the terminal set and the unique action homomorphism $G \to \Si_1$.\dqed
  \end{enumerate}
\end{definition}

\begin{definition}[Bicategory of Spans]\label{definition:Span}
  Suppose $\C$ is a small category with pullbacks, equipped with a choice of pullbacks for each pair of morphisms having a common codomain.
  The \emph{bicategory of spans in $\C$}\index{bicategory of spans}\index{spans!bicategory of} is denoted $\Span(\C)$\label{not:SpanC} and consists of the following.
  \begin{description}
  \item[0-Cells] The 0-cells are objects $X \in \C$.
  \item[1-Cells] The 1-cells with domain $X$ and codomain $Y$ are triples $(A,f,g)$ that are spans
    \begin{equation}\label{eq:Span1}
      X \fot{f} A \fto{g} Y \inspace \C.
    \end{equation}
    Since the object $A$ is determined by the two morphisms, a span is sometimes denoted by its pair of morphisms, $(f,g)$.
  \item[2-Cells] The 2-cells $(A,f,g) \to (A',f',g')$ are morphisms $w\cn A \to A'$ in $\C$ that make the following diagram commute in $\C$.
    \[
      \begin{tikzpicture}[x=13ex,y=7ex]
        \draw[0cell] 
        (0,0) node (x) {X}
        (1,.5) node (a) {A}
        (1,-.5) node (a') {A'}
        (2,0) node (y) {Y}
        ;
        \draw[1cell] 
        (a) edge['] node {f} (x)
        (a') edge node {f'} (x)
        (a) edge node {g} (y)
        (a') edge['] node {g'} (y)
        (a) edge node {w} (a')
        ;
      \end{tikzpicture}
    \]
  \item[Identities] The identity 1-cell on a 0-cell $X$ is the triple $\De_X = (X,1_X,1_X)$ given by the identity morphisms in $\C$.
    The identity 2-cell on a 1-cell $(A,f,g)$ is the identity morphism $1_A$ in $\C$.
  \item[Composition] For objects $X$, $Y$, and $Z$ in $\C$, the composition functor
    \[
      \Span(\C)(Y,Z)\times \Span(\C)(X,Y) \to \Span(\C)(X,Z)
    \]
    sends a composable pair to the span given by their chosen pullback, as shown below.
    \[
      \begin{tikzpicture}[x=20ex,y=10ex]
        \draw[0cell] 
        (0,0) node (x) {X}
        (x)+(1,0) node (y) {Y}
        (y)+(1,0) node (z) {Z}
        (x)+(.5,.5) node (a) {A}
        (y)+(.5,.5) node (b) {B}
        (a)+(.5,.5) node (p) {A \times_Y B}
        ;
        \draw[1cell] 
        (a) edge['] node {f} (x)
        (a) edge node {g} (y)
        (b) edge['] node {h} (y)
        (b) edge node {k} (z)
        (p) edge node {} (a)
        (p) edge node {} (b)
        ;
      \end{tikzpicture}
    \]
  \end{description}
  Having a chosen pullback for each pair of morphisms with a common codomain makes the composition of 1-cells well defined.
  Universality of pullbacks makes it associative and unital up to isomorphisms that satisfy the axioms of bicategorical composition.
  See \cite[Example 2.1.22]{johnson-yau} for further details of this construction.
\end{definition}

Now we use \cref{definition:FinG,definition:Span} to define the Burnside category and its specialization, the Burnside ring.
In \cref{definition:Burnside2} below we generalize further to a Burnside 2-category.
\begin{definition}[Burnside Category and Burnside Ring]\label{definition:Burnside1}
  The \emph{Burnside category}\index{Burnside!category}\index{category!Burnside} of a finite group $G$, denoted $\GB$\label{not:GB}, is an $\Ab$-enriched category defined as follows.
  The objects of $\GB$ are the finite $G$-sets $X \in \FinG$.
  The Abelian group $\GB(X,Y)$ for $X,Y \in \FinG$ is the Grothendieck group of isomorphism classes of spans
  \[
    X \leftarrow A \to Y \inspace \FinG.
  \]
  Thus, $\GB$ is the category obtained from $\Span(\FinG)$ by taking isomorphism classes of 1-cells and then group-completing each set of morphisms with respect to the Abelian monoid structure given by disjoint union.

  The \emph{Burnside ring}\index{Burnside!ring}\index{ring!Burnside} of $G$, denoted $\GA$\label{not:GA}, is obtained by taking isomorphism classes of objects in $\GB$.
  Equivalently, the additive group of elements is given by the Grothendieck group of isomorphism classes of finite $G$-sets, with addition given by disjoint union.
  Its multiplication is induced by Cartesian product.
\end{definition}

\begin{lemma}[Self-Duality of $\GB$]\label{lemma:GB-self-dual}
  There is an isomorphism of $\Ab$-enriched categories
  \begin{equation}\label{eq:GB-self-dual}
    \GB \fto{\iso} \GB^\op
  \end{equation} 
  that is the identity on objects and is induced on hom Abelian groups by the isomorphism
  \[
    \Span(\FinG)(X,Y) \fto{\iso} \Span(\FinG)(Y,X)
  \]
  that sends a span $(f,g)$ to its reverse, $(g,f)$.
\end{lemma}
\begin{proof}
  Functoriality of the indicated isomorphism follows from universality of the pullbacks defining composition.
\end{proof}

We warn the reader that the 2-categorical analog of the self-duality \cref{eq:GB-self-dual} does \emph{not} hold for the Burnside 2-category $\GE$ in \cref{definition:Burnside2} below.
Sending $(f,g)$ to its reverse $(g,f)$ does not define a 2-functor in that context; see \cref{remark:GE-non-self-dual}.

\begin{remark}[Self-Duality and Stable Orbit Spectra]\label{remark:self-duality}
  Self-duality of the Burnside category $\GB$ \cref{eq:GB-self-dual} is nearly transparent in its simplicity, but it is an algebraic artifact of a much deeper topological phenomenon.
  Each orbit $G/H$ has an equivariant suspension spectrum, $\Si^\infty G/H_+$, and there is an equivalent definition of $\GB$ with morphisms given by stable equivariant morphisms $\Si^\infty G/H_+ \to \Si^\infty G/K_+$; see \cite[Section~XIX.3]{alaska-notes}.
  The stable orbit spectra\index{stable orbit spectra}\index{orbit!stable - spectra} $\Si^\infty G/H_+$\label{not:suspGH} satisfy an equivariant self-duality (\cite[Section~XVI.7]{alaska-notes} or \cite[Section~8.0C]{hhrbook}) that implies that of \cref{lemma:GB-self-dual}.

  The definition of the Burnside category\index{Burnside!category}\index{category!Burnside} in terms of stable orbit spectra is the more general one, with origins in work of tom~Dieck \cite{tomdieck}; see \cite[Section~XVII.2]{alaska-notes}.
  The proofs that this definition can be given equivalently by spans of finite $G$-sets, as in \cref{definition:Burnside1}, depend on the assumption that $G$ is finite.
  In more general cases, the definition of the Burnside category in terms of stable orbit spectra is necessary.
\end{remark}

\begin{definition}\label{definition:Ab-G-Mackey}
  An \emph{Abelian $G$-Mackey functor}\index{Abelian Mackey functor}\index{Mackey functor!Abelian} is an $\Ab$-enriched presheaf
  \[
    \GB^\op \to \Ab.
  \]
  Because $\GB$ is isomorphic to $\GB^\op$ (\cref{lemma:GB-self-dual}), an Abelian $G$-Mackey functor is equivalently defined as a functor $\GB \to \Ab$.
\end{definition}

\begin{remark}\label{remark:HM}
  Each Abelian $G$-Mackey functor $M$ has an associated Eilenberg-Mac~Lane $G$-spectrum, $HM$.
  See \cite[Section~V.4]{alaska-notes} or \cite[Theorem~8.8.4]{hhrbook} for constructions via Elmendorf's \cref{theorem:elmendorf}.
  Such Mackey functors $M$, and their associated $G$-spectra $HM$, are the coefficient systems for Bredon cohomology of $G$-spectra.\index{Bredon cohomology}\index{cohomology!Bredon}
\end{remark}

\begin{explanation}\label{explanation:Mackey-axioms}
  An Abelian $G$-Mackey functor $M$ can be defined equivalently
  as a pair of functors\index{Abelian Mackey functor}\index{Mackey functor!Abelian}
  \[
    M_*\cn \FinG \to \Ab
    \andspace
    M^*\cn \FinG^\op \to \Ab
  \]
  that agree on objects and are subject to the following two axioms, where
  \[
    MX = M_*X = M^*X\label{not:Mstar}
  \]
  denotes the common value on objects.
  \begin{enumerate}
  \item For each pair of objects $X$ and $Y$ in $\FinG$, applying $M_*$ to the structure morphisms of the coproduct
    \[
      X \to X \bincoprod Y \leftarrow Y
    \]
    induces a universal morphism with domain $MX \oplus MY$ that is an isomorphism
    \[
      MX \oplus MY \fto{\iso} M(X \bincoprod Y).
    \]
  \item For each pullback diagram in $\FinG$,
    \[
      \begin{tikzpicture}[x=15ex,y=8ex]
        \draw[0cell] 
        (0,0) node (z) {Z}
        (z)+(0,1) node (x) {X}
        (z)+(-1,0) node (y) {Y}
        (y)+(0,1) node (w) {W}
        ;
        \draw[1cell] 
        (w) edge node {p} (x)
        (w) edge['] node {q} (y)
        (y) edge node {g} (z)
        (x) edge node {f} (z)
        ;
      \end{tikzpicture}
    \]
    the following equality of composite morphisms holds in $\Ab$:
    \[
      (M^*f)(M_*g) = (M_*p)(M^*q).
    \]
  \end{enumerate}
  See \cite[Section~2]{webb-mackey} and \cite[Definition 8.2.3]{hhrbook} for further discussion of this perspective, explanation of the equivalence with \cref{definition:Burnside1}, and several compelling examples.
\end{explanation}

\section{Equivariant Spectra and Presheaves on the Burnside 2-Category}\label{sec:Burnside2}

For the category $\C = \FinG$, there is a choice of pullbacks that makes $\Span(\FinG)$ nearly a 2-category.
Following Guillou-May \cite[Remark~1.8 and Definition~6.2]{guillou_may}, the following will be used in the definition of the Burnside 2-category (\cref{definition:Burnside2}) below.
A more general approach to such strictifications can be found in \cite{guillou}.
\begin{explanation}[Choices of Pullbacks in $\FinG$]\label{explanation:SpanFinG}
  Recall the lexicographic ordering\index{lexicographic ordering} of products from \cref{eq:lexord}.
  We use this to determine choices of pullbacks in $\FinG$, as follows.
  Suppose given the following composable pair of spans in $\FinG$,
  \[
    \begin{tikzpicture}[x=20ex,y=10ex]
      \draw[0cell] 
      (0,0) node (x) {X}
      (x)+(1,0) node (y) {Y}
      (y)+(1,0) node (z) {Z,}
      (x)+(.5,.5) node (a) {A}
      (y)+(.5,.5) node (b) {B}
      ;
      \draw[1cell] 
      (a) edge['] node {f} (x)
      (a) edge node {g} (y)
      (b) edge['] node {h} (y)
      (b) edge node {k} (z)
      ;
    \end{tikzpicture}
  \]
  where
  \[
    X = (\ufs{n}_X,\al_X),
    \quad
    A = (\ufs{n}_A,\al_A),
    \quad
    Y = (\ufs{n}_Y,\al_Y),
    \quad
    B = (\ufs{n}_B,\al_B),
    \text{\ and\ \ }
    Z = (\ufs{n}_Z,\al_Z).
  \]
  Let
  \[
    A \times_Y B = \{ (a,b) \in \ufs{n}_A \times \ufs{n}_B \ |\ g(a)=h(b) \}
  \]
  denote the pullback of $G$-sets, with its ordering induced by the lexicographic ordering on $\ufs{n}_A \times \ufs{n}_B$.
  This determines a unique order-preserving isomorphism of finite $G$-sets
  \begin{equation}\label{eq:p-rho-AYB}
    (\ufs{p},\rho) \fto{\iso} A \times_Y B
  \end{equation}
  with $(\ufs{p},\rho) \in \FinG$.

  We write $A \circ B = (\ufs{p},\rho)$ to denote this choice of pullback in $\FinG$ and let $\pi_A$ and $\pi_B$ denote the indicated composites below, where the unlabeled isomorphism is that of \cref{eq:p-rho-AYB}.
  \begin{equation}\label{eq:AcircB}
    \begin{tikzpicture}[x=20ex,y=10ex]
      \draw[0cell] 
      (0,0) node (x) {X}
      (x)+(1,0) node (y) {Y}
      (y)+(1,0) node (z) {Z}
      (x)+(.5,.5) node (a) {A}
      (y)+(.5,.5) node (b) {B}
      (a)+(.5,.3) node (p) {A \times_Y B}
      (a)+(.5,1) node (q) {A \circ B}
      ;
      \draw[1cell] 
      (a) edge['] node {f} (x)
      (a) edge node {g} (y)
      (b) edge['] node {h} (y)
      (b) edge node {k} (z)
      (p) edge node {} (a)
      (p) edge node {} (b)
      (q) edge['] node {\pi_A} (a)
      (q) edge node {\pi_B} (b)
      (q) edge node {\iso} (p)
      ;
    \end{tikzpicture}
  \end{equation} 
  We note three consequences of these choices via lexicographic ordering.
  \begin{enumerate}
  \item These choices for pullbacks make composition in $\Span(\FinG)$ strictly associative.
  \item The morphism $\pi_A$ is always order-preserving.
  \item The morphism $\pi_B$ is generally not order-preserving.
  \end{enumerate}

  For each $Y = (\ufs{n}_Y,\al_Y)$ in $\FinG$, let $\De_Y$ denote the unit 1-cell for $Y$ in $\Span(\FinG)$:
  \[
    \De_Y = \Big( Y \fot{1_Y} Y \fto{1_Y} Y \Big).
  \]
  In \cref{eq:AcircB} above, if the span $(h,k)$ is the unit $\De_Y$, then $B = Y$ and we have
  \[
    A \circ B = A,\qquad \pi_A = 1_A,\andspace \pi_B = g.
  \]
  Thus, $\De_Y$ is a strict right unit.

  Now suppose, instead, that the span $(f,g)$ in \cref{eq:AcircB} is the unit $\De_Y$.
  Then $A = Y$, but $\pi_B = g$ if and only if $h$ is an order-preserving $G$-map.
  In general, $\pi_B$ is an isomorphism of finite $G$-sets determined by the re-ordering of $\ufs{n}_B$ that is induced by the fibers of $h$.
\end{explanation}

To construct a 2-category from $\Span(\FinG)$, the identity 1-cells $\De_X$ are augmented by new strict identities via the following construction.
\begin{definition}[Whiskering a Category]\label{definition:Ddag}
  Suppose given a small category $\D$ with a distinguished object $\De \in \D$.
  Define the \index{whiskering!of a category}\index{category!whiskering}\emph{whiskering at $\De$}, denoted $\D^\dagger$, as a category whose objects consist of those of $\D$, together with a new object $I$ and an isomorphism
  \[
    I \fto[\iso]{\ze_{\De}} \De.
  \]
  The morphisms in $\D^\dagger$ are generated by those of $\D$ and composition with $\ze_\De$ and its inverse.
  Thus, $\D^\dagger$ is the pushout in $\Cat$ of the two inclusions
  \[
    \D \leftarrow \{ \De \} \to \{ \ze_\De^{\pm 1} \}
  \]
  where $\{ \De \}$ denotes the discrete category on $\De$ and the right hand side denotes the category generated by the isomorphism $\ze_\De$ and its inverse.
  A further elaboration of the whiskering construction is given in \cite[Definition~6.1]{guillou_may}.
\end{definition}

\begin{definition}[The Burnside 2-Category]\label{definition:Burnside2}
  The \index{Burnside!2-category}\index{2-category!Burnside}\emph{Burnside 2-category} of a finite group $G$ is a $\permcatsu$-enriched category (\cref{expl:perm-enr-cat}) denoted $\GE$\label{not:GE} and defined as follows.
  Its objects are the finite $G$-sets $X = (\ufs{n},\al)$ of $\FinG$ (\cref{definition:FinG}).
  For each pair of objects
  \[
    X = (\ufs{n},\al)
    \andspace
    Y = (\ufs{p},\be)
    \inspace
    \FinG,
  \]
  the category of 1- and 2-cells is given by
  \begin{equation}\label{eq:GEXY}
    \GE(X,Y) = 
    \begin{cases}
      \Span(\FinG)(X,Y) & \ifspace X \ne Y \orspace |X| \leq 1,\\
      \Span(\FinG)(X,X)^\dagger & \ifspace X = Y \andspace |X| \ge 2
    \end{cases}
  \end{equation}
  where $\Span(\FinG)$ is the bicategory of spans (\cref{definition:Span}) with the lexicographic choice of pullbacks from \cref{explanation:SpanFinG} and $\Span(\FinG)(X,X)^\dagger$ is the whiskering of the category $\Span(\FinG)(X,X)$ as in \cref{definition:Ddag} at the unit 1-cell $\De_{X}$.

  The horizontal composition of $\Span(\FinG)$ extends uniquely to $\GE$ such that the 1-cells $I_{\De_X} \in \Span(\FinG)(X,X)^\dag$ are strictly unital.
  The permutative structure of each $\Span(\FinG)(X,Y)$ given by disjoint union (\cref{definition:FinG}~\cref{it:FinG-disj}) also extends uniquely such that $(\ufs{0},!)$ remains its unit and for $Y \ne (\ufs{0},!)$ we have
  \[
    I_{\De_X} \bincoprod Y = X \bincoprod Y
    \andspace
    Y \bincoprod I_{\De_X} = Y \bincoprod X.
  \]
  For further explanation of this structure, see \cite[Definition~6.2]{guillou_may}, where our $\GE$ is denoted $\GE'$.
\end{definition}

\begin{remark}[Non-Self-Duality of $\GE$]\label{remark:GE-non-self-dual}
  Recall that the Burnside 1-category, $\GB$ in \cref{definition:Burnside1} is self-dual (\cref{lemma:GB-self-dual}).
  However, the assignment that sends a span $(f,g)$ as in \cref{eq:Span1} to its reverse $(g,f)$ does not define a 2-functor
  \[
    \GE \to \GE^\op
  \]
  because it does not preserve composition strictly.
  It is natural to consider the generalization from 2-functors to pseudofunctors, but the latter structure does not provide a $\permcatsu$-enriched functor in the sense of \cref{expl:perm-enr-functors}.
  This subtlety has further implications to be noted in \cref{rk:BO7.5} below.
\end{remark}

The following is a special case of more general enriched Mackey functors introduced in \cref{def:enr-diag-cat}.
\begin{definition}\label{definition:Sp-G-Mackey}
  Suppose given a (possibly non-symmetric) \index{K-theory@$K$-theory}$K$-theory multifunctor
  \[
    K\cn \permcatsu \to \Sp
  \]
  from permutative categories to spectra, and let $(-)_K$ denote the corresponding change of enrichment (\cref{def:mult-change-enr}).
  The category of \emph{spectral $G$-Mackey functors}\index{spectral!Mackey functor}\index{Mackey functor!spectral} for $K$ is the enriched presheaf category
  \[
    \Sp\mh\Cat\big( (\GE_{K})^\op, \Sp \big),
  \]
  consisting of $\Sp$-enriched functors and transformations, as in \cref{mcat-copm}.
\end{definition}

Note that if $K$ is a multifunctor in the symmetric sense---for example, if $K$ is the Elmendorf-Mandell $K$-theory, $\Kem$, in \cref{Kem}---then $(\GE_K)^\op$ and $(\GE^\op)_K$ are equal as $\Sp$-categories by \cref{dF-opposite}.
In such a case, the category of spectral $G$-Mackey functors is equal to $\Sp\mh\Cat\big( (\GE^\op)_K , \Sp \big)$.
However, if $K$ is not symmetric, then there is no such identification.
See, e.g., \cref{thm:Kemdg,rk:BO7.5} for particular uses of these details.

For further development of both theory and applications of spectral Mackey functors in equivariant algebraic $K$-theory, the reader is referred to \cite{bohmann_osorno-mackey,barwick,merling-malkiewich-Athy,barwick-glasman-shah,merling-malkiewich,guillou_may,gmmo}.
The key result for our purposes is the following from Guillou-May \cite{guillou_may}, which is a stable analog of Elmendorf's \cref{theorem:elmendorf}.
Here, $\bbK$\label{not:bbK} denotes the non-symmetric $K$-theory multifunctor in \cite{guillou_may,gmmo}.
\begin{theorem}[{\cite[Theorem~0.1]{guillou_may}}]\label{theorem:GM}
  There is a zigzag of Quillen equivalences
  \[
    \GSp \hty_Q \Sp\mh\Cat\big((\GE_{\bbK})^\op,\Sp\big)
  \]
  where $\GSp$\label{not:GSp} is the category of $G$-spectra.
\end{theorem}

Thus, the Guillou-May theorem shows that the homotopy theory of $G$-spectra is equivalent to that of spectral $G$-Mackey functors\index{Mackey functor!spectral}\index{spectral!Mackey functor} for $\bbK$.

\section{Stable Model Categories as Spectral Presheaf Categories}\label{sec:stmodcat}\index{model category!stable}\index{stable model category}

\begin{definition}\label{definition:model-variants}
  Suppose given a model category $\M$.
  We recall the following terms briefly and refer the reader to \cite{hovey,hirschhorn} for more detailed descriptions. 
  \begin{enumerate}
  \item We say $\M$ is \emph{simplicial}\index{simplicial model category}\index{model category!simplicial} if it is enriched, tensored, and cotensored over simplicial sets, such that the following \emph{pullback powering}\index{pullback powering condition} condition holds.
    For each cofibration $i\cn A \to B$ and fibration $p \cn X \to Y$ in $\M$, the universal morphism induced by $\M(i,X)$ and $\M(B,p)$,
    \[
      \M(B,X) \to \M(A,X) \times_{\M(A,Y)} \M(B,Y),
    \]
    is a Kan fibration that is acyclic whenever either $i$ or $p$ is acyclic.
  \item We say $\M$ is \emph{cofibrantly generated}\index{cofibrantly generated model category}\index{model category!cofibrantly generated} if it is equipped with two sets of morphisms, $\cI$ and $\cJ$, such that the following three statements hold.
    \begin{itemize}
    \item Both $\cI$ and $\cJ$ permit the small object argument.
    \item A morphism of $\M$ is a fibration if and only if it has the right lifting property with respect to every element of $\cJ$
    \item A morphism of $\M$ is an acyclic fibration if and only if it has the right lifting property with respect to every element of $\cI$.
    \end{itemize}
  \item We say that $\M$ is \emph{proper}\index{proper model category}\index{model category!proper} if the following two conditions hold.
    \begin{itemize}
    \item Every pushout of a weak equivalence along a cofibration is a weak equivalence.
    \item Every pullback of a weak equivalence along a fibration is a weak equivalence.
    \end{itemize}
  \item We say that $\M$ is \emph{stable}\index{stable model category}\index{model category!stable} if the suspension and loop functors on its homotopy category are inverse equivalences.\dqed
  \end{enumerate}
\end{definition}

For the remainder of this section we suppose that $\M$ is a simplicial, cofibrantly generated, proper, and stable model category.
The category of symmetric spectra over $\M$ \cite[Definition~3.6.1]{schwede-shipley_stable} is denoted $\Spm$.\label{not:Spm}
The following, from \cite[Definition 3.7.5]{schwede-shipley_stable}, describes an $\Sp$-enriched category generalizing the endomorphism spectrum associated to an object of $\M$.
\begin{definition}\label{definition:EP}
  Suppose $P$ is a set of cofibrant objects in $\M$.
  The \emph{spectral endomorphism category}\index{spectral!endomorphism category} $\EP$\label{not:EP} is the full $\Sp$-subcategory of $\Spm$ with objects given by the fibrant replacements, relative to the level model structure on $\Spm$, of the symmetric suspension spectra of the objects in $P$.
\end{definition}

The following result of Schwede-Shipley gives a characterization of $\M$ via $\Sp$-enriched presheaves.
In this result, $\Spcat\big(\EP^\op, \Sp\big)$ denotes the $\EP$-presheaf category of $\Sp$ as in \cref{mcat-copm}.
\begin{theorem}[{\cite[Theorem~3.3.3]{schwede-shipley_stable}}]\label{theorem:schwede-shipley}
  Suppose $P$ is a set of compact generators of a simplicial, cofibrantly generated, proper, and stable model category $\M$.
  Then there is a chain of simplicial Quillen equivalences
  \[\M \hty_Q \Spcat\big(\EP^\op, \Sp\big).\]
\end{theorem}
The work of Schwede-Shipley goes on to give a number of applications in (derived) Morita theory and equivariant stable homotopy.
In each case, their work characterizes the relevant stable model category as a category of spectral presheaves, also called enriched Mackey functors (see \cref{def:enr-diag-cat}).\index{Mackey functor!enriched}

\mainmatter

\part{Background on Multicategories and \texorpdfstring{$K$}{K}-Theory Functors}
\label{part:background}

\chapter{Categorically Enriched Multicategories}
\label{ch:ptmulticat}
In this chapter we discuss the following four categorically-enriched multicategories that are central to this work:
\begin{itemize}
\item $\Multicat$ of small multicategories (\cref{sec:multicatclosed}),
\item $\pMulticat$ of small pointed multicategories (\cref{sec:ptmulticatclosed}),
\item $\MoneMod$ of left $\Mone$-modules (\cref{sec:monemodules}), and
\item $\permcatsu$ of small permutative categories (\cref{sec:multpermcat}).
\end{itemize}
Each of the $\Cat$-multicategories,
\[\Multicat, \quad \pMulticat, \andspace \MoneMod,\]
is induced by a corresponding symmetric monoidal $\Cat$-category structure.  On the other hand, the $\Cat$-multicategory structure on $\permcatsu$ is not induced by the monoidal structure on $\MoneMod$ or $\pMulticat$.   The next table summarizes the various structures of these categories.
\smallskip
\begin{center}
\resizebox{.8\textwidth}{!}{
{\renewcommand{\arraystretch}{1.35}%
{\setlength{\tabcolsep}{1ex}
\begin{tabular}{|c|c|c|c|c|}\hline
& $\Multicat$ & $\pMulticat$ & $\MoneMod$ & $\permcatsu$ \\ \hline
2-category & \ref{v-multicat-2cat} & \ref{thm:pmulticat} & \ref{definition:MoneMod-prelim} & \ref{def:permcat} \\ \hline
symmetric monoidal $\Cat$-category & \ref{theorem:multicat-symmon} & \ref{thm:pmulticat-smclosed} & \ref{definition:MoneMod} & --- \\ \hline
$\Cat$-multicategory & \ref{expl:multicatcatmulticat} & \ref{expl:ptmulticatcatmulticat} & \ref{expl:monemodcatmulticat} & \ref{thm:permcatmulticat} \\ \hline
\end{tabular}}}}
\end{center}
\smallskip

These four $\Cat$-multicategories are related by several $\Cat$-multifunctors
\begin{equation}\label{EndU-intro}
\begin{tikzpicture}[xscale=3.5,yscale=1.5,vcenter]
\draw[0cell=.85]
(0,0) node (a) {\permcatsu}
(a)++(1,0) node (b) {\Multicat}
(a)++(0,-1) node (c) {\MoneMod}
(c)++(1,0) node (d) {\pMulticat}
;
\draw[1cell=.9]  
(a) edge node {\End} (b)
(a) edge node {\Endst} (d)
(a) edge node[swap] {\Endm} (c)
(c) edge node[pos=.4] {\Um} (d)
(d) edge node[swap] {\Ust} (b)
;
\end{tikzpicture}
\end{equation}
that we will explain in \cref{endfactorization}.  Here is a summary table.
\smallskip
\begin{center}
\resizebox{.8\width}{!}{
{\renewcommand{\arraystretch}{1.4}%
{\setlength{\tabcolsep}{1ex}
\begin{tabular}{|c|c|c|c|c|c|}\hline
& $\End$ & $\Endst$ & $\Endm$ & $\Ust$ & $\Um$ \\ \hline
2-functor & \ref{endtwofunctor} & \ref{endsttwofunctor} & \ref{endmtwofunctor} & \ref{Usttwofunctor} & \ref{proposition:EM2-5-1} \eqref{it:EM251-uniqueness} and \eqref{it:EM251-3} \\ \hline
symmetric monoidal $\Cat$-functor & --- & --- & --- & \ref{expl:Ust} & \ref{expl:moneinclusion} \\ \hline
$\Cat$-multifunctor & \ref{expl:end-catmulti} & \ref{expl:endst-catmulti} & \ref{expl:endm-catmulti} & \ref{expl:Ust-catmulti} & \ref{expl:Um-catmulti} \\ \hline
\end{tabular}}}}
\end{center}
\smallskip

\subsection*{Connection with Other Chapters}\

\subsubsection*{Infinite Loop Space Machines}

As summarized in \cref{eq:Ksummary}, the $\Cat$-multicategory $\permcatsu$ is connected to several categories in Segal $K$-theory and Elmendorf-Mandell $K$-theory via enriched multifunctors, including $\Endm$.

\subsubsection*{Equivalences of Homotopy Theories}

In \cref{ch:multperm} we discuss the fact that the endomorphism multicategory construction $\End$ in \cref{EndU-intro} is an equivalence of homotopy theories.   Moreover, in \cref{ch:ptmulticat-sp,ch:ptmulticat-alg} we extend this observation to the pointed setting by showing that each of $\Endst$, $\Endm$, and $\Um$ is an equivalence of homotopy theories.  In \cref{part:homotopy-mackey} we further extend these equivalences of homotopy theories to the respective categories of enriched diagrams and enriched Mackey functors.  See \cref{mackey-xiv-pmulticat,mackey-xiv-mone,mackey-pmulti-mone}.

\subsection*{Background}

Definitions about permutative categories, enriched multicategories, and pointed multicategories are in \cref{sec:monoidalcat,ch:prelim_multicat}.  Definitions for 2-categories and enriched categories are in \cref{sec:twocategories,sec:enrichedcat}.  Symmetric monoidal enriched categories are discussed in \cref{sec:enrmonoidalcat,sec:selfenr-smc,sec:change-enrichment}.  

\subsection*{Chapter Summary}

The following table lists the main content in this chapter.
\begin{center}
\resizebox{\textwidth}{!}{
{\renewcommand{\arraystretch}{1.4}%
{\setlength{\tabcolsep}{1ex}
\begin{tabular}{|c|c|}\hline
\multicolumn{2}{|c|}{\secnumname{sec:multicatclosed}} \\ \hline
multicategories as monadic algebras & \ref{thm:multigraphcat-monadic} \\ \hline
Boardman-Vogt tensor product & \ref{definition:multicat-tensor} and \ref{theorem:multicat-symmon} \\ \hline
$\Multicat$ as a $\Cat$-multicategory & \ref{expl:multicatcatmulticat} \\ \hline
internal hom multicategories & \ref{definition:multicat-hom} and \ref{theorem:multicat-sm-closed} \\ \hline\hline
\multicolumn{2}{|c|}{\secnumname{sec:ptmulticatclosed}} \\ \hline
smash product and pointed internal hom & \ref{def:multicat-wedge-smash} \\ \hline
pointed transformations & \ref{ex:pointedhommulticat} \\
\hline
\end{tabular}}}}
\end{center}
\medskip
The material in this chapter is adapted from \cite[Chapters 5, 6, 8, and 10]{cerberusIII}, which has all the detailed proofs.  We remind the reader of \cref{conv:universe} about universes and \cref{expl:leftbracketing} about left normalized bracketing for iterated products.

\section{Multicategories}
\label{sec:multicatclosed}

There is a 2-category $\Multicat$ of small multicategories, multifunctors, and multinatural transformations (\cref{v-multicat-2cat}).  In this section we review
\begin{enumerate}
\item the symmetric monoidal closed structure (\cref{theorem:multicat-symmon,theorem:multicat-sm-closed}) and
\item the $\Cat$-multicategory structure (\cref{expl:multicatcatmulticat})
\end{enumerate}
on $\Multicat$. 
\begin{itemize}
\item The monoidal product is in \cref{definition:multicat-tensor} after some preliminary constructions.  This monoidal product is often called the \emph{Boardman-Vogt tensor product} in the literature because of its origin in \cite{boardman-vogt}. 
\item The closed structure is given by the internal hom multicategory in \cref{definition:multicat-hom}.
\end{itemize} 

The material in this section is adapted from \cite[Chapters 5 and 6]{cerberusIII}.

\subsection*{Multicategories as Monadic Algebras}

The tensor product on small multicategories requires some preliminary constructions, which we recall first.  The first fact we need is that small multicategories are algebras over a monad.  We use the following notation for input profiles and output.  Recall the class of profiles $\Prof$ (\cref{def:profile}).

\begin{definition}\label{definition:xiotimes}
Suppose $C$ and $D$ are two classes.  Given profiles 
\[\ang{c} = \ang{c_i}_{i=1}^m \in \Prof(C) \andspace \ang{d} = \ang{d_j}_{j=1}^n \in \Prof(D),\] 
we define the following in $\Prof(C \times D)$:\label{not:angtensor}
\begin{align*}
\ang{c} \times d_j & = \ang{(c_i,d_j)}_{i=1}^m\\
c_i \times \ang{d} & = \ang{(c_i,d_j)}_{j=1}^n\\
\ang{c} \otimes \ang{d} & = \ang{\ang{(c_i,d_j)}_{i=1}^m}_{j=1}^n\\
\ang{c} \otimes^\transp \ang{d} & = \ang{\ang{(c_i,d_j)}_{j=1}^n}_{i=1}^m
\end{align*}
Denote by 
\begin{equation}\label{xitimesmn}
\xitimes = \xitimes_{m,n} \cn \ang{c} \otimes \ang{d} \fto{\iso} \ang{c} \otimes^\transp \ang{d}
\end{equation}
the \index{transpose permutation}\index{permutation!transpose}\emph{transpose permutation} induced by changing order of indexing.
\end{definition}

\begin{definition}\label{definition:multigraph}  
A \index{multigraph}\index{graph!multi-}\emph{multigraph} $X$ consists of
\begin{itemize}
\item a class $\Vt X$ of \emph{vertices} and
\item a set $X\mmap{x';\ang{x}}$ for each tuple of vertices $\ang{x}$ and $x'$.
\end{itemize}  
We refer to the elements of $X\mmap{x';\ang{x}}$ as \index{multiedges}\emph{multiedges}, with \index{source}\emph{source} $\ang{x}$ and \index{target}\emph{target} $x'$.  We let $\Prof(X)$ denote $\Prof(\Vt X)$.

A \index{multigraph!morphism}\index{morphism!multigraph}\emph{morphism} of multigraphs
\[f \cn X \to Y\]
consists of
\begin{itemize}
\item a function 
\[f\cn\Vt X \to \Vt Y\]
on vertices and
\item a function 
\[f\cn X\mmap{x';\ang{x}} \to Y\mmap{f(x');f\ang{x}}\]
on multiedges for each $\mmap{x';\ang{x}} \in \Prof(X) \times \Vt X$, with $f\ang{x} = \ang{fx_j}_{j=1}^n$ if $\angx = \ang{x_j}_{j=1}^n$.
\end{itemize} 
Moreover, we define the following.
\begin{itemize}
\item A multigraph is \index{small!multigraph}\index{multigraph!small}\emph{small} if its class of vertices is a set.
\item The collection of small multigraphs and their morphisms form a category, denoted \label{not:mgraph}$\MGraph$.
\end{itemize}
This finishes the definition.
\end{definition}

The following result combines \cite[5.5.9 and 5.5.11]{cerberusIII}.

\begin{theorem}\label{thm:multigraphcat-monadic}
There is a free-forgetful adjunction
\[\begin{tikzpicture}[x=30mm,y=25mm]
    \draw[0cell] 
    (0,0) node (x) {\MGraph}
    (1,0) node (y) {\Multicat}
    ;
    \draw[1cell=.9]
    (x) edge[bend left=11,transform canvas={yshift=0mm}] node {L} (y) 
    (y) edge[bend left=11,transform canvas={yshift=0mm}] node {U} (x) 
    ;
    \draw[2cell] 
    (.5,.01mm) node {\bot}
    ;
\end{tikzpicture}\]
that is strictly monadic.  
\end{theorem}

\subsection*{Two Auxiliary Products}

\begin{definition}\label{definition:amp-product}\index{multigraph!internal product}
For multigraphs $X$ and $Y$ with vertex classes $C$ and $D$, respectively, we define a multigraph $X \amtimes Y$\label{not:amtimes} with vertex class $C \times D$ as follows.  Given
\[\ang{c,d} = \bang{(c_j,d_j)}_{j=1}^n \in \Prof(C \times D) 
\andspace (c',d') \in C \times D,\]
the set of multiedges with source $\ang{c,d}$ and target $(c',d')$ is given by the coproduct
\begin{equation}\label{eq:XamtimesY-coprod}
(X \amtimes Y)\mmap{(c',d');\ang{c,d}} =
\coprod_{\ang{c''} \otimes \ang{d''} \,=\, \ang{c,d}} 
X\mmap{c';\ang{c''}} \times Y\mmap{d';\ang{d''}}.
\end{equation}
The coproduct is indexed by pairs 
\[\big(\ang{c''}, \ang{d''}\big) \in \Prof(C) \times \Prof(D) \stspace \ang{c''} \otimes \ang{d''} = \ang{c,d}\]
with the tensor product of profiles in \cref{definition:xiotimes}.
\end{definition}

\begin{definition}\label{definition:sharp-prod}
For small multicategories $\M$ and $\N$, we define the \index{multicategory!sharp product}\index{product!sharp}\emph{sharp product} $\M\shtimes\N$ as the pushout in $\Multicat$
\begin{equation}\label{sharp-pushout}
\begin{tikzpicture}[x=30mm,y=15mm,vcenter]
  \draw[0cell=.9]
  (0,0) node (a) {\Ob\M \times \Ob\N}
  (1,0) node (b) {\txcoprod_{d \in \Ob\N}\, \M}
  (0,-1) node (c) {\txcoprod_{c \in \Ob\M}\, \N}
  (1,-1) node (d) {\M \shtimes \N}
  ;
  \draw[1cell] 
  (a) edge node {} (b)
  (c) edge node {} (d)
  (a) edge node {} (c)
  (b) edge node {} (d)
  ;
\end{tikzpicture}
\end{equation}
along morphisms induced by the inclusions  \[\Ob\M \hookrightarrow \M \andspace \Ob\N \hookrightarrow \N.\dqed\]
\end{definition}

\begin{explanation}[Sharp Product]\label{explanation:unpacking-sharp}
Restricting \cref{definition:sharp-prod} to objects, there is a canonical bijection
\[\Ob(\M\shtimes\N) \iso \Ob\M \times \Ob\N.\]
The operations of $\M\shtimes\N$ are generated by operations of the form
\[\phi \times d \in \M \times \{d\} \andspace c \times \psi \in \{c\}\times \N\]
subject to the axioms \cref{sharp-i,sharp-ii,sharp-iii,sharp-iv,sharp-v} below, which are determined by the pushout \cref{sharp-pushout}.
\begin{romenumerate}
\item\label{sharp-i} For $(c,d) \in \M\shtimes\N$, there are equalities
\[\opu_c \times d = \opu_{(c,d)} = c \times \opu_d.\]
\item\label{sharp-ii} For operations $\phi, \phi_1, \ldots, \phi_n$  in $\M$ such that the composite below is defined, there is an equality
\[\ga\scmap{\phi \times d; \ang{\phi_j \times d}_{j=1}^n} =
\ga\scmap{\phi; \ang{\phi_j}_{j=1}^n} \times d.\]
\item\label{sharp-iii} For $\si \in \Si_n$, there is an equality
\[(\phi \times d)\cdot\si = (\phi\cdot\si) \times d.\]
\item\label{sharp-iv} For operations $\psi, \psi_1, \ldots, \psi_m$ in $\N$ such that the composite below is defined, there is an equality
\[\ga\scmap{c \times \psi; \ang{c \times \psi_i}_{i=1}^m} =
c \times \ga\scmap{\psi; \ang{\psi_i}_{i=1}^m}.\]
\item\label{sharp-v} For $\si \in \Si_m$, there is an equality
\[(c \times \psi)\cdot\si = c \times (\psi\cdot\si).\]
\end{romenumerate}
These conditions are equivalent to the requirement that a multifunctor
\[F\cn\M\shtimes\N \to \P\]
consists of an assignment on objects,
\[F(c,d) \in \Ob\P \forspace (c,d) \in \Ob\M \times \Ob\N,\]
such that each of
\[F(c,-)\cn\N \to \P \andspace F(-,d)\cn\M \to \P\]
is a multifunctor.
\end{explanation}

\subsection*{The Boardman-Vogt Tensor Product of Multicategories}

\begin{definition}\label{definition:tensor-operations}
Suppose given small multicategories $\M$ and $\N$ along with operations
\[\phi \in \M\mmap{c';\ang{c}} \andspace \psi \in \N\mmap{d';\ang{d}}.\]
We define the following:
\begin{align*}
\phi \times \ang{d} & = \ang{\phi\times d_j}_j \in \txprod_j\,
\M\mmap{c';\ang{c}} \times \{d_j\} \\
\ang{c} \times \psi & = \ang{c_i \times \psi}_i \in \txprod_i\, \{c_i\}
\times \N\mmap{d';\ang{d}}\\
\phi \otimes \psi & = \ga\scmap{c' \times \psi; \phi \times \ang{d}}
\in (\M\shtimes\N)\mmap{(c',d');\ang{c}\otimes\ang{d}}\\
\phi \otimes^\transp \psi & = \ga\scmap{\phi \times d'; \ang{c} \times
\psi} \in (\M \shtimes \N)\mmap{(c',d');\ang{c} \otimes^\transp \ang{d}}
\end{align*}
Denote by $\xitimes$ the bijection
\begin{equation}\label{eq:xitimes}
  (\M\shtimes\N)\mmap{(c',d');\ang{c} \otimes^\transp \ang{d}} \fto{\iso}
  (\M\shtimes\N)\mmap{(c',d');\ang{c} \otimes \ang{d}}
\end{equation} 
induced by the transpose permutation $\xitimes$ in \cref{xitimesmn} that interchanges order of indexing.
\end{definition}

\begin{definition}[Tensor Product of Multicategories]\label{definition:multicat-tensor}
For small multicategories $\M$ and $\N$, the tensor products of \cref{definition:tensor-operations} give two canonical morphisms of multigraphs
\[\begin{tikzpicture}[x=40mm,y=20mm,baseline={(ux.base)}]
    \draw[0cell] 
    (0,0) node (ux) {(U\M) \amtimes (U\N)}
    (1,0) node (uy) {U(\M \shtimes \N).}
    ;
    \draw[1cell] 
    (ux) edge[transform canvas={yshift=.8mm}] node {\otimes} (uy)
    (ux) edge[transform canvas={yshift=-.8mm}] node['] {\xitimes
      \circ \otimes^\transp} (uy)
    ;
\end{tikzpicture}\]
Taking adjoints, we obtain the two morphisms in $\Multicat$ below.  We define $\M \otimes \N$ to be their \index{multicategory!Boardman-Vogt tensor product}\index{multicategory!tensor product}\index{product!Boardman-Vogt tensor}\index{Boardman-Vogt tensor product}coequalizer in $\Multicat$:
\begin{equation}\label{bvtensor}
\begin{tikzpicture}[x=30mm,y=20mm,vcenter]
    \draw[0cell] 
    (0,0) node (ux) {L\big((U\M) \amtimes (U\N)\big)}
    (1,0) node (uy) {\M \shtimes \N}
    (1.8,0) node (w) {\M \otimes \N}
    ;
    \draw[1cell] 
    (ux) edge[transform canvas={yshift=.8mm}] node {} (uy)
    (ux) edge[transform canvas={yshift=-.8mm}] node['] {} (uy)
    (uy) edge[dashed] node {} (w)
    ;
\end{tikzpicture}
\end{equation}
For an object $(c,d) \in \M\shtimes\N$, we denote by $c \otimes d$ its image in $\M\otimes\N$.  Moreover, the tensor product $\otimes$ extends naturally to multifunctors.
\end{definition}

\begin{explanation}[Unpacking the Tensor Product]\label{explanation:unpacking-tensor}
Restricting \cref{definition:multicat-tensor} to objects, there are canonical bijections
\begin{equation}\label{ObMtensorN}
\Ob(\M \otimes \N) \iso \Ob(\M\shtimes\N) \iso \Ob\M \times \Ob\N.
\end{equation}
The operations of $\M\otimes\N$ are generated by
\[\phi \otimes d \in \M\mmap{c';\ang{c}}\times\{d\} \andspace
c \otimes \psi \in \{c\} \times \N\mmap{d';\ang{d}}\]
subject to the relations of $\M\shtimes\N$ in \cref{explanation:unpacking-sharp} along with one additional
\index{relation!interchange}\index{interchange relation}\emph{interchange relation}
\begin{equation}\label{eq:interchange-relation}
\phi \otimes \psi = (\phi \otimes^\transp \psi)\cdot\xitimes.
\end{equation}
A multifunctor
\[F \cn \M \otimes \N \to \P\]
consists of an assignment on objects, 
\[F(c,d) \in \Ob\P \forspace (c,d) \in \Ob\M \times \Ob\N,\]
such that the following two conditions hold.
\begin{itemize}
\item Each of 
\[F(c,-)\cn\N \to \P \andspace F(-,d)\cn\M \to \P\]
is a multifunctor.
\item There is an equality
\begin{equation}\label{eq:F-interchange-relation}
F(\phi \otimes \psi) = F(\phi \otimes^\transp \psi)\cdot\xitimes
\end{equation}
for each $\phi \in \M\mmap{c';\ang{c}}$ and $\psi \in \N\mmap{d';\ang{d}}$.\defmark
\end{itemize} 
\end{explanation}

\begin{definition}[Braiding on Multicategories]\label{definition:sharp-symm}
For small multicategories $\M$ and $\N$, suppose
\[\beta\cn\M\shtimes\N \to \N\shtimes\M\]
is the multifunctor given
\begin{itemize}
\item on objects by $\beta(c,d) = (d,c)$ and
\item on generating operations by
\[\beta(\phi\times d) = d \times \phi \andspace \beta(c \times \psi) = \psi \times c.\]
\end{itemize} 
Define the \emph{braiding}
\[\beta\cn\M\otimes\N \to \N\otimes\M\]
as the induced multifunctor on tensor products.
\end{definition}

Recall the following.
\begin{itemize}
\item A symmetric monoidal $\V$-category (\cref{definition:symm-monoidal-vcat}) is a symmetric monoidal category in the $\V$-enriched sense.  \cref{theorem:multicat-symmon} below involves the case $\V = (\Cat, \times, \boldone)$.
\item There is a 2-category $\Multicat$, which is \cref{v-multicat-2cat} with $\V = \Set$.
\item The initial operad $\Mtu$ in \cref{ex:vmulticatinitialterminal} \cref{ex:initialoperad} has a single object $*$ and a single unit operation $\opu_* \in \Mtu\smscmap{*;*}$.  
\end{itemize}  
The following result combines \cite[5.6.18 and 6.4.3]{cerberusIII}.

\begin{theorem}\label{theorem:multicat-symmon}
The following quadruple is a symmetric monoidal category, with the associativity and unit isomorphisms induced by those of the sharp product, $\shtimes$.  
\[\big(\Multicat, \otimes, \Mtu, \beta\big)\]
Moreover, the tensor product $\otimes$ extends componentwise to multinatural transformations such that the quadruple above becomes a symmetric monoidal $\Cat$-category.
\end{theorem}

\begin{explanation}[$\Multicat$ as a $\Cat$-Multicategory]\label{expl:multicatcatmulticat}
Since $\Multicat$ is a symmetric monoidal $\Cat$-category, it has the structure of a $\Cat$-multicategory by \cref{proposition:monoidal-v-cat-v-multicat}, with the following data.
\begin{itemize}
\item Its objects are small multicategories.
\item For small multicategories $\ang{\M_j}_{j=1}^n$ and $\N$, the $n$-ary multimorphism category is
\begin{equation}\label{multicathomcats}
\Multicat\scmap{\ang{\M_j}_{j=1}^n ; \N} = 
\begin{cases}
\Multicat\left(\txotimes_{j=1}^n \M_j \scs \N\right) & \text{if $n>0$ and}\\
\Multicat(\Mtu,\N) & \text{if $n=0$}.
\end{cases}
\end{equation}
If $n>0$, then this category has
\begin{itemize}
\item multifunctors 
\[\txotimes_{j=1}^n \M_j \to \N\]
as objects and
\item multinatural transformations between such multifunctors as morphisms.
\end{itemize}
If $n=0$, then the objects in 
\[\Multicat\scmap{\ang{};\N} = \Multicat\left(\Mtu, \N\right)\] 
are multifunctors 
\[\Mtu \to \N.\]  
Each such multifunctor is determined by a choice of an object in $\N$.  Thus, $\Multicat\scmap{\ang{};\N}$ is canonically isomorphic to the underlying category of $\N$ as in \cref{ex:unarycategory}.
\item The symmetric group action is induced by the braiding $\beta$ of the tensor product (\cref{definition:sharp-symm}). 
\item The multicategorical composition is given by tensor product and composition of multifunctors, and likewise for multinatural transformations.
\end{itemize}
This finishes the description of the $\Cat$-multicategory $\Multicat$.
\end{explanation}

\subsection*{Internal Hom for Multicategories}

We use the following notation for a tuple of multifunctors $\ang{F}$.

\begin{definition}\label{definition:multicat-hom-notation}
Suppose given multicategories $\M$ and $\N$ together with a tuple of multifunctors $\ang{F} = \ang{F_i}_{i=1}^m$ with each $F_i \cn \M \to \N$.  Then we use the following notation.
\begin{itemize}
\item For $c\in \Ob\M$, denote by\label{not:angFc}
\[\ang{F}c = \ang{F_ic}_{i=1}^m.\]
\item For $\ang{c} = \ang{c_j}_{j=1}^n \in \Prof(\Ob\M)$, denote by\label{not:angofFc}
\[\ang{Fc} = \ang{\ang{F_i c_j}_{i=1}^m}_{j=1}^n \andspace 
\ang{Fc}^{\transp} = \ang{\ang{F_i c_j}_{j=1}^n}_{i=1}^m.\]
\item For an $n$-ary operation $\phi \in \M\mmap{c';\ang{c}}$, denote by
\[\ang{F}\phi = \ang{F_i \phi}_{i=1}^m \in \txprod_{i=1}^m \N\mmap{F_i c'; F_i\ang{c}}.\]
\end{itemize}
This finishes the definition.
\end{definition}

\begin{definition}\label{definition:multicat-hom}
For small multicategories $\M$ and $\N$, the \index{hom!internal - multicategory}\index{multicategory!internal hom}\index{internal hom!multicategory}\emph{internal hom multicategory} 
\[\Hom(\M,\N)\]
is defined as follows.  
\begin{itemize}
\item The objects of $\Hom(\M,\N)$ are multifunctors $\M \to \N$. 
\item The $m$-ary operations
\begin{equation}\label{transformation-theta}
\theta \cn\ang{F} = \ang{F_i}_{i=1}^m \to G
\end{equation}
in $\Hom(\M,\N)$ are called \index{transformation!internal hom multicategory}\emph{transformations} and are given by component $m$-ary operations
\[\theta_c\in \N\mmap{Gc;\ang{F}c} \forspace c \in \ObM.\]
For each operation $\phi \in \M\mmap{c';\ang{c}}$ with $\ang{c} = \ang{c_j}_{j=1}^n$ and $\theta_{\ang{c}} = \ang{\theta_{c_j}}_{j=1}^n$, the following \index{naturality condition!internal hom multicategory}\emph{naturality condition} is required to hold:
\begin{equation}\label{eq:multinat-equality}
\ga\scmap{G\phi; \theta_{\ang{c}}} = \ga\scmap{\theta_{c'} ; \ang{F}\phi)} \cdot\xitimes_{m,n}.
\end{equation}
\item The unit operation 
\[1_G \cn G \to G\]
is given by the identity multinatural transformation whose component at $c$ is the $(Gc)$-colored unit $1_{Gc}$ in $\N$.
\item The composition and symmetric group action of $\Hom(\M,\N)$ are given componentwise by those of $\N$.
\end{itemize}
This finishes the definition of the internal hom multicategory $\Hom(\M,\N)$.  Moreover, $\Hom(-,-)$ extends naturally to multifunctors, contravariantly in the first argument and covariantly in the second argument.
\end{definition}

Recall from \cref{def:closedcat} that a symmetric monoidal category is \emph{closed} if, for each object $x$, the functor $- \otimes x$ admits a right adjoint.  The following is \cite[5.7.14]{cerberusIII}.

\begin{theorem}\label{theorem:multicat-sm-closed}\index{Boardman-Vogt tensor product!symmetric monoidal closed}\index{multicategory!symmetric monoidal closed}
Equipped with the internal hom of \cref{definition:multicat-hom}, the symmetric monoidal category $\Multicat$ in \cref{theorem:multicat-symmon} is closed.
\end{theorem}

\section{Pointed Multicategories}
\label{sec:ptmulticatclosed}

There is a 2-category $\pMulticat$ of small pointed multicategories, pointed multifunctors, and pointed multinatural transformations (\cref{sec:ptmulticatiicat}).  In this section we review
\begin{enumerate}
\item the symmetric monoidal closed structure (\cref{thm:pmulticat-smclosed}) and 
\item the $\Cat$-multicategory structure (\cref{expl:ptmulticatcatmulticat})
\end{enumerate} 
on $\pMulticat$.  The material in this section is adapted from \cite[Chapters 5 and 6]{cerberusIII}. 

\begin{definition}\label{def:multicat-wedge-smash}
Suppose $(\M,i^\M)$ and $(\N,i^\N)$ are small pointed multicategories.  We define the following small pointed multicategories.
\begin{enumerate}
\item The \index{pointed!multicategory!wedge}\index{multicategory!pointed!wedge}\index{product!wedge}\index{wedge product}\emph{wedge product}, also called the \index{wedge sum}\index{sum!wedge}\emph{wedge sum}, is the pointed multicategory defined by the following coequalizer in $\Multicat$.
\begin{equation}\label{multicat-wedge}
\begin{tikzpicture}[x=30mm,y=20mm,vcenter]
    \draw[0cell] 
    (0,0) node (ux) {\Mterm}
    (1,0) node (uy) {\M \bincoprod \N}
    (1.8,0) node (w) {\M \wed \N}
    ;
    \draw[1cell] 
    (ux) edge[transform canvas={yshift=.8mm}] node {i^\M} (uy)
    (ux) edge[transform canvas={yshift=-.8mm}] node['] {i^\N} (uy)
    (uy) edge[dashed] node {} (w)
    ;
\end{tikzpicture}
\end{equation}
\item The \index{smash product!pointed multicategories}\index{multicategory!pointed!smash product}\index{pointed!multicategory!smash product}\emph{smash product} is the pointed multicategory defined by the following pushout in $\Multicat$.
\begin{equation}\label{eq:multicat-smash-pushout}
\begin{tikzpicture}[x=30mm,y=15mm,vcenter]
  \draw[0cell] 
  (0,0) node (a) {(\M\otimes\Mterm) \bincoprod (\Mterm \otimes \N)}
  (1,0) node (b) {\M \otimes \N}
  (0,-1) node (c) {\Mterm}
  (1,-1) node (d) {\M \sma \N}
  ;
  \draw[1cell] 
  (a) edge node {} (b)
  (c) edge node {i^{\M \sma \N}} (d)
  (a) edge node {} (c)
  (b) edge node {\vpi_{\M,\N}} (d)
  ;
\end{tikzpicture}
\end{equation}
In \cref{eq:multicat-smash-pushout}, $\otimes$ is the tensor product in \cref{bvtensor}, and $\Mterm$ is the terminal multicategory in \cref{definition:terminal-operad-comm}.  Moreover, the smash product extends naturally to pointed multifunctors.
\item The \index{pointed!smash unit!pointed multicategory}\index{smash unit!pointed multicategory}\index{pointed!multicategory!smash unit}\index{multicategory!pointed!smash unit}\emph{smash unit} $\Mtup$ is the pointed multicategory
\begin{equation}\label{eq:smashunit}
\Mtup = \Mtu \bincoprod \Mterm
\end{equation}
with $\Mtu$ the initial operad in \cref{ex:vmulticatinitialterminal} \cref{ex:initialoperad}.  The pointed structure $\Mterm \to \Mtup$ is given by the $\Mterm$ summand in $\Mtup$.
\item The \index{hom!pointed internal - multicategory}\index{multicategory!pointed!hom}\index{pointed!multicategory!hom}\index{internal hom!pointed multicategory}\emph{pointed internal hom} is the pointed multicategory defined by the following pullback in $\Multicat$, where $\Hom(-,-)$ is the internal hom multicategory in \cref{definition:multicat-hom}.
\begin{equation}\label{eq:multicat-pHom}
\begin{tikzpicture}[x=30mm,y=15mm,vcenter]
    \draw[0cell] 
    (0,0) node (a) {\pHom(\M,\N)}
    (1,0) node (b) {\Mterm}
    (0,-1) node (c) {\Hom(\M,\N)}
    (1,-1) node (d) {\Hom(\Mterm,\N)}
    ;
    \draw[1cell] 
    (a) edge node {} (b)
    (c) edge node {} (d)
    (a) edge node {} (c)
    (b) edge node {} (d)
    ;
\end{tikzpicture}
\end{equation}
The pointed structure $\Mterm \to \pHom(\M,\N)$ is induced by the composite
\[\Mterm \iso \Hom(\M,\Mterm) \fto{(i^\N)_*} \Hom(\M,\N) \fto{(i^\M)^*} \Hom(\Mterm,\N),\]
which is equal to the right vertical morphism in \cref{eq:multicat-pHom}.  Moreover, $\pHom(-,-)$ extends naturally to pointed multifunctors, contravariantly in the first argument and covariantly in the second argument.
\end{enumerate}
This finishes the definition.
\end{definition}

\begin{explanation}[Pointed Internal Hom Multicategory]\label{ex:pointedhommulticat}
Unpacking the pullback \cref{eq:multicat-pHom}, we describe the pointed multicategory $\pHom(\M,\N)$ explicitly as follows.
\begin{itemize}
\item Its objects are pointed multifunctors $(\M,i^\M) \to (\N,i^\N)$.
\item For pointed multifunctors 
\[\ang{F} =  \ang{F_i}_{i=1}^m \scs G \cn \M \to \N,\] 
an $m$-ary operation in $\pHom(\M,\N)\scmap{\angF;G}$ is a transformation
\[\theta = \{\theta_c\}_{c \in \ObM} \cn\ang{F} \to G\]
as in \cref{transformation-theta} such that the component at the basepoint object $* \in \ObM$,
\begin{equation}\label{ptdtransformation}
\theta_* = \iota_m \in \N\scmap{\angF *;G(*)} = \N\scmap{\ang{*}_{i=1}^m;*},
\end{equation}
is equal to the $m$-ary basepoint operation in $\N$.
\end{itemize}
A transformation that satisfies the basepoint condition \cref{ptdtransformation} is called a \index{pointed!transformation}\index{transformation!pointed}\emph{pointed transformation}.
\end{explanation}

The following result combines \cite[5.7.22 and 6.4.4]{cerberusIII} and is the pointed variant of \cref{theorem:multicat-symmon,theorem:multicat-sm-closed}.

\begin{theorem}\label{thm:pmulticat-smclosed}
In the context of \cref{def:multicat-wedge-smash}, the quadruple
\[\big(\pMulticat, \sma, \Mtup, \pHom\big)\]
is a complete and cocomplete symmetric monoidal closed category, with the associativity and unit isomorphisms induced by those of $\Multicat$ in \cref{theorem:multicat-symmon}.

Moreover, the smash product $\sma$ extends componentwise to pointed multinatural transformations such that $\pMulticat$ becomes a symmetric monoidal $\Cat$-category.
\end{theorem}

Next is the pointed variant of \cref{expl:multicatcatmulticat}.

\begin{explanation}[$\pMulticat$ as a $\Cat$-Multicategory]\label{expl:ptmulticatcatmulticat}
Since $\pMulticat$ is a symmetric monoidal $\Cat$-category, it has the structure of a $\Cat$-multicategory by \cref{proposition:monoidal-v-cat-v-multicat}, with the following data.
\begin{itemize}
\item Its objects are small pointed multicategories.
\item For small pointed multicategories $\bang{(\M_j, i^{\M_j})}_{j=1}^n$ and $(\N,i^\N)$, the $n$-ary multimorphism category
\[\pMulticat\lrscmap{\bang{(\M_j, i^{\M_j})}_{j=1}^n ; (\N,i^\N)} 
= \pMulticat\brb{\txsma_{j=1}^n \M_j , \N}\]
has
\begin{itemize}
\item pointed multifunctors 
\[\txsma_{j=1}^n \M_j \to \N\]
as objects and
\item pointed multinatural transformations between such pointed multifunctors as morphisms.
\end{itemize}
If $n=0$, then the objects in 
\[\pMulticat\scmap{\ang{};(\N,i^\N)} = \pMulticat(\Mtup, \N)\]
are pointed multifunctors 
\[\Mtup = \Mtu \bincoprod \Mterm \to (\N,i^\N).\]  
Each such pointed multifunctor is determined by a choice of an object in $\N$.  Thus, $\pMulticat\scmap{\ang{};(\N,i^\N)}$ is canonically isomorphic to the underlying category of $\N$ as in \cref{ex:unarycategory}.
\item The symmetric group action is induced by the braiding of the smash product, which, in turn, is induced by the braiding of the tensor product (\cref{definition:sharp-symm}). 
\item The multicategorical composition is given by smash product and composition of pointed multifunctors, and likewise for pointed multinatural transformations.
\end{itemize}
This finishes the description of the $\Cat$-multicategory $\pMulticat$.
\end{explanation}

\begin{explanation}[Forgetting Basepoints]\label{expl:Ust}
The forgetful 2-functor\index{pointed!multicategory!forgetful functor}
\begin{equation}\label{Ust}
\Ust \cn \big(\pMulticat, \sma, \Mtup\big) \to \big(\Multicat, \otimes, \Mtu\big)
\end{equation}
in \cref{Usttwofunctor} is a symmetric monoidal $\Cat$-functor (\cref{definition:braided-monoidal-vfunctor}) with the following structure morphisms:
\begin{description}
\item[Unit Constraint] It is the multifunctor
\begin{equation}\label{Ust-unit}
i \cn \Mtu \to \Mtup = \Mtu \bincoprod \Mterm
\end{equation}
given by the inclusion of the $\Mtu$ summand in $\Mtup$.
\item[Monoidal Constraint] Its component for small pointed multicategories $\M$ and $\N$ is the multifunctor  
\begin{equation}\label{Ust-monconstraint}
\vpi_{\M,\N} \cn \M \otimes \N \to \M \sma \N
\end{equation}
given by the right vertical arrow in the pushout \cref{eq:multicat-smash-pushout} that defines the smash product.\defmark
\end{description}
\end{explanation}

\begin{explanation}[$\Ust$ as a $\Cat$-Multifunctor]\label{expl:Ust-catmulti}
Regarding $\Multicat$ and $\pMulticat$ as $\Cat$-multicategories as in \cref{expl:multicatcatmulticat,expl:ptmulticatcatmulticat}, respectively, the symmetric monoidal $\Cat$-functor $\Ust$ in \cref{Ust} induces a $\Cat$-multifunctor
\begin{equation}\label{Ust-multifunctor}
\Ust \cn \pMulticat \to \Multicat
\end{equation}
in the sense of \cref{def:enr-multicategory-functor} with the following structure:
\begin{description}
\item[Object Assignment] 
$\Ust$ sends a small pointed multicategory $(\M,i)$ to the multicategory $\M$.
\item[Multimorphism Functors] 
Suppose given small pointed multicategories $\ang{\M} = \ang{\M_j}_{j=1}^n$ and $\N$.  The $n$-ary multimorphism functor
\[\Ust \cn \pMulticat\brb{\txsma_{j=1}^n \M_j , \N} \to \Multicat\brb{\txotimes_{j=1}^n \M_j , \N}\]
sends a pointed multifunctor
\[P \cn \txsma_{j=1}^n \M_j \to \N\]
to the composite multifunctor
\[\begin{tikzpicture}[xscale=2.5,yscale=1,vcenter]
\draw[0cell=.9]
(0,0) node (a) {\txotimes_{j=1}^n \M_j}
(a)++(1,0) node (b) {\txsma_{j=1}^n \M_j}
(b)++(1,0) node (c) {\N}
;
\draw[1cell=.9]  
(a) edge node {\varpi} (b)
(b) edge node {P} (c)
;
\end{tikzpicture}\]
if $n > 0$, where $\varpi$ is an iterate of the monoidal constraint in \cref{Ust-monconstraint}.  If $n=0$, then $\Ust$ sends $P$ to the composite multifunctor
\[\begin{tikzpicture}[xscale=2,yscale=1,vcenter]
\draw[0cell=.9]
(0,0) node (a) {\Mtu}
(a)++(1,0) node (b) {\Mtup}
(b)++(1,0) node (c) {\N}
;
\draw[1cell=.9]  
(a) edge node {i} (b)
(b) edge node {P} (c)
;
\end{tikzpicture}\]
where $i$ is the unit constraint in \cref{Ust-unit}.  For a pointed multinatural transformation, $\Ust$ is defined similarly by whiskering with $\varpi$ if $n>0$ and with $i$ if $n=0$.\defmark
\end{description}
\end{explanation}

\section{\texorpdfstring{$\Mone$}{M1}-Modules}
\label{sec:monemodules}

In this section we review a full sub-2-category of $\pMulticat$ given by the left modules over a small pointed multicategory $\Mone$ (\cref{ex:mofone}).
\begin{itemize}
\item In \cref{definition:calM,definition:part-prod} we discuss the partition multicategory $\cM a$ of a pointed finite set $a$ and a pairing called the partition product.
\item The partition multicategory $\cM \ord{1}$ is equipped with the structure of a commutative monoid in $(\pMulticat, \sma, \Mtup)$ in \cref{def:Monecommonoid}.
\item The 2-category of left $\Mone$-modules is in \cref{definition:MoneMod-prelim}.  Its main properties are summarized in \cref{proposition:EM2-5-1}.
\item The symmetric monoidal $\Cat$-category of left $\Mone$-modules is in \cref{definition:MoneMod}.  Its induced $\Cat$-multicategory is discussed in \cref{expl:monemodcatmulticat}.
\end{itemize} 
The material in this section is adapted from \cite[Chapters 8 and 10]{cerberusIII}.

\subsection*{Partition Multicategories}

Recall from \cref{def:ordn} the permutative category
\[\big(\Fskel, \sma, \ord{1}, \xi\big)\]
of pointed finite sets and pointed functions with the smash product as the monoidal product.

\begin{definition}\label{definition:calM}
For a pointed finite set $a$ with basepoint $*$, the unpointed finite set $a^\punc$ is obtained from $a$ by removing its basepoint:
\begin{equation}\label{aflat}
a^\punc = a \setminus \{*\}.
\end{equation}
The \index{partition multicategory}\index{multicategory!partition}\emph{partition multicategory}, denoted $\cM a$, is the pointed multicategory defined as follows.
\begin{description}
\item[Objects] $\Ob(\cM a) = \Pset(a^\punc)$, the set of basepoint-free subsets of $a$. 
\item[Multimorphisms] For an $n$-tuple $\ang{s} = \ang{s_j}_{j=1}^n$ with each $s_j \in \Pset(a^\punc)$ and $t \in \Pset(a^\punc)$, the set of $n$-ary operations is\index{partition}
\[(\cM a)\mmap{t;\ang{s}} = \begin{cases}
\left\{\iota_{\ang{s}}\right\} & \text{if $\ang{s}$ is a partition of $t$ and}\\
\varnothing & \text{otherwise}. 
\end{cases}\]
In the first case, $\left\{\iota_{\ang{s}}\right\}$ is a one-element set.
\item[Pointed Structure]
It is the multifunctor $i \cn \Mterm \to \cM a$ given by 
\begin{itemize}
\item the empty set $\varnothing \in \Pset(a^\punc)$ and
\item the unique operations 
\[\iota^n = \iota_{\ang{\varnothing}_{j=1}^n} \in (\cM a)\mmap{\varnothing; \ang{\varnothing}_{j=1}^n} \forspace n \ge 0.\]
\end{itemize} 
\item[Other Structure] 
The colored units, symmetric group action, and composition are uniquely defined by the terminal property of a one-element set.
\end{description}  
This finishes the definition of the partition multicategory $\cM a$.
\end{definition}

An important pointed multicategory for this work is the partition multicategory $\Mone$ of the pointed finite set $\ord{1} = \{0,1\}$.

\begin{example}\label{ex:mofone}
We explain in detail the \index{partition multicategory!M1@$\Mone$}\index{multicategory!partition!M1@$\Mone$}partition multicategory $\Mone$.
\begin{itemize}
\item Its object set is
\[\Ob(\Mone) = \Pset(\ord{1}^\punc) = \big\{ \emptyset, \{1\} \big\}.\]
\item Its nonempty sets of operations are\label{not:moneoperations}
\[\begin{split}
\cM\ord{1} \mmap{\emptyset; \ang{\emptyset}_{j=1}^n} &= \{\iota^n\}\\
\cM\ord{1} \mmap{\{1\}; (\emptyset,\ldots,\{1\},\ldots,\emptyset)} &= \{\pi^n_j\}
\end{split}\]  
for $n \geq 0$ and $j \in \{1,\ldots,n\}$.  In the definition of $\pi^n_j$ above, 
\[(\emptyset,\ldots,\{1\},\ldots,\emptyset)\] 
has length $n$ with $\{1\}$ in the $j$th entry and $\emptyset$ in other entries.  The operations $\iota^n$ for $n \geq 0$ are closed under the symmetric group action and composition.
\item  The $\{1\}$-colored unit is 
\[\pi^1_1 \in \cM\ord{1}\mmap{\{1\}; \{1\}}.\]
\item The right $\Sigma_n$-action on $\pi^n_j$ is given by
\[\pi^n_j \cdot \sigma = \pi^n_{\sigmainv(j)} \forspace \sigma \in \Sigma_n.\]
\item The composition involving $\pi^n_j$ is given by
\[\ga\lrscmap{\pi^n_j ; \big( \ang{\iota^{k_i}}_{i=1}^{j-1}, \pi^{k_j}_p, \ang{\iota^{k_i}}_{i=j+1}^n\big)}
= \pi^{k_1+\cdots+k_n}_{k_1+\cdots+k_{j-1}+p}\]
for $k_1,\ldots,k_n \geq 0$ and $p \in \{1,\ldots,k_j\}$. 
\end{itemize}  
This finishes the description of the partition multicategory $\cM\ord{1}$.
\end{example}

\subsection*{The Commutative Monoid $\Mone$}

The partition multicategory $\Mone$ becomes a commutative monoid (\cref{def:commonoid}) via the following multiplication.  Recall the smash product defined in \cref{eq:multicat-smash-pushout}.

\begin{definition}\label{definition:part-prod}
For a pair of pointed finite sets $a$ and $b$, the \index{partition product}\index{product!partition}\emph{partition product} is the pointed multifunctor 
\begin{equation}\label{eq:part-prod}
\txprod_{a,b} \cn \cM a \sma \cM b \to \cM(a \sma b)
\end{equation}
defined as follows.
\begin{description}
\item[Object Assignment] This is given by the Cartesian product of subsets, noting that
\[s \times t \subset (a^\punc \times b^\punc) \iso (a \sma b)^\punc
\forspace s \subset a^\punc \andspace t \subset b^\punc.\]
\item[Multimorphism Assignment] By \cref{explanation:unpacking-tensor}, the generating operations of the tensor product
\[\cM a \otimes \cM b\]
are of the form
\[\iota_{\ang{s}} \otimes t \in \cM a\mmap{s';\ang{s}} \times \{t\}
\andspace 
s \otimes \iota_{\ang{t}}  \in \{s\} \times \cM b \mmap{t';\ang{t}}\]
with
\begin{itemize}
\item $s \subset a^\punc$, $t \subset b^\punc$,
\item $\ang{s}$ a partition of $s'$ in $\cM a$, and
\item $\ang{t}$ a partition of $t'$ in $\cM b$.
\end{itemize}
We define partitions of $s' \times t$ and $s \times t'$ by, respectively,
\[\ang{s} \times t = \ang{s_i \times t}_i
\andspace
s \times \ang{t} = \ang{s \times t_j}_j.\]
Then we define $\txprod_{a,b}$ on generating operations by
\[\iota_{\ang{s}} \otimes t \mapsto \iota_{\ang{s} \times t}
\andspace
s \otimes \iota_{\ang{t}} \mapsto \iota_{s \times \ang{t}}.\]
\end{description} 
The definition of $\txprod_{a,b}$ descends to the smash product $\cM a \sma \cM b$ because the Cartesian product of any set with the empty set is empty.  Therefore, each generating operation of the form
\[\iota_{\ang{s}} \times \varnothing
\orspace
\varnothing \times \iota_{\ang{t}}\]
is sent to a partition of the empty set, which is a basepoint operation of $\cM (a \sma b)$.  This finishes the definition of the partition product $\txprod_{a,b}$.
\end{definition}

\begin{lemma}\label{lemma:part-prod-1}
For each pointed finite set $b$, the partition products for $\ord{1}$ and $b$ are isomorphisms
\begin{align}
  \txprod_{\ord{1},b} \cn\ & \cM\ord{1} \sma \cM b \fto{\ \iso\ } \cM(\ord{1} \sma b) \iso \cM{b}\label{eq:part-prod-1}\\
  \txprod_{b,\ord{1}} \cn\ & \cM{b} \sma \cM\ord{1} \fto{\ \iso\ } \cM(b \sma \ord{1}) \iso \cM{b}.\nonumber
  \end{align}
\end{lemma}

Recall from \cref{thm:pmulticat-smclosed} the symmetric monoidal category $\pMulticat$.

\begin{proposition}\label{proposition:cM-symm-mon}\index{partition multicategory!symmetric monoidal functor}\index{multicategory!partition!symmetric monoidal functor}\index{symmetric monoidal functor!partition multicategory}
The partition multicategory $\cM$ defines a symmetric monoidal functor
\[\cM \cn (\Fskel^\op, \sma, \ord{1}) \to (\pMulticat, \sma, \Mtup)\]
with the following structure morphisms.
\begin{description}
\item[Unit Constraint] It is the pointed multifunctor 
\begin{equation}\label{eq:cM0}
\cM^0\cn \Mtup = \Mtu \bincoprod \Mterm \to \Mone
\end{equation}
determined by sending the unique object of $\Mtu$ to $\{1\} \in \Pset(\ord{1}^\punc)$.
\item[Monoidal Constraint] It is the composite of the partition product $\txprod$ with the lexicographic isomorphisms
\begin{equation}\label{eq:cM2}
\cM^2_{\ord{m},\ord{n}} \cn \cM\ord{m} \sma \cM\ord{n}
\fto{\txprod_{\ord{m},\ord{n}}}
\cM(\ord{m} \sma \ord{n}) \iso \cM(\ord{mn}).
\end{equation}
\end{description}
\end{proposition}

\begin{explanation}\label{expl:Mnotstrong}
The symmetric monoidal functor $\cM$ is neither strong nor strictly unital because the unit and monoidal constraints, $\cM^0$ and $\cM^2$ in \cref{eq:cM0,eq:cM2}, are not isomorphisms.
\end{explanation}

Recall from \cref{def:commonoid} the notion of a \index{commutative monoid}\index{monoid!commutative}\emph{commutative monoid} in a symmetric monoidal category.  The following definition uses the symmetric monoidal category $(\pMulticat, \sma, \Mtup)$ in \cref{thm:pmulticat-smclosed}.

\begin{definition}\label{def:Monecommonoid}\index{multicategory!partition!M1@$\Mone$}\index{partition multicategory!M1@$\Mone$}\index{M1@$\Mone$}
We define the commutative monoid
\[\big(\Mone, \txprod_{\ord{1},\ord{1}} \scs \cM^0\big)\]
in the symmetric monoidal category $(\pMulticat, \sma, \Mtup)$ as follows.
\begin{description}
\item[Object] It is the partition multicategory $\Mone$ (\cref{ex:mofone}).
\item[Multiplication] It is the partition product for $\ord{1}$ and $\ord{1}$ in \cref{eq:part-prod-1}:
\[\txprod_{\ord{1},\ord{1}} \cn \Mone \sma \Mone \fto{\iso} \Mone.\]
\item[Unit] It is the unit constraint $\cM^0 \cn \Mtup \to \Mone$ in \cref{eq:cM0}.\defmark
\end{description}
\end{definition}

\subsection*{The 2-Category of $\Mone$-Modules}

Recall from \cref{def:modules} that each monoid in a monoidal category has an associated category of left modules.  The following definition uses the commutative monoid $\Mone$ in $\pMulticat$ in \cref{def:Monecommonoid}.

\begin{definition}\label{definition:MoneMod-prelim}\index{2-category!of M1-modules@of $\Mone$-modules}\index{partition multicategory!M1-modules@$\Mone$-modules}\index{multicategory!partition!M1-modules@$\Mone$-modules}\index{M1-modules@$\Mone$-modules}\index{module!over M1@over $\Mone$}
We define the 2-category of left $\Mone$-modules, denoted $\MoneMod$, as follows.
\begin{itemize}
\item It has objects and 1-cells given by, respectively, left $\Mone$-modules and their morphisms as in \cref{def:modules}.
\item For left $\Mone$-module morphisms
\[F,F'\cn \N \to \N',\]
the set of \emph{left $\Mone$-module 2-cells} consists of pointed multinatural transformations (\cref{def:ptd-multicat})
\[\theta \cn F \to F'\]
such that the two whiskerings in \cref{eq:M1-mult-2nat} below are equal, where $\mu$ and $\mu'$ denote the left $\Mone$-module structures for $\N$ and $\N'$, respectively.
\begin{equation}\label{eq:M1-mult-2nat}
    \begin{tikzpicture}[x=35mm,y=15mm,vcenter]
    \def\w{15}
      \draw[0cell=.85]
      (0,0) node (a) {\Mone \sma \N}
      (1,0) node (b) {\Mone \sma \N'}
      (0,-1) node (c) {\N}
      (1,-1) node (d) {\N'}
      ;
      \draw[1cell=.85]
      (a) edge[bend left=\w] node[pos=.4] {1 F} (b)
      (c) edge[bend left=\w] node[pos=.4] {F} (d)
      (a) edge['] node {\mu} (c)
      (b) edge node {\mu'} (d)
      (a) edge[bend right=\w,'] node[pos=.6] {1 F'} (b)
      (c) edge[bend right=\w,'] node[pos=.6] {F'} (d)
      ;
      \draw[2cell]
      node[between=a and b at .5, rotate=-90, 2label={below,1 \theta}] {\Rightarrow}
      node[between=c and d at .5, rotate=-90, 2label={below,\theta}] {\Rightarrow}
      ;
    \end{tikzpicture}
\end{equation}
\item Identities and compositions in $\MoneMod$ are given by those of $\Multicat_*$ in \cref{thm:pmulticat}.
\end{itemize}
We use the same notation for the underlying 1-category of $\MoneMod$.
\end{definition}

\begin{example}[Endomorphism $\Mone$-Modules]\label{ex:endmc}
Each small permutative category $(\C,\oplus,\pu,\xi)$ has an associated left $\Mone$-module
\[\Endm(\C) = \big(\Endst(\C), \mu\big)\]
defined as follows.
\begin{itemize}
\item $\Endst(\C)$ is the pointed endomorphism multicategory in \cref{ex:endstc}.
\item The left $\Mone$-module structure
\[\mu \cn \Mone \sma \Endst(\C) \to \Endst(\C)\] 
is given by the following assignments for $a \in \ObC$ and multimorphisms $f$ in $\Endst(\C)$.
  \begin{align*}
    (\varnothing,a) & \mapsto \pu
    & \iota^n \sma a & \mapsto \iota^n = 1_{\pu}
    & \varnothing \sma f & \mapsto 1_{\pu}\\
    (\{1\},a) & \mapsto a
    & \pi^n_j \sma a & \mapsto 1_a
    & \{1\} \sma f & \mapsto f
  \end{align*}
\end{itemize}
As in \cref{ex:endc,ex:endstc}, it is enough to assume that $(\C,\otimes,\tu)$ is a symmetric monoidal category.  In this more general case,
\begin{itemize}
\item the image of $\iota^n \sma a$ is an iterate of the right unit isomorphism $\rho$ in $\C$, and
\item the image of $\pi^n_j \sma a$ is an iterate of the left unit isomorphism $\lambda$ and the right unit isomorphism $\rho$ in $\C$.
\end{itemize} 

Moreover, the following statements hold:
\begin{enumerate}
\item Each \emph{strictly unital} symmetric monoidal functor between symmetric monoidal categories
\[(P,P^2,P^0=1) \cn \C \to \D\]
induces a left $\Mone$-module morphism
\[\Endm(P) \cn \Endm(\C) \to \Endm(\D)\]
given by the pointed multifunctor $\Endst(P)$ in \cref{ex:endstc}.  The latter is given by $\End(P)$ in \cref{EndP}.
\item Each monoidal natural transformation between strictly unital symmetric monoidal functors between symmetric monoidal categories
\[\theta \cn (P,P^2,P^0=1) \to (Q,Q^2,Q^0=1) \cn \C \to \D\]
induces a left $\Mone$-module 2-cell
\[\begin{tikzpicture}
\def\h{3} \def\t{27}
\draw[0cell]
(0,0) node (a) {\phantom{X}}
(a)+(-.6,0) node (a') {\Endm(\C)}
(a)+(\h,0) node (b) {\phantom{X}}
(b)+(.6,0) node (b') {\Endm(\D)}
;
\draw[1cell=.8]
(a) edge[bend left=\t] node {\Endm(P)} (b)
(a) edge[bend right=\t] node[swap] {\Endm(Q)} (b)
;
\draw[2cell]
node[between=a and b at .3, rotate=-90, 2label={above,\Endm(\theta)}] {\Rightarrow}
;
\end{tikzpicture}\]
given by the pointed multinatural transformation $\Endst(\theta)$ in \cref{ex:endstc}.
\end{enumerate}
If there is no danger of confusion, we denote $\Endm(\C)$ by $\C$.
\end{example}

\begin{proposition}\label{endmtwofunctor}
The endomorphism left $\Mone$-module in \cref{ex:endmc} defines a 2-functor
\[\Endm \cn \permcatsu \to \MoneMod.\]
\end{proposition}

We also denote by $\Endm$ the restriction of the 2-functor in \cref{endmtwofunctor} to the locally-full sub-2-category $\permcatst$ in \cref{def:permcat}.

The following result from \cite[10.1.14, 10.1.28, 10.2.5, and 10.2.22]{cerberusIII} summarizes some of the main properties of $\Mone$-modules.

\begin{proposition}\label{proposition:EM2-5-1}\
\begin{enumerate}
\item\label{it:EM251-1} Suppose $\N$ is a left $\cM{\ord{1}}$-module in $\Multicat_*$.  Then the structure morphism
    \[
    \mu\cn \cM\ord{1} \sma \N \to \N
    \]
    is an isomorphism.  Its inverse is given by the unit
\begin{equation}\label{Monemodinv}
\N \fto[\iso]{\la^\inv} \Mtup \sma \N \fto{\cM^0 \sma 1} \cM\ord{1} \sma \N
\end{equation}    
with
\begin{itemize}
\item $\la$ the left unit isomorphism for $\sma$ and
\item $\cM^0 \cn \Mtup \to \Mone$ the unit constraint in \cref{eq:cM0}.
\end{itemize}
\item\label{it:EM251-2} The structure morphism for a
    right $\Mone$-module is an isomorphism with inverse given by $(1 \sma
    \eta) \circ \rho^\inv$, where $\rho$ is the right unit isomorphism
    for $\sma$.
\item\label{it:EM251-uniqueness} Each small pointed multicategory $\N$ admits at most one left $\Mone$-module structure and at most one right $\Mone$-module structure.
\item\label{it:EM251-3} The 2-category of left, respectively right, $\Mone$-modules is a full sub-2-category of $\Multicat_*$ in \cref{thm:pmulticat}.
\item\label{it:EM251-4} For a left $\cM\ord{1}$-module $(\N,\mu^\N)$
    and a right $\cM\ord{1}$-module $(\P,\mu^\P)$, the two morphisms
    in $\Multicat_*$
   \[ \begin{tikzpicture}[x=40mm,y=20mm,vcenter]
      \draw[0cell] 
      (0,0) node (ux) {(\P \sma \cM\ord{1}) \sma \N}
      (1,0) node (uy) {\P \sma \N}
      ;
      \draw[1cell] 
      (ux) edge[transform canvas={yshift=.8mm}] node {\mu^\P \sma 1} (uy)
      (ux) edge[transform canvas={yshift=-.8mm}] node['] {(1 \sma
        \mu^\N) \circ \al} (uy)
      ;
    \end{tikzpicture}\]
    are equal.  Therefore, the canonical morphism to the coequalizer
    \[
    \P \sma \N \fto{\iso} \P \sma_{\cM\ord{1}} \N
    \]
    is an isomorphism in $\Multicat_*$.
\item\label{monetoendc}
For each small permutative category $\C$, there are natural isomorphisms of pointed categories
\[\MoneMod\brb{\Mone,\Endm(\C)} \iso \pMulticat\brb{\Mtup, \Endst(\C)} \iso (\C,\pu).\]
\item\label{monebicomplete} There is a complete, cocomplete, symmetric monoidal, and \index{M1-modules@$\Mone$-modules!closed structure}\index{module!over M1@over $\Mone$!closed structure}closed category
\[\big(\MoneMod, \sma, \Mone, \pHom\big).\] 
\begin{itemize}
\item The monoidal product and internal hom are those of $(\pMulticat,\sma,\pHom)$ in \cref{thm:pmulticat-smclosed}. 
\item The monoidal unit is $\Mone$.
\item The unit isomorphisms are given by the left $\Mone$-module structure and braiding for $\sma$.
\end{itemize}
\end{enumerate}
\end{proposition}

Recall from \cref{def:twoadjunction} the notion of a 2-adjunction.

\begin{proposition}\label{MonesmaUmadj}
There is a free-forgetful 2-adjunction
\[\begin{tikzpicture}
\draw[0cell]
(0,0) node (a) {\pMulticat}
(a)+(3.25,0) node (b) {\MoneMod.}
(a)+(1.55,0) node (x) {\bot}
;
\draw[1cell=.9]
(a) edge[bend left=15,transform canvas={yshift=-.8ex}] node {\Monesma} (b)
(b) edge[bend left=15,transform canvas={yshift=.7ex}] node {\Um} (a)
;
\end{tikzpicture}\]
\end{proposition}

\begin{explanation}[Unit and Counit]\label{expl:MonesmaUmadj}
For the 2-adjunction $(\Monesma) \dashv \Um$, the unit 
\begin{equation}\label{etam-def}
\etahat \cn 1_{\pMulticat} \to \Um \circ~ (\Monesma)
\end{equation}
has component pointed multifunctor given by the composite
\[\begin{tikzpicture}[baseline={(a.base)}]
\def\u{.6}
\draw[0cell]
(0,0) node (a) {\M}
(a)+(2,0) node (b) {\Mtup \sma \M}
(b)+(3,0) node (c) {\Mone \sma \M}
;
\draw[1cell=.9]
(a) edge node {\lambda^\inv} node[swap] {\iso} (b)
(b) edge node {\cM^0 \sma 1} (c)
;
\draw[1cell=.9]
(a) [rounded corners=3pt] |- ($(b)+(-1,\u)$)
-- node {\etahat_\M} ($(b)+(1,\u)$) -| (c)
;
\end{tikzpicture}\]
in \cref{Monemodinv} for each small pointed multicategory $\M$.  This component is, in general, \emph{not} an isomorphism because $\M$ need not be a left $\Mone$-module.

The counit
\begin{equation}\label{epzm-def}
\epzhat \cn (\Monesma) \circ \Um \to 1_{\MoneMod}
\end{equation}
of the 2-adjunction $(\Monesma) \dashv \Um$ has component 
\[\epzhat_\N = \mu \cn \Mone \sma \N \fto{\iso} \N\]
given by the left $\Mone$-module structure morphism for each left $\Mone$-module $(\N,\mu)$.  This component is an isomorphism by \cref{proposition:EM2-5-1} \eqref{it:EM251-1}.
\end{explanation}

\subsection*{The Symmetric Monoidal $\Cat$-Category of $\Mone$-Modules}

Recall the notion of a symmetric monoidal $\V$-category (\cref{definition:symm-monoidal-vcat}).  The following definition involves the case $\V = (\Cat,\times,\boldone)$.

\begin{definition}\label{definition:MoneMod}\index{2-category!of M1-modules@of $\Mone$-modules}\index{partition multicategory!M1-modules@$\Mone$-modules}\index{multicategory!partition!M1-modules@$\Mone$-modules}\index{M1-modules@$\Mone$-modules}\index{module!over M1@over $\Mone$}
Define the symmetric monoidal $\Cat$-category
\[\big(\MoneMod, \sma ,\Mone\big)\]
with the following data.
\begin{itemize}
\item The base $\Cat$-category is the 2-category of left $\Mone$-modules in \cref{definition:MoneMod-prelim}.  By \cref{proposition:EM2-5-1} \eqref{it:EM251-3}, it has hom categories
\[\MoneMod(\N,\N') = \Multicat_*(\N,\N')\]
for left $\Mone$-modules $\N$ and $\N'$.
\item The monoidal composition is given by the smash product, $\sma$, in $\pMulticat$.  This is well defined by \cref{proposition:EM2-5-1} \eqref{it:EM251-4}.
\item The identity object is $\Mone$ (\cref{ex:mofone}).
\item The monoidal unitors and monoidal associator are given by those of $(\pMulticat, \sma)$ in \cref{thm:pmulticat-smclosed}.
\end{itemize}
This finishes the definition.
\end{definition}

\begin{explanation}[$\MoneMod$ as a $\Cat$-Multicategory]\label{expl:monemodcatmulticat}
Since $\MoneMod$ is a symmetric monoidal $\Cat$-category, it has the structure of a $\Cat$-multicategory by \cref{proposition:monoidal-v-cat-v-multicat}, with the following data.
\begin{itemize}
\item Its objects are left $\Mone$-modules.
\item For left $\Mone$-modules $\ang{\N_j}_{j=1}^n$ and $\N'$, the $n$-ary multimorphism category is
\[\begin{split}
\MoneMod\lrscmap{\ang{\N_j}_{j=1}^n ; \N'}
&= \MoneMod\left(\txsma_{j=1}^n \N_j \scs \N' \right)\\
&= \begin{cases}
\pMulticat\left(\txsma_{j=1}^n \N_j \scs \N' \right) & \text{if $n>0$ and}\\
\pMulticat\left(\Mone \scs \N'\right) & \text{if $n=0$}.
\end{cases}
\end{split}\]
If $n>0$, then this category has
\begin{itemize}
\item pointed multifunctors
\begin{equation}\label{monemodnarycat}
\txsma_{j=1}^n \N_j \to \N'
\end{equation}
as objects and
\item pointed multinatural transformations between such pointed multifunctors as morphisms.
\end{itemize}
If $n=0$, then an empty $\sma$ in \cref{monemodnarycat} means $\Mone$, the monoidal unit in $\MoneMod$.
\item The symmetric group action is induced by the braiding of the smash product in $\pMulticat$. 
\item The multicategorical composition is given by smash product and composition of pointed multifunctors, and likewise for pointed multinatural transformations.
\end{itemize}
This finishes the description of the $\Cat$-multicategory $\MoneMod$.  We note the subtle difference between the $\Cat$-multicategories $\MoneMod$ and $\pMulticat$ in \cref{expl:ptmulticatcatmulticat}, especially in arity 0.
\end{explanation}

\begin{proposition}\label{Monesma-CatSM}
  The 2-functor
  \[
    \Monesma \cn \pMulticat \to \MoneMod
  \]
  is a strong symmetric $\Cat$-monoidal functor, hence also a $\Cat$-multifunctor.
\end{proposition}
\begin{proof}
  The unit constraint for $\Monesma$ is the isomorphism
  \[
    \Mone \fto[\iso]{\rho^\inv} \Mone \sma \Mtup 
  \]
  where $\rho$ is the right unit isomorphism for $\sma$.
  The monoidal constraint is the composite isomorphism for $\M,\N \in \pMulticat$
  \[
    (\Mone \sma \M) \sma (\Mone \sma \N) \fto{\iso}
    (\Mone \sma \Mone) \sma (\M \sma \N) \fto{\txprod_{\ord{1},\ord{1}}}
    \Mone \sma (\M \sma \N),
  \]
  where the first isomorphism permutes the factors and the second isomorphism is the partition product from \cref{lemma:part-prod-1} with $b = 1$.
  The symmetric $\Cat$-monoidal functor axioms of \cref{definition:monoidal-V-fun,definition:braided-monoidal-vfunctor} then follow because $\Mone$ is a commutative monoid in $\pMulticat$.
\end{proof}

\begin{explanation}[Forgetting $\Mone$-Module Structure]\label{expl:moneinclusion}
The forgetful 2-functor
\begin{equation}\label{Umonemod}
\Um \cn \big(\MoneMod, \sma, \Mone\big) \to \big(\pMulticat, \sma, \Mtup\big)
\end{equation}
is a symmetric monoidal $\Cat$-functor (\cref{definition:braided-monoidal-vfunctor}) with the following structure morphisms:
\begin{description}
\item[Monoidal Constraint] It is the identity.
\item[Unit Constraint] It is the pointed multifunctor in \cref{eq:cM0},
\[\cM^0 \cn \Mtup \to \Mone.\]
\end{description}
By \cref{proposition:EM2-5-1} \eqref{it:EM251-3} and \eqref{it:EM251-4}, the underlying forgetful 2-functor 
\[\Um \cn \MoneMod \hookrightarrow \pMulticat\]
is an inclusion between the underlying 2-categories.  However, as a symmetric monoidal $\Cat$-functor, it is neither unital nor strong, since $\cM^0$ is not an isomorphism.  Therefore, $\MoneMod$ is \emph{not} a symmetric monoidal $\Cat$-subcategory of $\pMulticat$.
\end{explanation}

\begin{explanation}[$\Um$ as a $\Cat$-Multifunctor]\label{expl:Um-catmulti}
Regarding $\MoneMod$ and $\pMulticat$ as $\Cat$-multicategories using \cref{expl:ptmulticatcatmulticat,expl:monemodcatmulticat}, there is an induced $\Cat$-multifunctor (\cref{def:enr-multicategory-functor})
\begin{equation}\label{Um-multifunctor}
\Um \cn \MoneMod \to \pMulticat
\end{equation}
with the following structure:
\begin{description}
\item[Object Assignment]
$\Um$ sends a left $\Mone$-module $(\M,\mu)$ to the pointed multicategory $\M$.
\item[Multimorphism Functors]
Suppose given left $\Mone$-modules $\angM = \ang{\M_j}_{j=1}^n$ and $\N$.  The $n$-ary multimorphism functor
\[\Um \cn \MoneMod\brb{\txsma_{j=1}^n \M_j , \N} \to \pMulticat\brb{\txsma_{j=1}^n \M_j , \N}\]
is an isomorphism if $n > 0$ by \cref{proposition:EM2-5-1} \eqref{it:EM251-3} and \eqref{it:EM251-4}.  

If $n=0$, then $\Um$ sends a left $\Mone$-module morphism 
\[P \cn \Mone \to \N\]
to the composite pointed multifunctor
\[\begin{tikzpicture}[xscale=2,yscale=1,vcenter]
\draw[0cell=.9]
(0,0) node (a) {\Mtup}
(a)++(1,0) node (b) {\Mone}
(b)++(1,0) node (c) {\N}
;
\draw[1cell=.9]  
(a) edge node {\cM^0} (b)
(b) edge node {P} (c)
;
\end{tikzpicture}\]
where $\cM^0$ is the unit constraint in \cref{eq:cM0}.  For a left $\Mone$-module 2-cell (\cref{definition:MoneMod-prelim}), $\Um$ is defined similarly by whiskering with $\cM^0$.\defmark
\end{description}
\end{explanation}

\begin{proposition}\label{etahat-epzhat-monCatnat}
  The unit $\etahat$ in \cref{etam-def} and the counit $\epzhat$ in \cref{epzm-def} are monoidal $\Cat$-natural transformations, hence also $\Cat$-multinatural transformations.
\end{proposition}
\begin{proof}
  As in \cref{Monesma-CatSM}, the assertions about $\etahat$ and $\epzhat$ follow from the commutative monoid structure of $\Mone$.
\end{proof}

\section{Permutative Categories}
\label{sec:multpermcat}

There is a 2-category $\permcatsu$ of small permutative categories, strictly unital symmetric monoidal functors, and monoidal natural transformations (\cref{def:permcat}).  In this section we extend the 2-category $\permcatsu$ to a $\Cat$-multicategory.
\begin{itemize}
\item In \cref{def:nlinearfunctor,def:nlineartransformation} we define $n$-linear functors and $n$-linear transformations.  They generalize strictly unital symmetric monoidal functors and monoidal natural transformations, respectively.
\item The $\Cat$-multicategory $\permcatsu$ is in \cref{definition:permcatsus-homcat,definition:permcat-action,definition:permcat-comp}.
\item \cref{endfactorization} shows that the $\Cat$-multicategories $\permcatsu$, $\Multicat$, $\pMulticat$, and $\MoneMod$ are related by the various endomorphism multicategory constructions and forgetful functors.
\end{itemize}
The material in this section is adapted from \cite[Sections~6.5 and~6.6]{cerberusIII}.

Throughout this section, suppose $\angC = \ang{\C_j}_{j=1}^n$ and $\D$ are permutative categories.  In each of these permutative categories, the monoidal product, monoidal unit, and braiding are denoted by $\oplus$, $\pu$, and $\xi$, respectively.

\subsection*{Multilinear Functors and Transformations}

\begin{notation}\label{notation:compk}
Suppose $\ang{x} = \ang{x_j}_{j=1}^n$ is an $n$-tuple of symbols, and $y$ is a symbol with $k \in \{1,\ldots,n\}$.  We denote by\label{not:compk}
\[\ang{x \compk y} = \ang{x} \compk y 
= \big(\underbracket[0.5pt]{x_1, \ldots, x_{k-1}}_{\text{empty if $k=1$}}, y, \underbracket[0.5pt]{x_{k+1}, \ldots, x_n}_{\text{empty if $k=n$}}\big)\]
the $n$-tuple obtained from $\ang{x}$ by replacing its $k$-th entry by $y$.  Similarly, for $k \neq \ell \in \{1,\ldots,n\}$ and a symbol $z$, we denote by
\[\ang{x \compk y \comp_\ell z} = \ang{x} \compk y \comp_\ell z\]
the $n$-tuple obtained from $\ang{x \compk y}$ by replacing its $\ell$-th entry by $z$.
\end{notation}

\begin{definition}\label{def:nlinearfunctor}
An \index{n-linear@$n$-linear!functor}\index{functor!n-linear@$n$-linear}\emph{$n$-linear functor}
\[\begin{tikzcd}[column sep=2.3cm]
\txprod_{j=1}^n \C_j \ar{r}{\big(P,\, \{P^2_j\}_{j=1}^n\big)} & \D
\end{tikzcd}\]
consists of the following data.
\begin{itemize}
\item $P \cn \prod_{j=1}^n \C_j \to \D$ is a functor.
\item For each $j\in \{1,\ldots,n\}$, $P^2_j$ is a natural transformation, called the \index{linearity constraint}\index{constraint!linearity}\emph{$j$-th linearity constraint}, with component morphisms
\begin{equation}\label{laxlinearityconstraints}
\begin{tikzcd}[column sep=large]
P\ang{x \compj x_j} \oplus P\ang{x \compj x_j'} \ar{r}{P^2_j}
& P\bang{x \compj (x_j \oplus x_j')} \in \D
\end{tikzcd}
\end{equation}
for objects $\ang{x} \in \txprod_{j=1}^n \C_j$ and $x_j' \in \C_j$.
\end{itemize}
These data are required to satisfy the axioms \cref{nlinearunity,constraintunity,eq:ml-f2-assoc,eq:ml-f2-symm,eq:f2-2by2} below.
\begin{description}
\item[Unity]
For objects $\angx$ and morphisms $\angf$ in $\txprod_{j=1}^n \C_j$, the following object and morphism equalities hold for each $j \in \{1,\ldots,n\}$.
\begin{equation}\label{nlinearunity}
    \left\{\begin{aligned}
        P \ang{x \compj \pu} &= \pu\\
        P \ang{f \compj 1_{\pu}} &= 1_{\pu}
      \end{aligned}\right.
\end{equation}
\item[Constraint Unity]\index{constraint!- axiom!unity}
  \begin{equation}\label{constraintunity}
    P^2_j = 1 \qquad \text{if any $x_i = \pu$ or if $x_j' = \pu$}.
  \end{equation}
\item[Constraint Associativity]\index{constraint!- axiom!associativity} The following diagram commutes for each $i\in \{1,\ldots,n\}$ and objects $\ang{x} \in \txprod_{j=1}^n \C_j$, with $x_i', x_i'' \in \C_i$.
    \begin{equation}\label{eq:ml-f2-assoc}

    \end{equation}
  \item[Constraint Symmetry]\index{constraint!- axiom!symmetry} The following diagram commutes
    for each $i\in \{1,\ldots,n\}$ and objects $\ang{x} \in \txprod_{j=1}^n \C_j$, with
    $x_i' \in \C_i$.
    \begin{equation}\label{eq:ml-f2-symm}
    %
    \end{equation}
  \item[Constraint 2-By-2]\index{constraint!- axiom!2-by-2} The following diagram commutes
    for each
    \[
    i \neq k\in \{1,\ldots,n\}, 
    \quad
    \ang{x} \in \txprod_{j=1}^n \C_j,
    \quad
    x_i' \in \C_i,
    \andspace
    x_k' \in \C_k.
    \]
\begin{equation}\label{eq:f2-2by2}
%
\end{equation}
\end{description}
This finishes the definition of an $n$-linear functor.  

Moreover, we define the following.
\begin{itemize}
\item If $n=0$, then a 0-linear functor is a functor $\boldone \to \D$, which is also regarded as a choice of an object in $\D$.
\item An $n$-linear functor $(P, \{P^2_j\})$ is
\begin{itemize}
\item \index{multilinear functor!strong}\index{strong multilinear functor}\emph{strong} if each $P^2_j$ is a natural isomorphism and
\item \index{multilinear functor!strict}\index{strict multilinear functor}\emph{strict} if each $P^2_j$ is an identity natural transformation.
\end{itemize} 
\item A \index{multilinear functor}\index{functor!multilinear}\emph{multilinear functor} is an $n$-linear functor for some $n \geq 0$.\defmark
\end{itemize}
\end{definition}

\begin{example}\label{ex:onelinearfunctor}
A 1-linear functor $\C \to \D$ is precisely a \index{functor!symmetric monoidal!strictly unital}\index{strictly unital!monoidal functor}strictly unital symmetric monoidal functor (\cref{def:monoidalfunctor}).
\end{example}

\begin{definition}\label{def:nlineartransformation}
Suppose $P,Q$ are $n$-linear functors as displayed below.
\begin{equation}\label{nlineartransformation}
\begin{tikzpicture}[xscale=3,yscale=2.5,baseline={(a.base)}]
\draw[0cell=.9]
(0,0) node (a) {\prod_{j=1}^n \C_j}
(a)++(1,0) node (b) {\D}
;
\draw[1cell=.9]  
(a) edge[bend left=25] node[pos=.45] {\big( P, \{P^2_j\}\big)} (b)
(a) edge[bend right=25] node[swap,pos=.45] {\big( Q, \{Q^2_j\}\big)} (b)
;
\draw[2cell] 
node[between=a and b at .47, shift={(0,0)}, rotate=-90, 2label={above,\theta}] {\Rightarrow}
;
\end{tikzpicture}
\end{equation}
An \index{n-linear@$n$-linear!transformation}\emph{$n$-linear transformation} $\theta \cn P \to Q$ is a natural transformation of underlying functors that satisfies the following two \index{multilinearity conditions}\emph{multilinearity conditions}.
\begin{description}
\item[Unity] 
\begin{equation}\label{niitransformationunity}
\theta_{\ang{x}} = 1_{\pu} \qquad \text{if any $x_i = \pu \in \C_i$}.
\end{equation}
\item[Constraint Compatibility] The following diagram commutes for each $\ang{x} \in \txprod_{j=1}^n \C_j$ and $x_i' \in \C_i$ with $i \in \{1,\ldots,n\}$.
\begin{equation}\label{eq:monoidal-in-each-variable}
\begin{tikzpicture}[x=45mm,y=15mm,vcenter]
\tikzset{0cell/.append style={nodes={scale=.9}}}
\tikzset{1cell/.append style={nodes={scale=.9}}}
      \draw[0cell] 
      (0,0) node (a) {
        P\ang{x\compi x_i}
        \oplus P\ang{x\compi x_i'}
      }
      (1,0) node (b) {
        P\bang{x\compi (x_i \oplus x_i')}
      }
      (0,-1) node (c) {
        Q\ang{x\compi x_i}
        \oplus Q\ang{x\compi x_i'}
      }
      (1,-1) node (d) {
        Q\bang{x\compi (x_i \oplus x_i')}
      }
      ;
      \draw[1cell] 
      (a) edge node {P^2_i} (b)
      (a) edge['] node {\theta \oplus \theta} (c)
      (b) edge node {\theta} (d)
      (c) edge node {Q^2_i} (d)
      ;
  \end{tikzpicture}
\end{equation}
\end{description}
This finishes the definition of an $n$-linear transformation.  Moreover, we define the following.
\begin{itemize}
\item A \index{multilinear transformation}\index{transformation!multilinear}\emph{multilinear transformation} is an $n$-linear transformation for some $n \geq 0$. 
\item Identities and compositions of multilinear transformations are defined componentwise.\defmark
\end{itemize}
\end{definition}

\begin{example}\label{ex:onelineartr}
A 1-linear transformation between 1-linear functors is precisely a \index{monoidal natural transformation}\index{natural transformation!monoidal}monoidal natural transformation (\cref{def:monoidalnattr}) between corresponding strictly unital symmetric monoidal functors.
\end{example}

\subsection*{$\Cat$-Multicategory Structure}

Next we define the $\Cat$-multicategory (\cref{def:enr-multicategory}) $\permcatsu$ whose objects are small permutative categories.  For the rest of this section, $\angC = \ang{\C_j}_{j=1}^n$ and $\D$ are small permutative categories.  The notation in the following definition is chosen to match with the notation in \cref{def:permcat} in the 1-linear case; see \cref{ex:onelinearfunctor,ex:onelineartr}.

\begin{definition}[Multimorphism Categories]\label{definition:permcatsus-homcat}\index{multimorphism!category}\index{n-linear@$n$-linear!multimorphism category}\index{multilinear!multimorphism category}
We define the following categories of $n$-linear functors and transformations.
\begin{itemize}
\item $\permcatsu\scmap{\ang{\C};\D}$ is the category with
\begin{itemize}
\item $n$-linear functors $\ang{\C} \to \D$ as objects and
\item $n$-linear transformations between them as morphisms.
\end{itemize} 
\item $\permcatsus\scmap{\ang{\C};\D}$ is the full subcategory of strong $n$-linear functors.
\item $\permcatst\scmap{\ang{\C};\D}$ is the full subcategory of strict $n$-linear functors.\defmark
\end{itemize}
\end{definition}

\begin{definition}[Symmetric Group Action]\label{definition:permcat-action}
Suppose given $n$-linear functors $P$ and $Q$ together with an $n$-linear transformation $\theta$ as displayed below.
\begin{equation}\label{permiicatcddata}
\begin{tikzpicture}[xscale=3,yscale=2.5,vcenter]
\draw[0cell=.9]
(0,0) node (a) {\prod_{j=1}^n \C_j}
(a)++(1,0) node (b) {\D}
;
\draw[1cell=.9]  
(a) edge[bend left=25] node[pos=.45] {(P, \{P^2_j\})} (b)
(a) edge[bend right=25] node[swap,pos=.45] {(Q, \{Q^2_j\})} (b)
;
\draw[2cell] 
node[between=a and b at .5, shift={(0,0)}, rotate=-90, 2label={below,\theta}] {\Rightarrow}
;
\end{tikzpicture}
\end{equation}
For a permutation $\sigma \in \Sigma_n$, the \index{symmetric group!action!multimorphism category}symmetric group action 
\begin{equation}\label{permiicatsymgroupaction}
\begin{tikzcd}[column sep=large]
\permcatsu\mmap{\D; \ang{\C}} \ar{r}{\sigma}[swap]{\iso} & \permcatsu\mmap{\D; \ang{\C}\sigma}
\end{tikzcd}
\end{equation}
sends the data \cref{permiicatcddata} to the following composites and whiskerings, where $\si$ permutes the coordinates according to $\sigma$.
\begin{equation}\label{permiicatsigmaaction}
\begin{tikzpicture}[xscale=3,yscale=2.5,vcenter]
\draw[0cell=.9]
(0,0) node (a) {\prod_{j=1}^n \C_j}
(a)++(1,0) node (b) {\D}
(a)++(-.75,0) node (c) {\prod_{j=1}^n \C_{\sigma(j)}}
;
\draw[1cell=.9]  
(c) edge node {\sigma} (a)
(a) edge[bend left=25] node[pos=.45] {(P, \{P^2_j\})} (b)
(a) edge[bend right=25] node[swap,pos=.45] {(Q, \{Q^2_j\})} (b)
;
\draw[2cell] 
node[between=a and b at .5, shift={(0,0)}, rotate=-90, 2label={below,\theta}] {\Rightarrow}
;
\end{tikzpicture}
\end{equation}
For objects 
\[\ang{a} = \ang{a_j}_{j=1}^n \in \prod_{j=1}^n \C_{\sigma(j)} \andspace a_j' \in \C_{\sigma(j)},\] 
the $j$-th linearity constraint of $P^\sigma = P \circ \sigma$ has component given by the following composite in $\D$.
\begin{equation}\label{fsigmatwoj}
\begin{tikzcd}[column sep=large, row sep=small]
P^{\sigma} \ang{a} \oplus P^{\sigma} \ang{a \compj a_j'} \ar{r}{(P^{\sigma})^2_j} \ar[equal]{d} & P^{\sigma} \bang{a \compj (a_j \oplus a_j')} \ar[equal]{d}\\
P\big( \sigma\ang{a} \big) \oplus P\big( \sigma\ang{a} \comp_{\sigma(j)} a_j' \big) \ar{r}{P^2_{\sigma(j)}} & P\big( \sigma\ang{a} \comp_{\sigma(j)} (a_j \oplus a_j') \big)
\end{tikzcd}
\end{equation}
If $P$ is strong, respectively strict, with each $P^2_j$ a natural isomorphism, respectively identity, then $P^\sigma$ is also strong, respectively strict.
\end{definition}

\begin{definition}[Multicategorical Composition]\label{definition:permcat-comp}
Suppose given, for each $j \in \{1,\ldots,n\}$,
\begin{itemize}
\item permutative categories $\ang{\B_j} = \ang{\B_{j,i}}_{i=1}^{k_j}$,
\item $k_j$-linear functors $P'_j, Q'_j \cn \ang{\B_j} \to \C_j$, and
\item a $k_j$-linear transformation $\theta_j \cn P'_j \to Q'_j$ as follows.
\end{itemize} 
  \begin{equation}\label{permiicatbcdata}
    \begin{tikzpicture}[xscale=3,yscale=3,vcenter]
      \draw[0cell=.9]
      (0,0) node (a) {\prod_{i=1}^{k_j} \B_{j,i}}
      (a)++(1,0) node (b) {\C_j}
      ;
      \draw[1cell=.9]  
      (a) edge[bend left=20] node[pos=.45] {P'_j} (b)
      (a) edge[bend right=20] node[swap,pos=.45] {Q'_j} (b)
      ;
      \draw[2cell] 
      node[between=a and b at .5, shift={(0,0)}, rotate=-90, 2label={below,\theta_j}] {\Rightarrow}
      ;
    \end{tikzpicture}
  \end{equation}
With $\ang{\B} = \bang{\ang{\B_j}}_{j=1}^n$, the multicategorical composition functor
  \begin{equation}\label{permiicatgamma}
    \begin{tikzcd}[cells={nodes={scale=.85}}]
      \permcatsu\mmap{\D;\ang{\C}} \times \txprod_{j=1}^n \permcatsu\mmap{\C_j;\ang{\B_j}} \ar{r}{\gamma}
      & \permcatsu\mmap{\D;\ang{\B}}
    \end{tikzcd}
  \end{equation}
  sends the data \cref{permiicatcddata,permiicatbcdata} to the composites
  \begin{equation}\label{compositemodification}
    \begin{tikzpicture}[xscale=5,yscale=4,vcenter]
      \draw[0cell=.9]
      (0,0) node (a) {\prod_{j=1}^n \prod_{i=1}^{k_j} \B_{j,i}}
      (a)++(1,0) node (b) {\D}
      ;
      \draw[1cell=.9]  
      (a) edge[bend left=20] node[pos=.45] {P \circ \txprod_j P'_j} (b)
      (a) edge[bend right=20] node[swap,pos=.45] {Q \circ \txprod_j Q'_j} (b)
      ;
      \draw[2cell=.9] 
      node[between=a and b at .6, shift={(0,0)}, rotate=-90, 2label={below,\theta \otimes (\textstyle\prod_j \theta_j)}] {\Longrightarrow}
      ;
    \end{tikzpicture}
  \end{equation}
  defined as follows.

  \begin{description}
  \item[Composite Multilinear Functor]\index{composition!multilinear functor}\index{multilinear functor!composition}
  Suppose given tuples of objects 
  \begin{equation}\label{angwj}
    \begin{aligned}
      \ang{w_j} &= \ang{w_{j,i}}_{i=1}^{k_j} \in \txprod_{i=1}^{k_j} \B_{j,i} \qquad \text{for $j \in \{1,\ldots,n\}$ and}\\
      \ang{w} &= \ang{\ang{w_j}}_{j=1}^n \in \txprod_{j=1}^n \txprod_{i=1}^{k_j} \B_{j,i}.
    \end{aligned}
  \end{equation} 
  Then we have the object
  \begin{equation}\label{compositefangw}
    \begin{aligned}
      \textstyle \big(P \comp \prod_j P_j'\big) \ang{w} &= P \bang{P_j' \ang{w_j}}_{j=1}^n\\
      &= P\big(P_1'\ang{w_1}, \ldots, P_n'\ang{w_n}\big) \inspace \D.
    \end{aligned}
  \end{equation}

For the linearity constraints of the composite $P \circ \prod_j P'_j$ in \cref{compositemodification}, in addition to the objects in \cref{angwj}, consider
  \begin{itemize}
  \item an object $w_{j,i}' \in \B_{j,i}$ for some choice of $(j,i)$ with 
  \[\ell = k_1 + \cdots + k_{j-1} + i\]
  and
  \item $\ang{P'w} = \bang{P'_j\ang{w_j}}_{j=1}^n \in \prod_{j=1}^n \C_j$.
  \end{itemize} 
The following objects appear in \cref{ffjlinearity} below.
  \[\begin{split}
  \ang{w_j \compi w_{j,i}'} &= \big(\overbracket[.5pt]{w_{j,1},\ldots, w_{j,i-1}}^{\text{empty if $i=1$}},\, w_{j,i}' \scs \overbracket[.5pt]{w_{j,i+1},\ldots, w_{j,k_j}}^{\text{empty if $i=k_j$}}\big)\\
  \bang{w_j \compi (w_{j,i} \oplus w_{j,i}')} &= \big(\underbracket[.5pt]{w_{j,1},\ldots, w_{j,i-1}}_{\text{empty if $i=1$}},\, w_{j,i} \oplus w_{j,i}' \scs \underbracket[.5pt]{w_{j,i+1},\ldots, w_{j,k_j}}_{\text{empty if $i=k_j$}}\big).
  \end{split}\]
  The $\ell$-th linearity constraint $\big(P \circ \txprod_j P'_j\big)^2_\ell$ is defined as the following composite in $\D$.
  \begin{equation}\label{ffjlinearity}
    \begin{tikzpicture}[xscale=3,yscale=1,vcenter]
      \def\h{1} \def\v{-1}
      \draw[0cell=.75] 
      (0,0) node (x0) {P\ang{P'w} \oplus P\bang{P'w \compj P'_j\ang{w_j \compi w_{j,i}'}}}
      (x0)++(-\h,\v) node (x1) {(P \circ \txprod_j P'_j)\ang{w} \oplus (P \circ \txprod_j P'_j)\ang{w \comp_\ell w_{j,i}'}}
      (x1)++(1.3*\h,\v) node (x2) {P\bang{P'w \compj (P'_j\ang{w_j} \oplus P'_j\ang{w_j \compi w_{j,i}'})}} 
      (x1)++(0,2*\v) node (x3) {(P \circ \txprod_j P'_j)\bang{w \comp_\ell (w_{j,i} \oplus w_{j,i}')}}
      (x3)++(\h,\v) node (x4) {P\bang{P'w \compj P'_j\bang{w_j \compi (w_{j,i} \oplus w_{j,i}')}}}
      ;
      \draw[1cell=.75] 
      (x1) edge[-,double equal sign distance,transform canvas={xshift=-1em}] (x0)
      (x0) edge node[pos=.7] {P^2_j} (x2)
      (x2) edge node[pos=.3] {P\ang{1 \compj (P'_j)^2_i}} (x4)
      (x1) edge node[swap] {\big(P \circ \txprod_j P'_j\big)^2_\ell} (x3)
      (x3) edge[-,double equal sign distance,transform canvas={xshift=-1em}] (x4)
      ; 
    \end{tikzpicture}
  \end{equation}
If $P$ and each $P'_j$ are strong, respectively strict, then each linearity constraint $\big(P \circ \txprod_j P'_j\big)^2_\ell$ is componentwise invertible, respective identity, and, therefore, the composite $P \circ \txprod_j P'_j$ is also strong, respective strict.

\item[Composite Multinatural Transformation]\index{horizontal composition!multinatural transformation}\index{multinatural transformation!composition}
  The $n$-linear transformation $\theta \otimes \big(\prod_j \theta_j\big)$ in \cref{compositemodification} is the horizontal composite of the natural transformations $\prod_j \theta_j$ and $\theta$.  More explicitly, the component morphism $\big(\theta \otimes (\prod_j \theta_j)\big)_{\angw}$ is the following composite in $\D$.
  \begin{equation}\label{thetaprodthetaw}
    \begin{tikzpicture}[xscale=3.9,yscale=3.5,baseline={(a.base)}]
      \draw[0cell=.9]
      (0,0) node (a) {P \bang{P_j'\angwj}_{j=1}^n}
      (a)++(1.1,0) node (b) {P \bang{Q_j'\angwj}_{j=1}^n}
      (b)++(1,0) node (c) {Q \bang{Q_j'\angwj}_{j=1}^n}
      ;
      \draw[1cell=.9]  
      (a) edge node {P\ang{(\theta_j)_{\angwj}}_{j=1}^n} (b)
      (b) edge node {\theta_{\ang{Q_j'\angwj}_{j=1}^n}} (c)
      ;
    \end{tikzpicture}
  \end{equation}
  \end{description}
The finishes the definition of the multicategorical composition in $\permcatsu$.
\end{definition}

\begin{theorem}\label{thm:permcatmulticat}
There is a $\Cat$-multicategory\index{multicategory!- of permutative categories}\index{permutative category!multicategory of}
\[\permcatsu\]
defined by the following data.
\begin{itemize}
\item The objects are small permutative categories.
\item The multimorphism categories are in \cref{definition:permcatsus-homcat}.
\item The colored units are identity symmetric monoidal functors.
\item The symmetric group action is in \cref{definition:permcat-action}.
\item The multicategorical composition is in \cref{definition:permcat-comp}.
\end{itemize}
Moreover, there are sub-$\Cat$-multicategories
\[\permcatst \hookrightarrow \permcatsus \hookrightarrow \permcatsu\]
with the multimorphism categories in \cref{definition:permcatsus-homcat}.
\end{theorem}

\begin{explanation}\label{expl:permcatsutwocat}
The underlying 2-categories, in the sense of \cref{ex:unarycategory} with $\V = (\Cat, \times, \boldone)$, of the $\Cat$-multicategories 
\[\permcatsu, \quad \permcatsus, \andspace \permcatst\]
are the corresponding 2-categories in \cref{def:permcat}.
\end{explanation}

The following result combines \cite[6.5.10 and 6.5.13]{cerberusIII}.

\begin{proposition}\label{proposition:n-lin-equiv}\index{endomorphism!2-functor}
For small permutative categories $\ang{\C_j}_{j=1}^n$ and $\D$, the 2-functor 
\[\Endst \cn \permcatsu \to \pMulticat\]
in \cref{endsttwofunctor} induces an isomorphism of multimorphism categories
  \begin{align*}\label{eq:multilin-via-sma}
    \permcatsu\scmap{\ang{\C};\D} \fto[\iso]{\Endst}
    & \pMulticat\scmap{\ang{\Endst(\C)};\Endst(\D)}\\
    & = \pMulticat\big(\txsma_{j=1}^n \Endst(\C_j) \scs \Endst(\D)\big)
  \end{align*}
  between
\begin{itemize}
\item the category of $n$-linear functors and transformations $\ang{\C} \to \D$ and
\item the category of pointed multifunctors 
\[\txsma_{j=1}^n \Endst(\C_j) \to \Endst(\D)\]
and pointed multinatural transformations.
\end{itemize} 
Therefore, $\Endst$ is a $\Cat$-multifunctor.
\end{proposition}

\begin{explanation}[$\Endst$ as a $\Cat$-Multifunctor]\label{expl:endst-catmulti}
The $\Cat$-multifunctor $\Endst$ in \cref{proposition:n-lin-equiv} is given explicitly as follows.  

\emph{Object Assignment}.
It sends a small permutative category $\C$ to the pointed endomorphism multicategory 
\[\Endst(\C) = \big(\End(\C), i\big)\]
in \cref{ex:endstc}.

\emph{Multimorphism Functor on Objects}. 
An object in the $n$-ary multimorphism category $\permcatsu\scmap{\ang{\C};\D}$ is an $n$-linear functor (\cref{def:nlinearfunctor})
\[\begin{tikzcd}[column sep=2.3cm]
\txprod_{j=1}^n \C_j \ar{r}{\big(P,\, \{P^2_j\}_{j=1}^n\big)} & \D.
\end{tikzcd}\]
The image pointed multifunctor
\begin{equation}\label{EndstofP}
\Endst(P) \cn \txsma_{j=1}^n \Endst(\C_j) \to \Endst(\D)
\end{equation}
has object assignment induced by the object assignment of $P$.  This makes sense because, by the pushout \cref{eq:multicat-smash-pushout} that defines the smash product $\txsma_{j=1}^n \Endst(\C_j)$, its objects are represented by elements in
\[\Ob\left(\txotimes_{j=1}^n \End(\C_j)\right) 
= \txprod_{j=1}^n \Ob\big(\End(\C_j)\big) 
= \txprod_{j=1}^n \Ob(\C_j).\]
The object unity axiom \cref{nlinearunity} ensures that it descends to the objects of the smash product.

By the pushout \cref{eq:multicat-smash-pushout} again, each multimorphism in $\txsma_{j=1}^n \Endst(\C_j)$ is represented by a multimorphism in $\txotimes_{j=1}^n \End(\C_j)$.  By 
\begin{itemize}
\item \cref{explanation:unpacking-tensor} of the tensor product and
\item the definition of $\End$ (\cref{ex:endc}),
\end{itemize}  
the multimorphisms in $\txotimes_{j=1}^n \End(\C_j)$ are generated by
\begin{equation}\label{angccompjpsi}
c_1 \otimes \cdots \otimes c_{j-1} \otimes \psi \otimes c_{j+1} \otimes \cdots \otimes c_n
\end{equation}
for some $j \in \{1,\ldots,n\}$, objects $c_i \in \C_i$ for $i \neq j$, and multimorphism 
\[\psi \in \End(\C_j)\scmap{\ang{x_k}_{k=1}^r;y} 
= \C_j \brb{ \oplus_{k=1}^r x_k , y}.\]
We use the following notation.
\begin{equation}\label{eq:angccompj}
  \begin{aligned}
    \angc ~\compj~ ? &= \big(c_1, \ldots, c_{j-1}, ?, c_{j+1}, \ldots, c_n\big)\\
    1_{\angc} ~\compj~ ? &= \big(1_{c_1}, \ldots, 1_{c_{j-1}}, ?, 1_{c_{j+1}}, \ldots, 1_{c_n}\big)
  \end{aligned}
\end{equation} 
Then $\Endst(P)$ sends the multimorphism in \cref{angccompjpsi} to the following composite morphism in $\D$.
\begin{equation}\label{eq:EstP1compjpsi}
  \begin{tikzpicture}[xscale=3,yscale=1.3,vcenter]
    \draw[0cell=.9]
    (0,0) node (a) {\txoplus_{k=1}^r P\left(\angc \compj x_k\right)}
    (a)++(.5,-1) node (b) {P\left(\angc \compj \left(\txoplus_{k=1}^r x_k\right)\right)}
    (a)++(1,0) node (c) {P\left(\angc \compj y\right)}
    ;
    \draw[1cell=.9]  
    (a) edge node[swap,pos=.1] {P^2_j} (b)
    (b) edge node[swap,pos=.8] {P\left(1_{\angc} \compj \psi\right)} (c)
    ;
  \end{tikzpicture}
\end{equation}
This is an $r$-ary multimorphism in 
\[\End(\D)\lrscmap{\bang{ P\big(\angc \compj x_k\big) }_{k=1}^r ; P\big(\angc \compj y \big)}.\]
This assignment descends to the smash product by the unity axioms \cref{nlinearunity,constraintunity}.  These object and multimorphism assignments yield a pointed multifunctor $\Endst(P)$ as in \cref{EndstofP} by the other axioms of an $n$-linear functor (\cref{def:nlinearfunctor}).

\emph{Multimorphism Functor on Morphisms}. 
A morphism 
\[\theta \cn P \to Q \inspace \permcatsu\scmap{\ang{\C};\D}\]
is an $n$-linear transformation (\cref{def:nlineartransformation}).  A morphism in
\[\pMulticat\big(\txsma_{j=1}^n \Endst(\C_j) \scs \Endst(\D)\big)\]
is a pointed multinatural transformation (\cref{def:ptd-multicat}).  For each object $\angc \in \prod_{j=1}^n \C_j$, the component of 
\[\Endst(\theta) \cn \Endst(P) \to \Endst(Q)\]
at the object of $\txsma_{j=1}^n \Endst(\C_j)$ represented by $\angc$ is the component morphism
\begin{equation}\label{eq:Endst-theta-component}
  \theta_{\angc} \cn P\angc \to Q\angc \inspace \D.
\end{equation}
This defines a pointed multinatural transformation $\Endst(\theta)$ by the multilinearity conditions \cref{niitransformationunity,eq:monoidal-in-each-variable}.
\end{explanation}

The following result combines \cite[5.3.6, 5.3.9, 6.5.1, and 10.2.14]{cerberusIII}.

\begin{theorem}\label{endfactorization}\index{endomorphism!2-functor}\index{pointed!multicategory!forgetful functor}
There is a commutative diagram of $\Cat$-multifunctors
\begin{equation}\label{endufactor}
\begin{tikzpicture}[xscale=3.5,yscale=1.5,vcenter]
\draw[0cell=.85]
(0,0) node (a) {\permcatsu}
(a)++(1,0) node (b) {\Multicat}
(a)++(0,-1) node (c) {\MoneMod}
(c)++(1,0) node (d) {\pMulticat}
;
\draw[1cell=.9]  
(a) edge node {\End} (b)
(a) edge node {\Endst} (d)
(a) edge node[swap] {\Endm} (c)
(c) edge node[pos=.4] {\Um} (d)
(d) edge node[swap] {\Ust} (b)
;
\end{tikzpicture}
\end{equation}
defined as follows.
\begin{itemize}
\item $\Ust$, $\Um$, and $\Endst$ are the $\Cat$-multifunctors in \cref{expl:Ust-catmulti,expl:Um-catmulti,expl:endst-catmulti}, respectively.
\item $\End$ is the composite $\Cat$-multifunctor $\Ust \circ \Endst$, which restricts to the 2-functor in \cref{endtwofunctor}.
\item $\Endm$ is defined on objects in \cref{ex:endmc}.  It extends to a $\Cat$-multifunctor satisfying
\[\Um \circ \Endm = \Endst\]
by \cref{proposition:EM2-5-1} \cref{it:EM251-3,it:EM251-4,monetoendc} and \cref{proposition:n-lin-equiv}.
\end{itemize}
\end{theorem}

\begin{explanation}[$\End$ as a $\Cat$-Multifunctor]\label{expl:end-catmulti}
The $\Cat$-multifunctor\index{endomorphism!Cat-multifunctor@$\Cat$-multifunctor}
\[\End = \Ust \Endst \cn \permcatsu \to \Multicat\]
in \cref{endufactor} sends a small permutative category $\C$ to the endomorphism multicategory $\End(\C)$ in \cref{ex:endc}.  For small permutative categories $\angC$ and $\D$, the composite multimorphism functor
\[\begin{tikzpicture}[xscale=5,yscale=1.3,vcenter]
\draw[0cell=.85]
(0,0) node (a) {\permcatsu\scmap{\angC;\D}}
(a)++(.5,-1) node (b) {\pMulticat\brb{\txsma_{j=1}^n \Endst(\C_j) , \Endst(\D)}}
(a)++(1,0) node (c) {\Multicat\brb{\txotimes_{j=1}^n \End(\C_j) , \End(\D)}}
;
\draw[1cell=.9]  
(a) edge node {\End} (c)
(a) edge node[pos=.6] {\iso} node[swap,pos=.2] {\Endst} (b)
(b) edge node[swap,pos=.8] {\Ust} (c)
;
\end{tikzpicture}\]
is as described in \cref{expl:endst-catmulti} before descending to the smash product.
\end{explanation}

\begin{explanation}[$\Endm$ as a $\Cat$-Multifunctor]\label{expl:endm-catmulti}
The $\Cat$-multifunctor
\[\Endm \cn \permcatsu \to \MoneMod\]
in \cref{endufactor} sends a small permutative category $\C$ to the endomorphism left $\Mone$-module 
\[\Endm(\C) = \big(\Endst(\C),\mu\big)\] 
in \cref{ex:endmc}.  For small permutative categories $\angC$ and $\D$, the multimorphism functor $\Endm$ is the following composite isomorphism.
\[\begin{tikzpicture}[xscale=5,yscale=1.3,vcenter]
\draw[0cell=.85]
(0,0) node (a) {\permcatsu\scmap{\angC;\D}}
(a)++(.5,-1) node (b) {\MoneMod\brb{\txsma_{j=1}^n \Endm(\C_j) , \Endm(\D)}}
(a)++(1,0) node (c) {\pMulticat\brb{\txsma_{j=1}^n \Endst(\C_j) , \Endst(\D)}}
;
\draw[1cell=.9]  
(a) edge node {\Endst} node[swap] {\iso} (c)
(a) edge node[swap,pos=.15] {\Endm} (b)
(c) edge node[pos=.1] {(\Um)^{-1}} node[swap,pos=.6] {\iso} (b)
;
\end{tikzpicture}\]
\begin{itemize}
\item $\Endst$ is an isomorphism by \cref{proposition:n-lin-equiv}.
\item $\Um$ is an isomorphism by
\begin{itemize}
\item \cref{proposition:EM2-5-1} \eqref{it:EM251-3} and \eqref{it:EM251-4} if $n>0$ and
\item \cref{proposition:EM2-5-1} \eqref{monetoendc} if $n=0$.\defmark
\end{itemize} 
\end{itemize}
\end{explanation}

\chapter{Infinite Loop Space Machines}
\label{ch:Kspectra}
In this chapter we review two $K$-theory functors,
\begin{itemize}
\item $\Kse$ due to Segal \cite{segal} and
\item $\Kem$ due to Elmendorf-Mandell \cite{elmendorf-mandell,elmendorf-mandell-perm},
\end{itemize}  
from small permutative categories, $\permcatsu$, to connective symmetric spectra, $\Spc$:
\[
\]
Each of $\Kse$ and $\Kem$ is an equivalence of homotopy theories.  This implies that each of them induces an equivalence between the respective stable homotopy categories, which are obtained by inverting the respective classes of stable equivalences.  

Moreover, there are decompositions
\[\begin{split}
\Kse &= \Kf \circ~ \Ner_* \circ~ \Jse \andspace\\
\Kem &= \Kg \circ~ \Ner_* \circ~ \Jt \circ \Endm
\end{split}\]
as in the following diagram, which we explain in more detail in \cref{sec:segalEMK}.
\begin{equation}\label{KseKem-diagram}
%
\end{equation}  
Each category in \cref{KseKem-diagram} is an enriched symmetric monoidal category, except for $\permcatsu$, which is a $\Cat$-multicategory.  However, Segal $K$-theory, $\Kse$, is \emph{not} a multifunctor because $\Jse$ is not a multifunctor.  On the other hand, each constituent functor that comprises Elmendorf-Mandell $K$-theory, $\Kem$, is an enriched multifunctor, so $\Kem$ itself is an enriched multifunctor.  

Each functor in \cref{KseKem-diagram}, \emph{except} $\Jt$ and $\Kg$, is an equivalence of homotopy theories.  Each of the three pairs along the top row of \cref{KseKem-diagram},
\[(\cP, \Jse), \quad (S_*, \Ner_*), \andspace (\bA,\Kf),\]
induces mutually inverse equivalences between the respective stable homotopy categories.  Among the three homotopy inverses $\cP$, $S_*$, and $\bA$, only $\cP$ is a non-symmetric $\Cat$-multifunctor.  Each of $S_*$ and $\bA$ is incompatible with the multiplicative structures of its domain and codomain.

\subsection*{Connection with Other Chapters}\

\subsubsection*{Equivalences of Homotopy Theories}

We use equivalences of homotopy theories (\cref{sec:hty-thy}) in several subsequent chapters.
\begin{itemize}
\item In \cref{ch:multperm} we observe that there are equivalences of homotopy theories between small multicategories, $\Multicat$, and small permutative categories, $\permcatsu$.  Together with Segal $K$-theory, $\Kse$ (\cref{sec:segalEMK}), we obtain an equivalence of homotopy theories from $\Multicat$ to $\Spc$.
\item In \cref{ch:ptmulticat-sp,ch:ptmulticat-alg} we extend the equivalences of homotopy theories between $\Multicat$ and $\permcatsu$ first to \emph{pointed} multicategories, $\pMulticat$, and then further to left $\Mone$-modules, $\MoneMod$.  Together with Segal $K$-theory, this implies that there is an equivalence of homotopy theories from each of $\pMulticat$ and $\MoneMod$ to $\Spc$.
\end{itemize}

\subsubsection*{Enriched Mackey Functors}\

\begin{itemize}
\item In \cref{sec:factor-Kemse,sec:presheaf-K,sec:mult-mackey-spectra} we apply our general results about multicategorical standard enrichment, enriched diagrams, and enriched presheaves to the $\Cat$-multifunctors that constitute Elmendorf-Mandell $K$-theory, $\Kem$.  In particular, \cref{thm:Kemdg,gspectra-thm-xiv} prove that $\Kem$ induces a change-of-enrichment functor $\Kemdg$ that produces spectral Mackey functors from permutative Mackey functors based on any small category $\C$ enriched in $\permcatsu$.  Moreover, this functor factors through categories of enriched Mackey functors based on left $\Mone$-modules, $\Gstar$-categories, and $\Gstar$-simplicial sets.
\item In \cref{part:homotopy-mackey} we further extend the equivalences of homotopy theories between $\pMulticat$, $\MoneMod$, and $\permcatsu$ to their respective categories of enriched diagrams and enriched Mackey functors.  See \cref{mackey-xiv-pmulticat,mackey-xiv-mone,mackey-pmulti-mone}.
\end{itemize}

\subsection*{Background}

We use the notions of enriched (monoidal) categories and multicategories in \cref{sec:enrichedcat,sec:enrmonoidalcat,sec:enrmulticat}.  We also use the symmetric monoidal $\Cat$-category $\MoneMod$ and the $\Cat$-multicategory $\permcatsu$ in \cref{sec:monemodules,sec:multpermcat}, respectively.

\subsection*{Chapter Summary}

In \cref{sec:hty-thy} we review homotopy theories and their equivalences in the context of complete Segal spaces.  In \cref{sec:pointed-diagrams} we discuss pointed diagram categories and their symmetric monoidal closed structure.  In \cref{sec:Gamma-cat} we review $\Ga$-objects, which are pointed diagrams on the category $\Fskel$ of pointed finite sets.  In \cref{sec:Gstar-cat} we review $\Gstar$-objects, which are pointed diagrams on the category $\Gskel$ with finite tuples of pointed finite sets as objects.  In \cref{sec:segalEMK} we review the categories and functors that constitute Segal $K$-theory, its homotopy inverse, and Elmendorf-Mandell $K$-theory.  Here is a summary table.
\reftable{.9}{
\multicolumn{2}{|c|}{\secnumname{sec:hty-thy}} \\ \hline
complete Segal space model structure & \ref{theorem:css-fibrant} \\ \hline
relative categories, functors, and natural transformations & \ref{definition:rel-cat} \\ \hline
(criteria for) equivalences of homotopy theories & \ref{definition:rel-cat-pow} (\ref{gjo29}) \\ \hline\hline
\multicolumn{2}{|c|}{\secnumname{sec:pointed-diagrams}} \\ \hline
symmetric monoidal closed category of pointed objects & \ref{theorem:pC-sm-closed} \\ \hline
pointed unitary enrichment & \ref{def:unitary-enrichment} \\ \hline
pointed Day convolution and monoidal unit & \ref{definition:Dgm-pV-convolution-hom} \\ \hline
symmetric monoidal closed category of pointed diagrams & \ref{thm:Dgm-pv-convolution-hom} \\ \hline\hline
\multicolumn{2}{|c|}{\secnumname{sec:Gamma-cat}} \\ \hline
category $\Fskel$ and $\Ga$-objects & \ref{def:ordn} and \ref{def:Gamma-object} \\ \hline
symmetric monoidal closed category of $\Ga$-objects & \ref{expl:Gamma-object} \\ \hline\hline
\multicolumn{2}{|c|}{\secnumname{sec:Gstar-cat}} \\ \hline
category $\Gskel$ and $\Gstar$-objects & \ref{definition:Gstar}, \ref{definition:concatenation-product}, and \ref{def:Gstar-object} \\ \hline
symmetric monoidal closed category of $\Gstar$-objects & \ref{expl:Gstar-object} \\ \hline
functors $i \cn \Fskel \lradj \Gskel \cn \sma$ & \ref{i-incl} and \ref{smagf} \\ \hline\hline
\multicolumn{2}{|c|}{\secnumname{sec:segalEMK}} \\ \hline
Segal $K$-theory $\Kse$ & \ref{Kse} \\ \hline
homotopy inverses $\cP$, $S_*$, and $\bA$ & \ref{PinvK}, \ref{Sstar-mandell}, and \ref{A-segal} \\ \hline
Elmendorf-Mandell $K$-theory $\Kem$ & \ref{Kem} \\ \hline
}
Most of the material in this chapter is adapted from \cite[Chapters 4 and 8--10]{cerberusIII} and \cite{johnson-yau-invK,johnson-yau-multiK}.  We provide more references below.  We remind the reader of \cref{conv:universe,expl:leftbracketing}.

\section{Homotopy Theories via Complete Segal Spaces}\label{sec:hty-thy}

In this section we review equivalences of homotopy theories in terms of complete Segal spaces in the sense of \cite{rezk-homotopy-theory}.  See \cite{barwick-kan,dwyer-kan,toen-axiomatisation} for further development.  Practical criteria for checking equivalences of homotopy theories are discussed in \cite{gjo-extending,gjo1,johnson-yau-multiK}.  General references for model category theory are \cite{hirschhorn,hovey}.

\subsection*{Complete Segal Spaces}

The category \index{simplicial sets}$\sSet$ of simplicial sets is equipped with the standard \index{model structure!Kan}\index{Kan model structure}\emph{Kan model structure}.  The \index{nerve functor}\index{functor!nerve}\emph{nerve} functor is denoted
\begin{equation}\label{nerve}
\Ner \cn \Cat \to \sSet.
\end{equation}
See \cite[Section 7.2]{cerberusIII} for an elementary discussion of the nerve.  For $n \geq 0$, we denote by
\begin{equation}\label{ordn}
\ord{n} = \{0,1,\ldots,n\}
\end{equation}
the pointed finite set with $n+1$ elements and basepoint $0$.

\begin{definition}\label{def:bisimplicial}\
\begin{enumerate}
\item Denote by $\boldtwo$ the nerve of the category consisting of two isomorphic objects.
\item For $n \geq 2$ and $j \in \{1,\ldots,n\}$, the \emph{$j$-th characteristic map} is the pointed function
\[\chi_j \cn \ord{1} \to \ord{n} \stspace \chi_j(1) = j.\]
\item A \index{bisimplicial set}\emph{bisimplicial set} is a simplicial object in the category of simplicial sets.  
\item A bisimplicial set is \index{Reedy fibrant}\emph{Reedy fibrant} if it is a fibrant object in the Reedy model structure.
\item For a bisimplicial set $A$, the \index{Segal!morphism}\emph{$n$-th Segal morphism} is the simplicial map
\[A_n \to \overbracket[.5pt]{A_1 \times_{A_0} \cdots \times_{A_0} A_1}^{\text{$n$ copies of $A_1$}}\]
whose composite with the $j$-th coordinate projection is $A(\chi_j)$ for each $j \in \{1,\ldots,n\}$.
\item We say that a bisimplicial set \index{Segal!condition}\emph{satisfies the Segal condition} if the $n$-th Segal morphism is a weak equivalence of simplicial sets for each $n \geq 2$.\defmark
\end{enumerate}
\end{definition}

\begin{definition}\label{definition:css}
  A \index{complete Segal space}\index{Segal!complete - space}\emph{complete Segal space} is a bisimplicial set $A$ that satisfies the following three conditions.
\begin{description}
\item[Fibrancy] $A$ is Reedy fibrant.
\item[Segal Condition] $A$ satisfies the Segal condition.
\item[Path Condition] The morphism
\[A_0 \iso \Map(\De[0],A) \to \Map(\boldtwo,A)\]
induced by the unique morphism $\boldtwo \to \De[0]$ is a weak equivalence of simplicial sets.\defmark
\end{description}
\end{definition}

The following is \cite[Theorem 7.2]{rezk-homotopy-theory}.

\begin{theorem}\label{theorem:css-fibrant}
  There is a simplicial model structure on the category of bisimplicial sets, called the \index{model structure!complete Segal space}\index{complete Segal space!model structure}\index{Segal!complete - space!model structure}\emph{complete Segal space model structure}, that is given as a left Bousfield localization of the Reedy model structure and in which the fibrant objects are precisely the complete Segal spaces.
\end{theorem}

A weak equivalence in the complete Segal space model structure of bisimplicial sets is called a \index{weak equivalence!Rezk}\index{Rezk weak equivalence}\emph{Rezk weak equivalence}.

\subsection*{Homotopy Theories of Relative Categories}

A \index{subcategory!wide}\index{wide subcategory}\emph{wide subcategory} of a category $\C$ is a subcategory that contains all of the objects of $\C$.

\begin{definition}\label{definition:rel-cat}\
\begin{enumerate}
\item\label{relcat-i} A \index{category!relative}\index{relative!category}\emph{relative category} is a pair $(\C,\cW)$ consisting of
\begin{itemize}
\item a category $\C$ and
\item a wide subcategory $\cW$ of $\C$.
\end{itemize} 
We refer to morphisms in $\cW$ as \index{equivalence!stable}\index{stable equivalence}\emph{stable equivalences}.  If there is no danger of confusion, we denote a relative category $(\C,\cW)$ by $\C$.
\item\label{relcat-ii} A \index{functor!relative}\index{relative!functor}\emph{relative functor} between relative categories 
\[F \cn (\C,\cW) \to (\D,\cX)\]
is a functor $F \cn \C \to \D$ that restricts to a functor $\cW \to \cX$.
\item\label{relcat-iii} Suppose $(\C,\cW)$ and $(\D,\cX)$ are relative categories.  A functor $F \cn \C \to \D$ \emph{creates morphisms in $\cW$} if $\cW = F^\inv(\cX)$.
\item\label{relcat-iv} For relative functors
\[F, G \cn (\C,\cW) \to (\D,\cX),\]
a \index{natural transformation!relative}\index{relative!natural transformation}\emph{relative natural transformation}
\[\theta \cn F \to G\]
is a natural transformation such that each component of $\theta$ is a morphism in $\cX$.
\item\label{relcat-v} 
  A subclass $\cW \bigsubset \C$ of morphisms in a category $\C$ has the \index{2-out-of-3 property}\emph{2-out-of-3 property} if, for each pair of morphisms $f,g \in \C$ with $gf$ defined, whenever any two of 
\[f, \quad g, \andspace gf\]
belong to $\cW$, then so does the third morphism.
\item\label{relcat-vi} A \index{weak equivalence!category with}\index{category!with weak equivalences}\emph{category with weak equivalences} is a relative category $(\C, \cW)$ such that $\cW$
\begin{itemize}
\item contains all the isomorphisms and
\item has the 2-out-of-3 property.\defmark
\end{itemize}
\end{enumerate}
\end{definition}

For example, the class of weak equivalences in each model category contains all the isomorphisms and has the 2-out-of-3 property.

\begin{definition}\label{definition:rel-cat-pow}
Suppose $(\C,\cW)$ is a relative category.
\begin{enumerate}
\item\label{def:relcat-i} For a small category $\Dgm$, the \index{category!relative diagram -}\index{relative!diagram category}\emph{relative diagram category}
\[(\C,\cW)^\Dgm\]
is the wide subcategory of the diagram category $\C^\Dgm = [\Dgm,\C]$ in which a morphism is a natural transformation with each component in $\cW$.
\item\label{def:relcat-ii} The \index{diagram!classification}\index{classification diagram}\emph{classification diagram} of $(\C,\cW)$ is the bisimplicial set
\[\Nde(\C,\cW) = \Ner\big((\C,\cW)^{\De[?]}\big)\]
with $\De[n]$ denoting the category with $n$ composable arrows for $n \geq 0$.
\item\label{def:relcat-iii} A \index{relative!category!homotopy theory}\index{homotopy theory}\emph{homotopy theory of $(\C,\cW)$} is a fibrant replacement of $\Nde(\C,\cW)$ in the complete Segal space model structure.
\item\label{def:relcat-iv} A relative functor
\[F \cn (\C,\cW) \to (\D,\cX)\]
is an \index{equivalence!of homotopy theories}\index{homotopy theory!equivalence of}\emph{equivalence of homotopy theories} if the bisimplicial map $\R\big(\Nde(F)\big)$ is a \index{weak equivalence!Rezk}\index{Rezk weak equivalence}Rezk weak equivalence, where $\R$ denotes fibrant replacement in the complete Segal space model structure.  We sometimes denote an equivalence of homotopy theories by $\fto{\sim}$.\defmark
\end{enumerate} 
\end{definition}

\begin{definition}\label{def:inverse-heq}
Suppose given
\begin{itemize}
\item relative categories $(\C,\cW)$ and $(\D,\cX)$ and
\item functors 
\[L \cn \C \lradj \D \cn R.\]
\end{itemize}
Then we say that $L$ and $R$ are \index{equivalence!inverse - of homotopy theories}\index{homotopy theory!inverse equivalence of}\index{inverse equivalence!of homotopy theories}\emph{inverse equivalences of homotopy theories} if the following three conditions hold:
\begin{romenumerate}
\item\label{def:iheq-i} 
$L$ and $R$ are relative functors.
\item\label{def:iheq-ii}
$RL$ and $1_\C$ are connected by a zigzag of relative natural transformations.
\item\label{def:iheq-iii}
$LR$ and $1_\D$ are connected by a zigzag of relative natural transformations.
\end{romenumerate}
This finishes the definition.
\end{definition}

We emphasize that the functors $L$ and $R$ in \cref{def:inverse-heq} do \emph{not} have to form an adjunction.  The following result from \cite[Corollary 2.9]{gjo1} is our main tool for proving equivalences of homotopy theories.  

\begin{proposition}\label{gjo29}\index{homotopy theory!inverse equivalence of}\index{inverse equivalence!of homotopy theories}
Suppose given inverse equivalences of homotopy theories (\cref{def:inverse-heq})
\[L \cn (\C, \cW) \lradj (\D,\cX) \cn R.\]
Then $L$ and $R$ are equivalences of homotopy theories in the sense of \cref{definition:rel-cat-pow} \eqref{def:relcat-iv}.
\end{proposition}

The following definition from \cite[1.8]{gjo-extending} strengthens \cref{def:inverse-heq}.

\begin{definition}\label{def:heq}
Suppose $(\C,\cW)$ and $(\D,\cX)$ are relative categories, and suppose
\[L \cn \C \lradj \D \cn R\]
is an adjunction of categories with left adjoint $L$.  We call $L \dashv R$ an \index{equivalence!adjoint - of homotopy theories}\index{homotopy theory!adjoint equivalence of}\index{adjoint equivalence!of homotopy theories}\emph{adjoint equivalence of homotopy theories} if the following three conditions hold:
\begin{romenumerate}
\item\label{def:heq-i} $L$ and $R$ are relative functors.
\item\label{def:heq-ii} The unit of $L \dashv R$ is a relative natural transformation $1_\C \to RL$.
\item\label{def:heq-iii} The counit of $L \dashv R$ is a relative natural transformation $LR \to 1_\D$.
\end{romenumerate}
This finishes the definition.
\end{definition}

If $L \dashv R$ is an adjoint equivalence of homotopy theories, then $L$ and $R$ are inverse equivalences of homotopy theories in the sense of \cref{def:inverse-heq}.  \cref{gjo29} implies that $L$ and $R$ are equivalences of homotopy theories.

\section{Category of Pointed Diagrams}
\label{sec:pointed-diagrams}

To prepare for the discussion of $\Ga$-objects and $\Gstar$-objects in \cref{sec:Gamma-cat,sec:Gstar-cat}, in this section we review pointed diagrams.  The main results, \cref{thm:Dgm-pv-convolution-hom,theorem:diagram-omnibus}, state that the category of pointed functors and natural transformations is complete, cocomplete, symmetric monoidal closed, enriched, tensored, and cotensored.  The material in this section is adapted from \cite[Chapter 4]{cerberusIII}.

\subsection*{Symmetric Monoidal Closed Category of Pointed Objects}

\begin{definition}[Pointed Objects]\label{def:pointed-objects}
Suppose $\C$ is a category with a chosen terminal object $\term$.  We denote by $\pC$ the category under $\term$, which is defined as follows. 
\begin{itemize}
\item An object in $\pC$ is a pair $(a,i^a)$\label{not:aia} consisting of
\begin{itemize}
\item an object $a$ in $\C$ and
\item a morphism
\[i^a \cn \term \to a \in \C.\]
\end{itemize} 
We call $(a,i^a)$ a \index{object!pointed}\index{pointed!object}\emph{pointed object} with \index{pointed!structure}\emph{pointed structure} $i^a$.
\item For pointed objects $(a,i^a)$ and $(b,i^b)$, a morphism
\[f \cn (a,i^a) \to (b,i^b) \inspace \pC\]
is a morphism $f \cn a \to b$ in $\C$ such that the following diagram in $\C$ commutes.
\begin{equation}\label{ptmorphism-diagram}
\begin{tikzpicture}[xscale=1,yscale=1,vcenter]
\draw[0cell]
(0,0) node (t) {\term}
(t)++(2,.7) node (a) {a}
(t)++(2,-.7) node (b) {b}
;
\draw[1cell]
(t) edge node[pos=.6] {i^a} (a)
(t) edge node[swap,pos=.6] {i^b} (b)
(a) edge node {f} (b)
;
\end{tikzpicture}
\end{equation} 
We call a morphism in $\pC$ a \index{morphism!pointed}\index{pointed!morphism}\emph{pointed morphism}.
\end{itemize} 
This finishes the definition of $\pC$.
\end{definition}

Note that if $\C$ is complete and cocomplete, then so is $\pC$.

\begin{definition}[Smash Product and Pointed Hom]\label{def:wedge-smash-phom}
Suppose $(\C,\otimes,\tu,\Hom)$ is a complete and cocomplete symmetric monoidal closed category (\cref{def:closedcat}) with a chosen terminal object $\term$.  Suppose $(a,i^a)$ and $(b,i^b)$ are pointed objects.  We define the following pointed objects.
\begin{enumerate}
\item The \index{wedge}\emph{wedge} $a \wed b$ is the pushout in $\C$ of the span\label{not:awedgeb}
\[a \xleftarrow{i^a} \term \fto{i^b} b\]
with pointed structure given by the composite 
\[\term \fto{i^a} a \to a \wed b.\]
\item The \index{product!smash}\index{smash product}\emph{smash product} $a \sma b$ is the following pushout in $\C$.
\begin{equation}\label{eq:smash}
\begin{tikzpicture}[x=50mm,y=15mm,vcenter]
      \draw[0cell] 
      (0,0) node (a) {(a \otimes \term) \bincoprod (\term \otimes b)}
      (1,0) node (b) {a \otimes b}
      (0,-1) node (c) {\term}
      (1,-1) node (d) {a \sma b}
      ;
      \draw[1cell] 
      (a) edge node {(1_a \otimes i^b) \bincoprod (i^a \otimes 1_b)} (b)
      (c) edge node {i^{a \sma b}} (d)
      (a) edge node {} (c)
      (b) edge node {\pi_{a,b}} (d)
      ;
\end{tikzpicture}
\end{equation}
\item The \index{smash unit}\emph{smash unit} $\stu$ is the coproduct in $\C$
\begin{equation}\label{smash-unit-object}
\stu = \tu \bincoprod \term
\end{equation}
with pointed structure given by the inclusion of the $\term$ summand.
\item The \index{hom!pointed}\index{pointed!hom}\emph{pointed hom} $\pHom(a,b)$ is the following pullback in $\C$.
  \begin{equation}\label{eq:pHom-pullback}
  \begin{tikzpicture}[x=30mm,y=15mm,vcenter]
    \draw[0cell] 
    (0,0) node (a) {\pHom(a,b)}
    (1,0) node (b) {\term}
    (0,-1) node (c) {\Hom(a,b)}
    (1,-1) node (d) {\Hom(\term,b)}
    ;
    \draw[1cell] 
    (a) edge node {} (b)
    (c) edge node {(i^a)^*} (d)
    (a) edge node {} (c)
    (b) edge node {} (d)
    ;
  \end{tikzpicture}
  \end{equation}
The right vertical morphism in \cref{eq:pHom-pullback} is equal to the composite
\[\term \iso \Hom(a,\term) \fto{(i^b)_*} \Hom(a,b) \fto{(i^a)^*} \Hom(\term,b).\]
This induces the pointed structure of $\pHom(a,b)$.
\end{enumerate}
This finishes the definition.
\end{definition}

The following is \cite[4.20]{elmendorf-mandell-perm}, as presented in \cite[4.2.3]{cerberusIII}.

\begin{theorem}\label{theorem:pC-sm-closed}
In the context of \cref{def:wedge-smash-phom}, the quadruple
\[\left(\pC,\sma,\stu,\pHom\right)\]
is a complete and cocomplete symmetric monoidal closed category.
\end{theorem}

We mainly use \cref{theorem:pC-sm-closed} when $(\C,\otimes,\tu,\term)$ is
\begin{itemize}
\item $(\Cat,\times,\boldone,\boldone)$, the category of small categories (\cref{ex:cat}), and
\item $(\sSet,\times,*,*)$, the category of simplicial sets.
\end{itemize}

\subsection*{Pointed Unitary Enrichment}

\begin{definition}[Zero Objects and Zero Morphisms]\label{definition:zero}
  A \index{object!zero}\index{zero!object}\emph{zero object} in a category $\C$ is an object $\zob$ that is both initial and terminal.  Moreover, for a given zero object $\zob$ in $\C$, we define the following.
\begin{itemize}
\item A \index{morphism!zero}\index{zero!morphism}\emph{zero morphism} in $\C$ is a morphism that factors through $\zob$.  In other words, for objects $x,y \in \C$, the zero morphism is the composite
\[x \to \zob \to y.\]
\item A \emph{nonzero morphism} is a morphism that does not factor through $\zob$.
\item For objects $x,y \in \C$, we denote by 
\begin{equation}\label{Cpuncxy}
\C^\punc(x,y) = \C(x,y) \setminus \{\zob\}
\end{equation}
the set of nonzero morphisms $x \to y$.\defmark
\end{itemize}
\end{definition}

\begin{definition}\label{def:null-object}
Suppose $(\Dgm,\Dtimes,\pu)$ is a symmetric monoidal category.  A zero object $\termd$ in $\Dgm$ is called a \index{null object}\index{object!null}\emph{null object} if there are natural isomorphisms
\[a \Dtimes \termd \iso \termd \iso \termd \Dtimes \; a\]
for objects $a$ in $\Dgm$.
\end{definition}

The following definition uses \cref{theorem:pC-sm-closed} on $\V$ and categories enriched in $\pV$ as in \cref{def:enriched-category}.

\begin{definition}[Pointed Unitary Enrichment]\label{def:unitary-enrichment}
Suppose
\begin{itemize}
\item $(\Dgm,\Dtimes,\pu,\termd)$ is a small symmetric monoidal category with a chosen null object $\termd$ and
\item $(\V,\otimes,\tu,[,],\termv)$ is a complete and cocomplete symmetric monoidal closed category with a chosen terminal object $\termv$. 
\end{itemize} 
The \index{pointed!unitary enrichment}\index{unitary enrichment!pointed}\emph{pointed unitary enrichment} of $\Dgm$ over $\pV$, denoted $\Dhat$, is the $\pV$-category with
\begin{itemize}
\item the same class of objects as $\Dgm$ and
\item for any pair of objects $a,b \in \Dgm$, the morphism object
\begin{equation}\label{eq:Dhat-ptd-unitary-enr}
\Dhat(a,b) = \bigvee_{f \in \Dgm^\punc(a,b)} \stu \inspace \pV
\end{equation}
with the notation as follows.
\begin{itemize}
\item The wedge and the smash unit, $\stu = \tu \bincoprod \termv$, are as in \cref{def:wedge-smash-phom}.
\item An empty wedge means $\termv \in \pV$.  
\item $\Dgm^\punc(a,b)$ is the set of nonzero morphisms as in \cref{Cpuncxy}.
\end{itemize}
\end{itemize}
Moreover, we denote by
\[\Dhat \sma \Dhat\]
the tensor product $\pV$-category as in \cref{definition:vtensor-0}.
\end{definition}

By \cite[2.4.10]{cerberusIII}, the symmetric monoidal structure on $\Dgm$ induces a symmetric monoidal $\pV$-category structure on $\Dhat$ (\cref{definition:symm-monoidal-vcat}).

\subsection*{Pointed Diagrams}

\begin{definition}[Pointed Diagram Categories]\label{def:pointed-category}\
\begin{enumerate}
\item A \index{category!pointed}\index{pointed!category}\emph{pointed category} is a pair $(\C,*)$ consisting of
\begin{itemize}
\item a category $\C$ and
\item a chosen object $* \in \C$, which is called the \index{basepoint}\emph{basepoint}.
\end{itemize} 
\item For pointed categories $(\C,*)$ and $(\D,*)$, a \index{functor!pointed}\index{pointed!functor}\emph{pointed functor}
\[F \cn (\C,*) \to (\D,*)\]
is a functor $F \cn \C \to \D$ such that $F(*) = *$.  
\item For pointed functors $F,G \cn (\C,*) \to (\D,*)$, a \index{natural transformation!pointed}\index{pointed!natural transformation}\emph{pointed natural transformation}
\[\theta \cn F \to G\]
is a natural transformation such that
\[\theta_* = 1_* \cn F(*) = * \to * = G(*) \inspace \D.\]
\end{enumerate}
Moreover, suppose $\C$ is small.  We define the \index{diagram!category!pointed}\index{pointed!diagram category}\emph{pointed diagram category} $\CstarD$ with
\begin{itemize}
\item pointed functors $(\C,*) \to (\D,*)$ as objects,
\item pointed natural transformations between them as morphisms,
\item identity functors as identity morphisms, and
\item vertical composition of natural transformations as composition.
\end{itemize}
This finishes the definition.
\end{definition}

In what follows, as in \cref{def:unitary-enrichment}, an empty wedge means the chosen terminal object $\termv \in \pV$.
  
\begin{definition}[Pointed Day Convolution]\label{definition:Dgm-pV-convolution-hom}
In the context of \cref{def:unitary-enrichment}, suppose given pointed functors
\[A, B \cn (\Dgm,\termd) \to (\pV,\termv).\]
We define the following pointed functors $(\Dgm,\termd) \to (\pV,\termv)$.
\begin{enumerate}
\item The \index{pointed!diagrams!unit diagram}\index{unit diagram!pointed}\emph{monoidal unit diagram} is the pointed functor
\begin{equation}\label{eq:Dgm-pV-unit}
\du = \bigvee_{\Dgm^\punc(\pu,-)} \stu \cn (\Dgm,\termd) \to (\pV,\termv)
\end{equation}
with $\Dgm^\punc(-,-)$ the set of nonzero morphisms as in \cref{Cpuncxy}.
\item The \index{pointed!diagrams!Day convolution}\index{Day convolution!pointed diagrams}\emph{pointed Day convolution} is the $\pV$-coend
\begin{equation}\label{eq:Dgm-pV-convolution}
A \sma B = \ecint^{(a,b) \in \Dhat \sma \Dhat}
\bigvee_{\Dgm^\punc(a \Dtimes b, -)} (Aa \sma Bb)
\end{equation}
with $\Dhat$ the pointed unitary enrichment in \cref{def:unitary-enrichment}.  If the input object is $\termd$ in \cref{eq:Dgm-pV-convolution}, we choose $\termv$ for the coend.
\item The \index{pointed!diagrams!hom diagram}\index{hom diagram!pointed diagrams}\emph{pointed hom diagram} is the $\pV$-end
\begin{equation}\label{eq:Dgm-pV-hom}
\Homdstar(A,B) = \ecint_{b \in \Dhat}\, \big[Ab \scs B(- \Dtimes b)\big]_*
\end{equation}
where $[,]_*$ is the pointed hom \cref{eq:pHom-pullback} for $\pV$.  If the input object is $\termd$ in \cref{eq:Dgm-pV-hom}, we choose $\termv$ for the end.
\item The \index{mapping object!pointed diagrams}\index{pointed!diagrams!mapping object}\emph{pointed mapping object} is the $\pV$-end
\begin{equation}\label{eq:Dgm-pV-map}
\Mapdstar(A,B) = \ecint_{b \in \Dhat}\, [Ab, Bb]_*
\iso \left(\Homdstar(A,B)\right)(\pu).
\end{equation}
\end{enumerate}
Each of \cref{eq:Dgm-pV-convolution,eq:Dgm-pV-hom,eq:Dgm-pV-map} extends componentwise to pointed natural transformations.
\end{definition}

\begin{theorem}\label{thm:Dgm-pv-convolution-hom}
In the context of \cref{def:unitary-enrichment,def:pointed-category,definition:Dgm-pV-convolution-hom}, the quadruple
\[\left( \DstarV, \sma, \du, \Homdstar \right)\]
is a \index{pointed!diagram category!complete and cocomplete}complete and cocomplete \index{pointed!diagram category!symmetric monoidal closed}symmetric monoidal closed category.
\end{theorem}

It follows from \cref{thm:Dgm-pv-convolution-hom,theorem:v-closed-v-sm} that $\DstarV$ is a symmetric monoidal $(\DstarV)$-category.

\begin{definition}\label{def:evale}
In the context of \cref{thm:Dgm-pv-convolution-hom}, evaluation at the monoidal unit of $\Dgm$ defines a symmetric monoidal functor\index{evaluation!at $\pu$}
\[\ev_{\pu} \cn \DstarV \to \pV.\]  
It admits a strong symmetric monoidal left adjoint, denoted \index{evaluation!at $\pu$!left adjoint}$L_{\pu}$.  
\end{definition}

\begin{theorem}\label{theorem:diagram-omnibus}
In the context of \cref{thm:Dgm-pv-convolution-hom,def:evale}, the adjunction 
\begin{equation}\label{eq:pV-DstarV-adj}
  \begin{tikzpicture}[xscale=2.5,yscale=1,baseline={(x.base)}]
  \def\a{35}
    \draw[0cell] 
    (0,0) node (x) {\pV}
    (1,0) node (y) {\DstarV}
    (.45,0) node (b) {\bot}
    ;
    \draw[1cell]
    (x) edge[bend left=\a] node {L_{\pu}} (y) 
    (y) edge[bend left=\a] node {\ev_{\pu}} (x) 
    ;
  \end{tikzpicture}
\end{equation}
makes the pointed diagram category $\DstarV$ \index{pointed!diagram category!enriched}enriched, \index{pointed!diagram category!tensored and cotensored}tensored, and cotensored over $\pV$, with mapping objects given by $\Mapdstar$ in \cref{eq:Dgm-pV-map}.  In particular, $\DstarV$ is
\begin{itemize}
\item a symmetric monoidal $\V$-category (\cref{definition:symm-monoidal-vcat}) and
\item a $\V$-multicategory.
\end{itemize} 
\end{theorem}

\begin{explanation}\label{expl:DstarV-enr-multicat}
In \cref{theorem:diagram-omnibus}, the assertion about symmetric monoidal $\V$-category follows from
\begin{itemize}
\item the fact that $\DstarV$ is a symmetric monoidal $(\DstarV)$-category (\cref{thm:Dgm-pv-convolution-hom,theorem:v-closed-v-sm}) and
\item \cref{VseU} applied to the symmetric monoidal functor $\ev_{\pu}$ (\cref{def:evale}).
\end{itemize} 
The assertion about $\V$-multicategory follows from \cref{proposition:monoidal-v-cat-v-multicat}.
\end{explanation}

\section{\texorpdfstring{$\Gamma$}{Gamma}-Objects}
\label{sec:Gamma-cat}

In this section we review the symmetric monoidal closed category of $\Ga$-objects.  The material in this section is adapted from \cite[Section 8.1]{cerberusIII}.

Recall from \cref{def:symmoncat} that a \index{category!permutative}\index{permutative category}\emph{permutative category} is a strict symmetric monoidal category.
Recall from \cref{def:ordn} the permutative category
\[\big(\Fskel, \sma, \ord{1}, \xi\big)\]
of \index{finite sets!pointed}\index{pointed!finite sets}pointed finite sets and pointed functions with the smash product as the monoidal product.
Note that $\ord{0} \in \Fskel$ is a null object (\cref{def:null-object}).  For $\V = \sSet$, the following definition is due to Segal \cite{segal}.

\begin{definition}[$\Ga$-Objects]\label{def:Gamma-object}
  Suppose $(\V,\otimes,\tu,[,],\termv)$ is a complete and cocomplete symmetric monoidal closed category with a chosen terminal object $\termv$.  With the small pointed category $(\Fskel,\ord{0})$, we define the \index{object!9Gamma@$\Ga$-}\index{99Gamma@$\Ga$!-object}\index{category!of $\Ga$-objects}\emph{category of $\Ga$-objects} in $\V$ as the pointed diagram category\label{not:GaV}
\[\GaV = \FstarV\]
as in \cref{def:pointed-category}.  Moreover, we define the following.
\begin{itemize}
\item With
\[\big(\V,\otimes,\tu,\termv\big) = \big(\Cat, \times, \boldone, \boldone\big),\]
we call an object in $\Gacat$ a \index{category!9Gamma@$\Ga$-}\index{99Gamma@$\Ga$!-category}\emph{$\Ga$-category}.
\item With
\[\big(\V,\otimes,\tu,\termv\big) = \big(\sSet, \times, *, *\big),\]
we call an object in $\Gasset$ a \index{simplicial set!9Gamma@$\Ga$-}\index{99Gamma@$\Ga$!-simplicial set}\emph{$\Ga$-simplicial set}.\defmark
\end{itemize}
\end{definition}

\begin{explanation}[Symmetric Monoidal Closed Structure]\label{expl:Gamma-object}
A $\Ga$-object in $\V$ is a pointed functor
\[(\Fskel,\ord{0}) \to (\pV,\termv).\]
Morphisms of $\Ga$-objects are pointed natural transformations between such pointed functors.  With the indexing permutative category
\[\big(\Dgm, \Dtimes, \pu, \termd\big) = \big(\Fskel, \sma, \ord{1}, \ord{0}\big),\]
by \cref{thm:Dgm-pv-convolution-hom} there is a complete and cocomplete symmetric monoidal closed category
\begin{equation}\label{GammaV}
\big(\GaV, \sma, \du, \Homfstar\!\big).
\end{equation}
Moreover, by \cref{theorem:diagram-omnibus}, $\GaV$ is enriched, tensored, and cotensored over $\pV$.  In particular, $\GaV$ is
\begin{itemize}
\item a symmetric monoidal $\V$-category (\cref{definition:symm-monoidal-vcat}) and
\item a $\V$-multicategory (\cref{proposition:monoidal-v-cat-v-multicat}).
\end{itemize} 
Therefore, 
\begin{itemize}
\item $\Gacat$ is a $\Cat$-multicategory, and
\item $\Gasset$ is an $\sSet$-multicategory.
\end{itemize}
 
Specifying \cref{definition:Dgm-pV-convolution-hom} to the case $\Dgm = \Fskel$, the symmetric monoidal closed structure and $\pV$-enrichment of $\GaV$ is given as follows.  An empty wedge means the chosen terminal object $\termv \in \pV$.
\begin{enumerate}
\item The \index{pointed!diagrams!unit diagram}\index{unit diagram!pointed}\emph{monoidal unit diagram} is the pointed functor
\begin{equation}\label{eq:GaV-unit}
\du = \bigvee_{\Fskel^\punc(\ord{1},-)} \stu \cn (\Fskel,\ord{0}) \to (\pV,\termv).
\end{equation}
\item The \index{pointed!diagrams!Day convolution}\index{Day convolution!pointed diagrams}\emph{pointed Day convolution} is the $\pV$-coend
\begin{equation}\label{eq:GaV-convolution}
A \sma B = \ecint^{(\ord{m},\ord{n}) \in \Fhat \sma \Fhat}
\bigvee_{\Fskel^\punc(\ord{m} \sma \ord{n}, -)} (A\ord{m} \sma B\ord{n}).
\end{equation}
If the input object is $\ord{0}$, we choose $\termv$ for the coend.
\item The \index{pointed!diagrams!hom diagram}\index{hom diagram!pointed diagrams}\emph{pointed hom diagram} is the $\pV$-end
\begin{equation}\label{eq:GaV-hom}
\Homfstar(A,B) = \ecint_{\ord{n} \in \Fhat}\, \big[A\ord{n} \scs B(- \sma \ord{n})\big]_*.
\end{equation}
If the input object is $\ord{0}$, we choose $\termv$ for the end.
\item The \index{mapping object!pointed diagrams}\index{pointed!diagrams!mapping object}\emph{pointed mapping object} is the $\pV$-end
\begin{equation}\label{eq:GaV-map}
\Mapfstar(A,B) = \ecint_{\ord{n} \in \Fhat}\, [A\ord{n}, B\ord{n}]_*
\iso \left(\Homfstar(A,B)\right)(\ord{1}).
\end{equation}
\end{enumerate}
To understand \cref{eq:GaV-unit}, we observe that there is a canonical bijection
\[\Fskel^\punc(\ord{1}, \ord{n}) \iso \ord{n}^\punc = \{1,\ldots,n\}\]
for each $n \geq 0$, where $\Fskel^\punc(\ord{1},\ord{0}) = \emptyset$.
\end{explanation}

\section{\texorpdfstring{$\Gstar$}{Gstar}-Objects}
\label{sec:Gstar-cat}

In this section we review the symmetric monoidal closed category of $\Gstar$-objects.  The material in this section is adapted from \cite[Sections 9.1 and 9.2]{cerberusIII}.  

\subsection*{Smash Powers of $\Fskel$}

The indexing category $\Gskel$ involves the category $\Fskel$ in \cref{def:ordn}.  To define $\Gskel$, first we need some preliminary definitions.

\begin{definition}[Injections]\label{definition:ufs}
We define the category $\Inj$ as follows.
\begin{itemize}
\item Its objects are \index{unpointed finite sets}\index{finite sets!unpointed}\emph{unpointed finite sets}
\begin{equation}\label{ufsn}
\ufs{n} = 
\begin{cases}
\{1, \ldots, n\} & \ifspace n > 0,\\
\varnothing & \ifspace n = 0,
\end{cases}
\end{equation}
for $n \ge 0$.
\item Its morphisms are \index{injection}injections.
\end{itemize} 
Suppose $f \cn \ufs{q} \hookrightarrow \ufs{p}$ is an injection.  We define a functor
\[f_*\cn \Fskel^q \to \Fskel^p\]
called the \emph{reindexing injection} as follows.  Suppose given $q$-tuples of pointed finite sets or pointed functions
\begin{equation}\label{ordtunpsi}
\ordtu{n} = \ang{\ord{n}_k}_{k=1}^q
\orspace
\ang{\psi} = \ang{\psi_k}_{k=1}^q
\in \Fskel^q,
\end{equation}
respectively.  We define\label{not:fstarordtun}
\[f_*\ordtu{n} = \bang{\ord{n}_{f^\inv(j)}}_{j=1}^p \andspace
f_*\ang{\psi} = \bang{\psi_{f^\inv(j)}}_{j=1}^p \in \Fskel^p,\]
where
\[\ord{n}_\varnothing = \ord{1} \andspace
\psi_{\varnothing} = 1_{\ord{1}}.\]
This finishes the definition.
\end{definition}

\begin{definition}[Smash Powers of $\Fskel$]\label{definition:Fskel-roundsma}
For $q \geq 0$, we define the pointed category $\Fskel^{(q)}$, called the \emph{$q$-th smash power of $\Fskel$}, as follows.
\begin{description}
\item[The Case $q = 0$] 
We define\label{not:ObFzero}
\[\Ob\,\Fskel^{(0)} = \big\{\gst, \ang{}\big\},\]
which consists of the basepoint object $\gst$ and the empty tuple $\ang{}$.  We define the morphisms of $\Fskel^{(0)}$ such that
\begin{itemize}
\item $\gst$ is both initial and terminal and
\item the only nonzero morphism is the identity morphism of $\ang{}$.
\end{itemize} 
\item[The Case $q > 0$]
With $\Fskel$ having basepoint $\ord{0}$, we define the $q$-fold smash power of pointed categories
\[\Fskel^{(q)} = \Fskel^{\sma q}\]
as in \cref{eq:smash} applied to
\[(\C, \otimes, \tu) = (\Cat, \times, \boldone).\]
We denote the objects and morphisms of $\Fskel^{(q)}$ as $q$-tuples as in \cref{ordtunpsi}. 
\begin{itemize}
\item A $q$-tuple $\ordtu{n}$ is identified with the basepoint of $\Fskel^{(q)}$ if any $\ord{n}_k = \ord{0}$.  We call $q$ the \emph{length} of a $q$-tuple $\ordtu{n}$.
\item A $q$-tuple $\ang{\psi}$ is a \index{morphism!zero}\index{zero!morphism}\emph{zero morphism} if any $\psi_k$ factors through $\ord{0} \in \Fskel$.
\end{itemize} 
\end{description}
This finishes the definition.
\end{definition}

\subsection*{The Permutative Category $\Gskel$}

Next we define the indexing category $\Gskel$.

\begin{definition}[Tuples of Pointed Finite Sets]\label{definition:Gstar}
We define a small pointed category\index{99G@$\Gskel$}
\[(\Gskel,\gst)\]
as follows.
\begin{description}
\item[Objects] 
The set of objects is the wedge of pointed sets
\[\Ob\,\Gskel = \bigvee_{q \ge 0} \Ob \big(\Fskel^{(q)}\big)\]
as in \cref{def:wedge-smash-phom} applied to $(\C,\term) = (\Set, *)$.  The basepoint object $\gst$ is both initial and terminal in $\Gskel$.  
\item[Morphisms] 
For a pair of objects
\[\ordtu{n} = \ang{\ord{n}_k}_{k=1}^q \andspace 
\ordtu{m} = \ang{\ord{m}_j}_{j=1}^p,\]
the set of morphisms in $\Gskel$ is
\begin{align}
\Gskel\brb{\ordtu{n},\ordtu{m}} \label{eq:Gskelnm}
    & = \bigvee_{f\in \Inj(\ufs{q}\,,\,\ufs{p})}\;\Big(\,
    \Fskel^{(p)}\brb{f_*\ordtu{n}, \ordtu{m}}
    \,\Big)\\
    & = \bigvee_{f\in \Inj(\ufs{q}\,,\,\ufs{p})}\;\Big(\,
    \Opsma_{j=1}^p \Fskel\brb{\ord{n}_{f^\inv(j)}, \ord{m}_j}
    \,\Big). \nonumber
\end{align}
In \eqref{eq:Gskelnm} for $p > 0$, we denote a morphism by a pair $(f,\ang{\psi})$ with
\[\begin{tikzcd}
\ufs{q} \ar{r}{f} & \ufs{p} \in \Inj
\end{tikzcd}
\andspace
\begin{tikzcd}
f_*\ordtu{n} \ar{r}{\ang{\psi}} & \ordtu{m} \in \Fskel^{(p)}.
\end{tikzcd}\]
A morphism $( f, \ang{\psi} )$ is identified with the \index{morphism!zero}\index{zero!morphism}\emph{zero morphism} in $\Gskel(\ordtu{n}, \ordtu{m})$ if there exists a component morphism
\[\psi_j \cn \ord{n}_{f^\inv (j)} \to \ord{m}_j \forsomespace j \in \{1,\ldots,p\}\]
that factors through $\ord{0} \in \Fskel$.
\item[Identities]
The identity morphism of a $q$-tuple $\ordtu{n}$ is the pair $(1_{\ufs{q}},1_{\ordtu{n}})$.
\item[Composition] 
The composite of morphisms
\[\ordtu{n} \fto{(f,\ang{\psi})} \ordtu{m} \fto{(g,\ang{\phi})} \ang{\ord{\ell}}\]
in $\Gskel$ is the pair 
\[\big(gf, \ang{\phi}\circ g_*\ang{\psi}\big).\]
\end{description}
This finishes the definition of $(\Gskel,\gst)$.
\end{definition}

\begin{definition}[Permutative Structure on $\Gskel$]\label{definition:concatenation-product}
We define a permutative category
\[\big(\Gskel, \oplus, \ang{}, \xi\big)\]
as follows.
\begin{description}
\item[Monoidal Unit] 
It is the empty tuple $\ang{}$.
\item[Monoidal Product] 
The \index{concatenation!product}\index{product!concatenation}\emph{concatenation product}
\[\Gskel \times \Gskel \fto{\oplus} \Gskel\]
is the concatenation of tuples
\[\ordtu{n} \oplus \ordtu{n'} \andspace \ang{\psi} \oplus \ang{\psi'}\]
for tuples of pointed finite sets and pointed functions, respectively. 
\begin{itemize}
\item The basepoint $\gst \in \Gskel$ is defined as a null object for $\oplus$ (\cref{def:null-object}).
\item The concatenation product of any morphism with a morphism from, respectively to, $\gst$ is uniquely determined because $\gst$ is a null object.  
\end{itemize} 
Given morphisms
\[(f,\ang{\psi}) \in \Gskel\brb{\ordtu{n},\ordtu{m}}
\andspace
(f',\ang{\psi'}) \in \Gskel\brb{\ang{\ord{n}'},\ang{\ord{m}'}},\]
the concatenation of injections
\[f\cn \ufs{q} \hookrightarrow \ufs{p} \andspace
f'\cn \ufs{q'} \hookrightarrow \ufs{p'}\]
is the injection
\[(f \oplus f')(i) = \begin{cases}
    f(i) & \forspace i \in \{1,\ldots,q\}\\
    p+f'(i-q) & \forspace i \in \{q+1,\ldots,q+q'\}.
  \end{cases}\]
Then we define the morphism
\[(f,\ang{\psi}) \oplus (f',\ang{\psi'}) =
\brb{f \oplus f', \ang{\psi} \oplus \ang{\psi'}}.\]
\item[Symmetry]
For non-basepoint objects $\ordtu{n}$ and $\ordtu{n'}$ in $\Gskel$, the symmetry component
\begin{equation}\label{eq:Gskel-xi}
\xi_{\ordtu{n},\ordtu{n'}}\cn \ordtu{n} \oplus \ordtu{n'} \fto{\iso} \ordtu{n'} \oplus \ordtu{n}
\end{equation}
is the pair $(\tau_{q,q'},1)$ defined as follows.
\begin{itemize}
\item The first entry
\begin{equation}\label{tauqqprime}
\tau_{q,q'}\cn \ufs{q+q'} \fto{\iso} \ufs{q'+q}
\end{equation}
is the block permutation that swaps the first $q$ elements with the last $q'$ elements.
\item The second entry is the identity morphism on
\[(\tau_{q,q'})_* \big(\ordtu{n}\oplus\ordtu{n'}\big) =
\ordtu{n'}\oplus\ordtu{n}.\]
\end{itemize}
Each component of $\xi$ involving the null basepoint $\gst$ is $1_\gst$. 
\end{description}
This finishes the definition of the permutative structure on $\Gskel$.
\end{definition}

\subsection*{$\Gstar$-Objects}

For $\V = \Cat$ and $\sSet$, the following definition is due to Elmendorf and Mandell \cite{elmendorf-mandell}.

\begin{definition}[$\Gstar$-Objects]\label{def:Gstar-object}
  Suppose $(\V,\otimes,\tu,[,],\termv)$ is a complete and cocomplete symmetric monoidal closed category with a chosen terminal object $\termv$.  With the small pointed category $(\Gskel,\gst)$, we define the \index{99G@$\Gstar$!-objects}\index{category!of $\Gstar$-objects}\emph{category of $\Gstar$-objects} in $\V$ as the pointed diagram category\label{not:GstarV}
\[\GstarV\]
as in \cref{def:pointed-category}.  Moreover, we define the following.
\begin{itemize}
\item With
\[\big(\V,\otimes,\tu,\termv\big) = \big(\Cat, \times, \boldone, \boldone\big),\]
we call an object in $\Gstarcat$ a \index{category!99G@$\Gstar$-}\index{99G@$\Gstar$!-category}\emph{$\Gstar$-category}.
\item With
\[\big(\V,\otimes,\tu,\termv\big) = \big(\sSet, \times, *, *\big),\]
we call an object in $\Gstarsset$ a \index{simplicial set!99G@$\Gstar$-}\index{99G@$\Gstar$!-simplicial set}\emph{$\Gstar$-simplicial set}.\defmark
\end{itemize}
\end{definition}

\begin{explanation}[Symmetric Monoidal Closed Structure]\label{expl:Gstar-object}
A $\Gstar$-object in $\V$ is a pointed functor
\[(\Gskel,\gst) \to (\pV,\termv).\]
Morphisms of $\Gstar$-objects are pointed natural transformations between such pointed functors.  With the indexing permutative category
\[\big(\Dgm, \Dtimes, \pu, \termd\big) = \big(\Gskel, \oplus, \ang{}, \gst\big),\]
by \cref{thm:Dgm-pv-convolution-hom} there is a complete and cocomplete symmetric monoidal closed category
\begin{equation}\label{Gstar-V}
\big(\GstarV, \sma, \du, \Homgstar\!\big).
\end{equation}
Moreover, by \cref{theorem:diagram-omnibus}, $\GstarV$ is enriched, tensored, and cotensored over $\pV$.  In particular, $\GstarV$ is
\begin{itemize}
\item a symmetric monoidal $\V$-category (\cref{definition:symm-monoidal-vcat}) and
\item a $\V$-multicategory (\cref{proposition:monoidal-v-cat-v-multicat}).
\end{itemize} 
Therefore, 
\begin{itemize}
\item $\Gstarcat$ is a $\Cat$-multicategory, and
\item $\Gstarsset$ is an $\sSet$-multicategory.
\end{itemize}
 
Specifying \cref{definition:Dgm-pV-convolution-hom} to the case $\Dgm = \Gskel$, the symmetric monoidal closed structure and $\pV$-enrichment of $\GstarV$ is given as follows.  An empty wedge means the chosen terminal object $\termv \in \pV$.
\begin{enumerate}
\item The \index{pointed!diagrams!unit diagram}\index{unit diagram!pointed}\emph{monoidal unit diagram} is the pointed functor
\begin{equation}\label{eq:GstarV-unit}
\du = \bigvee_{\Gskel^\punc(\ang{},-)} \stu \cn (\Gskel,\gst) \to (\pV,\termv).
\end{equation}
\item The \index{pointed!diagrams!Day convolution}\index{Day convolution!pointed diagrams}\emph{pointed Day convolution} is the $\pV$-coend
\begin{equation}\label{eq:GstarV-convolution}
A \sma B = \ecint^{(\ordtu{m},\ordtu{n}) \in \Ghat \sma \Ghat}
\bigvee_{\Gskel^\punc(\ordtu{m} \oplus \ordtu{n}, -)} \big(A\ordtu{m} \sma B\ordtu{n}\big).
\end{equation}
If the input object is $\gst$, we choose $\termv$ for the coend.
\item The \index{pointed!diagrams!hom diagram}\index{hom diagram!pointed diagrams}\emph{pointed hom diagram} is the $\pV$-end
\begin{equation}\label{eq:GstarV-hom}
\Homgstar(A,B) = \ecint_{\ordtu{n} \in \Ghat}\, \big[A\ordtu{n} \scs B(- \oplus \ordtu{n})\big]_*.
\end{equation}
If the input object is $\gst$, we choose $\termv$ for the end.
\item The \index{mapping object!pointed diagrams}\index{pointed!diagrams!mapping object}\emph{pointed mapping object} is the $\pV$-end
\begin{equation}\label{eq:GstarV-map}
\Mapgstar(A,B) = \ecint_{\ordtu{n} \in \Ghat}\, \big[A\ordtu{n}, B\ordtu{n}\big]_*
\iso \left(\Homgstar(A,B)\right)\ang{}.
\end{equation}
\end{enumerate}
To understand \cref{eq:GstarV-unit}, we observe that there are canonical bijections
\[\Gskel^\punc\big(\ang{}, \ordtu{n}\big)
\iso \txprod_{k=1}^q \Fskel^\punc(\ord{1}, \ord{n}_k) 
\iso \txprod_{k=1}^q \ord{n}_k^\punc
= \txprod_{k=1}^q \ufs{n}_k\]
if $\ordtu{n} = \ang{\ord{n}_k}_{k=1}^q$ with each $\ord{n}_k \in \Fskel$.
\end{explanation}

\subsection*{Functors between $\Fskel$ and $\Gskel$}

The indexing categories $\Fskel$ and $\Gskel$ in, respectively, \cref{def:ordn,definition:Gstar} are related by the following functors.

\begin{definition}\label{definition:sma-functor}
We define the following pointed functors.
\begin{description}
\item[Length-One Inclusion]
The pointed functor
\begin{equation}\label{i-incl}
i\cn (\Fskel, \ord{0}) \to (\Gskel, \gst)
\end{equation}
sends
\begin{itemize}
\item each pointed finite set $\ord{n} \in \Fskel$ to the length-one tuple $(\ord{n}) \in \Gskel$ and
\item each morphism $\psi$ in $\Fskel$ to the pair $(1_{\ufs{1}},(\psi))$.
\end{itemize} 
\item[Smash Product]
The strict symmetric monoidal pointed functor
\begin{equation}\label{smagf}
\sma\cn \big(\Gskel,\oplus,\ang{},\gst\big) \to \big(\Fskel,\sma,\ord{1},\ord{0}\big)
\end{equation}
is defined on objects as follows:
\[\left\{\begin{aligned}
\sma \gst & = \ord{0},\\
\sma \ang{} & = \ord{1}, \andspace\\
\sma \ang{\ord{n}_k}_{k=1}^q & = \sma_{k=1}^q \ord{n}_k = \ord{n_1 \cdots n_q} \forspace q > 0.
\end{aligned}\right.\]
The image of a morphism 
\[(f,\ang{\psi}) \cn \ang{\ord{n}_k}_{k=1}^q \to \ang{\ord{m}_j}_{j=1}^p \inspace \Gskel\] 
under $\sma$ is the following composite in $\Fskel$.
\[\begin{tikzcd}[column sep=large]
\sma_{k=1}^q \ord{n}_k \ar{r}{f_*}[swap]{\iso} &
\sma_{k=1}^q \ord{n}_{f^\inv(k)} \ar{r}{\iso} &
\sma_{j=1}^p \ord{n}_{f^\inv(j)} \ar{r}{\sma_{j=1}^p \psi_j} &
\sma_{j=1}^p \ord{m}_j
\end{tikzcd}\]
\item[Induced Functors on Diagrams]
In the context of \cref{def:Gamma-object,def:Gstar-object}, $\sma$ in \cref{smagf} and $i$ in \cref{i-incl} induce the functors
\begin{equation}\label{smastar-istar}
\begin{tikzcd}[column sep=large]
\GaV \ar{r}{\sma^*} & \GstarV \ar{r}{i^*} & \GaV
\end{tikzcd}
\end{equation} 
given by precomposition and whiskering with $\sma$ and $i$, respectively.\defmark
\end{description}
\end{definition}

\begin{explanation}\label{expl:sma-i}
Consider \cref{definition:sma-functor}.
\begin{enumerate}
\item The functor $i$ is fully faithful, and 
\[\sma \circ i = 1_{\Fskel} \cn \Fskel \to \Gskel \to \Fskel.\]
This implies that the composite in \cref{smastar-istar}, $i^* \circ \sma^*$, is the identity functor on $\GaV$.
\item While $\sma$ is strict symmetric monoidal,  $i$ is neither monoidal nor oplax monoidal.\defmark
\end{enumerate}
\end{explanation}

\section{Segal and Elmendorf-Mandell \texorpdfstring{$K$}{K}-theory}
\label{sec:segalEMK}

In this section we briefly review 
\begin{itemize}
\item Segal $K$-theory 
\[\Kse = \Kf \circ~ \Ner_* \circ~ \Jse \cn \permcatsu \to \Spc \; ;\]
\item the homotopy inverse functors $\cP$, $S_*$, and $\bA$; and
\item Elmendorf-Mandell $K$-theory 
\[\Kem = \Kg \circ~ \Ner_* \circ~ \Jem \cn \permcatsu \to \Spc.\]
\end{itemize}  
They are summarized in the diagram \cref{eq:Ksummary} below.
\begin{equation}\label{eq:Ksummary}

\end{equation}
Each arrow in \cref{eq:Ksummary}, \emph{except} $\Jt$ and $\Kg$, is an equivalence of homotopy theories.  The following table summarizes their multicategorical and symmetric monoidal properties.
\smallskip
\begin{center}
\resizebox{.95\textwidth}{!}{
{\renewcommand{\arraystretch}{1.4}%
{\setlength{\tabcolsep}{1ex}
\begin{tabular}{|c|c|c|}\hline
$\permcatsu$ & $\Cat$-multicategory & \ref{thm:permcatmulticat} \\ \hline
$\MoneMod$, $\Gacat$, $\Gstarcat$ & symmetric monoidal $\Cat$-categories & \ref{definition:MoneMod}, \ref{expl:Gamma-object}, \ref{expl:Gstar-object} \\ \hline
$\Gasset$, $\Gstarsset$, $\Spc$ & symmetric monoidal $\sSet$-categories & \ref{expl:Gamma-object}, \ref{expl:Gstar-object}, \ref{SymSp} \\ \hline\hline
$\cP$ & \emph{non-symmetric} $\Cat$-multifunctor & \ref{PinvK} \\ \hline
$\Kem$ & $\sSet$-multifunctor & \ref{Kem} \\ \hline
$\Jem$, $\Endm$ & $\Cat$-multifunctors & \ref{Jem}, \ref{expl:endm-catmulti} \\ \hline
$\sma^*$ & enriched symmetric monoidal functors & \ref{smastar-istar} \\ \hline
$\Jt$ & symmetric monoidal $\Cat$-functor & \ref{Jt-smcat} \\ \hline
$\Ner_*$, $\Kf$, $\Kg$ & symmetric monoidal $\sSet$-functors & \ref{NerGa}, \ref{NerGs}, \ref{Kf-sm}, \ref{Kg-sm} \\ \hline\hline
$\Kse$, $\Jse$, $i^*$, $L$, $S_*$, $\bA$ & \emph{not} multifunctors & \ref{Kse}, \ref{Jse}, \ref{smastar-istar}, \ref{leftadjL}, \ref{Sstar-mandell}, \ref{A-segal} \\ \hline
\end{tabular}}}}
\end{center}
\smallskip
Later in this work, when we apply our general results to some of these functors, we only need to use some of their categorical properties, which we recall below.  Detailed discussion of these functors are in \cite[Chapters 8--10]{cerberusIII}, \cite{johnson-yau-invK,johnson-yau-multiK}, and the references cited below.

\subsection*{Categories}
\label{subsec:Kcat}

The categories in \cref{eq:Ksummary} are defined as follows.  Recall symmetric monoidal $\V$-category and $\V$-multicategory from \cref{definition:symm-monoidal-vcat,def:enr-multicategory}, respectively.  
\begin{itemize}
\item $\permcatsu$ is the category of small permutative categories and strictly unital symmetric monoidal functors in \cref{def:permcat}.  By \cref{thm:permcatmulticat} it is a $\Cat$-multicategory.  
\item $\MoneMod$ is the symmetric monoidal $\Cat$-category of left $\Mone$-modules in \cref{definition:MoneMod}.
\item $\Gacat$ and $\Gasset$ are the categories of $\Ga$-categories and $\Ga$-simplicial sets, respectively, in \cref{def:Gamma-object}.  They are symmetric monoidal categories enriched in $\Cat$ and $\sSet$, respectively, by \cref{expl:Gamma-object}.
\item $\Gstarcat$ and $\Gstarsset$ are the categories of $\Gstar$-categories and $\Gstar$-simplicial sets, respectively, in \cref{def:Gstar-object}.  They are symmetric monoidal categories enriched in $\Cat$ and $\sSet$, respectively, by \cref{expl:Gstar-object}.
\item $\Spc$ is the category of connective symmetric spectra.  It is the full subcategory of the category 
\begin{equation}\label{SymSp}
\Sp
\end{equation} 
of all symmetric spectra \cite{hss}.  Connective means that the negative degree homotopy groups are trivial.  The category $\Sp$ is complete, cocomplete, and symmetric monoidal closed.  Moreover, it is enriched, tensored, and cotensored over $\psset$.  In particular, $\Spc$ is a symmetric monoidal $\sSet$-category.  See \cite[Chapter 7]{cerberusIII} for an elementary discussion of symmetric spectra.
\end{itemize}

\subsection*{Segal $K$-Theory}
\label{subsec:segalK}

The left-to-right composite functor along the top row of \cref{eq:Ksummary},
\begin{equation}\label{Kse}
\Kse = \Kf \circ \Ner_* \circ~ \Jse \cn \permcatsu \to \Spc,
\end{equation}
is called \index{K-theory@$K$-theory!Segal}\index{Segal!K-theory@$K$-theory}\emph{Segal $K$-theory}.  We also denote by $\Kse$ its composite with the subcategory inclusion $\Spc \hookrightarrow \Sp$.  The constituent functors are as follows.
\begin{itemize}
\item The first functor is \index{J-theory@$J$-theory!Segal}\index{Segal!J-theory@$J$-theory}\emph{Segal $J$-theory} \cite{segal,may-permutative}
\begin{equation}\label{Jse}
\Jse \cn \permcatsu \to \Gacat
\end{equation}
that sends each small permutative category $\C$ to the $\Ga$-category
\[(\Jse\C)(-) = \pMulticat\big(\cM(-), \Endst\C\big) \cn \Fskel \to \pCat.\]
In this definition,
\begin{itemize}
\item $\cM(-)$ is the partition multicategory in \cref{definition:calM}, and
\item $\Endst\C$ is the pointed endomorphism multicategory in \cref{ex:endstc}.
\end{itemize} 
While both its domain and codomain are $\Cat$-multicategories, $\Jse$ is \emph{not} a multifunctor because it is incompatible with the multiplicative structures.  See \cite[Section 8.5]{cerberusIII} for a thorough discussion.
\item The second functor,
\begin{equation}\label{NerGa}
\Ner_* \cn \Gacat \to \Gasset,
\end{equation}
is induced by precomposing and whiskering with the nerve functor, $\Ner$, in \cref{nerve}.  Since the nerve is a right adjoint, it preserves all small limits, in particular, terminal objects and finite products.  Therefore, $\Ner_*$ is a symmetric monoidal $\sSet$-functor by \cite[3.7.28]{cerberusIII}.
\item The third functor \cite{bousfield_friedlander,segal},
\begin{equation}\label{Kf-sm}
\Kf \cn \Gasset \to \Spc,
\end{equation}
sends each $\Ga$-simplicial set $X$ to the connective symmetric spectrum 
\begin{equation}\label{KfX}
\Kf X = \Big\{ (\Kf X)k = \big|X(S^k)\big| \Big\}_{k \geq 0}
\end{equation}
with $|-|$ denoting the diagonal and $S^k = (S^1)^{\sma k}$ denoting the standard simplicial $k$-sphere.  Thus $(\Kf X)k $ is the pointed simplicial set whose set of $n$-simplices is given by
\[\big((\Kf X)k\big)_n = \left( X\big(S^k_n\big) \right)_n = \big( X\ord{n^k} \big)_n.\]
For $k \geq 1$ the structure morphism 
\[\Sigma \big|X(S^{k-1})\big| \to \big|X(S^k)\big|\]
is induced by the inclusions
\[S^{k-1}_n \iso \{0,i\} \sma S^{k-1}_n \hookrightarrow S^1_n \sma S^{k-1}_n \iso S^k_n\]
for $i \in \{1,\ldots,n\}$.  The $\Sigma_k$-action on $\left|X(S^k)\right|$ is induced by the $\Sigma_k$-action on $S^k = (S^1)^{\sma k}$ that permutes the $k$ smash factors.  See \cite[Section 8.2]{cerberusIII} for a thorough discussion.  Moreover, with $\Sp$ as the codomain, $\Kf$ is a symmetric monoidal $\sSet$-functor by \cite[9.4.18]{cerberusIII}.
\end{itemize}
We discuss the functors $\cP$, $S_*$, and $\bA$ in \cref{PinvK,Sstar-mandell,A-segal} below when we discuss equivalences of homotopy theories.

\subsection*{Elmendorf-Mandell $K$-Theory}
\label{subsec:EMK}

The composite along the bottom of \cref{eq:Ksummary},
\begin{equation}\label{Kem}
\begin{aligned}
\Kem &= \Kg \circ \Ner_* \circ~ \Jem\\
&= \Kg \circ \Ner_* \circ~ \Jt \circ \Endm \cn \permcatsu \to \Spc,
\end{aligned}
\end{equation}
is called \index{K-theory@$K$-theory!Elmendorf-Mandell}\index{Elmendorf-Mandell!K-theory@$K$-theory}\emph{Elmendorf-Mandell $K$-theory} \cite{elmendorf-mandell,elmendorf-mandell-perm}.  We also denote by $\Kem$ its composite with the subcategory inclusion $\Spc \hookrightarrow \Sp$.  The constituent enriched multifunctors are as follows.
\begin{itemize}
\item The endomorphism left $\Mone$-module $\Cat$-multifunctor 
\[\Endm \cn \permcatsu \to \MoneMod\]
is as in \cref{expl:endm-catmulti}.
\item The symmetric monoidal $\Cat$-functor
\begin{equation}\label{Jt-smcat}
\Jt = \pMulticat\big(\cT,-\big) \cn \MoneMod \to \Gstarcat
\end{equation}
sends each left $\Mone$-module $\M$ to the $\Gstar$-category
\[(\Jt\M)(-) = \pMulticat\big(\cT(-), \M\big) \cn \Gskel \to \pCat.\]
Here $\cT$ is defined by
\[\cT\ang{\ord{m}_j}_{j=1}^p = \txsma_{j=1}^p \cM\ord{m}_j \forspace \ang{\ord{m}_j}_{j=1}^p \in \Gskel\]
with $\cM(-)$ the partition multicategory in \cref{definition:calM}.  See \cite[10.3.17]{cerberusIII} for a detailed discussion.
\item \index{J-theory@$J$-theory!Elmendorf-Mandell}\index{Elmendorf-Mandell!J-theory@$J$-theory}\emph{Elmendorf-Mandell $J$-theory} is the composite $\Cat$-multifunctor
\begin{equation}\label{Jem}
\Jem = \Jt \circ \Endm \cn \permcatsu \to \Gstarcat,
\end{equation}
which associates to each small permutative category a $\Gstar$-category.  See \cite[Section 10.3]{cerberusIII} for a thorough discussion.  
\item The symmetric monoidal $\sSet$-functor \cite[9.2.19]{cerberusIII}
\begin{equation}\label{NerGs}
\Ner_* \cn \Gstarcat \to \Gstarsset
\end{equation}
is induced by precomposing and whiskering with the nerve functor, $\Ner$.
\item The symmetric monoidal $\sSet$-functor \cite[Definition 4.5]{elmendorf-mandell}
\begin{equation}\label{Kg-sm}
\Kg \cn \Gstarsset \to \Sp
\end{equation}
sends each $\Gstar$-simplicial set $X$ to the connective symmetric spectrum
\[\Kg X = \Big\{ (\Kg X)k = \big|X\big(\underbrace{S^1, \ldots, S^1}_{k}\big)\big| \Big\}_{k \geq 0}.\]
The structure morphisms and symmetric group action are analogous to those for $\Kf(-)$ in \cref{KfX}.  See \cite[Sections 9.3 and 9.4]{cerberusIII} for a thorough discussion.
\end{itemize}
It follows that $\Kem$ is an $\sSet$-multifunctor.

\subsection*{Relating Segal and Elmendorf-Mandell $K$-Theory}
\label{subsec:segalEMK}

The other functors in \cref{eq:Ksummary} are defined as follows.
\begin{itemize}
\item The functors
\[\begin{tikzcd}[column sep=large, row sep=0ex]
\Gacat \ar{r}{\sma^*} & \Gstarcat \ar{r}{i^*} & \Gacat\\
\Gasset \ar{r}{\sma^*} & \Gstarsset \ar{r}{i^*} & \Gasset
\end{tikzcd}\]
are induced by
\begin{itemize}
\item the smash product, $\sma \cn \Gskel \to \Fskel$, and
\item the length-one inclusion functor, $i \cn \Fskel \to \Gskel$, as in \cref{smastar-istar}.
\end{itemize}
Since $\sma$ is a strict symmetric monoidal functor, each $\sma^*$ is a symmetric monoidal functor in the enriched sense by \cite[9.4.18]{cerberusIII}.  However, neither $i^*$ is a multifunctor because $i$ is not compatible with the permutative structures of its domain and codomain.
\item In the middle square in \cref{eq:Ksummary}, associativity of composition of functors and whiskering implies the following equalities.
\[\begin{split}
\Ner_* \circ~ \sma^* = \sma^* \circ \Ner_* & \cn \Gacat \to \Gstarsset\\
\Ner_* \circ~ i^* = i^* \circ \Ner_* & \cn \Gstarcat \to \Gasset.
\end{split}\]
\item The right region in \cref{eq:Ksummary} commutes,
\[\Kg \circ \sma^* = \Kf \cn \Gasset \to \Spc.\]
See \cite[9.3.16]{cerberusIII} for a proof. 
\item In the left region in \cref{eq:Ksummary}, there is an equality
\[i^* \circ \Jem = \Jse \cn \permcatsu \to \Gacat\]
by the definitions of the functors involved; see \cite[8.5.1, 10.3.1, and 10.3.27]{cerberusIII}.  Moreover, there is a natural transformation \cite[Theorem 4.6]{elmendorf-mandell}\label{not:Pistar}
\[\Pi^* \cn \sma^* \circ \Jse \to \Jem \cn \permcatsu \to \Gstarcat.\]
See \cite[Section 10.6]{cerberusIII} for a detailed discussion.
\item Each of the functors
\begin{equation}\label{leftadjL}
\begin{tikzcd}[column sep=normal]
\Gacat \ar{r}{L} & \Gstarcat
\end{tikzcd}
\andspace
\begin{tikzcd}[column sep=normal]
\Gasset \ar{r}{L} & \Gstarsset
\end{tikzcd}
\end{equation}
is the left adjoint of the respective functor $i^*$.  See \cite[Section 3]{johnson-yau-multiK} for a detailed construction of $L$.  By the explicit formulas of the pointed \index{pointed!diagrams!Day convolution}\index{Day convolution!pointed diagrams}Day convolution, \cref{eq:GaV-convolution,eq:GstarV-convolution}, neither $L$ is compatible with the multiplicative structures of its domain and codomain.
\end{itemize}

\subsection*{Stable Equivalences}
\label{subsec:Kstableeq}

Each category in \cref{eq:Ksummary} is equipped with the structure of a relative category (\cref{definition:rel-cat}) as follows.
\begin{itemize}
\item The pair 
\[\brb{\Spc,\cS}\]
is a relative category, where $\cS$ is the wide subcategory of \index{equivalence!stable}\index{stable equivalence}\emph{stable equivalences} of connective symmetric spectra.
\item For each of $\permcatsu$, $\Gacat$, and $\Gasset$, we denote by $\cS$ the wide subcategory of morphisms created by the indicated functor:
\begin{equation}\label{perm-steq}
\brb{\permcatsu,\cS} \fto{\Jse}
\brb{\Gacat,\cS} \fto{\Ner_*}
\brb{\Gasset,\cS} \fto{\Kf}
\brb{\Spc,\cS}.
\end{equation}
In each case, we call morphisms in $\cS$ \emph{stable equivalences}.  In particular, stable equivalences in $\permcatsu$ are created by Segal $K$-theory $\Kse$ \cref{Kse}.
\item We denote by\label{not:csi}
\[\csi \bigsubset \Gstarcat \andspace \csi \bigsubset \Gstarsset\]
the wide subcategories created by, respectively, the functors
\[\begin{tikzcd}
\Gstarcat \ar{r}{i^*} & \Gacat
\end{tikzcd}
\andspace
\begin{tikzcd}
\Gstarsset \ar{r}{i^*} & \Gasset.
\end{tikzcd}\]
We refer to morphisms in $\csi$ as \emph{$i^*$-stable equivalences}.
\item The natural transformation $\Pi^*$ is componentwise an $i^*$-stable equivalence in $\Gstarcat$; see \cite[4.10]{johnson-yau-multiK} for an explanation.
\end{itemize}
In \cref{def:ptmulti-stableeq} we equip $\MoneMod$ with the structure of a relative category.

\begin{remark}[Stable Equivalences]\label{remark:steq}
  The following two key properties of stable equivalences will be used repeatedly below.
  \begin{enumerate}
  \item\label{it:steq-1} Each class of stable equivalences includes isomorphisms, is closed under composition, and has the 2-out-of-3 property.
    These follow from the stronger statement that there is a Quillen model structure on $\Sp$ whose weak equivalences are the stable equivalences \cite[3.4.4 and~5.3.8]{hss}.
  \item\label{it:steq-2} Suppose given
    \[
      P\cn \C \to \D \inspace \permcatsu.
    \]
    If the underlying functor of $P$ is a left or right adjoint, then $P$ is a stable equivalence.
    This follows from the observations that (i) an adjunction of categories induces a homotopy equivalence on nerves \cite[7.2.5]{cerberusIII} and (ii) the stable equivalences of symmetric spectra contain the level equivalences \cite[7.8.8]{cerberusIII}.\dqed
  \end{enumerate}
\end{remark}

\subsection*{Equivalences of Homotopy Theories}
\label{subsec:K-heq}

Equipped with the relative category structures above, each arrow in \cref{eq:Ksummary}, \emph{except} $\Jt$ and $\Kg$, is an equivalence of homotopy theories in the sense of \cref{definition:rel-cat-pow}.

\subsubsection*{Segal $K$-Theory}\

\begin{itemize}
\item Each of the two relative functors
\[\begin{tikzcd}[column sep=large]
\permcatsu \ar{r}{\Jse} & \Gacat \ar{r}{\Ner_*} & \Gasset
\end{tikzcd}\]
is an equivalence of homotopy theories by the work of Mandell \cite{mandell_inverseK}, which sharpens earlier work of Thomason \cite{thomason}.
\item The relative functor
\[\Kf \cn \Gasset \to \Spc\]
is an equivalence of homotopy theories by the work of Segal \cite{segal} and Bousfield-Friedlander \cite{bousfield_friedlander}.  Therefore, Segal $K$-theory
\begin{equation}\label{Kse-heq}
\Kse = \Kf \circ~ \Ner_* \circ~ \Jse \cn \permcatsu \to \Spc
\end{equation}
is an equivalence of homotopy theories.
\end{itemize}

\subsubsection*{Homotopy Inverses}

Each constituent functor in Segal $K$-theory admits a homotopy inverse functor given by the right-to-left functors along the top row of \cref{eq:Ksummary}.
\begin{itemize}
\item The work of Mandell \cite{mandell_inverseK} constructs the relative functor $\cP$ in
\begin{equation}\label{PinvK}
\begin{tikzcd}[column sep=large]
\permcatsu \ar[shift right]{r}[swap]{\Jse} & \Gacat \ar[shift right]{l}[swap]{\cP}
\end{tikzcd}
\end{equation}
and shows that the pair $(\cP,\Jse)$ forms inverse equivalences of homotopy theories as in \cref{def:inverse-heq}.  Thus $\cP$ is an equivalence of homotopy theories by \cref{gjo29}.  

While $\Jse$ is not a multifunctor even in the non-symmetric sense, $\cP$ \emph{is} a non-symmetric $\Cat$-multifunctor by \cite[1.3]{johnson-yau-invK}.  Moreover, $\cP$ is a \index{multifunctor!pseudo symmetric $\Cat$-}\index{pseudo symmetric!Cat-multifunctor@$\Cat$-multifunctor}\emph{pseudo symmetric} $\Cat$-multifunctor by \cite[10.12]{yau-multigro}.  This means that $\cP$ preserves the symmetric group action up to natural isomorphisms that satisfy further coherence axioms.  Since we will not use the pseudo symmetry of $\cP$ in this work, we refer the reader to \cite{yau-multigro} for detailed definitions and discussion.

\item The work of Mandell \cite{mandell_inverseK} also constructs the relative functor $S_*$ in
\begin{equation}\label{Sstar-mandell}
\begin{tikzcd}[column sep=large]
\Gacat \ar[shift right]{r}[swap]{\Ner_*} & \Gasset \ar[shift right]{l}[swap]{S_*}
\end{tikzcd}
\end{equation}
and shows that the pair $(S_*,\Ner_*)$ forms inverse equivalences of homotopy theories as in \cref{def:inverse-heq}.  Thus $S_*$ is an equivalence of homotopy theories by \cref{gjo29}.  In contrast to the symmetric monoidal $\sSet$-functor $\Ner_*$, the functor $S_*$ is \emph{not} a multifunctor even in the non-symmetric sense.  See the introduction of \cite{johnson-yau-multiK} for an explanation.

\item The relative functor $\bA$ in 
\begin{equation}\label{A-segal}
\begin{tikzcd}[column sep=large]
\Gasset \ar[shift right]{r}[swap]{\Kf} & \Spc \ar[shift right]{l}[swap]{\bA}
\end{tikzcd}
\end{equation}
is constructed in \cite[Def.\! 3.1]{segal}.  The work of Bousfield-Friedlander \cite[Theorem 5.8]{bousfield_friedlander} shows that the adjoint pair $\Kf \dashv \bA$ is a Quillen equivalence.  In contrast to the symmetric monoidal $\sSet$-functor $\Kf$, the functor $\bA$ is \emph{not} compatible with the multiplicative structures of its domain and codomain.  Thus $\bA$ is not a multifunctor even in the non-symmetric sense.
\end{itemize}

\subsubsection*{Elmendorf-Mandell $K$-Theory}\

\begin{itemize}
\item The authors show in \cite{johnson-yau-multiK} that each of the three relative functors
\begin{equation}\label{Kem-heq}
\left\{\begin{aligned}
\Jem & \cn \permcatsu \to \Gstarcat,\\
\Ner_* & \cn \Gstarcat \to \Gstarsset,\andspace\\
\Kem = \Kg \circ~ \Ner_* \circ~ \Jem & \cn \permcatsu \to \Spc
\end{aligned}\right.
\end{equation}
is an equivalence of homotopy theories.  Therefore, the composite
\[\Ner_* \circ~ \Jem \cn \permcatsu \to \Gstarsset\]
is also an equivalence of homotopy theories.  
\item The work of \cite{johnson-yau-multiK} also shows that each of the six relative functors
\[\begin{tikzpicture}[xscale=1,yscale=1]
\def\v{-1.5}
\draw[0cell=.9]
(0,0) node (a) {\Gacat}
(a)+(0,\v) node (b) {\Gstarcat}
(a)+(3,0) node (c) {\Gasset}
(c)+(0,\v) node (d) {\Gstarsset}
;
\draw[1cell=.9]
(b) edge[transform canvas={xshift=-.5ex}] node[swap] {i^*} (a)
(a) edge[transform canvas={xshift=-1.8ex}] node[swap] {\sma^*} (b)
(a) edge[bend left=30, transform canvas={xshift=.8em}] node {L} node[swap] {\vdash} (b)
(d) edge[transform canvas={xshift=.5ex}] node {i^*} (c)
(c) edge[transform canvas={xshift=1.8ex}] node {\sma^*} (d)
(c) edge[bend right=30, transform canvas={xshift=-.9em}] node[swap] {L} node {\dashv} (d)
;
\end{tikzpicture}\]
is an equivalence of homotopy theories.
\end{itemize}
In \cref{ptmulticat-thm-xi,ptmulticat-xxv} we show that
\[\Endm \cn \permcatsu \to \MoneMod\]
is an equivalence of homotopy theories.  We emphasize that the functors
\[\begin{split}
\Jt & \cn \MoneMod \to \Gstarcat \andspace\\
\Kg & \cn \Gstarsset \to \Spc
\end{split}\]
are probably not equivalences of homotopy theories because they are not known to be relative functors.

\chapter{Homotopy Theory of Multicategories}
\label{ch:multperm}
Recall the following 2-categories from \cref{def:permcat,v-multicat-2cat}.
\begin{itemize}
\item $\permcatst$ is the 2-category of small permutative categories, \emph{strict} symmetric monoidal functors, and monoidal natural transformations.
\item $\permcatsu$ is the larger 2-category with \emph{strictly unital} symmetric monoidal functors as 1-cells.
\item $\Multicat$ is the 2-category of small multicategories, multifunctors, and multinatural transformations.
\end{itemize}
To prepare for \cref{ch:ptmulticat-sp,ch:ptmulticat-alg}, in this chapter we review equivalences of homotopy theories between these three 2-categories.  Here is a summary diagram.
\begin{equation}\label{FEKse-heq-intro}

\end{equation} 
Each arrow in the diagram \cref{FEKse-heq-intro} is an equivalence of homotopy theories, where
\[\Kse \cn \brb{\permcatsu,\cS} \fto{\sim} \brb{\Spc,\cS}\]
is Segal $K$-theory.  Thus, there is a composite equivalence of homotopy theories
\[%
\]
Moreover, for each small non-symmetric $\Cat$-multicategory $\Q$, there are equivalences of homotopy theories
\begin{equation}\label{FEsu-intro}
%
\end{equation}
between categories of non-symmetric $\Q$-algebras.  See \cref{thm:alg-hty-equiv}.

\subsection*{Connection with Other Chapters}\

\subsubsection*{Pointed Extension}
 
In \cref{ch:ptmulticat-sp} we extend the 2-adjunction and adjoint equivalence of homotopy theories
\begin{equation}\label{FE-unpt-intro}
%
\end{equation}
to small \emph{pointed} multicategories, $\pMulticat$, and left $\Mone$-modules, $\MoneMod$.  One subtlety of this pointed extension is that it is \emph{not} achieved through the free-forgetful adjunction between $\Multicat$ and $\pMulticat$ in \cref{dplustwofunctor}.  Instead, the key ingredient is a detailed analysis of the pointed variant of $\Fr$, denoted $\Fst$ in \cref{ptmulticat-thm-i}.  Note that $\Fr$ in \cref{FE-unpt-intro} is \emph{not} a multifunctor because a multifunctor structure on $\Fr$ requires using strictly unital, but generally non-strict, symmetric monoidal functors.  See \cref{theorem:F-multi}.

\medskip
\subsubsection*{Pointed Multifunctorial Extension}

With $\permcatsu$ in place of $\permcatst$, in \cref{ch:ptmulticat-alg} we extend
\begin{itemize}
\item the $\Cat$-multifunctors $\Fr$ and $\End$ (non-symmetric in the case of $\Fr$) and
\item the equivalences of homotopy theories in \cref{FEsu-intro}
\end{itemize} 
to small pointed multicategories and left $\Mone$-modules.  This requires a nontrivial extension of the machinery in \cref{sec:Fmultifunctor} and key results in \cref{ch:ptmulticat-sp}. 

\medskip
\subsubsection*{Enriched Mackey Functors}

In \cref{part:homotopy-mackey} we further extend the main results in \cref{ch:ptmulticat-alg} to  equivalences of homotopy theories between
\begin{itemize}
\item enriched Mackey functors based on permutative categories,
\item enriched Mackey functors based on small pointed multicategories, and
\item enriched Mackey functors based on left $\Mone$-modules.
\end{itemize} 
See \cref{mackey-xiv-pmulticat,mackey-xiv-mone,mackey-pmulti-mone}.

\subsection*{Background}

We use 2-adjunctions in \cref{def:twoadjunction}, enriched multifunctors and multinatural transformations in \cref{sec:enrmulticat}, multilinear functors in \cref{def:nlinearfunctor}, and equivalences of homotopy theories in \cref{sec:hty-thy}.  Discussion of Segal $K$-theory is in \cref{sec:segalEMK}.  The material in this chapter is adapted from \cite{johnson-yau-Fmulti,johnson-yau-permmult}, where we refer the reader for detailed proofs.

\subsection*{Chapter Summary}

In \cref{sec:freepermcat} we discuss the 2-functor $\Fr$ in \cref{FE-unpt-intro}, which we call the \emph{free permutative category} construction.  In \cref{sec:free-twoadjoint} we discuss the 2-adjunction $\Fr \dashv \End$ in \cref{FE-unpt-intro}.  In \cref{sec:componentwiseadjoint} we discuss a componentwise right adjoint $\vrho$ of the counit $\epz$ of $\Fr \dashv \End$.  This componentwise right adjoint is, furthermore, a symmetric monoidal functor.  In \cref{sec:Fmultifunctor} we extend the 2-functor $\Fr$ to a non-symmetric $\Cat$-multifunctor with codomain $\permcatsu$.  In \cref{sec:multicat-model} we define the subcategories $\cSI$ and $\cSF$ and review the equivalences of homotopy theories in \cref{FEKse-heq-intro,FEsu-intro}.  Here is a summary table.
\reftable{.98}{
  $\Fr$ on multicategories, multifunctors, and multinatural transformations
  & \ref{definition:free-perm}, \ref{definition:free-smfun}, and \ref{definition:free-perm-multinat} \\ \hline
  unit $\eta \cn 1 \to \End\,\Fr$ and counit $\epz \cn \Fr\,\End \to 1$ for $\Fr \dashv \End$
  & \ref{definition:eta} and \ref{definition:epz} \\ \hline
  componentwise right adjoint $\vrho_\C$ of $\epz_\C$
  & \ref{def:epzrho-adjunction} \\ \hline
  symmetric monoidal functor $\left(\vrho_\C, \vrho_\C^2, \vrho_\C^0\right)$
  & \ref{def:rhoc-monoidal} \\ \hline
  multilinear functor $\Frn \cn \ang{\Fr\M_i} \to \Fr(\txotimes_i \M_i)$
  & \ref{def:S-multi} \\ \hline
  non-symmetric $\Cat$-multifunctor $\Fr$
  & \ref{definition:F-multi} \\ \hline
  non-symmetric $\Cat$-multinatural unit $\eta \cn 1 \to \End\, \Fr$
  & \ref{lemma:eta-mnat} \\ \hline
  stable equivalences $\cSI$ and $\cSF$
  & \ref{def:mult-stableeq} \\ \hline
  equivalences of homotopy theories in \cref{FEKse-heq-intro,FEsu-intro}
  & \ref{thm:F-heq}, \ref{thm:alg-hty-equiv}, \ref{thm:Fsu-heq}, and \ref{cor:I-heq} \\ \hline
} 
We remind the reader of \cref{conv:universe,expl:leftbracketing}.

\section{Free Permutative Categories}
\label{sec:freepermcat}

In this section we describe a 2-functor
\[\Fr \cn \Multicat \to \permcatst.\]
In \cref{sec:free-twoadjoint} we discuss unit and counit that make the pair $(\Fr,\End)$ into a 2-adjunction.
\[\begin{tikzpicture}
\draw[0cell]
(0,0) node (a) {\Multicat}
(a)+(3.25,0) node (b) {\permcatst}
(a)+(1.55,0) node (x) {\bot}
;
\draw[1cell=.9]
(a) edge[bend left=15,transform canvas={yshift=-.8ex}] node {\Fr} (b)
(b) edge[bend left=15,transform canvas={yshift=.7ex}] node {\End} (a)
;
\end{tikzpicture}\]

\subsection*{Free Permutative Category of a Multicategory}

Recall from \cref{ufsn} that 
\[\ufs{n} =\{1,\ldots,n\}\] 
denotes the unpointed finite set with $n$ elements, where $\ufs{0} = \emptyset$.  Concatenation of tuples is denoted by $\oplus$.  The definition of $\Fr$ uses the following notation for sub-tuples of objects and morphisms.  

\begin{definition}[Sub-tuples]\label{definition:free-perm-helper}
Suppose $\M$ is a multicategory (\cref{def:enr-multicategory}), and $\ang{x}=\ang{x_i}_{i=1}^r$ is an $r$-tuple of objects in $\M$.  Suppose
\[\begin{tikzcd}
\ufs{r} \ar{r}{f} & \ufs{s} \ar{r}{g} & \ufs{t}
\end{tikzcd}\]
are functions of unpointed finite sets with $r,s,t \ge 0$.  We define the following.
\begin{enumerate}
\item For $j \in \ufs{s}$, we define the sub-tuple of $\angx$,
\begin{equation}\label{eq:x-finv}
\ang{x}_{f^\inv(j)} = \begin{cases}
\ang{x_i}_{i \in f^\inv(j)} & \text{if $f^\inv(j) \neq \emptyset$ and}\\
\ang{} & \text{if $f^\inv(j) = \emptyset$,}
\end{cases}
\end{equation}
consisting of those objects $x_i$ such that $f(i) = j$.
\item For an $s$-tuple $\ang{\phi} = \ang{\phi_j}_{j=1}^s$ of multimorphisms in $\M$ and $k \in \ufs{t}$, we define the sub-tuple of $\ang{\phi}$,
\begin{equation}\label{angphiginv}
\ang{\phi}_{g^\inv(k)} = \begin{cases}
\ang{\phi_j}_{j \in g^\inv(k)} & \text{if $g^\inv(k) \neq \emptyset$ and}\\
\ang{} & \text{if $g^\inv(k) = \emptyset.$}
\end{cases}
\end{equation}
\item For $k \in \ufs{t}$, we define $\si^k_{g,f} \in \Si_t$ to be the unique permutation determined by the equality
\begin{equation}\label{eq:sigma-kgf}
\bigg[\bigoplus_{j \in g^\inv(k)} \ang{x}_{f^\inv(j)}\bigg] \cdot \sigma^k_{g,f} 
= \ang{x}_{(gf)^\inv(k)}.
\end{equation}
\begin{itemize}
\item The tuple $\bigoplus_{j \in g^\inv(k)} \ang{x}_{f^\inv(j)}$ in \cref{eq:sigma-kgf} is the concatenation of the tuples $\ang{x}_{f^\inv(j)}$ in the order specified by $j \in g^\inv(k)$. 
\item The right-hand side of \cref{eq:sigma-kgf} is a sub-tuple of $\angx$, defined as in \cref{eq:x-finv}.
\end{itemize}
\end{enumerate}
This finishes the definition.
\end{definition}

Recall from \cref{def:profile} that $\Prof(S)$ means the class of finite tuples in $S$.  The free permutative category in the next definition is sketched in \cite[Theorem~4.2]{elmendorf-mandell-perm}.  

\begin{definition}[Free Permutative Category]\label{definition:free-perm}
Given a multicategory $(\M,\ga,\opu)$, we define a permutative category 
\[\big(\Fr\M, \oplus, \ang{}, \xi\big),\] 
which is called the \index{category!free permutative}\index{permutative category!free}\index{free!permutative category}\emph{free permutative category on $\M$}, as follows.
\begin{description}
\item[Objects] $\Ob(\Fr\M) = \Prof(\M)$, the class of finite tuples $\ang{x} = \ang{x_i}_{i=1}^r$ with each $x_i \in \Ob\M$ and $r \ge 0$.
\item[Morphisms] Given finite tuples $\ang{x} = \ang{x_i}_{i=1}^r$ and $\ang{y} = \ang{y_j}_{j=1}^s$, a morphism
\begin{equation}\label{FMmorphism}
(f,\ang{\phi}) \cn \ang{x} \to \ang{y} \inspace \Fr\M
\end{equation}
is a pair consisting of
\begin{itemize}
\item a function
\[f \cn \ufs{r} \to \ufs{s},\]
called the \index{index map}\emph{index map}, and
\item an $s$-tuple of multimorphisms
\[\ang{\phi} = \ang{\phi_j}_{j=1}^s \withspace \phi_j \in \M\lrscmap{\ang{x_i}_{i \in f^\inv(j)}; y_j}.\]
\end{itemize}
\item[Identities] 
The identity morphism for an object $\angx = \ang{x_i}_{i=1}^r$ in $\Fr\M$ is the pair
\[1_{\angx} = \left(1_{\ufs{r}} \scs \ang{\opu_{x_i}}_{i=1}^r\right).\]
\item[Composition] Given a pair of morphisms
\[\ang{x} = \ang{x_i}_{i=1}^r \fto{(f,\ang{\phi})} 
\ang{y} = \ang{y_j}_{j=1}^s \fto{(g,\ang{\psi})} 
\ang{z} = \ang{z_k}_{k=1}^t\]
their composite is the morphism
\begin{equation}\label{eq:FM-comp}
\left(gf \scs \bang{\theta_k \cdot \si^k_{g,f}}_{k=1}^t \right) \cn \angx \to \angz
\end{equation}
with
\begin{equation}\label{eq:thetak}
\theta_k = \ga\lrscmap{\psi_k;\ang{\phi_j}_{j \in g^\inv(k)}} 
\in \M\lrscmap{\;\bigoplus_{j \in g^\inv(k)} \ang{x}_{f^\inv(j)} ; z_k}
\end{equation}
for each $k \in \ufs{t}$ and $\si^k_{g,f}$ as in \cref{eq:sigma-kgf}.
\item[Monoidal Product on Objects]
The monoidal product 
\begin{equation}\label{FMproduct}
\oplus \cn \Fr\M \times \Fr\M \to \Fr\M
\end{equation}
is given by concatenation of finite tuples on objects:
\begin{equation}\label{FMproduct-obj}
\ang{x_i}_{i=1}^r \oplus \ang{y_j}_{j=1}^s = \brb{\ang{x_i}_{i=1}^r, \ang{y_j}_{j=1}^s}.
\end{equation}
\item[Monoidal Product on Morphisms]
Given a pair of morphisms
\[\big(f,\ang{\phi_j}_{j=1}^s\big) \cn \ang{x_i}_{i=1}^r \to \ang{y_j}_{j=1}^s
\andspace
\big(f',\ang{\phi'_j}_{j=1}^{s'}\big) \cn \ang{x'_i}_{i=1}^{r'} \to \ang{y'_j}_{j=1}^{s'}\]
in $\Fr\M$, their monoidal product is the morphism
\begin{equation}\label{FMproduct-mor}
\brb{f \oplus f', \ang{\phi} \oplus \ang{\phi'}} \cn \ang{x} \oplus \ang{x'} \to \ang{y} \oplus \ang{y'}.
\end{equation}
In \cref{FMproduct-mor} the index map is the composite
\[\begin{tikzpicture}
\def\h{1.75} \def\u{.6}
\draw[0cell]
(0,0) node (a) {\ufs{r + r'}}
(a)+(\h,0) node (b) {\ufs{r} \bincoprod \ufs{r'}}
(b)+(2.7,0) node (c) {\ufs{s} \bincoprod \ufs{s'}}
(c)+(\h,0) node (d) {\ufs{s + s'}}
;
\draw[1cell=.8]
(a) edge node {\iso} (b)
(b) edge node {f \bincoprod f'} (c)
(c) edge node {\iso} (d)
;
\draw[1cell=.9]
(a) [rounded corners=3pt] |- ($(b)+(0,\u)$)
-- node {f \oplus f'} ($(c)+(0,\u)$) -| (d)
;
\end{tikzpicture}\]
given by
\begin{itemize}
\item the canonical order-preserving isomorphisms and
\item the disjoint union of $f$ with $f'$.
\end{itemize} 
\item[Monoidal Unit]
The strict monoidal unit is the empty sequence $\ang{}$.  The associativity and unit isomorphisms for $\oplus$ are identity natural transformations. 
\item[Braiding] 
The braiding for objects $\ang{x_i}_{i=1}^r$ and $\ang{y_j}_{j=1}^s$ is 
\begin{equation}\label{FMbraiding}
\xi_{\ang{x},\ang{y}} = \brb{\tau_{r,s} , \ang{\opu}} 
\cn \angx \oplus \angy \fto{\iso} \angy \oplus \angx.
\end{equation}
In \cref{FMbraiding} the index map is the composite
\[\begin{tikzpicture}
\def\h{1.75} \def\u{.6}
\draw[0cell]
(0,0) node (a) {\ufs{r + s}}
(a)+(\h,0) node (b) {\ufs{r} \bincoprod \ufs{s}}
(b)+(2.7,0) node (c) {\ufs{s} \bincoprod \ufs{r}}
(c)+(\h,0) node (d) {\ufs{s + r}}
;
\draw[1cell=.8]
(a) edge node {\iso} (b)
(b) edge node {\iso} node[swap] {\text{swap}} (c)
(c) edge node {\iso} (d)
;
\draw[1cell=.9]
(a) [rounded corners=3pt] |- ($(b)+(0,\u)$)
-- node {\tau_{r,s}} ($(c)+(0,\u)$) -| (d)
;
\end{tikzpicture}\]
given by
\begin{itemize}
\item the canonical order-preserving isomorphisms and
\item the block permutation that swaps $\ufs{r}$ and $\ufs{s}$, keeping the relative order within each block unchanged.
\end{itemize} 
Each entry in the $(r+s)$-tuple $\ang{\opu}$ in \cref{FMbraiding} is a colored unit of some $x_i$ or $y_j$.
\end{description}
This finishes the definition of $\big(\Fr\M, \oplus, \ang{}, \xi\big)$.
\end{definition}

A detailed proof of the following is in \cite[5.7]{johnson-yau-permmult}.

\begin{proposition}\label{proposition:free-perm}
For each multicategory $\M$, the quadruple in \cref{definition:free-perm}
\[\big(\Fr\M, \oplus, \ang{}, \xi\big)\]
is a permutative category.
\end{proposition}

\begin{example}[Free Permutative Category of the Initial Operad]\label{example:free-Mtu}
  The \index{operad!initial}\index{initial operad}initial operad $\Mtu$ in \cref{ex:vmulticatinitialterminal}~\cref{ex:initialoperad} has a single object $*$ and a single operation $\opu_* \in \Mtu_1$.
  The free permutative category $\Fr(\Mtu)$ is isomorphic to the permutation category defined as follows.
  \begin{itemize}
  \item Its objects are natural numbers, $n \geq 0$, corresponding to length-$n$ sequences of the object $* \in \Mtu$.
  \item Its morphisms are permutations
    \[
      \Fr(\Mtu)(p,q) =
      \begin{cases}
        \Si_p & \ifspace p = q,\\
        \varnothing & \ifspace p \neq q
      \end{cases}
    \]
    for $p,q \geq 0$.
  \item The permutative structure $\oplus$ is given by addition on objects and block sums on morphisms.
  \item The monoidal unit is the object 0.
  \item The braiding
    \[
      \xi_{p,q} \cn p+q \fto{\iso} q+p
    \]
    is the block permutation in $\Si_{p+q}$ that swaps the first $p$ elements with the last $q$ elements.
    This is denoted $\tau_{p,q}$ in \cref{tauqqprime}.\defmark
  \end{itemize} 
\end{example}

\begin{example}[Free Permutative Category of the Terminal Multicategory]\label{example:free-Mterm}
  The \index{multicategory!terminal}\index{terminal!multicategory}terminal multicategory $\Mterm$ in \cref{definition:terminal-operad-comm} has a single object and a unique $n$-ary operation for each $n$.
  The free permutative category $\Fr\Mterm$ is isomorphic to the natural number category $\mathbf{N}$, whose objects are natural numbers and whose morphisms are given by morphisms of finite sets
  \[
    \mathbf{N}(r,s) = \Set(\ufs{r},\ufs{s}).
  \]
  The natural number $r \in \mathbf{N}$ corresponds to the sequence of length $r$ where each term is the unique object of $\Mterm$.
  A morphism $f\cn \ufs{r} \to \ufs{s}$ in $\mathbf{N}$ corresponds to the morphism
  \[
    (f,\ang{\phi}) \in \Fr\Mterm
  \]
  where $\phi_j$ is the unique operation in $\Mterm$ of arity $|f^\inv(j)|$.
\end{example}

\subsection*{Free Permutative Category as a 2-Functor}

Next we define $\Fr$ on multifunctors (\cref{def:enr-multicategory-functor}) and multinatural transformations (\cref{def:enr-multicat-natural-transformation}).

\begin{definition}[$\Fr$ on Multifunctors]\label{definition:free-smfun}
Given a multifunctor $H\cn \M \to \N$, we define a \emph{strict} symmetric monoidal functor
\[\Fr H \cn \Fr\M \to \Fr\N\]
as follows.  
\begin{description}
\item[Object Assignment]
For an object $\ang{x_i}_{i=1}^r$ in $\Fr\M$, we define the object
\begin{equation}\label{FHobjects}
(\Fr H)\ang{x_i}_{i=1}^r = \ang{Hx_i}_{i=1}^r \inspace \Fr\N.
\end{equation}
\item[Morphism Assignment]
For a morphism 
\[\brb{f, \ang{\phi_j}_{j=1}^s} \cn \ang{x_i}_{i=1}^r \to \ang{y_j}_{j=1}^s \inspace \Fr\M\]
as in \cref{FMmorphism}, we define the morphism
\begin{equation}\label{eq:FH-fphi}
(\Fr H)(f,\ang{\phi}) = \brb{f, \ang{H\phi_j}_{j=1}^s} \cn \ang{Hx_i}_{i=1}^r \to \ang{Hy_j}_{j=1}^s
\end{equation}
in $\Fr\N$.
\item[Constraints]
The unit and monoidal constraints for $\Fr H$ are identities.
\end{description}
This finishes the definition of $\Fr H$.
\end{definition}

\begin{definition}[$\Fr$ on Multinatural Transformations]\label{definition:free-perm-multinat}
Suppose $H, K \cn \M \to \N$ are multifunctors.  Given a multinatural transformation
\[\omega \cn H \to K,\]
we define a monoidal natural transformation
\[\Fr\omega \cn \Fr H \to \Fr K\]
with component morphism
\begin{equation}\label{eq:Fka-x}
(\Fr\omega)_{\ang{x}} = \left(1_{\ufs{r}} \scs \ang{\omega_{x_i}}_{i=1}^r \right) 
\cn \ang{Hx_i}_{i=1}^r \to \ang{Kx_i}_{i=1}^r \inspace \Fr\N 
\end{equation}
for each object $\ang{x_i}_{i=1}^r$ in $\Fr\M$.
\end{definition}

The following result is \cite[5.13]{johnson-yau-permmult}.

\begin{proposition}\label{proposition:free-perm-functor}\index{category!free permutative - 2-functor}\index{permutative category!free!2-functor}\index{free!permutative category 2-functor}
The constructions in \cref{definition:free-perm,definition:free-smfun,definition:free-perm-multinat} provide a 2-functor
\[\Fr \cn \Multicat \to \permcatst.\]
\end{proposition}

We also use $\Fr$ to denote the composite of the 2-functor in \cref{proposition:free-perm-functor} with any one of the inclusion 2-functors in \cref{permcatinclusion}.

\section{Free Permutative Category as a Left 2-Adjoint}
\label{sec:free-twoadjoint}

In this section we recall the fact that the 2-functor $\Fr$ in \cref{proposition:free-perm-functor} is a left 2-adjoint of the \index{endomorphism!2-functor}endomorphism multicategory 2-functor in \cref{endtwofunctor} 
\[\End \cn \permcatst \to \Multicat.\]
Recall from \cref{ex:endc} that, for each permutative category $\C$, the endomorphism multicategory $\End(\C)$ has the same objects as $\C$.  Next we define the unit and counit for the 2-adjunction $(\Fr,\End)$ in the sense of \cref{def:twoadjunction}.

\begin{definition}[Unit]\label{definition:eta}
Given a multicategory $\M$, we define a multifunctor
\[\eta_{\M} \cn \M \to \End \;\Fr\M\]
as follows.
\begin{description}
\item[Object Assignment]
For an object $y$ in $\M$, we define the object
\begin{equation}\label{Funit-obj}
\eta_{\M}(y) = (y) \inspace \End\; \Fr\M,
\end{equation}
where $(y)$ on the right-hand side is the length-1 tuple consisting of the object $y$.
\item[Multimorphism Assignment]
For an $r$-ary multimorphism 
\[\phi \cn \angx = \ang{x_i}_{i=1}^r \to y \inspace \M,\] 
we define the $r$-ary multimorphism
\begin{equation}\label{Funit-operation}
\eta_{\M}(\phi) = \brb{\iota_r, (\phi)} \cn \ang{x} \to (y)
\end{equation}
in
\[\big(\End\; \Fr\M\big) \lrscmap{\bang{(x_i)}_{i=1}^r; (y)} 
= (\Fr\M)\brb{\angx,(y)}.\]
On the right-hand side of \cref{Funit-operation},
\begin{itemize}
\item $\iota_r \cn \ufs{r} \to \ufs{1}$ is the unique function, and
\item $(\phi)$ is the length-1 tuple consisting of $\phi$.
\end{itemize} 
\end{description}
This finishes the definition of $\eta_{\M}$.
Multifunctoriality of $\eta_\M$ follows from the definitions of monoidal sum and composition in $\Fr\M$.
The 2-naturality of $\eta$ follows from the termwise definitions of $\Fr H$ and $\Fr \om$ in \cref{definition:free-smfun,definition:free-perm-multinat}, respectively.
\end{definition}

\begin{definition}[Counit]\label{definition:epz}
Given a permutative category $(\C, \oplus, \pu, \xi)$, we define a \emph{strict} symmetric monoidal functor
\[\epz_\C \cn \Fr\; \End(\C) \to \C\]
as follows.
\begin{description}
\item[Object Assignment]
For an $r$-tuple $\ang{x_i}_{i=1}^r$ of objects in $\C$, we define the object
\begin{equation}\label{Fcounit-obj}
\epz_\C\ang{x} = \txoplus_{i=1}^r x_i \inspace \C,
\end{equation} 
where an empty $\oplus$ means the monoidal unit $\pu$.
\item[Morphism Assignment]
Suppose given a morphism
\begin{equation}\label{FEndC-mor}
\brb{f, \ang{\phi_j}_{j=1}^s} \cn \ang{x_i}_{i=1}^r \to \ang{y_j}_{j=1}^s \inspace \Fr\; \End(\C)
\end{equation}
with each
\[\phi_j \in \End(\C)\lrscmap{\ang{x}_{f^\inv(j)}; y_j} 
= \C\left(\txoplus_{i \in f^\inv(j)} x_i \scs y_j\right).\]
We define the morphism
\[\epz_\C\brb{f,\ang{\phi}} \cn \epz_\C \angx \to \epz_\C \angy \inspace \C\]
as the following composite.
\begin{equation}\label{Fcounit-mor}
\begin{tikzpicture}[baseline={(a.base)}]
\def\u{.7}
\draw[0cell]
(0,0) node (a) {\txoplus_{i=1}^r x_i}
(a)+(3,0) node (b) {\txoplus_{j=1}^s \txoplus_{i \in f^\inv(j)} x_i}
(b)+(3.6,0) node (c) {\txoplus_{j=1}^s y_j}
;
\draw[1cell=.8]
(a) edge node {\xi_f} node[swap] {\iso} (b)
(b) edge node {\txoplus_{j=1}^s \phi_j} (c)
;
\draw[1cell=.9]
(a) [rounded corners=3pt] |- ($(b)+(-1,\u)$)
-- node {\epz_\C\brb{f,\ang{\phi}}} ($(b)+(1,\u)$) -| (c)
;
\end{tikzpicture}
\end{equation}
In \cref{Fcounit-mor} $\xi_f$ is the unique coherence isomorphism in $\C$ that permutes the terms of the sum.  Existence and uniqueness of this coherence isomorphism follows from the symmetric monoidal Coherence Theorem \cite[XI.1, Th.\ 1]{maclane}.
\item[Constraints]
The unit and monoidal constraints of $\epz_\C$ are defined as the identities.
\end{description}
This finishes the definition of $\epz_\C$.
Verification that $\epz_\C$ is strict symmetric monoidal follows from strictness of concatenation and uniqueness of the coherence isomorphisms $\xi_f$.
The 2-naturality of $\epz$ follows because strict symmetric monoidal functors preserve the monoidal sums $\oplus_i x_i$ and coherence isomorphisms $\xi_f$ in \cref{Fcounit-obj,Fcounit-mor}.
\end{definition}

The following result combines \cite[6.2, 6.8, and 6.11]{johnson-yau-permmult}.  Recall the notion of a 2-adjunction from \cref{def:twoadjunction}.

\begin{theorem}\label{theorem:FE-adj}
There is a 2-adjunction
\[\begin{tikzpicture}
\draw[0cell]
(0,0) node (a) {\Multicat}
(a)+(3.25,0) node (b) {\permcatst}
(a)+(1.55,0) node (x) {\bot}
;
\draw[1cell=.9]
(a) edge[bend left=15,transform canvas={yshift=-.8ex}] node {\Fr} (b)
(b) edge[bend left=15,transform canvas={yshift=.7ex}] node {\End} (a)
;
\end{tikzpicture}\]
consisting of the following data.
\begin{itemize}
\item $\permcatst$ is the 2-category in \cref{def:permcat}. 
\item $\Multicat$ is the 2-category in \cref{v-multicat-2cat}.
\item The left adjoint is $\Fr$ in \cref{proposition:free-perm-functor}.
\item The right adjoint is $\End$ in \cref{endtwofunctor} restricted to $\permcatst$.
\item The unit 
\[\eta \cn 1_{\Multicat} \to \End\; \Fr\]
has components in \cref{definition:eta}.
\item The counit
\[\epz \cn \Fr\; \End \to 1_{\permcatst}\]
has components in \cref{definition:epz}.
\end{itemize}
\end{theorem}

\begin{remark}[Naturality of the Counit]\label{rk:epznaturality}
The counit $\epz$ in \cref{theorem:FE-adj} is only natural with respect to \emph{strict} symmetric monoidal functors.  This implies that the 2-adjunction $\Fr \dashv \End$ does \emph{not} extend to a 2-adjunction, or even an adjunction of underlying categories, to any of the larger 2-categories in \cref{def:permcat} in which the 1-cells are, in general, not strict symmetric monoidal.
\end{remark}

\section{Componentwise Right Adjoint of the Counit}
\label{sec:componentwiseadjoint}

In this section we discuss a componentwise right adjoint of the counit of the 2-adjunction $\Fr \dashv \End$ in \cref{theorem:FE-adj}.  We also discuss a symmetric monoidal structure on this componentwise right adjoint; see \cref{rhoc-monoidal}.

\begin{definition}\label{def:epzrho-adjunction}
For each permutative category $(\C,\oplus)$, we define an adjunction
\[\begin{tikzpicture}
\draw[0cell]
(0,0) node (a') {\Fr\; \End(\C)}
(a')+(.6,0) node (a) {\phantom{\C}}
(a)+(2,0) node (b) {\C}
(a)+(1,0) node (x) {\bot}
;
\draw[1cell=.9]
(a) edge[bend left=15] node {\epz_\C} (b)
(b) edge[bend left=15] node {\vrho_\C} (a)
;
\end{tikzpicture}\]
as follows, where $\epz_\C$ is the strict symmetric monoidal functor in \cref{definition:epz}.
\begin{description}
\item[Right Adjoint]
The functor $\vrho_\C$ is defined by the object and morphism assignments
\begin{equation}\label{rhoc-assignments}
\left\{\begin{aligned}
\vrho_\C (x) & = (x) \forspace x \in \ObC \andspace\\
\vrho_\C(\phi) & = \brb{1_{\ufs{1}},(\phi)} \cn (x) \to (y) \forspace \phi \in \C(x,y).
\end{aligned}\right.
\end{equation}
On the right-hand side of \cref{rhoc-assignments}, $(x)$, $(y)$, and $(\phi)$ are tuples of length 1.
\item[Counit]
The following composite is the identity functor.  
\begin{equation}\label{vrho-epz}
\C \fto{\vrho_\C} \Fr\; \End(\C) \fto{\epz_\C} \C
\end{equation}
We define the counit for the adjunction $(\epz_\C,\vrho_\C)$ as the identity functor,
\[\epsilon = 1_{\C}.\]
\item[Unit]
The composite
\begin{equation}\label{epz-vrho}
\Fr\; \End(\C) \fto{\epz_\C} \C \fto{\vrho_\C} \Fr\; \End(\C)
\end{equation}
is given by the following assignments for each object $\ang{x_i}_{i=1}^r$ and morphism $\brb{f,\ang{\phi_j}_{j=1}^s}$ in $\Fr\, \End(\C)$ as in \cref{FEndC-mor}.
\begin{equation}\label{epzrho-unit}
\left\{\begin{aligned}
\ang{x} & \mapsto \left(\txoplus_{i=1}^r x_i\,\right)\\
\brb{f,\ang{\phi}} & \mapsto \left(1_{\ufs{1}} \scs \big(\txoplus_{j=1}^s \phi_j\big)\circ \xi_f\,\right) \cn
\left(\txoplus_{i=1}^r x_i\,\right) \to \left(\txoplus_{j=1}^s y_j\,\right)
\end{aligned}\right.
\end{equation}
We define the unit for the adjunction $(\epz_\C,\vrho_\C)$ as the natural transformation
\[\ups \cn 1_{\Fr\,\End(\C)} \to \vrho_\C \epz_\C\]
with component morphism
\begin{equation}\label{eq:ups-x}
\ups_{\ang{x}} = \big( \iota_r \scs 1_{\oplus_{i=1}^r x_i} \big) \cn \ang{x_i}_{i=1}^r \to \big(\txoplus_{i=1}^r x_i\, \big) \inspace \Fr\,\End(\C)
\end{equation}
for each length-$r$ tuple $\ang{x_i}_{i=1}^r$ of objects in $\C$.  In \eqref{eq:ups-x},
\begin{itemize}
\item $\iota_r\cn \ufs{r} \to \ufs{1}$ is the unique function, and
\item the identity morphism
\[1_{\oplus_{i=1}^r x_i} \in \big(\End(\C)\big)\lrscmap{\ang{x_i}_{i=1}^r; \txoplus_{i=1}^r x_i}
= \C\left(\txoplus_{i=1}^r x_i \scs \txoplus_{i=1}^r x_i\right)\]
is an $r$-ary multimorphism in $\End(\C)$.
\end{itemize} 
\end{description}
This finishes the definition.
\end{definition}

The following result is \cite[6.13]{johnson-yau-permmult}.

\begin{proposition}\label{proposition:epz-rho-adj}
In the context of \cref{def:epzrho-adjunction}, there is an adjunction of categories
\begin{equation}\label{epzc-rhoc}
\begin{tikzpicture}[baseline={(a'.base)}]
\draw[0cell]
(0,0) node (a') {\Fr\; \End(\C)}
(a')+(.6,0) node (a) {\phantom{\C}}
(a)+(2,0) node (b) {\C}
(a)+(1,0) node (x) {\bot}
;
\draw[1cell=.9]
(a) edge[bend left=15] node {\epz_\C} (b)
(b) edge[bend left=15] node {\vrho_\C} (a)
;
\end{tikzpicture}
\end{equation}
with the following data.
\begin{itemize}
\item The left adjoint is $\epz_\C$ in \cref{definition:epz}.
\item The right adjoint is $\vrho_\C$ in \cref{rhoc-assignments}.
\item The counit $\epsilon$ is the identity functor on $\C$.
\item The unit $\ups$ has components in \cref{eq:ups-x}.
\end{itemize}
\end{proposition}

\subsection*{Symmetric Monoidal Structure}

The right adjoint $\vrho_\C$ in \cref{epzc-rhoc} is a symmetric monoidal functor with the following structure morphisms.  The next definition uses the permutative category structure on $\Fr\,\End(\C)$ in \cref{definition:free-perm}, with $\M = \End(\C)$ in \cref{ex:endc}.

\begin{definition}\label{def:rhoc-monoidal}
For each permutative category $(\C,\oplus,\pu,\xi)$, we define unit constraint $\vrho_\C^0$ and monoidal constraint $\vrho_\C^2$ for the functor in \cref{epzc-rhoc}
\[\vrho_\C \cn \C \to \Fr\; \End(\C)\]
as follows.
\begin{description}
\item[Unit Constraint] It is the morphism
\begin{equation}\label{rhoc-unit}
\vrho_\C^0 = \srb{\iota_0, 1_\pu} \cn \ang{} \to (\pu) = \vrho_\C(\pu) \inspace \Fr\; \End(\C)
\end{equation}
defined as follows.
\begin{itemize}
\item $\iota_0 \cn \ufs{0} = \emptyset \to \ufs{1}$ is the unique function.
\item The identity morphism
\[1_\pu \in \big(\End(\C)\big)\scmap{\ang{}; \pu} = \C(\pu, \pu)\]
is a nullary multimorphism in $\End(\C)$.
\end{itemize}
\item[Monoidal Constraint]
For each pair of objects $x,y \in \C$, the monoidal constraint has a component morphism
\begin{equation}\label{rhoc-monconstraint}
(\vrho_\C^2)_{x,y} = \brb{\iota_2, 1_{x \oplus y}} \cn (x) \oplus (y) = (x,y) \to (x \oplus y) \inspace \Fr\; \End(\C)
\end{equation}
defined as follows.
\begin{itemize}
\item $\iota_2 \cn \ufs{2} \to \ufs{1}$ is the unique function.
\item The identity morphism
\[1_{x \oplus y} \in \big(\End(\C)\big)\scmap{x,y; x \oplus y} = \C\big(x \oplus y, x \oplus y\big)\]
is a binary multimorphism in $\End(\C)$.
\end{itemize}
\end{description}
This finishes the definition of $\vrho_\C^0$ and $\vrho_\C^2$.
\end{definition}

The terminal property of $\ufs{1} = \{1\}$ and the permutative category axioms of $\C$ imply the following.

\begin{lemma}\label{rhoc-monoidal}
In the context of \cref{def:rhoc-monoidal}, the triple
\[\big(\vrho_\C, \vrho_C^2, \vrho_\C^0\big) \cn \C \to \Fr\; \End(\C)\]
is a symmetric monoidal functor.
\end{lemma}

\begin{remark}[Not Strictly Unital]\label{rk:rhonotunital}
The symmetric monoidal functor $\left(\vrho_\C, \vrho_C^2, \vrho_\C^0\right)$ in \cref{rhoc-monoidal} is \emph{neither} strictly unital nor strong because the unit constraint, which is the morphism 
\[\vrho_\C^0 = \srb{\iota_0, 1_\pu} \cn \ang{} \to (\pu)\]
in \cref{rhoc-unit}, is not an isomorphism in $\Fr\, \End(\C)$.  This stands in stark contrast with its left adjoint, $\epz_\C$ in \cref{definition:epz}, which is a \emph{strict} symmetric monoidal functor.
\end{remark}

\section[Non-Symmetric Multifunctor]{Free Permutative Category as a Non-Symmetric \texorpdfstring{$\Cat$}{Cat}-Multifunctor}
\label{sec:Fmultifunctor}

In this section we extend the 2-functor in \cref{proposition:free-perm-functor} 
\[\Fr \cn \Multicat \to \permcatst\] 
to a \emph{non-symmetric} $\Cat$-multifunctor
\[\Fr \cn \Multicat \to \permcatsu\]
in the sense of \cref{def:enr-multicategory-functor} with $(\V,\otimes) = (\Cat,\times)$.  This is one of the main results in \cite{johnson-yau-Fmulti}.
\begin{itemize}
\item The domain $\Multicat$ is the $\Cat$-multicategory in \cref{expl:multicatcatmulticat}.  It is induced by a symmetric monoidal $\Cat$-category structure (\cref{theorem:multicat-symmon}), with the tensor product $\otimes$ in \cref{definition:multicat-tensor}.
\item The codomain $\permcatsu$ is the $\Cat$-multicategory in \cref{thm:permcatmulticat}.  Its multimorphism categories have multilinear functors as objects and multilinear transformations as morphisms (\cref{definition:permcatsus-homcat}).
\end{itemize}
Here is an outline of this section.
\begin{itemize}
\item The non-symmetric $\Cat$-multifunctor $\Fr$ involves strong $n$-linear functors $\Frn$ in \cref{def:S-multi}.  The definition of each $\Frn$, in turn, requires some auxiliary constructions involving tuples in \cref{definition:xonen,definition:phionen}.
\item The multimorphism functors of $\Fr$ are in \cref{definition:F-multi}.
\end{itemize}

\subsection*{Tensor Product of Tuples}

Suppose $\ang{\M_i}_{i=1}^n$ and $\N$ are small multicategories for some $n \geq 0$.    Recall from \cref{definition:free-perm} that objects in the free permutative category $\Fr\N$ are finite tuples of objects in $\N$.

\begin{definition}[Tensor Product of Tuples of Objects]\label{definition:xonen}
Suppose given objects
\begin{equation}\label{angxi-oneton}
\ang{x^i} = \ang{x^i_j}_{j=1}^{r_i} \in \Fr\M_i \forspace i \in \{1,\ldots,n\}
\end{equation}
with each $x^i_j$ an object in $\M_i$. 
\begin{itemize}
\item For each $n$-tuple of indices $\ang{j_i}_{i=1}^n$ with each $j_i \in \{1,\ldots,r_i\}$, we define the object
\begin{equation}\label{xonenjonen-def}
x^{\onen}_{\jonejn} = \bang{x^i_{j_i}}_{i=1}^n \inspace \txotimes_{i=1}^n \M_i
\end{equation}
using the canonical bijection \cref{ObMtensorN} for objects in the tensor product.
\item We define an object
\begin{equation}\label{angxonen-def}
\ang{x^{\onen}} = \txotimes_{i=1}^n \ang{x^i} \inspace \Fr\big(\txotimes_{i=1}^n \M_i \big)
\end{equation}
using the tensor product of tuples in \cref{definition:xiotimes}.
\end{itemize} 
In other words, with
\begin{equation}\label{ronen}
r_{\onen} = \txprod_{i=1}^n r_i,
\end{equation}
the object $\ang{x^{\onen}}$ in \cref{angxonen-def} is the $r_{\onen}$-tuple 
\begin{equation}\label{eq:xonen}
\ang{x^{\onen}} = \lrang{ \,\cdots\, \lrang{x^{\onen}_{\jonejn}}_{j_1=1}^{r_1} \,\cdots\, }_{j_n = 1}^{r_n}
\end{equation}
with each entry as in \cref{xonenjonen-def}. 
\end{definition}

The following definition extends the construction $\ang{x^{\onen}}$ to morphisms.

\begin{definition}[Tensor Product of Morphisms]\label{definition:phionen}
Suppose given morphisms
\begin{equation}\label{fiangphi-FMi}
\brb{f^i,\ang{\phi^i}} \cn 
\ang{x^i} = \ang{x^i_j}_{j=1}^{r_i} \to \ang{y^i} = \ang{y^i_k}_{k=1}^{s_i} 
\inspace \Fr\M_i
\end{equation}
for $i \in \{1,\ldots,n\}$ with
\begin{itemize}
\item each $f^i \cn \ufs{r_i} \to \ufs{s_i}$ an index map,
\item each $\ang{\phi^i} = \ang{\phi^i_k}_{k=1}^{s_i}$ an $s_i$-tuple, and
\item each
\[\phi^i_k \in \M_i\scmap{\ang{x^i_j}_{j \in (f^i)^\inv(k)}; y^i_k}\]
an $|(f^i)^\inv(k)|$-ary multimorphism in $\M_i$.
\end{itemize}
We define the following objects and (multi)morphisms.
\begin{itemize}
\item First we define an index map $f^{\onen}$ as the composite function
\begin{equation}\label{eq:fonen}
\begin{tikzpicture}[baseline={(a.base)}]
\def\h{2} \def\v{.6}
\draw[0cell]
(0,0) node (a) {\ufs{r_{\onen}}}
(a)+(\h,0) node (b) {\binprod_{i=1}^n \, \ufs{r_i}}
(b)+(3,0) node (c) {\binprod_{i=1}^n \, \ufs{s_i}}
(c)+(\h,0) node (d) {\ufs{s_{\onen}}.}
;
\draw[1cell=.9]
(a) edge node {\iso} (b)
(c) edge node {\iso} (d)
(b) edge node {\binprod_{i=1}^n \, f^i} (c)
;
\draw[1cell]
(a) [rounded corners=3pt] |- ($(b)+(0,\v)$)
-- node {f^{\onen}} ($(c)+(0,\v)$) -| (d)
;
\end{tikzpicture}
\end{equation} 
In \cref{eq:fonen},
\begin{itemize}
\item the integers 
\[r_{\onen} = \txprod_{i=1}^n r_i \andspace 
s_{\onen} = \txprod_{i=1}^n s_i\]
are as in \cref{ronen}, and
\item the two unlabeled isomorphisms are given by the reverse lexicographic ordering of the products.
\end{itemize} 
\item For each $n$-tuple of indices $\ang{k_i}_{i=1}^n$ with each $k_i \in \{1,\ldots,s_i\}$, we define the object in $\Fr\big(\txotimes_{i=1}^n \M_i \big)$
\begin{equation}\label{angxonen-fkonen}
\begin{aligned}
\ang{x^{\onen}}_{f;\,\konekn}
&= \txotimes_{i=1}^n \ang{x^i_j}_{j \in (f^i)^\inv(k_i)} \\
&= \lrang{\,\cdots\, \lrang{ x^{\onen}_{\jonejn} }_{j_1 \in (f^1)^\inv(k_1)} \,\cdots\, }_{j_n \in (f^n)^\inv(k_n)}
\end{aligned}
\end{equation}
as a sub-tuple of $\ang{x^{\onen}}$ in \cref{eq:xonen}.  Then we define the multimorphism
\begin{equation}\label{eq:phionenkonekn}
\phi^{\onen}_{\konekn} = \txotimes_{i=1}^n \phi^i_{k_i} 
\cn \ang{x^{\onen}}_{f;\,\konekn} \to y^{\onen}_{\konekn} 
\inspace \txotimes_{i=1}^n \M_i.
\end{equation} 
\item We define the $s_{\onen}$-tuple of multimorphisms
\begin{equation}\label{eq:phionen}
\begin{aligned}
\ang{\phi^{\onen}} 
&= \txotimes_{i=1}^n \ang{\phi^i}\\
&= \lrang{\,\cdots\, \lrang{ \phi^{\onen}_{\konekn} }_{k_1 = 1}^{s_1} \,\cdots\, }_{k_n = 1}^{s_n}
\end{aligned}
\end{equation}
with each entry as in \cref{eq:phionenkonekn}.
\item Finally, we define the morphism
\begin{equation}\label{fonen-angphionen}
\brb{f^{\onen},\ang{\phi^{\onen}}} \cn \ang{x^{\onen}} \to \ang{y^{\onen}}
\inspace \Fr\big(\txotimes_{i=1}^n\M_i\big)
\end{equation}
with
\begin{itemize}
\item the objects $\ang{x^{\onen}}$ and $\ang{y^{\onen}}$ as in \cref{angxonen-def},
\item the index map $f^{\onen}$ in \cref{eq:fonen}, and
\item $\ang{\phi^{\onen}}$ the $s_{\onen}$-tuple in \cref{eq:phionen}. 
\end{itemize} 
\end{itemize}
This finishes the definition of $\brb{f^{\onen},\ang{\phi^\onen}}$.
\end{definition}

\subsection*{The Strong Multilinear Functor $\Frn$}
Recall 
\begin{itemize}
\item from \cref{def:nlinearfunctor} the notion of an \emph{$n$-linear functor} and
\item from \cref{example:free-Mtu} the free permutative category $\Fr(\Mtu)$ of the initial operad $\Mtu$.
\end{itemize}

\begin{definition}[Multilinear Functor $\Frn$]\label{def:S-multi}
Suppose $\ang{\M_i}_{i=1}^n$ are small multicategories.  We define the data of an $n$-linear functor
\[\brb{\Frn, \bang{(\Frn)^2_p}_{p=1}^n} \cn 
\txprod_{i=1}^n \Fr\M_i \to \Fr\big( \txotimes_{i=1}^n \M_i \big)\]
as follows.  For $n = 0$, we define the 0-linear functor\label{not:Frzero}
\[\Fr^0 \cn \boldone \to \Fr(\Mtu)\]
by the choice of the length-one tuple $(*) \in \Fr(\Mtu)$.

Suppose $n > 0$ for the rest of this definition.  Suppose given morphisms 
\[\brb{f^i,\ang{\phi^i}} \cn 
\ang{x^i} = \ang{x^i_j}_{j=1}^{r_i} \to \ang{y^i} = \ang{y^i_k}_{k=1}^{s_i} 
\inspace \Fr\M_i\]
for $i \in \{1,\ldots,n\}$ as in \cref{fiangphi-FMi}.
\begin{description}
\item[Object Assignment] 
We define the object 
\begin{equation}\label{eq:Sxi}
\Frn \bang{\ang{x^i}}_{i=1}^n = \ang{x^\onen} \inspace \Fr\big( \txotimes_{i=1}^n \M_i \big)
\end{equation}
using \cref{angxonen-def}.
\item[Morphism Assignment]
We define the morphism
\begin{equation}\label{eq:Sfi}
\Frn\lrang{\brb{f^i,\ang{\phi^i}}}_{i=1}^n = \brb{f^\onen,\ang{\phi^\onen}} \inspace \Fr\big( \txotimes_{i=1}^n \M_i \big)
\end{equation}
using \cref{fonen-angphionen}.
\item[Linearity Constraints]
For $p \in \{1,\ldots,n\}$, suppose given an object
\[\ang{\hat{x}^p} = \bang{\hat{x}^p_j}_{j=1}^{\hat{r}_p} \inspace \Fr\M_p\]
with each $\hat{x}^p_j$ an object in $\M_p$.  We first define the object
\[\ang{\tilde{x}^p} = \ang{x^p} \oplus \ang{\hat{x}^p} \inspace \Fr\M_p\]
with length $r_p + \hat{r}_p$.  Then we define the objects
\begin{equation}\label{eq:xhatonen-xtilonen}
\ang{\hat{x}^\onen} \andspace \ang{\tilde{x}^\onen} 
\inspace \Fr\big( \txotimes_{i=1}^n \M_i \big)
\end{equation}
as in \cref{angxonen-def}, using $\ang{\hat{x}^p}$ and $\ang{\tilde{x}^p}$, respectively, in place of $\ang{x^p}$.

The $p$th linearity constraint, $(\Frn)^2_p$, is defined by component isomorphisms in $\Fr\big( \txotimes_{i=1}^n \M_i \big)$
\begin{equation}\label{eq:S2b}
(\Frn)^2_p = \brb{\rho_{r_p,\hat{r}_p}, \ang{1}} \cn 
\ang{x^\onen} \oplus \ang{\hat{x}^\onen} \fto{\iso} \ang{\tilde{x}^\onen}.
\end{equation} 
In \cref{eq:S2b} the two components of $(\Frn)^2_p$ are as follows.
\begin{itemize}
\item The first component 
\begin{equation}\label{rhorbrhatb}
\rho_{r_p,\hat{r}_p} \in \Sigma_{r_1 \cdots (r_p + \hat{r}_p) \cdots r_n}
\end{equation}
is the unique permutation of entries determined by the domain and codomain of $(\Frn)^2_p$.
\item Each entry in $\ang{1}$ is a colored unit of an entry in either $\ang{x^\onen}$ or $\ang{\hat{x}^\onen}$.
\end{itemize}
\end{description}
This finishes the definition of $\brb{\Frn, \bang{(\Frn)^2_p}_{p=1}^n}$.
\end{definition}

\begin{remark}\label{rk:Snequalsone}\
\begin{enumerate}
\item\label{Snequalsone-i}
  For $n=1$, the 1-linear functor $\Fr^1$ is the identity symmetric monoidal functor on $\Fr\M_1$.
\item\label{Snequalsone-ii}
  For $n \geq 2$, the permutation $\rho_{r_p,\hat{r}_p}$ in \cref{rhorbrhatb} is the identity if $p=n$, but it is not the identity in general.
\item\label{Snequalsone-iii}
  In \cite[Section~7]{johnson-yau-Fmulti}, the multilinear functors $(\Frn,(\Frn)^2_p)$ are denoted $(\Sr,\Sr^2_p)$. \defmark
\end{enumerate}
\end{remark}

The following result combines \cite[7.12, 7.14, and 7.16]{johnson-yau-Fmulti}.

\begin{proposition}\label{S-multifunctorial}
Suppose $\ang{\M_i}_{i=1}^n$ are small multicategories.  Then the following statements hold.
\begin{enumerate}
\item\label{S-multi-i} 
The data in \cref{def:S-multi},
\[\brb{\Frn, \bang{(\Frn)^2_p}_{p=1}^n} \cn 
\txprod_{i=1}^n \Fr\M_i \to \Fr\big( \txotimes_{i=1}^n \M_i \big),\]
form a strong $n$-linear functor.
\item\label{S-multi-ii}
Each $\Frn$ is 2-natural with respect to multifunctors and multinatural transformations.
\end{enumerate}
\end{proposition}

\begin{explanation}[2-Naturality of $\Frn$]\label{expl:S-twonatural}
In \cref{S-multifunctorial} \eqref{S-multi-ii}, the 2-naturality of $\Frn$ with respect to multifunctors between small multicategories
\[H_i \cn \M_i \to \N_i \forspace i \in \{1,\ldots,n\}\]
means that the following two composite $n$-linear functors are equal.
\begin{equation}\label{S-nat-multif}
\begin{tikzpicture}[x=40mm,y=15mm,vcenter]
\draw[0cell=.9]
(0,0) node (a) {\txprod_{i=1}^n \Fr\M_i}
(0,-1) node (b) {\Fr\big(\txotimes_{i=1}^n \M_i \big)}
(1,0) node (a') {\txprod_{i=1}^n \Fr\N_i}
(1,-1) node (b') {\Fr\big(\txotimes_{i=1}^n \N_i \big)}
;
\draw[1cell=.9] 
(a) edge node {\txprod_{i=1}^n \Fr H_i} (a')
(b) edge node {\Fr(\txotimes_{i=1}^n H_i)} (b')
(a) edge['] node {\Frn} (b)
(a') edge node {\Frn} (b')
;
\end{tikzpicture}
\end{equation}
The 2-naturality of $\Frn$ with respect to multinatural transformations
\[\theta_i \cn H_i \to K_i \forspace i \in \{1,\ldots,n\}\]
means the following equality of $n$-linear transformations.  
\begin{equation}\label{S-nat-tr}
1_{(\Frn)} * \big( \txprod_{i=1}^n \Fr \theta_i \big) 
= \Fr \big( \txotimes_{i=1}^n \theta_i \big) * 1_{(\Frn)}
\end{equation}
This equality is obtained from the diagram \cref{S-nat-multif} by replacing each $H_i$ with $\theta_i$.
\end{explanation}

\subsection*{The Non-Symmetric $\Cat$-Multifunctor $\Fr$}

Next we extend the 2-functor $\Fr$ in \cref{proposition:free-perm-functor} to multimorphism categories.

\begin{convention}\label{convention:Fun}
To avoid confusion in \cref{definition:F-multi} below, for small multicategories $\M$ and $\N$, we denote by
\[\Frbar \cn \Multicat(\M,\N) \to \permcatst(\Fr\M,\Fr\N)\] 
the assignment of $\Fr$ on multifunctors and multinatural transformations as in \cref{definition:free-smfun,definition:free-perm-multinat}, respectively.
\end{convention}

In \cref{eq:FbarS} below, we use the multilinear functor $\Frn$ (\cref{S-multifunctorial}).

\begin{definition}\label{definition:F-multi}
Suppose $\angM = \ang{\M_i}_{i=1}^n$ and $\N$ are small multicategories.  We define a functor between multimorphism categories
\begin{equation}\label{F-multimorphism-cat}
\Fr\cn \Multicat\scmap{\ang{\M};\N} \to \permcatsu\scmap{\ang{\Fr\M};\Fr\N}
\end{equation}
as follows.  Suppose given multifunctors $H$ and $K$ and a multinatural transformation $\theta$ as in the diagram below.
\[\begin{tikzpicture}[baseline={(a.base)}]
\def\t{22}
\draw[0cell]
(0,0) node (a) {\ang{\M}}
(a)++(2,0) node (b) {\N}
;
\draw[1cell=.9]  
(a) edge[bend left=\t] node[pos=.43] {H} (b)
(a) edge[bend right=\t] node[swap,pos=.43] {K} (b)
;
\draw[2cell] 
node[between=a and b at .47, rotate=-90, 2label={above,\,\theta}] {\Rightarrow}
;
\end{tikzpicture}\]
Then $\Fr$ sends these data to the following composite $n$-linear functors and whiskering.
\begin{equation}\label{eq:FbarS}
\begin{tikzpicture}[baseline={(a.base)}]
\def\t{20}
\draw[0cell]
(-2.75,0) node (z) {\ang{\Fr\M}}
(0,0) node (a') {\Fr\big(\txotimes_{i=1}^n \M_i\big)}
(a')+(.65,0) node (a) {\phantom{\Fr\N}}
(a)++(2.5,0) node (b) {\Fr\N}
;
\draw[1cell=.9]  
(z) edge node {\Frn} (a')
(a) edge[bend left=\t] node[pos=.5] {\Frbar H} (b)
(a) edge[bend right=\t] node[swap,pos=.5] {\Frbar K} (b)
;
\draw[2cell] 
node[between=a and b at .42, rotate=-90, 2label={above,\,\Frbar \theta}] {\Rightarrow}
;
\end{tikzpicture}
\end{equation}
This finishes the definition of the multimorphism functor $\Fr$.
\end{definition}

\begin{explanation}\label{expl:Frnequalsone}
In the case $n=1$, the 1-linear functor $\Fr^1$ is the identity symmetric monoidal functor (\cref{rk:Snequalsone} \eqref{Snequalsone-i}).  The multimorphism functor $\Fr$ in \cref{F-multimorphism-cat} reduces to the hom functor of the 2-functor $\Fr$ in \cref{proposition:free-perm-functor}.  Therefore, there is no ambiguity in reusing the notation $\Fr$ in \cref{definition:F-multi}.
\end{explanation}

\begin{explanation}\label{expl:Frnequalszero}
  In the case $n=0$, recall from \cref{expl:multicatcatmulticat} that a nullary operation $H \in \Multicat\scmap{\ang{};\N}$ consists of a choice of object $H* \in \N$, where $*$ is the unique object of $\Mtu$.
  The 0-linear functor $\Fr^0 \cn \boldone \to \Fr(\Mtu)$ chooses the length-one tuple $(*) \in \Fr(\Mtu)$, and hence $\Fr H = (\ol{\Fr}H) \circ \Fr^0$ is determined by the length-one tuple $(H*) \in \Fr\N$.
\end{explanation}

Recall from \cref{def:enr-multicategory-functor} that a \emph{non-symmetric $\Cat$-multifunctor} between $\Cat$-multicategories preserve colored units and composition, but it is \emph{not} required to preserve the symmetric group action as in \cref{enr-multifunctor-equivariance}.  The following result is \cite[8.1]{johnson-yau-Fmulti}.

\begin{theorem}\label{theorem:F-multi}\index{category!free permutative - multifunctor}\index{permutative category!free - multifunctor}\index{free!permutative category multifunctor}
There is a non-symmetric $\Cat$-multifunctor
\begin{equation}\label{eq:thm-F-multi}
  \Fr \cn \Multicat \to \permcatsu
\end{equation} 
defined by the following data.
\begin{itemize}
\item $\Multicat$ is the $\Cat$-multicategory in \cref{expl:multicatcatmulticat}.
\item $\permcatsu$ is the $\Cat$-multicategory in \cref{thm:permcatmulticat}.
\item The object assignment of $\Fr$ is the free permutative category in \cref{proposition:free-perm}. 
\item The multimorphism functors of $\Fr$ are in \cref{F-multimorphism-cat}.
\end{itemize}
\end{theorem}

\begin{explanation}\label{expl:Fnotsymmetric}
Consider \cref{theorem:F-multi}.
\begin{enumerate}
\item\label{Fnotsym-i}
As explained in \cite[8.2]{johnson-yau-Fmulti}, the non-symmetry of $\Fr$ in \cref{eq:FbarS} is due to the incompatibility of $\Frn$ with permutations.  
\item\label{Fnotsym-ii}
As in \cref{definition:free-smfun}, $\Frbar$ sends each multifunctor to a \emph{strict} symmetric monoidal functor.  Moreover, by \cref{S-multifunctorial} \eqref{S-multi-i}, each $\Frn$ is a \emph{strong} multilinear functor.  Therefore, the composite in \cref{eq:FbarS}
\[\Fr H = \Frbar H \circ \Frn\]
is a \emph{strong} $n$-linear functor by definition \cref{ffjlinearity}.\defmark
\end{enumerate}
\end{explanation}

Recall from \cref{expl:end-catmulti} the $\Cat$-multifunctor
\[\End \cn \permcatsu \to \Multicat.\]
Also recall from \cref{expl:catmultitransformation} an explicit description of a (non-symmetric) $\Cat$-multinatural transformation.  The following result is \cite[9.2]{johnson-yau-Fmulti}, where $\Fr$ is the non-symmetric $\Cat$-multifunctor in \cref{theorem:F-multi}.

\begin{lemma}\label{lemma:eta-mnat}
The unit in \cref{theorem:FE-adj}
\[\begin{tikzpicture}[baseline={(a.base)}]
\def\t{22}
\draw[0cell]
(0,0) node (a') {\Multicat}
(a')+(.5,0) node (a) {\phantom{X}}
(a)+(2.5,0) node (b) {\phantom{X}}
(b)+(.5,0) node (b') {\Multicat}
;
\draw[1cell=.9]  
(a) edge[bend left=\t] node[pos=.5] {1_{\Multicat}} (b)
(a) edge[bend right=\t] node[swap,pos=.5] {\End\; \Fr} (b)
;
\draw[2cell] 
node[between=a and b at .45, rotate=-90, 2label={above,\eta}] {\Rightarrow}
;
\end{tikzpicture}\]
is a non-symmetric $\Cat$-multinatural transformation.
\end{lemma}

\begin{remark}[Non-Existence of Counit Analog]\label{rk:epz-not-mnat}
The counit in \cref{theorem:FE-adj} does \emph{not} yield a non-symmetric $\Cat$-multinatural transformation
\[\epz \cn \Fr\; \End \to 1_{\permcatsu}.\]
As we mentioned in \cref{rk:epznaturality}, that counit $\epz$ is only natural with respect to \emph{strict} symmetric monoidal functors but not strictly unital symmetric monoidal functors in general.  Thus, the analog of \cref{lemma:eta-mnat} does not hold for the counit.
\end{remark}

\section[Homotopy Equivalences]{Homotopy Equivalences between Multicategories and Permutative Categories}
\label{sec:multicat-model}

In this section we review equivalences of homotopy theories between the categories $\Multicat$, $\permcatst$, and $\permcatsu$.
\begin{itemize}
\item \cref{thm:F-heq}, which relates $\Multicat$ and $\permcatst$, is the main result of \cite{johnson-yau-permmult}.  The morphisms in the category $\permcatst$ are \emph{strict} symmetric monoidal functors.  The proof of this theorem relies on the componentwise right adjoint of the counit $\epz$ in \cref{proposition:epz-rho-adj}.  In \cref{ch:ptmulticat-sp} we use a pointed version of this componentwise right adjoint to extend the equivalence of homotopy theories to \emph{pointed} multicategories and left $\Mone$-modules.
\item \cref{thm:alg-hty-equiv} relates $\Q$-algebras in $\Multicat$ and $\Q$-algebras in $\permcatsu$ for a small non-symmetric $\Cat$-multicategory $\Q$.  It is the main result of \cite{johnson-yau-Fmulti}.  In the larger category $\permcatsu$, the morphisms are \emph{strictly unital} symmetric monoidal functors.  In \cref{ch:ptmulticat-alg} we extend this equivalence of homotopy theories between categories of non-symmetric algebras, as well as the $\Cat$-multifunctoriality of $\Fr$ (non-symmetric in the case of $\Fr$) and $\End$, to pointed multicategories and left $\Mone$-modules.
\item \cref{thm:Fsu-heq}, which relates $\Multicat$ and $\permcatsu$, is the special case of \cref{thm:alg-hty-equiv} when $\Q$ is the initial operad.  
\item A consequence of \cref{thm:F-heq,thm:Fsu-heq} is that the inclusion functor from $\permcatst$ to $\permcatsu$ is an equivalence of homotopy theories.  See \cref{cor:I-heq}.
\end{itemize}

\subsection*{Stable Equivalences}

Recall from \cref{sec:segalEMK} that Segal $K$-theory 
\[\Kse \cn \brb{\permcatsu,\cS} \fto{\sim} \brb{\Spc,\cS}\]
is an equivalence of homotopy theories with the following relative category structures.  Recall that a subcategory is \emph{wide} if it contains all the objects of the larger category.
\begin{itemize}
\item In the codomain, the wide subcategory 
\[\cS \bigsubset \Spc\]
consists of stable equivalences of connective symmetric spectra.
\item In the domain, the wide subcategory of stable equivalences 
\[\cS \bigsubset \permcatsu\] 
is created by $\Kse$.  In other words, a \index{permutative category!stable equivalence}\index{equivalence!stable!of permutative categories}\index{stable equivalence!of permutative categories}\emph{stable equivalence} in $\permcatsu$ is a strictly unital symmetric monoidal functor $P$ between small permutative categories such that $\Kse P \in \Spc$ is a stable equivalence of connective symmetric spectra.
\end{itemize}

To relate the homotopy theories of $\Multicat$, $\permcatst$, and $\permcatsu$, first we specify classes of stable equivalences in $\Multicat$ and $\permcatst$.  

\begin{definition}[Stable Equivalences]\label{def:mult-stableeq}
We define the wide subcategories
\[\begin{split}
\cSI & = I^\inv(\cS) \bigsubset \permcatst \andspace \\
\cSF & = \Fr^\inv(\cSI) \bigsubset \Multicat
\end{split}\]
as the subcategories created by the indicated functors below.
\begin{equation}\label{FI}
\begin{tikzcd}[column sep=normal]
\brb{\Multicat,\cSF} \ar{r}{\Fr} & \brb{\permcatst,\cSI} \ar[hookrightarrow]{r}{I} & \brb{\permcatsu,\cS}
\end{tikzcd}
\end{equation}
\begin{itemize}
\item $\Fr$ is the underlying functor of the 2-functor in \cref{proposition:free-perm-functor}.
\item $I$ is the underlying functor of the inclusion 2-functor in \cref{permcatinclusion}. 
\end{itemize}
We call morphisms in $\cSI$ and $\cSF$ \index{permutative category!stable equivalence}\index{equivalence!stable!of permutative categories}\index{stable equivalence!of permutative categories}\emph{stable equivalences} and \emph{$\Fr$-stable equivalences}, respectively.
\end{definition}

\subsection*{Equivalences of Homotopy Theories}

Recall from \cref{def:heq} the notion of an \index{equivalence!adjoint - of homotopy theories}\index{homotopy theory!adjoint equivalence of}\index{adjoint equivalence!of homotopy theories}\emph{adjoint equivalence of homotopy theories}.  The following result is \cite[7.3]{johnson-yau-permmult}.  Its proof makes crucial use of \cref{proposition:epz-rho-adj} on the componentwise right adjoint of the counit.

\begin{theorem}\label{thm:F-heq}
The adjunction in \cref{theorem:FE-adj}
\[\begin{tikzpicture}
\draw[0cell]
(0,0) node (a') {\brb{\Multicat, \cSF}}
(a')+(1,0) node (a) {\phantom{X}}
(a)+(1,0) node (x) {\bot}
(a)+(2,0) node (b) {\phantom{X}}
(b)+(1.1,0) node (b') {\brb{\permcatst, \cSI}}
;
\draw[1cell=.9]
(a) edge[bend left=18] node {\Fr} (b)
(b) edge[bend left=18] node {\End} (a)
;
\end{tikzpicture}\]
is an adjoint equivalence of homotopy theories.
\end{theorem}

In particular, each of the unit and the counit of $\Fr \dashv \End$ in \cref{thm:F-heq} is a relative natural transformation with respect to $\cSF$ and $\cSI$, respectively.

The equivalences of homotopy theories---but not the adjunction---in \cref{thm:F-heq} can be extended to the larger category $\permcatsu$, together with algebraic structures in the following sense.  

\begin{definition}[Symmetric and Non-Symmetric Algebras]\label{def:nonsymalgebra}
Suppose $\P$ is a small $\Cat$-multicategory, and $\N$ is a $\Cat$-multicategory.
\begin{itemize}
\item A \index{algebra!over a multicategory}\emph{$\P$-algebra} in $\N$ is a $\Cat$-multifunctor
\[\P \to \N.\]
\item A \emph{morphism} of $\P$-algebras in $\N$ is a $\Cat$-multinatural transformation (\cref{expl:catmultitransformation}).
\item The category of $\P$-algebras and their morphisms in $\N$ is denoted $\N^\P$.\label{not:NP}
\end{itemize}

There is also a non-symmetric analog, for which we use the same notation and similar terminology.
Suppose $\Q$ is a small non-symmetric $\Cat$-multicategory, and $\N$ as above is a $\Cat$-multicategory.
\begin{itemize}
\item A \index{non-symmetric!algebra!over a multicategory}\index{algebra!non-symmetric - over a multicategory}\emph{non-symmetric $\Q$-algebra} in $\N$ is a non-symmetric $\Cat$-multifunctor
  \[\Q \to \N.\]
\item A \emph{morphism} of non-symmetric $\Q$-algebras in $\N$ is a non-symmetric $\Cat$-multinatural transformation.
\item The category of non-symmetric $\Q$-algebras and their morphisms in $\N$ is denoted $\N^\Q$.\label{not:NQ}
\end{itemize}

If the underlying 1-category of $\N$ is a relative category $(\N,\cW)$, then the categories of $\P$-algebras, respectively non-symmetric $\Q$-algebras, have an induced relative structure.
\begin{itemize}
\item We define the wide subcategories\label{not:cWP}
  \[
    \cW^\P \bigsubset \N^\P \andspace
    \cW^\Q \bigsubset \N^\Q
  \]
  to be those that contain all the morphisms with each component in $\cW$. 
\item We consider the pairs $(\N^\P, \cW^\P)$ and $(\N^\Q, \cW^\Q)$ as relative categories.
\end{itemize} 
This finishes the definition.
\end{definition}

Recall the $\Cat$-multifunctors (non-symmetric for $\Fr$) 
\[\Fr \cn \Multicat \lradj \permcatsu \cn \End\]
in \cref{theorem:F-multi,expl:end-catmulti}.  \cref{thm:alg-hty-equiv} below extends the equivalences of homotopy theories in \cref{thm:F-heq} to non-symmetric algebras.  It is the main result in \cite[1.1]{johnson-yau-Fmulti}.  

\begin{theorem}\label{thm:alg-hty-equiv}
Suppose $\Q$ is a small non-symmetric $\Cat$-multicategory.  Then the functors  
\[\Fr^\Q \cn \brb{\Multicat^\Q, (\cSF)^\Q} \lrsimadj \brb{(\permcatsu)^\Q, \cS^\Q} 
\cn \End^\Q,\]
induced by post-composition and whiskering with, respectively, $\Fr$ and $\End$, are inverse equivalences of homotopy theories in the sense of \cref{def:inverse-heq}.
\end{theorem}

Thus, by \cref{gjo29}, $\Fr^\Q$ and $\End^\Q$ are equivalences of homotopy theories between categories of non-symmetric $\Q$-algebras.

\begin{remark}[Subtleties]\label{rk:FE-not-adjunction}
There are two main subtleties of \cref{thm:alg-hty-equiv}.
\begin{enumerate}
\item\label{rk:FE-not-adjunction-i} 
Unlike \cref{thm:F-heq}, the two functors in \cref{thm:alg-hty-equiv} do \emph{not} form an adjunction in general, as discussed in \cref{rk:epznaturality}.  Nevertheless, the unit
\[\eta \cn 1_{\Multicat} \to \End\; \Fr\]
in \cref{theorem:FE-adj} is
\begin{itemize}
\item componentwise an $\Fr$-stable equivalence by \cref{thm:F-heq} and
\item a non-symmetric $\Cat$-multinatural transformation by \cref{lemma:eta-mnat}.
\end{itemize}
\item\label{rk:FE-not-adjunction-ii} 
For a permutative category $\C$, the functor in \cref{epzc-rhoc} 
\[\vrho_\C \cn \C \to \Fr \; \End(\C)\]
is \emph{not} strictly unital, as discussed in \cref{rk:rhonotunital}.  Therefore, we cannot use $\vrho$ directly to compare the identity functor on $(\permcatsu)^\Q$ and the composite $\Fr^\Q\, \End^\Q$.  To use \cref{gjo29}, the proof in \cite{johnson-yau-Fmulti} compares these two functors via a zigzag of relative natural transformations.
\end{enumerate}
As we discuss in \cref{ch:ptmulticat-sp}, the pointed variant of $\vrho_\C$ \emph{is} strictly unital.  In this sense, the pointed variant of \cref{thm:alg-hty-equiv} is more natural than the unpointed version here.
\end{remark}

Taking $\Q$ as the initial operad in \cref{ex:vmulticatinitialterminal} \cref{ex:initialoperad}, \cref{thm:alg-hty-equiv} yields the following special case.

\begin{theorem}\label{thm:Fsu-heq}
The relative functors
\[\Fr \cn \brb{\Multicat,\cSF} \lrsimadj \brb{\permcatsu,\cS} \cn \End\]
are equivalences of homotopy theories.
\end{theorem}

\cref{cor:I-heq} below follows from
\begin{itemize}
\item the commutative diagram 
\begin{equation}\label{Fpermcatsu}
\begin{tikzpicture}[baseline={(a.base)}]
\def\u{.6}
\draw[0cell]
(0,0) node (a) {\Multicat}
(a)+(3,0) node (b) {\permcatst}
(b)+(3,0) node (c) {\permcatsu,}
;
\draw[1cell=.9]
(a) edge node {\Fr} (b)
(b) edge[right hook->] node {I} (c)
;
\draw[1cell=.9]
(a) [rounded corners=3pt] |- ($(b)+(-1,\u)$)
-- node {\Fr} ($(b)+(1,\u)$) -| (c)
;
\end{tikzpicture}
\end{equation}
\item \cref{thm:F-heq,thm:Fsu-heq}, and 
\item the 2-out-of-3 property of Rezk weak equivalences (\cref{theorem:css-fibrant}).
\end{itemize} 

\begin{corollary}\label{cor:I-heq}
The inclusion relative functor in \cref{FI}
\[\begin{tikzcd}[column sep=large]
\brb{\permcatst,\cSI} \ar[hookrightarrow]{r}{I}[swap]{\sim} & \brb{\permcatsu,\cS}
\end{tikzcd}\]
is an equivalence of homotopy theories.
\end{corollary}

Combining \cref{Kse-heq,thm:F-heq,thm:Fsu-heq,cor:I-heq}, we see that each arrow in the diagram
\begin{equation}\label{FEKse-heq}
\begin{tikzpicture}[baseline={(a.base)}]
\def\u{.7}
\draw[0cell=.8]
(0,0) node (a') {\brb{\Multicat, \cSF}}
(a')+(.8,0) node (a) {\phantom{X}}
(a)+(.7,0) node (x) {\bot}
(a)+(1.4,0) node (b) {\phantom{X}}
(b)+(.9,0) node (b') {\brb{\permcatst, \cSI}}
(b')+(3,0) node (c) {\brb{\permcatsu, \cS}}
(c)+(2.7,0) node (d) {\brb{\Spc,\cS}}
;
\draw[1cell=.8]
(a) edge[bend left=18] node {\Fr} (b)
(b) edge[bend left=18] node {\End} (a)
(b') edge[right hook->] node {I} (c)
(c) edge node {\Kse} (d)
;
\draw[1cell=.8]
(a') [rounded corners=3pt] |- ($(b')+(-1,\u)$)
-- node {\Fr} ($(b')+(1,\u)$) -| (c)
;
\draw[1cell=.8]
(c) [rounded corners=3pt] |- ($(b')+(1,-\u)$)
-- node[swap] {\End} ($(b')+(-1,-\u)$) -| (a')
;
\end{tikzpicture}
\end{equation}
\bigskip
is an equivalence of homotopy theories.

\part{Homotopy Theory of Pointed Multicategories, \texorpdfstring{$\Mone$}{M1}-Modules, and Permutative Categories}
\label{part:multicat}

\chapter{Pointed Multicategories and \texorpdfstring{$\Mone$}{M1}-Modules Model All Connective Spectra}
\label{ch:ptmulticat-sp}
The title of this chapter refers to equivalences of homotopy theories
\[
  \begin{aligned}
    \Kse \circ I \circ \Fst & \cn \brb{\pMulticat,\cSst} \fto{\sim} \brb{\Spc,\cS} \andspace\\
    \Kse \circ I \circ \Fm & \cn \brb{\MoneMod,\cSM} \fto{\sim} \brb{\Spc,\cS}
  \end{aligned}
\]
that are established below.
The main results and their context are summarized in the following diagram of 2-adjunctions and functors.
\begin{equation}\label{ptmulticat-page-v-summary}
\begin{tikzpicture}[baseline={(a.base)}]
\def\s{10} \def\t{18} \def\h{3.5}
\draw[0cell=.8]
(0,0) node (a) {\brb{\pMulticat,\cSst}}
(a)+(\h,0) node (b) {\brb{\permcatst,\cSI}}
(a)+(\h/2,1.75) node (c) {\brb{\Multicat,\cSF}}
(a)+(\h/2,-1.75) node (d) {\brb{\MoneMod,\cSM}}
(b)+(3,0) node (e) {\brb{\permcatsu,\cS}}
(e)+(2.5,0) node (f) {\brb{\Spc,\cS}}
;
\draw[1cell=.8]
(a) edge[transform canvas={yshift=-.7ex},bend left=\s] node {\Fst} (b)
(b) edge[transform canvas={yshift=-.5ex}] node {\Endst} (a)
(c) edge[transform canvas={xshift=2ex},bend left=\t] node[pos=.4] {\Fr} (b)
(b) edge[transform canvas={xshift=1.5ex}] node[pos=.6] {\End} (c)
(c) edge[bend right=\t, transform canvas={xshift=-1.6ex}] node[swap] {\dplus} (a)
(a) edge[transform canvas={xshift=-1ex}] node[swap,pos=.6] {\Ust} (c)
;
\draw[1cell=.8]
(d) edge[transform canvas={xshift=2ex},bend right=\t] node[swap] {\Fm} (b)
(b) edge[transform canvas={xshift=1ex}] node[swap,pos=.6] {\Endm} (d)
(d) edge[transform canvas={xshift=-1ex}] node[swap,pos=.4] {\Um} (a)
(a) edge[transform canvas={xshift=-2ex},bend right=\t] node[swap] {\Monesma} (d)
;
\draw[1cell=.8]
(b) edge[right hook->] node {I} (e)
(e) edge node {\Kse} (f)
;
\end{tikzpicture}
\end{equation}
Each of the arrows in \cref{ptmulticat-page-v-summary}, \emph{except $\Ust$}, is a relative functor.
\begin{itemize}
\item $\Kse$ is Segal $K$-theory in \cref{eq:Ksummary,Kse}.
\item $I$ is the inclusion functor in \cref{permcatinclusion}, which is the identity on objects.
\item The 2-adjunctions
\[\dplus \dashv \Ust, \quad \Fr \dashv \End, \andspace (\Monesma) \dashv \Um\] 
are in \cref{dplustwofunctor}, \cref{theorem:FE-adj}, and \cref{MonesmaUmadj}, respectively.
\item We establish the 2-adjunctions
\[\Fst \dashv \Endst \andspace \Fm \dashv \Endm\]
in \cref{ptmulticat-thm-v,ptmulticat-thm-ii}, respectively, below.
\item The diagram involving the right adjoints 
\[\End = \Ust \circ \Endst = \Ust \circ \Um \circ \Endm\] 
is the restriction of the commutative diagram \cref{endufactor} to $\permcatst$ and underlying 2-functors.
\item The diagram involving the left adjoints commutes up to 2-natural isomorphisms:
\[\Fr \iso \Fst \circ \dplus \iso \Fm \circ (\Monesma) \circ \dplus.\]
\end{itemize}
Each arrow in \cref{ptmulticat-page-v-summary}, \emph{except $\Ust$}, is an equivalence of homotopy theories in the sense of \cref{definition:rel-cat-pow}~\cref{def:relcat-iv}.
\begin{itemize}
\item The inclusion functor $I$ and Segal $K$-theory $\Kse$ are equivalences of homotopy theories by \cref{cor:I-heq,Kse-heq}, respectively.
\item $\Fr \dashv \End$ is an adjoint equivalence of homotopy theories by \cref{thm:F-heq}. 
\item We show that $\dplus$ is an equivalence of homotopy theories in \cref{ptmulticat-cor-page-vi} below.
\item We establish the adjoint equivalences of homotopy theories
\[\Fst \dashv \Endst, \quad (\Monesma) \dashv \Um, \andspace \Fm \dashv \Endm\]
in \cref{ptmulticat-thm-x,ptmulticat-thm-vi,ptmulticat-thm-xi}, respectively, below.
\item We do not know whether $\Ust$ is a relative functor with respect to the subcategories $\cSst$ and $\cSF$, and thus cannot conclude that it is an equivalence of homotopy theories.
  See \cref{question:Ust-rel} for further discussion of this point.
\end{itemize}
The definition of $\Fst$ and some of its properties follow formally from the adjunction $\Fr \dashv \End$.
In particular, \cref{proposition:FstM-pushout} shows that $\Fst$ can be constructed as a 2-categorical pushout that collapses basepoint operations.
However, the adjoint equivalences of homotopy theories are slightly more subtle because our arguments depend on particular details of $\Fst$ described in \cref{sec:Fst-ex}.
\Cref{expl:rhoc-not-su,explanation:upsst} contain further comments on this point.

\subsection*{Connection with Other Chapters}
\cref{ch:ptmulticat-alg} extends $\Fst$ and $\Fm$ to $\Cat$-enriched multifunctors and extends the equivalences of homotopy theories in this chapter to categories of algebras.
This is further extended in \cref{ch:mackey_eq} to enriched diagrams and Mackey functors.
See \cref{remark:ptmulticat-prop-viii} for more detailed technical comments about these extensions.

\subsection*{Background}
The material in \cref{ch:multperm} describes the equivalences of homotopy theories given by the unpointed free construction $\Fr$ and its right 2-adjoint $\End$.
Pointed multicategories are reviewed in \cref{sec:ptmulticatclosed}.
The symmetric monoidal $\Cat$-category of $\Mone$-modules, and its relation to pointed multicategories, is described in \cref{sec:monemodules}.
Equivalences of homotopy theories are described in \cref{sec:hty-thy}.

\subsection*{Chapter Summary}
In \cref{sec:Fsttwofunctor} we define the \emph{pointed free permutative category} construction, $\Fst$, via certain equivalence relations on $\Fr$.
\Cref{sec:relating-FFst} further develops the relationships between $\Fr$ and $\Fst$.
\Cref{sec:FEst-adj} shows that $\Fst$ has a right 2-adjoint $\Endst$.
\Cref{sec:FEm-adj} explains how $\Fst$ and related constructions restrict to the sub 2-category of $\Mone$-modules.
In \cref{sec:Fst-ex} we compute several examples $\Fst\M$ for certain pointed multicategories $\M$.
\Cref{sec:epzst-right-adj} gives a componentwise right adjoint, providing a pointed analog of \cref{sec:componentwiseadjoint}.
\Cref{sec:ptmultperm-heq,sec:moneperm-heq} then establish the equivalences of homotopy theories described above.
Here is a summary table.
\reftable{.9}{
  Underlying category and permutative structure on $\Fst\M$
  & \ref{def:Fst-object} and \ref{def:Fst-permutative}
  \\ \hline
  $\Fst$ on pointed multifunctors and multinatural transformations
  & \ref{def:Fst-onecells} and \ref{def:Fst-twocells}
  \\ \hline
  $\pst_\M \cn \Fr\M \to \Fst\M$
  & \ref{def:FFst-projection}
  \\ \hline
  $\Fst$ as a 2-pushout
  & \ref{proposition:FstM-pushout}
  \\ \hline
  unit $\etast \cn 1 \to \Endst \Fst$ and counit
  $\epzst \cn \Fst \Endst \to 1$
  & \ref{def:etast} and \ref{def:epzst}
  \\ \hline
  isomorphism $\Fr \M \iso \Fst(\M_+)$
  & \ref{explanation:FFst+iso}
  \\ \hline
  2-adjunction $\Fm \dashv \Endm$
  & \ref{ptmulticat-thm-ii}
  \\ \hline
  symmetric monoidal functor $\vrhost_\C \cn \C \to \Fst\Endst\C$
  & \ref{def:vrhostC}
  \\ \hline
  stable equivalences $\cSM$ and $\cSst$
  & \ref{def:ptmulti-stableeq}
  \\ \hline
  equivalences of homotopy theories
  & \ref{ptmulticat-thm-x},
    \ref{ptmulticat-cor-page-vi},
    \ref{ptmulticat-thm-vi},
    \ref{ptmulticat-thm-xi}
  \\ \hline
}
We remind the reader of \cref{conv:universe} about universes.

\section{Pointed Free Permutative Categories}
\label{sec:Fsttwofunctor}

This section defines the pointed free construction $\Fst$ in details that will be useful below.
\Cref{sec:relating-FFst} shows that $\Fst$ is a certain pushout constructed from $\Fr$.

\subsection*{Underlying Category $\Fst\M$}
The objects and morphisms of $\Fst\M$ are determined by the following equivalence relations on objects and morphisms of $\Fr$.
\begin{definition}[Removing Basepoints]\label{definition:removing-basepts}
  Suppose $(\M,i^\M)$ is a pointed multicategory with basepoint $* = i^\M(*)$ and $n$-ary basepoint operations $\iota^n \in \M\scmap{\ang{*};*}$.
  \begin{enumerate}
  \item For each tuple of objects $\ang{x_i}_{i=1}^r \in \Fr\M$, with $x_i \in \M$, let $\ang{x}^\wedge$ be the sub-tuple consisting of non-basepoint objects, $x_i \ne * = i^\M(*)$.
  \item For each morphism
    \[
      \ang{x_i}_{i=1}^r \fto{(f, \ang{\phi})} \ang{y_j}_{j=1}^s, \inspace \Fr\M
    \]
    define\label{not:fphiprime}
    \[
      \ang{x}' \fto{(f', \ang{\phi}')} \ang{y}' \inspace \Fr\M
    \]
    as follows.
    An index $j \in s$ is called \emph{removable} if $\phi_j$ is a basepoint operation in $\M$.  That is, if
    \begin{itemize}
    \item $y_j = *$,
    \item $x_i = *$ for each $i \in f^\inv(j)$, and
    \item $\phi_j = \iota^{|f^\inv(j)|} \in \M\scmap{\ang{*};*}$.
    \end{itemize}
    An index $j$ is \emph{irremovable} if it is not removable.
    Let $\ufs{s}' \subset \ufs{s}$ be the subset consisting of irremovable $j$.
    Let $\ufs{r}' = f^\inv(\ufs{s}')$ and then define
    \begin{equation}\label{remove-bpt-ops}
      \ang{y}' = \ang{y_j}_{j \in \ufs{s}'},\quad
      \ang{x}' = \ang{x_i}_{i \in \ufs{r}'},\quad
      \ang{\phi}' = \ang{\phi_j}_{j \in \ufs{s}'},\andspace
      f' = \left.f\right|_{\ufs{r}'}.
      \dqed
    \end{equation}
 \end{enumerate}
\end{definition}

\begin{definition}[Up-To-Basepoint Equivalence]\label{definition:up-to-basept-equiv}
  Suppose $(\M,i^\M)$ is a pointed multicategory.
  Define the following equivalence relations, called \emph{up-to-basepoint equivalence} on objects and morphisms of $\Fr\M$.
  \begin{description}
  \item[Objects]
    Up-to-basepoint equivalence is denoted $\obsim$ on $\Ob(\Fr \M)$ and defined by
    \[
      \ang{x} \obsim \ang{y} \iffspace \ang{x}^\wedge = \ang{y}^\wedge.
    \]
    Thus, $\ang{x} \obsim \ang{y}$ if the tuples $\ang{x}$ and $\ang{y}$ agree up to insertion or deletion of basepoints.
    The equivalence class of $\ang{x}$ is denoted \label{not:obsimclass}$[\ang{x}]$.

  \item[Morphisms]
    An $m$-tuple of morphisms $\bigl((f_1,\ang{\phi_1}),\ldots,(f_m,\ang{\phi_m})\bigr)$ is \emph{$\obsim$-composable} if 
    \[
      \codom(f_i,\ang{\phi_i}) \obsim \dom(f_{i+1},\ang{\phi_{i+1}}) \foreachspace 1 \leq i \leq m-1.
    \]
    Let \label{not:wtmor}$\wt{\Mor}(\Fr\M) \subset \coprod_{m \ge 1} \Mor(\Fr\M)^{\times m}$ denote the collection of $\obsim$-composable tuples of morphisms.

    For $\obsim$-composable tuples of morphisms
    \[
      \uf = \bigl((f_1,\ang{\phi_1}),\ldots,(f_m,\ang{\phi_m})\bigr)
      \andspace
      \ug=\bigl((g_1,\ang{\psi_1}),\ldots,(g_n,\ang{\psi_n})\bigr),
    \]
    up-to-basepoint equivalence is denoted $\simeq$ and is generated by $\isim$ and $\iisim$ as follows.
    \begin{enumerate} 
    \item For relation one, \label{not:isim}$\uf \isim \ug$ if $m = n+1$ and there is some $1 \leq i \leq m-1$ such that
      \begin{itemize}
      \item for $j \not\in \{i,i+1\}$:
        \begin{align*}
          (f_j,\ang{\phi_j}) & = (g_j,\ang{\psi_j}) \ifspace 1 \leq j < i\\
          (f_j,\ang{\phi_j}) & = (g_{j-1},\ang{\psi_{j-1}}) \ifspace i+1 < j \leq m
        \end{align*}
      \item $\codom(f_i,\ang{\phi_i}) = \dom(f_{i+1},\ang{\phi_{i+1}})$ in $\Fr\M$, and
        \[
          (g_i, \ang{\psi_i}) = (f_{i+1},\ang{\phi_{i+1}})\circ(f_i,\ang{\phi_i}),
        \]
        the composite in $\Fr\M$ \cref{eq:FM-comp}.
      \end{itemize}

    \item For relation two, \label{not:iisim}$\uf \iisim \ug$ if $m = n$ and there is some $1 \leq i \leq m$ such that
      \begin{itemize}
      \item $(f_j,\ang{\phi_j}) = (g_j,\ang{\psi_j})$ for $j \ne i$, and
      \item $(f_i,\ang{\phi_i})' = (g_i,\ang{\psi_i})'$ as in \cref{remove-bpt-ops} above.
      \end{itemize}
    \end{enumerate}
    Thus, two $\obsim$-composable tuples of morphisms are equivalent if they differ by either
    \begin{itemize}
    \item composition in $\Fr\M$ of adjacent entries, or
    \item insertion or deletion of basepoint operations.
    \end{itemize}
    The equivalence class of $\uf$ is denoted \label{not:ufclass}$[\uf]$.
    \dqed
  \end{description} 
\end{definition}

\begin{definition}[Pointed Free Permutative Category]\label{def:Fst-object}
  Suppose $(\M,i^\M)$ is a pointed multicategory.
  Define the data of the \index{category!free permutative!pointed}\index{permutative category!free!pointed}\index{free!permutative category!pointed}\index{pointed!free permutative category}\emph{pointed free permutative category} $\Fst\M$ as follows.
  \begin{description}
  \item[Objects] The objects are given by $\obsim$-equivalence classes:
    \[
      \Ob(\Fst\M) = \Ob(\Fr\M)/\obsim.
    \]
  \item[Morphisms] The morphisms are given by $\simeq$-equivalence classes of $\obsim$-composable tuples:
    \[
      \Mor(\Fst\M) = \wt{\Mor}(\Fr\M)/\simeq.
    \]
    For $\uf = \bigl( (f_1,\ang{\phi_1}), \ldots, (f_m,\ang{\phi_m}) \bigr)$ in $\wt{\Mor}(\Fr\M)$, define
    \[
      \dom([\uf]) = [\dom(f_1,\ang{\phi_1})]
      \andspace
      \codom([\uf]) = [\codom(f_m,\ang{\phi_m})].
    \]
    Note that these are well defined since both relations $\isim$ and $\iisim$ preserve $\obsim$-equivalence classes of (co)domain.

  \item[Identities] For an $\obsim$-equivalence class $[\ang{x}]$ in $\Fst\M$, define
    \[
      1_{[\ang{x}]} = [1_{\ang{x}}].
    \]
    This is well defined by relation $\iisim$ for $1_{\ang{x}} = (1,\ang{1_{x_i}}_i)$.

  \item[Composition] Composition of equivalence classes
    \[
      [\ang{x}] \fto{[\uf]} [\ang{y}] \fto{[\ug]} [\ang{z}] \in \Fst\M
    \]
    is given by concatenation of representative tuples, 
   \[
     [\ug] \circ [\uf] = [\uf,\ug] = \bigl[ \bigl( (f_1,\ang{\phi_1}), \ldots, (f_m,\ang{\phi_m}), (g_1,\ang{\psi_1}), \ldots, (g_n,\ang{\psi_n}) \bigr)\bigr],
   \]
   where
   \[
     \uf = \bigl( (f_1,\ang{\phi_1}), \ldots, (f_m,\ang{\phi_m}) \bigr) \andspace
     \ug = \bigl( (g_1,\ang{\psi_1}), \ldots, (g_n,\ang{\psi_n}) \bigr).
   \]
   This concatenation is $\obsim$-composable because
   \[
     \codom([\uf]) = [\codom(f_m,\ang{\phi_m})]
     \andspace
     \dom([\ug]) = [\dom(g_1,\ang{\phi_1})].
   \]
   The composite $[\ug]\circ[\uf]$ is well defined by the definitions of $\isim$ and $\iisim$.
 \end{description}
 This finishes the definition of $\Fst\M$.
\end{definition}

\begin{remark}
  The description of $\Fst$ in \cref{definition:up-to-basept-equiv,def:Fst-object} is an extension of the explicit description for coequalizers in $\Cat$ from \cite[Section~1.4]{yau-involutive}.
  \Cref{proposition:FstM-pushout} shows that $\Fst\M$ can equivalently be constructed as a pushout in $\permcatsu$.
\end{remark}

\begin{lemma}\label{FstM-category}
  In the context of \cref{def:Fst-object}, $\Fst\M$ is a category.
  If $\M$ is a small multicategory, then $\Fst\M$ is a small category.
\end{lemma}
\begin{proof}
  Associativity of composition holds because concatenation of sequences is associative.
  The composition with identities is unital by relation $\isim$.
  If $\M$ is small, then $\Fr\M$, and hence also $\Fst\M$, is small because its objects are finite tuples of objects of $\M$ (\cref{definition:free-perm}).
\end{proof}

\subsection*{Permutative Structure for $\Fst\M$}
Recall from \cref{definition:free-perm} that the unpointed free construction $\Fr\M$ has monoidal sum $\oplus$ given by concatenation and monoidal unit the empty tuple, $\ang{}$.
There is an induced permutative structure on $\Fst\M$, described as follows.
\begin{definition}\label{definition:uf1-1uf'}
  In the context of \cref{def:Fst-object}, suppose given morphisms
  \[
  [\uf]\cn[\ang{x}] \to [\ang{y}] \andspace [\uf']\cn [\ang{x'}] \to [\ang{y'}]
  \]
  in $\Fst\M$, with
  \begin{align*}
    \uf & = \bigl( (f_1,\ang{\phi_1}), \ldots, (f_m,\ang{\phi_m}) \bigr)
          \andspace \\
    \uf' & = \bigl( (f_1',\ang{\phi_1'}), \ldots, (f_{n}',\ang{\phi_{n}'}) \bigr).
  \end{align*}
  Define
  \begin{align}
    \uf \oplus 1_{\ang{x'}} & = \bigl( (f_1,\ang{\phi_1})\oplus 1_{\ang{x'}}, \ldots, (f_m,\ang{\phi_m})\oplus 1_{\ang{x'}} \bigr)\label{eq:uf1}
                              \andspace \\
    1_{\ang{y}} \oplus \uf' & = \bigl( 1_{\ang{y}} \oplus (f_1',\ang{\phi_1'}), \ldots, 1_{\ang{y}} \oplus (f_{n}',\ang{\phi_{n}'}) \bigr),\label{eq:1uf'}
  \end{align}
  where the sums at right are those of $\Fr\M$.
\end{definition}
\begin{explanation}\label{explanation:uf1-1uf'}
  Observe that the $\simeq$-equivalence classes of \cref{eq:uf1,eq:1uf'} are well defined in the following senses.
  \begin{itemize}
  \item If $\ang{x'} \obsim \ang{w'}$, then $\uf \oplus 1_{\ang{x'}} \iisim \uf \oplus 1_{\ang{w'}}$.
  \item If $\uf \isim \ug$, then $\uf \oplus 1_{\ang{x'}} \isim \ug \oplus 1_{\ang{x'}}$ by functoriality of $\oplus$.
  \item If $\uf \iisim \ug$, then $\uf \oplus 1_{\ang{x'}} \iisim \ug \oplus 1_{\ang{x'}}$.
  \end{itemize}
  The corresponding statements for \cref{eq:1uf'} hold likewise.
\end{explanation}

\begin{definition}[Permutative Structure on $\Fst\M$]\label{def:Fst-permutative}\index{category!free permutative!pointed}\index{permutative category!free!pointed}\index{free!permutative category!pointed}\index{pointed!free permutative category}
  Suppose $(\M,i^\M)$ is a pointed multicategory.
  Define the data of a permutative category
  \[
    \big(\Fst\M, \oplus, [\ang{}], \xi\big)
  \]
  as follows.
  \begin{description}
  \item[Monoidal Sum] The monoidal sum functor
    \[
      \Fst\M \times \Fst\M \fto{\oplus} \Fst\M
    \]
    is defined by that of $\Fr\M$ on representative objects and morphisms, as follows.
    \begin{itemize}
    \item For objects $[\ang{x}]$ and $[\ang{x'}]$ in $\Fst\M$, let
      \[
        [\ang{x}] \oplus [\ang{x'}] = [\ang{x} \oplus \ang{x'}],
      \]
      where the monoidal sum at right is that of $\Fr\M$.
      This is well defined because concatenation of representatives in $\Fr\M$ preserves $\obsim$-equivalence classes.
    \item For morphisms $[\uf]\cn[\ang{x}] \to [\ang{y}]$ and $[\uf']\cn [\ang{x'}] \to [\ang{y'}]$ in $\Fst\M$, use \cref{definition:uf1-1uf'} and let
      \begin{align*}
        [\uf] \oplus [\uf']
        & = [1_{\ang{y}} \oplus \uf'] \circ [\uf \oplus 1_{\ang{x'}}]\\
        & = [\uf \oplus 1_{\ang{y'}}] \circ [1_{\ang{x}} \oplus \uf'].
      \end{align*}
      The $\simeq$-equivalence classes at right are well defined by \cref{explanation:uf1-1uf'}.
      The two indicated composites are equal as morphisms in $\Fst\M$ by $\isim$-equivalence and functoriality of $\oplus$ in $\Fr\M$.
    \end{itemize}
  \item[Monoidal Unit] The monoidal unit is $[\ang{}]$, the equivalence class of the empty tuple of objects.
  \item[Symmetry] The symmetry isomorphism for $\Fst\M$ is given by the equivalence class of the symmetry for $\Fr\M$, as below for objects $[\ang{x}]$ and $[\ang{x'}]$.
    \[
      \begin{tikzpicture}[x=40mm,y=8mm]
        \draw[0cell] 
        (0,0) node (a) {[\ang{x}] \oplus [\ang{x'}]}
        (a)+(0,-1) node (b) {[\ang{x} \oplus \ang{x'}]}
        (a)+(1,0) node (a') {[\ang{x'}] \oplus [\ang{x}]}
        (a')+(0,-1) node (b') {[\ang{x'} \oplus \ang{x}]}
        ;
        \draw[1cell] 
        (a) edge node {\xi_{[\ang{x}],[\ang{x'}]}} (a')
        (b) edge node {[\xi_{\ang{x},\ang{x'}}]} (b')
        (a) edge[equal] node {} (b)
        (a') edge[equal] node {} (b')
        ;
      \end{tikzpicture}
    \]
    This is well defined because $\xi_{\ang{w},\ang{w'}} \iisim \xi_{\ang{x},\ang{x'}}$ for $\ang{x} \obsim \ang{w}$ and $\ang{x'} \obsim \ang{w'}$.
  \end{description}
  The monoidal sum is strictly associative and unital because concatenation of tuples is so.
  The symmetry and hexagon axioms \cref{symmoncatsymhexagon} follow from those of $\Fr\M$.
\end{definition}

\subsection*{The 2-Functor $\Fst$}
Recall the descriptions of $\Fr$ on multifunctors and multinatural transformations from \cref{definition:free-smfun,definition:free-perm-multinat}.
These induce pointed variants for $\Fst$, which we now describe.
\begin{definition}[$\Fst$ on 1-Cells]\label{def:Fst-onecells}\index{category!free permutative!pointed}\index{permutative category!free!pointed}\index{free!permutative category!pointed}\index{pointed!free permutative category}
  Suppose given a pointed multifunctor between pointed multicategories
  \[
    H \cn (\M,i^\M) \to (\N,i^\N).
  \]
  Define a strict symmetric monoidal functor
  \[
    \Fst H \cn \Fst\M \to \Fst\N
  \]
  induced by $\Fr H$ as follows.
  \begin{description}
  \item[Object Assignment] For an object $[\ang{x_{i}}_{i=1}^r]$ in $\Fst\M$, define the object
    \begin{equation}\label{eq:FstHx}
      (\Fst H)[\ang{x_i}_{i=1}^r] = [\ang{Hx_i}_{i=1}^r] \inspace \Fst\N.
    \end{equation}
    The assumption that $H$ is a pointed multifunctor ensures that \cref{eq:FstHx} is well defined with respect to $\obsim$-equivalence. 
  \item[Morphism Assignment] For a morphism
    \[
      [\uf] \cn [\ang{x}] \to [\ang{y}] \inspace \Fst\M
    \]
    with
    \[
      \uf = \bigl( (f_1,\ang{\phi_1}), \ldots, (f_m,\ang{\phi_m}) \bigr),
    \]
    define
    \begin{align}
      (\Fst H)[\uf]
      & = \big[ (\Fr H)(f_1,\ang{\phi_1}), \ldots, (\Fr H)(f_m,\ang{\phi_m}) \big],\label{eq:FstHfphi}\\
      & = \big[ (f_1,\ang{H\phi_{1,j_1}}_{j_1=1}^{s_1}), \ldots, (f_m,\ang{H\phi_{m,j_m}}_{j_m=1}^{s_m}) \big],\nonumber
    \end{align}
    where each $\ang{\phi_i}$ has length $s_i$.
    Multifunctoriality of $H$ ensures that \cref{eq:FstHfphi} is well defined with respect to $\isim$-equivalence.
    The assumption that $H$ is pointed, and therefore preserves basepoint operations, ensures that \cref{eq:FstHfphi} is well defined with respect to $\iisim$-equivalence.
  \item[Constraints] The unit and monoidal constraints for $\Fst H$ are identities.
  \end{description}
  Functoriality and the strict symmetric monoidal conditions for $\Fst H$ follow from those of $\Fr H$ on representatives.
  This finishes the definition of $\Fst H$.
\end{definition}

\begin{definition}[$\Fst$ on 2-Cells]\label{def:Fst-twocells}\index{category!free permutative!pointed}\index{permutative category!free!pointed}\index{free!permutative category!pointed}\index{pointed!free permutative category}
  Suppose given a pointed multinatural transformation $\theta$ between pointed multifunctors between pointed multicategories
  \[\begin{tikzpicture}
      \def\h{2} \def\t{23}
      \draw[0cell]
      (0,0) node (a) {\phantom{X}}
      (a)+(\h,0) node (b) {\phantom{X}}
      (a)+(-.4,0) node (a') {(\M,i^\M)}
      (b)+(.4,0) node (b') {(\N,i^\N).}
      ;
      \draw[1cell=.9]
      (a) edge[bend left=\t,transform canvas={yshift=0ex}] node {H} (b)
      (a) edge[bend right=\t,transform canvas={yshift=0ex}] node[swap] {K} (b)
      ;
      \draw[2cell]
      node[between=a and b at .45, rotate=-90, 2label={above,\,\theta}] {\Rightarrow}
      ;
    \end{tikzpicture}\]
  Define a monoidal natural transformation
  \[\begin{tikzpicture}
      \def\h{2.5} \def\t{18}
      \draw[0cell]
      (0,0) node (a) {\Fst\M}
      (a)+(\h,0) node (b) {\Fst\N}
      ;
      \draw[1cell=.9]
      (a) edge[bend left=\t,transform canvas={yshift=0ex}] node {\Fst H} (b)
      (a) edge[bend right=\t,transform canvas={yshift=0ex}] node[swap] {\Fst K} (b)
      ;
      \draw[2cell]
      node[between=a and b at .42, rotate=-90, 2label={above,\,\Fst\theta}] {\Rightarrow}
      ;
    \end{tikzpicture}\]
  induced by $\Fr \theta$ with component morphism in $\Fst\N$
  \begin{equation}\label{eq:Fst-theta-x}
    (\Fst\theta)_{[\ang{x}]}
    = [(\Fr\theta)_{\ang{x}}]
    = \left[\left(1_{\ufs{r}} \scs \ang{\theta_{x_i}}_{i=1}^r \right) \right]
    \cn [\ang{Hx_i}_{i=1}^r] \to [\ang{Kx_i}_{i=1}^r]
  \end{equation}
  for each object $[\ang{x_i}_{i=1}^r]$ in $\Fst\M$.
  The assumption that $\theta$ is a pointed multifunctor, and hence its basepoint component is an identity operation, ensures that \cref{eq:Fst-theta-x} is well defined with respect to $\obsim$-equivalence.
  The monoidal naturality conditions for $\Fst\theta$ follow from those of $\Fr\theta$ on representatives.
  This finishes the definition of $\Fst\theta$.
\end{definition}

The 2-functoriality of $\Fr$ (\cref{proposition:free-perm-functor}) implies that the assignments above determine a 2-functor.
\begin{theorem}\label{ptmulticat-thm-i}\index{category!free permutative - 2-functor!pointed}\index{permutative category!free!2-functor!pointed}\index{free!permutative category 2-functor!pointed}
  The assignments on objects, 1-cells, and 2-cells given in \cref{def:Fst-permutative,def:Fst-onecells,def:Fst-twocells}, respectively,
  determine a 2-functor
  \[
    \Fst \cn \pMulticat \to \permcatst.
  \]
\end{theorem}

\section{Relating Unpointed and Pointed Free Permutative Categories}
\label{sec:relating-FFst}

In this section we define a 2-natural transformation
\[
  \begin{tikzpicture}
    \def\h{2} \def\t{23}
    \draw[0cell]
    (0,0) node (a) {\phantom{X}}
    (a)+(\h,0) node (b) {\phantom{X}}
    (a)+(-.55,0) node (a') {\pMulticat}
    (b)+(.6,.06) node (b') {\permcatst.}
    ;
    \draw[1cell=.9]
    (a) edge[bend left=\t,transform canvas={yshift=0ex}] node {\Fr \circ \Ust} (b)
    (a) edge[bend right=\t,transform canvas={yshift=0ex}] node[swap] {\Fst} (b)
    ;
    \draw[2cell]
    node[between=a and b at .45, rotate=-90, 2label={above,\,\pst}] {\Rightarrow}
    ;
  \end{tikzpicture}
\]
We usually suppress $\Ust$ and abbreviate $\Fr \circ \Ust$ to $\Fr$.
Using $\pst$, we show in \cref{proposition:FstM-pushout} that $\Fst\M$ is the pushout in $\permcatst$ of the span
\[
  \boldone \leftarrow \Fr\Mterm \fto{\Fr i^\M} \Fr\M,
\]
where $\boldone$ is the terminal permutative category, $\Mterm$ is the terminal multicategory (\cref{definition:terminal-operad-comm}), and $(\M,i^\M)$ is a small pointed multicategory.

\begin{definition}\label{def:FFst-projection}
  For a small pointed multicategory $(\M,i^\M)$, define a strict symmetric monoidal functor
  \[
    \pst_{\M} \cn \Fr\M \to \Fst\M
  \]
  on objects $\ang{x}$ and morphisms $(f,\ang{\phi})$ by the assignments
  \begin{equation}\label{eq:pst}
    \left\{
      \begin{aligned}
        \ang{x} & \mapsto [\ang{x}] \andspace \\
        \ (f,\ang{\phi}) & \mapsto [(f , \ang{\phi})],
      \end{aligned}
    \right.
  \end{equation}
  where a morphism of $\Fr\M$ is regarded as an $\obsim$-composable tuple of length one.
  These assignments are functorial by definition of $\isim$.
  Recall that the monoidal sum in $\Fst\M$ is given by concatenation, and the monoidal sum in $\Fr\M$ is functorial.
  These imply that $\pst_\M$ is a strict symmetric monoidal functor.
\end{definition}

\begin{proposition}\label{pst-iinatural}
  The strict symmetric monoidal functors $\pst_\M$ of \cref{def:FFst-projection} are components of a 2-natural transformation
  \begin{equation}\label{eq:pst-iinatural}
    \pst \cn \Fr \to \Fst.
  \end{equation}
\end{proposition}
\begin{proof}
  Suppose given a pointed multinatural transformation $\theta$ between pointed multifunctors between small pointed multicategories
  \[\begin{tikzpicture}
      \def\h{2} \def\t{23}
      \draw[0cell]
      (0,0) node (a) {\phantom{X}}
      (a)+(\h,0) node (b) {\phantom{X}}
      (a)+(-.4,0) node (a') {(\M,i^\M)}
      (b)+(.4,0) node (b') {(\N,i^\N).}
      ;
      \draw[1cell=.9]
      (a) edge[bend left=\t,transform canvas={yshift=0ex}] node {H} (b)
      (a) edge[bend right=\t,transform canvas={yshift=0ex}] node[swap] {K} (b)
      ;
      \draw[2cell]
      node[between=a and b at .45, rotate=-90, 2label={above,\,\theta}] {\Rightarrow}
      ;
    \end{tikzpicture}\]
  In the following diagram we have
  \[
    \pst_\N \circ (\Fr H) = (\Fst H) \circ \pst_\M
    \andspace
    \pst_\N \circ (\Fr K) = (\Fst K) \circ \pst_\M
  \]
  by definition of $\pst$ and \cref{def:Fst-onecells} for $\Fst$ on 1-cells.
  \[
    \begin{tikzpicture}
      \def\h{3} \def\t{23}
      \def\v{2}
      \draw[0cell=.9] 
      (0,0) node (a) {\Fr\M}
      (a)+(\h,0) node (b) {\Fr\N}
      (a)+(0,-\v) node (c) {\Fst\M}
      (c)+(\h,0) node (d) {\Fst\N}
      ;
      \draw[1cell=.9] 
      (a) edge[bend left=\t] node {\Fr H} (b)
      (c) edge[bend left=\t] node {\Fst H} (d)
      (a) edge[',bend right=\t] node {\Fr K} (b)
      (c) edge[',bend right=\t] node {\Fst K} (d)
      (a) edge['] node {\pst_\M} (c)
      (b) edge node {\pst_\N} (d)
      ;
      \draw[2cell]
      node[between=a and b at .45, rotate=-90, 2label={above,\,\Fr\theta}] {\Rightarrow}
      node[between=c and d at .45, rotate=-90, 2label={above,\,\Fst\theta}] {\Rightarrow}
      ;
    \end{tikzpicture}
  \]
  By \cref{def:Fst-twocells}, we also have the equality of whiskerings
  \[
    (\Fst \theta) * \pst_\M = \pst_\N * (\Fr \theta).
  \]
  This completes the proof that $\pst$ is 2-natural.
\end{proof}

Recall from \cref{example:free-Mterm} that the free permutative category \index{multicategory!terminal}\index{terminal!multicategory}$\Fr\Mterm$ is the natural number category $\mathbf{N}$ whose objects are natural numbers and morphisms are given by morphisms of finite sets.
\begin{proposition}\label{proposition:FstM-pushout}
  For each small pointed multicategory $(\M,i^\M)$, the diagram
  \begin{equation}\label{eq:FstM-pushout}
    \begin{tikzpicture}[x=23mm,y=13mm,vcenter]
      \draw[0cell=.9] 
      (0,0) node (a) {\Fr\Mterm}
      (a)+(1,0) node (b) {\Fr\M}
      (a)+(0,-1) node (c) {\boldone}
      (c)+(1,0) node (d) {\Fst\M}
      ;
      \draw[1cell=.9] 
      (a) edge node {\Fr i^\M} (b)
      (a) edge['] node {} (c)
      (c) edge node {} (d)
      (b) edge node {\pst_\M} (d)
      ;
    \end{tikzpicture}
  \end{equation}
  is a 2-pushout in $\permcatst$.
\end{proposition}
\begin{proof}
  Suppose given a permutative category $\C$ and a strict symmetric monoidal functor $Q$ such that the outer diagram below commutes.
  \begin{equation}\label{eq:FstM-pushout-Q}
    \begin{tikzpicture}[x=23mm,y=13mm,vcenter]
      \draw[0cell=.9] 
      (0,0) node (a) {\Fr\Mterm}
      (a)+(1,0) node (b) {\Fr\M}
      (a)+(0,-1) node (c) {\boldone}
      (c)+(1,0) node (d) {\Fst\M}
      (c)+(1.5,-.5) node (q) {\C}
      ;
      \draw[1cell=.9] 
      (a) edge node {\Fr i^\M} (b)
      (a) edge['] node {} (c)
      (c) edge node {} (d)
      (c) edge[bend right=15] node {} (q)
      (b) edge['] node {\pst_\M} (d)
      (b) edge[bend left=25] node {Q} (q)
      (d) edge[dashed] node[pos=.3] {\ol{Q}} (q)
      ;
    \end{tikzpicture}
  \end{equation}
  Commutativity of the outer diagram implies that $Q$ sends the objects and morphisms of $\Fr\Mterm$ to, respectively, the monoidal unit $\pu \in \C$ and identities on that unit.

  If $\ang{x} \obsim \ang{y}$ in $\Fr\M$ then, because $Q$ is strict symmetric monoidal, we have
  \[
    Q\ang{x} = \oplus_i Qx_i = \oplus_j Qy_j = Q\ang{y},
  \]
  where the middle equality holds because $\ang{x} \obsim \ang{y}$ and $\pu \in \C$ is a strict unit.
  Therefore, $Q$ uniquely determines a well-defined assignment on objects $\ol{Q}$ that commutes with $\pst_\M$.
  
  Now consider
  \[
    [\uf] = \big[ (f_1,\ang{\phi_1}), \ldots, (f_m,\ang{\phi_m}) \big]
    \inspace \Mor(\Fst\M) = \wt{\Mor}(\Fr\M)/\simeq.
  \]
  For $\ol{Q}$ to be functorial and commute with $\pst_\M$, we must have
  \[
    \ol{Q} [\uf] = Q(f_m,\ang{\phi_m}) \circ \cdots \circ Q(f_1,\ang{\phi_1}).
  \]
  This is well defined with respect to $\isim$-equivalence classes by functoriality of $Q$.
  The definition of $\ol{Q}$ on $\iisim$-equivalence classes is well defined because the morphisms in $\Fr\M$ that are induced by basepoint operations will be sent by $Q$ to identity morphisms $1_\pu$ in $\C$.  

  Thus, there is a unique strict symmetric monoidal functor $\ol{Q}$ in \cref{eq:FstM-pushout-Q} that commutes with $\pst_\M$.
  This proves that $\Fst\M$ satisfies the desired 1-dimensional universal property.
  
  For the 2-dimensional aspect of the universal property, suppose given a monoidal natural transformation $\ka$ as in the diagram below, such that the whiskering $\ka * (\Fr i^\M)$ is equal to the identity 2-cell of the constant functor given by the left-bottom composite.
  \begin{equation}\label{eq:FstM-pushout-2}
    \begin{tikzpicture}[x=23mm,y=13mm,vcenter]
      \draw[0cell=.9] 
      (0,0) node (a) {\Fr\Mterm}
      (a)+(1,0) node (b) {\Fr\M}
      (a)+(0,-1) node (c) {\boldone}
      (c)+(1,0) node (d) {\Fst\M}
      (c)+(2,-.5) node (q) {\C}
      ;
      \draw[1cell=.9] 
      (a) edge node {\Fr i^\M} (b)
      (a) edge['] node {} (c)
      (c) edge node {} (d)
      (c) edge[bend right=25,in=220] node {} (q)
      (b) edge['] node {\pst_\M} (d)
      (b) edge[bend left=25] node[pos=.3] (Q) {Q} (q)
      (b) edge[',bend right=5] node[pos=.3] (R) {R} (q)
      (d) edge[dashed,bend left=5] node[pos=.3,scale=.7] (Qb) {\ol{Q}} (q)
      (d) edge[',dashed,bend right=30] node[pos=.3,scale=.7] (Rb) {\ol{R}} (q)
      ;
      \draw[2cell]
      node[between=Q and R at .5, rotate=-135, 2label={above,\ka}] {\Rightarrow}
      node[between=Qb and Rb at .5, rotate=-125, 2label={above,\ol{\ka}}, scale=.8] {\Rightarrow}
      ;
    \end{tikzpicture}
  \end{equation}
  Hence, the component of $\ka$ at the basepoint object is an identity morphism.
  Because $\ka$ is monoidal natural, this implies that
  \[
    \ka_{\ang{x}} = \ka_{\ang{y}}
  \]
  whenever $\ang{x} \obsim \ang{y}$.
  Therefore, $\ol{\ka}$ is uniquely determined and well defined with components
  \[
    \ol{\ka}_{[\ang{x}]} = \ka_{\ang{x}}.
  \]
  This completes the proof.
\end{proof}

\begin{definition}\label{definition:bpt-collapsing}
  In the context of \cref{proposition:FstM-pushout}, say that a strict symmetric monoidal functor
  \[
    Q\cn \Fr\M \to \C
  \]
  is \emph{basepoint-collapsing} if the outer diagram \cref{eq:FstM-pushout-Q} commutes.
  Say that a monoidal natural transformation
  \[
    \ka\cn Q \to R
  \]
  is \emph{basepoint-collapsing} if $Q$ and $R$ are basepoint-collapsing and 
  \[
  \ka * \Fr i^\M = 1_{1_\pu}
  \]
  as in \cref{eq:FstM-pushout-2}.
  Thus, $Q$ is basepoint-collapsing if and only if $Q\iota_n = 1_\pu$ for each basepoint operation $\iota_n$.
  Similarly, $\ka$ is basepoint-collapsing if and only if $\ka_{*} = 1_\pu$.
  With these terms, \cref{proposition:FstM-pushout} asserts that basepoint-collapsing 1- and 2-cells of $\permcatst$ extend uniquely along $\pst_\M$.
\end{definition}

\section{Pointed Free Permutative Category as a Left 2-Adjoint}
\label{sec:FEst-adj}

In this section we show that the adjunction $\Fr \dashv \End$ from \cref{theorem:FE-adj} induces an adjunction $\Fst \dashv \Endst$.
As above, we usually suppress the forgetful $\Ust$.

\begin{definition}[Unit]\label{def:etast}
  Suppose $\M$ is a small pointed multicategory.
  We define the composite multifunctor
  \begin{equation}\label{etast-def}
    \begin{tikzpicture}[vcenter]
      \def\u{.6}
      \draw[0cell]
      (0,0) node (a) {\M}
      (a)+(3,0) node (b) {\End\,\Fr\M = \Endst\, \Fr\M}
      (b)+(4.5,0) node (c) {\Endst\, \Fst\M}
      ;
      \draw[1cell]
      (a) edge node {\eta_\M} (b)
      (b) edge node {\Endst(\pst_\M)} (c)
      ;
      \draw[1cell]
      (a) [rounded corners=3pt] |- ($(b)+(-1,\u)$)
      -- node {\etast_\M} ($(b)+(1,\u)$) -| (c)
      ;
    \end{tikzpicture}
  \end{equation}
  as follows.
  \begin{itemize}
  \item $\eta_\M$ is the unit of the 2-adjunction $\Fr \dashv \End$ in \cref{definition:eta}.
  \item $\Endst$ is the $\Cat$-multifunctor in \cref{expl:endst-catmulti}.
  \item $\pst_\M \cn \Fr \M \to \Fst\M$ is the strict symmetric monoidal functor in \cref{def:FFst-projection}.
  \end{itemize}
  This finishes the definition of $\etast_\M$.
\end{definition}

\begin{explanation}\label{expl:etast-explicit}
  The multifunctor $\etast_\M$ in \cref{etast-def} is given explicitly by the assignments
  \begin{align*}
    y & \mapsto [(y)], \forspace y \in \M, \andspace\\
    \phi & \mapsto \big[ \big( \iota_r,(\phi) \big) \big] \forspace \phi \in \M\scmap{\ang{x_i}_{i=1}^r;y},
  \end{align*}
  where $\iota_r \cn \ufs{r} \to \ufs{1}$ is the unique index map.
  Note, in particular, that $\etast_\M$ sends each basepoint operation of $\M$ to the equivalence class of the identity $1_{[\ang{}]}$ by relation $\iisim$.
\end{explanation}

\begin{explanation}\label{expl:eta-not-ptd}
  In \cref{etast-def}, even though $\M$ and $\Endst\, \Fr\M$ are pointed multicategories, we emphasize that $\eta_\M$ is \emph{not} pointed because
  $\eta_\M$ sends the basepoint object $* \in \M$ to the length-one tuple $(*) \ne \ang{}$.
  The composite, $\etast_\M$ in \cref{etast-def} \emph{is} pointed because $[(*)] = [\ang{}]$ in $\Fst\M$.
\end{explanation}

\begin{lemma}\label{etast-pointed}
  In the context of \cref{def:etast},
  \[
    \etast_\M \cn \M \to \Endst\, \Fst\M
  \]
  is a pointed multifunctor that is, moreover, 2-natural in $\M$.
\end{lemma}
\begin{proof}
  Multifunctoriality of $\etast_\M$ follows from that of $\eta_\M$.
  As noted in \cref{expl:etast-explicit,expl:eta-not-ptd}, $\etast_\M$ is pointed by relations $\obsim$ on objects and $\iisim$ on morphisms in $\Fst\M$.
  The 2-naturality of $\etast$ follows from that of $\eta_\M$ and $\pst_\M$, together with $\Cat$-multifunctoriality of $\Endst$.
\end{proof}

For a small permutative category $\C$, recall the counit
\[
  \epz_\C \cn \Fr\, \End \C \to \C
\]
from \cref{definition:epz}.
Note that $\epz_\C$ is basepoint-collapsing (\cref{definition:bpt-collapsing}) by its definition on morphisms \cref{Fcounit-mor}.
\begin{definition}[Counit]\label{def:epzst}
  Suppose $(\C,\oplus,e)$ is a small permutative category.
  We define a strict symmetric monoidal functor $\epzst_\C$ as the unique dashed arrow
  \begin{equation}\label{epzst-def}
    \begin{tikzpicture}[vcenter]
      \def\u{.6}
      \draw[0cell]
      (0,0) node (a) {\Fr\, \End\C = \Fr\, \Endst\C}
      (a)+(4,0) node (b) {\Fst\, \Endst\C}
      (b)+(3.5,0) node (c) {\C}
      ;
      \draw[1cell]
      (a) edge node {\pst_{\Endst \C}} (b)
      (b) edge[dashed] node {\epzst_\C} (c)
      ;
      \draw[1cell]
      (a.160) [rounded corners=3pt] |- ($(b)+(-1,\u)$)
      -- node {\epz_\C} ($(b)+(1,\u)$) -| (c)
      ;
    \end{tikzpicture}
  \end{equation}
  induced by the pushout \cref{eq:FstM-pushout-Q}, where
  \begin{itemize}
  \item $\M = \Endst \C$,
  \item $\pst_{\Endst\C}$ is the strict symmetric monoidal functor in \cref{def:FFst-projection}, and
  \item $Q = \epz_\C$ is the counit of the 2-adjunction $\Fr \dashv \End$ in \cref{definition:epz}.
  \end{itemize}
  This finishes the definition of $\epzst_\C$.
\end{definition}

\begin{explanation}\label{expl:epzst-explicit}
  The proof of \cref{proposition:FstM-pushout} explains that the multifunctor $\epzst_\C$ in \cref{epzst-def} is given by $\epz_\C$ on representative objects and operations of $\Fr\Endst\C$.
\end{explanation}

\begin{lemma}\label{epzst-strict}
  In the context of \cref{def:epzst}, $\epzst_\C$ is 2-natural in $\C$.
\end{lemma}
\begin{proof}
  The 2-naturality of $\epzst$ follows from the universality of the pushouts in \cref{proposition:FstM-pushout}.
  Indeed, given a strict symmetric monoidal functor $R\cn \C \to \D$,
  consider the following.
  \begin{equation}\label{eq:epzst-nat}
    \begin{tikzpicture}[vcenter]
      \def\u{.6}
      \draw[0cell=.9]
      (0,0) node (a) {\Fr\, \Endst\C}
      (a)+(3.5,0) node (b) {\Fst\, \Endst\C}
      (b)+(3.5,0) node (c) {\C}
      (0,-2) node (a') {\Fr\, \Endst\D}
      (a')+(3.5,0) node (b') {\Fst\, \Endst\D}
      (b')+(3.5,0) node (c') {\D}
      ;
      \draw[1cell=.9]
      (a) edge node {\pst_{\Endst \C}} (b)
      (b) edge node {\epzst_\C} (c)
      (a') edge node {\pst_{\Endst \D}} (b')
      (b') edge node {\epzst_\D} (c')
      (a) edge node {} (a')
      (b) edge node {} (b')
      (c) edge node {} (c')
      (b) edge[dashed] node {} (c')
      ;
      \draw[1cell=.9]
      (a) [rounded corners=3pt] |- ($(b)+(-1,\u)$)
      -- node {\epz_\C} ($(b)+(1,\u)$) -| (c)
      ;
      \draw[1cell]
      (a') [rounded corners=3pt] |- ($(b')+(-1,-\u)$)
      -- node {\epz_\D} ($(b')+(1,-\u)$) -| (c')
      ;
    \end{tikzpicture}
  \end{equation}
  In the above diagram, the rectangle at left and the outer rectangle commute by naturality of $\pst$ and $\epz$, respectively.
  This implies that both composites around the rectangle at right are basepoint-collapsing.
  Therefore, the two composites are equal by uniqueness of the universal dashed functor in \cref{eq:epzst-nat}.
  Similar reasoning, using the 2-dimensional aspect of the pushout, verifies 2-naturality of $\epzst$ with respect to monoidal natural transformations.
\end{proof}

The diagram \cref{FEFEstdplusUst} below uses the 2-adjunctions 
\[
  \dplus \dashv \Ust \andspace \Fr \dashv \End
\]
in \cref{dplustwofunctor,theorem:FE-adj}, respectively.

\begin{theorem}\label{ptmulticat-thm-v}
  There is a 2-adjunction
  \[\begin{tikzpicture}
      \def\t{17}
      \draw[0cell]
      (0,0) node (a) {\phantom{X}}
      (a)+(-.55,0) node (a') {\pMulticat}
      (a)+(1,0) node (x) {\bot}
      (a)+(2,0) node (b) {\phantom{X}}
      (b)+(.65,.06) node (b') {\permcatst}
      ;
      \draw[1cell=.9]
      (a) edge[bend left=\t] node {\Fst} (b)
      (b) edge[bend left=\t] node {\Endst} (a)
      ;
    \end{tikzpicture}\]
  with unit $\etast$ and counit $\epzst$ in \cref{def:etast,def:epzst}, respectively.
  Moreover, \cref{FEst-adj-i,FEst-adj-ii} below hold regarding the following diagram.
  \begin{equation}\label{FEFEstdplusUst}
    \begin{tikzpicture}[vcenter]
      \def\s{12} \def\t{18} \def\h{2.5}
      \draw[0cell=.9]
      (0,0) node (a) {\pMulticat}
      (a)+(.6,0) node (x) {}
      (x)+(\h,0) node (y) {}
      (y)+(.7,.06) node (b) {\permcatst}
      (x)+(\h/2,1.75) node (c) {\Multicat}
      ;
      \draw[1cell=.85]
      (x) edge[bend left=\s,transform canvas={yshift=.5ex}] node {\Fst} (y)
      (y) edge[transform canvas={yshift=-.5ex}] node {\Endst} (x)
      (c) edge[bend left=\t,transform canvas={xshift=1ex}] node[pos=.6] {\Fr} (b)
      (b) edge[transform canvas={xshift=.5ex}] node[pos=.6] {\End} (c)
      (c) edge[bend right=\t, transform canvas={xshift=-1.6ex}] node[swap,pos=.6] {\dplus} (a)
      (a) edge[transform canvas={xshift=-1ex}] node[swap,pos=.6] {\Ust} (c)
      ;
    \end{tikzpicture}
  \end{equation}
  \begin{romenumerate}
  \item\label{FEst-adj-i}
    There is an equality of 2-functors
    \begin{equation}\label{EUstEst}
      \End = \Ust \circ \Endst \cn \permcatst \to \Multicat
    \end{equation}
    given by restricting the diagram \cref{endufactor} to $\permcatst$.
  \item\label{FEst-adj-ii}
    There is a 2-natural isomorphism
    \begin{equation}\label{FFstdplus}
      \Fr \iso \Fst \circ \dplus \cn \Multicat \to \permcatst.
    \end{equation}
  \end{romenumerate}
\end{theorem}
\begin{proof}
  Throughout this proof we write $\Er = \End$ and $\Est = \Endst$.
  The unit $\etast$ and counit $\epzst$ are shown to be 2-natural in \cref{etast-pointed,epzst-strict}, respectively.
  We now verify that the two triangle identities for $\etast$ and $\epzst$ follow from those of $\eta$ and $\epz$.

  Consider the following diagram, for a small pointed multicategory $\M$.
  As above, we suppress the forgetful $\Ust$.
  \begin{equation}\label{eq:ltriang-etast-epzst}
    \begin{tikzpicture}[x=35mm,y=13mm,vcenter]
      \draw[0cell=.9] 
      (0,0) node (a) {\Fr\M}
      (a)+(1,0) node (b) {\Fr\Er\Fr\M}
      (b)+(1,0) node (c) {\Fr\M}
      (b)+(0,-7mm) node (b') {\Fr\Est\Fr\M}
      (b)+(-.3,-1.2) node (x) {\Fr\Est\Fst\M}
      (b)+(.3,-1.2) node (y) {\Fst\Est\Fr\M}
      (a)+(0,-2) node (d) {\Fst\M}
      (d)+(1,0) node (e) {\Fst\Est\Fst\M}
      (e)+(1,0) node (f) {\Fst\M}
      ;
      \draw[1cell=.9] 
      (a) edge node {\Fr\eta_\M} (b)
      (b) edge node {\epz_{\Fr\M}} (c)
      (d) edge['] node {\Fst\etast_\M} (e)
      (e) edge['] node {\epzst_{\Fst\M}} (f)
      (a) edge['] node {\pst_\M} (d)
      (c) edge node {\pst_\M} (f)
      (b) edge[equal] node {} (b')
      (b') edge['] node {\Fr\Est\pst_\M} (x)
      (b') edge node {\pst_{\Est\Fr\M}} (y)
      (x) edge['] node[pos=.3] {\pst_{\Est\Fst\M}} (e)
      (y) edge node[pos=.3] {\Fst\Est\pst_{\M}} (e)
      (a) edge['] node {\Fr\etast_\M} (x)
      (y) edge['] node {\epzst_{\Fr\M}} (c)
      ;
    \end{tikzpicture}
  \end{equation}
  In the above diagram, the two upper regions commute by definition of $\etast$ and $\epzst$.
  The lower left quadrilateral commutes by 2-naturality of $\pst$ with respect to $\etast_\M$.
  The middle quadrilateral commutes by 2-naturality of $\pst$ with respect to $\Est\pst_\M$.
  The lower right quadrilateral commutes by 2-naturality of $\epzst$ with respect to $\pst_\M$.

  The composite along the top of \cref{eq:ltriang-etast-epzst} is equal to $1_{\Fr\M}$ by the triangle identity for $\eta$ and $\epz$.
  Thus, the composite along the top and right, which is equal to $\pst_\M$, is basepoint-collapsing in the sense of \cref{definition:bpt-collapsing}.
  Therefore, by the universal property of the pushout (\cref{proposition:FstM-pushout}), the bottom composite of \cref{eq:ltriang-etast-epzst} is equal to $1_{\Fst\M}$.
  
  The other triangle identity is simpler.
  For each small permutative category $\C$, the following diagram commutes by definition of $\etast$ and $\epzst$.
  \[
    \begin{tikzpicture}[x=35mm,y=13mm]
      \def\veq{6.5mm}
      \draw[0cell=.9] 
      (0,0) node (a) {\Er\C}
      (a)+(1,.5) node (b) {\Er\Fr\Er\C}
      (b)+(1,-.5) node (c) {\Er\C}
      (a)+(0,-\veq) node (d) {\Est\C}
      (d)+(1,-.5) node (e) {\Est\Fst\Est\C}
      (c)+(0,-\veq) node (f) {\Est\C}
      (b)+(0,-\veq) node (b') {\Est\Fr\Est\C}
      ;
      \draw[1cell=.9] 
      (a) edge node[pos=.6] {\eta_{\Er\C}} (b)
      (b) edge node[pos=.4] {\Er\epz_\C} (c)
      (d) edge node[pos=.4] {\etast_{\Est\C}} (e)
      (e) edge node[pos=.65] {\Est\epzst_\C} (f)
      (a) edge[equal] node {} (d)
      (c) edge[equal] node {} (f)
      (b) edge[equal] node {} (b')
      (b') edge node {\Est\pst_{\Est\C}} (e)
      ;
    \end{tikzpicture}
  \]
  Commutativity of the above diagram, together with the respective triangle identity for $\eta$ and $\epz$, implies that the bottom composite is equal to $1_{\Est\C}$.

  The equality \cref{FEst-adj-i} holds by definition.
  Uniqueness of left adjoints implies the existence of an isomorphism as in \cref{FEst-adj-ii}.
\end{proof}

\begin{explanation}\label{explanation:FFst+iso}
  An explicit description of the isomorphism
  \[
    \Fr \iso \Fst \circ \dplus \cn \Multicat \to \permcatst
  \]
  in \cref{ptmulticat-thm-v}~\cref{FEst-adj-ii} can be given as follows.
  For a small multicategory $\M$, the canonical inclusion into $\M_+$ induces a strict symmetric monoidal functor
  \[
    \Fr(\M) \to \Fr(\M_+) \fto{\pst_\M} \Fst(\M_+).
  \]
  For the reverse direction, there is a strict symmetric monoidal functor
  \begin{equation}\label{eq:FFst+reverse}
    \Fr(\M_+) \to \Fr\M
  \end{equation}
  that is given on length-one tuples by
  \[
    (x) \mapsto 
    \begin{cases}
      (x) & \ifspace x \ne *\\
      \ang{} & \ifspace x = *.
    \end{cases}
  \]
  On generating morphisms $(\iota_r,(\phi))$, where $\phi$ is an $r$-ary operation and $\iota_r\cn \ufs{r} \to \ufs{1}$ is the unique index map, \cref{eq:FFst+reverse} is given by
  \[
    (\iota_r,(\phi)) \mapsto
    \begin{cases}
      (\iota_r,(\phi)) & \ifspace \phi \in \M\scmap{\ang{x};y}\\
      1_{\ang{}} & \text{\quad otherwise,}
    \end{cases}
  \]
  where the second case applies if $\phi$ is a basepoint operation of $\M_+$.
  These assignments uniquely determine a strict symmetric monoidal functor that is basepoint-collapsing, and hence induce a unique
  \begin{equation}\label{eq:FstM+-to-FM}
    \Fst(\M_+) \to \Fr\M.
  \end{equation}
  The composite
  \[
    \Fr\M \to \Fr(\M_+) \to \Fst(\M_+) \to \Fr\M
  \]
  is the identity by definition.
  The other composite is the identity by uniqueness of \cref{eq:FstM+-to-FM}.
  The 2-naturality of this isomorphism follows from 2-naturality of the canonical inclusion $\M \to \M_+$ and the 2-dimensional universality of the pushout $\Fst(\M_+)$.
\end{explanation}

\section{Free Permutative Categories of \texorpdfstring{$\Mone$}{M1}-Modules}
\label{sec:FEm-adj}

In this section we use the 2-adjunctions
\[
  (\Monesma) \dashv \Um \andspace \Fst \dashv \Endst
\]
from \cref{MonesmaUmadj,ptmulticat-thm-v}, respectively, to describe an induced adjunction for $\Mone$-modules.

\begin{theorem}\label{ptmulticat-thm-ii}
  There is a 2-adjunction
  \[
    \begin{tikzpicture}
      \draw[0cell]
      (0,0) node (a) {\MoneMod}
      (a)+(3.25,0) node (b) {\permcatst}
      (a)+(1.55,0) node (x) {\bot}
      ;
      \draw[1cell=.9]
      (a) edge[bend left=15,transform canvas={yshift=-.8ex}] node {\Fm} (b)
      (b) edge[bend left=15,transform canvas={yshift=.7ex}] node {\Endm} (a)
      ;
    \end{tikzpicture}
  \]
  such that \cref{FEm-adj-i,FEm-adj-ii,FEm-adj-iii} below hold regarding the following diagram.
  \begin{equation}\label{FEmFEstMonesmaUm}
    \begin{tikzpicture}[vcenter]
      \def\s{12} \def\t{18} \def\h{2.5}
      \draw[0cell=.9]
      (0,0) node (a) {\pMulticat}
      (a)+(.6,0) node (x) {}
      (x)+(\h,0) node (y) {}
      (y)+(.7,.06) node (b) {\permcatst}
      (x)+(\h/2,-1.75) node (c) {\MoneMod}
      ;
      \draw[1cell=.85]
      (x) edge[transform canvas={yshift=.8ex},bend left=\s] node {\Fst} (y)
      (y) edge[transform canvas={yshift=0ex}] node {\Endst} (x)
      (c) edge[transform canvas={xshift=1ex},bend right=\t] node[swap] {\Fm} (b)
      (b) edge[transform canvas={xshift=.5ex}] node[swap,pos=.6] {\Endm} (c)
      (c) edge[transform canvas={xshift=-1ex}] node[swap,pos=.4] {\Um} (a)
      (a) edge[transform canvas={xshift=-1.5ex},bend right=\t] node[swap] {\Monesma} (c)
      ;
    \end{tikzpicture}
  \end{equation}
  \begin{romenumerate}
  \item\label{FEm-adj-i} 
    $\Fm$ is defined as the following composite 2-functor.
    \begin{equation}\label{Fm-def}
      \begin{tikzpicture}[vcenter]
        \def\u{.6}
        \draw[0cell]
        (0,0) node (a) {\MoneMod}
        (a)+(3,0) node (b) {\pMulticat}
        (b)+(3,0) node (c) {\permcatst}
        ;
        \draw[1cell]
        (a) edge node {\Um} (b)
        (b) edge node {\Fst} (c)
        ;
        \draw[1cell]
        (a) [rounded corners=3pt] |- ($(b)+(-1,\u)$)
        -- node {\Fm} ($(b)+(1,\u)$) -| (c)
        ;
      \end{tikzpicture}
    \end{equation}
  \item\label{FEm-adj-ii}
    There is an equality of 2-functors
    \begin{equation}\label{EstUmEm}
      \Endst = \Um \circ \Endm
    \end{equation}
    given by restricting \cref{endufactor} to $\permcatst$.
  \item\label{FEm-adj-iii}
    There is a 2-natural isomorphism
    \begin{equation}\label{FstFmMonesma}
      \Fst \iso \Fm \circ (\Monesma).
    \end{equation}
  \end{romenumerate}
\end{theorem}
\begin{proof}
  The 2-adjunction
  \[
    \Fm \dashv \Endm
  \]
  is given by the following 2-natural isomorphisms, explained below, for each $\N \in \MoneMod$ and $\Q \in \permcatst$:
  \begin{align*}
    \permcatst \big( \Fm \N \scs \Q \big)
    & = \permcatst \big( \Fst \Um \N \scs \Q \big)\\
    & \iso \pMulticat \big( \Um \N \scs \Endst \Q \big)\\
    & = \pMulticat \big( \Um \N \scs \Um\Endm \Q \big)\\
    & \iso \MoneMod \big( \N \scs \Endm \Q \big).
  \end{align*}
  In the above computation, the two equalities follow from \cref{Fm-def,EstUmEm}.
  The first isomorphism is given by the 2-adjunction $\Fst \dashv \Endst$ of \cref{ptmulticat-thm-v}.
  The last isomorphism holds because $\MoneMod$ is a full sub-2-category of $\pMulticat$ (\cref{proposition:EM2-5-1}~\cref{it:EM251-3}).

  The 2-natural isomorphism \cref{FstFmMonesma} follows from uniqueness of left adjoints.
  This completes the proof.
\end{proof}

\begin{explanation}\label{expl:etam-epzm}
  The unit and counit of the 2-adjunction $\Fm \dashv \Endm$ in \cref{ptmulticat-thm-ii} are given as follows.
  \begin{description}
  \item[Unit] The unit is the 2-natural transformation
    \begin{equation}\label{etam-define}
      \etam \cn 1_{\MoneMod} \to \Endm \Fm
    \end{equation}
    with component at each left $\Mone$-module $\M$ given by the morphism
    \[
      \begin{aligned}
        \etam_\M \cn \M \to
        & \Endm \Fm \M\\
        & = \Endm \Fst\, \Um \M \inspace \MoneMod.
      \end{aligned}
    \]
    This morphism is uniquely determined by the image
    \[
      \begin{aligned}
        \Um\etam_\M = \etast_{\Um\M} \cn \Um\M \to
        & \Um \Endm \Fst\, \Um \M\\
        & = \Endst\, \Fst\, \Um \M
      \end{aligned}
    \]
    in $\pMulticat$ using \cref{proposition:EM2-5-1}~\cref{it:EM251-3}, with $\etast$ the unit of $\Fst \dashv \Endst$ in \cref{etast-def}.
  \item[Counit] The counit is the 2-natural transformation
    \begin{equation}\label{epzm-define}
      \epzm \cn \Fm \Endm \to 1_{\permcatst}
    \end{equation}
    with components defined by $\epzst$ as follows:
    \begin{equation}\label{eq:epzst-epzm}
      \begin{tikzpicture}[x=30mm,y=20mm,vcenter]
        \draw[0cell] 
        (0,0) node (a) {
          \Fm\Endm\C
        }
        (a.east)+(7.5ex,0) node (b) {
          \Fst\Um\Em\C
        }
        (b.east)+(6ex,0) node (c) {
          \Fst\Endst\C
        }
        (c)+(1,0) node (d) {
          \C
        }
        ;
        \draw[1cell] 
        (c) edge node {\epzst_\C} (d)
        (a) edge[equal] node {} (b)
        (b) edge[equal] node {} (c)
        ;
        \draw[1cell,rounded corners=3pt]
        (a) |- (8ex,.35) -- node {\epzm_\C} ++(24ex,0) -| (d)
        ;
      \end{tikzpicture}
    \end{equation}
    for each $\C \in \permcatst$. \dqed
\end{description} 
\end{explanation}

\section{Examples of Pointed Free Permutative Categories} \label{sec:Fst-ex}

In this section we describe some examples of $\Fst\M$ for pointed multicategories $\M$.
Our first two examples make use of the isomorphism \cref{FFstdplus}:
\[
  \Fst(\M_+) \iso \Fr(\M)
\]
for unpointed multicategories $\M$.
Recall, from \cref{example:free-Mtu}, that \index{operad!initial}\index{initial operad}$\Fr(\Mtu)$ is isomorphic to the permutation category with objects given by natural numbers and morphisms given by symmetric groups.
\begin{example}[The Smash Unit]\label{example:Fst-Mtu}
  Recall from \cref{thm:pmulticat-smclosed} that $\Mtup$ is the \index{pointed!smash unit!pointed multicategory}\index{smash unit!pointed multicategory}\index{pointed!multicategory!smash unit}\index{multicategory!pointed!smash unit}unit for the smash product of pointed multicategories.
  Using \cref{FFstdplus} we have
  \begin{equation}
    \Fst(\Mtup) \iso \Fr(\Mtu).
  \end{equation}
  Each object of $\Fst(\Mtup)$ is an $\obsim$-equivalence class of tuples of objects of $\Mtup$ and, therefore, is identified with a natural number for its non-basepoint entries.
  Each morphism of $\Fr(\Mtup)$ is $\iisim$-equivalent to one with only non-basepoint operations.
  Because the only such operations are the unit operation in $\Mtu$, each morphism of $\Fr(\Mtup)$ is $\iisim$-equivalent to one given by a permutation of non-basepoint entries.
\end{example}

\begin{example}[The Terminal Multicategory]\label{example:Fst-Mterm}\ 
  Let $\varnothing$ denote the empty multicategory, with no objects and no operations.
  Then \index{multicategory!terminal}\index{terminal!multicategory}$\Mterm = \varnothing_+$ and we have
  \[
    \Fst(\Mterm) \iso \Fr(\varnothing) \iso \boldone,
  \]
  where the last isomorphism identifies the unique object of $\boldone$ with the empty tuple $\ang{} \in \Fr(\varnothing)$.
\end{example}

For further examples, we will make use of the following observation for general pointed multicategories.
\begin{explanation}[Morphisms Arising from Nullary Operations]\label{expl:w-angx}
  For each pointed multicategory $\M$, with basepoint object $*$, there is a nullary operation
  \[
    \iota_0 \in \M\scmap{\ang{};*}.
  \]
  For each object $x \in \M$, this gives a morphism in $\Fr\M$
  \[
    \brb{ f^{(x)}, {(1_x,\iota_0)} } \,\cn\,(x) \to (x,*),
  \]
  where the index map $f^{(x)}$ is the standard inclusion $\ufs{1} \hookrightarrow \ufs{2}$.

  More generally, suppose given an object $\ang{x} \in \Fr\M$.
  Recall that $\ang{x}^\wedge$ denotes the sub-tuple consisting of non-basepoint objects, $x_i \ne *$.
  The nullary operations give rise to a morphism in $\Fr\M$
  \begin{equation}\label{eq:w-angx}
    w^{\ang{x}} = \brb{f^{\ang{x}},\ang{\de^{x_i}}_i} \cn \ang{x}^\wedge \to \ang{x}
  \end{equation} 
  defined as follows.
  \begin{itemize}
  \item The index map $f^{\ang{x}}$ is the inclusion to the subset that indexes non-basepoint objects of $\ang{x}$, preserving their order.
  \item For each entry $x_i$ in $\ang{x}$,
    \[
      \de^{x_i} =
      \begin{cases}
        1_{x_i} & \ifspace x_i \ne *\\
        \iota_0 & \ifspace x_i = *.
      \end{cases}
    \]
  \end{itemize} 

  Note that we have
  \[
    w^{\ang{x}} \iisim 1_{\ang{x}}.
  \]
  Therefore, for any morphism
  \[
    (f,\ang{\phi})\cn \ang{x} \to \ang{x'} \inspace \Fr\M,
  \]
  we have
  \begin{equation}\label{eq:fphi-simeq-fphiw}
    (f,\ang{\phi})
    \isim \big( 1_{\ang{x}} \scs (f,\ang{\phi}) \big)
    \iisim \big( w^{\ang{x}} \scs (f,\ang{\phi}) \big)
    \inspace \wt{\Mor}(\Fr\M).
  \end{equation}
  We will make use of these equivalences now in the following examples, and later in the proof of \cref{ptmulticat-prop-viii}.
\end{explanation}

We next consider $\Fst\M = \Fm\M$ for $\Mone$-modules $\M$.
We use the following characterization, from \cite{cerberusIII}, to construct a partial inverse for $w^{\ang{x}}$.
Recall the interchange permutation $\xi^\otimes$ from \cref{eq:interchange-relation}.
\begin{proposition}[{\cite[10.2.7]{cerberusIII}}]\label{M1-mod-str}
  Suppose $\M$ is a pointed multicategory with basepoint $*$.
  Then a left $\Mone$-module structure on $\M$
  \[
    \mu\cn \Mone \sma \M \to \M.
  \]
  determines and is uniquely determined by operations
  \[
    \pi^2_1(x) \in \M\scmap{x,*;x} \forspace x \in \M
  \]
  such that the following basepoint, unit, and interchange conditions hold for objects $x$ in $\M$ and operations $\phi$ in $\M\scmap{\ang{y};x}$ with $\ang{y} = (y_1,\ldots,y_m)$.
  \begin{align}
    \pi^2_1(*) & = \iota^2 & \inspace & \M\scmap{*,*;*}\label{eq:pi2cond-1}\\
    \ga\scmap{\pi^2_1(x);1_x,\iota^0} & = 1_x & \inspace & \M\scmap{x;x}\label{eq:pi2cond-2}\\
    \ga\scmap{\phi;\ang{\pi^2_1(y_j)}_j} & = \ga\scmap{\pi^2_1(x);\phi,\iota^m} \cdot \xi^\otimes & \inspace & \M\scmap{y_1,*,\ldots,y_m,*;x}\label{eq:pi2cond-3}
  \end{align}
\end{proposition}
\begin{remark}
  Using equivariance on operations in $\M$ and induction, the operations $\pi^2_1$ determine a larger family with similar properties
  \[
    \pi^n_j(x) \in \M\scmap{*,\ldots,x,\ldots,*;x},
  \]
  where the input profile has $x$ in position $j$ and the $n-1$ other entries are basepoints.
  See \cite[8.4.5,10.2.4]{cerberusIII} for further discussion of these operations.
\end{remark}
\begin{example}[Operations Arising from Strict Units]
  If $\M = \cM a$ for a pointed finite set $a$, then the operations $\pi^2_1(x)$ arise because $* = \varnothing$ is a strict unit for the disjoint union in $\cP(a^\punc)$.
  If $\M = \Endst\C$ for a permutative category $\C$, then the operations $\pi^2_1(x)$ arise because $* = \pu$ is strict unit for the monoidal sum in $\C$.
\end{example}

We will show that the operations $\pi^2_1(x_i)$, for entries $x_i$ in a tuple $\ang{x}$, provide a partial inverse to $w^{\ang{x}}$.
Note, however, that the only morphism in $\Fr\M$ with codomain $\ang{}$ is the identity $1_{\ang{}}$.
Therefore, cases where $\ang{x}$ consists entirely of basepoint objects must be treated separately, as in the following.
\begin{definition}\label{definition:angxtil}
  Suppose $\M$ is a pointed multicategory with basepoint object $*$.
  Let $(*^r)$ denote the $r$-tuple consisting entirely of basepoints. Define the following for $\ang{x} \in \Fr\M$:
  \begin{equation}\label{eq:angxtil}
    \ang{x}^{\sim} =
    \begin{cases}
      (*) & \text{if $\ang{x} = (*^r)$ for $r > 0$}\\
      \ang{} & \text{if $\ang{x} = \ang{}$}\\
      \ang{x}^\wedge & \text{if $\ang{x}$ has any non-basepoint entries.}
    \end{cases}
  \end{equation}
  Thus, $\ang{x}^\sim$ is equal to $\ang{x}^\wedge$ unless $\ang{x} = (*^r)$ for $r > 0$.
\end{definition}

\begin{lemma}\label{lemma:c-angx}
  Suppose $\M$ is a left $\Mone$-module.
  Then, for each $\ang{x} \in \Fr\M$ there is a morphism
  \[
    c^{\ang{x}} \cn \ang{x} \to \ang{x}^\sim \inspace \Fr\M
  \]
  such that the composite
  \[
    \ang{x}^\wedge \fto{w^{\ang{x}}} \ang{x} \fto{c^{\ang{x}}} \ang{x}^\sim
  \]
  is $\simeq$-equivalent to the identity $1_{\ang{x}}$ in $\wt{\Mor}(\Fr\M)$.
\end{lemma} 
\begin{proof}
  Suppose $\ang{x}$ has length $r$.
  We use the operations $\pi^2_1(x)$ described in \cref{M1-mod-str} to show that there is a morphism $c^{\ang{x}}$ such that
  \begin{equation}\label{eq:c-angx-prop}
    c^{\ang{x}} \circ w^{\ang{x}} =
    \begin{cases}
      (\varnothing \fto{!} \ufs{1} \scs \iota_0)
      & \text{if $\ang{x} = (*^r)$ for $r > 0$},\\
      1_{\ang{x}^\wedge}
      & \text{otherwise},
    \end{cases}
  \end{equation}
  where $!$ denotes the unique index map.
  In either of the two cases \cref{eq:c-angx-prop}, we have the conclusion
  \[
    [c^{\ang{x}} \circ w^{\ang{x}}] = 1_{[\ang{x}]} \inspace \Fst\M = \Fm\M.
  \]

  First we discuss three trivial cases, where $\ang{x}$ consists of either $r > 0$ basepoint entries, or is empty, or consists of $r > 0$ non-basepoint entries.
  \begin{itemize}
  \item If $\ang{x} = (*^r)$ for $r > 0$, then the $r$-ary basepoint operation $\iota_r$ defines a morphism
  \[
    c^{(*^r)} = (\ufs{r} \fto{!} \ufs{1}, \iota_r) \cn (*^r) \to (*)
  \]
  as desired.
  \item If $r = 0$, note $w^{\ang{}} = 1_{\ang{}}$ and choose
    \[
      c^{\ang{}} = 1_{\ang{}}.
    \]
  \item If $r > 0$ and $\ang{x}$ has no basepoint entries, then $\ang{x} = \ang{x}^\wedge = \ang{x}^\sim$ and we choose
    \[
      c^{\ang{x}} = 1_{\ang{x}} = 1_{\ang{x}^\wedge}.
    \]
  \end{itemize}
  Each of these choices satisfies \cref{eq:c-angx-prop}.

  For the remainder of the proof we suppose $r \ge 2$ and $\ang{x}$ has at least one basepoint entry that is adjacent to one non-basepoint entry $x_i \ne *$.
  Note, therefore, that $\ang{x}^\sim = \ang{x}^\wedge$.
  If $r = 2$ and $\ang{x} = (x_1,*)$, then there is an operation
  \[
    \pi = \pi^2_1(x_1) \in \M\scmap{(x_1,*);x_1}
  \]
  determined by \cref{M1-mod-str}.
  In this case, choose
  \[
    c^{(x_1,*)} = (\ufs{2} \fto{!} \ufs{1} \scs \pi).
  \]
  If $r = 2$ and $\ang{x} = (*,x_2)$ then choose $c^{(*,x_2)}$ similarly, using the operation
  \[
    \pi = \pi^2_1(x_2)\tau \in \M\scmap{(*,x_2);x_2}
  \]
  given by the right action of the transposition $\tau\in \Si_2$.
  Note that
  \begin{equation}\label{eq:c-angx-ii}
    c^{(x_1,*)} \circ w^{(x_1,*)} = 1_{(x_1)}
    \andspace
    c^{(*,x_2)} \circ w^{(*,x_2)} = 1_{(x_2)}
    \inspace \wt{\Mor}(\Fr\M)
  \end{equation}
  by property \cref{eq:pi2cond-2} of the operations $\pi$ and hence these choices satisfy \cref{eq:c-angx-prop}.

  If $r > 2$, then $\ang{x}$ decomposes as a sum
  \[
    \ang{x} = \ang{x'} \oplus \ang{x''} \oplus \ang{x'''},
  \]
  where the first and last summands are possibly empty and the middle summand $\ang{x''}$ is either $(x_i,*)$ or $(*,x_i)$.
  Then choose $c^{\ang{x}}$ inductively as the composite along the top and right of the following diagram, where $\ang{y} = \ang{x'} \oplus (x_i) \oplus \ang{x'''}$, so $\ang{y}$ has length $r-1$, and $c^{\ang{y}}$ is any choice of morphism satisfying \cref{eq:c-angx-prop}.
  \[
    \begin{tikzpicture}[x=40mm,y=17mm]
      \draw[0cell] 
      (0,0) node (a) {
        \ang{x'} \oplus \ang{x''} \oplus \ang{x'''} 
      }
      (a)+(-6em,0) node (a') {
        \ang{x}
      }
      (a)+(1.5,0) node (b) {
        \ang{x'} \oplus (x_i) \oplus \ang{x'''}
      }
      (b)+(0,-2) node (c) {
        \ang{x}^\wedge
      }
      (a)+(0,-1) node (p) {
        \ang{x'} \oplus (x_i) \oplus \ang{x'''}
      }
      (a')+(0,-2) node (q) {
        \ang{x}^\wedge
      }
      (q)+(3.3em,0) node (q') {
        \ang{y}^\wedge
      }
      ;
      \draw[1cell] 
      (a') edge[equal] node {} (a)
      (q') edge[equal] node {} (q)
      (a) edge node {1 \oplus c^{\ang{x''}} \oplus 1} (b)
      (b) edge node {c^{\ang{y}}} (c)
      (q') edge['] node {w^{\ang{y}}} (p)
      (p) edge['] node {1 \oplus w^{(x_i)} \oplus 1} (a)
      (q) edge node {w^{\ang{x}}} (a')
      (q') edge node {1} (c)
      (p) edge['] node {1} (b)
      ;
    \end{tikzpicture}
  \]
  In the above diagram, the region at left commutes by definition of $w^{\ang{x}}$ and composition in $\Fr\M$ and the lower right region commutes by inductive hypothesis on $c^{\ang{y}}$.
  The upper triangle commutes by \cref{eq:c-angx-ii} for $c^{\ang{x''}}$, and hence the composite $c^{\ang{x}}$ also satisfies \cref{eq:c-angx-prop}.
  This completes the proof.
\end{proof}
Combining \cref{lemma:c-angx} with \cref{eq:fphi-simeq-fphiw} gives an application that we explain in \cref{lemma:nonsing-reps} below.
\begin{explanation}[Non-Uniqueness of $c^{\ang{x}}$]\label{expl:c-angx-nonunq}
  Note that the morphisms $c^{\ang{x}}$ constructed in the proof of \cref{lemma:c-angx} generally depend on a choice of non-basepoint entry $x_i$ in $\ang{x}$ and are not necessarily unique.
  However, the condition
  \[
    c^{\ang{x}} \circ w^{\ang{x}} \simeq 1_{\ang{x}} \inspace \wt{\Mor}(\Fr\M),
  \]
  together with the observation $w^{\ang{x}} \iisim 1_{\ang{x}}$ from \cref{expl:w-angx}, implies that $c^{\ang{x}} \simeq 1_{\ang{x}}$ for any choice of $c^{\ang{x}}$.
  Therefore, all such choices result in identity morphisms in $\Fst\M$.
  This is the key property of the morphisms $c^{\ang{x}}$.
\end{explanation}

The following result uses the notation of \cref{expl:w-angx,definition:angxtil,lemma:c-angx} above.
\begin{lemma}\label{lemma:nonsing-reps}
  Suppose $\M$ is a pointed multicategory and suppose that, for each object $\ang{x} \in \Fr\M$, there is a morphism
  \[
    c^{\ang{x}} \cn \ang{x} \to \ang{x}^\sim \inspace \Fr\M
  \]
  such that the composite
  \[
    \ang{x}^\wedge \fto{w^{\ang{x}}} \ang{x} \fto{c^{\ang{x}}} \ang{x}^\sim
  \]
  is $\simeq$-equivalent to the identity $1_{\ang{x}}$ in $\wt{\Mor}(\Fr\M)$.
  Then each morphism $[\ang{x}] \to [\ang{y}]$ in $\Fst\M$ is represented by a length-one sequence consisting of a morphism
  \[
    (f',\ang{\phi'}) \cn \ang{x}^\wedge \to \ang{y}^\sim \inspace \Fr\M.
  \]
\end{lemma}
\begin{proof}
  Applying the observation \cref{eq:fphi-simeq-fphiw} to $c^{\ang{x}}$, we have
  the following equivalences in $\wt{\Mor}(\Fr\M)$:
  \[
    \big( c^{\ang{x}} \big) \simeq \brb{w^{\ang{x}},c^{\ang{x}}} \isim \big( c^{\ang{x}}\circ w^{\ang{x}} \big) \simeq \big( 1_{\ang{x}} \big),
  \]
  where the final relation holds by hypothesis on $c^{\ang{x}}$.
  Now suppose given a morphism
  \[
    \brb{f,\ang{\phi}} \cn \ang{x} \to \ang{y} \inspace \Fr\M.
  \]
  Let $(f',\ang{\phi'})$ and $(f'',\ang{\phi''})$ denote the composites indicated below.
  \begin{equation}\label{eq:f'f''}
    \begin{tikzpicture}[x=30mm,y=15mm,vcenter]
      \draw[0cell] 
      (0,0) node (a) {\ang{x}}
      (1,0) node (b) {\ang{y}}
      (0,-1) node (c') {\ang{x}^\sim}
      (0,-2) node (c) {\ang{x}^\wedge}
      (1,-1) node (d) {\ang{y}^\sim}
      ;
      \draw[1cell] 
      (a) edge node {(f,\ang{\phi})} (b)
      (c') edge[dashed] node {(f',\ang{\phi'})} (d)
      (c) edge[dashed,'] node {(f'',\ang{\phi''})} (d)
      (c) edge[bend left=40] node {w^{\ang{x}}} (a)
      (b) edge node {c^{\ang{y}}} (d)
      (c) edge['] node {} (c')
      (c') edge['] node {\wt{w}^{\ang{x}}} (a)
      ;
    \end{tikzpicture}
  \end{equation}
  In the above diagram, either
  \begin{itemize}
  \item $\ang{x}$ has some non-basepoint operation, in which case the unlabeled morphism $\ang{x}^\wedge \to \ang{x}^\sim$ is the identity and $\wt{w}^{\ang{x}}$ is equal to $w^{\ang{x}}$, or
  \item $\ang{x} = (*^r)$ consists of all basepoint entries, in which case the unlabeled morphism, $\wt{w}^{\ang{x}}$, and $w^{\ang{x}}$, are all determined by basepoint operations and a choice of index map $\ufs{1} \to \ufs{r}$ for $\wt{w}^{\ang{x}}$.
  \end{itemize}
  In either case, note that $\wt{w}^{\ang{x}}$ is also $\simeq$-equivalent to an identity.
  This conclusion holds independently of which index map $\ufs{1} \to \ufs{r}$ we choose for $\wt{w}^{\ang{x}}$ in the case $\ang{x} = (*^r)$.

  Since each of $w^{\ang{x}}$, $\wt{w}^{\ang{x}}$, and $c^{\ang{y}}$ is $\simeq$-equivalent to an identity, we have
  \[
    \brb{f',\ang{\phi'}} \simeq \brb{f,\ang{\phi}} \simeq \brb{f'',\ang{\phi''}}
    \inspace \wt{\Mor}(\Fr\M).
  \]
  This proves the assertion for morphisms of the form $\big[ \brb{f,\ang{\phi}} \big]$ in $\Fst\M$.

  For a pair of $\obsim$-composable morphisms,
  \[
    \ang{x} \fto{(f,\ang{\phi})} \ang{y} \obsim \ang{u} \fto{(g,\ang{\psi})} \ang{v},
  \]
  note that we have $\ang{y}^\wedge = \ang{u}^\wedge$ and $\ang{y}^\sim = \ang{u}^\sim$ by definition of $\obsim$ and \cref{eq:angxtil}.
  Considering the following diagram,
  \[
    \begin{tikzpicture}[x=30mm,y=15mm]
      \draw[0cell] 
      (0,0) node (a) {\ang{x}}
      (1,0) node (b) {\ang{y}}
      (0,-1) node (c') {\ang{x}^\sim}
      (0,-2) node (c) {\ang{x}^\wedge}
      (1,-1) node (d) {\ang{y}^\sim}
      ;
      \draw[1cell] 
      (a) edge node {(f,\ang{\phi})} (b)
      (c') edge node {(f',\ang{\phi'})} (d)
      (c) edge['] node {(f'',\ang{\phi''})} (d)
      (c) edge[bend left=40] node {w^{\ang{x}}} (a)
      (b) edge node {c^{\ang{y}}} (d)
      (c) edge['] node {} (c')
      (c') edge['] node {\wt{w}^{\ang{x}}} (a)
      ;
      \draw[0cell] 
      (b)+(.4,0) node (u) {\ang{u}}
      (u)+(1,0) node (v) {\ang{v}}
      (u)+(0,-1) node (u') {\ang{u}^\sim}
      (u')+(1,0) node (v') {\ang{v}^\sim}
      ;
      \draw[1cell] 
      (u) edge node {(g,\ang{\psi})} (v)
      (u') edge['] node {(g',\ang{\psi'})} (v')
      (v) edge node {c^{\ang{v}}} (v')
      ;
      \draw[1cell]
      (u') edge['] node {\wt{w}^{\ang{u}}} (u)
      ;
      \draw[1cell]
      node[between=b and u at .5] {\obsim}
      node[between=d and u' at .5] {=}
      ;
    \end{tikzpicture}
  \]
  we conclude that
  \[
    \big[ (f,\ang{\phi})\scs (g,\ang{\psi}) \big] \simeq
    \big[ (f'',\ang{\phi''})\scs (g',\ang{\psi'}) \big] \isim
    \big[ (g',\ang{\psi'}) \circ (f'',\ang{\phi''}) \big].
  \]
  Therefore, by induction, each tuple of $\obsim$-composable morphisms in $\wt{\Mor}(\Fr\M)$ is equivalent to a single morphism in $\Fr\M$, with domain of the form $\ang{x}^\wedge$ and codomain of the form $\ang{y}^\sim$.
  This completes the proof.
\end{proof}

Combining \cref{lemma:c-angx,lemma:nonsing-reps} yields the following result.
\begin{corollary}\label{corollary:FmM-reps}
  Suppose $\M$ is a left $\Mone$-module.
  Then each morphism $[\ang{x}] \to [\ang{y}]$ in $\Fm\M = \Fst\Um\M$ is represented by a length-one sequence consisting of a morphism
  \[
    (f',\ang{\phi'}) \cn \ang{x}^\wedge \to \ang{y}^\sim \inspace \Fr\M.
  \]
\end{corollary}

\begin{example}[Partition Multicategories]\label{example:FstMa}
  Applying \cref{corollary:FmM-reps} for a partition multicategory $\M = \cM a$, note that the only partitions of the empty subset $\varnothing \subset a^\punc$ are those consisting of empty subsets.
  Therefore, the only morphisms with codomain $(\varnothing) \in \Fr(\cM a)$ are those that are $\iisim$-equivalent to $1_{\ang{}}$.
  Thus \cref{corollary:FmM-reps} simplifies slightly to give the following description of the permutative category $\Fm(\cM a) = \Fst\Um(\cM a)$.
  \begin{description}
  \item[Objects] Each object of $\Fm(\cM a)$ is uniquely represented by a tuple $\ang{s}$ such that each entry $s_i$ is a nonempty subset of $a^\punc$.
  \item[Morphisms] Each morphism of $\Fm(\cM a)$ is uniquely represented by a morphism
    \[
      (f,\ang{\phi})\cn \ang{s} \to \ang{t}
    \]
    where each $s_i$ and $t_j$ is nonempty and each $\ang{s}_{f^\inv(j)}$ is a partition of $t_j$, with corresponding operation $\phi_j = \iota_{\ang{s}_{f^\inv(j)}}$ in $\cM a$.
  \end{description}

  As a special case, with $a = \ord{1}$, the discussion above shows that $\Fm(\Mone)$ is isomorphic to the permutation category with
  \begin{itemize}
  \item objects given by natural numbers $n \ge 0$, corresponding to length-$n$ tuples of the subset $s = \{1\}$, and
  \item morphisms given by symmetric groups, whose elements permute the entries of a tuple.
  \end{itemize} 
  Thus, we conclude
  \begin{equation}\label{eq:FstMone-FrMtu}
    \Fm(\Mone) \iso \Fr(\Mtu).
  \end{equation}
  In fact, this isomorphism is a special case $\M = \Mtu$ combining the general isomorphisms
  \[
    \Fr\M \iso \Fst(\M_+) \iso \Fm(\Mone \sma \M_+)
  \]
  from \cref{FFstdplus,FstFmMonesma} above.
\end{example}

\Cref{corollary:FmM-reps} can be used for an application similar to that of \cref{example:FstMa} in the case $\M = \Endst\C$, although there may be nontrivial morphisms $x \to \pu$ in $\C$ and hence this case does not admit the same simplification as in \cref{example:FstMa}.
An alternative application is given in \cref{ptmulticat-prop-viii} below, where we use \cref{expl:w-angx} to show that there is an adjunction in $\permcatsu$ between $\C$ and $\Fst\Endst\C$.

\section{Componentwise Right Adjoint of the Pointed Adjunction}
\label{sec:epzst-right-adj}

In this section we extend the componentwise right adjoint $\vrho_\C$ from \cref{proposition:epz-rho-adj} to the pointed context.

\begin{definition}\label{def:vrhostC}
  For each small permutative category $\C$, we define the symmetric monoidal functor $\vrhost_\C$ as the composite
  \begin{equation}\label{vrhostC-def}
    \begin{tikzpicture}[vcenter]
      \def\u{.6}
      \draw[0cell]
      (0,0) node (a) {\C}
      (a)+(3.2,0) node (b) {\Fr\,\End(\C) = \Fr\,\Endst\C}
      (b)+(4.4,0) node (c) {\Fst\, \Endst\C}
      ;
      \draw[1cell]
      (a) edge node {\vrho_\C} (b)
      (b) edge node {\pst_{\Endst\C}} (c)
      ;
      \draw[1cell]
      (a) [rounded corners=3pt] |- ($(b)+(-1,\u)$)
      -- node {\vrhost_\C} ($(b)+(1,\u)$) -| (c)
      ;
    \end{tikzpicture}
  \end{equation}
  of
  \begin{itemize}
  \item the symmetric monoidal functor $\vrho_\C$ in \cref{rhoc-monoidal} and
  \item the strict symmetric monoidal functor $\pst_{\Endst\C}$ in \cref{def:FFst-projection}.\defmark
  \end{itemize}
\end{definition}

\begin{explanation}\label{expl:rhoc-not-su}
  The symmetric monoidal functor $\vrhost_\C$ in \cref{vrhostC-def} is given explicitly by combining the definition of $\vrho_\C$ from \cref{rhoc-assignments} with that of $\pst$ from \cref{eq:pst}:
  \begin{equation}\label{rhostc-assignments}
    \left\{
      \begin{aligned}
        \vrhost_\C (x) & = [(x)]
        & \forspace & x \in \ObC \andspace\\
        \vrhost_\C(\phi) & = \big[\brb{1_{\ufs{1}},(\phi)}\big] \cn [(x)] \to [(y)]
        & \forspace & \phi \in \C(x,y).
      \end{aligned}
    \right.
  \end{equation}
  We emphasize that $\vrho_\C$ is \emph{not} strictly unital, as discussed in \cref{rk:rhonotunital}. 
  Despite this, \cref{ptmulticat-xxi} below shows that $\vrhost_\C$ \emph{is} strictly unital.
\end{explanation}

\begin{lemma}\label{ptmulticat-xxi}
  For each small permutative category $\C$, the symmetric monoidal functor 
  \[
    \vrhost_\C \cn \C \to \Fst\, \Endst\C
  \]
  in \cref{vrhostC-def} is strictly unital. 
\end{lemma}
\begin{proof}
  The monoidal and unit constraints for $\vrho_\C$ are given in \cref{def:rhoc-monoidal}.
  Since $\pst_{\Endst\C}$ is strict symmetric monoidal, the monoidal and unit constraints of $\vrhost_\C$ are given by $\pst_{\Endst\C}(\vrho_\C^2)$ and $\pst_{\Endst\C}(\vrho_\C^0)$, respectively.
  By \cref{rhoc-unit}, 
  \[
    \vrho_\C^0 = (\iota_0,1_\pu)\cn \ang{} \to (\pu).
  \]
  Therefore $\pst_{\End\C}(\vrho_\C^0) = \big[ (\iota_0,1_\pu) \big] = 1_{[\ang{}]}$, by relation $\iisim$.
\end{proof}

\begin{proposition}\label{ptmulticat-prop-viii}
  For each small permutative category $\C$, the adjunction $\epz_\C \dashv \vrho_\C$ in \cref{epzc-rhoc} extends along $\pst_{\Endst(\C)}$ to an adjunction in $\permcatsu$.
  \begin{equation}\label{epzstc-rhostc}
    \begin{tikzpicture}[baseline={(a'.base)}]
      \draw[0cell]
      (0,0) node (a') {\Fst\, \Endst\C}
      (a')+(.55,0) node (a) {\phantom{\C}}
      (a)+(2,.03) node (b) {\C}
      (a)+(1,0) node (x) {\bot}
      ;
      \draw[1cell=.9]
      (a) edge[bend left=15] node {\epzst_\C} (b)
      (b) edge[bend left=15,transform canvas={yshift=-.2ex}] node {\vrhost_\C} (a)
      ;
    \end{tikzpicture}
  \end{equation}
  with
  \begin{itemize}
  \item $\epzst_\C$ the counit in \cref{epzst-def} and
  \item $\vrhost_\C$ the strictly unital symmetric monoidal functor in \cref{ptmulticat-xxi}.
  \end{itemize}
\end{proposition}
\begin{proof}
  Throughout this proof we write
  \[
  \Er = \End,\quad \Est = \Endst, \andspace \pst = \pst_{\Est\C}.
  \]
  Recalling \cref{proposition:epz-rho-adj}, the adjunction $\epz_\C \dashv \vrho_\C$ has unit and counit
  \[
    1_{\Fr\Er\C} \fto{\ups} \vrho_\C\epz_\C
    \andspace
    \epz_\C\vrho_\C \fto{1} 1_\C,
  \]
  respectively, where $\ups$ is defined in \cref{eq:ups-x}.
  \Cref{rhoc-monoidal,rk:rhonotunital} explain that $\rho_\C$ is a symmetric monoidal functor that is generally not strong monoidal.

  To show that $\ups$ gives a well defined a unit\label{not:upsst}
  \[
    \upsst \cn 1_{\Fst\Est\C} \to \vrhost_\C\epzst_\C,
  \]
  first recall from \cref{eq:fphi-simeq-fphiw} with $\M = \Est\C$ that we have
  \begin{equation}\label{eq:fphi-simeq-fphiw-ii}
    (f,\ang{\phi}) \simeq \big( w^{\ang{x}} \scs (f,\ang{\phi}) \big)
    \inspace \wt{\Mor}(\Fr\Est\C),
  \end{equation}
  where $(f,\ang{\phi})\cn \ang{x} \to \ang{z}$ is any morphism in $\Fr\Est\C$ and
  \[
    w^{\ang{x}} \cn \ang{x}^\wedge \to \ang{x}
  \]
  is determined by the nullary basepoint operations (\cref{expl:w-angx}).

  Now suppose $\ang{x} \obsim \ang{y}$ in $\Fr\Est\C$.
  Then
  \begin{itemize}
  \item $\ang{x}^\wedge = \ang{y}^\wedge$ by definition of $\obsim$,
  \item $\oplus_i x_i = \oplus_j y_j$ by strictness of the monoidal unit in $\C$, and
  \item the following equalities hold in $\Fr\Est\C$:
    \[
      \ups_{\ang{x}} \circ w^{\ang{x}}
      = \ups_{\ang{x}^\wedge}
      = \ups_{\ang{y}^\wedge}
      = \ups_{\ang{y}} \circ w^{\ang{y}}.
    \]
  \end{itemize}
  Therefore, by \cref{eq:fphi-simeq-fphiw-ii} with $\M = \Est\C$ we have
  \[
    \big( \ups_{\ang{x}} \big)
    \simeq
    \big( w^{\ang{x}} \scs \ups_{\ang{x}} \big)
    \simeq
    \big( w^{\ang{y}} \scs \ups_{\ang{y}} \big)
    \simeq
    \big( \ups_{\ang{y}} \big).
  \]
  Hence, the components
  \begin{equation}\label{eq:upsst}
    \upsst_{[\ang{x}]} = [\ups_{\ang{x}}] \cn [\ang{x}] \to [(\oplus_i x_i)]
    \inspace \Fst\Est\C
  \end{equation}
  are well defined with respect to $\obsim$-equivalence classes of objects.
  With this definition, naturality of $\ups$ implies that of $\upsst$ and we have
  \[
    \pst * \ups = \upsst * \pst
  \]
  in the following diagram in $\permcatsu$.
  \begin{equation}\label{eq:upsst-diagram}
    \begin{tikzpicture}[x=35mm,y=16mm,vcenter]
      \def\veq{6.5mm}
      \draw[0cell=.9] 
      (0,0) node (a) {\Fr\Er\C}
      (a)+(0,-\veq) node (b) {\Fr\Est\C}
      (a)+(.5,-.25) node (c) {\C}
      (a)+(0,-1) node (d) {\Fst\Est\C}
      (1,0) node (a') {\Fr\Er\C}
      (a')+(0,-\veq) node (b') {\Fr\Est\C}
      (d)+(.5,-.25) node (c') {\C}
      (a')+(0,-1) node (d') {\Fst\Est\C}
      ;
      \draw[1cell=.9] 
      (a) edge[equal] node {} (b)
      (b) edge['] node {\pst} (d)
      (a') edge[equal] node {} (b')
      (b') edge node {\pst} (d')
      (a) edge[bend left=15] node (T) {1} (a')
      (a) edge['] node {\epz_\C} (c)
      (c) edge['] node {\vrho_\C} (a')
      (d) edge['] node {\epzst_\C} (c')
      (c') edge['] node {\vrhost_\C} (d')
      (d) edge[bend left=15] node (T') {1} (d')
      ;
      \draw[2cell]
      node[between=T and c at .5, rotate=-90, 2label={above,\ups}]
      {\Rightarrow}
      node[between=T' and c' at .5, rotate=-90, 2label={above,\upsst}] {\Rightarrow}
      ;
    \end{tikzpicture}
  \end{equation}

  Since $\epzst_\C \pst = \epz_\C$, we also have
  \[
    \epzst_\C\vrhost_\C = 1_\C.
  \]
  One of the triangle identities for $(\epzst_\C,\vrhost_\C,\upsst,1_\C)$ is, therefore, immediate.
  The other follows from the corresponding identity for $(\epz_\C,\rho_\C,\ups,1_\C)$ on representative objects and morphisms.
\end{proof}

\begin{explanation}\label{explanation:upsst}
  Note that our proof of \cref{ptmulticat-prop-viii} depends on specific details from \cref{expl:w-angx} instead of the 2-categorical pushout description of $\Fst\M$ in \cref{proposition:FstM-pushout}.
  This is because \cref{proposition:FstM-pushout} describes a pushout in $\permcatst$, while \cref{eq:upsst-diagram} is a diagram in $\permcatsu$.
  Thus, a more general proof of \cref{ptmulticat-prop-viii} using \cref{proposition:FstM-pushout} would require a comparison of 2-dimensional pushouts.
  The strictification theory of \cite{bkp} provides some comparisons of this sort, and may be one way of approaching more general versions of the results here and below.
\end{explanation}

\begin{remark}\label{remark:ptmulticat-prop-viii}
  Our first main application of \cref{ptmulticat-prop-viii} is \cref{ptmulticat-thm-x}, showing that the components of $\epzst$ are stable equivalences.
  The further applications of \cref{ptmulticat-prop-viii}, in \cref{ptmulticat-xxiii,ptmulticat-xxv,mackey-xiv-pmulticat,mackey-xiv-mone} use the same result for $\vrhost$ instead.
  The reason for this change is that $\vrhost$ is shown to be $\Cat$-multinatural in \cref{ptmulticat-xxii}, while the corresponding result for $\epzst$ does not hold (see \cite[Remark~10.4]{johnson-yau-Fmulti}).
\end{remark}

Recall from \cref{ex:endstc}~\cref{EndstP} that each \emph{strictly unital} symmetric monoidal functor $P$ induces a \emph{pointed} multifunctor $\Endst(P)$.
Since the components of $\vrhost$ are strictly unital, by \cref{ptmulticat-xxi}, we may consider $\Endst\vrhost$.
For comparison, recall the description of $\etast$ from \cref{expl:etast-explicit}.
The following result is used in the proof of \cref{mackey-xiv-pmulticat}, step
\cref{mackey-xiv-pm-v}.

\begin{lemma}\label{etaEEvrho}
  Suppose $\C$ is a small permutative category.
  Then the two pointed multifunctors below are equal.
  \[
    \begin{tikzpicture}[baseline={(a.base)}]
      \draw[0cell]
      (0,0) node (a) {\Endst \C}
      (a)+(3.5,0) node (b) {\Endst\Fst\Endst \C}
      ;
      \draw[1cell=.9]
      (a) edge[transform canvas={yshift=.6ex}] node {\etast_{\Endst\C}} (b)
      (a) edge[transform canvas={yshift=-.5ex}] node[swap] {\Endst \vrhost_\C} (b)
      ;
    \end{tikzpicture}
  \]
\end{lemma}
\begin{proof}
  We need to show that $\etast_{\Endst\C}$ and $\Endst \vrhost_\C$ have (i) the same object assignment and (ii) the same multimorphism assignment.
  For each object $x$ in $\C$, there are object equalities
  \[
    \etast_{\Endst\C}(x) = [(x)] = (\Endst \vrhost_\C)(x) \inspace \Endst\Fst\Endst \C.
  \]
  This proves (i).

  To prove (ii), consider an $r$-ary multimorphism
  \[
    \psi \in (\Endst\C)\scmap{\ang{x_i}_{i=1}^r; y} = \C\brb{\txoplus_{i=1}^r x_i , y}.
  \]
  Then there are equalities
  \begin{equation}\label{etaEEvrho-morphism}
    \etast_{\Endst\C}(\psi) 
    = \left[\big(\iota_r \scs (\psi)\big)\right] 
    = (\Endst \vrhost_\C)(\psi)
  \end{equation}
  in
  \[
    (\Endst\Fst\Endst \C)\scmap{\ang{[(x_i)]}_{i=1}^r; [(y)]} 
    = (\Fst\Endst \C)\brb{[\ang{x_i}_{i=1}^r], [(y)]}
  \]
  for the following reasons.
  \begin{itemize}
  \item The first equality in \cref{etaEEvrho-morphism} follows from definition \cref{Funit-operation} of 
    \[
      \eta_{\End\C}(\psi) = \big(\iota_r \scs (\psi)\big) \cn \ang{x_i}_{i=1}^r \to (y)
    \]
    and the fact that each component of $\pst$ is a strict symmetric monoidal functor (\cref{def:FFst-projection}).
  \item By definition \cref{rhoc-assignments}, we have the equality
    \[
      \vrho_\C(\psi) = \big(1_{\ufs{1}} \scs (\psi)\big) \cn \left(\txoplus_{i=1}^r x_i\right) \to (y).
    \]
    The symmetric monoidal structure of $\vrho_\C$ in \cref{rhoc-monoidal} implies that $(\Endst \vrhost_\C)(\psi)$ is the following composite.
    \[
      [\ang{x_i}_{i=1}^r] 
      \fto{[(\iota_r \scs 1_{\oplus_i x_i})]} [(\txoplus_{i=1}^r x_i)] 
      \fto{[(1_{\ufs{1}} \scs (\psi))]} [(y)]
    \]
    This composite is equal to the middle entry in \cref{etaEEvrho-morphism} because
    \[
      \begin{aligned}
        1_{\ufs{1}} \circ \iota_r &= \iota_r \cn \ufs{r} \to \ufs{1} \andspace\\
        \psi \circ 1_{\oplus_i x_i} &= \psi \cn \txoplus_{i=1}^r x_i \to y.
      \end{aligned}
    \]
  \end{itemize}
  This proves \cref{etaEEvrho-morphism}.
\end{proof}

\section{Homotopy Theory of Pointed Multicategories}
\label{sec:ptmultperm-heq}

Recall the notions about relative categories in \cref{definition:rel-cat}.
We extend the relative category structure on $\permcatst$ to $\pMulticat$ and $\MoneMod$ as follows.  

\begin{definition}[Stable Equivalences]\label{def:ptmulti-stableeq}
  We define the wide subcategories
  \[
    \begin{aligned}
      \cSM &= (\Fm)^\inv(\cSI) \bigsubset \MoneMod \andspace\\
      \cSst &= \Fst^\inv(\cSI) \bigsubset \pMulticat
    \end{aligned}
  \]
  as the subcategories created by the indicated functors below.
  \begin{equation}\label{SM-Sst-def}
    \begin{tikzpicture}[vcenter]
      \draw[0cell]
      (0,0) node (a) {\brb{\MoneMod, \cSM}}
      (a)+(0,-1) node (b) {\brb{\pMulticat, \cSst}}
      (a)+(3.5,-.5) node (c) {\phantom{X}}
      (c)+(1.1,0) node (d) {\brb{\permcatst,\cSI}}
      ;
      \draw[1cell]
      (a) edge node[pos=.3] {\Fm} (c)
      (b) edge node[swap,pos=.3] {\Fst} (c)
      ;
    \end{tikzpicture}
  \end{equation}
  \begin{itemize}
  \item $\cSI = I^\inv(\cS) \bigsubset \permcatst$ is the wide subcategory in \cref{def:mult-stableeq}.
  \item $\Fm$ is the underlying functor of the 2-functor in \cref{Fm-def}.
  \item $\Fst$ is the underlying functor of the 2-functor in \cref{ptmulticat-thm-i}.
  \end{itemize}
  We refer to morphisms in $\cSM$ and $\cSst$ as \index{multicategory!pointed - stable equivalence}\index{equivalence!stable!of pointed multicategories}\index{stable equivalence!of pointed multicategories}\emph{$\Fm$-stable equivalences} and \emph{$\Fst$-stable equivalences}, respectively.
\end{definition}

Note that $\Fst$-stable equivalences are the preimages of the stable equivalences in $\cS \bigsubset \permcatsu$ \cref{perm-steq} under the following composite.
\[
  \begin{tikzcd}[column sep=large]
    \pMulticat \ar{r}{\Fst} & \permcatst \ar[hookrightarrow]{r}{I} & \permcatsu
  \end{tikzcd}
\]
Recall the notion of an \index{equivalence!adjoint - of homotopy theories}\index{homotopy theory!adjoint equivalence of}\index{adjoint equivalence!of homotopy theories}\emph{adjoint equivalence of homotopy theories} in \cref{def:heq}.

\begin{theorem}\label{ptmulticat-thm-x}
  The adjunction in \cref{ptmulticat-thm-v}
  \[
    \begin{tikzpicture}[vcenter]
      \def\t{17}
      \draw[0cell]
      (0,0) node (a) {\phantom{X}}
      (a)+(2,0) node (b) {\phantom{X}}
      (a)+(1,0) node (x) {\bot}
      (a)+(-1.05,0) node (a') {\brb{\pMulticat,\cSst}}
      (b)+(1.1,0) node (b') {\brb{\permcatst,\cSI}}
      ;
      \draw[1cell=.9]
      (a) edge[bend left=\t] node {\Fst} (b)
      (b) edge[bend left=\t] node {\Endst} (a)
      ;
    \end{tikzpicture}
  \]
  is an adjoint equivalence of homotopy theories.
\end{theorem}
\begin{proof}
  The left adjoint $\Fst$ is a relative functor by definition of $\cSst$.
  To see that $\Endst$ is a relative functor, first recall that each $\epzst_\C$ is a left adjoint in $\permcatsu$ by \cref{ptmulticat-prop-viii}.
  Thus, the components of $\epzst$ are stable equivalences by \cref{remark:steq}~\cref{it:steq-2}.
  Naturality of $\epzst$ and the 2-out-of-3 property for stable equivalences (\cref{remark:steq}~\cref{it:steq-1}) then imply that
  \[
    \Fst\Endst P \cn \Fst\Endst\C  \to \Fst\Endst\D
  \]
  is a stable equivalence whenever $P$ is a stable equivalence.
  This, in turn, implies that $\Endst P$ is an $\Fst$-stable equivalence.
  Hence $\Endst$ is a relative functor.

  The triangle identities for $\etast$ and $\epzst$, together with the 2-out-of-3 property, imply that the components of $\etast$ are also stable equivalences.
  This completes the proof.
\end{proof}

Recall from \cref{def:mult-stableeq} the wide subcategory 
\[\cSF = \Fr^\inv(\cSI) \bigsubset \Multicat\]
of $\Fr$-stable equivalences.  Also recall the 2-functor in \cref{dplus-def}
\[\dplus \cn \Multicat \to \pMulticat\]
given by adjoining a basepoint.  \cref{ptmulticat-cor-page-vi} below uses the underlying functor of this 2-functor.

\begin{corollary}\label{ptmulticat-cor-page-vi}
The functor
\[\dplus \cn \brb{\Multicat,\cSF} \fto{\sim} \brb{\pMulticat,\cSst}\]
is an equivalence of homotopy theories, with $\cSst$ as in \cref{def:ptmulti-stableeq}.
\end{corollary}
\begin{proof}
  The result follows by the 2-out-of-3 property for equivalences of homotopy theories, using \cref{thm:F-heq,ptmulticat-thm-x} together with \cref{ptmulticat-thm-v}~\cref{FEst-adj-ii}.
\end{proof}

\section{Homotopy Theory of \texorpdfstring{$\Mone$}{M1}-Modules}
\label{sec:moneperm-heq}

\begin{theorem}\label{ptmulticat-thm-vi}
  The adjunction in \cref{MonesmaUmadj}
  \begin{equation}\label{eq:MonesmaUmadj-steq}
    \begin{tikzpicture}[vcenter]
      \def\t{17}
      \draw[0cell]
      (0,0) node (a) {\phantom{X}}
      (a)+(2.4,0) node (b) {\phantom{X}}
      (a)+(1.2,0) node (x) {\bot}
      (a)+(-1,0) node (a') {\brb{\pMulticat,\cSst}}
      (b)+(1.05,0) node (b') {\brb{\MoneMod,\cSM}}
      ;
      \draw[1cell=.9]
      (a) edge[bend left=\t] node {\Monesma} (b)
      (b) edge[bend left=\t] node {\Um} (a)
      ;
    \end{tikzpicture}
  \end{equation}
  is an adjoint equivalence of homotopy theories, with $\cSst$ and $\cSM$ as in \cref{def:ptmulti-stableeq}.
\end{theorem}
\begin{proof}
  The left adjoint, $\Mone \sma -$, creates $\Fst$-stable equivalences because, by \cref{ptmulticat-thm-ii}~\cref{FEm-adj-iii}, there is a 2-natural isomorphism
  \[
    \Fst \iso \Fm \circ (\Mone \sma -).
  \]
  The right adjoint, $\Um$, creates $\Fm$-stable equivalences by definition of $\Fm$ \cref{Fm-def}.

  Recalling \cref{expl:MonesmaUmadj}, the counit of the adjunction \cref{eq:MonesmaUmadj-steq} is a componentwise isomorphism and, therefore, a componentwise $\Fm$-stable equivalence.
  Using the 2-out-of-3 property and the fact that $\Mone \sma -$ creates stable equivalences, the left triangle identity then shows that the unit of \cref{eq:MonesmaUmadj-steq} is a componentwise $\Fst$-stable equivalence.
\end{proof}

\begin{theorem}\label{ptmulticat-thm-xi}
  The adjunction in \cref{ptmulticat-thm-ii}
  \[\begin{tikzpicture}
      \def\t{17}
      \draw[0cell]
      (0,0) node (a) {\phantom{X}}
      (a)+(2,0) node (b) {\phantom{X}}
      (a)+(1,0) node (x) {\bot}
      (a)+(-1.05,0) node (a') {\brb{\MoneMod,\cSM}}
      (b)+(1.1,0) node (b') {\brb{\permcatst,\cSI}}
      ;
      \draw[1cell=.9]
      (a) edge[bend left=\t] node {\Fm} (b)
      (b) edge[bend left=\t] node {\Endm} (a)
      ;
    \end{tikzpicture}\]
  is an adjoint equivalence of homotopy theories, with $\cSM$ as in \cref{def:ptmulti-stableeq}.
\end{theorem}
\begin{proof}
  Recall from \cref{ptmulticat-thm-x} that the components of $\etast$ and $\epzst$ are stable equivalences.
  By \cref{expl:etam-epzm} we have
  \[
    \Um \etam_{\M} = \etast_{\Um\,\M}.
  \]
  This shows that $\etam$ is a componentwise $\Fm$-stable equivalence
  because $\Um$ creates $\Fm$-stable equivalences \cref{Fm-def}.
  Likewise, $\epzm$ is a componentwise stable equivalence by \cref{eq:epzst-epzm}.
\end{proof}

\chapter[Multiplicative Homotopy Theory]{Multiplicative Homotopy Theory of Pointed Multicategories and \texorpdfstring{$\Mone$}{M1}-Modules}
\label{ch:ptmulticat-alg}

The main results of this chapter extend the equivalences of homotopy theories between categories of non-symmetric $\Q$-algebras, for $\Q$ a small non-symmetric $\Cat$-multicategory, in \cref{thm:alg-hty-equiv},
\[\Fr^\Q \cn \brb{\Multicat^\Q, (\cSF)^\Q} \lrsimadj \brb{(\permcatsu)^\Q, \cS^\Q} 
\cn \End^\Q,\]
from $\Multicat^\Q$ to $\pMulticat^\Q$ and $(\MoneMod)^\Q$.
The three pairs of functors in the diagram below are shown to be inverse equivalences of homotopy theories in \cref{ptmulticat-xxiii,ptmulticat-xxv,Monesma-Um-algebra}.
\begin{equation}\label{ptmulticatalg-summary}
\begin{tikzpicture}[x=12mm,y=12mm,vcenter]
\def\s{10} \def\t{18} \def\h{3.5}
\draw[0cell=.8]
(0,0) node (a) {\brb{\pMulticat^\Q,\cSst^\Q}}
(a)+(\h,0) node (b) {\brb{(\permcatsu)^\Q,\cS^\Q}}
(a)+(\h/2,-1.75) node (d) {\brb{(\MoneMod)^\Q,(\cSM)^\Q}}
;
\draw[1cell=.8]
(a) edge[transform canvas={yshift=-.7ex},bend left=\s] node {\Fst^\Q} (b)
(b) edge[transform canvas={yshift=-.5ex}] node {\Endst^\Q} (a)
;
\draw[1cell=.8]
(d) edge[transform canvas={xshift=2ex},bend right=\t] node[swap] {\Fm^\Q} (b)
(b) edge[transform canvas={xshift=1ex}] node[swap,pos=.6] {\Endm^\Q} (d)
(d) edge[transform canvas={xshift=-1ex}] node[swap,pos=.4] {\Um^\Q} (a)
(a) edge[transform canvas={xshift=-2ex},bend right=\t] node[swap] {(\Monesma)^\Q} (d)
;
\end{tikzpicture}
\end{equation}

\subsection*{Connection with Other Chapters}

In \cref{ch:mackey_eq} we show, after developing the relevant basic theory in the intervening chapters, that the multifunctors
$\Fst$ and $\Fm$ also induce equivalences of homotopy theories between enriched diagram categories.

\subsection*{Background}

Recall \cref{def:inverse-heq} for inverse equivalences of homotopy theories.
\Cref{ch:ptmulticat-sp} describes the underlying adjoint pairs
$\Fst \dashv \Endst$ and $\Fm \dashv \Endm$ in \cref{ptmulticat-thm-v,ptmulticat-thm-ii}, respectively.
Note that these adjunctions are restricted to $\permcatst$ in \cref{ch:ptmulticat-sp}, but the natural domain of $\Endst$, and hence also $\Endm$, is $\permcatsu$ (see \cref{expl:endst-catmulti}).
The necessity of expanding to $\permcatsu$ for the codomain of $\Fst$ as a non-symmetric $\Cat$-multifunctor is described in \cref{expl:su-vs-st}.
Also recall \cref{remark:ptmulticat-prop-viii} for further technical remarks about how \cref{ptmulticat-prop-viii} is used below.
The pair $(\Mone \sma - ) \dashv \Um$ is described in \cref{MonesmaUmadj}.

\subsection*{Chapter Summary}

\Cref{sec:Sst-multilinear} defines the strong multilinear functors $\Fstn$ with their linearity constraints $(\Fstn)^2_p$.
\Cref{sec:Fstmultifunctor} uses $\Fstn$ to define $\Fst$ on multimorphism categories, showing that $\Fst$ is a non-symmetric $\Cat$-multifunctor.
\Cref{sec:etast-rhost} develops the (non-symmetric) $\Cat$-multinaturality of $\etast$ and $\vrhost$.
\Cref{sec:pmulticat-alg,sec:monemod-alg} develop the three inverse equivalences of homotopy theories shown in the diagram \cref{ptmulticatalg-summary} above.
Here is a summary table.
\reftable{.55}{
  $n$-linear functors $\brb{\Fstn,(\Fstn)^2_p}$
  & \ref{def:Sst}, \ref{def:Sst-constraints}
  \\ \hline
  non-symmetric $\Cat$ multifunctor $\Fst$
  & \ref{def:Fst-multi}, \ref{ptmulticat-xvii}
  \\ \hline \hline
  \multicolumn{2}{|c|}{non-symmetric $\Cat$-multinatural transformations}
  \\ \hline
   $\pst\cn \Fr\Ust \to \Fst$
  & \ref{pst-mnat}
  \\ \hline
  $\etast\cn 1_{\pMulticat} \to \Endst\Fst$
  & \ref{ptmulticat-xx}
  \\ \hline
  $\vrhost\cn 1_{\permcatsu} \to \Fst\Endst$
  & \ref{ptmulticat-xxii}
  \\ \hline \hline
  \multicolumn{2}{|c|}{inverse equivalences of homotopy theories}
  \\ \hline
  $\brb{\Fst^\Q,\Endst^\Q}$
  & \ref{ptmulticat-xxiii}
  \\ \hline
  $\brb{\Fm^\Q,\Endm^\Q}$
  & \ref{ptmulticat-xxv}
  \\ \hline
  $\brb{(\Mone\sma -)^\Q,\Um^\Q}$
  & \ref{Monesma-Um-algebra}
  \\ \hline
}
We remind the reader of \cref{conv:universe} about universes and \cref{expl:leftbracketing} about left normalized bracketing for iterated products.

\section{The Strong Multilinear Functor \texorpdfstring{$\Fstn$}{Fn•}}
\label{sec:Sst-multilinear}

This section defines pointed variants of the strong multilinear functors $\Frn$ from \cref{S-multifunctorial}.
These are determined by commutativity with the strict symmetric monoidal functors
\[
  \pst_\M \cn \Fr\M \to \Fst\M
\]
for small pointed multicategories $\M$, from \cref{def:FFst-projection}.

To begin, the definition of underlying functors $\Fstn$ depends on 
the tensor products of objects \cref{angxonen-def} and morphisms \cref{fonen-angphionen} defining $\Frn$.
Also recall the universal morphism $\vpi$ \cref{eq:multicat-smash-pushout} from a tensor product of small pointed multicategories to their corresponding smash product.
For the case $n=0$, recall from \cref{example:Fst-Mtu} that $\Fst(\Mtup) \iso \Fr(\Mtu)$ is isomorphic to the permutation category.

\begin{definition}[The Functors $\Fstn$]\label{def:Sst}
  Suppose $\angM = \ang{\M_i}_{i=1}^n$ are small pointed multicategories.  We define the data of a functor $\Fstn$ such that the following diagram commutes.
  \begin{equation}\label{eq:Frn-Fstn}
    \begin{tikzpicture}[x=40mm,y=13mm,vcenter]
      \draw[0cell=.9] 
      (0,0) node (a) {\txprod_{i=1}^n \Fr\M_i}
      (a)+(1,0) node (b) {\Fr\big( \txotimes_{i=1}^n \M_i \big)}
      (a)+(0,-2) node (c) {\txprod_{i=1}^n \Fst\M_i}
      (b)+(0,-1) node (d) {\Fr\big( \txsma_{i=1}^n \M_i \big)}
      (d)+(0,-1) node (e) {\Fst\big( \txsma_{i=1}^n \M_i \big)}
      ;
      \draw[1cell=.9] 
      (a) edge node {\Frn} (b)
      (b) edge node {\Fr \vpi} (d)
      (a) edge['] node {\txprod_i \pst_{\M_i}} (c)
      (c) edge[dashed] node {\Fstn} (e)
      (d) edge node {\pst_{(\sma_i \M_i)}} (e)
      ;
    \end{tikzpicture}
  \end{equation}

  If $n = 0$, then recall $\Mterm$ is the empty product in $\Multicat$ and $\Mtu_+$ is the empty smash product in $\pMulticat$.
  The functor\label{not:Fstzero}
  \[
    \Fst^0 \cn \boldone \to \Fst(\Mtup)
  \]
  is defined by the choice of object $[(1)] \in \Fst(\Mtup)$, where $(1)$ is the length-one tuple consisting of the unique object of the initial operad, $\Mtu$.

  Now suppose $n > 0$.
  \begin{description}
  \item[Assignment on Objects]
    Define $\Fstn$ on objects by
    \begin{equation}\label{eq:Sst-obj}
      \Fstn \big( [\ang{x^1}] , \ldots , [\ang{x^n}] \big)
      = [\ang{x^\onen}] \in \Fst\big( \txsma_{i=1}^n \M_i\big),
    \end{equation}
    using the tensor product of tuples 
    \begin{equation}\label{eq:xonen-st}
      \ang{x^{\onen}}
      = \lrang{ \,\cdots\, \lrang{x^{\onen}_{\jonejn}}_{j_1=1}^{r_1} \,\cdots\, }_{j_n = 1}^{r_n} 
    \end{equation}
    from \cref{angxonen-def}.
    This assignment is well defined on $\obsim$-equivalence classes because the term
    \[
      x^\onen_{\jonejn} = \ang{x^i_{j_i}}_{i = 1}^n \inspace \txotimes_{i=1}^n \M_i
    \]
    is sent to the basepoint object of $\txsma_{i} \M_i$ whenever any $x^i_{j_i} = * \in \M_i$.
    Thus, any basepoint terms of $\ang{x^i}$ produce corresponding basepoint terms of $\ang{x^\onen}$.
  \item[Assignment on Morphisms]
    For morphisms, suppose given length-one sequences
    \[
      \big[ \brb{f^i,\ang{\phi^i}} \big] \cn [\ang{x^i}] \to [\ang{y^i}] \inspace \Fst\M_i.
    \]
    Define
    \begin{equation}\label{eq:Sst-mor}
      \Fstn \big( [(f^1,\ang{\phi^1})] , \ldots , [(f^n,\ang{\phi^n})] \big)
      = [(f^\onen,\ang{\phi^\onen})],
    \end{equation}
    using the index map $f^\onen$ and the tensor product of morphisms
    \begin{equation}\label{eq:phionenkonekn-st}
      \ang{\phi^{\onen}}
      = \txotimes_{i=1}^n \ang{\phi^i}
    \end{equation}
    from \cref{fonen-angphionen}.
    This assignment is well defined on $\iisim$-equivalence classes because the operation
    \[
      \txotimes_{i=1}^n \phi^i_{k_i} \inspace \txotimes_{i=1}^n \M_i
    \]
    is sent to a basepoint operation of $\txsma_{i} \M_i$ whenever any $\phi^i_{k_i}$ is a basepoint operation in $\M_i$.

    Composition in $\Fst \M_i$ is given by concatenation of representative sequences, so the definition of $\Fstn$ on morphisms with representative sequences of length greater than one is determined by the desired functoriality of $\Fstn$.
    These assignments are well defined on $\isim$-equivalence classes by functoriality of $\Frn$, $\pst$, and $\Fst(\vpi)$.
  \end{description}

  The assignments above show that \cref{eq:Frn-Fstn} commutes on objects and morphisms.
  The definition of $\Fstn$ on morphisms implies that it is functorial.
\end{definition}

\begin{definition}[Linearity Constraints of $\Fstn$]\label{def:Sst-constraints}
  Suppose $\angM = \ang{\M_i}_{i=1}^n$ is a tuple of small pointed multicategories.
  For each $p \in \{1,\ldots,n\}$, we define the data of a natural transformation $(\Fstn)^2_p$ with components determined on representatives by $(\Frn)^2_p$.
  That is, we define
  \begin{equation}\label{eq:Fstn2p}
    (\Fstn)^2_p = \big[ (\Frn)^2_p \big] = 
    \big[ \brb{\rho_{r_p,\hat{r}_p}, \ang{1}} \big] \cn 
    [\ang{x^\onen}] \oplus [\ang{\hat{x}^\onen}] \fto{\iso} [\ang{\tilde{x}^\onen}]
  \end{equation}
  for
  \[\ang{\hat{x}^p} \in \Fr\M_p, \andspace
    \ang{x^\onen},\ \ang{\hat{x}^\onen},\ \ang{\tilde{x}^\onen}
    \in \Fr\big( \txotimes_{i=1}^n \M_i \big)
  \]
  as in \cref{eq:xhatonen-xtilonen}.
  The components \cref{eq:Fstn2p} are well defined because the operations determining $(\Frn)^2_p$ are colored units and, therefore, $\obsim$-equivalent representatives in the domain or codomain of \cref{eq:Fstn2p} will result in $\iisim$-equivalent components $(\Frn)^2_p$.

  Naturality of \cref{eq:Fstn2p} follows, for length-one morphism sequences, from that of $(\Frn)^2_p$.
  This implies naturality with respect to general $\obsim$-composable morphism sequences because composition in $\Fst\big( \txsma_{i=1}^n \M_i \big)$ is given by concatenation.
\end{definition}

Now we show that the constructions above assemble to form multilinear functors (\cref{def:nlinearfunctor}).

\begin{proposition}\label{ptmulticat-xv}
  Suppose $\angM = \ang{\M_i}_{i=1}^n$ is a tuple of small pointed multicategories for $n \ge 0$.
  The data in \cref{def:Sst,def:Sst-constraints}
  \begin{equation}\label{Sst-nlinear}
    \left(\Fstn \scs \bang{(\Fstn)^2_p}_{p=1}^n \right) \cn 
    \txprod_{i=1}^n \Fst\M_i \to \Fst\big( \txsma_{i=1}^n \M_i\big)
  \end{equation}
  form a strong $n$-linear functor.
\end{proposition}
\begin{proof}
  Each of the multilinearity axioms in \cref{def:nlinearfunctor} follows from the corresponding axiom for $\left(\Frn \scs \bang{(\Frn)^2_p}_{p=1}^n \right)$ on representatives.
  The components of each $(\Fstn)^2_p$ are isomorphisms because their representatives $(\Frn)^2_p$ are so.
\end{proof}

\begin{proposition}\label{ptmulticat-xvi}
  For $n \ge 0$, the strong $n$-linear functor in \cref{Sst-nlinear},
  \[
    \left(\Fstn \scs \bang{(\Fstn)^2_p}_{p=1}^n \right) \cn 
    \txprod_{i=1}^n \Fst\M_i \to \Fst\big( \txsma_{i=1}^n \M_i\big),
  \]
  is 2-natural with respect to pointed multifunctors and pointed multinatural transformations.
\end{proposition}
\begin{proof}
  Suppose given a tuple $\ang{H} = \ang{H_i}_{i=1}^n$ of pointed multifunctors
  \[
    H_i \cn \M_i \to \N_i \forspace 1 \leq i \leq n.
  \]
  Naturality of $\Fstn$ with respect to $\ang{H_i}$ is verified by commutativity of the inner rectangle in the following diagram.
  \begin{equation}\label{eq:nat-Fstn}
    \begin{tikzpicture}[x=35mm,y=33mm,vcenter]
      \draw[0cell=.9] 
      (0,0) node (a) {\txprod_{i=1}^n \Fst\M_i}
      (a)+(1,0) node (b) {\txprod_{i=1}^n \Fst\N_i}
      (a)+(0,-.5) node (c) {\Fst\big( \txsma_{i=1}^n \M_i \big)}
      (b)+(0,-.5) node (d) {\Fst\big( \txsma_{i=1}^n \N_i \big)}
      (a)+(-.707,.4) node (x) {\txprod_{i=1}^n \Fr\M_i}
      (b)+(.707,.4) node (y) {\txprod_{i=1}^n \Fr\N_i}
      (c)+(225:1) node (z) {\Fr\big( \txotimes_{i=1}^n \M_i \big)}
      (d)+(315:1) node (w) {\Fr\big( \txotimes_{i=1}^n \N_i \big)}
      (c)+(225:.5) node (z') {\Fr\big( \txsma_{i=1}^n \M_i \big)}
      (d)+(315:.5) node (w') {\Fr\big( \txsma_{i=1}^n \N_i \big)}
      ;
      \draw[1cell=.9] 
      (x) edge['] node {\txprod_i \pst_{\M_i}} (a)
      (y) edge node {\txprod_i \pst_{\N_i}} (b)
      (z) edge node {\Fr\vpi} (z')
      (z') edge node {\pst} (c)
      (w) edge['] node {\Fr\vpi} (w')
      (w') edge['] node {\pst} (d)
      (x) edge node {\txprod_i \Fr H_i} (y)
      (z) edge node {\Fr (\otimes_i H_i)} (w)
      (z') edge node {\Fr (\sma_i H_i)} (w')
      (x) edge['] node {\Frn} (z)
      (y) edge node {\Frn} (w)
      (a) edge node {\txprod_i \Fst H_i} (b)
      (c) edge node {\Fst (\sma_i H_i)} (d)
      (a) edge['] node {\Fstn} (c)
      (b) edge node {\Fstn} (d)
      ;
    \end{tikzpicture}
  \end{equation}
  The functors $\pst$ are surjective on objects and on length-one tuples of morphisms.
  Since morphisms in each $\Fst\M_i$ are generated under composition (concatenation) by those of length one, it suffices to verify
  \[
    \Fstn \circ \big( \txprod_i \Fst H_i \big) \circ \big( \txprod_i \pst_{\M_i} \big) 
    = 
    \Fst \big( \txsma_i H_i \big) \circ \Fstn \circ \big( \txprod_i \pst_{\M_i} \big).
  \]
  The equality above holds by commutativity of the following in \cref{eq:nat-Fstn}.
  \begin{itemize}
  \item The trapezoids at left and right commute by \cref{eq:Frn-Fstn}.
  \item The bottom trapezoid commutes by naturality of $\vpi$ \cref{eq:multicat-smash-pushout}.
  \item The remaining two trapezoid regions commute by naturality of $\pst$ \cref{eq:pst-iinatural}.
  \item The outer rectangle commutes by naturality of $\Frn$.
  \end{itemize}

  This shows that $\Fstn$ is natural with respect to $\ang{H}$.
  A similar analysis for pointed multinatural transformations shows that $\Fstn$ is 2-natural.
\end{proof}

\section[Non-Symmetric \texorpdfstring{$\Cat$}{Cat}-Multifunctor]{Pointed Free Permutative Category as a Non-Symmetric \texorpdfstring{$\Cat$}{Cat}-Multifunctor}
\label{sec:Fstmultifunctor}

In this section we extend the 2-functor 
\[\Fst \cn \pMulticat \to \permcatst\] 
in \cref{ptmulticat-thm-i} to multimorphism categories.

\begin{convention}\label{conv:Fbst}
To avoid confusion in \cref{def:Fst-multi} below, for small pointed multicategories $\M$ and $\N$, we denote by
\[\Fbst \cn \pMulticat(\M,\N) \to \permcatst(\Fst\M,\Fst\N)\] 
the assignment of $\Fst$ on pointed multifunctors and pointed multinatural transformations as in \cref{def:Fst-onecells,def:Fst-twocells}, respectively.
This is the pointed analog of \cref{convention:Fun} for $\Fr$.
\end{convention}

In \cref{Sst-Fbst} below, we use the multilinear functor $\Fstn$ (\cref{ptmulticat-xv}).

\begin{definition}\label{def:Fst-multi}
Suppose $\angM = \ang{\M_i}_{i=1}^n$ and $\N$ are small pointed multicategories.  We define a functor between multimorphism categories
\begin{equation}\label{Fst-multimorphism-cat}
\Fst\cn \pMulticat\scmap{\ang{\M};\N} \to \permcatsu\scmap{\ang{\Fst\M};\Fst\N}
\end{equation}
as follows.  Suppose given pointed multifunctors $H$ and $K$ and a pointed multinatural transformation $\theta$ as in the diagram below.
\[\begin{tikzpicture}[baseline={(a.base)}]
\def\t{22}
\draw[0cell=.9]
(0,0) node (a) {\ang{\M}}
(a)++(2,0) node (b) {\N}
;
\draw[1cell=.85]  
(a) edge[bend left=\t] node[pos=.43] {H} (b)
(a) edge[bend right=\t] node[swap,pos=.43] {K} (b)
;
\draw[2cell=.9] 
node[between=a and b at .47, rotate=-90, 2label={above,\theta}] {\Rightarrow}
;
\end{tikzpicture}\]
Then $\Fst$ sends these data to the following composite $n$-linear functors and whiskering.
\begin{equation}\label{Sst-Fbst}
\begin{tikzpicture}[baseline={(a.base)}]
\def\t{20}
\draw[0cell=.9]
(-2.75,0) node (z) {\ang{\Fst\M}}
(0,0) node (a') {\Fst\big(\txsma_{i=1}^n \M_i\big)}
(a')+(.65,0) node (a) {\phantom{\Fr\N}}
(a)++(2.5,0) node (b) {\Fst\N}
;
\draw[1cell=.85]  
(z) edge node {\Fstn} (a')
(a) edge[bend left=\t] node[pos=.5] {\Fbst H} (b)
(a) edge[bend right=\t] node[swap,pos=.5] {\Fbst K} (b)
;
\draw[2cell=.9] 
node[between=a and b at .39, rotate=-90, 2label={above,\,\Fbst \theta}] {\Rightarrow}
;
\end{tikzpicture}
\end{equation}
This finishes the definition of the multimorphism functor $\Fst$.
\end{definition}

\begin{explanation}[Codomain Not Strict]\label{expl:su-vs-st}
In \cref{Fst-multimorphism-cat} above and \cref{ptmulticat-xvii} below, we stress that the codomain uses $\permcatsu$ and \emph{not} $\permcatst$, which is the codomain of $\Fbst$ in \cref{conv:Fbst}.  The reason is that the definition of $\Fst$ in \cref{Sst-Fbst} involves the strong $n$-linear functor $\Fstn$ in \cref{Sst-nlinear}.  The latter is \emph{not} strict because the components of its linearity constraints $(\Fstn)^2_p$ in \cref{def:Sst-constraints} are not identities in general.
\end{explanation}

\begin{theorem}\label{ptmulticat-xvii}\index{category!free permutative - multifunctor!pointed}\index{permutative category!free - multifunctor!pointed}\index{free!permutative category multifunctor!pointed}
There is a non-symmetric $\Cat$-multifunctor
\begin{equation}\label{eq:ptmulticat-xvii}
  \Fst \cn \pMulticat \to \permcatsu
\end{equation}
with
\begin{itemize}
\item object assignment in \cref{def:Fst-permutative} and
\item multimorphism functors in \cref{def:Fst-multi}.
\end{itemize}
Moreover, $\Fst$ extends the 2-functor in \cref{ptmulticat-thm-i}.
\end{theorem}
\begin{proof}
  To see that $\Fst$ preserves units, note that $\Fst^1$ is the identity monoidal functor.
  Since $\Fbst$ is functorial, we have
  \[
    \Fst(1_\M) = 1_{\Fst \M}
  \]
  for each small pointed multicategory $\M$.
  
  To see that $\Fst$ preserves composition, suppose given
  \begin{align*}
    H_i & \in \pMulticat\scmap{\ang{\M_i};\M'_i} \forspace 1 \leq i \leq n, \andspace\\
    H' & \in \pMulticat\scmap{\ang{\M'};\M''}.
  \end{align*}
  Let $k_i$ denote the length of $\ang{\M_i}$.
  The two multilinear functors
  \[
    \Fst\big( \ga\scmap{H';\ang{H}} \big)
    \andspace
    \ga\scmap{\Fst H'; \ang{\Fst H}}
  \]
  are given by the two composites around the boundary in the following diagram, where the unlabeled isomorphisms are given by reordering terms.
  \[
    \begin{tikzpicture}[x=27mm,y=18mm]
      \draw[0cell=.9] 
      (0,0) node (t1) {\txprod_{i,j} \Fst\M_{i,j}}
      (t1)++(0,-2) node (t2) {\Fst\big(\txsma_{i,j} \M_{i,j}\big)}
      (t2)++(1,0) node (t3) {\Fst\big(\txsma_{i} \txsma_{j} \M_{i,j}\big)}
      (t3)++(1.3,0) node (t4) {\Fst\big(\txsma_{i} \M'_i \big)}
      (t4)++(1,0) node (t5) {\Fst \M''}
      (t1)++(1,0) node (b1) {\txprod_i \txprod_j \Fst\M_{i,j}}
      (b1)++(0,-1) node (b2) {\txprod_i \Fst\big(\txsma_j \M_{i,j}\big)}
      (b2)++(1.3,0) node (b3) {\txprod_i \Fst \M'_{i}}
      ;
      \draw[1cell=.9] 
      (t1) edge['] node {\Fst^{k_1+\cdots+k_n}} (t2)
      (t2) edge node {\iso} (t3)
      (t3) edge node {\Fbst(\sma_i H_i)} (t4)
      (t4) edge node {\Fbst H'} (t5)
      (t1) edge node {\iso} (b1)
      (b1) edge node {\binprod_i \Fst^{k_i}} (b2)
      (b2) edge node {\binprod_i \Fbst(H_i)} (b3)
      (b3) edge node {\Fstn} (t4)
      (b2) edge node {\Fstn} (t3)
      ;
    \end{tikzpicture}
  \]
  
  In the above diagram, the two composites around the middle rectangle are equal as multilinear functors by naturality of $\Fstn$ (\cref{ptmulticat-xvi}) with respect to the multifunctors $H_i$.
  The rectangle at left commutes on objects and length-one sequences of morphisms because the corresponding diagram for $\Fr$ and $(\Frn,(\Frn)^2_p)$ commutes (\cref{theorem:F-multi} c.f., proof of \cite[8.1]{johnson-yau-Fmulti}).
  The commutativity for general morphisms then follows from functoriality of $\Fstn$.
  A similar diagram for multinatural transformations $H_i \to K_i$ and $H'_i \to K'_i$ commutes by the 2-naturality of $\Fstn$.
\end{proof}

Recall (non-symmetric) $\Cat$-multinatural transformation from \cref{expl:catmultitransformation}.
\begin{lemma}\label{pst-mnat}
  The 2-natural transformation $\pst$ of \cref{pst-iinatural} extends to a non-symmetric $\Cat$-multinatural transformation
  \[
    \begin{tikzpicture}
      \def\h{2} \def\t{23}
      \draw[0cell]
      (0,0) node (a) {\phantom{X}}
      (a)+(\h,0) node (b) {\phantom{X}}
      (a)+(-.55,0) node (a') {\pMulticat}
      (b)+(.6,.06) node (b') {\permcatsu.}
      ;
      \draw[1cell=.9]
      (a) edge[bend left=\t,transform canvas={yshift=0ex}] node {\Fr \circ \Ust} (b)
      (a) edge[bend right=\t,transform canvas={yshift=0ex}] node[swap] {\Fst} (b)
      ;
      \draw[2cell]
      node[between=a and b at .45, rotate=-90, 2label={above,\,\pst}] {\Rightarrow}
      ;
    \end{tikzpicture}
  \]
\end{lemma}
\begin{proof}
  Recall from \cref{expl:Ust-catmulti} that the forgetful
  \[
    \Ust \cn \pMulticat \to \Multicat
  \]
  is a $\Cat$-multifunctor with $n$-ary multimorphism functors
  \[
    \pMulticat\scmap{\ang{\M};\N} \to \Multicat\scmap{\Ust\ang{\M};\Ust\N}
  \]
  given by composition and whiskering with $\varpi\cn \otimes_i \M_i \to \sma_i \M_i$ \cref{Ust-monconstraint}.
  The $\Cat$-multifunctoriality of $\Fr$ is described in \cref{theorem:F-multi} and its multimorphism functors are given by precomposition and whiskering with the $n$-linear functors $\Frn$ from \cref{def:S-multi}.

  For objects of $\pMulticat\scmap{\ang{\M};\N}$, i.e., multifunctors
  \[
    H\cn \txsma_i \M_i \to \N,
  \]
  the object $\Cat$-naturality diagram \cref{catmultinaturalitydiagram} for $\pst$ is the following.
  In this diagram, we suppress $\Ust$ except in the lower left arrow.
  \[
    \begin{tikzpicture}[x=27mm,y=13mm]
      \draw[0cell=.9] 
      (0,0) node (a) {\txprod_i \Fr\M_i}
      (a)+(0,-1) node (b) {\Fr(\txotimes_i \M_i)}
      (b)+(0,-1) node (c) {\Fr\N}
      (a)+(2,0) node (a') {\txprod_i \Fst\M_i}
      (a')+(0,-1) node (b') {\Fst(\txsma_i \M_i)}
      (b')+(0,-1) node (c') {\Fst\N}
      (b)+(1,0) node (x) {\Fr(\txsma_i \M_i)}
      ;
      \draw[1cell=.9] 
      (a) edge['] node {\Frn} (b)
      (b) edge['] node {\Frbar (\Ust H)} (c)
      (a') edge node {\Fstn} (b')
      (b') edge node {\Fbst H} (c')
      (a) edge node {\txprod_i \pst_{\M_i}} (a')
      (b) edge node {\Frbar\varpi} (x)
      (x) edge node {\pst_{\sma_i \M_i}} (b')
      (c) edge node {\pst_{\N}} (c')
      (x) edge['] node {\Frbar H} (c)
      ;
    \end{tikzpicture}
  \]

  In the above diagram, the top rectangle commutes by construction \cref{eq:Frn-Fstn}, the lower left triangle commutes by definition of $\Ust H$, and the lower right trapezoid commutes by naturality of $\pst$.
  Therefore, the object $\Cat$-naturality condition \cref{catmultinaturalitydiagram} holds for each $H \in \pMulticat\scmap{\ang{\M};\N}$.
  For morphisms, i.e., multinatural transformations
  \[
    \theta \cn H \to H' \cn \txsma_i \M_i \to \N,
  \]
  a similar argument applies to verify the morphism $\Cat$-naturality condition \cref{catmultinatiicellpasting}, using 2-naturality of $\pst$.
\end{proof}

\begin{example}[Partition Products]\label{example:FstPiab}
  Recall from \cref{definition:part-prod} the \index{partition product}\index{product!partition}partition product multifunctor
  \[
    \Pi_{a,b}\cn \cM a \sma \cM b \to \cM(a \sma b)
  \]
  for finite pointed sets $a$ and $b$.
  It is given on objects by the Cartesian product, and preserves operations because a pairwise Cartesian product of partitions $\ang{s}$ and $\ang{t}$ provides a partition of $(\cup_i s_i) \times (\cup_j t_j)$.
  Using the description from \cref{example:FstMa}, and omitting the forgetful $\Um$, the composite
  \[
    \Fst(\cM a) \times \Fst(\cM b)
    \fto{\Fst^2} \Fst(\cM a \sma \cM b)
    \fto{\Fst(\Pi_{a,b})} \Fst(\cM(a \sma b))
  \]
  is given as follows for two tuples of nonempty subsets $s^1_i \subset a^\punc$ and $s^2_j \subset b^\punc$:
  \[
    [\ang{s^1}],[\ang{s^2}] \mapsto [\ang{s^{12}}] \mapsto [\ang{\ang{s^1_i \times s^2_j}_i}_j].\dqed
  \]
\end{example}

\section{Comparison Transformations}
\label{sec:etast-rhost}

Consider the diagram
\begin{equation}\label{FstEst-context}
\Fst \cn \pMulticat \lradj \permcatsu \cn \Endst
\end{equation}
consisting of the following data.
\begin{itemize}
\item $\pMulticat$ is the $\Cat$-multicategory in \cref{expl:ptmulticatcatmulticat}.
\item $\permcatsu$ is the $\Cat$-multicategory in \cref{thm:permcatmulticat}.
\item $\Endst$ is the $\Cat$-multifunctor in \cref{expl:endst-catmulti}.
\item $\Fst$ is the non-symmetric $\Cat$-multifunctor in \cref{ptmulticat-xvii}.
\end{itemize}
Recall from \cref{ptmulticat-thm-v} that the 2-adjunction
\[\begin{tikzpicture}
\def\t{17}
\draw[0cell]
(0,0) node (a) {\phantom{X}}
(a)+(-.55,0) node (a') {\pMulticat}
(a)+(1,0) node (x) {\bot}
(a)+(2,0) node (b) {\phantom{X}}
(b)+(.65,.06) node (b') {\permcatst}
;
\draw[1cell=.9]
(a) edge[bend left=\t] node {\Fst} (b)
(b) edge[bend left=\t] node {\Endst} (a)
;
\end{tikzpicture}\]
has unit
\[\etast \cn 1_{\pMulticat} \to \Endst\Fst\]
in \cref{def:etast}.  Recall \emph{(non-symmetric) $\Cat$-multinatural transformation} from \cref{expl:catmultitransformation}.

\begin{lemma}\label{ptmulticat-xx}
  In the context of \cref{FstEst-context}, 
  \[
    \etast \cn 1_{\pMulticat} \to \Endst\Fst
  \]
  is a non-symmetric $\Cat$-multinatural transformation.
\end{lemma}
\begin{proof}
  In this proof we will omit the adjective non-symmetric for the $\Cat$-multinatural transformations under discussion.
  Recall from \cref{etast-def} that $\etast$ is the composite of $\Endst(\pst) = \Endst * \pst$ with $\eta$.
  Thus, $\Cat$-multinaturality of $\etast$ follows from that of $\eta$ (\cref{lemma:eta-mnat}) and that of $\pst * \Endst$ (\cref{pst-mnat,expl:endst-catmulti}).
  Indeed, as a 2-natural transformation, $\etast$ satisfies the following equality of pasting diagrams.
  \[
    \begin{tikzpicture}[x=22mm,y=19mm]
      \def\wl{2} 
      \def\wr{0} 
      \def\h{.5} 
      \def\m{.5} 
      \draw[font=\Large] (0,0) node (eq) {=}; 
      \newcommand{\boundary}{
        \draw[0cell=.8] 
        (0,0) node (a) {\pMulticat}
        (a)+(1,-1) node (b) {\permcatsu}
        (b)+(1,.2) node (a') {\pMulticat}
        (b)+(1,1) node (c') {\Multicat}
        (a)+(1,.3) node (c'') {\pMulticat}
        ; 
        \draw[1cell=.8] 
        (a) edge['] node (Fst) {\Fst} (b)
        (b) edge['] node (Est) {\Endst} (a')
        (a') edge['] node {\Ust} (c')
        (a) edge node {1} (c'')
        (c'') edge node {\Ust} (c')
        ;
      }
      \begin{scope}[shift={(-\wl-\m,\h)}]
        \boundary
        \draw[0cell=.8]
        ;
        \draw[1cell=.8] 
        (a) edge node (V) {1} (a')
        ;
        \draw[2cell] 
        node[between=V and b at .5, rotate=-90, 2label={above,\etast}] {\Rightarrow}
        ;
      \end{scope}
      \begin{scope}[shift={(\wr+\m,\h)}]
        \boundary
        \draw[0cell=.8]
        (a)+(1,-.3) node (c) {\Multicat}
        ;
        \draw[1cell=.8] 
        (a) edge node {\Ust} (c)
        (c) edge node {\Fr} (b)
        (c) edge node (W) {1} (c')
        (b) edge['] node (E) {\End} (c')
        ;
        \draw[2cell] 
        node[between=c and Fst at .5, rotate=-90, 2label={above,\pst}] {\Rightarrow}
        node[between=c and E at .5, rotate=-90, 2label={above,\eta}] {\Rightarrow}
        ;
      \end{scope}
    \end{tikzpicture}
  \]
  The right hand side is also a diagram of $\Cat$-multinatural transformations, and thus $\Ust * \etast$ is $\Cat$-multinatural.
  Now since
  \begin{romenumerate}
  \item the two conditions for $\Cat$-multinaturality in \cref{expl:catmultitransformation} consist of certain equalities \cref{catmultinaturalitydiagram,catmultinatiicellpasting} involving composition of operations, and
  \item and such equalities are detected on underlying multicategories, by applying $\Ust$,
  \end{romenumerate}
  we conclude that $\etast$ also satisfies the conditions for $\Cat$-multinaturality.
\end{proof}

Recall from \cref{ptmulticat-xxi} the strictly unital symmetric monoidal functor
\[\vrhost_\C \cn \C \to \Fst\, \Endst\C\]
for each small permutative category $\C$.

\begin{lemma}\label{ptmulticat-xxii}
In the context of \cref{FstEst-context}, 
\[\vrhost \cn 1_{\permcatsu} \to \Fst\, \Endst\]
is a non-symmetric $\Cat$-multinatural transformation.
\end{lemma}
\begin{proof}
  For objects of $\permcatsu\scmap{\ang{\C};\D}$, i.e., multilinear functors
  \[
    Q\cn \txprod_{i = 1}^n \C_i \to \D,
  \]
  the object $\Cat$-naturality diagram \cref{catmultinaturalitydiagram} for $\vrhost$ is the following, where we use the abbreviation $\Est = \Endst$.
  \begin{equation}\label{eq:vrhost-mfun}
    \begin{tikzpicture}[x=35mm,y=13mm,vcenter]
      \draw[0cell=.9] 
      (0,0) node (a) {\txprod_i \C_i}
      (a)+(1,0) node (b) {\txprod_i\Fr \Est \C_i}
      (b)+(1,0) node (c) {\txprod_i \Fst \Est \C_i}
      (a)+(0,-2) node (a') {\D}
      (a')+(1,0) node (b') {\Fr\Est \D}
      (b')+(1,0) node (c') {\Fst\Est \D}
      (b)+(0,-1) node (x) {\Fr(\txotimes_i \Est \C_i)}
      (c)+(0,-1) node (y) {\Fst(\txsma_i \Est \C_i)}
      ;
      \draw[1cell=.9] 
      (a) edge node {\txprod_i \vrho_{\C_i}} (b)
      (b) edge node {\txprod_i \pst_{\Est\C_i}} (c)
      (a') edge node {\vrho_\D} (b')
      (b') edge node {\pst_{\Est\D}} (c')
      (a) edge['] node {Q} (a')
      (b) edge node {\Frn} (x)
      (x) edge node {\Fr(\Est Q)} (b')
      (c) edge node {\Fstn} (y)
      (y) edge node {\Fst(\Est Q)} (c')
      ;
    \end{tikzpicture}
  \end{equation}
  In the above diagram, the rectangle at right commutes by multifunctoriality of $\pst$ (\cref{pst-mnat}).
  To verify that the rectangle at left commutes, we first consider the underlying functors and then check linearity constraints.
  
  For morphisms
  \[
    \phi^i \cn b^i \to c^i \inspace \C_i,
  \]
  the composite to $\Fr\Est\D$ along the top and vertical middle functors of \cref{eq:vrhost-mfun} is given as follows:
  \begin{align*}
    (\phi^1, \cdots, \phi^n)
    & \mapsto
      \Big(\, \big(1_{\ord{1}},(\phi^1)\big), \ldots \big(1_{\ord{1}},(\phi^n)\big) \,\Big)
    & & \text{by \cref{rhoc-assignments} for $\vrho$,} 
    \\
    & \mapsto
      \Big(\, 1_{\ord{1}} \scs \big(\phi^\onen\big) \,\Big)
    & & \text{by \cref{eq:Sfi} for $\Frn$,} 
    \\
    & \mapsto
      \Big(\, 1_{\ord{1}} \scs \big(H\phi^\onen\big) \,\Big)
    & & \text{by \cref{eq:FH-fphi} for $H = \Est Q$,} 
  \end{align*}
  where $\big( \phi^\onen \big)$ is the length-one tuple consisting of the morphism $\phi^\onen = \otimes_i \phi^i$ in $\otimes_i \C_i$.
  By \cref{eq:EstP1compjpsi} we have
  \[
    H\phi^\onen = (\Est Q)(\phi^1,\ldots,\phi^n) = Q(\phi^1,\ldots,\phi^n);
  \]
  the linear constraints $Q^2$ do not appear because each $\phi^i$ is a unary operation in $\Est\C$.
  This verifies that the left side of \cref{eq:vrhost-mfun} commutes as a diagram of underlying functors.

  Next we consider the $j$th linearity constraints on the left side of \cref{eq:vrhost-mfun}, for $1 \leq j \leq n$.  
  Suppose given
  \[
    c^i \in \C_i \foreachspace 1 \leq i \leq n,  \andspace \hat{c}^j \in \C_j.
  \]
  Recall the notation $\circ_j$ from \cref{notation:compk}.
  We will use the following:
  \begin{align*}
    \tilde{c}^j & = c^j \oplus \hat{c}^j \inspace \C_j,\\
    \{ c \} & = (c^1,\ldots,c^n) \inspace \txprod_i \C_i,\\
    \{ \hat{c} \} & = \{ c \circ_j \hat{c}^j \} = (c^1,\ldots,\hat{c}^j,\ldots,c^n) \inspace \txprod_i \C_i, \andspace\\
    \{ \tilde{c} \} & = \{ c \circ_j \tilde{c}^j \} = (c^1,\ldots,\tilde{c}^j,\ldots,c^n) \inspace \txprod_i \C_i.
  \end{align*}
  So each of $\{ \hat{c} \}$ and $\{ \tilde{c} \}$ has the same entries as $\{ c \}$ in all but the $j$th position.
  Next we introduce further notation:
  \begin{align*}
    \txprod\vrho & = \txprod_i \vrho_{\C_i},\\
    (c^\onen) & = \txotimes_i (c^i) = \Frn\big( \txprod\vrho \{ c \} \big)
                \inspace \Fr\big(\txotimes_i \Est \C_i),\\
    (\hat{c}^\onen) & = \Frn\big( \txprod\vrho \{ \hat{c} \} \big)
                      \inspace \Fr\big(\txotimes_i \Est \C_i),\\
    (\tilde{c}^\onen) & = \Frn\big( \txprod\vrho \{ \tilde{c} \} \big)
                        \inspace \Fr\big(\txotimes_i \Est \C_i), \andspace\\
    \ol{c} & = \txprod\vrho\{ c \} \circ_j (c^j,\hat{c}^j)
             \inspace \txprod_i \Fr\Est \C_i.
  \end{align*}
  So each of $(c^\onen)$, $(\hat{c}^\onen)$, and $(\tilde{c}^\onen)$ is a length-one tuple and
  \[
    \ol{c} = \big( (c^1), \ldots, (c^j,\hat{c}^j), \ldots, (c^n) \big)
  \]
  has the same entries as $\txprod\vrho\{c\}$ except entry $j$, which is $\txprod\vrho(c) \oplus \txprod\vrho(\hat{c}^j) = (c^j,\hat{c}^j)$.
  With this notation, note that $\Frn\ol{c}$ is the length-two tuple $(c^\onen,\hat{c}^\onen)$.

  Applying \cref{ffjlinearity,eq:S2b}, the $j$th linearity constraint of $\Frn \circ \txprod\vrho$ at the objects $\{ c \}$ and $\{ \hat{c} \}$ is the composite in $\Frn(\otimes_i \Est \C_i)$
  \begin{equation}\label{eq:Frntxprodvrho2j}
    (c^\onen,\hat{c}^\onen) \fto{(\Frn)^2_j = 1}
    \Frn\ol{c} \fto{\Frn\big(\ang{1 \circ_j \vrho^2_{\C_j}} \big)}
    (\tilde{c}^\onen),
  \end{equation}
  where we note
  \begin{align*}
    \Frn(\txprod\vrho\{c\}) \oplus \Frn(\txprod\vrho\{\hat{c}\})
    & = (c^\onen,\hat{c}^\onen) \andspace \\
    \Frn(\txprod\vrho \{ \tilde{c} \})
    & = (\tilde{c}^\onen).
  \end{align*}

  Recall \cref{rhoc-monconstraint} for the monoidal constraint $\rho^2_{\C_j}$. 
  Since $\Fr(\Est Q)$ is strict symmetric monoidal, the $j$th linearity constraint of the top and vertical middle composite of \cref{eq:vrhost-mfun} is given by applying $\Fr(\Est Q)$ to \cref{eq:Frntxprodvrho2j}, resulting in the following:
  \[
    Q(c^\onen) \oplus Q(\hat{c}^\onen)
    \fto{1}
    Q(c^\onen) \oplus Q(\hat{c}^\onen)
    \fto{\brb{\iota_2,(Q^2_j)}}
    Q(\tilde{c}^\onen).
  \]
  Now since
  \[
    \brb{\iota_2,(Q^2_j)} = \vrho_\D\big(Q^2_j),
  \]
  we conclude that the $j$th linearity constraints of the two composites in the left half of \cref{eq:vrhost-mfun} are equal.
  This completes verification of the object $\Cat$-naturality condition \cref{catmultinaturality} for $\vrho$.

  Verification of the morphism $\Cat$-naturality condition \cref{catmultinaturalityiicell} is similar.
  Given a multilinear transformation
  \[
    \theta \cn Q \to Q',
  \]
  one verifies
  \[
    \Fst(\Est \theta) * (\Fstn \circ \txprod_i \pst_{\Est \C_i} \circ \txprod_i \vrho_{\C_i})
    = 
    \pst_{\Est \D} * \Fr(\Est \theta) * (\Frn \circ \txprod_i \vrho_{\C_i})
    =
    (\pst_{\Est \D} \circ \vrho_{\D}) * \theta,
  \]
  using \cref{pst-mnat} for the first equality and \cref{eq:Endst-theta-component} with \cref{eq:Fka-x} for the second.
  This completes the proof.
\end{proof}

\section{Multiplicative Homotopy Theory of Pointed Multicategories}
\label{sec:pmulticat-alg}

Recall from \cref{def:ptmulti-stableeq} the wide subcategory
\[
  \cSst \bigsubset \pMulticat
\]
of $\Fst$-stable equivalences created by 
\[
  \Fst \cn \brb{\pMulticat,\cSst} \to \brb{\permcatst,\cSI}.
\]  
Recall from \cref{def:nonsymalgebra} the notion of \emph{non-symmetric algebras}.

\cref{ptmulticat-xxiii} below considers the $\Cat$-multifunctors (non-symmetric for $\Fst$)
\[\Fst \cn \pMulticat \lradj \permcatsu \cn \Endst\]
in \cref{FstEst-context}.
This result simultaneously extends
\begin{itemize}
\item \cref{thm:alg-hty-equiv} from $\Multicat$ to $\pMulticat$ and
\item \cref{ptmulticat-thm-x} to non-symmetric algebras.
\end{itemize} 
As in \cref{thm:alg-hty-equiv}, each of the two induced functors $\Fst^\Q$ and $\Endst^\Q$ is given by post-composition and whiskering with the respective functor.

\begin{theorem}\label{ptmulticat-xxiii}
Suppose $\Q$ is a small non-symmetric $\Cat$-multicategory.  In the context of \cref{FstEst-context}, the induced functors  
\[\Fst^\Q \cn \brb{\pMulticat^\Q, \cSst^\Q} \lrsimadj \brb{(\permcatsu)^\Q, \cS^\Q} 
\cn \Endst^\Q\]
are inverse equivalences of homotopy theories in the sense of \cref{def:inverse-heq}.
\end{theorem}
\begin{proof}
  We first verify that $\Fst^\Q$ and $\Endst^\Q$ are relative functors.
  The functor $\Fst^\Q$ is a relative functor because the stable equivalences $\cSst^\Q$ are determined componentwise and $\Fst$ creates the stable equivalences $\cSst$.

  In \cref{ptmulticat-thm-x}, $\Endst$ is shown to be a relative functor with respect to \emph{strict} symmetric monoidal functors.
  For more general strictly unital symmetric monoidal functors, recall from \cref{ptmulticat-xxi} the strictly unital
  \[
    \vrhost_\C \cn \C \to \Fst\Endst\C \forspace \C \in \permcatst.
  \]
  \Cref{ptmulticat-prop-viii} shows that each $\vrhost_\C$ is a right adjoint and, hence, a stable equivalence of permutative categories.
  \Cref{ptmulticat-xxii} shows that $\vrhost$ is $\Cat$-multinatural and, in particular, 2-natural with respect to strictly unital symmetric monoidal functors.
  Therefore, if $P \cn \C \to \D$ is a stable equivalence in $\permcatsu$, then the naturality diagram
  \begin{equation}\label{eq:vrhost-nat-steq}
    \begin{tikzpicture}[x=27mm,y=13mm,vcenter]
      \draw[0cell=.9] 
      (0,0) node (c) {\C}
      (1,0) node (fec) {\Fst\Endst\C}
      (0,-1) node (d) {\D}
      (1,-1) node (fed) {\Fst\Endst\D}
      ;
      \draw[1cell=.9] 
      (c) edge node {\vrhost_\C} node['] {\sim} (fec)
      (d) edge node {\vrhost_\D} node['] {\sim} (fed)
      (c) edge['] node {P} node['] {\sim} (d)
      (fec) edge node {\Fst\Endst P} (fed)
      ;
    \end{tikzpicture}
  \end{equation}
  shows that $\Fst\Endst P$ is a stable equivalence by the 2-out-of-3 property.
  Since $\Fst$ creates stable equivalences, this implies that $\Endst P$ is a stable equivalence.
  Therefore, $\Endst$ is a relative functor with respect to the strictly unital stable equivalences.
  Thus it follows that $\Endst^\Q$ is also a relative functor with respect to $\cS^\Q$.

  Since $\vrhost$ is 2-natural and componentwise a stable equivalence, the induced $(\vrhost)^\Q$ is a natural stable equivalence
  \[
    1_{(\permcatsu)^\Q} \fto[\sim]{(\vrhost)^\Q} \Fst^\Q \Endst^\Q.
  \]
  For the reverse composite, recall from \cref{ptmulticat-xx} that $\etast$ is $\Cat$-multinatural and, in particular, 2-natural with respect to pointed multifunctors.
  The components of $\etast$ are shown to be stable equivalences in the proof of \cref{ptmulticat-thm-x}, and hence the induced $(\etast)^\Q$ is a natural stable equivalence
  \[
    1_{\pMulticat^\Q} \fto[\sim]{(\etast)^\Q} \Endst^\Q \Fst^\Q.
  \]
  Thus, by \cref{gjo29}, $\Fst^\Q$ and $\Endst^\Q$ are equivalences of homotopy theories between categories of non-symmetric $\Q$-algebras.
  This completes the proof.
\end{proof}

\section{Multiplicative Homotopy Theory of \texorpdfstring{$\Mone$}{M1}-Modules}
\label{sec:monemod-alg}


Consider the diagram
\begin{equation}\label{FmEmcontext}
\Fm \cn \MoneMod \lradj \permcatsu \cn \Endm
\end{equation}
consisting of the following.
\begin{itemize}
\item $\MoneMod$ is the $\Cat$-multicategory in \cref{expl:monemodcatmulticat}.
\item $\permcatsu$ is the $\Cat$-multicategory in \cref{thm:permcatmulticat}.
\item $\Endm$ is the $\Cat$-multifunctor in \cref{expl:endm-catmulti}.
\item $\Fm$ is the non-symmetric $\Cat$-multifunctor given by the following composite.
\begin{equation}\label{Fm-multi-def}
\begin{tikzpicture}[vcenter]
\def\u{.6}
\draw[0cell]
(0,0) node (a) {\MoneMod}
(a)+(3,0) node (b) {\pMulticat}
(b)+(3,0) node (c) {\permcatsu}
;
\draw[1cell]
(a) edge node {\Um} (b)
(b) edge node {\Fst} (c)
;
\draw[1cell]
(a) [rounded corners=3pt] |- ($(b)+(-1,\u)$)
-- node {\Fm} ($(b)+(1,\u)$) -| (c)
;
\end{tikzpicture}
\end{equation}
In \cref{Fm-multi-def},
\begin{itemize}
\item $\Um$ is the $\Cat$-multifunctor in \cref{expl:Um-catmulti}, and
\item $\Fst$ is the non-symmetric $\Cat$-multifunctor in \cref{ptmulticat-xvii}.
\end{itemize}
\end{itemize}

\begin{explanation}[2-Functors]\label{expl:Fm-multi}
Considering the underlying 2-functors, the diagram \cref{Fm-multi-def} factors as follows.
\begin{equation}\label{Fm-composite}
\begin{tikzpicture}[vcenter]
\def\u{-.6}
\draw[0cell]
(0,0) node (a) {\pMulticat}
(a)+(1.5,1.3) node (m) {\MoneMod}
(a)+(3,0) node (b) {\permcatst}
(b)+(3,0) node (c) {\permcatsu}
;
\draw[1cell=.85]
(a) edge node {\Fst} (b)
(b) edge[right hook->] node {I} (c)
(m) edge node[swap,pos=.7] {\Um} (a)
(m) edge node[pos=.7] {\Fm} (b)
(m) edge[bend left=15] node {\Fm} (c)
;
\draw[1cell=.85]
(a) [rounded corners=3pt] |- ($(b)+(-2,\u)$)
-- node[pos=.3] {\Fst} ($(b)+(1,\u)$) -| (c)
;
\end{tikzpicture}
\end{equation}
\smallskip
\begin{itemize}
\item The interior $\Fm$ is the 2-functor in \cref{Fm-def}.
\item The interior $\Fst$ is the 2-functor in \cref{ptmulticat-thm-i}.
\item $I$ is the inclusion 2-functor in \cref{permcatinclusion}.
\end{itemize}
Thus the non-symmetric $\Cat$-multifunctor $\Fm$ in \cref{Fm-multi-def} extends the 2-functor $\Fm$ in \cref{Fm-def}.
\end{explanation}

\subsection*{Comparison Transformations}

\begin{definition}[Comparing $\Endm \Fm$ and the Identity]\label{def:etam-multi}
In the context of \cref{FmEmcontext}, we define the non-symmetric $\Cat$-multinatural transformation
\begin{equation}\label{etam-multi-def}
\etam \cn 1_{\MoneMod} \to \Endm \Fm
\end{equation}
with, for each left $\Mone$-module $\M$, the component morphism
\[\begin{split}
\etam_\M \cn \M \to & \Endm \Fm \M\\
& = \Endm \Fst\, \Um \M \inspace \MoneMod
\end{split}\]
given by the unit in \cref{etam-define}.
\end{definition}  

\begin{explanation}\label{expl:etam-multi}
As explained in \cref{etam-define}, the left $\Mone$-module morphism $\etam_\M$ is uniquely determined by its image in $\pMulticat$:
\[\begin{split}
\Um\etam_\M = \etast_{\Um\M} \cn \Um\M \to & \Um \Endm \Fst\, \Um \M\\
& = \Endst\, \Fst\, \Um \M.
\end{split}\]
The components $\etam_\M$ define a non-symmetric $\Cat$-multinatural transformation $\etam$ as in \cref{etam-multi-def} by the following facts.
\begin{itemize}
\item $\MoneMod$ is a full sub-2-category of $\pMulticat$ (\cref{proposition:EM2-5-1} \eqref{it:EM251-3}).
\item $\Um$ is a $\Cat$-multifunctor (\cref{expl:Um-catmulti}).
\item $\etast$ is a non-symmetric $\Cat$-multinatural transformation (\cref{ptmulticat-xx}).\defmark
\end{itemize} 
\end{explanation}

\begin{definition}[Comparing $\Fm\Endm$ and the Identity]\label{def:vrhom-multi}
In the context of \cref{FmEmcontext}, we define the non-symmetric $\Cat$-multinatural transformation $\vrhom$,
\begin{equation}\label{vrhom-multi-def}
\begin{aligned}
\vrhom = \vrhost \cn 1_{\permcatsu} \to & \Fm\Endm\\
&= \Fst\, \Um\Endm = \Fst\, \Endst,
\end{aligned}
\end{equation}
as $\vrhost$ in \cref{ptmulticat-xxii}.
\end{definition}

\begin{explanation}\label{expl:vrhom-multi}
For each small permutative category $\C$, the component strictly unital symmetric monoidal functor
\[\vrhom_\C = \vrhost_\C \cn \C \to \Fm\Endm\C = \Fst\, \Endst\C\]
is in \cref{def:vrhostC,expl:rhoc-not-su}.
\end{explanation}

The following result is the $\MoneMod$ analog of \cref{etaEEvrho}.
It is used in \cref{mackey-xiv-mone}, which is the $\MoneMod$ analog of \cref{mackey-xiv-pmulticat}.

\begin{lemma}\label{etavrhoMone}
Suppose $\C$ is a small permutative category.  Then the two left $\Mone$-module morphisms below are equal.
\[\begin{tikzpicture}[baseline={(a.base)}]
\draw[0cell]
(0,0) node (a) {\Endm \C}
(a)+(4.5,0) node (b) {\Endm\Fm\Endm \C}
;
\draw[1cell=.9]
(a) edge[transform canvas={yshift=.7ex}] node {\etam_{\Endm\C}} (b)
(a) edge[transform canvas={yshift=-.4ex}] node[swap] {\Endm \vrhom_\C} (b)
;
\end{tikzpicture}\]
\end{lemma}

\begin{proof}
Left $\Mone$-module morphisms are determined by their underlying multifunctors.  Thus it suffices to show that the two arrows in question are equal as multifunctors.  In this case, the equality between them follows from \cref{etaEEvrho}, which gives the equality
\[\etast_{\Endst\C} = \Endst \vrhost_\C \scs\]
along with \cref{ex:endmc} about $\Endm$, \cref{Fm-def} about $\Fm$,  \cref{etam-define} about $\etam$, and \cref{vrhom-multi-def} about $\vrhom$.
\end{proof}

\subsection*{Equivalences of Homotopy Theories}

Recall from \cref{def:ptmulti-stableeq} the wide subcategory
\[\cSM \bigsubset \MoneMod\]
of $\Fm$-stable equivalences created by 
\[\Fm \cn \brb{\MoneMod,\cSM} \to \brb{\permcatst,\cSI}.\]  
Recall from \cref{def:nonsymalgebra} the notion of \emph{non-symmetric algebras}.  \cref{ptmulticat-xxv} below simultaneously extends
\begin{itemize}
\item \cref{ptmulticat-xxiii} from $\pMulticat$ to $\MoneMod$ and
\item \cref{ptmulticat-thm-xi} to non-symmetric algebras.
\end{itemize}
As in \cref{thm:alg-hty-equiv,ptmulticat-xxiii}, each of the two induced functors $\Fm^\Q$ and $\Endm^\Q$ is given by post-composition and whiskering with the respective functor.

\begin{theorem}\label{ptmulticat-xxv}
  Suppose $\Q$ is a small non-symmetric $\Cat$-multicategory.  In the context of \cref{FmEmcontext}, the induced functors  
  \[
    \Fm^\Q \cn \brb{(\MoneMod)^\Q, (\cSM)^\Q} \lrsimadj \brb{(\permcatsu)^\Q, \cS^\Q} 
    \cn \Endm^\Q
  \]
  are inverse equivalences of homotopy theories in the sense of \cref{def:inverse-heq}.
\end{theorem}
\begin{proof}
  This proof is similar to that of \cref{ptmulticat-xxiii}.
  The functor $\Fm$ creates stable equivalences by definition of $\cSM$, and thus the forgetful $\Um$ also creates stable equivalences, because $\Fm = \Fst\Um$.
  Also recall from \cref{proposition:EM2-5-1}~\cref{it:EM251-3} that $\MoneMod$ is a full sub-2-category of $\pMulticat$.
  Together these imply that
  \begin{itemize}
    \item the components of $\vrhom$ in \cref{vrhom-multi-def} are natural 
      stable equivalences and
    \item the components of $\etam$ in \cref{etam-multi-def} are natural stable equivalences.
  \end{itemize}
  The naturality diagram \cref{eq:vrhost-nat-steq} then shows that $\Endm$
  is a relative functor.
  
  Therefore, we conclude that $(\Fm)^\Q$ and $(\Endm)^\Q$ are both relative functors and the induced $(\vrhom)^\Q$ and $(\etam)^\Q$ are natural stable equivalences
  \[
    1_{(\permcatsu)^\Q} \fto[\sim]{(\vrhom)^\Q} \Fm^\Q\Endm^\Q
    \andspace
    1_{(\MoneMod)^\Q} \fto[\sim]{(\etam)^\Q} \Endm^\Q\Fm^\Q.
  \]

  Thus, by \cref{gjo29}, $\Fm^\Q$ and $\Endm^\Q$ are inverse equivalences of homotopy theories between categories of non-symmetric $\Q$-algebras.
  This completes the proof.
\end{proof}

Recall that there is an adjoint equivalence of homotopy theories
\begin{equation}\label{eq:Monesma-UM-Cat-lradj}
  \Monesma \cn \pMulticat \lrsimadj \MoneMod \cn \Um
\end{equation}
by \cref{ptmulticat-thm-vi}
Recall, furthermore, that
\begin{itemize}
\item $\Monesma$ is $\Cat$-multifunctorial by \cref{Monesma-CatSM} and
\item $\Um$ is $\Cat$-multifunctorial by \cref{expl:Um-catmulti}.
\end{itemize}

Since $\Monesma$ and $\Um$ are $\Cat$-multifunctors in the symmetric sense, induced by symmetric $\Cat$-monoidal functors, the following result holds for both symmetric and non-symmetric $\Q$.
Recall from \cref{def:nonsymalgebra} that, for a $\Cat$-multicategory $\N$, the notation $\N^\Q$ denotes the category of $\Q$-algebras (in the symmetric sense) if $\Q$ is a small $\Cat$-multicategory in the symmetric sense.
\begin{theorem}\label{Monesma-Um-algebra}
  Suppose $\Q$ is either a symmetric or non-symmetric small multicategory.
  The induced functors between categories of algebras (non-symmetric if $\Q$ is non-symmetric)
  \begin{equation}\label{eq:MonesmaQ-UMQ-hteq}
    (\Monesma)^\Q
    \cn \brb{\pMulticat^\Q, \cSst^\Q} \lrsimadj \brb{(\MoneMod)^\Q, (\cSM)^\Q} \cn
    \Um^\Q
  \end{equation}
  are inverse equivalences of homotopy theories in the sense of \cref{def:inverse-heq}.
\end{theorem}
\begin{proof}
  Recall from \cref{expl:MonesmaUmadj,etahat-epzhat-monCatnat} the unit $\etahat$ and counit $\epzhat$ for the adjunction $(\Mone \sma - ) \dashv \Um$ are monoidal $\Cat$-natural transformations, and hence $\Cat$-multinatural.
  It follows, therefore, that $(\Mone \sma -)^\Q$ and $\Um^\Q$ are inverse equivalences of homotopy theories because $\etahat$ and $\epzhat$ are shown to be componentwise stable equivalences in \cref{ptmulticat-thm-vi}.
\end{proof}

\part{Enrichment of Diagrams and Mackey Functors in Closed Multicategories}
\label{part:enrpresheave}

\chapter{Multicategorically Enriched Categories}
\label{ch:menriched}

This chapter defines and develops the basic properties of enrichment in a non-symmetric multicategory $\M$.
This is similar to, but more general than, the concept of enrichment in a monoidal category $\V$ from \cref{ch:prelim_enriched}.
For the special case $\M = \End\,\V$, \cref{EndV-enriched} shows that the two notions of enrichment agree.
\cref{ex:EndV-enriched} lists a number of symmetric monoidal closed categories for which enrichment over $\End\,\V$ applies.

The main application for our further work is $\M = \permcatsu$.
\cref{permcat-selfenr} shows that $\permcatsu$ is a $\permcatsu$-category.
\cref{ev-bilinear,ev-comp} provide further details that will be used in \cref{ch:gspectra,ch:std_enrich} to discuss a corresponding closed structure for $\permcatsu$.

In the case that $\M$ has a symmetric group action making it a multicategory, there is a notion of opposite $\M$-category described in \cref{sec:opp-enr-cat}.
For the special case $\M = \End\,\V$ with $\V$ a symmetric monoidal category, \cref{v-opposite-mcat} shows that the notions of opposite $\M$-category and opposite $\V$-category agree.
This will be used in \cref{ch:gspectra_Kem} for the discussion of enriched diagrams and enriched Mackey functors.

The abstract theory of multicategorical enrichment developed here and in the remaining chapters of \cref{part:enrpresheave} will find homotopy-theoretic applications in \cref{part:homotopy-mackey}.
\cref{ch:mackey} develops conditions for change of enrichment functors, between enriched diagram and Mackey functor categories, to induce equivalences of homotopy theories.
\cref{ch:mackey_eq} applies this to enrichment over, and diagrams in, $\permcatsu$, $\pMulticat$, and $\MoneMod$.

\subsection*{Connection with Other Chapters}

The remaining chapters in this work depend on the multicategorical enrichment developed here.
Change of enrichment is discussed in \cref{ch:change_enr}.
\cref{ch:gspectra} develops the basic theory of closed multicategories and extends the results about $\permcatsu$ in 
\cref{sec:perm-self-enr,sec:eval-perm} to a closed multicategory structure.
In \cref{ch:std_enrich} the theory of self-enriched multicategories is developed further, with additional applications to $\permcatsu$.
Opposite enriched categories (\cref{sec:opp-enr-cat}) are important in \cref{ch:gspectra_Kem,ch:mackey,ch:mackey_eq}, where they are the domains of enriched Mackey functor categories \cref{mcat-copm}.

\subsection*{Background}

The content of this chapter depends only on that of \cref{ch:ptmulticat,ch:prelim,ch:prelim_enriched,ch:prelim_multicat}.

\subsection*{Chapter Summary}

\cref{sec:menriched-cat} defines enrichment in a non-symmetric multicategory.
\cref{sec:enr-endomorphism} shows that enrichment in an endomorphism multicategory agrees with enrichment in the underlying monoidal category.
\cref{sec:enriched-in-perm} describes enrichment in $\permcatsu$, the $\Cat$-multicategory of permutative categories.
As an important special case, \cref{sec:perm-self-enr} describes the self-enrichment of $\permcatsu$.
\cref{sec:eval-perm} explains the bilinear evaluation for $\permcatsu$.
\cref{sec:opp-enr-cat} develops the basic theory for opposites of multicategorically-enriched categories.
Here is a summary table.
\reftable{.9}{
  2-category of $\M$-categories
  & \ref{mcat-iicat}
  \\ \hline
  $\V$-categories and $\End(\V)$ categories
  & \ref{EndV-enriched} and
   \ref{ex:EndV-enriched}
  \\ \hline
  2-category of $\permcatsu$-categories
  & \ref{expl:perm-enr-cat},
   \ref{expl:perm-enr-functors}, and
   \ref{expl:perm-nattr}
  \\ \hline
  self-enrichment of $\permcatsu$
  & \ref{def:perm-selfenr} and
  \ref{permcat-selfenr}
  \\ \hline
  bilinear evaluation for $\permcatsu$
  & \ref{def:perm-evaluation},
   \ref{ev-bilinear}, and
   \ref{ev-comp}
   \\ \hline
  opposite $\M$-categories
  &
  \ref{def:opposite-mcat} and
  \ref{v-opposite-mcat}
  \\
}

We remind the reader of \cref{conv:universe} about universes and \cref{expl:leftbracketing} about left normalized bracketing for iterated products.

\section{Enrichment in a Multicategory}
\label{sec:menriched-cat}

Throughout this section, we assume that 
\[(\M,\ga,\opu)\]
is a non-symmetric multicategory (\cref{def:enr-multicategory}).  This means that $\M$ is a non-symmetric $\Set$-multicategory, with $(\Set,\times,*)$ the symmetric monoidal category of sets and functions with the Cartesian product as the monoidal product.  In this section we define categories, functors, and natural transformations enriched in $\M$.  This section is organized as follows.
\begin{itemize}
\item $\M$-categories, $\M$-functors, and $\M$-natural transformations are in \cref{def:menriched-cat,def:mfunctor,def:mnaturaltr}.
\item Vertical composition of $\M$-natural transformations are discussed in \cref{def:mnaturaltr-vcomp,mnat-vcomp-welldef}.
\item Horizontal composition of $\M$-natural transformations are discussed in \cref{def:mnaturaltr-hcomp,mnat-hcomp-welldef}.
\item \cref{mcat-iicat} proves the existence of a 2-category with small $\M$-categories as objects.
\end{itemize}
We discuss enrichment in the endomorphism multicategory of a monoidal category in \cref{sec:enr-endomorphism}.  In particular, \cref{EndV-enriched} proves that, for a monoidal category $\V$, the 2-category of small $\V$-categories and the 2-category of small $(\End\,\V)$-categories are the same.  Thus the two notions of enrichment over a monoidal category coincide.

\subsection*{$\M$-Categories}

\begin{definition}\label{def:menriched-cat}
An \emph{$\M$-category} $(\C,\mcomp,i)$, which is also called a \index{multicategory!enrichment in a -}\index{enriched!category - in a multicategory}\emph{category enriched in $\M$}, consists of the following data.
\begin{description}
\item[Objects]
$\C$ is equipped with a class $\Ob\C$ of \emph{objects}.  We usually write $x \in \C$ instead of $x \in \Ob\C$.
\item[Hom Objects]
For each pair of objects $x,y \in \C$, $\C$ is equipped with an object\label{not:Cxyhom}
\[\C(x,y) \in \Ob\M,\]
which is called the \index{object!hom}\index{hom object}\emph{hom object} with \emph{domain} $x$ and \emph{codomain} $y$.  We sometimes abbreviate $\C(x,y)$ to $\C_{x,y}$ to save space.
\item[Composition]
For objects $x,y,z \in \C$, $\C$ is equipped with a binary multimorphism
\begin{equation}\label{menriched-cat-comp}
\mcomp_{x,y,z} \cn \big(\C(y,z) \scs \C(x,y)\big) \to \C(x,z) \inspace \M,
\end{equation}
which is called the \index{composition!enriched category}\emph{composition}.
\item[Identities]
Each object $x \in \C$ is equipped with a nullary multimorphism
\begin{equation}\label{menriched-cat-identity}
i_x \cn \ang{} \to \C(x,x) \inspace \M,
\end{equation}
which is called the \emph{identity} of $x$.  Here $\ang{}$ denotes the empty sequence.
\end{description}
These data are required to make the following \index{associativity!enriched category}associativity and \index{unity!enriched category}unity diagrams in $\M$ commute for objects $w,x,y,z \in \C$.
\begin{equation}\label{menriched-cat-assoc}
\begin{tikzpicture}[vcenter]
\def\v{-1.4}
\draw[0cell=.9]
(0,0) node (a) {\big(\C(y,z) \scs \C(x,y) \scs \C(w,x)\big)}
(a)+(5.5,0) node (b) {\big(\C(x,z) \scs \C(w,x)\big)}
(a)+(0,\v) node (c) {\big(\C(y,z) \scs \C(w,y)\big)}
(b)+(0,\v) node (d) {\C(w,z)}
;
\draw[1cell=.9]
(a) edge node {(\mcomp_{x,y,z} \scs \opu)} (b)
(b) edge node {\mcomp_{w,x,z}} (d)
(a) edge node[swap] {(\opu \scs \mcomp_{w,x,y})} (c)
(c) edge node {\mcomp_{w,y,z}} (d)
;
\end{tikzpicture}
\end{equation}
\begin{equation}\label{menriched-cat-unity}
\begin{tikzpicture}[vcenter]
\def\h{3.2} \def\v{-1.4}
\draw[0cell=.9]
(0,0) node (a1) {\big(\C(x,y) \scs \ang{}\big)}
(a1)+(\h,0) node (a2) {\C(x,y)}
(a2)+(\h,0) node (a3) {\big(\ang{} \scs \C(x,y)\big)}
(a1)+(0,\v) node (b1) {\big(\C(x,y) \scs \C(x,x)\big)}
(a2)+(0,\v) node (b2) {\C(x,y)}
(a3)+(0,\v) node (b3) {\big(\C(y,y) \scs \C(x,y)\big)}
;
\draw[1cell=.9]
(a1) edge[-,double equal sign distance] (a2)
(a2) edge[-,double equal sign distance] (a3)
(a1) edge node[swap] {(\opu \scs i_x)} (b1)
(a2) edge node {\opu} (b2)
(a3) edge node {(i_y \scs \opu)} (b3)
(b1) edge node {\mcomp_{x,x,y}} (b2)
(b3) edge node[swap] {\mcomp_{x,y,y}} (b2)
;
\end{tikzpicture}
\end{equation}
This finishes the definition of an $\M$-category.  An $\M$-category is \index{small!enriched category}\index{enriched category!small}\emph{small} if its class of objects is a set.
\end{definition}

\begin{explanation}[$\M$-Categories]\label{expl:mcat}
The composition $\mcomp_{x,y,z}$ in \cref{menriched-cat-comp} is an element in 
\[\M\scmap{\C_{y,z} \scs \C_{x,y}; \C_{x,z}}.\]  
The identity $i_x$ of $x$ in \cref{menriched-cat-identity} is an element in 
\[\M\scmap{\ang{}; \C_{x,x}}.\]
The associativity diagram \cref{menriched-cat-assoc} means the equality of 3-ary multimorphisms
\[\ga\scmap{\mcomp_{w,y,z}; \opu_{\C_{y,z}} \scs \mcomp_{w,x,y}} 
= \ga\scmap{\mcomp_{w,x,z}; \mcomp_{x,y,z} \scs \opu_{\C_{w,x}}}\]
in $\M\scmap{\C_{y,z} \scs \C_{x,y} \scs \C_{w,x}; \C_{w,z}}$.  The unity diagram \cref{menriched-cat-unity} means the equalities of unary multimorphisms
\[\ga\scmap{\mcomp_{x,x,y}; \opu_{\C_{x,y}} \scs i_x} = \opu_{\C_{x,y}} 
= \ga\scmap{\mcomp_{x,y,y}; i_y \scs \opu_{\C_{x,y}}}\]
in $\M\scmap{\C_{x,y};\C_{x,y}}$.  Other commutative diagrams in $\M$ below are interpreted similarly.
\end{explanation}

\subsection*{$\M$-Functors}

\begin{definition}\label{def:mfunctor}
Suppose $(\C,\mcomp,i)$ and $(\D,\mcomp,i)$ are $\M$-categories.  An \emph{$\M$-functor}
\[F \cn \C \to \D,\]
which is also called a \index{functor!enriched}\index{enriched!functor}\emph{functor enriched in $\M$}, consists of the following data.
\begin{description}
\item[Object Assignment]
$F$ is equipped with a function
\[F \cn \Ob\C \to \Ob\D,\]
which is called the \emph{object assignment}.
\item[Component Morphisms]
For each pair of objects $x,y \in \C$, $F$ is equipped with a unary multimorphism
\begin{equation}\label{Fxy-component}
F_{x,y} \cn \C(x,y) \to \D(Fx,Fy) \inspace \M,
\end{equation}
which is called the \emph{$(x,y)$-component} of $F$.  We sometimes abbreviate $F_{x,y}$ to $F$.
\end{description}
These data are required to make the following two diagrams in $\M$ commute for objects $x,y,z \in \C$.
\begin{equation}\label{mfunctor-diagrams}
\begin{tikzpicture}[vcenter]
\def\v{-1.4} \def\h{3.8}
\draw[0cell=.8]
(0,0) node (a) {\big(\C(y,z) \scs \C(x,y)\big)}
(a)+(\h,0) node (b) {\C(x,z)}
(a)+(0,\v) node (c) {\big(\D(Fy,Fz) \scs \D(Fx,Fy)\big)}
(b)+(0,\v) node (d) {\D(Fx,Fz)}
;
\draw[1cell=.8]
(a) edge node {\mcomp_{x,y,z}} (b)
(b) edge node {F_{x,z}} (d)
(a) edge node[swap] {(F_{y,z} \scs F_{x,y})} (c)
(c) edge node {\mcomp_{Fx,Fy,Fz}} (d)
;
\begin{scope}[shift={(\h+1.2,0)}]
\draw[0cell=.8]
(0,0) node (a) {\ang{}}
(a)+(1.7,0) node (b) {\C(x,x)}
(b)+(0,\v) node (c) {\D(Fx,Fx)}
;
\draw[1cell=.8]
(a) edge node {i_x} (b)
(b) edge node {F_{x,x}} (c)
(a) edge node[swap] {i_{Fx}} (c)
;
\end{scope}
\end{tikzpicture}
\end{equation}
This finishes the definition of an $\M$-functor.

Moreover, we define the following.
\begin{itemize}
\item The \index{identity!enriched functor}\emph{identity $\M$-functor} 
\begin{equation}\label{id-mfunctor}
1_\C \cn \C \to \C
\end{equation}
is defined by the identity object assignment and unit $(x,y)$-component
\[(1_\C)_{x,y} = \opu_{\C(x,y)} \cn \C(x,y) \to \C(x,y) \inspace \M\]
for each pair of objects $x,y \in \C$.
\item For $\M$-functors $F \cn \C \to \D$ and $G \cn \D \to \E$, the \index{composition!enriched functor}\emph{composite} $\M$-functor
\begin{equation}\label{mfunctor-composition}
GF \cn \C \to \E
\end{equation}
is defined by composing the object assignments.  The $(x,y)$-component of $GF$ is the composite unary multimorphism 
\begin{equation}\label{Fxy-GFxFy}
\C(x,y) \fto{F_{x,y}} \D(Fx,Fy) \fto{G_{Fx,Fy}} \E(GFx,GFy)
\end{equation}
in $\M$.  The commutativity of the diagrams \cref{mfunctor-diagrams} for $GF$ follows from the commutativity of the corresponding diagrams for $F$ and $G$.  Composition of $\M$-functors is associative and unital with respect to identity $\M$-functors by the associativity and unity axioms of $\M$ in \cref{enr-multicategory-associativity,enr-multicategory-right-unity,enr-multicategory-left-unity}. 
\end{itemize}
This finishes the definition.
\end{definition}

\begin{remark}[History]\label{rk:lambek-mcat}
Enrichment in a non-symmetric multicategory goes back to the beginning of multicategory theory.  It is mentioned in \cite[page 106]{lambek}, immediately after the definition of a non-symmetric multicategory.  The explicit data and axioms of an $\M$-category and an $\M$-functor in \cref{def:menriched-cat,def:mfunctor} are from \cite[Section 2]{bohmann_osorno-mackey}.
\end{remark}

\subsection*{$\M$-Natural Transformations}

\begin{definition}\label{def:mnaturaltr}
  Suppose $F,G \cn \C \to \D$ are $\M$-functors between $\M$-categories.  An \index{enriched!natural transformation}\index{natural transformation!enriched}\emph{$\M$-natural transformation} 
\[\theta \cn F \to G\] 
consists of, for each object $x$ in $\C$, a nullary multimorphism
\begin{equation}\label{mnatural-component}
\theta_x \cn \ang{} \to \D(Fx,Gx) \inspace \M,
\end{equation}
which is called the \index{component}\emph{$x$-component of $\theta$}.   The following \index{naturality diagram}\emph{naturality diagram} in $\M$ is required to commute for each pair of objects $x,y \in \C$, with $\mcomp$ denoting the composition of $\D$.
\begin{equation}\label{mnaturality-diag}
\begin{tikzpicture}[vcenter]
\def\h{4.6} \def\e{1}
\draw[0cell=.9]
(0,0) node (a) {\big(\ang{} \scs \C(x,y)\big)}
(a)+(\h,0) node (b) {\big(\D(Fy,Gy) \scs \D(Fx,Fy)\big)}
(a)+(-\e,-\e) node (a') {\big(\C(x,y) \scs \ang{}\big)}
(a')+(0,-1.3) node (b') {\big(\D(Gx,Gy) \scs \D(Fx,Gx)\big)}
(b')+(\e+\h,0) node (d) {\D(Fx,Gy)}
;
\draw[1cell=.9]
(a) edge node {(\theta_y \scs F_{x,y})} (b)
(b) edge node {\mcomp} (d)
(a) edge[-,double equal sign distance, transform canvas={xshift=-1ex}] (a')
(a') edge node[swap] {(G_{x,y} \scs \theta_x)} (b')
(b') edge node {\mcomp} (d)
;
\end{tikzpicture}
\end{equation}
We also call $\theta$ a \emph{natural transformation enriched in $\M$}.

Moreover, the \index{identity!enriched natural transformation}\emph{identity $\M$-natural transformation} 
\begin{equation}\label{id-mnat}
1_F \cn F \to F
\end{equation}
is defined by the identity components:
\[(1_F)_x = i_{Fx} \cn \ang{} \to \D(Fx,Fx) \forspace x \in \C.\]
We use the 2-cell notation \cref{twocellnotation} for $\M$-natural transformations.
\end{definition}

Next we define vertical and horizontal compositions for $\M$-natural transformations.

\begin{definition}\label{def:mnaturaltr-vcomp}
Suppose $\theta$ and $\psi$ are $\M$-natural transformations between $\M$-functors between $\M$-categories as in the left diagram below.
\begin{equation}\label{mnat-vcomp-diag}
\begin{tikzpicture}[baseline={(a.base)}]
\def\h{2.5} \def\t{40} \def\s{30}
\draw[0cell]
(0,0) node (a) {\C}
(a)+(\h,0) node (b) {\D}
;
\draw[1cell=.8]
(a) edge[bend left=\t] node {F} (b)
(a) edge node[pos=.2] {G} (b)
(a) edge[bend right=\t] node[swap] {H} (b)
;
\draw[2cell]
node[between=a and b at .45, shift={(0,.3)}, rotate=-90, 2label={above,\theta}] {\Rightarrow}
node[between=a and b at .45, shift={(0,-.3)}, rotate=-90, 2label={above,\psi}] {\Rightarrow}
;
\begin{scope}[shift={(4,0)}]
\draw[0cell]
(0,0) node (a) {\C}
(a)+(2,0) node (b) {\D}
;
\draw[1cell=.9]
(a) edge[bend left=\s] node {F} (b)
(a) edge[bend right=\s] node[swap] {H} (b)
;
\draw[2cell]
node[between=a and b at .4, rotate=-90, 2label={above,\psi\theta}] {\Rightarrow}
;
\end{scope}
\end{tikzpicture}
\end{equation}
The \index{vertical composition!enriched natural transformation}\emph{vertical composite} $\M$-natural transformation $\psi\theta$, as in the right diagram above, has, for each object $x$ in $\C$, $x$-component given by the following composite in $\M$, with $\mcomp$ denoting the composition of $\D$.
\begin{equation}\label{mnat-vcomp-component}
\begin{tikzpicture}[baseline={(a.base)}]
\def\w{.6}
\draw[0cell=.9]
(0,0) node (a) {\ang{}}
(a)+(1,0) node (a') {\big(\ang{} \scs \ang{}\big)}
(a')+(4,0) node (b) {\big(\D(Gx,Hx) \scs \D(Fx,Gx)\big)}
(b)+(3.5,0) node (c) {\D(Fx,Hx)}
;
\draw[1cell=1.3]
(a) edge[-,double equal sign distance] (a')
;
\draw[1cell=.9]
(a') edge node {(\psi_x \scs \theta_x)} (b)
(b) edge node {\mcomp} (c)
;
\draw[1cell=.9]
(a) [rounded corners=3pt] |- ($(b)+(-1,\w)$)
-- node[pos=0] {(\psi\theta)_x} ($(b)+(1,\w)$) -| (c)
;
\end{tikzpicture}
\end{equation}
This finishes the definition of $\psi\theta$.
\end{definition}

\begin{lemma}\label{mnat-vcomp-welldef}
In the context of \cref{def:mnaturaltr-vcomp}, the following statements hold.
\begin{enumerate}
\item\label{mnat-vcomp-i} The vertical composite $\psi\theta$ is a well-defined $\M$-natural transformation.
\item\label{mnat-vcomp-ii} Vertical composition of $\M$-natural transformations is associative.
\item\label{mnat-vcomp-iii} Identity $\M$-natural transformations \cref{id-mnat} are two-sided units for vertical composition.
\end{enumerate} 
\end{lemma}

\begin{proof}
\emph{Assertion \cref{mnat-vcomp-i}}.  Suppose $x$ and $y$ are objects in $\C$.  By the definition \cref{mnat-vcomp-component} of $(\psi\theta)_x$, the naturality diagram \cref{mnaturality-diag} for $\psi\theta$ is the boundary of the following diagram in $\M$, where we abbreviate $\C(x,y)$ to $\C_{x,y}$ and likewise for $\D$.
\[\begin{tikzpicture}
\def\h{8} \def\t{-1.3} \def\g{1.2}
\draw[0cell=.8]
(0,0) node (a1) {(\C_{x,y} \scs \ang{} \scs \ang{})}
(a1)+(\h/2,0) node (a2) {(\ang{} \scs \C_{x,y} \scs \ang{})}
(a2)+(\h/2,0) node (a3) {(\ang{} \scs \ang{} \scs \C_{x,y})}
(a1)+(0,\t) node (b1) {(\D_{Hx,Hy} \scs \D_{Gx,Hx} \scs \D_{Fx,Gx})}
(a3)+(0,\t) node (b2) {(\D_{Gy,Hy} \scs \D_{Fy,Gy} \scs \D_{Fx,Fy})}
(b1)+(0,2*\t) node (c1) {(\D_{Hx,Hy} \scs \D_{Fx,Hx})}
(b2)+(0,2*\t) node (c2) {(\D_{Fy,Hy} \scs \D_{Fx,Fy})}
(a2)+(0,\t) node (d1) {(\D_{Gy,Hy} \scs \D_{Gx,Gy} \scs \D_{Fx,Gx})}
(d1)+(-\h/4,\t) node (d2) {(\D_{Gx,Hy} \scs D_{Fx,Gx})}
(d1)+(\h/4,\t) node (d3) {(\D_{Gy,Hy} \scs \D_{Fx,Gy})}
(d1)+(0,2*\t) node (d4) {\D_{Fx,Hy}}
;
\draw[1cell=.8]
(a1) edge[-,double equal sign distance] (a2)
(a2) edge[-,double equal sign distance] (a3)
(a1) edge node[swap] {(H_{x,y} \scs \psi_x \scs \theta_x)} (b1)
(b1) edge node[swap] {(\opu \scs \mcomp)} (c1)
(c1) edge node {\mcomp} (d4)
(a3) edge node {(\psi_y \scs \theta_y \scs F_{x,y})} (b2)
(b2) edge node {(\mcomp \scs \opu)} (c2)
(c2) edge node[swap] {\mcomp} (d4)
(a2) edge node {(\psi_y \scs G_{x,y} \scs \theta_x)} (d1)
(b1) edge node {(\mcomp \scs \opu)} (d2)
(d1) edge node[swap] {(\mcomp \scs \opu)} (d2)
(d1) edge node {(\opu \scs \mcomp)} (d3)
(b2) edge node[swap] {(\opu \scs \mcomp)} (d3)
(d2) edge node {\mcomp} (d4)
(d3) edge node[swap] {\mcomp} (d4)
;
\end{tikzpicture}\]
The diagram above is commutative for the following reasons.
\begin{itemize}
\item The top left and right pentagons are commutative by the naturality \cref{mnaturality-diag} of $\psi$ and $\theta$, respectively.
\item The other three sub-regions are commutative by the associativity \cref{menriched-cat-assoc} of $\D$.
\end{itemize}
This proves that $\psi\theta$ is an $\M$-natural transformation.

\medskip
\emph{Assertion \cref{mnat-vcomp-ii}}.  Suppose $\vphi \cn H \to I$ is an $\M$-natural transformation for an $\M$-functor $I \cn \C \to \D$.  We must show that, for each object $x$ in $\C$, the $x$-components of $(\vphi\psi)\theta$ and $\vphi(\psi\theta)$ are equal.  Consider the following diagram in $\M$.
\[\begin{tikzpicture}
\def\h{1.8} \def\v{1}
\draw[0cell=.85]
(0,0) node (a) {(\ang{} \scs \ang{} \scs \ang{})}
(a)+(4.5,0) node (b) {(\D_{Hx,Ix} \scs \D_{Gx,Hx} \scs \D_{Fx,Gx})}
(b)+(\h,\v) node (c1) {(\D_{Gx,Ix} \scs \D_{Fx,Gx})}
(b)+(\h,-\v) node (c2) {(\D_{Hx,Ix} \scs \D_{Fx,Hx})}
(b)+(2*\h,0) node (d) {\D_{Fx,Ix}}
;
\draw[1cell=.85]
(a) edge node {(\vphi_x \scs \psi_x \scs \theta_x)} (b)
(b) edge[transform canvas={xshift=-1em}] node[pos=.2] {(\mcomp \scs \opu)} (c1)
(b) edge[transform canvas={xshift=-1em}] node[swap,pos=.2] {(\opu \scs \mcomp)} (c2)
(c1) edge[transform canvas={xshift=1em}] node[pos=.6] {\mcomp} (d)
(c2) edge[transform canvas={xshift=1em}] node[swap,pos=.6] {\mcomp} (d)
;
\end{tikzpicture}\]
The following statements hold for the diagram above.
\begin{itemize}
\item The composite along the top is the $x$-component of $(\vphi\psi)\theta$.
\item The composite along the bottom is the $x$-component of $\vphi(\psi\theta)$.
\item The right quadrilateral commutes by the associativity \cref{menriched-cat-assoc} of $\D$.
\end{itemize}
This proves that the $\M$-natural transformations $(\vphi\psi)\theta$ and $\vphi(\psi\theta)$ are equal.

\medskip
\emph{Assertion \cref{mnat-vcomp-iii}}.  We must show that $\theta 1_F$ and $1_G \theta$ have, for each object $x$ in $\C$, the same $x$-component as $\theta$.  We consider the following diagram in $\M$.
\[\begin{tikzpicture}
\def\h{4} \def\v{-1.3}
\draw[0cell=.85]
(0,0) node (a1) {(\ang{} \scs \ang{})}
(a1)+(\h,0) node (a2) {(\D_{Fx,Gx} \scs \ang{}) = (\ang{} \scs \D_{Fx,Gx})}
(a2)+(\h,0) node (a3) {(\ang{} \scs \ang{})}
(a1)+(0,\v) node (b1) {(\D_{Fx,Gx} \scs \D_{Fx,Fx})}
(a2)+(0,\v) node (b2) {\D_{Fx,Gx}}
(a3)+(0,\v) node (b3) {(\D_{Gx,Gx} \scs \D_{Fx,Gx})}
;
\draw[1cell=.85]
(a1) edge node {\theta_x} (a2)
(a3) edge node[swap] {\theta_x} (a2)
(a1) edge node[swap] {(\theta_x \scs i_{Fx})} (b1)
(a2) edge node[swap,pos=.6] {(\opu \scs i_{Fx})} (b1)
(a2) edge node {\opu} (b2)
(a2) edge node[pos=.6] {(i_{Gx} \scs \opu)} (b3)
(a3) edge node {(i_{Gx} \scs \theta_x)} (b3)
(b1) edge node {\mcomp} (b2)
(b3) edge node[swap] {\mcomp} (b2)
;
\end{tikzpicture}\]
The following statements hold for the diagram above.
\begin{itemize}
\item The left-bottom composite, $\ga\scmap{\mcomp; \theta_x, i_{Fx}}$, is the $x$-component of $\theta 1_F$.  This composite is equal to the Z-shaped composite by the associativity \cref{enr-multicategory-associativity} and left unity \cref{enr-multicategory-left-unity} of $\M$.
\item The right-bottom composite, $\ga\scmap{\mcomp; i_{Gx},\theta_x}$, is the $x$-component of $1_G \theta$.  This composite is equal to the S-shaped composite by the associativity \cref{enr-multicategory-associativity} and left unity \cref{enr-multicategory-left-unity} of $\M$.
\item The composite $\opu \theta_x$ is equal to $\theta_x$ by the left unity \cref{enr-multicategory-left-unity} of $\M$.
\item The middle two triangles are commutative by the unity \cref{menriched-cat-unity} of $\D$.
\end{itemize}
This proves that $\theta 1_F$ and $1_G \theta$ are both equal to $\theta$.
\end{proof}

\begin{definition}\label{def:mnaturaltr-hcomp}
Suppose $\theta$ and $\theta'$ are $\M$-natural transformations between $\M$-functors between $\M$-categories as in the left diagram below.
\[
\]
The \index{horizontal composition!enriched natural transformation}\emph{horizontal composite} $\M$-natural transformation $\theta' \!\ast \theta$, as in the right diagram above, has, for each object $x$ in $\C$, $x$-component given by the following composite in $\M$, with $\mcomp$ denoting the composition of $\E$.
\begin{equation}\label{mnat-hcomp-component}
%
\end{equation}
This finishes the definition of $\theta' \!\ast \theta$.
\end{definition}

\begin{lemma}\label{mnat-hcomp-welldef}
In the context of \cref{def:mnaturaltr-hcomp}, the following statements hold.
\begin{enumerate}
\item\label{mnat-hcomp-i} 
$(\theta' \!\ast \theta)_x$ in \cref{mnat-hcomp-component} is equal to the following composite in $\M$.
\begin{equation}\label{mnat-hcomp-comp}
%
\end{equation}
\item\label{mnat-hcomp-ii} 
$\theta' \!\ast \theta$ is a well-defined $\M$-natural transformation.
\item\label{mnat-hcomp-iii}
Horizontal composition of $\M$-natural transformations is associative.
\item\label{mnat-hcomp-iv}
Identity $\M$-natural transformations \cref{id-mnat} of identity $\M$-functors \cref{id-mfunctor} are two-sided units for horizontal composition.
\item\label{mnat-hcomp-v}
Horizontal composition preserves identity $\M$-natural transformations \cref{id-mnat}.
\item\label{mnat-hcomp-vi}
Horizontal composition preserves vertical composition of $\M$-natural transformations \cref{mnat-vcomp-diag}.
\end{enumerate} 
\end{lemma}

\begin{proof}
\emph{Assertion \cref{mnat-hcomp-i}}.  This follows from the naturality of $\theta'$ \cref{mnaturality-diag}, which implies that the right composites in \cref{mnat-hcomp-component,mnat-hcomp-comp} are equal.

\medskip
\emph{Assertion \cref{mnat-hcomp-ii}}.  By \cref{Fxy-GFxFy,mnat-hcomp-component}, for objects $x,y \in \C$ the naturality diagram \cref{mnaturality-diag} for $\theta' \!* \theta$ is the boundary of the following diagram in $\M$, with $\C(x,y)$ abbreviated to $\C_{x,y}$ and likewise for $\D$ and $\E$.
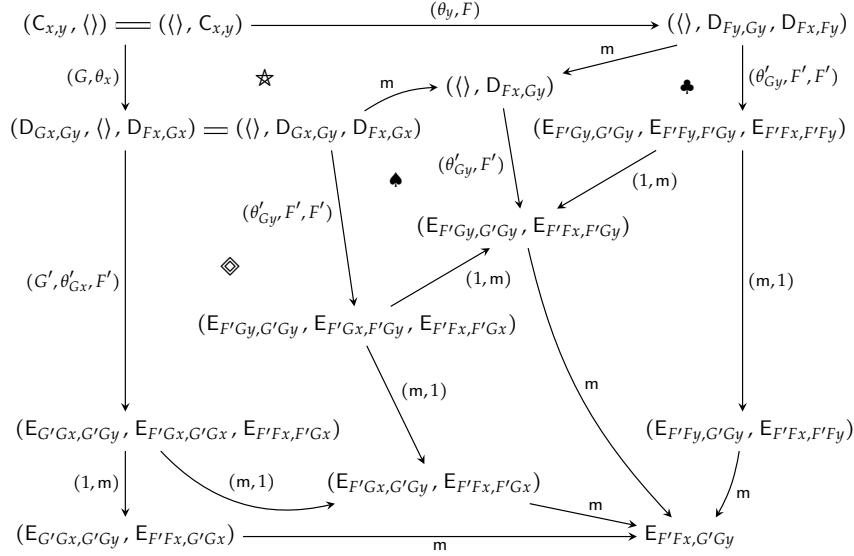
\[
\]
The following statements hold for the diagram above.
\begin{itemize}
\item The sub-regions labeled $\pentagram$ and $\diamonddiamond$ are commutative by the naturality \cref{mnaturality-diag} of $\theta$ and $\theta'$, respectively.
\item The sub-regions labeled $\spadesuit$ and $\clubsuit$ are commutative by the compatibility of $F'$ with composition \cref{mfunctor-diagrams}.
\item The three unlabeled sub-regions are commutative by the associativity of the composition of $\E$ \cref{menriched-cat-assoc}.
\end{itemize}
This proves that $\theta' \!* \theta$ is an $\M$-natural transformation.

\medskip
\emph{Assertion \cref{mnat-hcomp-iii}}.  To show that horizontal composition is associative, consider horizontally composable $\M$-natural transformations as follows.  
\[
\]
To show that, for each object $x$ in $\B$, $\theta'' \!* (\theta' \!* \theta)$ and $(\theta'' \!* \theta') * \theta$ have the same $x$-components, we consider the following diagram in $\M$.
\[%
\]
The following statements hold for the diagram above.
\begin{itemize}
\item The left-bottom composite is $\big(\theta'' \!* (\theta' \!* \theta)\big)_x$. 
\item The other composite is $\big((\theta'' \!* \theta') * \theta\big)_x$.
\item The left quadrilateral is commutative by the compatibility of $F''$ with composition \cref{mfunctor-diagrams}. 
\item The right quadrilateral is commutative by the associativity of the composition of $\E$ \cref{menriched-cat-assoc}.
\end{itemize}
This proves that $\theta'' \!* (\theta' \!* \theta)$ is equal to $(\theta'' \!* \theta') * \theta$.

\medskip
\emph{Assertion \cref{mnat-hcomp-iv}}.  Consider the following $\M$-natural transformations.
\[
\]
To show that $1_{1_\D} \!* \theta$ is equal to $\theta$, we consider, for each object $x \in \C$, the following diagram in $\M$.
\[%
\]
The following statements hold for the diagram above.
\begin{itemize}
\item The left-bottom-right composite is the $x$-component of $1_{1_\D} \!* \theta$.
\item The upper left triangle is commutative by the left unity of $\M$ \cref{enr-multicategory-left-unity}.
\item The lower right triangle is commutative by the unity of $\D$ \cref{menriched-cat-unity}.
\end{itemize}
This proves that $1_{1_\D} \!* \theta$ is equal to $\theta$.

To show that $\theta * 1_{1_\C}$ is equal to $\theta$, we consider, for each object $x \in \C$, the following diagram in $\M$.
\[%
\]
The following statements hold for the diagram above.
\begin{itemize}
\item The left-bottom-right composite is the $x$-component of $\theta * 1_{1_{\C}}$.
\item The bottom sub-region and the middle quadrilateral are commutative by definition.
\item The top sub-region is commutative by the left unity of $\M$ \cref{enr-multicategory-left-unity}.
\item The left triangle is commutative by the compatibility of $F$ with identities \cref{mfunctor-diagrams}.
\item The right triangle is commutative by the unity of $\D$ \cref{menriched-cat-unity}.
\end{itemize}
This proves that $\theta * 1_{1_\C}$ is equal to $\theta$.

\medskip
\emph{Assertion \cref{mnat-hcomp-v}}.  Consider two horizontally composable identity $\M$-natural transformations as follows.
\[\begin{tikzpicture}
\def\h{2} \def\t{25}
\draw[0cell]
(0,0) node (a) {\C}
(a)+(\h,0) node (b) {\D}
(b)+(\h,0) node (c) {\E}
;
\draw[1cell=.8]
(a) edge[bend left=\t] node {F} (b)
(a) edge[bend right=\t] node[swap] {F} (b)
(b) edge[bend left=\t] node {F'} (c)
(b) edge[bend right=\t] node[swap] {F'} (c)
;
\draw[2cell]
node[between=a and b at .45, rotate=-90, 2label={above,1_F}] {\Rightarrow}
node[between=b and c at .45, rotate=-90, 2label={above,1_{F'}}] {\Rightarrow}
;
\end{tikzpicture}\]
To show that $1_{F'} * 1_F$ is equal to $1_{F'F}$, we consider, for each object $x \in \C$, the following diagram in $\M$.
\[\begin{tikzpicture}
\def\h{3.5} \def\v{-1.4} \def\w{-.6}
\draw[0cell=.85]
(0,0) node (a1) {\ang{}}
(a1)+(\h,0) node (a2) {\E_{F'Fx,F'Fx}}
(a1)+(-\h/2,\v) node (b1) {(\ang{} \scs \D_{Fx,Fx})}
(b1)+(\h,0) node (b2) {(\ang{} \scs \E_{F'Fx,F'Fx})}
(b2)+(\h,0) node (b3) {\phantom{(\E_{F'Fx,F'Fx} \scs \E_{F'Fx,F'Fx})}}
(b3)+(.8,0) node (b3') {(\E_{F'Fx,F'Fx} \scs \E_{F'Fx,F'Fx})}
;
\draw[1cell=.9]
(a1) edge node {i_{F'Fx}} (a2)
(a1) edge node[pos=.6] {i_{F'Fx}} (b2)
(b2) edge node {\opu} (a2)
(a1) edge node[swap] {i_{Fx}} (b1)
(b1) edge node {F'} (b2)
(b2) edge node {(i_{F'Fx} \scs \opu)} (b3')
(b3) edge node[swap] {\mcomp} (a2)
;
\draw[1cell=.9]
(b1) [rounded corners=3pt] |- ($(b2)+(-1,\w)$)
-- node[swap] {(i_{F'Fx} \scs F')} ($(b2)+(1,\w)$) -| (b3)
;
\end{tikzpicture}\]
The following statements hold for the diagram above.
\begin{itemize}
\item The left-bottom-right composite is the $x$-component of $1_{F'} * 1_{F}$.
\item The bottom sub-region is commutative by definition.
\item The left triangle is commutative by the compatibility of $F'$ with identities \cref{mfunctor-diagrams}.
\item The middle triangle is commutative by the left unity of $\M$ \cref{enr-multicategory-left-unity}.
\item The right triangle is commutative by the unity of $\E$ \cref{menriched-cat-unity}.
\end{itemize}
This proves that horizontal composition preserves identity $\M$-natural transformations.

\medskip
\emph{Assertion \cref{mnat-hcomp-vi}}.  To show that horizontal composition preserves vertical composition, consider the following $\M$-natural transformations.
\[\begin{tikzpicture}[baseline={(a.base)}]
\def\h{2.5} \def\t{40} \def\s{30} 
\draw[0cell]
(0,0) node (a) {\C}
(a)+(\h,0) node (b) {\D}
(b)+(\h,0) node (c) {\E}
;
\draw[1cell=.8]
(a) edge[bend left=\t] node {F} (b)
(a) edge node[pos=.2] {G} (b)
(a) edge[bend right=\t] node[swap] {H} (b)
(b) edge[bend left=\t] node {F'} (c)
(b) edge node[pos=.2] {G'} (c)
(b) edge[bend right=\t] node[swap] {H'} (c)
;
\draw[2cell]
node[between=a and b at .45, shift={(0,.3)}, rotate=-90, 2label={above,\theta}] {\Rightarrow}
node[between=a and b at .45, shift={(0,-.3)}, rotate=-90, 2label={above,\psi}] {\Rightarrow}
node[between=b and c at .45, shift={(0,.3)}, rotate=-90, 2label={above,\theta'}] {\Rightarrow}
node[between=b and c at .45, shift={(0,-.3)}, rotate=-90, 2label={above,\psi'}] {\Rightarrow}
;
\end{tikzpicture}\]
We must show that, for each object $x$ in $\C$, the following equality holds in $\M$.
\begin{equation}\label{mcat-middle-two}
\big((\psi' \!* \psi)(\theta' \!* \theta)\big)_x = \big((\psi'\theta') * (\psi\theta)\big)_x
\end{equation}
To prove \cref{mcat-middle-two} we use the notation
\[\begin{split}
X &= (\E_{G'Hx,H'Hx} \scs \E_{G'Gx,G'Hx} \scs \E_{F'Gx,G'Gx} \scs \E_{F'Fx,F'Gx}) \scs\\
Y &= (\E_{G'Gx,H'Hx} \scs \E_{F'Gx,G'Gx} \scs \E_{F'Fx,F'Gx}) \scs \andspace\\
Z &= (\E_{G'Hx,H'Hx} \scs \E_{F'Hx,G'Hx} \scs \E_{F'Fx,F'Hx}) 
\end{split}\]
and consider the following diagram in $\M$.
\[\begin{tikzpicture}
\def\h{4} \def\u{1.5} \def\v{3} \def\w{.6}
\draw[0cell=.75]
(0,0) node (a1) {(\ang{} \scs \ang{})}
(a1)+(\h,\u) node (a2) {(\ang{} \scs \D_{Gx,Hx} \scs \D_{Fx,Gx})}
(a2)+(\h,-\u) node (a3) {(\ang{} \scs \ang{} \scs \D_{Fx,Hx})}
(a1)+(0,-\u) node (b1) {(\ang{} \scs \D_{Gx,Hx} \scs \ang{} \scs \D_{Fx,Gx})}
(a3)+(0,-\u) node (b2) {Z}
(b1)+(0,-\u) node (c1) {X}
(b2)+(0,-\u) node (c2) {(\E_{F'Hx,H'Hx} \scs \E_{F'Fx,F'Hx})}
(c1)+(0,-\v) node (d1) {(\E_{G'Gx,H'Hx} \scs \E_{F'Fx,G'Gx})}
(c2)+(0,-\v) node (d2) {\E_{F'Fx,H'Hx}}
(a2)+(0,-\u) node (e1) {(\ang{} \scs \E_{F'Hx,G'Hx} \scs \E_{F'Gx,F'Hx} \scs \E_{F'Fx,F'Gx})}
(e1)+(0,-\u-.7) node (e2) {(\E_{G'Hx,H'Hx} \scs \E_{F'Gx,G'Hx} \scs \E_{F'Fx,F'Gx})}
(c1)+(\h/3,-\v/2) node (e3) {Y}
(c2)+(-\h/2,-\v/2) node (e4) {(\E_{G'Hx,H'Hx} \scs \E_{F'Fx,G'Hx})}
(d1)+(\h,0) node (e5) {(\E_{F'Gx,H'Hx} \scs \E_{F'Fx,F'Gx})}
;
\draw[1cell=.75]
(a1) edge node[swap] {(\psi_x,\theta_x)} (b1)
(b1) edge node[swap] {(\psi'_{Hx},G',\theta'_{Gx},F')} (c1)
(c1) edge node[swap] {(\mcomp,\mcomp)} (d1)
(a1) edge node {(\psi_x,\theta_x)} (a2)
(a2) edge node {\mcomp} (a3)
(a3) edge node {(\psi'_{Hx},\theta'_{Hx},F')} (b2)
(b2) edge node {(\mcomp, \opu)} (c2)
(c2) edge node {\mcomp} (d2)
;
\draw[1cell=.8]
(d1) [rounded corners=3pt] |- ($(e5)+(-2,-\w)$)
-- node[pos=0] {\mcomp} ($(e5)+(1,-\w)$) -| (d2)
;
\draw[1cell=.75]
(a2) edge node[pos=.7] {(\theta'_{Hx},F',F')} (e1)
(e1) edge node {(\psi'_{Hx},\opu,\mcomp)} (b2)
(e1) edge node {(\psi'_{Hx},\mcomp,\opu)} (e2)
(c1) edge node {(\opu,\mcomp,\opu)} (e2)
(c1) edge node {(\mcomp,\opu,\opu)} (e3)
(e3) edge node[pos=.3] {(\opu,\mcomp)} (d1)
(e3) edge node[pos=.4] {(\mcomp,\opu)} (e5)
(e5) edge node {\mcomp} (d2)
(b2) edge[bend right=35] node[pos=.3] {(\opu,\mcomp)} (e4)
(e2) edge node[swap] {(\mcomp,\opu)} (e5)
(e2) edge node {(\opu,\mcomp)} (e4)
(e4) edge node {\mcomp} (d2)
;
\end{tikzpicture}\]
The following statements hold for the diagram above.
\begin{itemize}
\item The left-bottom composite is the left-hand side of \cref{mcat-middle-two}.
\item The top-right composite is the right-hand side of \cref{mcat-middle-two}. 
\item The upper left quadrilateral is commutative by the naturality of $\theta'$ \cref{mnaturality-diag}.
\item The upper right quadrilateral is commutative by the compatibility of $F'$ with composition \cref{mfunctor-diagrams}.
\item The lower left triangle is commutative by definition.
\item The other five sub-regions are commutative by the associativity of the composition of $\E$ \cref{menriched-cat-assoc}.
\end{itemize}
This proves the desired equality \cref{mcat-middle-two}.
\end{proof}

Recall the notion of a 2-category in \cref{def:twocategory}.

\begin{theorem}\label{mcat-iicat}\index{2-category!of small enriched categories!over a non-symmetric multicategory}
For each non-symmetric multicategory $(\M,\ga,\opu)$, there is a 2-category
\[\MCat\]
defined by the following data.
\begin{itemize}
\item The objects are small $\M$-categories (\cref{def:menriched-cat}).
\item The 1-cells are $\M$-functors (\cref{def:mfunctor}).
\item Identity 1-cells are identity $\M$-functors \cref{id-mfunctor}.
\item Horizontal composition of 1-cells is composition of $\M$-functors \cref{mfunctor-composition}.
\item The 2-cells are $\M$-natural transformations (\cref{def:mnaturaltr}).
\item Identity 2-cells are identity $\M$-natural transformations \cref{id-mnat}.
\item Vertical and horizontal compositions of 2-cells are those of $\M$-natural transformations (\cref{def:mnaturaltr-vcomp,def:mnaturaltr-hcomp}).
\end{itemize}
\end{theorem}

\begin{proof}
Axioms \cref{twocat-i,twocat-ii,twocat-iii,twocat-iv} in \cref{def:twocategory} of a 2-category hold for $\MCat$ by, respectively, 
\begin{romenumerate}
\item \cref{mnat-vcomp-welldef},
\item \cref{mnat-hcomp-welldef} \cref{mnat-hcomp-v,mnat-hcomp-vi},
\item \cref{def:mfunctor}, and
\item \cref{mnat-hcomp-welldef} \cref{mnat-hcomp-iii,mnat-hcomp-iv}.
\end{romenumerate}
This finishes the proof.
\end{proof}

\begin{example}\label{ex:permcatsu-enr-cat}
By \cref{thm:permcatmulticat} $\permcatsu$ is a $\Cat$-multicategory, in particular a multicategory.  By \cref{mcat-iicat} there is a 2-category $\permcatsucat$ of small categories, functors, and natural transformations enriched in the multicategory $\permcatsu$.  We describe this 2-category more explicitly in \cref{sec:enriched-in-perm}.
\end{example}

\section{Enrichment in an Endomorphism Multicategory}
\label{sec:enr-endomorphism}

For a monoidal category $\V$ (\cref{def:monoidalcategory}), there are two notions of enrichment over $\V$.
\begin{enumerate}
\item By \cref{ex:vcatastwocategory} there is a 2-category $\VCat$ of small $\V$-categories, $\V$-functors, and $\V$-natural transformations.
\item By \cref{ex:endc} there is a non-symmetric endomorphism multicategory $\End\, \V$, which is, furthermore, a multicategory if $\V$ is symmetric monoidal.  By \cref{mcat-iicat} there is a 2-category $\EndVCat$ of small $(\End\,\V)$-categories, $(\End\,\V)$-functors, and $(\End\,\V)$-natural transformations.
\end{enumerate} 
Next we observe that these two notions of enrichment are the same.  This result is stated in \cite[Remark 2.10]{bohmann_osorno-mackey}.

\begin{proposition}\label{EndV-enriched}\index{category!enriched}\index{enriched category}
For each monoidal category $(\V,\otimes,\tu)$, there is an equality of 2-categories
\[\VCat = \EndVCat.\]
\end{proposition}

\begin{proof}
The objects in $\VCat$ and $\EndVCat$ are small $\V$-categories and small $(\End\,\V)$-categories, respectively.  The identification of these objects follows by comparing
\begin{itemize}
\item \cref{def:menriched-cat} for $(\EndV)$-categories and
\item \cref{def:enriched-category} for $\V$-categories.
\end{itemize} 
More explicitly, suppose $\C$ is a $\V$-category.
\begin{itemize}
\item The identity of an object $x \in \C$ is a morphism 
\[i_x \cn \tu \to \C(x,x) \in \V.\]
By the definition \cref{endc-angxy} of $\EndV$, such a morphism $i_x$ is the same as a nullary multimorphism in 
\[(\End\,\V)\scmap{\ang{}; \C(x,x)} = \V\big(\tu, \C(x,x)\big).\]
\item For objects $x,y,z \in \C$, the multiplication of $\C$ is a morphism
\[\mcomp_{x,y,z} \cn \C(y,z) \otimes \C(x,y) \to \C(x,z) \inspace \V.\]
This is the same as a binary multimorphism in 
\[(\End\,\V)\scmap{\C(y,z), \C(x,y); \C(x,z)} = \V\big(\C(y,z) \otimes \C(x,y) \scs \C(x,z)\big).\]
\item Under the above identifications, the associativity axiom \cref{enriched-cat-associativity} and the unity axiom \cref{enriched-cat-unity} of a $\V$-category are equivalent to those of an $(\End\,\V)$-category in \cref{menriched-cat-assoc,menriched-cat-unity}.
\end{itemize}
Thus $\C$ is equivalently an $(\End\,\V)$-category.

A similar comparison of
\begin{itemize}
\item \cref{def:mfunctor,def:mnaturaltr,def:mnaturaltr-vcomp,def:mnaturaltr-hcomp} for $\EndVCat$ and
\item \cref{def:enriched-functor,def:enriched-natural-transformation} for $\VCat$
\end{itemize} 
proves that the rest of the 2-category structures in $\VCat$ and $\EndVCat$---namely, (identity) 1-cells, (identity) 2-cells, vertical composition of 2-cells, and horizontal composition of 1-cells and 2-cells---are the same.  For example, the two diagrams in \cref{eq:enriched-composition} for a $\V$-functor are equivalent to the diagrams in \cref{mfunctor-diagrams} for an $(\End\,\V)$-functor.  The naturality diagram \cref{enr-naturality} for a $\V$-natural transformation is equivalent to the naturality diagram \cref{mnaturality-diag} for an $(\End\,\V)$-natural transformation.
\end{proof}

For a monoidal category $\V$, \cref{EndV-enriched} identifies $\V$-enrichment and $(\End\,\V$)-enrichment, with $\End\,\V$ the non-symmetric endomorphism multicategory.  In what follows, we use them interchangeably.

\begin{example}[Self-Enrichment]\label{ex:V-EndV-enr}\index{canonical self-enrichment}\index{self-enrichment}\index{symmetric monoidal category!closed!self-enrichment}\index{category!symmetric monoidal!closed!self-enrichment}
Suppose $\V$ is a symmetric monoidal closed category (\cref{def:closedcat}).  Then $\V$ is also a symmetric monoidal $\V$-category with the canonical self-enrichment (\cref{theorem:v-closed-v-sm}).  By \cref{EndV-enriched} the canonical self-enrichment of $\V$ is equal to an $(\End\,\V)$-enrichment, making $\V$ into an $(\End\,\V)$-category.
\end{example}

\begin{example}\label{ex:EndV-enriched}
\cref{EndV-enriched,ex:V-EndV-enr} apply to the following symmetric monoidal closed categories:
\begin{itemize}
\item $\Multicat$ of small multicategories (\cref{theorem:multicat-sm-closed});
\item $\pMulticat$ of small pointed multicategories (\cref{thm:pmulticat-smclosed});
\item $\MoneMod$ of left $\Mone$-modules (\cref{proposition:EM2-5-1} \eqref{monebicomplete});
\item $\pC$ of pointed objects in a complete and cocomplete symmetric monoidal closed category $\C$ (\cref{theorem:pC-sm-closed});
\item $\DstarV$ of pointed diagrams in a complete and cocomplete symmetric monoidal closed category $\V$ with a chosen terminal object (\cref{thm:Dgm-pv-convolution-hom});
\item $\GaV$ of $\Ga$-objects in $\V$ \cref{GammaV};
\item $\GstarV$ of $\Gstar$-objects in $\V$ \cref{Gstar-V}; and
\item $\Sp$ of symmetric spectra \cref{SymSp}.
\end{itemize}
However, \cref{EndV-enriched,ex:V-EndV-enr} do not apply to $\permcatsu$ (\cref{thm:permcatmulticat}) because its multicategory structure is not induced by a monoidal structure.  See \cite[5.7.23 and 10.2.17]{cerberusIII}.  We discuss its self-enrichment in \cref{permcat-selfenr}.
\end{example}

\section[Enrichment in Permutative Categories]{Enrichment in the Multicategory of Permutative Categories}
\label{sec:enriched-in-perm}

Recall from \cref{thm:permcatmulticat} that $\permcatsu$ is a $\Cat$-multicategory, hence also a multicategory. 
\begin{itemize}
\item Its objects are small permutative categories (\cref{def:symmoncat}).
\item Its $n$-ary multimorphisms are $n$-linear functors (\cref{def:nlinearfunctor}).
\item Its colored units are identity symmetric monoidal functors.
\item Its symmetric group action and composition are in \cref{definition:permcat-action,definition:permcat-comp}, respectively.
\end{itemize}
In this section we explicitly describe categories, functors, and natural transformations enriched in $\permcatsu$ in, respectively, \cref{expl:perm-enr-cat,expl:perm-enr-functors,expl:perm-nattr}.  We discuss the self-enrichment of $\permcatsu$ in \cref{sec:perm-self-enr}.

We use the shortened notation
\begin{equation}\label{permcatsu-psu}
\psu = \permcatsu
\end{equation}
to simplify the presentation.

\begin{explanation}[$\psu$-Categories]\label{expl:perm-enr-cat}\index{permutative category!enrichment in the multicategory of}\index{category!permutative!enrichment in the multicategory of}\index{category!enriched in permutative categories}
Unpacking \cref{def:menriched-cat} with $\M = \psu$, a $\psu$-category $(\C,\mcomp,i)$ consists of the following data.
\begin{description}
\item[Objects] $\C$ is equipped with a class $\Ob\C$ of objects.
\item[Hom Permutative Categories] For each pair of objects $x,y \in \C$, $\C$ is equipped with a small permutative category (\cref{def:symmoncat})
\[\big(\C(x,y), \oplus, \pu, \xi\big),\]
which is also denoted $\C_{x,y}$.
\item[Composition] For each triple of objects $x,y,z \in \C$, $\C$ is equipped with a bilinear functor (\cref{def:nlinearfunctor})
\[\mcomp_{x,y,z} \cn \C(y,z) \times \C(x,y) \to \C(x,y).\]
Its first and second linearity constraints are denoted \label{not:mcompi}$\mcomp_{x,y,z}^1$ and $\mcomp_{x,y,z}^2$, respectively.  If there is no danger of confusion, we sometimes omit the subscripts.
\item[Identities] Recalling that a 0-linear functor is a choice of an object in the codomain category, each object $x \in \C$ is equipped with an object
\[i_x \in \C(x,x).\]
This is also regarded as a functor $\boldone \to \C(x,x)$ from the terminal category.  We emphasize that $i_x$ is \emph{not} required to be the monoidal unit $\pu$ of $\C(x,x)$.
\end{description}

Below, we describe some implications of this structure.
We abbreviate each $\mcomp(-,-)$ to concatenation, so
\begin{equation}\label{mxyzfg}
  \mcomp_{x,y,z}(f,g) = fg.
\end{equation}

The constraint 2-by-2 axiom \cref{eq:f2-2by2} for $\mcomp$ says that, for objects $x,y,z\in\C$, with
\begin{equation}\label{m2by2data}
  f,f' \in \C(y,z), \andspace g,g' \in \C(x,y), 
\end{equation}
the following diagram in $\C(x,z)$ commutes.
\begin{equation}\label{psucat-m-2by2}
  \begin{tikzpicture}[x=35mm,y=15mm,vcenter]
    \draw[0cell] 
    (0,0) node (a) {fg \oplus f'g \oplus fg' \oplus f'g'}
    (a)+(1,1) node (a') {(f \oplus f')g \oplus (f \oplus f')g'}
    (a')+(.5,-1.5) node (c) {(f \oplus f')(g \oplus g')}
    (a)+(0,-1) node (b) {fg \oplus fg' \oplus f'g \oplus f'g'}
    (b)+(1,-1) node (b') {f(g \oplus g') \oplus f'(g \oplus g')}
    ;
    \draw[1cell] 
    (a) edge['] node {1 \oplus \xi \oplus 1} (b)
    (a) edge node {\mcomp^1 \oplus \mcomp^1} (a')
    (b) edge['] node {\mcomp^2 \oplus \mcomp^2} (b')
    (a') edge node {\mcomp^2} (c)
    (b') edge['] node {\mcomp^1} (c)
    ;
  \end{tikzpicture}
\end{equation}

The associativity axiom \cref{menriched-cat-assoc} of a $\psu$-category says that, for objects $w,x,y,z \in \C$, the following two composite 3-linear functors are equal, with each 1 denoting the identity symmetric monoidal functor.
\begin{equation}\label{psucat-associativity}
\begin{tikzpicture}[vcenter]
\def\v{-1.4}
\draw[0cell=.9]
(0,0) node (a1) {\C(y,z) \times \C(x,y) \times \C(w,x)}
(a1)+(5,0) node (a2) {\C(x,z) \times \C(w,x)}
(a1)+(0,\v) node (b1) {\C(y,z) \times \C(w,y)}
(a2)+(0,\v) node (b2) {\C(w,z)} 
;
\draw[1cell=.9]
(a1) edge node {\mcomp_{x,y,z} \times 1} (a2)
(a2) edge node {\mcomp_{w,x,z}} (b2)
(a1) edge node[swap] {1 \times \mcomp_{w,x,y}} (b1)
(b1) edge node {\mcomp_{w,y,z}} (b2)
;
\end{tikzpicture}
\end{equation}
This means, first of all, that the two composites in \cref{psucat-associativity} are equal as functors.  Furthermore, their respective linearity constraints are equal.  To make this explicit, consider objects 
\begin{equation}\label{fghprime}
f,f' \in \C(y,z), \quad g,g' \in \C(x,y), \andspace h,h' \in \C(w,x).
\end{equation}
Linearity constraints of composite multilinear functors are defined in \cref{ffjlinearity}.  The equality of, respectively, the first, second, and third linearity constraints of the two composites in \cref{psucat-associativity} are the following three commutative diagrams in $\C(w,z)$.
\begin{equation}\label{psucat-con-i}

\end{equation}

The unity axiom \cref{menriched-cat-unity} of a $\psu$-category says that, for objects $x,y \in \C$, the following diagram of symmetric monoidal functors commutes.
\begin{equation}\label{psucat-unity}
%
\end{equation}
Both boundary composites in \cref{psucat-unity} are strictly unital by the unity axiom \cref{nlinearunity} of each bilinear functor $\mcomp$.  In terms of the monoidal constraints, the commutative diagram \cref{psucat-unity} means the following commutative diagram in $\C(x,y)$ for objects $g,g' \in \C(x,y)$, where we use concatenation to denote $\mcomp$ as in \cref{mxyzfg}.
\begin{equation}\label{psucat-unity-diag}
%
\end{equation}
This finishes the description of a $\psu$-category.
\end{explanation}

\begin{explanation}[$\psu$-Functors]\label{expl:perm-enr-functors}
Unpacking \cref{def:mfunctor} with $\M = \psu$, a $\psu$-functor 
\[F \cn (\C,\mcomp,i) \to (\D,\mcomp,i)\] 
between $\psu$-categories consists of 
\begin{itemize}
\item an object assignment $F \cn \Ob\C \to \Ob\D$ and
\item for each pair of objects $x,y \in \C$, a strictly unital symmetric monoidal functor (\cref{def:monoidalfunctor})
\[\big(F_{x,y} \scs F^2_{x,y} \scs F^0_{x,y} = 1\big) \cn \C(x,y) \to \D(Fx,Fy).\]
\end{itemize}

The compatibility of $F$ with composition \cref{mfunctor-diagrams} is the following commutative diagram of bilinear functors.
\begin{equation}\label{psufunctor-comp}

\end{equation}
This means, first of all, that \cref{psufunctor-comp} is a commutative diagram of functors.  Moreover, the equality of, respectively, the first and second linearity constraints of the two composites in \cref{psufunctor-comp} are the following two commutative diagrams in $\D(Fx,Fz)$ for objects $f,f' \in \C(y,z)$ and $g,g' \in \C(x,y)$, with $\mcomp$ shortened to concatenation as in \cref{mxyzfg}.
\begin{equation}\label{psufunctor-constraints}
%
\end{equation}

The compatibility of $F$ with identities \cref{mfunctor-diagrams} means the object equality
\begin{equation}\label{psufunctor-identities}
F_{x,x}(i_x) = i_{Fx} \inspace \D(Fx,Fx)
\end{equation}
for each object $x$ in $\C$.  This finishes the description of a $\psu$-functor.
\end{explanation}

\begin{explanation}[$\psu$-Natural Transformations]\label{expl:perm-nattr}
Unpacking \cref{def:mnaturaltr} with $\M = \psu$, a $\psu$-natural transformation
\[%
\]
between $\psu$-functors between $\psu$-categories consists of, for each object $x$ in $\C$, an object
\[\theta_x \in \D(Fx,Gx).\]
This is also regarded as a functor $\boldone \to \D(Fx,Gx)$.  

The naturality diagram \cref{mnaturality-diag} for a $\psu$-natural transformation is the following commutative diagram of symmetric monoidal functors for each pair of objects $x,y \in \C$.
\begin{equation}\label{psu-naturality}
%
\end{equation}
Both composites in \cref{psu-naturality} are strictly unital because 
\begin{itemize}
\item $F_{x,y}$ and $G_{x,y}$ are strictly unital and
\item both bilinear functors $\mcomp$ satisfy the unity axiom \cref{nlinearunity}.
\end{itemize} 
In addition to being a commutative diagram of functors, the equality of the monoidal constraints of the two composites in \cref{psu-naturality} means the following commutative diagram in $\D(Fx,Gy)$ for objects $g,g' \in \C(x,y)$, with $\mcomp$ denoted by concatenation as in \cref{mxyzfg}.
\begin{equation}\label{psu-nat-constraint}
\begin{tikzpicture}[vcenter]
\def\h{4} \def\g{.5} \def\v{.8} \def\u{-1.4}
\draw[0cell=.85]
(0,0) node (a1) {(Gg) \theta_x \oplus (Gg') \theta_x}
(a1)+(\g,\v) node (a2) {\theta_y (Fg) \oplus \theta_y (Fg')}
(a2)+(\h,0) node (a3) {\theta_y (Fg \oplus Fg')}
(a1)+(0,\u) node (b1) {(Gg \oplus Gg') \theta_x}
(b1)+(\g+\h,0) node (b2) {\big(G(g \oplus g')\big) \theta_x}
(b2)+(0,\v) node (b3) {\theta_y F(g \oplus g')}
;
\draw[1cell=.85]
(a1) edge[-,double equal sign distance] (a2)
(b2) edge[-,double equal sign distance] (b3)
(a2) edge node {\mcomp^2_2} (a3)
(a3) edge node {1_{\theta_y} F^2} (b3)
(a1) edge node[swap] {\mcomp^2_1} (b1)
(b1) edge node {G^2 1_{\theta_x}} (b2)
;
\end{tikzpicture}
\end{equation}
This finishes the description of a $\psu$-natural transformation.
\end{explanation}

\section{Self-Enrichment of the Multicategory of Permutative Categories}
\label{sec:perm-self-enr}

In this section we observe that the 2-category $\permcatsu$ (\cref{def:permcat}) of small permutative categories has the additional structure of a category enriched in the multicategory $\permcatsu$ (\cref{thm:permcatmulticat}) in the sense of \cref{def:menriched-cat}.  See \cref{permcat-selfenr}.  There are two ways to think about this result.
\begin{enumerate}
\item It is a permutative extension of the 2-category structure on $\permcatsu$.  So we have hom \emph{permutative} categories and composition \emph{bilinear} functors between permutative categories.  From this viewpoint, this section is an expanded version of \cite[Section 5]{bohmann_osorno-mackey}; see also \cite[Example 3.4]{guillou}.
\item \cref{permcat-selfenr} is a precursor of \cref{perm-closed-multicat}, which extends the $\Cat$-multicategory $\permcatsu$ and its self-enrichment to a \emph{closed} multicategory structure.  One may consider this section as a warm-up exercise for \cref{perm-closed-multicat}.  Each non-symmetric closed multicategory has a canonical self-enrichment (\cref{cl-multi-cl-cat}), so, in particular, $\permcatsu$ has a canonical self-enrichment.  As we discuss in more detail in \cref{ex:perm-closed-multicat}, the self-enrichment of $\permcatsu$ in \cref{permcat-selfenr} is the same as its canonical self-enrichment obtained from its closed multicategory structure.  However, a non-symmetric closed multicategory has more structure than its self-enrichment.
\end{enumerate}

This section is organized as follows.
\begin{itemize}
\item The hom permutative categories are constructed in \cref{def:permutative-homcat} and verified in \cref{psucd-hom-permcat}.
\item The composition bilinear functors are constructed in \cref{def:psu-bilinear-comp} and verified in \cref{psu-mBCD}.
\item The self-enrichment of $\permcatsu$ is constructed in \cref{def:perm-selfenr} and verified in \cref{permcat-selfenr}.
\end{itemize}
To simplify the presentation, we also use the shortened notation in \cref{permcatsu-psu}:
\[\psu = \permcatsu.\]
In a typical permutative category, the monoidal product, monoidal unit, and braiding are denoted $\oplus$, $\pu$, and $\xi$, respectively.  

\begin{definition}[Hom Permutative Categories]\label{def:permutative-homcat}
Given small permutative categories $\C$ and $\D$, we define a small permutative category\index{hom!permutative category}
\begin{equation}\label{psucd-permutative}
\big(\psu(\C,\D), \oplus, \pu, \xi\big)
\end{equation}
as follows.  The small category $\psu(\C,\D)$ is a hom category of the 2-category $\permcatsu$ in \cref{def:permcat}. 
\begin{itemize}
\item Its objects are strictly unital symmetric monoidal functors $\C \to \D$ (\cref{def:monoidalfunctor}).
\item Its morphisms are monoidal natural transformations (\cref{def:monoidalnattr}).
\item Identity morphisms and composition are those of monoidal natural transformations.
\end{itemize}

The monoidal product
\begin{equation}\label{psucd-oplus}
\oplus \cn \psu(\C,\D) \times \psu(\C,\D) \to \psu(\C,\D)
\end{equation}
is defined as follows.

\medskip
\emph{Objects}.  It sends a pair of strictly unital symmetric monoidal functors 
\[(F,F^2) \scs (G,G^2) \cn \C \to \D\]
to the following composite functor, with diag denoting the diagonal functor.
\begin{equation}\label{FoplusG-def}

\end{equation}
In other words, the monoidal product is defined pointwise in $\D$:
\begin{equation}\label{FoplusG-x}
(F \oplus G)(x) = Fx \oplus Gx \forspace x \in \C.
\end{equation}
The unit constraint and monoidal constraint of $F \oplus G$ are defined as follows.

\medskip
\emph{Unit Constraint}: It is the identity morphism 
\[1_\pu \cn \pu \to \pu = \pu \oplus \pu = (F \oplus G)(\pu) \inspace \D.\]  

\emph{Monoidal Constraint}: It is the following composite in $\D$ for objects $x,y \in \C$.
\begin{equation}\label{FoplusG-mon-constraint}
%
\end{equation}
The naturality of $(F \oplus G)^2$ follows from
\begin{itemize}
\item the naturality of the braiding $\xi$ of $\D$, $F^2$, and $G^2$, and
\item the functoriality of $\oplus$ for $\D$.
\end{itemize} 

\medskip
\emph{Morphisms}: For a pair of monoidal natural transformations
\[%
\]
between strictly unital symmetric monoidal functors, their monoidal product 
\[\theta \oplus \theta' \cn F \oplus F' \to G \oplus G'\] 
is the following whiskering.
\[%
\]
In other words, for each object $x \in \C$, the $x$-component of $\theta \oplus \theta'$ is
\begin{equation}\label{theta-oplus-thetap-x}
(\theta \oplus \theta')_x = \theta_x \oplus \theta'_x \cn Fx \oplus F'x \to Gx \oplus G'x.
\end{equation}
This defines a monoidal natural transformation for the following reasons.
\begin{itemize}
\item The naturality of $\theta \oplus \theta'$ follows from the naturality of $\theta$ and $\theta'$, together with the functoriality of $\oplus$ in $\D$.
\item $\theta \oplus \theta'$ satisfies the unity constraint axiom in \cref{monnattr} because, if $x = \pu$ in $\C$, then $\theta_\pu = 1_\pu = \theta'_\pu$.
\item $\theta \oplus \theta'$ satisfies the monoidal constraint axiom in \cref{monnattr} by
\begin{itemize}
\item the naturality of the braiding $\xi$ in $\D$ and
\item the monoidal constraint axiom \cref{monnattr} for $\theta$ and $\theta'$.
\end{itemize}
\end{itemize}
The functoriality of $\oplus$ for $\psu(\C,\D)$ follows from the functoriality of $\oplus$ for $\D$.

\medskip
\emph{Monoidal Unit}.  It is the constant functor 
\begin{equation}\label{psucd-unit}
\pu \cn \C \to \D
\end{equation}
at the monoidal unit of $\D$.  Its unit constraint and monoidal constraint are both given by the identity morphism $1_\pu$ in $\D$.  

\medskip
\emph{Braiding}.  For strictly unital symmetric monoidal functors $F$ and $G$ as in \cref{FoplusG-def}, the $(F,G)$-component of the braiding $\xi$ is the natural isomorphism
\begin{equation}\label{xiFG}
\begin{tikzpicture}[baseline={(a.base)}]
\def\t{22}
\draw[0cell]
(0,0) node (a) {\C}
(a)+(2.5,0) node (b) {\D}
;
\draw[1cell=.8]
(a) edge[bend left=\t] node {F \oplus G} (b)
(a) edge[bend right=\t] node[swap] {G \oplus F} (b)
;
\draw[2cell]
node[between=a and b at .38, rotate=-90, 2label={above,\xi_{F,G}}] {\Rightarrow}
;
\end{tikzpicture}
\end{equation}
with, for each object $x \in \C$, $x$-component given by the braiding
\begin{equation}\label{xiFxGx}
(F \oplus G)(x) = Fx \oplus Gx \fto[\iso]{\xi_{Fx,Gx}} Gx \oplus Fx = (G \oplus F)(x)
\end{equation}
in $\D$.  The naturality of the braiding in $\D$ implies the naturality of both
\begin{itemize}
\item $(\xi_{F,G})_x$ with respect to $x \in \C$ and
\item $\xi_{F,G}$ with respect to $F$ and $G$ in $\psu(\C,\D)$.
\end{itemize}  
This finishes the definition.
\end{definition}

\begin{lemma}\label{psucd-hom-permcat}
For small permutative categories $\C$ and $\D$, the quadruple in \cref{psucd-permutative} 
\[\big(\psu(\C,\D), \oplus, \pu, \xi\big)\]
is a small permutative category.
\end{lemma}

\begin{proof}
We already explained some of the required conditions in \cref{def:permutative-homcat}.  It remains to prove statements \cref{psucd-hom-i,psucd-hom-ii,psucd-hom-iii,psucd-hom-iv,psucd-hom-v} below.
\begin{romenumerate}
\item\label{psucd-hom-i} The pair defined in \cref{FoplusG-def,FoplusG-mon-constraint}
\[\big(F \oplus G, (F \oplus G)^2\big) \cn \C \to \D\]
is a strictly unital symmetric monoidal functor.
\item\label{psucd-hom-ii} The functor $\oplus$ in \cref{psucd-oplus} is associative.
\item\label{psucd-hom-iii} The constant functor $\pu$ in \cref{psucd-unit} is a strict two-sided unit for $\oplus$.
\item\label{psucd-hom-iv} The natural transformation $\xi_{F,G}$ in \cref{xiFG} is monoidal (\cref{def:monoidalnattr}).
\item\label{psucd-hom-v} $\psu(\C,\D)$ satisfies the symmetry and hexagon axioms \cref{symmoncatsymhexagon}.
\end{romenumerate}

\medskip
\emph{Statement \cref{psucd-hom-i}}.  We need to check the unity axiom \cref{monoidalfunctorunity}, the associativity axiom \cref{monoidalfunctorassoc}, and the braiding axiom \cref{monoidalfunctorbraiding} for $F \oplus G$.

Since its unit constraint is $1_{\pu}$, the unity axiom \cref{monoidalfunctorunity} for $F \oplus G$ means the equalities
\[(F \oplus G)^2_{\pu,?} = 1_{F? \oplus G?} = (F \oplus G)^2_{?,\pu} \inspace \D.\]
These equalities follow from the following equalities in $\D$.
\[\xi_{\pu,?} = 1_? = \xi_{?,\pu}\qquad 
F^2_{\pu,?} = 1_{F?} = F^2_{?,\pu}\qquad
G^2_{\pu,?} = 1_{G?} = G^2_{?,\pu}\]

For objects $x,y,z \in \C$, the associativity diagram \cref{monoidalfunctorassoc} for $F \oplus G$ is the boundary of the following diagram in $\D$.
\[\begin{tikzpicture}
\def\h{3.5} \def\v{-1.4} \def\u{-1} \def\t{15} \def\s{10} \def\q{30}
\draw[0cell=.8]
(0,0) node (a) {Fx \oplus Gx \oplus Fy \oplus Gy \oplus Fz \oplus Gz}
(a)+(-\h,\u) node (b1) {Fx \oplus Fy \oplus Gx \oplus Gy \oplus Fz \oplus Gz}
(a)+(\h,\u) node (b2) {Fx \oplus Gx \oplus Fy \oplus Fz \oplus Gy \oplus Gz}
(b1)+(0,\v) node (c1) {F(x \oplus y) \oplus G(x \oplus y) \oplus Fz \oplus Gz}
(b2)+(0,\v) node (c2) {Fx \oplus Gx \oplus F(y \oplus z) \oplus G(y \oplus z)}
(c1)+(0,\v) node (d1) {F(x \oplus y) \oplus Fz \oplus G(x \oplus y) \oplus Gz}
(c2)+(0,\v) node (d2) {Fx \oplus F(y \oplus z) \oplus Gx \oplus G(y \oplus z)}
(d1)+(\h,\u) node (e) {F(x \oplus y \oplus z) \oplus G(x \oplus y \oplus z)}
node[between=b1 and d2 at .5, shift={(0,.6)}] (f) {Fx \oplus Fy \oplus Fz \oplus Gx \oplus Gy \oplus Gz}
;
\draw[1cell=.8]
(a) edge[bend right=\t] node[swap] {1 \oplus \xi \oplus 1} (b1)
(a) edge[bend left=\t] node {1 \oplus \xi \oplus 1} (b2)
(b1) edge node[swap] {F^2 \oplus G^2 \oplus 1} (c1)
(b2) edge node {1 \oplus F^2 \oplus G^2} (c2)
(c1) edge node[swap] {1 \oplus \xi \oplus 1} (d1)
(c2) edge node {1 \oplus \xi \oplus 1} (d2)
(d1) edge[bend right=\t] node[swap,pos=.3] {F^2 \oplus G^2} (e)
(d2) edge[bend left=\t] node[pos=.3] {F^2 \oplus G^2} (e)
(b1) edge[bend left=\s] node[pos=.4] {1 \oplus \xi \oplus 1} (f)
(b2) edge[bend right=\s] node[swap,pos=.4] {1 \oplus \xi \oplus 1} (f)
;
\draw[1cell=.75]
(f) edge[bend left=\q] node[swap] {F^2 \oplus 1 \oplus G^2 \oplus 1} (d1)
(f) edge[bend right=\q] node {1 \oplus F^2 \oplus 1 \oplus G^2} (d2)
;
\end{tikzpicture}\]
The following statements hold for the diagram above.
\begin{itemize}
\item The top sub-region commutes by the coherence theorem for symmetric monoidal categories \cite[XI.1 Theorem 1]{maclane}.
\item The left and right sub-regions commute by the naturality of the braiding $\xi$ in $\D$.
\item The bottom sub-region commutes by the axiom \cref{monoidalfunctorassoc} for $F$ and $G$.
\end{itemize}
This proves the associativity axiom \cref{monoidalfunctorassoc} for $F \oplus G$.

For objects $x,y \in \C$, the braiding diagram \cref{monoidalfunctorbraiding} for $F \oplus G$ is the boundary of the following diagram in $\D$.
\[\begin{tikzpicture}
\def\v{-1.3} \def\u{-.8}
\draw[0cell=.85]
(0,0) node (a1) {(F \oplus G)(x) \oplus (F \oplus G)(y)}
(a1)+(5,0) node (a2) {(F \oplus G)(y) \oplus (F \oplus G)(x)}
(a1)+(0,\u) node (b1) {Fx \oplus Gx \oplus Fy \oplus Gy}
(a2)+(0,\u) node (b2) {Fy \oplus Gy \oplus Fx \oplus Gx}
(b1)+(0,\v) node (c1) {Fx \oplus Fy \oplus Gx \oplus Gy}
(b2)+(0,\v) node (c2) {Fy \oplus Fx \oplus Gy \oplus Gx}
(c1)+(0,\v) node (d1) {F(x \oplus y) \oplus G(x \oplus y)}
(c2)+(0,\v) node (d2) {F(y \oplus x) \oplus G(y \oplus x)}
;
\draw[1cell=.85]
(a1) edge[-,double equal sign distance] (b1)
(a2) edge[-,double equal sign distance] (b2)
(b1) edge node[swap] {1 \oplus \xi \oplus 1} (c1)
(b2) edge node {1 \oplus \xi \oplus 1} (c2)
(c1) edge node[swap] {F^2 \oplus G^2} (d1)
(c2) edge node {F^2 \oplus G^2} (d2)
(b1) edge node {\xi} (b2)
(c1) edge node {\xi \oplus \xi} (c2)
(d1) edge node {F\xi \oplus G\xi} (d2)
;
\end{tikzpicture}\]
The following statements hold for the diagram above.
\begin{itemize}
\item The top rectangle commutes by the coherence theorem for symmetric monoidal categories \cite[XI.1 Theorem 1]{maclane}.
\item The bottom rectangle commutes by the axiom \cref{monoidalfunctorbraiding} for $F$ and $G$.
\end{itemize}
This proves the braiding axiom \cref{monoidalfunctorbraiding} for $F \oplus G$.

\medskip
\emph{Statement \cref{psucd-hom-ii}}.  To verify that $\oplus$ in \cref{psucd-oplus} is associative on objects, we consider strictly unital symmetric monoidal functors
\[(F,F^2) \scs (G,G^2) \scs (H,H^2) \cn \C \to \D.\]
The definition \cref{FoplusG-x} of $(F \oplus G)(x)$ and the associativity of $\oplus$ in $\D$ imply that $(F \oplus G) \oplus H$ is equal to $F \oplus (G \oplus H)$ as functors.  Moreover, each of them has unit constraint given by $1_\pu$.  To show that their monoidal constraints are the same, we consider the following diagram in $\D$ for objects $x,y \in \C$.
\[\begin{tikzpicture}
\def\g{.5} \def\h{5.5} \def\v{-1.3} \def\u{.8}
\draw[0cell=.8]
(0,0) node (a1) {(Fx \oplus Gx) \oplus Hx \oplus (Fy \oplus Gy) \oplus Hy}
(a1)+(\g,\u) node (a2) {Fx \oplus (Gx \oplus Hx) \oplus Fy \oplus (Gy \oplus Hy)}
(a2)+(\h,0) node (a3) {Fx \oplus Fy \oplus Gx \oplus (Hx \oplus Gy) \oplus Hy}
(a3)+(0,\v-\u) node (a4) {Fx \oplus Fy \oplus Gx \oplus Gy \oplus Hx \oplus Hy}
(a1)+(0,\v) node (b1) {Fx \oplus (Gx \oplus Fy) \oplus Gy \oplus Hx \oplus Hy}
(a4)+(0,\v) node (b3) {F(x \oplus y) \oplus G(x \oplus y) \oplus H(x \oplus y)}
;
\draw[1cell=.8]
(a1) edge[-,double equal sign distance] (a2)
(a2) edge node {1 \oplus \xi \oplus 1} (a3)
(a3) edge node {1 \oplus \xi \oplus 1} (a4)
(a4) edge node {F^2 \oplus G^2 \oplus H^2} (b3)
(a1) edge node[swap] {1 \oplus \xi \oplus 1} (b1)
(b1) edge node {1 \oplus \xi \oplus 1} (a4)
;
\end{tikzpicture}\]
The following statements hold for the diagram above.
\begin{itemize}
\item The left-bottom composite is the monoidal constraint $\big((F \oplus G) \oplus H\big)^2_{x,y}$.
\item The top-right composite is the monoidal constraint $\big(F \oplus (G \oplus H)\big)^2_{x,y}$.
\item The top sub-region commutes by the coherence theorem for symmetric monoidal categories \cite[XI.1 Theorem 1]{maclane}.
\end{itemize}
This shows that $\oplus$ in \cref{psucd-oplus} is associative on objects.  Its associativity on morphisms follows from the definition \cref{theta-oplus-thetap-x} of $(\theta \oplus \theta')_x$ and the associativity of $\oplus$ in $\D$.  Thus the monoidal product of $\psu(\C,\D)$ is associative.

\medskip
\emph{Statement \cref{psucd-hom-iii}}.  
The constant functor $\pu$ in \cref{psucd-unit}, with $1_\pu$ as the unit and monoidal constraints, is a strict two-sided unit for $\oplus$ in $\psu(\C,\D)$ for the following reasons.
\begin{itemize}
\item $F \oplus \pu$ and $\pu \oplus F$ are both equal to $F$ as functors because the monoidal unit $\pu$ in $\D$ is strict.
\item The equalities \cref{braidedunity}
\[\xi_{\pu,?} = 1_? = \xi_{?,\pu} \inspace \D\]
imply that for $F \oplus \pu$ and $\pu \oplus F$, the morphism $1 \oplus \xi \oplus 1$ in \cref{FoplusG-mon-constraint} is the identity.  The second morphism in \cref{FoplusG-mon-constraint} is $F^2_{x,y} \oplus 1_\pu$ or $1_\pu \oplus F^2_{x,y}$, which are both equal to $F^2_{x,y}$ by the strict unity of $\pu$ in $\D$.
\end{itemize}
This shows that $\big(\psu(\C,\D),\oplus,\pu\big)$ is a strict monoidal category.

\medskip
\emph{Statement \cref{psucd-hom-iv}}.  
The natural transformation $\xi_{F,G}$ in \cref{xiFG} satisfies the unity axiom in \cref{monnattr} because its $\pu$-component is 
\[\xi_{F\pu,G\pu} = \xi_{\pu,\pu} = 1_\pu\]
by either unity properties in \cref{braidedunity}.   

To show that $\xi_{F,G}$ is compatible with the monoidal constraints of its domain and codomain in the sense of \cref{monnattr}, we consider the following diagram in $\D$ for objects $x,y \in \C$.
\[\begin{tikzpicture}
\def\v{-1.3} \def\u{-.8}
\draw[0cell=.85]
(0,0) node (a1) {(F \oplus G)(x) \oplus (F \oplus G)(y)}
(a1)+(5,0) node (a2) {(G \oplus F)(x) \oplus (G \oplus F)(y)}
(a1)+(0,\u) node (b1) {Fx \oplus Gx \oplus Fy \oplus Gy}
(a2)+(0,\u) node (b2) {Gx \oplus Fx \oplus Gy \oplus Fy}
(b1)+(0,\v) node (c1) {Fx \oplus Fy \oplus Gx \oplus Gy}
(b2)+(0,\v) node (c2) {Gx \oplus Gy \oplus Fx \oplus Fy}
(c1)+(0,\v) node (d1) {F(x \oplus y) \oplus G(x \oplus y)}
(c2)+(0,\v) node (d2) {G(x \oplus y) \oplus F(x \oplus y)}
;
\draw[1cell=.85]
(a1) edge[-,double equal sign distance] (b1)
(a2) edge[-,double equal sign distance] (b2)
(b1) edge node[swap] {1 \oplus \xi \oplus 1} (c1)
(b2) edge node {1 \oplus \xi \oplus 1} (c2)
(c1) edge node[swap] {F^2 \oplus G^2} (d1)
(c2) edge node {G^2 \oplus F^2} (d2)
(b1) edge node {\xi \oplus \xi} (b2)
(c1) edge node {\xi} (c2)
(d1) edge node {\xi} (d2)
;
\end{tikzpicture}\]
The following statements hold for the diagram above.
\begin{itemize}
\item The left vertical composite is the monoidal constraint $(F \oplus G)^2_{x,y}$.
\item The right vertical composite is the monoidal constraint $(G \oplus F)^2_{x,y}$.
\item The top rectangle commutes by the coherence theorem for symmetric monoidal categories \cite[XI.1 Theorem 1]{maclane}.
\item The bottom rectangle commutes by the naturality of the braiding in $\D$.
\end{itemize}
This shows that $\xi_{F,G}$ is a monoidal natural transformation (\cref{def:monoidalnattr}).

\medskip
\emph{Statement \cref{psucd-hom-v}}.  The symmetry and hexagon axioms \cref{symmoncatsymhexagon} hold in $\psu(\C,\D)$ by 
\begin{itemize}
\item the componentwise definitions \cref{FoplusG-x,theta-oplus-thetap-x,xiFxGx}, and
\item the corresponding axioms in $\D$.
\end{itemize}
This finishes the proof.
\end{proof}

From now on $\psu(\C,\D)$ is a permutative category as in \cref{psucd-hom-permcat}.

\begin{definition}[Composition]\label{def:psu-bilinear-comp}
For small permutative categories $\B$, $\C$, and $\D$, we define the data of a bilinear functor (\cref{def:nlinearfunctor})
\begin{equation}\label{mBCD}
\mcomp_{\B,\C,\D} \cn \psu(\C,\D) \times \psu(\B,\C) \to \psu(\B,\D)
\end{equation}
as follows, where we abbreviate $\mcomp_{\B,\C,\D}$ to $\mcomp$.
\begin{description}
\item[Objects and Morphisms] The underlying functor of $\mcomp$ is given by
\begin{itemize}
\item composition of strictly unital symmetric monoidal functors on objects and
\item horizontal composition of monoidal natural transformations on morphisms,
\end{itemize} 
as displayed below.
\begin{equation}\label{FGHI-theta-psi}
\begin{tikzpicture}[baseline={(b.base)}]
\def\a{.7} \def\h{1.8} \def\t{25} \def\s{22}
\draw[0cell=.9]
(0,0) node (b) {\B}
(b)+(\h,0) node (c) {\C}
(c)+(\h,0) node (d) {\D}
(d)+(\a,0) node (m1) {}
(m1)+(.8,0) node (m2) {}
(m2)+(\a,0) node (b') {\B}
(b')+(2,0) node (d') {\D}
;
\draw[1cell=.8]
(b) edge[bend left=\t] node {F} (c)
(b) edge[bend right=\t] node[swap] {G} (c)
(c) edge[bend left=\t] node {H} (d)
(c) edge[bend right=\t] node[swap] {I} (d)
(b') edge[bend left=\t] node {HF} (d')
(b') edge[bend right=\t] node[swap] {IG} (d')
;
\draw[1cell]
(m1) edge[|->] node {\mcomp} (m2)
;
\draw[2cell=.9]
node[between=b and c at .45, rotate=-90, 2label={above,\theta}] {\Rightarrow}
node[between=c and d at .43, rotate=-90, 2label={above,\psi}] {\Rightarrow}
node[between=b' and d' at .35, rotate=-90, 2label={above,\psi * \theta}] {\Rightarrow}
;
\end{tikzpicture}
\end{equation}
This is part of the 2-category structure on $\permcatsu$ (\cref{def:permcat}).
\item[First Monoidal Constraint] It is the identity natural transformation
\begin{equation}\label{HFIF}
HF \oplus IF \fto{\mcomp^2_1 = 1} (H \oplus I)F.
\end{equation}
Thus for each object $x \in \B$, its $x$-component is the identity morphism
\[HFx \oplus IFx \fto{1} (H \oplus I)(Fx).\]
\item[Second Monoidal Constraint]
It is the natural transformation 
\begin{equation}\label{HFHG}
HF \oplus HG \fto{\mcomp^2_2} H(F \oplus G)
\end{equation}
with, for each object $x \in \B$, $x$-component given by the monoidal constraint
\[HFx \oplus HGx \fto{H^2_{Fx,Gx}} H(Fx \oplus Gx).\]
The naturality of $\mcomp^2_2$ with respect to $x \in \B$ follows from the naturality of $H^2$.
\end{description}
This finishes the definition of $\mcomp_{\B,\C,\D}$.
\end{definition}

\begin{lemma}\label{psu-mBCD}
In the context of \cref{def:psu-bilinear-comp}, 
\[\big(\mcomp_{\B,\C,\D} \scs \mcomp^2_1=1 \scs \mcomp^2_2\big) \cn \psu(\C,\D) \times \psu(\B,\C) \to \psu(\B,\D)\]
is a bilinear functor.
\end{lemma}

\begin{proof}
We prove statements \cref{psu-m-i,psu-m-ii,psu-m-iii,psu-m-iv,psu-m-v} below.
\begin{romenumerate}
\item\label{psu-m-i} The natural transformation $\mcomp^2_1$ in \cref{HFIF} is monoidal.
\item\label{psu-m-ii} The natural transformation $\mcomp^2_2$ in \cref{HFHG} is monoidal.
\item\label{psu-m-iii} $\mcomp^2_1$ is natural in $F$, $H$, and $I$.
\item\label{psu-m-iv} $\mcomp^2_2$ is natural in $F$, $G$, and $H$.
\item\label{psu-m-v} $(\mcomp,\mcomp^2_1,\mcomp^2_2)$ satisfies the axioms of a bilinear functor (\cref{def:nlinearfunctor}).
\end{romenumerate}

\medskip
\emph{Statement \cref{psu-m-i}}.  
To check that the natural transformation $\mcomp^2_1$ in \cref{HFIF} is monoidal (\cref{def:monoidalnattr}), first recall that each of its components is an identity morphism.  Thus $\mcomp^2_1$ satisfies the unit axiom \cref{monnattr} because its $\pu$-component is $1_\pu$.  Moreover, for objects $x,y \in \B$, the $(x,y)$-component of the monoidal constraint of each of $HF \oplus IF$ and $(H \oplus I)F$ is the following composite in $\D$.
\[\begin{tikzpicture}
\def\g{0} \def\v{1.3}
\draw[0cell=.9]
(0,0) node (a) {HFx \oplus IFx \oplus HFy \oplus IFy}
(a)+(\g,-\v) node (b) {HFx \oplus HFy \oplus IFx \oplus IFy}
(b)+(5.5,0) node (c) {H(Fx \oplus Fy) \oplus I(Fx \oplus Fy)}
(c)+(\g,\v) node (d) {HF(x \oplus y) \oplus IF(x \oplus y)}
;
\draw[1cell=.9]
(a) edge node[swap] {1 \oplus \xi \oplus 1} (b)
(b) edge node {H^2 \oplus I^2} (c)
(c) edge node[swap] {H(F^2) \oplus I(F^2)} (d)
;
\end{tikzpicture}\]
This shows that $\mcomp^2_1$ in \cref{HFIF} is a monoidal natural transformation.

\medskip
\emph{Statement \cref{psu-m-ii}}.  
To check that the natural transformation $\mcomp^2_2$ in \cref{HFHG} is monoidal, first note that the unity axioms \cref{monoidalfunctorunity} of the strictly unital symmetric monoidal functor $(H,H^2,H^0 = 1_\pu)$ imply the equalities
\begin{equation}\label{Htwopu}
H^2_{\pu,?} = 1_{H?} = H^2_{?,\pu}.
\end{equation}
Since the domain and codomain of $\mcomp^2_2$ are strictly unital, the unit axiom \cref{monnattr} for $\mcomp^2_2$ is the equality
\[H^2_{F\pu,G\pu} = 1_\pu,\]
which holds by \cref{Htwopu} and the strict unity of $F$, $G$, and $H$.

Next we check the compatibility of $\mcomp^2_2$ with the monoidal constraints of its domain and codomain \cref{monnattr}.  For objects $x,y \in \B$, we use the notation
\[\begin{split}
X_1 &= HFx \oplus H(Gx \oplus Fy) \oplus HGy,\\
X_2 &= H(Fx \oplus Gx \oplus Fy) \oplus HGy,\\
Y_1 &= HFx \oplus H(Fy \oplus Gx) \oplus HGy, \andspace\\
Y_2 &= H(Fx \oplus Fy \oplus Gx) \oplus HGy,
\end{split}\]
and consider the following diagram in $\D$.
\[\begin{tikzpicture}
\def\h{7} \def\v{-2.7} 
\draw[0cell=.8]
(0,0) node (a1) {HFx \oplus HGx \oplus HFy \oplus HGy}
(a1)+(\h,0) node (a2) {H(Fx \oplus Gx) \oplus H(Fy \oplus Gy)}
(a1)+(0,\v) node (b1) {HFx \oplus HFy \oplus HGx \oplus HGy}
(a2)+(0,\v) node (b2) {H(Fx \oplus Gx \oplus Fy \oplus Gy)}
(b1)+(0,\v) node (c1) {H(Fx \oplus Fy) \oplus H(Gx \oplus Gy)}
(b2)+(0,\v) node (c2) {H(Fx \oplus Fy \oplus Gx \oplus Gy)}
(c1)+(0,-1.3) node (d1) {HF(x \oplus y) \oplus HG(x \oplus y)}
(c2)+(0,-1.3) node (d2) {H\big(F(x \oplus y) \oplus G(x \oplus y)\big)}
node[between=a1 and b2 at .4] (x1) {X_1}
node[between=a1 and b2 at .67] (x2) {X_2}
(x1)+(0,\v) node (y1) {Y_1}
(x2)+(0,\v) node (y2) {Y_2}
node[between=b1 and x1 at .25, shift={(0,.4)}] () {\spadesuit}
node[between=y1 and x2 at .5] () {\clubsuit}
node[between=y2 and b2 at .25] () {\clubsuit}
node[between=c1 and d2 at .5, shift={(0,.3)}] () {\clubsuit}
;
\draw[1cell=.8]
(a1) edge node[swap] {1 \oplus \xi \oplus 1} (b1)
(b1) edge node[swap] {H^2 \oplus H^2} (c1)
(c1) edge node[swap] {H(F^2) \oplus H(G^2)} (d1)
(d1) edge node {H^2} (d2)
(a1) edge node {H^2 \oplus H^2} (a2)
(a2) edge node {H^2} (b2)
(b2) edge node {H(1 \oplus \xi \oplus 1)} (c2)
(c2) edge node {H(F^2 \oplus G^2)} (d2)
(c1) edge node {H^2} (c2)
(x1) edge node[swap,pos=.25] {1 \oplus H(\xi) \oplus 1} (y1)
(x2) edge node[pos=.6] {H(1 \oplus \xi) \oplus 1} (y2)
(a1) edge node[swap,pos=.7] {1 \oplus H^2 \oplus 1} (x1)
(x1) edge node[pos=.2] {H^2 \oplus 1} (x2)
(x2) edge node[pos=.4] {H^2} (b2)
(b1) edge node[swap,pos=.7] {1 \oplus H^2 \oplus 1} (y1)
(y1) edge node[pos=.2] {H^2 \oplus 1} (y2)
(y2) edge node[pos=.4] {H^2} (c2)
;
\end{tikzpicture}\]
The following statements hold for the diagram above.
\begin{itemize}
\item The left vertical composite is the $(x,y)$-component of the monoidal constraint of $HF \oplus HG$, which is the domain of $\mcomp^2_2$.
\item The right vertical composite is the $(x,y)$-component of the monoidal constraint of $H(F \oplus G)$, which is the codomain of $\mcomp^2_2$.
\item The sub-region labeled $\spadesuit$ is commutative by the compatibility of $H$ with the braiding \cref{monoidalfunctorbraiding}.
\item The three sub-regions labeled $\clubsuit$ are commutative by the naturality of $H^2$.
\item The two unlabeled triangles are commutative by the associativity of $H^2$ \cref{monoidalfunctorassoc}.
\end{itemize}
This shows that $\mcomp^2_2$ in \cref{HFHG} is a monoidal natural transformation.

\medskip
\emph{Statement \cref{psu-m-iii}}.  To show that $\mcomp^2_1$ is natural in $F$, $H$, and $I$, we consider monoidal natural transformations
\[
\]
between strictly unital symmetric monoidal functors.  The corresponding naturality diagram for $\mcomp^2_1$ commutes because, for each object $x \in \B$, both $x$-components
\[\big((\psi * \theta) \oplus (\vphi * \theta)\big)_x \andspace \big((\psi \oplus \vphi) * \theta\big)_x\]
are given by the composite
\[%
\]
in $\D$.

\medskip
\emph{Statement \cref{psu-m-iv}}.  To show that $\mcomp^2_2$ is natural in $F$, $G$, and $H$, we consider monoidal natural transformations
\[%
\]
between strictly unital symmetric monoidal functors and the corresponding naturality diagram below.
\[%
\]
For each object $x \in \B$, the $x$-component of the diagram above is the boundary of the following diagram in $\D$.
\[%
\]
The top rectangle commutes by the naturality of $H^2$.  The bottom rectangle commutes by the monoidality of $\psi$ \cref{monnattr}.  This shows that $\mcomp^2_2$ is natural in $F$, $G$, and $H$.

\medskip
\emph{Statement \cref{psu-m-v}}.  Next we check the axioms of a bilinear functor for $(\mcomp,\mcomp^2_1,\mcomp^2_2)$.

\medskip
\emph{Unity \cref{nlinearunity}}.  It holds for $\mcomp$ in \cref{mBCD} for the following reasons. 
\begin{itemize}
\item The monoidal unit in each hom permutative category $\psu(\C,\D)$ is the constant functor at the monoidal unit of $\D$ \cref{psucd-unit}.
\item Objects in $\psu(\C,\D)$ are \emph{strictly unital} symmetric monoidal functors.
\item Each morphism in $\psu(\C,\D)$ has $\pu$-component given by $1_\pu$ by the unit axiom of a monoidal natural transformation \cref{monnattr}.
\end{itemize}

\medskip
\emph{Constraint Unity \cref{constraintunity}}.  It holds for $\mcomp^2_1$ in \cref{HFIF} because each of its components is an identity morphism.  For the same reason, $\mcomp^2_1$ also satisfies the next two axioms, \cref{eq:ml-f2-assoc,eq:ml-f2-symm}.

The constraint unity axiom holds for $\mcomp^2_2$ in \cref{HFHG} for the following reasons.
\begin{itemize}
\item If $H$ is the monoidal unit $\pu$ \cref{psucd-unit}, then its monoidal constraint $H^2$ is $1_\pu$ by definition.
\item If either $F$ or $G$ is the monoidal unit $\pu$, then \cref{Htwopu} yields the desired equalities 
\[H^2_{\pu,Gx} = 1_{HGx} \andspace H^2_{Fx,\pu} = 1_{HFx}.\]
\end{itemize}

\medskip
\emph{Constraint Associativity \cref{eq:ml-f2-assoc} and Symmetry \cref{eq:ml-f2-symm}}.  They hold for $\mcomp^2_2$ by the associativity \cref{monoidalfunctorassoc} and symmetry \cref{monoidalfunctorbraiding} of the monoidal constraint $H^2$.

\medskip
\emph{Constraint 2-By-2 \cref{eq:f2-2by2}}.  We consider strictly unital symmetric monoidal functors as follows.
\[
\]
In the constraint 2-by-2 axiom, if $i=1$ and $k=2$, then the diagram \cref{eq:f2-2by2} is the left pentagon below, which commutes by the definition of $(H \oplus I)^2$ in \cref{FoplusG-mon-constraint}.
\[%
\]
If $i=2$ and $k=1$, then the constraint 2-by-2 diagram \cref{eq:f2-2by2} is the right pentagon above.  For each object $x \in \B$, the $x$-component of this pentagon is the following diagram in $\D$.
\[%
\]
This diagram commutes by the symmetry axiom \cref{symmoncatsymhexagon} and the functoriality of $\oplus$ in the permutative category $\D$.

This finishes the proof that $(\mcomp,\mcomp^2_1,\mcomp^2_2)$ is a bilinear functor.
\end{proof}

Recall the description of a $\permcatsu$-category in \cref{expl:perm-enr-cat}.

\begin{definition}\label{def:perm-selfenr}\index{permutative category!self-enrichment of multicategory of}\index{category!permutative!self-enrichment of the multicategory of}
The \emph{self-enrichment} of $\psu$, which we also denote by $\psu$, is the $\psu$-category defined as follows.  \cref{permcat-selfenr} verifies that $\psu$ is a $\psu$-category.
\begin{description}
\item[Objects] The objects are small permutative categories (\cref{def:symmoncat}).
\item[Hom Permutative Categories] For small permutative categories $\C$ and $\D$, the hom permutative category $\psu(\C,\D)$ is the one in \cref{psucd-hom-permcat}.
\item[Composition] For small permutative categories $\B$, $\C$, and $\D$, the composition bilinear functor
\[\big(\mcomp_{\B,\C,\D} \scs \mcomp^2_1, \mcomp^2_2\big) \cn \psu(\C,\D) \times \psu(\B,\C) \to \psu(\B,\D)\]
is the one in \cref{psu-mBCD}.
\item[Identities] Each small permutative category $\B$ is equipped with the identity symmetric monoidal functor $1_\B$, which is also regarded as an object in $\psu(\B,\B)$.
\end{description}
This finishes the definition of the self-enrichment of $\psu$.
\end{definition}

\begin{theorem}\label{permcat-selfenr}
Equipped with the self-enrichment in \cref{def:perm-selfenr}, $\permcatsu$ is a $\permcatsu$-category.
\end{theorem}

\begin{proof}
We need to check that $\psu$ satisfies the associativity axiom \cref{menriched-cat-assoc} and the unity axiom \cref{menriched-cat-unity} of a $\psu$-category.

\medskip
\emph{Associativity \cref{menriched-cat-assoc}}.  For $\psu$ this axiom means the commutativity of the  following diagram of composite 3-linear functors for small permutative categories $\A$, $\B$, $\C$, and $\D$.
\begin{equation}\label{psuself-assoc}
\begin{tikzpicture}[vcenter]
\def\v{-1.4}
\draw[0cell=.85]
(0,0) node (a1) {\psu(\C,\D) \times \psu(\B,\C) \times \psu(\A,\B)}
(a1)+(5.5,0) node (a2) {\psu(\B,\D) \times \psu(\A,\B)}
(a1)+(0,\v) node (b1) {\psu(\C,\D) \times \psu(\A,\C)}
(a2)+(0,\v) node (b2) {\psu(\A,\D)}
;
\draw[1cell=.9]
(a1) edge node {\mcomp_{\B,\C,\D} \times 1} (a2)
(a2) edge node {\mcomp_{\A,\B,\D}} (b2)
(a1) edge node[swap] {1 \times \mcomp_{\A,\B,\C}} (b1)
(b1) edge node {\mcomp_{\A,\C,\D}} (b2)
;
\end{tikzpicture}
\end{equation}
As a diagram of functors, \cref{psuself-assoc} is commutative because 
\begin{itemize}
\item composition of strictly unital symmetric monoidal functors and
\item horizontal composition of monoidal natural transformations
\end{itemize}
are both associative.  

To show that the two composites in \cref{psuself-assoc} have the same linearity constraints, as in \cref{psucat-con-i,psucat-con-ii-iii}, we consider strictly unital symmetric monoidal functors as follows.
\[\begin{tikzpicture}
\def\h{2}
\draw[0cell]
(0,0) node (a) {\A}
(a)+(\h,0) node (b) {\B}
(b)+(\h,0) node (c) {\C}
(c)+(\h,0) node (d) {\D}
;
\draw[1cell=.85]
(a) edge[transform canvas={yshift=.5ex}] node {H} (b)
(a) edge[transform canvas={yshift=-.5ex}] node[swap] {H'} (b)
(b) edge[transform canvas={yshift=.5ex}] node {G} (c)
(b) edge[transform canvas={yshift=-.5ex}] node[swap] {G'} (c)
(c) edge[transform canvas={yshift=.5ex}] node {F} (d)
(c) edge[transform canvas={yshift=-.5ex}] node[swap] {F'} (d)
;
\end{tikzpicture}\]
We consider the three diagrams in \cref{psucat-con-i,psucat-con-ii-iii} in the current context of $\psu$.
\begin{itemize}
\item The diagram for the first linearity constraint \cref{psucat-con-i} is commutative because each arrow is the identity, since $\mcomp^2_1 = 1$ by definition \cref{HFIF}.
\item The diagram for the second linearity constraint (= left diagram in \cref{psucat-con-ii-iii}) is commutative because $\mcomp^2_1 = 1$ and each $\mcomp^2_2$ is given by the monoidal constraint $F^2$ of $F$.
\item The diagram for the third linearity constraint (= right diagram in \cref{psucat-con-ii-iii}) is as follows.
\[\begin{tikzpicture}
\def\h{3.8} \def\g{.5} \def\v{.8} \def\u{-1.3}
\draw[0cell=.85]
(0,0) node (a1) {F(GH) \oplus F(GH')}
(a1)+(\g,\v) node (a2) {(FG)H \oplus (FG)H'}
(a2)+(\h,0) node (a3) {(FG)(H \oplus H')}
(a1)+(0,\u) node (b1) {F(GH \oplus GH')}
(b1)+(\g+\h,0) node (b2) {F\big(G(H \oplus H')\big)}
;
\draw[1cell=.85]
(a1) edge[-,double equal sign distance] (a2)
(a3) edge[-,double equal sign distance] (b2)
(a2) edge node {(FG)^2} (a3)
(a1) edge node[swap] {F^2} (b1)
(b1) edge node {F(G^2)} (b2)
;
\end{tikzpicture}\]
This diagram commutes by the definition of the monoidal constraint $(FG)^2$ (\cref{def:mfunctor-comp}).
\end{itemize}
This proves the associativity axiom \cref{menriched-cat-assoc} for $\psu$.

\medskip
\emph{Unity \cref{menriched-cat-unity}}.  In the current context of $\psu$, the unity diagram \cref{psucat-unity} commutes as a diagram of functors because, for each small permutative category $\B$, the identity $i_\B$ is the identity symmetric monoidal functor $1_\B$ (\cref{def:perm-selfenr}).  Moreover, the diagram \cref{psucat-unity-diag} for the monoidal constraints is commutative for the following reasons.
\begin{itemize}
\item The left half of \cref{psucat-unity-diag} is commutative because $\mcomp^2_1 = 1$.
\item The right half of \cref{psucat-unity-diag} is commutative because the monoidal constraint of $1_\B$ is the identity.
\end{itemize}
This finishes the proof that $\psu$ is a $\psu$-category.
\end{proof}

\section{Bilinear Evaluation for Permutative Categories}
\label{sec:eval-perm}

In this section we discuss a bilinear evaluation for permutative categories.  The bilinear evaluation in \cref{evCD} below is an analog of the evaluation in a symmetric monoidal closed category \cref{evaluation}.  Along with the self-enrichment of $\permcatsu$ (\cref{permcat-selfenr}), the bilinear evaluation in this section is also a part of the \emph{closed} multicategory structure on $\permcatsu$, which we discuss in \cref{ch:gspectra}.  We will use the bilinear evaluation in \cref{expl:Kemse,expl:Endmse,expl:Endstse} below.

This section is organized as follows.
\begin{itemize}
\item The evaluation bilinear functor for small permutative categories is constructed in \cref{def:perm-evaluation} and verified in \cref{ev-bilinear}.
\item \cref{ev-comp} shows that evaluation is compatible with the composition bilinear functor $\mcomp$ in \cref{psu-mBCD}.
\end{itemize}
Recall \index{functor!multilinear}\index{multilinear functor}\index{functor!bilinear}\index{bilinear functor}\emph{bilinear functors} between permutative categories (\cref{def:nlinearfunctor}) and the permutative category $\permcatsu(\C,\D)$ in \cref{psucd-hom-permcat}.  To simplify the presentation, we use the shortened notation in \cref{permcatsu-psu}:
\[\psu = \permcatsu.\]

\begin{definition}\label{def:perm-evaluation}
For small permutative categories $\C$ and $\D$, we define the data of a bilinear functor
\begin{equation}\label{evCD}
\big(\evCD \scs (\evCD)^2_1 \scs (\evCD)^2_2\big) \cn \psu(\C,\D) \times \C \to \D,
\end{equation}
which is called the \index{bilinear!evaluation functor}\index{evaluation!bilinear functor}\emph{evaluation}, as follows.
\begin{description}
\item[Objects] For a strictly unital symmetric monoidal functor $(F,F^2) \cn \C \to \D$ and an object $x \in \C$, the object assignment is defined as
\begin{equation}\label{evCD-Fx}
\evCD(F,x) = Fx \inspace \D.
\end{equation}
\item[Morphisms] Consider
\begin{itemize}
\item a monoidal natural transformation $\theta \cn F \to G$ between strictly unital symmetric monoidal functors $F,G \cn \C \to \D$ and
\item a morphism $f \cn x \to y$ in $\C$.
\end{itemize} 
The morphism 
\[\evCD(F,x) = Fx \fto{\evCD(\theta,f)} Gy = \evCD(G,y)\]
is defined as either one of the following two composites in $\D$, which are equal by the naturality of $\theta$.
\begin{equation}\label{evCDthetaf}
\begin{tikzpicture}[vcenter]
\def\v{-1.3}
\draw[0cell]
(0,0) node (a1) {Fx}
(a1)+(2,0) node (a2) {Gx}
(a1)+(0,\v) node (b1) {Fy}
(a2)+(0,\v) node (b2) {Gy}
;
\draw[1cell=.9]
(a1) edge node {\theta_x} (a2)
(a2) edge node {Gf} (b2)
(a1) edge node[swap] {Ff} (b1)
(b1) edge node {\theta_y} (b2)
;
\end{tikzpicture}
\end{equation}
\item[First Linearity Constraint] It is the identity natural transformation with components as follows.
\begin{equation}\label{evCD-constraint-i}
\begin{tikzpicture}[vcenter]
\def\u{-.8}
\draw[0cell=.9]
(0,0) node (a1) {\evCD(F,x) \oplus \evCD(G,x)}
(a1)+(5,0) node (a2) {\evCD(F \oplus G,x)}
(a1)+(0,\u) node (b1) {Fx \oplus Gx}
(a2)+(0,\u) node (b2) {(F \oplus G)x}
;
\draw[1cell=.9]
(a1) edge node {(\evCD)^2_1} (a2)
(b1) edge node {1} (b2)
(a1) edge[-,double equal sign distance] (b1)
(a2) edge[-,double equal sign distance] (b2)
;
\end{tikzpicture}
\end{equation}
\item[Second Linearity Constraint] It is given componentwise by the monoidal constraint of the first variable as follows.
\begin{equation}\label{evCD-constraint-ii}
\begin{tikzpicture}[vcenter]
\def\u{-.8}
\draw[0cell=.9]
(0,0) node (a1) {\evCD(F,x) \oplus \evCD(F,y)}
(a1)+(5,0) node (a2) {\evCD(F,x \oplus y)}
(a1)+(0,\u) node (b1) {Fx \oplus Fy}
(a2)+(0,\u) node (b2) {F(x \oplus y)}
;
\draw[1cell=.9]
(a1) edge node {(\evCD)^2_2} (a2)
(b1) edge node {F^2_{x,y}} (b2)
(a1) edge[-,double equal sign distance] (b1)
(a2) edge[-,double equal sign distance] (b2)
;
\end{tikzpicture}
\end{equation}
\end{description}
This finishes the definition of the evaluation.  If there is no danger of confusion, we omit the subscripts in $\ev_{\C,\D}$.
\end{definition}

\cref{ev-bilinear,ev-comp} below show that evaluation is bilinear and has the expected property with respect to composition.  \cref{ev-bilinear} is analogous to the composition bilinear functor in \cref{psu-mBCD}.  In \cref{clp-evaluation} we extend it to a multilinear functor as part of the closed multicategory structure on $\permcatsu$.

\begin{proposition}\label{ev-bilinear}
For small permutative categories $\C$ and $\D$, the triple
\[\big(\evCD \scs (\evCD)^2_1 \scs (\evCD)^2_2\big) \cn \psu(\C,\D) \times \C \to \D\]
in \cref{evCD} is a bilinear functor.
\end{proposition}

\begin{proof}
We prove statements \cref{evCD-bilinear-i,evCD-bilinear-ii,evCD-bilinear-iii} below.
\begin{romenumerate}
\item\label{evCD-bilinear-i} $\evCD$ is a functor.
\item\label{evCD-bilinear-ii} $(\evCD)^2_1$ in \cref{evCD-constraint-i} and $(\evCD)^2_2$ in \cref{evCD-constraint-ii} are natural transformations.
\item\label{evCD-bilinear-iii} $\evCD$ satisfies the bilinear functor axioms.
\end{romenumerate}

\medskip
\emph{Statement \cref{evCD-bilinear-i}}.  The assignment $\evCD$ preserves identity morphisms because all four morphisms in \cref{evCDthetaf} are identity morphisms if $\theta = 1_F$ and $f = 1_x$.  

To see that $\evCD$ preserves composition, consider a monoidal natural transformation $\psi \cn G \to H$ and a morphism $g \cn y \to z$ in $\C$.  In the following diagram in $\D$, the rectangle is commutative by the naturality of $\psi$.
\[\begin{tikzpicture}[vcenter]
\def\h{2} \def\v{-1.3}
\draw[0cell]
(0,0) node (a1) {Fx}
(a1)+(\h,0) node (a2) {Gx}
(a2)+(\h,0) node (a3) {Hx}
(a2)+(0,\v) node (b2) {Gy}
(a3)+(0,\v) node (b3) {Hy}
(b3)+(0,\v) node (c) {Hz}
;
\draw[1cell=.9]
(a1) edge node {\theta_x} (a2)
(a2) edge node {\psi_x} (a3)
(b2) edge node {\psi_y} (b3)
(a2) edge node[swap] {Gf} (b2)
(a3) edge node {Hf} (b3)
(b3) edge node {Hg} (c)
;
\end{tikzpicture}\]
The top-right composite is 
\[\big(H(gf)\big) (\psi\theta)_x = \evCD(\psi\theta,gf).\]  
The other composite is 
\[(Hg)\psi_y(Gf)\theta_x = \evCD(\psi,g) \evCD(\theta,f).\]
This shows that $\evCD$ is a functor.

\medskip
\emph{Statement \cref{evCD-bilinear-ii}}.  The morphism $(\evCD)^2_1$ in \cref{evCD-constraint-i} is natural 
\begin{itemize}
\item with respect to $F$ and $G$ by \cref{theta-oplus-thetap-x} and
\item with respect to $x$ by \cref{FoplusG-x}.
\end{itemize}
The morphism $(\evCD)^2_2$ in \cref{evCD-constraint-ii} is natural 
\begin{itemize}
\item with respect to $F$ by the monoidality of $\theta$ \cref{monnattr} and
\item with respect to $x$ and $y$ by the naturality of the monoidal constraint $F^2$.
\end{itemize}

\medskip
\emph{Statement \cref{evCD-bilinear-iii}}.  The unity axiom \cref{nlinearunity} follows from the definitions \cref{evCD-Fx,evCDthetaf}.

The first linearity constraint $(\evCD)^2_1$ satisfies the constraint unity, associativity, and symmetry axioms, \cref{constraintunity,eq:ml-f2-assoc,eq:ml-f2-symm}, because its components are identity morphisms.  

The second linearity constraint $(\evCD)^2_2$ satisfies the constrain unity axiom \cref{constraintunity} for the following reasons.
\begin{itemize}
\item The monoidal constraint of the monoidal unit \cref{psucd-unit} in $\psu(\C,\D)$ is the identity morphism $1_\pu$ in $\D$.
\item If either $x$ or $y$ is the monoidal unit $\pu$ in $\C$, then 
\[F^2_{\pu,y} = 1_{Fy} \andspace F^2_{x,\pu} = 1_{Fx}\]
by the unity of $(F,F^2,F^0 = 1_\pu)$ in \cref{monoidalfunctorunity}.
\end{itemize}
The constraint associativity and symmetry axioms, \cref{eq:ml-f2-assoc,eq:ml-f2-symm}, hold for $(\evCD)^2_2$ by the associativity \cref{monoidalfunctorassoc} and symmetry \cref{monoidalfunctorbraiding} of the monoidal constraint $F^2$.

The constraint 2-by-2 axiom \cref{eq:f2-2by2} is
\begin{itemize}
\item the left pentagon below if $(i,k)=(1,2)$ and
\item the right pentagon below if $(i,k)=(2,1)$.
\end{itemize}
\[\begin{tikzpicture}
\def\g{.5} \def\h{2} \def\u{1} \def\v{1.3}
\draw[0cell=.8]
(0,0) node (a1) {Fx \oplus Gx \oplus Fy \oplus Gy}
(a1)+(\h,\u) node (a2) {(F \oplus G)x \oplus (F \oplus G)y}
(a1)+(\h+\g,-\v/2) node (a3) {(F \oplus G)(x \oplus y)}
(a1)+(0,-\v) node (b1) {Fx \oplus Fy \oplus Gx \oplus Gy}
(b1)+(\h,-\u) node (b2) {F(x \oplus y) \oplus G(x \oplus y)}
;
\draw[1cell=.8]
(a1) edge[transform canvas={xshift=-.5em}] node[pos=.3] {1 \oplus 1} (a2)
(a2) edge node[pos=.6] {(F \oplus G)^2} (a3)
(a1) edge node[swap] {1 \oplus \xi \oplus 1} (b1)
(b1) edge[transform canvas={xshift=-.5em}] node[swap,pos=.3] {F^2 \oplus G^2} (b2)
(b2) edge node[swap,pos=.6] {1} (a3)
;
\begin{scope}[shift={(5.5,0)}]
\draw[0cell=.8]
(0,0) node (a1) {Fx \oplus Fy \oplus Gx \oplus Gy}
(a1)+(\h,\u) node (a2) {F(x \oplus y) \oplus G(x \oplus y)}
(a1)+(\h+\g,-\v/2) node (a3) {(F \oplus G)(x \oplus y)}
(a1)+(0,-\v) node (b1) {Fx \oplus Gx \oplus Fy \oplus Gy}
(b1)+(\h,-\u) node (b2) {(F \oplus G)x \oplus (F \oplus G)y}
;
\draw[1cell=.8]
(a1) edge[transform canvas={xshift=-.5em}] node[pos=.3] {F^2 \oplus G^2} (a2)
(a2) edge node {1} (a3)
(a1) edge node[swap] {1 \oplus \xi \oplus 1} (b1)
(b1) edge[transform canvas={xshift=-.5em}] node[swap,pos=.3] {1 \oplus 1} (b2)
(b2) edge node[swap,pos=.7] {(F \oplus G)^2} (a3)
;
\end{scope}
\end{tikzpicture}\]
The left pentagon is commutative by the definition of $(F \oplus G)^2$ in \cref{FoplusG-mon-constraint}.  The right pentagon is commutative by the symmetry axiom \cref{symmoncatsymhexagon} and the functoriality of $\oplus$ in the permutative category $\D$.
\end{proof}

The following result is the permutative categorical analog of the left diagram in \cref{eq:adj-comp} for symmetric monoidal closed categories.  We will use this result in \cref{ex:perm-closed-multicat} as part of the identification of the two self-enrichment of $\permcatsu$; see \cref{perm-clmulti-v} in that proof.  Recall the composition of multilinear functors in \cref{definition:permcat-comp} and the composition bilinear functor $\mcomp_{\B,\C,\D}$ in \cref{psu-mBCD}.

\begin{proposition}\label{ev-comp}
For small permutative categories $\B$, $\C$, and $\D$, the following two composite 3-linear functors are equal.
\begin{equation}\label{perm-ev-comp}

\end{equation}
\end{proposition}

\begin{proof}
We prove statements \cref{perm-ev-comp-i,perm-ev-comp-ii} below.
\begin{romenumerate}
\item\label{perm-ev-comp-i} The two composites in \cref{perm-ev-comp} are equal as functors.
\item\label{perm-ev-comp-ii} Their three respective linearity constraints are equal.
\end{romenumerate}

\medskip
\emph{Statement \cref{perm-ev-comp-i}}.  Consider a morphism $f \cn x \to y$ in $\B$ and monoidal natural transformations between strictly unital symmetric monoidal functors as follows.
\[%
\]
Each composite in \cref{perm-ev-comp} sends the triple $(H,F,x)$ to the object $HFx \in \D$.  For morphisms we consider the following diagram in $\D$.
\[%
\]
The following statements hold for the diagram above.
\begin{itemize}
\item The top-right composite is 
\[(IGf)(\psi * \theta)_x = \ev_{\B,\D}(\psi * \theta,f) 
= \ev_{\B,\D}\big(\mcomp_{\B,\C,\D}(\psi,\theta), f\big).\]
\item The bottom composite is $\evCD\big(\psi,\ev_{\B,\C}(\theta,f)\big)$ by the left-bottom composite in \cref{evCDthetaf} applied to $\evCD(\psi,-)$.
\item The left triangle is commutative by the definition \cref{evCDthetaf} of $\ev_{\B,\C}(\theta,f)$.
\item The right rectangle is commutative by the naturality of $\psi$.
\end{itemize}
This shows that \cref{perm-ev-comp} is commutative as a diagram of functors.

\medskip
\emph{Statement \cref{perm-ev-comp-ii}}.  Recall from \cref{ffjlinearity} the definition of the linearity constraints of a composite multilinear functor.  Now we consider the three linearity constraints of the two composites in \cref{perm-ev-comp}.
\begin{itemize}
\item Their first linearity constraints are equal because both first linearity constraints $\mcomp^2_1$ in \cref{HFIF} and $\ev^2_1$ in \cref{evCD-constraint-i} are the identity.
\item Their second linearity constraints are equal because both second linearity constraints $\mcomp^2_2$ in \cref{HFHG} and $\ev^2_2$ in \cref{evCD-constraint-ii} are given by the monoidal constraint of the monoidal functor in question.
\item Their third linearity constraints are equal because the monoidal constraint of the composite monoidal functor $HF$ is $H(F^2) \circ H^2$ by \cref{def:mfunctor-comp}.
\end{itemize}
This proves that the two composite 3-linear functors in \cref{perm-ev-comp} are equal.
\end{proof}

\section{Opposite Enriched Categories}
\label{sec:opp-enr-cat}

In this section we discuss the opposite of an $\M$-category when $(\M,\ga,\opu)$ is a multicategory; see \cref{opposite-mcat}.  By \cref{def:enr-multicategory} $\M$ is a $\Set$-multicategory, where $(\Set,\times,*)$ is the symmetric monoidal category of sets and functions with the monoidal product given by the Cartesian product.  We emphasize that a multicategory has a symmetric group action \cref{rightsigmaaction}, which is necessary to define opposite enriched categories.  

When the enriching category is symmetric monoidal, we observe that the opposite enriched category in \cref{definition:vcat-op} is the same as the one in this section; see \cref{v-opposite-mcat}.  Looking ahead we consider change of enrichment of opposite enriched categories in \cref{dF-opposite}.  Moreover, opposite enriched categories are important in \cref{ch:gspectra_Kem,part:homotopy-mackey}, where they are the domains of enriched presheaf categories \cref{mcat-copm}.

Recall the notion of an $\M$-category in \cref{def:menriched-cat}.

\begin{definition}\label{def:opposite-mcat}
  Suppose $(\M,\ga,\opu)$ is a multicategory, and $(\C,\mcomp,i)$ is an $\M$-category.  The \index{opposite!enriched category}\index{enriched category!opposite}\emph{opposite $\M$-category} 
\[(\Cop, \mcompop, i)\]
is the $\M$-category defined as follows.
\begin{description}
\item[Objects] It has the same class of objects as $\C$.
\item[Hom Objects] For objects $x,y\in \Cop$, its hom object is the object
\[\Cop(x,y) = \C(y,x) \inspace \M.\]
\item[Composition] For objects $x,y,z \in \Cop$, its composition is the binary multimorphism
\begin{equation}\label{mcompop-def}
\mcompop_{x,y,z} \cn \big( \Cop(y,z) \scs \Cop(x,y)\big) \to \Cop(x,z) \inspace \M
\end{equation}
given by the image of the composition binary multimorphism
\[\mcomp_{z,y,x} \cn \big( \C(y,x) \scs \C(z,y)\big) \to \C(z,x)\]
under the symmetric group action of the nonidentity permutation $\tau \in \Sigma_2$:
\begin{equation}\label{symm-group-action-tau}
\M\big(\C(y,x) \scs \C(z,y) \sscs \C(z,x)\big) \fto[\iso]{\tau} 
\M\big(\C(z,y) \scs \C(y,x) \sscs \C(z,x)\big).
\end{equation}
\item[Identities] The identity of an object $x \in \Cop$ is the nullary multimorphism
\begin{equation}\label{Cop-identity}
i_x \cn \ang{} \to \C(x,x) = \Cop(x,x) \inspace \M.
\end{equation}
This is the same as the identity of $x$ as an object in $\C$.
\end{description}
This finishes the definition of $(\Cop,\mcompop,i)$.  \cref{opposite-mcat} proves that $\Cop$ is actually an $\M$-category.
\end{definition}

\begin{explanation}\label{expl:cop-opposite-mcat}
We denote the composition $\mcompop_{x,y,z}$ in \cref{mcompop-def} diagrammatically as follows.
\begin{equation}\label{mcompop-diagram}
\begin{tikzpicture}[vcenter]
\def\u{-.8} \def\h{4} \def\t{15}
\draw[0cell=.9]
(0,0) node (a1) {\big(\Cop(y,z) \scs \Cop(x,y)\big)}
(a1)+(\h,0) node (a2) {\Cop(x,z)}
(a1)+(0,\u) node (b1) {\big(\C(z,y) \scs \C(y,x)\big)}
(a2)+(0,\u) node (b2) {\C(z,x)}
(b1)+(\h/2,-1) node (c) {\big(\C(y,x) \scs \C(z,y)\big)}
;
\draw[1cell=.9]
(a1) edge[-,double equal sign distance, shorten <=-.3ex, shorten >=-.3ex] (b1)
(a2) edge[-,double equal sign distance] (b2)
(a1) edge node {\mcompop_{x,y,z}} (a2)
(b1) edge[bend right=\t, transform canvas={xshift=-1.5ex}, shorten >=-1.5ex] node[swap] {\tau} (c)
(c) edge[bend right=\t, transform canvas={xshift=1ex}, shorten <=-1ex] node[swap] {\mcomp_{z,y,x}} (b2)
;
\end{tikzpicture}
\end{equation}
We understand this composite as given by the symmetric group action \cref{symm-group-action-tau} in $\M$.  This diagram is an analog of the composite \cref{opposite-comp} that defines the composition of an opposite $\V$-category with $\V$ a braided monoidal category.
\end{explanation}

We now check that $\Cop$ satisfies the axioms of an $\M$-category.  The following observation is stated in \cite[Remark 2.8]{bohmann_osorno-mackey}.

\begin{proposition}\label{opposite-mcat}
In the context of \cref{def:opposite-mcat}, $(\Cop,\mcompop,i)$ is an $\M$-category.
\end{proposition}

\begin{proof}
We need to prove the associativity axiom \cref{menriched-cat-assoc} and the unity axiom \cref{menriched-cat-unity} for $\Cop$.  

For objects $w,x,y,z \in \Cop$, the associativity diagram \cref{menriched-cat-assoc} for $\Cop$ is the boundary of the following diagram in $\M$, with $\C_{x,y}$ denoting $\C(x,y)$ and $\tau \in \Sigma_2$ denoting the nonidentity permutation.
\[\begin{tikzpicture}
\def\g{3.7} \def\h{3.3} \def\v{-1.4} \def\w{.6} \def\x{1.7} \def\u{1.3}
\draw[0cell=.85]
(0,0) node (a1) {(\C_{z,y} \scs \C_{y,x} \scs \C_{x,w})}
(a1)+(\g,0) node (a2) {(\C_{y,x} \scs \C_{z,y} \scs \C_{x,w})}
(a2)+(\h,0) node (a3) {(\C_{z,x} \scs \C_{x,w})}
(a1)+(0,\v) node (b1) {(\C_{z,y} \scs \C_{x,w} \scs \C_{y,x})}
(b1)+(\g,0) node (b2) {(\C_{x,w} \scs \C_{y,x} \scs \C_{z,y})}
(b2)+(\h,0) node (b3) {(\C_{x,w} \scs \C_{z,x})}
(b1)+(0,\v) node (c1) {(\C_{z,y} \scs \C_{y,w})}
(c1)+(\g,0) node (c2) {(\C_{y,w} \scs \C_{z,y})}
(c2)+(\h,0) node (c3) {\C_{z,w}}
;
\draw[1cell=.85]
(a1) edge node {(\tau,\opu)} (a2)
(a2) edge node {(\mcomp,\opu)} (a3)
(b1) edge node {\tau\ang{1,2}} (b2)
(b2) edge node {(\opu,\mcomp)} (b3)
(c1) edge node {\tau} (c2)
(c2) edge node {\mcomp} (c3)
(a1) edge node {(\opu,\tau)} (b1)
(b1) edge node {(\opu,\mcomp)} (c1)
(a2) edge node {\tau\ang{2,1}} (b2)
(b2) edge node {(\mcomp,\opu)} (c2)
(a3) edge node[swap] {\tau} (b3)
(b3) edge node[swap] {\mcomp} (c3)
;
\draw[1cell=.85]
(a1) [rounded corners=3pt] |- ($(a2)+(-1,\w)$)
-- node[pos=0] {(\mcompop,\opu)} ($(a2)+(1,\w)$) -| (a3)
;
\draw[1cell=.85]
(c1) [rounded corners=3pt] |- ($(c2)+(-1,-\w)$)
-- node[pos=0] {\mcompop} ($(c2)+(1,-\w)$) -| (c3)
;
\draw[1cell=.85]
(a1) [rounded corners=3pt] -| ($(b1)+(-\x,.7)$)
-- node[pos=0] {(\opu,\mcompop)} ($(b1)+(-\x,-.7)$) |- (c1)
;
\draw[1cell=.85]
(a3) [rounded corners=3pt] -| ($(b3)+(\u,.7)$)
-- node[pos=0,swap] {\mcompop} ($(b3)+(\u,-.7)$) |- (c3)
;
\end{tikzpicture}\]
The following statements hold for the diagram above.
\begin{itemize}
\item The right and bottom strips are the definition of $\mcompop$ \cref{mcompop-diagram}.
\item The left and top strips are commutative by \cref{mcompop-diagram} and the right unity axiom \cref{enr-multicategory-right-unity} for $\M$.
\item Each $\tau\ang{r,s} \in \Sigma_3$ is the block permutation \cref{blockpermutation} that permutes a block of length $r$ with a block of length $s$.  There are equalities of permutations
\[\tau\ang{1,2} \circ (1 \times \tau) = (1,3) = \tau\ang{2,1} \circ (\tau \times 1) \inspace \Sigma_3,\]
where $(1,3)$ denotes the transposition of 1 and 3.  The upper left rectangle composed with either $\ga\scmap{\mcomp;\mcomp,\opu}$ or $\ga\scmap{\mcomp;\opu,\mcomp}$ is commutative by the symmetric group axiom \cref{enr-multicategory-symmetry} for $\M$.  
\item The lower left rectangle composed with the lower right horizontal $\mcomp$ is commutative by the top equivariance axiom \cref{enr-operadic-eq-1} for $\M$.
\item The lower right rectangle is commutative by the associativity axiom \cref{menriched-cat-assoc} for $\C$.
\item The upper right rectangle composed with the lower right vertical $\mcomp$ is commutative by the top equivariance axiom \cref{enr-operadic-eq-1} for $\M$.
\end{itemize}
This proves the associativity axiom for $\Cop$.

For objects $x,y \in \Cop$, the unity diagram \cref{menriched-cat-unity} for $\Cop$ is the boundary of the following diagram in $\M$.
\[\begin{tikzpicture}
\def\g{.8} \def\h{4} \def\f{1.5} \def\v{-1.3} \def\t{15}
\draw[0cell=.8]
(0,0) node (a1) {(\C_{y,x} \scs \ang{})}
(a1)+(\h,0) node (a2) {\phantom{(\ang{} \scs \C_{y,x})}}
(a2)+(-\g,0) node (a2') {(\ang{} \scs \C_{y,x})}
(a2)+(\g,0) node (a2'') {(\C_{y,x} \scs \ang{})}
(a2)+(\h,0) node (a3) {(\ang{} \scs \C_{y,x})}
(a1)+(0,\v) node (b1) {(\C_{y,x} \scs \C_{x,x})}
(a3)+(0,\v) node (b3) {(\C_{y,y} \scs \C_{y,x})}
(b1)+(\f,\v) node (c1) {(\C_{x,x} \scs \C_{y,x})}
(a2)+(0,2*\v) node (c2) {\C_{y,x}}
(b3)+(-\f,\v) node (c3) {(\C_{y,x} \scs \C_{y,y})}
;
\draw[1cell=.8]
(a1) edge node {\tau\ang{1,0} = 1} (a2')
(a2') edge[-,double equal sign distance] (a2'')
(a3) edge node[swap] {\tau\ang{0,1} = 1} (a2'')
(a1) edge node[swap] {(\opu,i_x)} (b1)
(b1) edge[bend right=\t, transform canvas={xshift=-.5ex}] node[swap] {\tau} (c1)
(c1) edge node {\mcomp} (c2)
(a3) edge node {(i_y,\opu)} (b3)
(b3) edge[bend left=\t, transform canvas={xshift=.5ex}] node {\tau} (c3)
(c3) edge node[swap] {\mcomp} (c2)
(a2') edge node[swap,pos=.4] {(i_x,\opu)} (c1)
(a2) edge node {\opu} (c2)
(a2'') edge node[pos=.4] {(\opu,i_y)} (c3)
;
\end{tikzpicture}\]
The following statements hold for the diagram above.
\begin{itemize}
\item The middle two sub-regions are commutative by the unity axiom \cref{menriched-cat-unity} for $\C$.
\item Since the block permutations $\tau\ang{1,0}$ and $\tau\ang{0,1} \in \Sigma_1$ are the identity, they act on $\M$ as the identity by the first part of the symmetric group axiom \cref{enr-multicategory-symmetry} for $\M$.
\item The left sub-region composed with the bottom left $\mcomp$ is commutative by the top equivariance axiom \cref{enr-operadic-eq-1} for $\M$.
\item The right sub-region composed with the bottom right $\mcomp$ is commutative by the top equivariance axiom \cref{enr-operadic-eq-1} for $\M$.
\end{itemize}
This proves the unity axiom for $\Cop$.
\end{proof}

For a monoidal category $\V$, \cref{EndV-enriched} shows that a $\V$-category is the same thing as an $(\End\,\V)$-category for the non-symmetric endomorphism multicategory $\EndV$.  Next we observe that, if $\V$ is symmetric monoidal, then we can also identify opposite categories enriched in $\V$ and in $\EndV$.

\begin{proposition}\label{v-opposite-mcat}\index{opposite!enriched category}\index{enriched category!opposite}
For each symmetric monoidal category $(\V,\otimes,\tu)$ and $\V$-category $\C$,
\begin{itemize}
\item the opposite $\V$-category $\Cop$ in \cref{definition:vcat-op} and
\item the opposite $(\End\, \V)$-category $\Cop$ in \cref{def:opposite-mcat}
\end{itemize}
are the same.
\end{proposition}

\begin{proof}
A comparison of \cref{definition:vcat-op,def:opposite-mcat} shows that $\Cop$ in these two definitions have
\begin{itemize}
\item the same objects, namely, the objects of $\C$;
\item the same hom objects, namely,
\[\Cop(x,y) = \C(y,x)\]
for objects $x,y \in \Cop$; and
\item the same identity for each object $x \in \Cop$, namely,
\[i_x \in (\End\,\V)\scmap{\ang{}; \C(x,x)} = \V\big(\tu,\C(x,x)\big).\]
\end{itemize} 
Their compositions, \cref{mcompop-diagram,opposite-comp}, are also the same because the symmetric group action of the multicategory $\End\,\V$ is induced by the braiding of the symmetric monoidal category $\V$.
\end{proof}

\begin{example}\label{ex:v-opposite-mcat}
\cref{v-opposite-mcat} applies to all the symmetric monoidal categories in \cref{ex:EndV-enriched}.  For example, consider the symmetric monoidal closed category $\MoneMod$ in \cref{proposition:EM2-5-1} \cref{monebicomplete}, with associated $\Cat$-multicategory in \cref{expl:monemodcatmulticat}.  For a $\MoneMod$-category $\C$, the opposite $\MoneMod$-category $\Cop$ in the sense of \cref{definition:vcat-op,def:opposite-mcat} are the same.
\end{example}

\chapter{Change of Multicategorical Enrichment}
\label{ch:change_enr}
This chapter defines and develops the basic properties of the change-of-enrichment 2-functor induced by a non-symmetric multifunctor
\[
  F \cn \M \to \N
\]
between non-symmetric multicategories $\M$ and $\N$.
If $\M$ and $\N$ are multicategories---equipped with the necessary symmetric group actions---then change of enrichment is shown to preserve opposites in \cref{dF-opposite}.
Composing change of enrichment functors is treated in \cref{func-change-enr}, and \cref{change-enr-twofunctor} extends this to show that there is a 2-functor
\[\Enr \cn \Multicatns \to \iicat\]
given by the assignments
\begin{itemize}
\item $\Enr\M = \MCat$ (\cref{mcat-iicat}) on objects,
\item $\Enr F = \dF$ (\cref{mult-change-enrichment}) on 1-cells, and
\item $\Enr\theta = \dtheta$ (\cref{dtheta-twonat}) on 2-cells.
\end{itemize}
Our main motivation for these results is for application to $K$-theory multifunctors, discussed in \cref{ex:dF-opposite,ex:mon-change-enr,ex:func-change-enr}.

\subsection*{Connection with Other Chapters}

The material in this chapter is used in each of \cref{ch:std_enrich,ch:gspectra_Kem,ch:mackey,ch:mackey_eq}.
Of these, \cref{ch:std_enrich,ch:gspectra_Kem} extend the theory here to that of self-enrichments and enriched diagrams, respectively.
The factorization of Elmendorf-Mandell $K$-theory \cref{Kem-diagram} is discussed further
\begin{itemize}
\item in \cref{gspectra-thm-xi} in the context of standard enrichment and
\item in \cref{gspectra-thm-xiv} in the context of presheaf change of enrichment.
\end{itemize}
\cref{ch:mackey,ch:mackey_eq} give homotopy-theoretic applications, with  \cref{ch:mackey_eq} focusing on diagrams and Mackey functors enriched in $\pMulticat$ and $\permcatsu$.

\subsection*{Background}

The content of this chapter depends on the multicategorical enrichment developed in \cref{ch:menriched}.

\subsection*{Chapter Summary}

\cref{sec:change-enr-mfunctor} defines the change-of-enrichment 2-functor along a non-symmetric multifunctor.
\cref{sec:pres-opposite} specializes to the symmetric case and shows that change of enrichment along a multifunctor preserves opposite enriched categories.
\cref{sec:change-enr-monfunctor} shows that the two notions of change of enrichment along a monoidal functor---monoidal or multicategorical---agree.
\cref{sec:functoriality-change-enr} describes compositionality for change of enrichment, and \cref{sec:twofunc-change-enr} extends this to show that change of enrichment is 2-functorial.
Here is a summary table.
\reftable{.9}{
  change of enrichment along a non-symmetric multifunctor
  & \ref{def:mult-change-enr} and \ref{mult-change-enrichment}
  \\ \hline
  examples of change of enrichment
  & \ref{ex:mult-change-enr} and \ref{ex:dFst}
  \\ \hline
  preservation of opposite enriched categories
  & \ref{dF-opposite} and \ref{ex:dF-opposite}
  \\ \hline
  change of enrichment along a monoidal functor
  & \ref{mon-change-enrichment}
  \\ \hline
  composition of change-of-enrichment
  & \ref{func-change-enr}
  \\ \hline
  2-functoriality of change of enrichment
  & \ref{change-enr-twofunctor}
  \\
}

We remind the reader of \cref{conv:universe} about universes and \cref{expl:leftbracketing} about left normalized bracketing for iterated products.

\section{Change of Enrichment along a Multifunctor}
\label{sec:change-enr-mfunctor}

Each monoidal functor between monoidal categories 
\[U \cn \V \to \W\]
induces a change-of-enrichment 2-functor (\cref{proposition:U-VCat-WCat})
\[\dU \cn \VCat \to \WCat.\]
In this section we generalize this construction from a monoidal functor to a non-symmetric multifunctor.  This section is organized as follows.
\begin{itemize}
\item The change of enrichment along a non-symmetric multifunctor $F$ is constructed in \cref{def:mult-change-enr} and is shown to be a 2-functor in \cref{mult-change-enrichment}.
\item \cref{ex:mult-change-enr} illustrates \cref{mult-change-enrichment} with $K$-theory multifunctors, some of which are non-symmetric.
\item As a further illustration of change of enrichment, in \cref{ex:dFst} we explicitly describe the change-of-enrichment 2-functor induced by the non-symmetric multifunctor (\cref{ptmulticat-xvii})
\[\Fst \cn \pMulticat \to \permcatsu\]
given by the pointed free permutative category construction.
\end{itemize}

\subsection*{Defining Change of Enrichment}

Recall
\begin{itemize}
\item non-symmetric multicategories (\cref{def:enr-multicategory}),
\item non-symmetric multifunctors (\cref{def:enr-multicategory-functor}), and
\item for a non-symmetric multicategory $\M$, the 2-category $\MCat$ of small $\M$-categories, $\M$-functors, and $\M$-natural transformations (\cref{mcat-iicat}).
\end{itemize}
We emphasize that a non-symmetric multifunctor preserves colored units \cref{enr-multifunctor-unit} and composition \cref{v-multifunctor-composition}, but not necessarily the symmetric group action even if its domain and codomain are multicategories.

First we define the object, 1-cell, and 2-cell assignments of change of enrichment.  Recall the notion of a 2-functor (\cref{def:twofunctor}).

\begin{definition}\label{def:mult-change-enr}
Suppose given a non-symmetric multifunctor between non-symmetric multicategories
\[F \cn (\M,\ga,\opu) \to (\N,\ga,\opu).\]
We define the data of a 2-functor\label{not:dF}
\[\dF \cn \MCat \to \NCat,\]
which is called the \index{change of enrichment}\index{enriched category!change of enrichment}\index{monoidal functor!change of enrichment}\index{symmetric monoidal functor!change of enrichment}\index{2-functor!change of enrichment}\emph{change of enrichment} or the \emph{change-of-enrichment 2-functor} along $F$, as follows.
\begin{description}
\item[Object Assignment]
The image of an $\M$-category $(\C,\mcomp,i)$ (\cref{def:menriched-cat}) under $\dF$ is the $\N$-category
\begin{equation}\label{mcat-change-enrichment}
\big(\C_F,\mcomp_F,i_F\big)
\end{equation}
consisting of the following data.
\begin{itemize}
\item The objects of $\C_F$ are those of $\C$.
\item For each pair of objects $x,y \in \C_F$, the hom object is
\begin{equation}\label{CFxy}
(\C_F)(x,y) = F\C(x,y) \inspace \N.
\end{equation}
\item For objects $x,y,z \in \C_F$, the composition binary multimorphism
\begin{equation}\label{mFxyz}
\big(F\C(y,z) \scs F\C(x,y)\big) \fto{(\mcomp_F)_{x,y,z} = F\left(\mcomp_{x,y,z}\right)} F\C(x,z) \inspace \N
\end{equation}
is the image under $F$ of the composition $\mcomp_{x,y,z}$ in \cref{menriched-cat-comp}. 
\item For each object $x \in \C$, the identity nullary multimorphism
\begin{equation}\label{iFx}
\ang{} \fto{(i_F)_x = F(i_x)} F\C(x,x) \inspace \N
\end{equation}
is the image under $F$ of the identity $i_x$ in \cref{menriched-cat-identity}.
\end{itemize}
This finishes the definition of the $\N$-category $(\C_F,\mcomp_F,i_F)$ in \cref{mcat-change-enrichment}.  Its associativity diagram \cref{menriched-cat-assoc} and unity diagram \cref{menriched-cat-unity} are obtained from those for $\C$ by applying $F$, which preserves colored units and composition.
\item[1-Cell Assignment]
The image of an $\M$-functor between $\M$-categories (\cref{def:mfunctor})
\[H \cn (\C,\mcomp,i) \to (\D,\mcomp,i)\]
under $\dF$ is the $\N$-functor
\begin{equation}\label{mfunctor-change-enrichment}
H_F \cn \big(\C_F,\mcomp_F,i_F\big) \to \big(\D_F,\mcomp_F,i_F\big)
\end{equation}
defined as follows.
\begin{itemize}
\item The object assignment of $H_F$ is the object assignment of $H$.
\item For each pair of objects $x,y \in \C_F$, the component unary multimorphism 
\begin{equation}\label{HFxy}
F\C(x,y) \fto{(H_F)_{x,y} = F(H_{x,y})} F\D(Hx,Hy) \inspace \N
\end{equation}
is the image under $F$ of the component $H_{x,y}$ in \cref{Fxy-component}.
\end{itemize}
This finishes the definition of the $\N$-functor $H_F$ in \cref{mfunctor-change-enrichment}.  Its compatibility axioms \cref{mfunctor-diagrams} are obtained from those for $H$ by applying the non-symmetric multifunctor $F$.
\item[2-Cell Assignment]
Suppose $\theta \cn H \to G$ is an $\M$-natural transformation (\cref{def:mnaturaltr}) as in the left diagram below.
\[\begin{tikzpicture}
\def\h{2} \def\t{25} \def\s{30}
\draw[0cell]
(0,0) node (a) {\C}
(a)+(\h,0) node (b) {\D}
;
\draw[1cell=.8]
(a) edge[bend left=\t] node {H} (b)
(a) edge[bend right=\t] node[swap] {G} (b)
;
\draw[2cell]
node[between=a and b at .45, rotate=-90, 2label={above,\theta}] {\Rightarrow}
;
\begin{scope}[shift={(4,0)}]
\draw[0cell]
(0,0) node (a) {\C_F}
(a)+(\h,0) node (b) {\D_F}
;
\draw[1cell=.8]
(a) edge[bend left=\t] node {H_F} (b)
(a) edge[bend right=\t] node[swap] {G_F} (b)
;
\draw[2cell]
node[between=a and b at .42, rotate=-90, 2label={above,\theta_F}] {\Rightarrow}
;
\end{scope}
\end{tikzpicture}\]
We define $\theta_F$ as the $\N$-natural transformation, as in the right diagram above, with, for each object $x \in \C_F$, $x$-component nullary multimorphism
\begin{equation}\label{mnat-change-enrichment}
\ang{} \fto{(\theta_F)_x = F(\theta_x)} F\D(Hx,Gx) \inspace \N.
\end{equation}
This is the image under $F$ of the $x$-component 
\[\theta_x \cn \ang{} \to \D(Hx,Gx)\]
of $\theta$ in \cref{mnatural-component}, which is a nullary multimorphism in $\M$.  The naturality diagram \cref{mnaturality-diag} for $\theta_F$ is obtained from the naturality diagram for $\theta$ by applying $F$.
\end{description}
This finishes the definition of $\dF$.
\end{definition}

The following result is stated in \cite[Proposition 2.11]{bohmann_osorno-mackey}.

\begin{proposition}\label{mult-change-enrichment}
For each non-symmetric multifunctor $F \cn \M \to \N$ between non-symmetric multicategories, the change of enrichment along $F$ in \cref{def:mult-change-enr} is a 2-functor
\[\dF \cn \MCat \to \NCat.\]
\end{proposition}

\begin{proof}
We need to check that the assignment $\dF$ preserves
\begin{itemize}
\item identity 1-cells \cref{id-mfunctor},
\item horizontal composition of 1-cells \cref{mfunctor-composition},
\item identity 2-cells \cref{id-mnat},
\item vertical composition of 2-cells \cref{mnat-vcomp-component}, and
\item horizontal composition of 2-cells \cref{mnat-hcomp-component}.
\end{itemize}
These preservation properties of $\dF$ follow from
\begin{enumerate}
\item the componentwise definitions \cref{CFxy,mFxyz,iFx,HFxy,mnat-change-enrichment}, and
\item the fact that $F$ preserves colored units and composition.
\end{enumerate} 
This finishes the proof.
\end{proof}

\subsection*{Examples of Change of Enrichment}

\begin{example}[$K$-Theory Multifunctors]\label{ex:mult-change-enr}\index{K-theory@$K$-theory!change of enrichment along}
The change-of-enrichment 2-functor exists for each of the multifunctors in the following diagram from \cref{endufactor,eq:Ksummary,eq:thm-F-multi} along with \cref{eq:ptmulticat-xvii,Fm-multi-def}.  Among these arrows, $\Fr$, $\Fst$, $\Fm$, and $\cP$ are non-symmetric.
\begin{equation}\label{mult-change-enr-diagram}
\begin{tikzpicture}[vcenter]
\def\h{3.2} \def\v{-1.4} \def\w{.6}
\draw[0cell=.9]
(0,0) node (p) {\permcatsu}
(p)+(-\h,0) node (m) {\Multicat}
(p)+(0,\v) node (mone) {\MoneMod}
(m)+(0,\v) node (pm) {\pMulticat}
(p)+(\h,0) node (ga) {\Gacat}
(ga)+(0,\v) node (gs) {\Gstarcat}
(ga)+(1.5,0) node (sp) {\Sp}
;
\draw[1cell=.8,shift left]
(p) edge node[pos=.7] {\End} (m)
(p) edge[shorten <=.5ex] node[pos=.7] {\Endst} (pm)
(p) edge node[pos=.7] {\Endm} (mone)
;
\draw[1cell=.8,shift left]
(m) edge node[pos=.3] {\Fr} (p)
(pm) edge[shorten >=1.5ex] node[pos=.2] {\Fst} (p)
(mone) edge node[pos=.3] {\Fm} (p)
;
\draw[1cell=.9]
(p) edge[shorten <=1ex] node {\Jem} (gs)
(ga) edge node[swap] {\cP} (p)
;
\draw[1cell=.9]
(p) [rounded corners=3pt] |- ($(ga)+(-1,\w)$)
-- node[pos=.2] {\Kem} ($(ga)+(1,\w)$) -| (sp)
;
\end{tikzpicture}
\end{equation}
In other words, each of these multifunctors (non-symmetric for $\Fr$, $\Fst$, $\Fm$, and $\cP$) induces a change-of-enrichment 2-functor as in \cref{mult-change-enrichment}.  We emphasize that \emph{none} of the arrows in \cref{mult-change-enr-diagram} is a monoidal functor because the multicategory structure on $\permcatsu$ is not induced by a monoidal structure.  Thus \cref{proposition:U-VCat-WCat} does not apply to the arrows in \cref{mult-change-enr-diagram}.  In \cref{ex:mon-change-enr} below we extend this example to include other $K$-theoretic symmetric monoidal functors. 
\end{example}

\begin{explanation}[Change of Enrichment along $\Fst$]\label{ex:dFst}\index{category!free permutative!pointed}\index{permutative category!free!pointed}\index{free!permutative category!pointed}\index{pointed!free permutative category}
To illustrate \cref{def:mult-change-enr}, we describe the change-of-enrichment 2-functor
\begin{equation}\label{dFst-explicit}
\dFst \cn \pMulticatcat \to \permcatsucat
\end{equation}
induced by the non-symmetric multifunctor (\cref{ptmulticat-xvii})
\[\Fst \cn \pMulticat \to \permcatsu.\]
To simplify the presentation, we use the shortened notation\label{not:pM}
\[\pM = \pMulticat \andspace \psu = \permcatsu.\]
Since $\pM$ is a symmetric monoidal closed category (\cref{thm:pmulticat-smclosed}), by \cref{EndV-enriched} the 2-category $\pMcat$ is the same regardless of whether we consider it in the sense of \cref{ex:vcatastwocategory} or \cref{mcat-iicat}.

\medskip
\emph{$\dFst$ on Objects}.  Suppose $(\C,\mcomp,i)$ is an $\pM$-category (\cref{def:enriched-category}).  Besides its class of objects, $\C$ consists of the following data.
\begin{itemize}
\item For each pair of objects $x,y\in \C$, the hom object $\C(x,y)$ is a small pointed multicategory (\cref{def:ptd-multicat}).
\item For objects $x,y,z \in \C$, the composition is a pointed multifunctor
\begin{equation}\label{mcomp-xyz-dFst}
\mcomp_{x,y,z} \cn \C(y,z) \sma \C(x,y) \to \C(x,z).
\end{equation}
\item The identity of each object $x \in \C$ is a pointed multifunctor
\begin{equation}\label{ix-dFst}
i_x \cn \Mtup = \Mtu \bincoprod \Mterm \to \C(x,x)
\end{equation}
from the smash unit $\Mtup$ in \cref{eq:smashunit}.  Since $i_x$ preserves the basepoint, it sends the unique object $* \in \Mterm$ to the basepoint of $\C(x,x)$.  Thus $i_x$ is determined by the object 
\[i_x(1) \in \C(x,x)\]
with $1$ denoting the unique object in the initial operad $\Mtu$ in \cref{ex:vmulticatinitialterminal} \cref{ex:initialoperad}.  To simplify the notation, we also denote the object $i_x(1)$ by $i_x$.
\end{itemize}
The associativity diagram \cref{enriched-cat-associativity} and the unity diagram \cref{enriched-cat-unity} are required to commute.

Applying \cref{CFxy,mFxyz,iFx} to $\Fst$, the $\psu$-category (\cref{expl:perm-enr-cat})
\[\big(\C_{\Fst} \scs \mcomp_{\Fst} \scs i_{\Fst}\big)\]
consists of the following data.
\begin{itemize}
\item $\Ob\big(\C_{\Fst}\big) = \Ob\C$.
\item For each pair of objects $x,y \in \C_{\Fst}$, the hom object
\[\C_{\Fst}(x,y) = \Fst \C(x,y)\]
is the pointed free permutative category (\cref{def:Fst-object,def:Fst-permutative}) of the small pointed multicategory $\C(x,y)$.
\item For objects $x,y,z \in \C_{\Fst}$, applying the $n=2$ case of \cref{Sst-Fbst} to $\mcomp_{x,y,z}$, the composition bilinear functor 
\[\begin{tikzpicture}
\def\w{.6}
\draw[0cell=.9]
(0,0) node (a) {\Fst\C(y,z) \times \Fst\C(x,y)}
(a)+(4,0) node (b) {\Fst\big(\C(y,z) \sma \C(x,y)\big)}
(b)+(4,0) node (c) {\Fst\C(x,z)}
;
\draw[1cell=.85]
(a) edge node{\Fst^2} (b)
(b) edge node {\Fst(\mcomp_{x,y,z})} (c)
;
\draw[1cell=.9]
(a) [rounded corners=3pt] |- ($(b)+(-1,\w)$)
-- node {(\mcomp_{\Fst})_{x,y,z}} ($(b)+(1,\w)$) -| (c)
;
\end{tikzpicture}\]
is the composite of
\begin{itemize}
\item the strong bilinear functor $\Fst^2$ in \cref{ptmulticat-xv} and
\item the strict symmetric monoidal functor $\Fst(\mcomp_{x,y,z})$ in \cref{def:Fst-onecells}.  
\end{itemize} 
\item For each object $x \in \C$, we also regard the pointed multifunctor $i_x$ in \cref{ix-dFst} as an object in the pointed multicategory $\C(x,x)$, given by the image of the unique object $1 \in \Mtu$.  The identity of an object $x \in \C_{\Fst}$ is the $\obsim$-equivalence class 
\[(i_{\Fst})_x = [(i_x)] \in \Fst\C(x,x)\]
of the length-one sequence $(i_x) \in \Fr\C(x,x)$ (\cref{def:Fst-object}).  More explicitly, applying the $n=0$ case of \cref{Sst-Fbst} to $i_x$, the identity $(i_{\Fst})_x$ is given by the 0-linear functor
\begin{equation}\label{FstzeroFstix}
\ang{} \fto{\Fst^0} \Fst(\Mtup) \fto{\Fst(i_x)} \Fst\C(x,x).
\end{equation}
By \cref{def:Sst} $\Fst^0$ is given by the $\obsim$-equivalence class $[(1)] \in \Fst(\Mtup)$.  Then $\Fst(i_x)$ sends $[(1)]$ to $[(i_x)]$ by definition \cref{eq:FstHx} applied to the pointed multifunctor $i_x$.
\end{itemize}

\medskip
\emph{$\dFst$ on 1-Cells}. 
Consider an $\pM$-functor (\cref{def:enriched-functor}) between $\pM$-categories
\[H \cn (\C,\mcomp,i) \to (\D,\mcomp,i).\]
Besides its object assignment, $H$ has, for each pair of objects $x,y \in \C$, a component pointed multifunctor
\begin{equation}\label{Hxy-ptmfunctor}
H_{x,y} \cn \C(x,y) \to \D(Hx,Hy).
\end{equation}
These pointed multifunctors are compatible with composition and identities in the sense of \cref{eq:enriched-composition}.  In particular, compatibility with identities means that, for each object $x \in \C$, the identity $i_x \in \C(x,x)$ is sent to the identity
\[H_{x,x}(i_x) = i_{Hx} \in \D(Hx,Hx).\]

Applying the change of enrichment along $\Fst$, the $\psu$-functor (\cref{expl:perm-enr-functors})
\[H_{\Fst} \cn \big(\C_{\Fst} \scs \mcomp_{\Fst} \scs i_{\Fst}\big) \to 
\big(\D_{\Fst} \scs \mcomp_{\Fst} \scs i_{\Fst}\big)\]
has the same object assignment as $H$.  For objects $x,y \in \C_{\Fst}$, $H_{\Fst}$ has the component strict symmetric monoidal functor \cref{HFxy}
\[\Fst\C(x,y) \fto{\Fst H_{x,y}} \Fst\D(Hx,Hy).\]
This is obtained from the pointed multifunctor $H_{x,y}$ in \cref{Hxy-ptmfunctor} by applying the 1-cell assignment of $\Fst$ in \cref{def:Fst-onecells}.

\medskip
\emph{$\dFst$ on 2-Cells}. 
Consider an $\pM$-natural transformation (\cref{def:enriched-natural-transformation})
\[\begin{tikzpicture}
\def\t{25}
\draw[0cell]
(0,0) node (a) {\phantom{C}}
(a)+(1.8,0) node (b) {\phantom{C}}
(a)+(-.5,0) node (a') {(\C,\mcomp,i)}
(b)+(.5,0) node (b') {(\D,\mcomp,i)}
;
\draw[1cell=.8]
(a) edge[bend left=\t] node {H} (b)
(a) edge[bend right=\t] node[swap] {G} (b)
;
\draw[2cell]
node[between=a and b at .45, rotate=-90, 2label={above,\theta}] {\Rightarrow}
;
\end{tikzpicture}\]
between $\pM$-functors between $\pM$-categories.  For each object $x \in \C$, the $x$-component of $\theta$ is a pointed multifunctor
\[\theta_x \cn \Mtup = \Mtu \bincoprod \Mterm \to \D(Hx,Gx).\]
Similar to \cref{ix-dFst}, the pointed multifunctor $\theta_x$ is determined by the object
\[\theta_x(1) \in \D(Hx,Gx)\]
with $1 \in \Mtu$ denoting the unique object.  As before we abbreviate $\theta_x(1)$ to $\theta_x$.  For objects $x,y \in \C$, the naturality diagram \cref{enr-naturality} for $\theta$ is the following commutative diagram of pointed multifunctors.
\begin{equation}\label{pMnt-naturality}
\begin{tikzpicture}[baseline={(a1.base)}]
\def\g{.3} \def\h{4.5} \def\v{1.2}
\draw[0cell=.9]
(0,0) node (a1) {\C(x,y)}
(a1)+(\g,\v) node (a2) {\Mtup \sma \C(x,y)}
(a2)+(\h,0) node (a3) {\D(Hy,Gy) \sma \D(Hx,Hy)}
(a3)+(\g,-\v) node (a4) {\D(Hx,Gy)}
(a1)+(\g,-\v) node (b1) {\C(x,y) \sma \Mtup}
(b1)+(\h,0) node (b2) {\D(Gx,Gy) \sma \D(Hx,Gx)}
;
\draw[1cell=.85]
(a1) edge node[pos=.2] {\lambda^\inv} node[swap,pos=.7] {\iso} (a2)
(a2) edge node {\theta_y \sma H_{x,y}} (a3)
(a3) edge node[pos=.7] {\mcomp} (a4)
(a1) edge node[pos=.7] {\iso} node[swap,pos=.1] {\rho^\inv} (b1)
(b1) edge node {G_{x,y} \sma \theta_x} (b2)
(b2) edge node[swap,pos=.7] {\mcomp} (a4)
;
\end{tikzpicture}
\end{equation}
The commutative diagram \cref{pMnt-naturality} is equivalent to the following equality in $\D(Hx,Gy)$ for objects and multimorphisms $f \in \C(x,y)$.
\[\mcomp\big(\theta_y \sma (H_{x,y}f)\big) = \mcomp\big((G_{x,y}f) \sma \theta_x\big)\]

Applying the change of enrichment along $\Fst$ to $\theta$ yields the following $\psu$-natural transformation (\cref{expl:perm-nattr}).
\[\begin{tikzpicture}
\def\t{28}
\draw[0cell]
(0,0) node (a) {\phantom{C}}
(a)+(1.8,0) node (b) {\phantom{C}}
(a)+(-1,0) node (a') {\big(\C_{\Fst} \scs \mcomp_{\Fst} \scs i_{\Fst}\big)}
(b)+(1,0) node (b') {\big(\D_{\Fst} \scs \mcomp_{\Fst} \scs i_{\Fst}\big)}
;
\draw[1cell=.8]
(a) edge[bend left=\t] node {H_{\Fst}} (b)
(a) edge[bend right=\t] node[swap] {G_{\Fst}} (b)
;
\draw[2cell]
node[between=a and b at .42, rotate=-90, 2label={above,\theta_{\Fst}}] {\Rightarrow}
;
\end{tikzpicture}\]
For each object $x \in \C_{\Fst}$, its $x$-component is the $\obsim$-equivalence class \cref{mnat-change-enrichment}
\[(\theta_{\Fst})_x = [(\theta_x)] \in \Fst\D(Hx,Gx)\]
of the length-one sequence $(\theta_x) \in \Fr\D(Hx,Gx)$.  The reasoning is the same as in \cref{FstzeroFstix}, with $i_x$ replaced by $\theta_x$.

This finishes our description of the change-of-enrichment 2-functor $\dFst$.  We will use $\dFst$ in \cref{ex:Fstse,expl:Fstdgr,expl:Fstdgr-object,expl:Fstdgr-morphism}.
\end{explanation}

\section{Preservation of Opposite Enriched Categories}
\label{sec:pres-opposite}

By \cref{mult-change-enrichment} each non-symmetric multifunctor $F$ yields a change-of-enrichment 2-functor $\dF$.  Recall that a multifunctor preserves colored units, composition, \emph{and} symmetric group action (\cref{def:enr-multicategory-functor}).  In this section we observe that change of enrichment along a multifunctor preserves opposite enriched categories (\cref{def:opposite-mcat}); see \cref{dF-opposite}.  We will use this opposite-preservation property in \cref{presheaf-change-enr} to define the presheaf change of enrichment of $F$.  In that context, we will use \cref{dF-opposite} in \cref{gspectra-thm-v-cor,gspectra-thm-vii-cor,rk:Kemdg-opposite,mackey-xiv-cor}.

\begin{proposition}\label{dF-opposite}\index{opposite!enriched category}\index{enriched category!opposite}
Suppose given a multifunctor between multicategories
\[F \cn (\M,\ga,\opu) \to (\N,\ga,\opu)\]
and an $\M$-category $\C$.  Then there is an equality of $\N$-categories
\[(\Cop)_F = (\C_F)^\op\]
with $\dF$ the change of enrichment in \cref{mcat-change-enrichment}.
\end{proposition}
\begin{proof}
For an $\M$-category $(\C,\mcomp,i)$, we need to show that the following two $\N$-categories are the same.
\begin{itemize}
\item $(\Cop)_F$ is the change of enrichment along $F$ \cref{mcat-change-enrichment} of the opposite $\M$-category $\Cop$ (\cref{def:opposite-mcat}).
\item $(\C_F)^\op$ is the opposite $\N$-category of $\C_F$, which is the change of enrichment of $\C$ along $F$.
\end{itemize}
In both $(\Cop)_F$ and $(\C_F)^\op$, the objects are those of $\C$.  For objects $x,y \in \C$, the hom objects are equal:
\[(\Cop)_F(x,y) = F\Cop(x,y) = F\C(y,x) = \C_F(y,x) = (\C_F)^\op(x,y).\]
By \cref{Cop-identity,iFx}, the identity of an object $x$ in $(\Cop)_F$ and $(\C_F)^\op$ is the nullary multimorphism
\[F(i_x) \cn \ang{} \to F\C(x,x) \inspace \N.\]
For objects $x,y,z \in (\Cop)_F$ and $(\C_F)^\op$, the compositions are equal:
\begin{equation}\label{CopF-composition}
\begin{aligned}
\big((\mcompop)_F\big)_{x,y,z} 
&= F(\mcomp_{z,y,x} \cdot \tau)\\
&= (F\mcomp_{z,y,x}) \cdot \tau\\
&= \big((\mcomp_F)^\op\big)_{x,y,z}.
\end{aligned}
\end{equation}
The first and third equalities above use the definitions \cref{mcompop-def,mFxyz}.  The second equality uses the fact that $F$ is a multifunctor, which preserves the symmetric group action \cref{enr-multifunctor-equivariance}.
\end{proof}

We emphasize that \cref{dF-opposite} does \emph{not} extend to non-symmetric multifunctors between multicategories because the second equality in \cref{CopF-composition} requires that $F$ preserves the symmetric group action in the strict sense.

\begin{example}[$K$-Theory Multifunctors]\label{ex:dF-opposite}\index{K-theory@$K$-theory!change of enrichment along}\index{Elmendorf-Mandell!K-theory@$K$-theory!change of enrichment along}\index{Elmendorf-Mandell!J-theory@$J$-theory!change of enrichment along}
\cref{dF-opposite} applies to the multifunctors in the following sub-diagram of \cref{mult-change-enr-diagram}.
\begin{equation}\label{ex:dF-opp-diagram}
\begin{tikzpicture}[vcenter]
\def\h{3} \def\v{-1.4} \def\w{.6}
\draw[0cell=.9]
(0,0) node (p) {\permcatsu}
(p)+(-\h,0) node (m) {\Multicat}
(p)+(0,\v) node (mone) {\MoneMod}
(m)+(0,\v) node (pm) {\pMulticat}
(p)+(\h,0) node (sp) {\Sp}
(sp)+(0,\v) node (gs) {\Gstarcat}
;
\draw[1cell=.9]
(p) edge node[swap] {\End} (m)
(p) edge[shorten <=.5ex] node[swap,pos=.6] {\Endst} (pm)
(p) edge node[pos=.7] {\Endm} (mone)
(p) edge[shorten <=1ex] node[pos=.6] {\Jem} (gs)
(p) edge node {\Kem} (sp)
;
\end{tikzpicture}
\end{equation}
In other words, the change of enrichment induced by each multifunctor in \cref{ex:dF-opp-diagram} preserves opposite enriched categories.  For example, consider a $\permcatsu$-category $\C$ (\cref{expl:perm-enr-cat}).  Then \cref{dF-opposite} yields the equality of $\Sp$-categories
\[(\Cop)_{\Kem} = (\C_{\Kem})^\op.\]
The opposite on the left-hand side is taken in $\permcatsucat$.  The opposite on the right-hand side is taken in $\Spcat$, after the change of enrichment along $\Kem$.  The above equality of $\Sp$-categories is important in the context of presheaf change of enrichment; see \cref{thm:Kemdg,rk:Kemdg-opposite,rk:BO7.5}.

Note that $\Fr$, $\Fst$, $\Fm$, and $\cP$ in \cref{mult-change-enr-diagram} are excluded from the diagram \cref{ex:dF-opp-diagram} because they are \emph{non-symmetric} multifunctors and do not preserve opposite enriched categories in general.  Specifically, each of these non-symmetric multifunctors fails to satisfy the second equality in \cref{CopF-composition}, which requires strict preservation of symmetric group action.
\end{example}

\section{Change of Enrichment along a Monoidal Functor}
\label{sec:change-enr-monfunctor}

For a monoidal functor between monoidal categories (\cref{def:monoidalfunctor})
\[U \cn \V \to \W,\]
there are two change-of-enrichment 2-functors as follows. 
\begin{enumerate}
\item By \cref{proposition:U-VCat-WCat} there is a change-of-enrichment 2-functor
\[\dU \cn \VCat \to \WCat.\]
\item By \cref{ex:endc} there is a non-symmetric multifunctor
\[\End\,U \cn \End\,\V \to \End\,\W.\]
By \cref{EndV-enriched,mult-change-enrichment} there is a change-of-enrichment 2-functor
\[\dEndU \cn \EndVCat = \VCat \to \EndWCat = \WCat.\]
\end{enumerate}
In this section we observe that these two change-of-enrichment 2-functors are equal.

\begin{proposition}\label{mon-change-enrichment}\index{change of enrichment}\index{enriched category!change of enrichment}\index{monoidal functor!change of enrichment}\index{symmetric monoidal functor!change of enrichment}\index{2-functor!change of enrichment}
For each monoidal functor between monoidal categories
\[(U,U^2,U^0) \cn (\V,\otimes,\tu) \to (\W,\otimes,\tu),\]
there is an equality of change-of-enrichment 2-functors
\[\dU = \dEndU \cn \VCat \to \WCat.\]
\end{proposition}
\begin{proof}
We need to show that $\dU$ and $\dEndU$ are equal on the objects, 1-cells, and 2-cells of the 2-category (\cref{EndV-enriched})
\[\VCat = \EndVCat.\]  
First we consider a $\V$-category $(\C,\mcomp,i)$ (\cref{def:enriched-category}).
\begin{itemize}
\item By \cref{def:mult-change-enr,definition:U-VCat-WCat}, both $\C_U$ and $\C_{\EndU}$ have
\begin{itemize}
\item the same objects as $\C$ and
\item hom objects $U\C(x,y)$ for objects $x,y \in \C$.
\end{itemize} 
\item The identity of an object $x \in \C_U$ is given by the composite \cref{eq:id-U}
\[\tu \fto{U^0} U\tu \fto{Ui_x} U\C(x,x).\]
Regarding $x$ as an object in $\C_{\EndU}$, \cref{iFx} implies that its identity is also given by the composite above because $\EndU$ on a nullary multimorphism \cref{PtwoPf} is the composite of the unit constraint $U^0$ followed by $U(-)$.
\item For objects $x,y,z \in \C_U$ the composition is given by the composite \cref{eq:comp-U} 
\[U\C(y,z) \otimes U\C(x,y) \fto{U^2} U\big(\C(y,z) \otimes \C(x,y)\big) \fto{U\mcomp} U\C(x,z).\]
Regarding $x,y$, and $z$ as objects in $\C_{\EndU}$, \cref{mFxyz} implies that the composition is also given by the composite above because $\EndU$ on a binary multimorphism \cref{PtwoPf} is the composite of the monoidal constraint $U^2$ followed by $U(-)$.
\end{itemize}
This proves that $\dU$ and $\dEndU$ are equal on objects.

Next we consider a $\V$-functor between $\V$-categories (\cref{def:enriched-functor})
\[H \cn \C \to \D.\] 
By \cref{HFxy,FUxy,PtwoPf}, for objects $x,y \in \C_U$ and $\C_{\EndU}$, there are equalities of component morphisms in $\W$:
\[(H_U)_{x,y} = U(H_{x,y}) = (\EndU)(H_{x,y}) = (H_{\EndU})_{x,y}.\]
Thus $\dU$ and $\dEndU$ are equal on 1-cells.

By \cref{mnat-change-enrichment,thetaUx}, the proof that $\dU$ and $\dEndU$ are equal on each $\V$-natural transformation $\theta$ is the same as the argument above for the identity of an object, with $i_x$ replaced by the component $\theta_x$ for objects $x \in \C$.
\end{proof}

\begin{example}[$K$-Theory Symmetric Monoidal Functors]\label{ex:mon-change-enr}
\cref{mon-change-enrichment} applies to the symmetric monoidal functors in the following commutative diagram from \cref{Ust,Umonemod,smastar-istar,eq:Ksummary}.
\begin{equation}\label{ex:mon-change-enr-diag}
\begin{tikzpicture}[baseline={(pm.base)}]
\def\v{-1.3} \def\h{2.5}
\draw[0cell=.9]
(0,0) node (m) {\Multicat}
(m)+(0,\v) node (pm) {\pMulticat}
(pm)+(0,\v) node (mone) {\MoneMod}
(pm)+(\h,0) node (gac) {\Gacat}
(gac)+(\h,0) node (gas) {\Gasset}
(gas)+(\h,0) node (sp) {\Sp}
(gac)+(0,\v) node (gsc) {\Gstarcat}
(gas)+(0,\v) node (gss) {\Gstarsset}
;
\draw[1cell=.9]
(mone) edge node {\Um} (pm)
(pm) edge node {\Ust} (m)
(mone) edge node {\Jt} (gsc)
(gsc) edge node {\Ner_*} (gss)
(gss) edge node[swap,pos=.4] {\Kg} (sp)
(gac) edge node {\Ner_*} (gas)
(gas) edge node {\Kf} (sp)
(gac) edge node[swap,pos=.4] {\smastar} (gsc)
(gas) edge node[swap,pos=.4] {\smastar} (gss)
;
\end{tikzpicture}
\end{equation}
In other words, for each of these symmetric monoidal functors, the induced change-of-enrichment 2-functors in \cref{proposition:U-VCat-WCat,mult-change-enrichment} are the same.  Moreover, each of these change-of-enrichment 2-functors preserves opposite enriched categories in the sense of \cref{dF-opposite}.  
\end{example}

\section{Composition of Change-of-Enrichment 2-Functors}
\label{sec:functoriality-change-enr}

By \cref{mult-change-enrichment} each non-symmetric multifunctor $F$ has an associated change-of-enrichment 2-functor $\dF$.  In this section we observe that change of enrichment respects composition of non-symmetric multifunctors; see \cref{func-change-enr}.  This is an analog of \cref{proposition:change-enr-horiz-comp} for change of enrichment along monoidal functors.  \cref{func-change-enr} is important in the context of standard enrichment, diagram change of enrichment, and Mackey functor change of enrichment; see \cref{gspectra-thm-iv,EndmseKNJ,gspectra-thm-vii,gspectra-thm-vii-cor,Fdg-Edg,EdgFdgr-A-diag,Fmdgr-factorization}.

\begin{proposition}\label{func-change-enr}\index{change of enrichment!composition of}\index{composition!- of change of enrichment}
For non-symmetric multifunctors between non-symmetric multicategories
\[\M \fto{F} \N \fto{G} \P,\]
the following diagram of change-of-enrichment 2-functors commutes.
\[\begin{tikzpicture}[vcenter]
\def\h{2.5} \def\w{.6}
\draw[0cell=1]
(0,0) node (a) {\MCat}
(a)+(\h,0) node (b) {\NCat}
(b)+(\h,0) node (c) {\PCat}
;
\draw[1cell=.9]
(a) edge node {\dF} (b)
(b) edge node {\dG} (c)
;
\draw[1cell=.9]
(a) [rounded corners=3pt] |- ($(b)+(-1,\w)$)
-- node {\dGF} ($(b)+(1,\w)$) -| (c)
;
\end{tikzpicture}\]
\end{proposition}

\begin{proof}
We need to show that $\dGF$ and $\dG \circ \dF$ are equal on the objects, 1-cells, and 2-cells of the 2-category $\MCat$ (\cref{mcat-iicat}).

First we consider an $\M$-category $(\C,\mcomp,i)$ (\cref{def:menriched-cat}).
\begin{itemize}
\item Both $(\C_F)_G$ and $\C_{GF}$ have the same objects as $\C$ and hom objects
\[GF\C(x,y) \in \P \forspace x,y \in \C.\]
\item For each object $x \in (\C_F)_G$ and $\C_{GF}$, by \cref{iFx,composite-multifunctor} the identity is the nullary multimorphism
\[\ang{} \fto{GF(i_x)} GF\C(x,x) \inspace \P.\]
\item For objects $x,y,z  \in (\C_F)_G$ and $\C_{GF}$, by \cref{mFxyz,composite-multifunctor} the composition is the binary multimorphism
\[\big(GF\C(y,z) \scs GF\C(x,y)\big) \fto{GF(\mcomp_{x,y,z})} GF\C(x,z) \inspace \P.\]
\end{itemize}
This shows that $\dGF$ and $\dG \circ \dF$ are equal on $\M$-categories.

Next we consider an $\M$-functor $H \cn \C \to \D$ between $\M$-categories (\cref{def:mfunctor}).
\begin{itemize}
\item Both $(H_F)_G$ and $H_{GF}$ have the same object assignment as $H$.
\item For objects $x,y \in (\C_F)_G$ and $\C_{GF}$, by \cref{HFxy} there are equalities of component unary multimorphisms
\[\big((H_F)_G\big)_{x,y} = GF(H_{x,y}) = (H_{GF})_{x,y} \inspace \P.\]
\end{itemize}
This shows that $\dGF$ and $\dG \circ \dF$ are equal on $\M$-functors.

For an $\M$-natural transformation $\theta$ and an object $x \in (\C_F)_G$ and $\C_{GF}$, by \cref{mnat-change-enrichment} there are equalities of $x$-component nullary multimorphisms
\[\big((\theta_F)_G\big)_x = GF(\theta_x) = (\theta_{GF})_x \inspace \P.\]
Thus $\dGF$ and $\dG \circ \dF$ are equal on $\M$-natural transformations.
\end{proof}

We emphasize that \cref{func-change-enr} does \emph{not} require $F$ and $G$ to preserve the symmetric group action even if $\M$, $\N$, and $\P$ are multicategories.

\begin{example}[$K$-Theory Multifunctors]\label{ex:func-change-enr}
\cref{func-change-enr} applies to all the (non-symmetric) multifunctors in \cref{mult-change-enr-diagram,ex:mon-change-enr-diag}.  For example, consider the following commutative sub-diagram of \cref{eq:Ksummary} consisting of multifunctors.
\begin{equation}\label{Kem-diagram}
\begin{tikzpicture}[xscale=1.2,yscale=1,vcenter]
\def\t{.35} \def\v{-1.5} \def\u{.6} \def\w{-.6}
\draw[0cell=.8]
(0,0) node (a1) {\permcatsu}
(a1)+(3,0) node (a2) {\Gacat}
(a2)+(2,0) node (a3) {\Gasset}
(a3)+(2,0) node (a4) {\Sp}
(a2)+(0,\v) node (b2) {\Gstarcat}
(a3)+(0,\v) node (b3) {\Gstarsset}
(b2)+(-2.3,0) node (b1) {\MoneMod}
;
\draw[1cell=.85]
(a2) edge node {\Ner_*} (a3)
(a3) edge node {\Kf} (a4)
(a1) edge node[pos=.4] {\Jem} (b2)
(a1) edge node[pos=.85] {\Endm} (b1)
(b1) edge node[pos=.4] {\Jt} (b2)
(b2) edge node {\Ner_*} (b3)
(b3) edge node[pos=.4,swap] {\Kg} (a4)
(a2) edge node[pos=.4] {\sma^*} (b2)
(a3) edge node[pos=.4] {\sma^*} (b3)
;
\draw[1cell=.85]
(a1) [rounded corners=3pt] |- ($(b2)+(0,\w)$)
-- node[pos=\t] {\Kem} ($(b3)+(0,\w)$) -| (a4)
;
\end{tikzpicture}
\end{equation}
By \cref{func-change-enr} the associated diagram of change-of-enrichment 2-functors is commutative.  In particular, the change-of-enrichment 2-functor along $\Kem$ is equal to the composite of the change-of-enrichment 2-functors along $\Endm$, $\Jt$, $\Ner_*$, and $\Kg$.  We discuss the factorization \cref{Kem-diagram} of $\Kem$ further
\begin{itemize}
\item in \cref{gspectra-thm-xi} in the context of standard enrichment and
\item in \cref{gspectra-thm-xiv} in the context of presheaf change of enrichment.\defmark
\end{itemize} 
\end{example}

\section{2-Functoriality of Change of Enrichment}
\label{sec:twofunc-change-enr}

In this section we show that the (change of) enrichment constructions
\[\M \mapsto \MCat \andspace F \mapsto \dF\]
in \cref{mcat-iicat,mult-change-enrichment}, respectively, are part of a 2-functor; see \cref{change-enr-twofunctor}.  To explain this precisely, first we construct the change of enrichment of multinatural transformations (\cref{def:multinat-change-enr}).  We show in \cref{dtheta-twonat} that this yields a well-defined 2-natural transformation.  

Beyond its immediate usage in \cref{change-enr-twofunctor}, we will use \cref{dtheta-twonat} in \cref{Fdgr-def} to construct a functor that goes in the backward direction as a diagram change of enrichment.  Then we use that functor to construct an equivalence of homotopy theories in that context; see \cref{mackey-gen-xiv}.

Recall 2-natural transformations and multinatural transformations in \cref{def:twonaturaltr,def:enr-multicat-natural-transformation}, respectively.

\begin{definition}\label{def:multinat-change-enr}
Suppose $\theta \cn F \to G$ is a multinatural transformation between non-symmetric multifunctors between non-symmetric multicategories, as in the left diagram below.  
\begin{equation}\label{dtheta-def}
\begin{tikzpicture}[baseline={(a.base)}]
\def\s{25} \def\t{24}
\draw[0cell]
(0,0) node (a) {\phantom{M}}
(a)+(1.7,0) node (b) {\phantom{M}}
(a)+(-.45,0) node (a') {(\M,\ga,\opu)}
(b)+(.4,0) node (b') {(\N,\ga,\opu)}
;
\draw[1cell=.8]
(a) edge[bend left=\s] node {F} (b)
(a) edge[bend right=\s] node[swap] {G} (b)
;
\draw[2cell]
node[between=a and b at .46, rotate=-90, 2label={above,\theta}] {\Rightarrow}
;
\begin{scope}[shift={(4.5,0)}]
\draw[0cell]
(0,0) node (a) {\phantom{\M}}
(a)+(2.3,0) node (b) {\phantom{\N}}
(a)+(-.3,0) node (a') {\MCat}
(b)+(.3,0) node (b') {\NCat}
;
\draw[1cell=.9]
(a) edge[bend left=\t] node[pos=.53] {\dF} (b)
(a) edge[bend right=\t] node[swap,pos=.53] {\dG} (b)
;
\draw[2cell]
node[between=a and b at .39, rotate=-90, 2label={above,\dtheta}] {\Rightarrow}
;
\end{scope}
\end{tikzpicture}
\end{equation}
We define the data of a 2-natural transformation $\dtheta$, as in the right diagram in \cref{dtheta-def}, as follows.  Here 
\begin{itemize}
\item $\MCat$ and $\NCat$ are the 2-categories in \cref{mcat-iicat}, and
\item $\dF$ and $\dG$ are the change-of-enrichment 2-functors in \cref{mult-change-enrichment}.\end{itemize}
For each small $\M$-category $(\C,\mcomp,i)$ (\cref{def:menriched-cat}), we define the data of an $\N$-functor (\cref{def:mfunctor})
\begin{equation}\label{C-sub-theta}
\C_\theta \cn (\C_F,\mcomp_F,i_F) \to (\C_G,\mcomp_G,i_G)
\end{equation}
as follows.
\begin{itemize}
\item $\C_\theta$ is the identity assignment on objects.  This is well defined since $\C_F$ and $\C_G$ both have the same objects as $\C$.
\item For each pair of objects $x,y \in \C$, the $(x,y)$-component of $\C_\theta$ is defined as the $\C(x,y)$-component of $\theta$:
\begin{equation}\label{Ctheta-ab}
\C_F(x,y) = F\C(x,y) \fto{(\C_\theta)_{x,y} = \theta_{\C(x,y)}} \C_G(x,y) = G\C(x,y).
\end{equation}
This is a unary multimorphism in $\N$.
\end{itemize}
This finishes the definition of $\dtheta$.
\end{definition}

Now we check that $\dtheta$ is well defined.

\begin{proposition}\label{dtheta-twonat}
In the context of \cref{def:multinat-change-enr}, $\dtheta$ is a 2-natural transformation.
\end{proposition}

\begin{proof}
We prove statements \cref{Ctheta-i,Ctheta-ii,Ctheta-iii} below.
\begin{romenumerate}
\item\label{Ctheta-i}
$\C_\theta$ in \cref{C-sub-theta} is an $\N$-functor (\cref{def:mfunctor}).
\item\label{Ctheta-ii}
$\dtheta$ is natural with respect to $\M$-functors \cref{onecellnaturality}.
\item\label{Ctheta-iii}
$\dtheta$ is natural with respect to $\M$-natural transformations \cref{twocellnaturality}.
\end{romenumerate}

\medskip
\emph{Statement \cref{Ctheta-i}}.
We abbreviate $\C(x,y)$ to $\C_{x,y}$.  For objects $x,y,z \in \C_F$, the two compatibility diagrams in \cref{mfunctor-diagrams} for $\C_\theta$ are as follows.
\[\begin{tikzpicture}[vcenter]
\def\v{-1.4} \def\h{3}
\draw[0cell=.9]
(0,0) node (a) {\big(F\C_{y,z} \scs F\C_{x,y}\big)}
(a)+(\h,0) node (b) {F\C_{x,z}}
(a)+(0,\v) node (c) {\big(G\C_{y,z} \scs G\C_{x,y}\big)}
(b)+(0,\v) node (d) {G\C_{x,z}}
;
\draw[1cell=.85]
(a) edge node {F\mcomp_{x,y,z}} (b)
(b) edge node {\theta_{\C_{x,z}}} (d)
(a) edge node[swap] {(\theta_{\C_{y,z}} \scs \theta_{\C_{x,y}})} (c)
(c) edge node {G\mcomp_{x,y,z}} (d)
;
\begin{scope}[shift={(4.5,0)}]
\draw[0cell=.9]
(0,0) node (a) {\ang{}}
(a)+(1.7,0) node (b) {F\C_{x,x}}
(b)+(0,\v) node (c) {G\C_{x,x}}
;
\draw[1cell=.85]
(a) edge node {Fi_x} (b)
(b) edge node {\theta_{\C_{x,x}}} (c)
(a) edge node[swap] {Gi_x} (c)
;
\end{scope}
\end{tikzpicture}\]
These two diagrams commute by the naturality of $\theta$ \cref{enr-multinat} for, respectively, the binary multimorphism $\mcomp_{x,y,z}$ and the nullary multimorphism $i_x$ in $\M$.

\medskip
\emph{Statement \cref{Ctheta-ii}}.
For an $\M$-functor $P \cn \C \to \D$ (\cref{def:mfunctor}), the 1-cell naturality diagram \cref{onecellnaturality} for $\dtheta$ is the left diagram of $\N$-functors below, where $P_F$ and $P_G$ are defined in \cref{mfunctor-change-enrichment}.
\[\begin{tikzpicture}
\def\v{-1.4}
\draw[0cell=.9]
(0,0) node (a1) {\C_F}
(a1)+(2,0) node (a2) {\C_G}
(a1)+(0,\v) node (b1) {\D_F}
(a2)+(0,\v) node (b2) {\D_G}
;
\draw[1cell=.85]
(a1) edge node {\C_\theta} (a2)
(a2) edge node {P_G} (b2)
(a1) edge node[swap] {P_F} (b1)
(b1) edge node {\D_\theta} (b2)
;
\begin{scope}[shift={(4,0)}]
\draw[0cell=.9]
(0,0) node (a1) {F\C_{x,y}}
(a1)+(3,0) node (a2) {G\C_{x,y}}
(a1)+(0,\v) node (b1) {F\D_{Px,Py}}
(a2)+(0,\v) node (b2) {G\D_{Px,Py}}
;
\draw[1cell=.85]
(a1) edge node {\theta_{\C_{x,y}}} (a2)
(a2) edge node {GP_{x,y}} (b2)
(a1) edge node[swap] {FP_{x,y}} (b1)
(b1) edge node {\theta_{\D_{Px,Py}}} (b2)
;
\end{scope}
\end{tikzpicture}\]
Both $P_G \C_\theta$ and $\D_\theta P_F$ have the same object assignment as $P$, since $\C_\theta$ and $\D_\theta$ are the identity functions on objects.  For objects $x,y \in \C_F$, the $(x,y)$-component of the left diagram above is the right diagram in $\N$.  The latter commutes by the naturality of $\theta$ \cref{enr-multinat} for the unary multimorphism 
\[P_{x,y} \cn \C(x,y) \to \D(Px,Py) \inspace \M.\]

\medskip
\emph{Statement \cref{Ctheta-iii}}.
Consider an $\M$-natural transformation $\psi \cn P \to Q$ (\cref{def:mnaturaltr}) for $\M$-functors $P,Q \cn \C \to \D$.  The 2-cell naturality diagram \cref{twocellnaturality} for $\dtheta$ is the left diagram of $\N$-natural transformations below, where $\psi_F$ and $\psi_G$ are defined in \cref{mnat-change-enrichment}.
\[\begin{tikzpicture}
\def\v{-1.4} \def\s{45}
\draw[0cell=.9]
(0,0) node (a1) {\C_F}
(a1)+(2.5,0) node (a2) {\C_G}
(a1)+(0,\v) node (b1) {\D_F}
(a2)+(0,\v) node (b2) {\D_G}
;
\draw[1cell=.8]
(a1) edge node {\C_\theta} (a2)
(b1) edge node[swap] {\D_\theta} (b2)
(a1) edge[bend right=\s] node[swap,pos=.4] {P_F} (b1)
(a1) edge[bend left=\s] node[pos=.6] {Q_F} (b1)
(a2) edge[bend right=\s] node[swap,pos=.4] {P_G} (b2)
(a2) edge[bend left=\s] node[pos=.6] {Q_G} (b2)
;
\draw[2cell]
node[between=a1 and b1 at .6, 2label={above,\psi_F}] {\Rightarrow}
node[between=a2 and b2 at .6, 2label={above,\psi_G}] {\Rightarrow}
;
\begin{scope}[shift={(5,0)}]
\draw[0cell=.9]
(0,0) node (a) {\ang{}}
(a)+(0,\v) node (b) {F\D_{Px,Qx}}
(b)+(2.5,0) node (c) {G\D_{Px,Qx}}
;
\draw[1cell=.85]
(a) edge node[swap] {F\psi_x} (b)
(b) edge node {\theta_{\D_{Px,Qx}}} (c)
(a) edge[bend left=20] node {G\psi_x} (c)
;
\end{scope}
\end{tikzpicture}\]
For each object $x \in \C_F$, since $\C_\theta$ is the identity on objects, the $x$-component of the left diagram above is the right diagram in $\N$.  The latter commutes by the naturality of $\theta$ \cref{enr-multinat} for the nullary multimorphism
\[\psi_x \cn \ang{} \to \D(Px,Qx) \inspace \M.\]
This finishes the proof.
\end{proof}

To state the main result of this section, recall the following 2-categories. 
\begin{itemize}
\item $\iicat$ is the 2-category of small 2-categories, 2-functors, and 2-natural transformations (\cref{ex:iicat}).
\item $\Multicatns$ is the 2-category of non-symmetric small multicategories, multifunctors, and multinatural transformations \cref{multicat-def}.
\end{itemize} 
We now observe that change of enrichment is a 2-functor between these 2-categories.  The following result is a multicategorical analog of \cite[2.2.7]{cerberusIII}, which is about enrichment in monoidal categories.  See \cref{expl:change-enr-size} for a discussion related to universes.

\begin{theorem}\label{change-enr-twofunctor}\index{change of enrichment!2-functoriality of}\index{2-functor!defined by change of enrichment}
There is a 2-functor
\[\Enr \cn \Multicatns \to \iicat\]
given by the assignments
\begin{itemize}
\item $\Enr\M = \MCat$ (\cref{mcat-iicat}) on objects,
\item $\Enr F = \dF$ (\cref{mult-change-enrichment}) on 1-cells, and
\item $\Enr\theta = \dtheta$ (\cref{dtheta-twonat}) on 2-cells.
\end{itemize}
\end{theorem}

\begin{proof}
We prove statements \cref{Enr-i,Enr-ii,Enr-iii,Enr-iv} below.
\begin{romenumerate}
\item\label{Enr-i}
$\Enr$ preserves identity 1-cells and horizontal composition of 1-cells.
\item\label{Enr-ii}
$\Enr$ preserves identity 2-cells.
\item\label{Enr-iii}
$\Enr$ preserves vertical composition of 2-cells.
\item\label{Enr-iv}
$\Enr$ preserves horizontal composition of 2-cells.
\end{romenumerate}

\medskip
\emph{Statement \cref{Enr-i}}. 
For a small non-symmetric multicategory $\M$, the identity non-symmetric multifunctor $1_\M$ (\cref{def:enr-multicategory-functor}) consists of the identity functions on objects and multimorphisms.  The change of enrichment along $1_\M$,
\[\dsub{1_{\M}} \cn \MCat \to \MCat,\]
is the identity 2-functor on $\MCat$ by the componentwise definitions \cref{CFxy,mFxyz,iFx,HFxy,mnat-change-enrichment}.  Moreover, $\Enr$ preserves horizontal composition of 1-cells by \cref{func-change-enr}.

\medskip
\emph{Statement \cref{Enr-ii}}. 
For a non-symmetric multifunctor $F \cn \M \to \N$, the identity multinatural transformation $1_F$ (\cref{def:enr-multicat-natural-transformation}) has, for each object $x \in \M$, $x$-component given by the colored unit $1_{Fx}$.  By definitions \cref{id-mfunctor,Ctheta-ab}, each component of the 2-natural transformation 
\[\dsub{1_{F}} \cn \dF \to \dF\] 
is an identity $\N$-functor, so $\dsub{1_{F}} = 1_{\dF}$.

\medskip
\emph{Statement \cref{Enr-iii}}. 
We consider vertically composable multinatural transformations between non-symmetric multifunctors, as in the left diagram below.
\[
\]
The desired equality of 2-natural transformations is
\[\dpsi \dtheta = \dpsitha,\]
where $\psi\theta \cn F \to H$ is the vertical composition \cref{multinatvcomp}.  We need to show that, for each small $\M$-category $\C$, the following diagram of $\N$-functors commutes.
\[%
\]
On objects of $\C_F$, each of the three arrows in the diagram above is the identity function.  For objects $x,y \in \C_F$, the $(x,y)$-component is the following diagram of unary multimorphisms in $\N$.
\[%
\]
This diagram commutes because $\psi\theta$ is defined componentwise \cref{multinat-vcomp-def}.

\medskip
\emph{Statement \cref{Enr-iv}}. 
We consider horizontally composable multinatural transformations $\theta$ and $\theta'$ between non-symmetric multifunctors \cref{multinathcomp}, as in the left diagram below.
\[%
\]
This yields the following 2-natural transformations.
\[%
\]
The desired equality of 2-natural transformations is
\[\dthap * \dtheta = \dsub{\theta' \!* \theta}.\]
For a small $\M$-category $\C$, the $\C$-component of $\dthap * \dtheta$ is the bottom composite $\P$-functor in the following diagram, while the $\C$-component of $\dsub{\theta' \!* \theta}$ is the top $\P$-functor.
\[%
\]
Each of these two $\P$-functors is the identity on objects.  For objects $x,y \in \C_{F'F}$, the $(x,y)$-components yield the following diagram of unary multimorphisms in $\P$.
\[%
\]
This diagram commutes by the definition \cref{multinat-hcomp-def} of each component of $\theta' \!* \theta$.
\end{proof}

\begin{explanation}[Universes]\label{expl:change-enr-size}\index{universe}\index{Axiom of Universes}\index{convention!universe}
In \cref{change-enr-twofunctor}, for a small non-symmetric multicategory $\M$ with respect to a given universe $\calu$, the 2-category $\MCat$ is, in general, not small with respect to $\calu$.  We implicitly use Grothendieck's Axiom of Universes (\cref{conv:universe}) to choose a larger universe $\calu'$ such that, for each $\calu$-small non-symmetric multicategory $\M$, $\Ob(\MCat)$ is a member of $\calu'$.  Then we consider the 2-category $\iicat'$, whose objects are $\calu'$-small 2-categories.  The precise codomain of $\Enr$ in \cref{change-enr-twofunctor} is $\iicat'$.  For related discussion in the context of enrichment in monoidal categories, the reader is referred to \cite[2.2.6 and 2.2.8]{cerberusIII}.
\end{explanation}

\chapter{The Closed Multicategory of Permutative Categories}
\label{ch:gspectra}
This chapter defines and develops the basic properties of closed multicategories $\M$.
This is similar to, but more general than, the concept of symmetric monoidal closed categories $\V$ from \cref{def:closedcat}.
For the special case $\M = \End\,\V$, \cref{smclosed-closed-multicat} shows that the two notions agree.
\cref{sec:internal-hom-permcat,sec:multiev-permcat,sec:permcatsu-clmulti} develop the special case of $\permcatsu$, the multicategory of permutative categories.
This, along with the closed symmetric monoidal categories $\pMulticat$ and $\MoneMod$, will be the main case of interest for applications in \cref{ch:mackey,ch:mackey_eq}.

More general closed structures for
$\permcat$, $\permcatst$, and $\permcatsus$ are discussed in \cref{sec:lax-strong-multicat}.
These structures involve notions of lax, respectively strong, respectively strictly unital strong, multilinear functors.

\subsection*{Connection with Other Chapters}

\cref{ch:std_enrich} develops the general theory of self-enrichment for closed multicategories, with $\permcatsu$ being one of the key examples for further applications.
\cref{ch:gspectra_Kem} develops the theory of enriched diagrams and enriched Mackey functors in a closed multicategory $\M$.
\cref{ch:mackey_eq} applies theory from \cref{ch:mackey} to $\permcatsu$-enriched categories and diagrams.

\subsection*{Background}
The $\Cat$-multicategory structure for $\permcatsu$ is discussed in \cref{sec:multpermcat}.
The self-enrichment of $\permcatsu$ is discussed in \cref{sec:perm-self-enr}.

\subsection*{Chapter Summary}

\cref{sec:closed-multicat} gives the basic definitions for closed multicategories.
\cref{sec:internal-hom-permcat} describes the internal hom for $\permcatsu$, and \cref{sec:multiev-permcat} describes the multicategorical evaluation.
In \cref{sec:permcatsu-clmulti} these are combined to give the closed multicategory structure for $\permcatsu$.
\cref{sec:lax-strong-multicat} discusses how the previous structures can be generalized to $\permcat$ and other multicategories of permutative categories.
Here is a summary table.
\reftable{.9}{
  definition of a closed multicategory
  & \ref{def:closed-multicat}
  \\ \hline
  internal hom for $\permcatsu$
  & \ref{def:clp-angcd} and \ref{clp-angcd-permutative}
  \\ \hline
  symmetric group action on internal hom
  & \ref{def:clp-ancd-sigma} and \ref{clp-angcd-strict}
  \\ \hline
  multicategorical $\ev$ and $\chi$
  & \ref{def:clp-eval}, \ref{clp-evaluation}, \ref{expl:domain-chi}, and \ref{expl:clp-partners}
  \\ \hline
  multicategorical evaluation axioms
  & \ref{clp-evbij-axiom} and \ref{clp-evbij-equiv}
  \\ \hline
  lax multilinear functors
  & \ref{def:lax-nlinearfunctor}, \ref{lemma:lax-nlinear-tensor}, \ref{proposition:laxmultilin-tensor}, and \ref{theorem:permcat-catmulticat}
  \\ \hline
  closed structure in lax case
  & \ref{definition:laxinternalhom}, \ref{def:lax-clp-ancd-sigma}, \ref{def:lax-clp-eval}, and \ref{lax-clp-partners}
  \\ \hline
  main results
  & \ref{perm-closed-multicat} and \ref{lax-perm-closed-multicat}
  \\
}

We remind the reader of \cref{conv:universe} about universes and \cref{expl:leftbracketing} about left normalized bracketing for iterated products.

\section{Closed Multicategories}
\label{sec:closed-multicat}

Recall from \cref{def:enr-multicategory} that a \emph{multicategory} means a $\Set$-multicategory for the symmetric monoidal category $(\Set,\times)$ of sets and functions with the Cartesian product as the monoidal product.  In this section we define \emph{closed multicategories}, which are multicategories equipped with internal hom objects and compatible evaluations.  This concept provides a common setting for 
\begin{itemize}
\item symmetric monoidal closed categories (\cref{smclosed-closed-multicat}) and
\item the closed structure on the multicategory $\permcatsu$ (\cref{perm-closed-multicat}).
\end{itemize}
We emphasize that, just like multicategories, closed multicategories have symmetric group action compatible with the closed structure.

In the absence of a symmetric monoidal closed structure, a closed multicategory structure is the closest substitute that is still sufficient for a robust theory of self-enrichment (\cref{sec:selfenr-clmulti}), standard enrichment (\cref{sec:std-enr-multifunctor}), enriched diagrams, enriched Mackey functors (\cref{sec:enr-diag-psh}), and change of enrichment in those contexts (\cref{sec:fun-std-enr-multi,sec:enr-diag-change-enr}).  In short, closed multicategories are the focus of the rest of this work.

After defining closed multicategories (\cref{def:closed-multicat}), in \cref{rk:closed-multicat,rk:clmulti-zero-case} we discuss the relationship between our definition and those in the literature \cite{lambek,manzyuk,zakharevich}.  \cref{smclosed-closed-multicat} shows that each symmetric monoidal closed category yields a closed multicategory via the endomorphism construction.  

\begin{definition}\label{def:closed-multicat}
  A \index{closed!multicategory}\index{multicategory!closed}\emph{closed multicategory} is a triple
\[\big(\M \scs \clM \scs \ev\big)\]
consisting of the following data.
\begin{description}
\item[Underlying Multicategory] $\M = (\M,\ga,\opu)$ is a multicategory (\cref{def:enr-multicategory}).
\item[Internal Hom Objects]
For $n \geq 0$ and each $(n+1)$-tuple of objects $\angx = \ang{x_i}_{i=1}^n, y$ in $\M$, it is equipped with an object
\begin{equation}\label{clMangxy}
\clM\scmap{\angx;y} \in \M,
\end{equation}
which is called an \index{internal hom!multicategory}\index{hom!internal - multicategory}\emph{$n$-ary internal hom object} and also denoted $\clM_{\angx;\, y}$. 
\item[Symmetric Group Action on Internal Hom]
For objects $\angx, y \in \M$ as above and each permutation $\sigma \in \Sigma_n$, it is equipped with an invertible unary multimorphism
\begin{equation}\label{sigma-clMxy}
\begin{tikzcd}[column sep=large]
\clM\scmap{\angx;y} \ar{r}{\sigma}[swap]{\iso} & \clM\scmap{\angx\sigma;y} 
\end{tikzcd}\inspace \M
\end{equation}
with $\angx\sigma = \ang{x_{\sigma(i)}}_{i=1}^n$.  It is called the \index{symmetric group!action}\index{right action}\emph{right symmetric group action} or the \emph{right $\sigma$-action} on internal hom objects.
\item[Multicategorical Evaluation] For objects $\angx, y \in \M$ as above, it is equipped with an $(n+1)$-ary multimorphism
\begin{equation}\label{evangxy}
\begin{tikzcd}[column sep=large]
\left(\clM\scmap{\angx;y} \scs \angx\right) \ar{r}{\ev_{\angx;\, y}} & y 
\end{tikzcd}
\inspace \M,
\end{equation}
which is called the \index{evaluation!multicategorical}\index{multicategorical evaluation}\emph{multicategorical evaluation} or the \emph{evaluation} at $(\angx; y)$.
\end{description}
The above data are required to satisfy the axioms \cref{right-id-action,hom-equivariance,eval-bijection,eval-bij-eq} below.
\begin{description}
\item[Equivariance of Internal Hom]
For the identity permutation $\id_n \in \Sigma_n$, the right $\id_n$-action
\begin{equation}\label{right-id-action}
\clM\scmap{\angx;y} \fto{\id_n = \opu} \clM\scmap{\angx;y}
\end{equation}
is the colored unit of the internal hom object $\clM\scmap{\angx;y}$.  Moreover, for $\sigma,\tau \in \Sigma_n$, the following diagram of unary multimorphisms in $\M$ commutes.
\begin{equation}\label{hom-equivariance}
\begin{tikzcd}[column sep=large]
\clM\scmap{\angx;y} \ar{r}{\sigma} \ar{dr}[swap]{\sigma\tau} 
& \clM\scmap{\angx\sigma;y} \ar{d}{\tau}\\
& \clM\scmap{\angx\sigma\tau;y} 
\end{tikzcd}
\end{equation}
\item[Evaluation Bijection]
For objects $\angx = \ang{x_i}_{i=1}^n$, $\angy = \ang{y_j}_{j=1}^p$, $z \in \M$, the function
\begin{equation}\label{eval-bijection}
\begin{tikzcd}[/tikz/column 1/.append style={anchor=base east},
/tikz/column 2/.append style={anchor=base west}, column sep=huge, row sep=0ex]
\M\left(\angx \sscs \clM\scmap{\angy;z} \right) \ar{r}{\chi_{\angx;\, \angy;\, z}}[swap]{\iso} 
& \M\scmap{\angx, \angy; z}\\
f \ar[maps to]{r} & \ga\big(\ev_{\angy;\, z} \sscs f, \ang{1_{y_j}}_{j=1}^p\big)
\end{tikzcd}
\end{equation}
is a bijection, which is called the \index{bijection!evaluation}\index{evaluation bijection}\emph{evaluation bijection}.  Two multimorphisms in $\M$ that correspond under this bijection are called \index{partner}\emph{partners}.  We write $\pnf$ for the partner of a multimorphism $f$, so\label{not:partner} 
\[\chi(f) = \pn{f} \andspace \chi^\inv(g) = \pn{g}\]
for $f \in \M\left(\angx \sscs \clM\scmap{\angy;z} \right)$ and $g \in \M\scmap{\angx, \angy; z}$.
\item[Equivariance of Evaluation Bijection]
For objects $\angx, \angy, z \in \M$ as above and permutations $\sigma \in \Sigma_n$ and $\vsi \in \Sigma_p$, the following diagram of bijections commutes.
\begin{equation}\label{eval-bij-eq}
\begin{tikzpicture}[vcenter]
\def\h{5} \def\v{-1.3}
\draw[0cell=.9]
(0,0) node (a) {\M\left(\angx \sscs \clM\scmap{\angy;z} \right)}
(a)+(\h,0) node (b) {\M\scmap{\angx, \angy; z}}
(a)+(0,\v) node (c) {\M\left(\angx\sigma \sscs \clM\scmap{\angy;z} \right)}
(c)+(0,\v) node (d) {\M\left(\angx\sigma \sscs \clM\scmap{\angy\vsi;z} \right)}
(d)+(\h,0) node (e) {\M\scmap{\angx\sigma, \angy\vsi; z}}
;
\draw[1cell]
(a) edge node {\chi_{\angx;\, \angy;\, z}} (b)
(b) edge node {\sigma \times \vsi} (e)
(a) edge node[swap] {\sigma} (c)
(c) edge node[swap] {\ga(\vsi \sscs -)} (d)
(d) edge node {\chi_{\angx\sigma;\, \angy\vsi;\, z}} (e)
;
\end{tikzpicture}
\end{equation}
In \cref{eval-bij-eq} the arrows are defined as follows.
\begin{itemize}
\item The top left arrow $\sigma$ and the right vertical arrow $\sigma\times\vsi$ are right symmetric group action of $\M$ \cref{rightsigmaaction}.
\item In the lower left arrow, $\vsi$ is the right $\vsi$-action on the internal hom object $\clM\scmap{\angy;z}$ in \cref{sigma-clMxy}.
\item $\ga(\vsi \sscs -)$ is composition with the unary multimorphism $\vsi$ in $\M$. 
\item The two horizontal arrows $\chi$ are the evaluation bijections in \cref{eval-bijection}.
\end{itemize} 
\end{description}
This finishes the definition of a closed multicategory.

Moreover, a \index{multicategory!non-symmetric closed}\index{closed!multicategory!non-symmetric}\index{non-symmetric!closed multicategory}\emph{non-symmetric closed multicategory} is a triple $\big(\M,\clM,\ev\big)$ as above with the changes \cref{nscm-i,nscm-ii,nscm-iii} below.
\begin{romenumerate}
\item\label{nscm-i} $(\M,\ga,\opu)$ is a non-symmetric multicategory.
\item\label{nscm-ii} The internal hom objects $\clM\scmap{\angx;y}$ are not equipped with the right symmetric group action \cref{sigma-clMxy}.
\item\label{nscm-iii} We do not require the equivariance of
\begin{itemize}
\item internal hom objects, \cref{right-id-action,hom-equivariance}, and
\item evaluation bijection, \cref{eval-bij-eq}.
\end{itemize}   
Thus the only axiom is the evaluation bijection axiom \cref{eval-bijection}.
\end{romenumerate}
This finishes the definition of a non-symmetric closed multicategory.
\end{definition}

\begin{example}[Waldhausen Categories]\label{ex:waldhausen}
The 2-category of small Waldhausen categories, exact functors, and natural transformations extends to a closed multicategory $\Wald$ by \cite[5.6]{zakharevich}.  Since we do not use that result in this work, we refer the reader to \cite{zakharevich} for further discussion of the closed multicategory $\Wald$.  All the results in this work about (non-symmetric) closed multicategories apply to $\Wald$.  See, for example, \cref{cl-multi-cl-cat}.
\end{example}

We discuss more examples below after some explanation and remarks.

\begin{explanation}[Closed Multicategories]\label{expl:closed-multi}
Suppose $(\M,\clM,\ev)$ is a non-symmetric closed multicategory.
\begin{enumerate}
\item\label{expl:clmulti-i}
If $\M$ is a closed multicategory, then the right $\sigma$-action on internal hom object \cref{sigma-clMxy} is an element
\[\sigma \in \M\Big(\clM\scmap{\angx;y} \sscs \clM\scmap{\angx\sigma;y}\Big).\]
It is a unary multimorphism in $\M$ regardless of the length of $\angx = \ang{x_i}_{i=1}^n$.
\item\label{expl:clmulti-ii}
The evaluation at $(\angx;y)$ in \cref{evangxy}, which is an $(n+1)$-ary multimorphism, is an element
\[\ev_{\angx;\, y} \in \M\Big(\clM\scmap{\angx;y} \scs \angx \sscs y \Big).\]
This is an analog of the evaluation \cref{evaluation} in a symmetric monoidal closed category.
\item\label{expl:clmulti-iii}
If $\angx = \ang{}$, then evaluation at $(\ang{}; y)$ is a unary multimorphism
\[\clM\scmap{\ang{};y} \fto{\ev_{\ang{};\, y}} y \inspace \M.\]
We do \emph{not} require this to be the colored unit of $y$.  We elaborate on this point in \cref{rk:clmulti-zero-case} below.
\item\label{expl:clmulti-iv}
In the evaluation bijection $\chi_{\angx;\, \angy;\, z}$ in \cref{eval-bijection} with $\angx = \ang{x_i}_{i=1}^n$ and $\angy = \ang{y_j}_{j=1}^p$,
\begin{itemize}
\item the domain $\M\left(\angx \sscs \clM\scmap{\angy;z} \right)$ is an $n$-ary multimorphism set in $\M$, and
\item the codomain $\M\scmap{\angx, \angy; z}$ is an $(n+p)$-ary multimorphism set in $\M$.
\end{itemize}
The evaluation bijection is an analog of the $\otimes$-$\Hom$ adjunction in a symmetric monoidal closed category (\cref{def:closedcat}).  Partners---which mean multimorphisms that correspond to each other under the evaluation bijection---are analogs of adjoints.
\item\label{expl:clmulti-v}
If $\angx = \ang{}$, then the evaluation bijection $\chi_{\ang{};\, \angy;\, z}$ is the following bijection.
\begin{equation}\label{eval-bij-zero}
\begin{tikzcd}[/tikz/column 1/.append style={anchor=base east},
/tikz/column 2/.append style={anchor=base west}, column sep=huge, row sep=0ex]
\M\left(\ang{} \sscs \clM\scmap{\angy;z} \right) \ar{r}{\chi_{\ang{};\, \angy;\, z}}[swap]{\iso} 
& \M\scmap{\angy; z}\\
f \ar[maps to]{r} & \ga\big(\ev_{\angy;\, z} \sscs f, \ang{1_{y_j}}_{j=1}^p\big)
\end{tikzcd}
\end{equation}
Thus, via this evaluation bijection, each multimorphism set $\M\scmap{\angy;z}$ is in bijection with the nullary multimorphism set with output given by the internal hom object $\clM\scmap{\angy;z}$.\defmark
\end{enumerate}
\end{explanation}

\begin{remark}[Variants]\label{rk:closed-multicat}
There are other variants of closed multicategories in the literature.  We briefly discuss some of them here.
\begin{enumerate}
\item A closed multicategory as in \cref{def:closed-multicat} is called a \index{closed!symmetric multicategory}\emph{closed symmetric multicategory} in \cite[1.2]{zakharevich}.  There are no other differences between the two definitions.  In particular, the main result in \cite{zakharevich}---that small Waldhausen categories form a closed multicategory $\Wald$ (\cref{ex:waldhausen})---still holds with our \cref{def:closed-multicat}.  The permutative analog of this observation is \cref{perm-closed-multicat}.
\item A \index{biclosed monoidal multicategory}\emph{biclosed monoidal multicategory} in \cite[page 106]{lambek} is analogous to a non-symmetric closed multicategory as in \cref{def:closed-multicat}.  However, Lambek's definition has more structure, denoted $\texttt{i}$ and $\texttt{m}$ there, which roughly correspond to a monoidal unit and a monoidal product. 
\item A \index{multicategory!closed}\index{closed!multicategory}\emph{closed multicategory} in \cite[3.6, 3.7]{manzyuk} is a more restrictive version of a non-symmetric closed multicategory as in \cref{def:closed-multicat}.  We discuss this nontrivial difference in more detail in \cref{rk:clmulti-zero-case}.
\end{enumerate}
Our terminology is in line with \cref{def:enr-multicategory}, where a $\V$-multicategory is equipped with a symmetric group action.  We add the adjective \emph{non-symmetric} to the variant without the symmetric group action. 
\end{remark}

\begin{remark}[Important Differences with Manzyuk's Definition]\label{rk:clmulti-zero-case}
A closed multicategory as in \cite[3.7]{manzyuk} is more restrictive than a non-symmetric closed multicategory as in \cref{def:closed-multicat}.  Specifically, the definition in \cite[3.7]{manzyuk} requires (i) the object equality
\begin{equation}\label{manzyuk-def-i}
\clM\scmap{\ang{};y} = y \forspace y \in \Ob\M
\end{equation}
and (ii) the unary multimorphism equality
\begin{equation}\label{manzyuk-def-ii}
\big( \clM\scmap{\ang{};y} \scs \ang{} \big) \fto{\ev_{\ang{};\, y} = \opu_y} y \inspace \M\scmap{y;y}.
\end{equation}

On the other hand, \cref{def:closed-multicat} does \emph{not} require these two equalities.  The reason that we do not impose the equalities \cref{manzyuk-def-i,manzyuk-def-ii} is that they are \emph{not} satisfied by our most basic examples of endomorphism multicategories in \cref{smclosed-closed-multicat} below.  As stated in \cref{EndV-closed-hom}, the nullary internal hom object is the internal hom object
\[\clEndV\scmap{\ang{}; y} = [\tu, y] \forspace y \in \V\]
because an empty $\otimes$ is, by definition, the monoidal unit $\tu$ in $\V$.  The object $[\tu,y]$ is, in general, not equal to $y$.  The nullary multicategorical evaluation $\ev_{\ang{};\,y}$ in \cref{EndV-ev-empty} is also not an identity in general.  Due to this nontrivial difference between \cref{def:closed-multicat} and the one in \cite{manzyuk}, in this work we will not use any of the results in \cite{manzyuk}.
\end{remark}

\cref{smclosed-closed-multicat} below says that the endomorphism multicategory (\cref{ex:endc}) of a symmetric monoidal closed category (\cref{def:closedcat}) is a closed multicategory, as one would expect.  It is briefly mentioned in \cite[Section 1]{zakharevich}; we provide a proof here for completeness.  We remind the reader that an iterated monoidal product $\txotimes_{i=1}^n$ is left normalized (\cref{expl:leftbracketing}), and an empty $\otimes$ is the monoidal unit $\tu$.

\begin{proposition}\label{smclosed-closed-multicat}
For each symmetric monoidal closed category $\big(\V,\otimes,\tu,\xi,[,]\big)$, the endomorphism multicategory $\End\,\V$ becomes a closed multicategory 
\[\big(\EndV \scs \clEndV \scs \ev\big)\]
when it is equipped with the following data.
\begin{itemize}
\item For objects $\ang{x_i}_{i=1}^n, y \in \V$, the $n$-ary internal hom object \cref{clMangxy} is defined as the object
\begin{equation}\label{EndV-closed-hom}
\clEndV\scmap{\ang{x_i}_{i=1}^n;y} = \big[\txotimes_{i=1}^n x_i \scs y \big] \inspace \V.
\end{equation}
\item For a permutation $\sigma \in \Sigma_n$, the right $\sigma$-action \cref{sigma-clMxy} on internal hom objects
\begin{equation}\label{EndV-closed-sym}
\big[\txotimes_{i=1}^n x_i \scs y \big] \fto[\iso]{\sigma} \big[\txotimes_{i=1}^n x_{\sigma(i)} \scs y \big] \inspace \V
\end{equation}
is the image under $[-,y]$ of the unique symmetry coherence isomorphism
\begin{equation}\label{permute-factors-sigma}
\txotimes_{i=1}^n x_{\sigma(i)} \fto[\iso]{\sigma} \txotimes_{i=1}^n x_i
\end{equation}
that permutes the $n$ factors according to $\sigma$.
\item The multicategorical evaluation $\ev_{\angx;y}$ \cref{evangxy} is defined as the following composite in $\V$ if $n > 0$.
\begin{equation}\label{EndV-ev}
\begin{tikzpicture}[vcenter]
\draw[0cell=.9]
(0,0) node (a) {\big[\txotimes_{i=1}^n x_i \scs y\big] \otimes x_1 \otimes \cdots \otimes x_n}
(a)+(0,-1.3) node (b) {\big[\txotimes_{i=1}^n x_i \scs y\big] \otimes \big(\txotimes_{i=1}^n x_i\big)}
(a)+(4,0) node (c) {y}
;
\draw[1cell=.9]
(a) edge node {\ev_{\angx;\, y}} (c)
(a) edge node[swap] {\alpha} node {\iso} (b)
(b) edge[bend right=15] node[swap,pos=.6] {\ev_{\otimes_{i=1}^n x_i,\, y}} (c)
;
\end{tikzpicture}
\end{equation}
In \cref{EndV-ev} $\alpha$ is the unique coherence isomorphism in $\V$ that moves parentheses, which is the identity if $n=1$, and $\ev_{\otimes_{i=1}^n x_i,\, y}$ is the evaluation in $\V$ \cref{evaluation}.  If $n=0$, then the multicategorical evaluation $\ev_{\ang{};\,y}$ is defined as the following composite, with $\rho$ the right unit isomorphism and $\ev_{\tu,y}$ the evaluation in $\V$.
\begin{equation}\label{EndV-ev-empty}
\begin{tikzpicture}[baseline={(a.base)}]
\def\h{2.5} \def\w{.6}
\draw[0cell=1]
(0,0) node (a) {[\tu,y]}
(a)+(2.5,0) node (b) {[\tu,y] \otimes \tu}
(b)+(2.2,0) node (c) {y}
;
\draw[1cell=.9]
(a) edge node {\rho^\inv} node[swap] {\iso} (b)
(b) edge node {\ev_{\tu,y}} (c)
;
\draw[1cell=.9]
(a) [rounded corners=3pt] |- ($(b)+(-1,\w)$)
-- node {\ev_{\ang{};\,y}} ($(b)+(1,\w)$) -| (c)
;
\end{tikzpicture}
\end{equation}
\end{itemize}
\end{proposition}

\begin{proof}
We check the axioms \cref{right-id-action,hom-equivariance,eval-bijection,eval-bij-eq} for $\EndV$.  For the rest of this proof, $\angx = \ang{x_i}_{i=1}^n$, $\angy = \ang{y_j}_{j=1}^p$, $y$, and $z$ are objects in $\V$.

\medskip
\emph{Equivariance of Internal Hom}.  For the identity permutation $\id_n \in \Sigma_n$, the symmetry coherence isomorphism \cref{permute-factors-sigma} is the identity morphism.  Thus the right $\id_n$-action in \cref{EndV-closed-sym} is the identity morphism, proving the axiom \cref{right-id-action}.

For permutations $\sigma, \tau \in \Sigma_n$, the axiom \cref{hom-equivariance} requires the  commutativity of the following diagram in $\V$.
\begin{equation}\label{EndV-hom-eq}
\begin{tikzpicture}[baseline={(a.base)}]
\def\h{3} \def\w{.6}
\draw[0cell=1]
(0,0) node (a) {\big[\txotimes_{i=1}^n x_i \scs y\big]}
(a)+(\h,0) node (b) {\big[\txotimes_{i=1}^n x_{\sigma(i)} \scs y\big]}
(b)+(\h+.4,0) node (c) {\big[\txotimes_{i=1}^n x_{\sigma\tau(i)} \scs y\big]}
;
\draw[1cell=.9]
(a) edge node {\sigma} (b)
(b) edge node {\tau} (c)
;
\draw[1cell=.9]
(a) [rounded corners=3pt, shorten >=-.5ex, shorten <=-.5ex] |- ($(b)+(-1,\w)$)
-- node {\sigma\tau} ($(b)+(1,\w)$) -| (c)
;
\end{tikzpicture}
\end{equation}
The diagram \cref{EndV-hom-eq} is the image under $[-,y]$ of the following diagram in $\V$.
\[\begin{tikzpicture}[baseline={(a.base)}]
\def\h{2.5} \def\w{.6}
\draw[0cell=1]
(0,0) node (a) {\txotimes_{i=1}^n x_{\sigma\tau(i)}}
(a)+(\h+.3,0) node (b) {\txotimes_{i=1}^n x_{\sigma(i)}}
(b)+(\h,0) node (c) {\txotimes_{i=1}^n x_i}
;
\draw[1cell=.9]
(a) edge node {\tau} (b)
(b) edge node {\sigma} (c)
;
\draw[1cell=.9]
(a) [rounded corners=3pt] |- ($(b)+(-1,\w)$)
-- node {\sigma\tau} ($(b)+(1,\w)$) -| (c)
;
\end{tikzpicture}\]
This diagram commutes by the uniqueness of the symmetry coherence isomorphism that permutes the $n$ factors according to $\sigma\tau$ \cite[XI.1 Theorem 1]{maclane}.  Thus the diagram \cref{EndV-hom-eq} is also commutative.

\medskip
\emph{Evaluation Bijection}.  Note that each of the two cases of $\ev_{\angx;y}$ in \cref{EndV-ev,EndV-ev-empty} consists of a coherence isomorphism followed by an instance of the evaluation in $\V$.  For the rest of this proof, the symbol $\iso$ denotes a coherence isomorphism.

The function $\chi = \chi_{\angx;\, \angy;\, z}$ in \cref{eval-bijection}---which we want to show is a bijection---sends a morphism
\begin{equation}\label{f-tensor-to-hom}
\txotimes_{i=1}^n x_i \fto{f} \big[\txotimes_{j=1}^p y_j \scs z\big] \inspace \V
\end{equation}
to the following composite.
\begin{equation}\label{chi-f}

\end{equation}
The desired inverse of $\chi$ sends a morphism
\[\big(\txotimes_{i=1}^n x_i\big) \otimes y_1 \otimes \cdots \otimes y_p \fto{g} z \inspace \V\]
to the adjoint
\[\txotimes_{i=1}^n x_i \fto{\gbar^\#} \big[\txotimes_{j=1}^p y_j \scs z\big]\]
of the following composite $\gbar$.
\[%
\]
The assignments, $f \mapsto \chi(f)$ and $g \mapsto \gbar^\#$, are inverses of each other by the uniqueness of adjoints and the fact that the evaluation in $\V$ \cref{evaluation} is the counit of the $\otimes$-$[,]$ adjunction.  Thus $\chi$ is a bijection.

\medskip
\emph{Equivariance of Evaluation Bijection}. 
For permutations $\sigma \in \Sigma_n$ and $\vsi \in \Sigma_p$ and a morphism $f$ as in \cref{f-tensor-to-hom}, the desired commutative diagram \cref{eval-bij-eq} is the boundary of the following diagram in $\V$.
\[%
\]
The following statements hold for the diagram above.
\begin{itemize}
\item The top triangle and the quadrilateral under it commute by the functoriality of $\otimes$.
\item The middle quadrilateral commutes by the coherence theorem for symmetric monoidal categories \cite[XI.1]{maclane}.
\item The lower left quadrilateral commutes by the naturality of coherence isomorphisms.
\item The lower right quadrilateral commutes because each of the two composites has adjoint
\[\big[\txotimes_{j=1}^p y_j \scs z\big] \to \big[\txotimes_{j=1}^p y_{\vsi(j)} \scs z\big]\]
given by the image under $[-,z]$ of the symmetry coherence isomorphism
\[\txotimes_{j=1}^p y_{\vsi(j)} \fto[\iso]{\vsi} \txotimes_{j=1}^p y_j.\]
\end{itemize}
This finishes the proof of the axioms \cref{right-id-action,hom-equivariance,eval-bijection,eval-bij-eq} for $\EndV$.
\end{proof}

\begin{example}\label{ex:smc-clmult}
\cref{smclosed-closed-multicat} applies to each of the symmetric monoidal closed categories listed in \cref{ex:EndV-enriched}.  In other words, via the endomorphism multicategory construction, each of those symmetric monoidal closed categories is a closed multicategory.  However, \cref{smclosed-closed-multicat} does not apply to $\permcatsu$ because it is not a monoidal category.
\end{example}

\section{Internal Hom Permutative Categories}
\label{sec:internal-hom-permcat}

Recall the following two results about permutative categories.
\begin{enumerate}
\item By \cref{thm:permcatmulticat} $\permcatsu$ is a $\Cat$-multicategory.  As a multicategory, it has 
\begin{itemize}
\item small permutative categories (\cref{def:symmoncat}) as objects and
\item multilinear functors (\cref{def:nlinearfunctor}) as multimorphisms.
\end{itemize}
\item By \cref{permcat-selfenr} the category $\permcatsu$ is a $\permcatsu$-category.
\end{enumerate}
In \cref{sec:internal-hom-permcat,sec:multiev-permcat,sec:permcatsu-clmulti}, we generalize these two facts by showing that $\permcatsu$ is a closed multicategory (\cref{def:closed-multicat}); see \cref{perm-closed-multicat}.  In this section we construct the internal hom objects and their symmetry group action in the closed multicategory structure on $\permcatsu$.  Here is an outline of this section.
\begin{itemize}
\item The internal hom objects are constructed in \cref{def:clp-angcd} and verified in \cref{clp-angcd-permutative}.
\item The symmetric group action on the internal hom objects are constructed in \cref{def:clp-ancd-sigma} and verified in \cref{clp-angcd-strict}.
\end{itemize}
We discuss the multicategorical evaluation and the closed multicategory axioms of $\permcatsu$ in subsequent sections.

\subsection*{Internal Hom Objects}

To construct the closed multicategory structure on $\permcatsu$, first we define the $n$-ary internal hom objects \cref{clMangxy}.  For a generic permutative category, we denote the monoidal product, monoidal unit, and braiding by $\oplus$, $\pu$, and $\xi$, respectively.

\begin{definition}\label{def:clp-angcd}
For small permutative categories $\D$ and $\angC = \ang{\C_i}_{i=1}^n$ for $n \geq 0$, we define the data of a small permutative category
\begin{equation}\label{clp-internal-hom}
\Big(\clp\scmap{\angC;\D} \scs \oplus \scs \clpu \scs \clxi \Big),
\end{equation}
which is called an \index{permutative category!internal hom -}\index{internal hom!permutative category}\emph{internal hom permutative category}, as follows.  We also use the shortened notation 
\[\psu = \permcatsu \andspace \clpsu = \clp.\]
\begin{description}
\item[Underlying Category]
It is the category $\psu\scmap{\angC;\D}$ in \cref{definition:permcatsus-homcat}.
\begin{itemize}
\item Its objects are $n$-linear functors $\angC \to \D$ (\cref{def:nlinearfunctor}).
\item Its morphisms are $n$-linear transformations (\cref{def:nlineartransformation}).
\end{itemize}
\item[Monoidal Product on Objects]
For two $n$-linear functors
\begin{equation}\label{PQnlinear}
\big(P, \{P^2_i\}_{i=1}^n\big) \scs \big(Q, \{Q^2_i\}_{i=1}^n\big) 
\cn \txprod_{i=1}^n \C_i \to \D,
\end{equation}
their monoidal product has underlying functor $P \oplus Q$ given by the composite functor below.
\begin{equation}\label{PplusQfunctor}
\begin{tikzpicture}[baseline={(a.base)}]
\def\h{2.5} \def\w{.6}
\draw[0cell=1]
(0,0) node (a) {\txprod_{i=1}^n \C_i}
(a)+(\h,0) node (b) {\D \times \D}
(b)+(\h-.5,0) node (c) {\D}
;
\draw[1cell=.9]
(a) edge node {(P,Q)} (b)
(b) edge node {\oplus} (c)
;
\draw[1cell=.9]
(a) [rounded corners=3pt] |- ($(b)+(-1,\w)$)
-- node {P \oplus Q} ($(b)+(1,\w)$) -| (c)
;
\end{tikzpicture}
\end{equation}
In other words, the functor $P \oplus Q$ is given by the objectwise monoidal product
\[(P \oplus Q)\angx = P\angx \oplus Q\angx\]
for objects and morphisms $\angx$ in $\txprod_{i=1}^n \C_i$.

For each $i \in \{1,\ldots,n\}$, the $i$-th linearity constraint of $P \oplus Q$, denoted $(P \oplus Q)^2_i$, is defined as follows.  For objects $\angx \in \txprod_{i=1}^n \C_i$ and $x_i' \in \C_i$, using \cref{notation:compk} we denote by
\begin{equation}\label{angxp-angxpp}
\angxp = \angx \compi x_i' \andspace \angxpp = \angx \compi (x_i \oplus x_i').
\end{equation}
The corresponding component of $(P \oplus Q)^2_i$ is the following composite in $\D$, with $\xi$ denoting the braiding in $\D$. 
\begin{equation}\label{PplusQtwoi}

\end{equation}
The naturality of $(P \oplus Q)^2_i$ follows from the naturality of $P^2_i$, $Q^2_i$, and $\xi$, together with the functoriality of $\oplus$ in $\D$.
\item[Monoidal Product on Morphisms]
For $n$-linear functors $P$ and $Q$ as above, suppose given $n$-linear transformations $\theta$ and $\psi$ as follows. 
\begin{equation}\label{theta-psi-tr}
%
\end{equation}
Their monoidal product 
\[\theta \oplus \psi \cn P \oplus Q \to P' \oplus Q'\]
is the natural transformation defined by the whiskering below.
\begin{equation}\label{theta-plus-psi}
%
\end{equation}
In other words, for each object $\angx \in \txprod_{i=1}^n \C_i$, $\theta \oplus \psi$ has $\angx$-component given by the monoidal product
\[(\theta \oplus \psi)_{\angx} = \theta_{\angx} \oplus \psi_{\angx} \cn P\angx \oplus Q\angx \to P'\angx \oplus Q'\angx.\]
This defines an $n$-linear transformation $\theta \oplus \psi$ for the following reasons.
\begin{itemize}
\item The naturality of $\theta \oplus \psi$ follows from the naturality of $\theta$ and $\psi$, together with the functoriality of $\oplus$ in $\D$.  
\item The unity axiom \cref{niitransformationunity} holds for $\theta \oplus \psi$ because, if any $x_i = \pu$ in $\C_i$, then $\theta_{\angx} = 1_{\pu} = \psi_{\angx}$.
\item The constraint compatibility axiom \cref{eq:monoidal-in-each-variable} holds for $\theta \oplus \psi$ by
\begin{itemize}
\item the naturality of the braiding $\xi$ in $\D$ and
\item the axiom \cref{eq:monoidal-in-each-variable} for $\theta$ and $\psi$.
\end{itemize} 
\end{itemize}
The construction $\oplus$ for $\clpsu\scmap{\angC;\D}$ preserves identities and composition of $n$-linear transformations by the functoriality of $\oplus$ in $\D$.  Moreover, $\oplus$ is associative on $n$-linear transformations because $\oplus$ in $\D$ is associative.
\item[Monoidal Unit]
The monoidal unit is the constant functor 
\begin{equation}\label{clpu}
\clpu \cn \txprod_{i=1}^n \C_i \to \D
\end{equation}
at the monoidal unit $\pu$ in $\D$.  Each of its $n$ linearity constraints is given componentwise by the identity morphism
\[1_\pu \cn \pu \oplus \pu = \pu \to \pu \inspace \D.\]
The axioms of an $n$-linear functor hold for $\clpu$ because
\begin{itemize}
\item $\pu$ is a strict monoidal unit in $\D$ and
\item there are morphism equalities 
\begin{equation}\label{braiding-identity}
\xi_{\pu,?} = 1_{?} = \xi_{?,\pu} \inspace \D.
\end{equation}
\end{itemize} 
These two facts also imply that $\clpu$ is a strict two-sided unit for $\oplus$.
\item[Braiding] 
For $n$-linear functors $P$ and $Q$ as above, the $(P,Q)$-component of the braiding\label{not:clxi}
\[\clxi_{P,Q} \cn P \oplus Q \to Q \oplus P\]
is the natural isomorphism given by the following pasting diagram, with $\tau$ swapping the two factors and $\xi$ denoting the braiding in $\D$.
\begin{equation}\label{clp-braiding}
\begin{tikzpicture}[baseline={(a.base)}]
\def\h{2.5} \def\v{.7}
\draw[0cell]
(0,0) node (a) {\txprod_{i=1}^n \C_i}
(a)+(\h,\v) node (b1) {\D \times \D}
(a)+(\h,-\v) node (b2) {\D \times \D}
(b1)+(\h,-\v) node (c) {\D}
;
\draw[1cell=.9]
(a) edge node[pos=.7] {(P,Q)} (b1)
(a) edge node[swap,pos=.7] {(Q,P)} (b2)
(b1) edge node[swap] (s) {\tau} (b2)
(b1) edge node[pos=.4] {\oplus} (c)
(b2) edge node[swap,pos=.4] {\oplus} (c)
;
\draw[2cell]
node[between=s and c at .4, rotate=-90, 2label={below,\xi}] {\Rightarrow}
;
\end{tikzpicture}
\end{equation}
In other words, for each object $\angx \in \txprod_{i=1}^n \C_i$, $\clxi_{P,Q}$ has $\angx$-component given by the braiding in $\D$
\[(\clxi_{P,Q})_{\angx} = \xi_{P\angx,Q\angx} \cn P\angx \oplus Q\angx \fto{\iso} Q\angx \oplus P\angx.\]
The naturality of the braiding $\xi$ in $\D$ implies the naturality of each of
\begin{itemize}
\item $(\clxi_{P,Q})_{\angx}$ with respect to $\angx$ and
\item $\clxi_{P,Q}$ with respect to $P$ and $Q$.
\end{itemize} 
\end{description}
This finishes the definition of $\clp\scmap{\angC;\D}$.  \cref{clp-angcd-permutative} proves that it is a permutative category.  
\end{definition}

\begin{explanation}\label{expl:clp-zero}
Consider \cref{def:clp-angcd}.
\begin{itemize}
\item If $n=0$, then $\angC$ is the empty sequence and
\[\clp\scmap{\ang{};\D} = \D\]
as permutative categories.
\item If $n=1$, then  
\[\clp(\C;\D) = \permcatsu(\C,\D),\]
the hom permutative category in \cref{psucd-hom-permcat}.\defmark
\end{itemize}
\end{explanation}

\begin{lemma}\label{clp-angcd-permutative}
For small permutative categories $\angC = \ang{\C_i}_{i=1}^n$ and $\D$, the quadruple
\[\Big(\clp\scmap{\angC;\D} \scs \oplus \scs \clpu \scs \clxi \Big)\] 
in \cref{def:clp-angcd} is a small permutative category.
\end{lemma}

\begin{proof}
By \cref{expl:clp-zero} and \cref{psucd-hom-permcat} for the case $n=1$, we only need to check the cases for $n>1$.  We already explained some of the required conditions in \cref{def:clp-angcd}.  It remains to check statements \cref{clp-angcd-perm-i,clp-angcd-perm-ii,clp-angcd-perm-iii} below.
\begin{romenumerate}
\item\label{clp-angcd-perm-i} The data defined in \cref{PplusQfunctor,PplusQtwoi}
\[\big(P \oplus Q \scs \big\{(P \oplus Q)^2_i\big\}_{i=1}^n\big) \cn \txprod_{i=1}^n \C_i \to \D\]
satisfy the axioms \cref{nlinearunity,constraintunity,eq:ml-f2-assoc,eq:ml-f2-symm,eq:f2-2by2} of an $n$-linear functor.
\item\label{clp-angcd-perm-ii} The construction $\oplus$ is associative on $n$-linear functors $\txprod_{i=1}^n \C_i \to \D$.
\item\label{clp-angcd-perm-iii} For $n$-linear functors $P$ and $Q$, the natural isomorphism defined in \cref{clp-braiding}
\[\clxi_{P,Q} \cn P \oplus Q \to Q \oplus P\] 
satisfies the axioms \cref{niitransformationunity} and \cref{eq:monoidal-in-each-variable} of an $n$-linear transformation.
\end{romenumerate}
Once we establish statements \cref{clp-angcd-perm-i,clp-angcd-perm-ii,clp-angcd-perm-iii} above, the symmetry and hexagon axioms \cref{symmoncatsymhexagon} for $\clpsu\scmap{\angC;\D}$ follow from those for $\D$.

\medskip
\emph{Statement \cref{clp-angcd-perm-i}}.  The unity axiom \cref{nlinearunity} for $P \oplus Q$ follows from the unity axiom for $P$ and $Q$, together with the strict unity of $\pu$ in $\D$.  The constraint unity, associativity, and symmetry axioms, \cref{constraintunity,eq:ml-f2-assoc,eq:ml-f2-symm}, are proved using the argument for statement \cref{psucd-hom-i} in the proof of \cref{psucd-hom-permcat}, with an appropriate change of notation.

For the constraint 2-by-2 axiom \cref{eq:f2-2by2}, suppose $i \neq k \in \{1,\ldots,n\}$ and consider objects
\[\angx = \ang{x_j^0}_{j=1}^n \in \txprod_{j=1}^n \C_j \scs \quad 
x_i^1 \in \C_i \scs \andspace x_k^1 \in \C_k.\]
We define the following objects for $\ell \in \{i,k\}$, $a,b \in \{0,1,2\}$, and $R \in \{P,Q\}$.
\[x_\ell^2 = x_\ell^0 \oplus x_\ell^1 \in \C_\ell \qquad 
R_{a,b} = R\bang{x \compi x_i^a \compk x_k^b} \in \D\]
For example, we have the objects
\[P_{0,0} = P\angx \andspace Q_{1,2} = Q\bang{x \compi x_i^1 \compk (x_k^0 \oplus x_k^1)}.\]
In the following diagram, we omit all the $\oplus$ symbols to save space.  With these conventions, the constraint 2-by-2 diagram \cref{eq:f2-2by2} for $P \oplus Q$ is the boundary of the following diagram in $\D$.
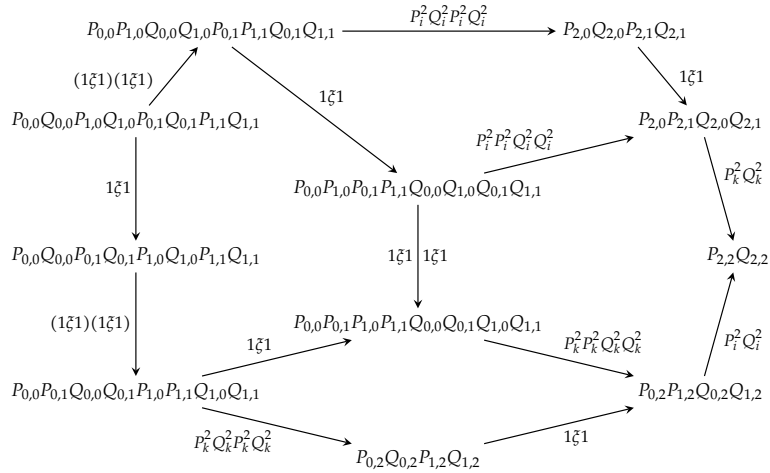
\[\begin{tikzpicture}
\def\g{1} \def\h{7.5} \def\u{1.2} \def\v{1.8} \def\t{0}
\draw[0cell=.7]
(0,0) node (a1) {P_{0,0} Q_{0,0} P_{1,0} Q_{1,0} P_{0,1} Q_{0,1} P_{1,1} Q_{1,1}}
(a1)+(\g,\u) node (a2) {P_{0,0} P_{1,0} Q_{0,0} Q_{1,0} P_{0,1} P_{1,1} Q_{0,1} Q_{1,1}}
(a1)+(\h,0) node (a4) {P_{2,0} P_{2,1} Q_{2,0} Q_{2,1}}
(a4)+(-\g,\u) node (a3) {P_{2,0} Q_{2,0} P_{2,1} Q_{2,1}}
(a4)+(\g/2,-\v) node (a5) {P_{2,2} Q_{2,2}}
(a1)+(0,-\v) node (b1) {P_{0,0} Q_{0,0} P_{0,1} Q_{0,1} P_{1,0} Q_{1,0} P_{1,1} Q_{1,1}}
(b1)+(0,-\v) node (b2) {P_{0,0} P_{0,1} Q_{0,0} Q_{0,1} P_{1,0} P_{1,1} Q_{1,0} Q_{1,1}}
(b2)+(\h/2,-\v/2) node (b3) {P_{0,2} Q_{0,2} P_{1,2} Q_{1,2}} 
(b2)+(\h,0) node (b4) {P_{0,2} P_{1,2} Q_{0,2} Q_{1,2}}
node[between=b1 and a4 at .5] (c1) {P_{0,0} P_{1,0} P_{0,1} P_{1,1} Q_{0,0} Q_{1,0} Q_{0,1} Q_{1,1}}
(c1)+(0,-\v) node (c2) {P_{0,0} P_{0,1} P_{1,0} P_{1,1} Q_{0,0} Q_{0,1} Q_{1,0} Q_{1,1}}
;
\draw[1cell=.7]
(a1) edge node[pos=.2] {(1\xi 1)(1\xi 1)} (a2)
(a2) edge node {P^2_i Q^2_i P^2_i Q^2_i} (a3)
(a3) edge node[pos=.75] {1\xi 1} (a4)
(a4) edge node {P^2_k Q^2_k} (a5)
(a1) edge node[swap] {1\xi 1} (b1)
(b1) edge node[swap] {(1\xi 1)(1\xi 1)} (b2)
(b2) edge[bend right=\t] node[swap] {P^2_k Q^2_k P^2_k Q^2_k} (b3)
(b3) edge[bend right=\t] node[swap] {1\xi 1} (b4)
(b4) edge node[swap] {P^2_i Q^2_i} (a5)
(a2) edge node {1\xi 1} (c1)
(c1) edge node {P^2_i P^2_i Q^2_i Q^2_i} (a4)
(c1) edge node {1\xi 1} node[swap] {1\xi 1} (c2)
(b2) edge[bend left=\t] node {1\xi 1} (c2)
(c2) edge[bend left=\t] node {P^2_k P^2_k Q^2_k Q^2_k} (b4)
;
\end{tikzpicture}\]
The following statements hold for the diagram above.
\begin{itemize}
\item The top and bottom quadrilaterals commute by the naturality of the braiding $\xi$ in $\D$.
\item The left sub-region commutes by the coherence theorem for symmetric monoidal categories \cite[XI.1 Theorem 1]{maclane}.
\item The right pentagon commutes by the constraint 2-by-2 axiom for $P$ and $Q$.
\end{itemize}
This proves that $P \oplus Q$ is an $n$-linear functor.

\medskip
\emph{Statement \cref{clp-angcd-perm-ii}}.  The associativity of $\oplus$ for $n$-linear functors
\[\big(P,\{P^2_i\}_{i=1}^n\big) \scs \big(Q,\{Q^2_i\}_{i=1}^n\big) 
\scs \big(R,\{R^2_i\}_{i=1}^n\big) \cn \txprod_{i=1}^n \C_i \to \D\]
is proved using the argument for statement \cref{psucd-hom-ii} in the proof of \cref{psucd-hom-permcat}, with an appropriate change of notation.

\medskip
\emph{Statement \cref{clp-angcd-perm-iii}}.  The axioms \cref{niitransformationunity} and \cref{eq:monoidal-in-each-variable} of an $n$-linear transformation for
\[\clxi_{P,Q} \cn P \oplus Q \to Q \oplus P\] 
are proved using the argument for statement \cref{psucd-hom-iv} in the proof of \cref{psucd-hom-permcat}, with an appropriate change of notation.
\end{proof}

From now on, $\clpsu\scmap{\angC;\D}$ is a permutative category as in \cref{clp-angcd-permutative}.

\subsection*{Symmetric Group Action on Internal Hom}

Next we define the right symmetric group action \cref{sigma-clMxy} on the internal hom permutative categories $\clpsu\scmap{\angC;\D}$.  Recall from \cref{def:clp-angcd} that the underlying category of $\clpsu\scmap{\angC;\D}$ is the category $\psu\scmap{\angC;\D}$ in \cref{definition:permcatsus-homcat}, with $n$-linear functors as objects and $n$-linear transformations as morphisms.

\begin{definition}\label{def:clp-ancd-sigma}
For small permutative categories $\D$ and $\angC = \ang{\C_i}_{i=1}^n$ for $n \geq 0$ and a permutation $\sigma \in \Sigma_n$, we define the functor
\begin{equation}\label{clp-sigma-action}
\clp\scmap{\angC;\D} \fto[\iso]{\sigma} \clp\scmap{\angC\sigma;\D}
\end{equation}
as the isomorphism
\[\permcatsu\scmap{\angC;\D} \fto[\iso]{\sigma} \permcatsu\scmap{\angC\sigma;\D}\]
in \cref{permiicatsymgroupaction}.  
\end{definition}

\begin{lemma}\label{clp-angcd-strict}
In the context of \cref{def:clp-ancd-sigma}, the following statements hold.
\begin{enumerate}
\item\label{clp-angcd-sigma-i} Equipped with identity monoidal and unit constraints, $\sigma$ in \cref{clp-sigma-action} is a strict symmetric monoidal isomorphism.  
\item\label{clp-angcd-sigma-ii}
The equivariance axioms for internal hom objects, \cref{right-id-action,hom-equivariance}, hold.
\end{enumerate}
\end{lemma}

\begin{proof}
Once statement \cref{clp-angcd-sigma-i} is proved, statement \cref{clp-angcd-sigma-ii} follows from (i) the definition \cref{permiicatsigmaaction} and (ii) the associativity of functor composition and horizontal composition of natural transformations.

For statement \cref{clp-angcd-sigma-i}, we need to show that $\sigma$ preserves (i) the monoidal unit and (ii) the monoidal product of its domain and codomain.  We denote $\sigma(-)$ by $(-)^\sigma$.

\medskip
\emph{Preservation of Monoidal Unit}.  The isomorphism $\sigma$ preserves the monoidal unit $\clpu$ in \cref{clpu} because the composite
\[\txprod_{i=1}^n \C_{\sigma(i)} \fto{\sigma} \txprod_{i=1}^n \C_i \fto{\clpu} \D\]
is the constant functor at the monoidal unit $\pu$ in $\D$, since $\clpu$ is constant at $\pu$.  By definition \cref{fsigmatwoj}, for each $i \in \{1,\ldots,n\}$ the $i$-th linearity constraint of $\clpu^\sigma$ is the $\sigma(i)$-th linearity constraint of $\clpu$, which is componentwise given by $1_\pu$ in $\D$.  Thus $\clpu^\sigma$ is the monoidal unit in $\clpsu\scmap{\angC\sigma;\D}$.

\medskip
\emph{Preservation of Monoidal Product}.  For $n$-linear functors $P,Q \cn \angC \to \D$ as in \cref{PQnlinear}, there is an equality of $n$-linear functors
\[P^\sigma \oplus Q^\sigma = (P \oplus Q)^\sigma \cn \angC\sigma \to \D\]
for the following reasons.  First, by definitions \cref{permiicatsigmaaction,PplusQfunctor}, each of the functors $P^\sigma \oplus Q^\sigma$ and $(P \oplus Q)^\sigma$ is given by the following composite.
\begin{equation}\label{PplusQsigma}
\txprod_{i=1}^n \C_{\sigma(i)} \fto{\sigma} \txprod_{i=1}^n \C_i 
\fto{(P,Q)} \D \times \D \fto{\oplus} \D
\end{equation}

Next, to show that $P^\sigma \oplus Q^\sigma$ and $(P \oplus Q)^\sigma$ have the same $i$-th linearity constraint for each $i \in \{1,\ldots,n\}$, we consider objects
\[\angx = \ang{x_j}_{j=1}^n \in \txprod_{j=1}^n \C_{\sigma(j)} \andspace x_i' \in \C_{\sigma(i)}\]
and use the notation $x_i'' = x_i \oplus x_i'$.  For these objects, by definition \cref{fsigmatwoj} and \cref{PplusQtwoi}, each of $(P^\sigma \oplus Q^\sigma)^2_i$ and $\big((P \oplus Q)^\sigma\big)^2_i$ is given by the following composite morphism in $\D$.
\[\begin{tikzpicture}
\def\h{5.5} \def\u{-.8} \def\t{15}
\draw[0cell=.7]
(0,0) node (a1) {(P \oplus Q)^\sigma \angx \oplus (P \oplus Q)^\sigma \ang{x \compi x_i'}}
(a1)+(\h,0) node (a2) {(P \oplus Q)^\sigma\ang{x \compi x_i''}}
(a1)+(0,\u) node (b1) {P(\sigma\angx) \oplus Q(\sigma\angx) \oplus P\big(\sigma\angx \comp_{\sigma(i)} x_i'\big) \oplus Q\big(\sigma\angx \comp_{\sigma(i)} x_i'\big)}
(b1)+(\h,0) node (b2) {P\big(\sigma\angx \comp_{\sigma(i)} x_i''\big) \oplus Q\big(\sigma\angx \comp_{\sigma(i)} x_i''\big)}
(b1)+(\h/2,-1) node[align=left] (b) {$\phantom{\oplus\,}P(\sigma\angx) \oplus P\big(\sigma\angx \comp_{\sigma(i)} x_i'\big)$\\ $\oplus\, Q(\sigma\angx) \oplus Q\big(\sigma\angx \comp_{\sigma(i)} x_i'\big)$}
;
\draw[1cell=.8]
(a1) edge[-,double equal sign distance] (b1)
(a2) edge[-,double equal sign distance] (b2)
(b1) edge[out=-60,in=180] node[swap,pos=.2] {1 \oplus \xi \oplus 1} (b)
(b) edge[out=0,in=-120] node[swap,pos=.8] {P^2_{\sigma(i)} \oplus Q^2_{\sigma(i)}} (b2)
;
\end{tikzpicture}\]
This proves that $P^\sigma \oplus Q^\sigma$ and $(P \oplus Q)^\sigma$ are equal as $n$-linear functors.

Finally, by \cref{permiicatsigmaaction,theta-plus-psi,PplusQsigma} with the $n$-linear functors $(P,Q)$ replaced by $n$-linear transformations $(\theta,\psi)$, the functor $(-)^\sigma$ preserves the monoidal product of morphisms.
\end{proof}

\section{Multicategorical Evaluation for Permutative Categories}
\label{sec:multiev-permcat}

In \cref{sec:internal-hom-permcat} we constructed internal hom permutative categories and their symmetric group action.  To continue the construction of the closed multicategory structure on $\permcatsu$, in this section we construct its multicategorical evaluation.
\begin{itemize}
\item We define the multicategorical evaluation for small permutative categories in \cref{def:clp-eval}.  We verify their multilinearity in \cref{clp-evaluation}.
\item \cref{expl:domain-chi,expl:clp-partners} provide a thorough description of the function $\chi$ \cref{eval-bijection} in the definition of a closed multicategory, in the context of small permutative categories.  We use this discussion in the next section to establish the evaluation bijection and its equivariance axiom; see \cref{clp-evbij-axiom,clp-evbij-equiv}.
\end{itemize}

\subsection*{Multicategorical Evaluation}

Now we define the evaluation \cref{evangxy} for the internal hom permutative categories $\clpsu\scmap{\angC;\D}$ in \cref{clp-angcd-permutative}.

\begin{definition}\label{def:clp-eval}\index{multilinear!evaluation functor}\index{evaluation!multilinear functor}\index{permutative category!multilinear evaluation}
For small permutative categories $\D$ and $\angC = \ang{\C_i}_{i=1}^n$ with $n \geq 0$, we define the data of an $(n+1)$-linear functor
\begin{equation}\label{clp-ev-angcd}
\clp\scmap{\angC;\D} \times \txprod_{i=1}^n \C_i \fto{\ev_{\angC;\,\D}} \D
\end{equation}
as follows.
\begin{description}
\item[Underlying Functor] For an $n$-linear functor
\[\big(P, \{P^2_i\}_{i=1}^n\big) \cn \txprod_{i=1}^n \C_i  \to \D\]
and an object $\angx \in \txprod_{i=1}^n \C_i$, we define the object
\begin{equation}\label{evangcd-object}
\ev_{\angC;\,\D} \big(P,\angx\big) = P\angx \inspace \D.
\end{equation}
For an $n$-linear transformation $\theta \cn P \to Q$ between $n$-linear functors
\[\big(P, \{P^2_i\}_{i=1}^n\big) \scs \big(Q, \{Q^2_i\}_{i=1}^n\big) \cn 
\txprod_{i=1}^n \C_i  \to \D\]
and a morphism $\angf \cn \angx \to \angy$ in $\txprod_{i=1}^n \C_i$, we define the morphism
\[\ev_{\angC;\,\D} \big(\theta,\angf\big) \cn P\angx \to Q\angy \inspace \D\]
as either one of the following two composites.
\begin{equation}\label{evangcd-thetaangf}

\end{equation}
The diagram \cref{evangcd-thetaangf} commutes by the naturality of $\theta$.
\item[Linearity Constraints]
Suppose $P$ and $Q$ are $n$-linear functors as above, and suppose $\angx \in \txprod_{i=1}^n \C_i$ and $x_i' \in \C_i$ are objects.
\begin{itemize}
\item The first linearity constraint of $\ev_{\angCD}$ is given by the following identity morphism in $\D$.
\begin{equation}\label{evangcd-first-linearity}
%
\end{equation}
\item For each $i \in \{1,\ldots,n\}$, the $(i+1)$-st linearity constraint of $\ev_{\angCD}$ is given by the following $i$-th linearity constraint of $P$.
\begin{equation}\label{evangcd-iplusone-linearity}
%
\end{equation}
\end{itemize}
\end{description}
This finishes the definition of $\ev_{\angC;\, \D}$.
\end{definition}

\begin{explanation}\label{expl:clp-eval-zero}
Consider \cref{def:clp-eval}.
\begin{itemize}
\item If $n=0$, then $\angC$ is empty, and
\[\ev_{\ang{};\,\D} \cn \clp\scmap{\ang{};\D} = \D \to \D\]
is the identity symmetric monoidal functor on $\D$.
\item If $n=1$, then 
\[\ev_{\C;\,\D} \cn \clpsu(\C;\D) \times \C \to \D\] 
is equal to the bilinear functor $\ev_{\C,\D}$ in \cref{ev-bilinear}.\defmark
\end{itemize}
\end{explanation}

\begin{lemma}\label{clp-evaluation}
For small permutative categories $\D$ and $\angC = \ang{\C_i}_{i=1}^n$, the data 
\[\big(\ev_{\angC;\, \D} \scs \big\{(\ev_{\angC;\, \D})^2_i\big\}_{i=1}^{n+1}\big) \cn 
\clp\scmap{\angC;\D} \times \txprod_{i=1}^n \C_i \to \D\] 
in \cref{def:clp-eval} form an $(n+1)$-linear functor.
\end{lemma}

\begin{proof}
By \cref{expl:clp-eval-zero,ev-bilinear} for the case $n=1$, we only need to check statements \cref{clp-evaluation-i,clp-evaluation-ii,clp-evaluation-iii} below for $n > 1$.
\begin{romenumerate}
\item\label{clp-evaluation-i} $\ev_{\angC;\, \D}$ is a functor.
\item\label{clp-evaluation-ii} $(\ev_{\angC;\, \D})^2_i$ is a natural transformation for each $i \in \{1,\ldots,n+1\}$.
\item\label{clp-evaluation-iii} $\ev_{\angC;\, \D}$ satisfies the axioms \cref{nlinearunity,constraintunity,eq:ml-f2-assoc,eq:ml-f2-symm,eq:f2-2by2} for an $(n+1)$-linear functor.
\end{romenumerate}

\medskip
\emph{Statements \cref{clp-evaluation-i,clp-evaluation-ii}}.  These assertions are proved using the proofs for statements \cref{evCD-bilinear-i,evCD-bilinear-ii}, respectively, in the proof of \cref{ev-bilinear}, with an appropriate change of notation.

\medskip
\emph{Statement \cref{clp-evaluation-iii}}.  The unity axiom \cref{nlinearunity} holds for $\ev_{\angC;\, \D}$ by the definitions \cref{clpu} of $\clpu$, \cref{evangcd-object} of $\ev_{\angC;\, \D}(P,\angx)$, and \cref{evangcd-thetaangf} of $\ev_{\angC;\,\D}(\theta,\angf)$.

The constraint unity, associativity, and symmetry axioms, \cref{constraintunity,eq:ml-f2-assoc,eq:ml-f2-symm}, hold for $(\evangCD)^2_1$ because it is the identity natural transformation.  For $(\evangCD)^2_{i+1}$ with $i \in \{1,\ldots,n\}$, these axioms hold by the definition \cref{clpu} of $\clpu$ and the corresponding axioms for $n$-linear functors.

For the constraint 2-by-2 axiom \cref{eq:f2-2by2} for $\evangCD$, we consider distinct indices $i \neq k \in \{1,\ldots,n+1\}$ as follows.
\begin{itemize}
\item If either $i=1$ or $k=1$, then the desired diagram \cref{eq:f2-2by2} commutes by the definition \cref{PplusQtwoi} of $(P \oplus Q)^2_i$, the symmetry axiom \cref{symmoncatsymhexagon}, and the functoriality of $\oplus$ in $\D$, as in the proof of statement \cref{evCD-bilinear-iii} in the proof of \cref{ev-bilinear}.  
\item If $i,k > 1$, then the desired diagram \cref{eq:f2-2by2} is the constraint 2-by-2 diagram for an $n$-linear functor $P$ and the indices $i-1 \neq k-1 \in \{1,\ldots,n\}$, which is commutative.
\end{itemize} 
This proves that $\evangCD$ is an $(n+1)$-linear functor.
\end{proof}

\subsection*{The Function $\chi$ for Permutative Categories}

To prepare for the proofs of the other two axioms of a closed multicategory, \cref{eval-bijection,eval-bij-eq}, for $\permcatsu$ in the next section, here we explain the function $\chi$ in \cref{eval-bijection} in the current context, starting with its domain.

\begin{explanation}[Domain of $\chi$]\label{expl:domain-chi}
We consider the following context:
\begin{itemize}
\item $\M = \psu$, the multicategory of small permutative categories and multilinear functors in \cref{thm:permcatmulticat};
\item $\clM(?;?) = \clpsu(?;?)$, the internal hom permutative categories in \cref{clp-angcd-permutative}; 
\item the right symmetric group action on internal hom in \cref{clp-angcd-strict}; and
\item the multilinear evaluation $\ev_{?;?}$ in \cref{clp-evaluation}.
\end{itemize} 
In this context, we consider small permutative categories
\[\B \scs \quad \angC = \ang{\C_i}_{i=1}^n \scs \andspace \angD = \ang{\D_j}_{j=1}^p;\]
morphisms
\begin{equation}\label{ang-wxyz}
\angw \fto{\angf} \angx \in \txprod_{i=1}^n \C_i \andspace 
\angy \fto{\angg} \angz \in \txprod_{j=1}^p \D_j;
\end{equation}
and an $n$-linear functor (\cref{def:nlinearfunctor})
\begin{equation}\label{PPtwoi-nlinear}
\big(P, \{P^2_i\}_{i=1}^n\big) \cn \txprod_{i=1}^n \C_i \to \clpsu\pangDB.
\end{equation}
We explain the $n$-linear functor $\big(P, \{P^2_i\}_{i=1}^n\big)$ and establish some notation.

\medskip
\emph{Underlying Functor}.  By \cref{def:clp-angcd}, the object $P\angw$ in $\clpsu\pangDB$ is a $p$-linear functor, and $P\angf$ is a $p$-linear transformation (\cref{def:nlineartransformation}), as in the left diagram below.
\begin{equation}\label{PangfPtwoi}
\begin{tikzpicture}[baseline={(a.base)}]
\def\t{24} \def\s{27}
\draw[0cell=.9]
(0,0) node (a') {\txprod_{j=1}^p \D_j}
(a')+(.4,.03) node (a) {\phantom{B}}
(a)+(2.2,0) node (b) {\B}
;
\draw[1cell=.8]
(a) edge[bend left=\t] node {\big(P\angw, \{P\angw^2_j\}_{j=1}^p\big)} (b)
(a) edge[bend right=\t] node[swap] {\big(P\angx, \{P\angx^2_j\}_{j=1}^p\big)} (b)
;
\draw[2cell=.9]
node[between=a and b at .37, rotate=-90, 2label={above,P\angf}] {\Rightarrow}
;
\begin{scope}[shift={(4.5,0)}]
\draw[0cell=.9]
(0,0) node (a') {\txprod_{j=1}^p \D_j}
(a')+(.4,.03) node (a) {\phantom{B}}
(a)+(2.2,0) node (b) {\B}
;
\draw[1cell=.8]
(a) edge[bend left=\s] node {P\angw \oplus P\angwp} (b)
(a) edge[bend right=\s] node[swap] {P\angwpp} (b)
;
\draw[2cell=.75]
node[between=a and b at .28, rotate=-90, 2label={above,(P^2_i)_{\angw;\, x_i}}] {\Rightarrow}
;
\end{scope}
\end{tikzpicture}
\end{equation}

\medskip
\emph{Linearity Constraints}.  For indices $i \in \{1,\ldots,n\}$ and $j \in \{1,\ldots,p\}$, we use the following notation.
\begin{equation}\label{angwp-angyp}
\begin{aligned}
\angwp &= \ang{w \compi x_i} & \angyp &= \ang{y \compj z_j}\\
\angwpp &= \bang{w \compi (w_i \oplus x_i)} & \angypp &= \bang{y \compj (y_j \oplus z_j)}\\
\big(P\angw\big)\angy &= P\angw\angy &&
\end{aligned}
\end{equation}
For each $i \in \{1,\ldots,n\}$, the $i$-th linearity constraint $P^2_i$ of $P$ is a natural transformation, with a typical component given by a $p$-linear transformation, as in the right diagram in \cref{PangfPtwoi} above.  Next we unpack its multilinearity axioms.

\medskip
\emph{Unity}.  The axiom \cref{niitransformationunity} for the $p$-linear transformation $(P^2_i)_{\angw;\, x_i}$ says that, if any $y_j = \pu$ in $\D_j$, then the $\angy$-component
\begin{equation}\label{Ptwoi-angwxi-unity}
P\angw\angy \oplus P\angwp\angy 
\fto{\big((P^2_i)_{\angw;x_i}\big)_{\angy}} P\angwpp\angy
\end{equation}
is equal to $1_\pu$ in $\B$.

\medskip
\emph{Constraint Compatibility}.  For the axiom \cref{eq:monoidal-in-each-variable} for the $p$-linear transformation $(P^2_i)_{\angw;\, x_i}$, first recall the linearity constraints of $P\angw \oplus P\angwp$ defined in \cref{PplusQtwoi}.  The diagram \cref{eq:monoidal-in-each-variable} for $(P^2_i)_{\angw;\, x_i}$ is the following commutative diagram in $\B$.
\begin{equation}\label{Ptwoi-angwxi-compat}
\begin{tikzpicture}[vcenter]
\def\h{3} \def\u{1} \def\v{1.5}
\draw[0cell=.75]
(0,0) node (a1) {P\angw\angy \oplus P\angwp\angy \oplus P\angw\angyp \oplus P\angwp\angyp}
(a1)+(\h,\u) node (a2) {P\angwpp\angy \oplus P\angwpp\angyp}
(a1)+(\h+.5,-\v/2) node (a3) {P\angwpp\angypp}
(a1)+(0,-\v) node (b1) {P\angw\angy \oplus P\angw\angyp \oplus P\angwp\angy \oplus P\angwp\angyp}
(b1)+(\h,-\u) node (b2) {P\angw\angypp \oplus P\angwp\angypp}
;
\draw[1cell=.75]
(a1) edge[transform canvas={xshift=-1em}] node {\big((P^2_i)_{\angw;\, x_i}\big)_{\angy} \oplus \big((P^2_i)_{\angw;\, x_i}\big)_{\angyp}} (a2)
(a2) edge node[pos=.75] {\big(P\angwpp^2_j\big)_{\angy;\, z_j}} (a3)
(a1) edge node[swap] {1 \oplus \xi \oplus 1} (b1)
(b1) edge[transform canvas={xshift=-1em}] node[swap] {\big(P\angw^2_j\big)_{\angy;\, z_j} \oplus \big(P\angwp^2_j\big)_{\angy;\, z_j}} (b2)
(b2) edge node[swap,pos=.75] {\big((P^2_i)_{\angw;\, x_i}\big)_{\angypp}} (a3)
;
\end{tikzpicture}
\end{equation}
This finishes our description of the $n$-linear functor $\big(P,\{P^2_i\}_{i=1}^n\big)$.
\end{explanation}

\begin{explanation}[The Function $\chi$]\label{expl:clp-partners}
In the context of \cref{expl:domain-chi}, we consider the function defined in \cref{eval-bijection}
\begin{equation}\label{chi-angCDB}
\psu\left( \angC \sscs \clpsu\pangDB\right) \fto{\chi_{\angC;\, \angD;\, \B}} \psu\scmap{\angC,\angD;\B},
\end{equation}
which we abbreviate to $\chi$.  This function sends an $n$-linear functor 
\[\big(P, \{P^2_i\}_{i=1}^n\big) \cn \txprod_{i=1}^n \C_i \to \clpsu\pangDB\]
to the $(n+p)$-linear functor $\chi P$ given by the following composite, where we suppress the associativity isomorphism for the Cartesian product.
\begin{equation}\label{chi-of-P}
\begin{tikzpicture}[baseline={(a.base)}]
\def\h{2.5} \def\w{.6}
\draw[0cell=.9]
(0,0) node (a) {\txprod_{i=1}^n \C_i \times \txprod_{j=1}^p \D_j}
(a)+(4,0) node (b) {\clpsu\pangDB \times \txprod_{j=1}^p \D_j}
(b)+(3.2,0) node (c) {\phantom{\B}}
(c)+(0,.06) node (c') {\B}
;
\draw[1cell=.9]
(a) edge node {P \times 1} (b)
(b) edge node {\ev_{\angDB}} (c)
;
\draw[1cell=.9]
(a) [rounded corners=3pt, shorten <=-.4ex] |- ($(b)+(-1,\w)$)
-- node[pos=.2] {\big(\chi P \scs \{(\chi P)^2_r\}_{r=1}^{n+p}\big)} ($(b)+(1,\w)$) -| (c')
;
\end{tikzpicture}
\end{equation}
Next we unpack this $(n+p)$-linear functor, using the notation in \cref{angwp-angyp}.

\medskip
\emph{Underlying Functor}.  By \cref{def:clp-eval}, on objects $\chi P$ is given by
\begin{equation}\label{chiP-objects}
\begin{aligned}
(\chi P)\big(\angw,\angy\big) &= \ev_{\angDB}\big(P\angw,\angy\big)\\
&= P\angw\angy \inspace \B.
\end{aligned}
\end{equation}
The morphism 
\[(\chi P)\big(\angf,\angg\big) = \ev_{\angDB}\big(P\angf,\angg\big)\]
is given by either one of the following two boundary composites in $\B$, which are equal by the naturality of $P\angf$.
\begin{equation}\label{chiP-morphisms}

\end{equation}

\emph{Linearity Constraints}.  The linearity constraints for composite multilinear functors are defined in \cref{ffjlinearity}.  By \cref{evangcd-first-linearity}, for each $i \in \{1,\ldots,n\}$, the $i$-th linearity constraint $(\chi P)^2_i$ is a natural transformation, with a typical component given by a component morphism in $\B$ of $(P^2_i)_{\angw;\, x_i}$---which is a component of the $i$-th linearity constraint $P^2_i$ in \cref{PangfPtwoi}---as follows.
\begin{equation}\label{chiPtwoi}
%
\end{equation}
By \cref{evangcd-iplusone-linearity}, for each $j \in \{1,\ldots,p\}$, the $(n+j)$-th linearity constraint $(\chi P)^2_{n+j}$ is a natural transformation, with a typical component given by a component morphism in $\B$ of $P\angw^2_j$---which is the $j$-th linearity constraint of $P\angw$ in \cref{PangfPtwoi}---as follows.
\begin{equation}\label{chiPtwonplusj}
%
\end{equation}

\medskip
\emph{Multilinearity Axioms}.  The axioms \cref{nlinearunity,constraintunity,eq:ml-f2-assoc,eq:ml-f2-symm,eq:f2-2by2} of an $(n+p)$-linear functor hold for $\chi P$ for the following reasons.
\begin{itemize}
\item The unity axiom \cref{nlinearunity} holds for $\chi P$ by the same axiom for the $n$-linear functor $P$ and each $p$-linear functor $P\angw$.
\item The constraint unity axiom \cref{constraintunity} holds for $\chi P$ by 
\begin{itemize}
\item the axiom \cref{constraintunity} for the $n$-linear functor $P$ and each $p$-linear functor $P\angw$, 
\item the unity axiom \cref{Ptwoi-angwxi-unity} for the $p$-linear transformation $(P^2_i)_{\angw;\, x_i}$, and
\item the unity axiom \cref{nlinearunity} for $P$.
\end{itemize}
\item The constraint associativity and symmetry axioms, \cref{eq:ml-f2-assoc,eq:ml-f2-symm}, hold for $\chi P$ by the same axioms for the $n$-linear functor $P$ and each $p$-linear functor $P\angw$.
\item The constraint 2-by-2 axiom \cref{eq:f2-2by2} for $\chi P$ and distinct indices 
\[r \neq s \in \{1,\ldots,n+p\}\]
holds for the following reasons.
\begin{itemize}
\item If both $r,s \in \{1,\ldots,n\}$, then, by \cref{chiPtwoi}, the axiom \cref{eq:f2-2by2} for $\chi P$ follows from the same axiom for the $n$-linear functor $P$. 
\item If both $r,s \in \{n+1,\ldots,n+p\}$, then, by \cref{chiPtwonplusj}, the axiom \cref{eq:f2-2by2} for $\chi P$ follows from the same axiom for the $p$-linear functor $P\angw$ and indices
\[r-n \neq s-n \in \{1,\ldots,p\}.\]
\item If $r \in \{1,\ldots,n\}$ and $s \in \{n+1,\ldots,n+p\}$, then the axiom \cref{eq:f2-2by2} for $\chi P$ follows from the commutative diagram \cref{Ptwoi-angwxi-compat}. 
\item If $r \in \{n+1,\ldots,n+p\}$ and $s \in \{1,\ldots,n\}$, then we again use the commutative diagram \cref{Ptwoi-angwxi-compat} but with the arrow $1 \oplus \xi \oplus 1$ reversed, which is possible by the symmetry axiom \cref{symmoncatsymhexagon} for $\B$.
\end{itemize}
\end{itemize}
This finishes our description of the function $\chi$ in \cref{chi-angCDB}.
\end{explanation}

\section{The Closed Multicategory Structure}
\label{sec:permcatsu-clmulti}

In this section we complete the construction of the closed multicategory structure on $\permcatsu$.
\begin{itemize}
\item In \cref{clp-evbij-axiom} we prove the evaluation bijection axiom.
\item In \cref{clp-evbij-equiv} we prove the equivariance axiom for evaluation bijection.
\item The closed multicategory structure on $\permcatsu$ is stated in \cref{perm-closed-multicat}.
\end{itemize}

\subsection*{Multicategorical Evaluation Axioms}

\begin{lemma}\label{clp-evbij-axiom}\index{permutative category!evaluation bijection}
The evaluation bijection axiom \cref{eval-bijection} holds for
\[\big(\permcatsu \scs \clp \scs \ev \big)\]
defined in \cref{thm:permcatmulticat,clp-internal-hom,clp-ev-angcd}
\end{lemma}

\begin{proof}
We need to show that, for small permutative categories $\B$, $\angC = \ang{\C_i}_{i=1}^n$, and $\angD = \ang{\D_j}_{j=1}^p$, the function $\chi = \chi_{\angC;\, \angD;\, \B}$ in \cref{chi-angCDB}, as displayed below, is a bijection.
\begin{equation}\label{chiPsi}
\begin{tikzpicture}[vcenter]
\draw[0cell]
(0,0) node (a) {\psu\left( \angC \sscs \clpsu\pangDB\right)}
(a)+(4.5,0) node (b) {\psu\scmap{\angC,\angD;\B}}
;
\draw[1cell=.9]
(a) edge[transform canvas={yshift=.6ex}] node {\chi} (b)
(b) edge[transform canvas={yshift=-.4ex}] node {\Psi} (a)
;
\end{tikzpicture}
\end{equation}
We show that $\chi$ is a bijection by constructing an explicit inverse $\Psi$ in the following two steps.
\begin{enumerate}
\item\label{chiinv-i}  We construct a function $\Psi$ that goes in the opposite direction as $\chi$.
\item\label{chiinv-ii}  We show that $\chi$ and $\Psi$ are inverses of each other.
\end{enumerate}

\stepproof{\cref{chiinv-i}}{The Function $\Psi$}

Given an $(n+p)$-linear functor
\begin{equation}\label{Rnplusp-linear}
\big(R, \{R^2_r\}_{r=1}^{n+p}\big) \cn \txprod_{i=1}^n \C_i \times \txprod_{j=1}^p \D_j \to \B,
\end{equation}
we define an $n$-linear functor
\begin{equation}\label{PsiR-nlinear}
\big( \Psi R, \{(\Psi R)^2_i\}_{i=1}^n \big) \cn \txprod_{i=1}^n \C_i \to \clpsu\scmap{\angD;\B}
\end{equation}
in steps \cref{chiinv-i-i,chiinv-i-ii,chiinv-i-iii,chiinv-i-iv} below.
\begin{romenumerate}
\item\label{chiinv-i-i} We define $\Psi R$ on objects in \cref{PsiRangw,Rtwonplusj}.
\item\label{chiinv-i-ii} We define $\Psi R$ on morphisms in \cref{PsiRangf-component}.
\item\label{chiinv-i-iii} We define the linearity constraints $(\Psi R)^2_i$ in \cref{PsiRtwoi-component}.
\item\label{chiinv-i-iv} We check the $n$-linear functor axioms for $\big( \Psi R, \{(\Psi R)^2_i\}_{i=1}^n \big)$.
\end{romenumerate}
In the rest of this proof, we use the notation in \cref{ang-wxyz,angwp-angyp} for objects and morphisms, so we ask the reader to briefly review them.

\medskip
\emph{Step \cref{chiinv-i}\cref{chiinv-i-i}: Objects}.  For an object $\angw$ in $\txprod_{i=1}^n \C_i$, we define the functor
\begin{equation}\label{PsiRangw}
\txprod_{j=1}^p \D_j \fto{(\Psi R)\angw = R\big(\angw,-\big)} \B.
\end{equation}
For each $j \in \{1,\ldots,p\}$, its $j$-th linearity constraint $(\Psi R)\angw^2_j$ is defined componentwise by the $(n+j)$-th linearity constraint of $R$, as indicated below.
\begin{equation}\label{Rtwonplusj}

\end{equation}
This is natural in $\angy$ and $z_j$ because $R^2_{n+j}$ is a natural transformation.  The $p$-linear functor axioms---\cref{nlinearunity,constraintunity,eq:ml-f2-assoc,eq:ml-f2-symm,eq:f2-2by2}---for $(\Psi R)\angw$ defined in \cref{PsiRangw,Rtwonplusj} follow from the corresponding axioms for the given $(n+p)$-linear functor $R$ in \cref{Rnplusp-linear}.

\medskip
\emph{Step \cref{chiinv-i}\cref{chiinv-i-ii}: Morphisms}.  For a morphism $\angf \cn \angw \to \angx$ in $\txprod_{i=1}^n \C_i$, we define the natural transformation
\[%
\]
with, for each object $\angy$ in $\txprod_{j=1}^p \D_j$, $\angy$-component given by the following morphism in $\B$.
\begin{equation}\label{PsiRangf-component}
%
\end{equation}
Its naturality in $\angy$ follows from the functoriality of $R$.  

Moreover, $(\Psi R)\angf$ is a $p$-linear transformation (\cref{def:nlineartransformation}) for the following reasons.
\begin{itemize}
\item By definition \cref{PsiRangf-component}, the unity axiom \cref{niitransformationunity} for $(\Psi R)\angf$ says that, if any $y_j = \pu$ in $\D_j$, then
\[R\big(\angf,1_{\angy}\big) = 1_\pu \inspace \B.\]
This is true by the unity axiom \cref{nlinearunity} for the $(n+p)$-linear functor $R$.
\item By definition \cref{Rtwonplusj}, the constraint compatibility axiom \cref{eq:monoidal-in-each-variable} for $(\Psi R)\angf$ is the following diagram in $\B$, which commutes by the naturality of $R^2_{n+j}$.
\[
\]
\end{itemize}
Thus $(\Psi R)\angf$ is a $p$-linear transformation.  The functoriality of $\Psi R$ follows from the definition \cref{PsiRangf-component} and the functoriality of $R$.

\medskip
\emph{Step \cref{chiinv-i}\cref{chiinv-i-iii}: Linearity Constraints}.  For each $i \in \{1,\ldots,n\}$ and objects $\angw \in \txprod_{i=1}^n \C_i$ and $x_i \in \C_i$, we define the $i$-th linearity constraint $(\Psi R)^2_i$ as having the following component $p$-linear transformation.
\[%
\]
For each object $\angy \in \txprod_{j=1}^p \D_j$, its $\angy$-component is the corresponding component morphism of the $i$-th linearity constraint $R^2_i$, as indicated below.
\begin{equation}\label{PsiRtwoi-component}
%
\end{equation}
The naturality of
\begin{itemize}
\item $\left(\big((\Psi R)^2_i\big)_{\angw;\, x_i}\right)_{\angy}$ with respect to $\angy$ and
\item $\big((\Psi R)^2_i\big)_{\angw;\, x_i}$ with respect to $\angw$ and $x_i$
\end{itemize}
follows from the naturality of $R^2_i$.  

The components \cref{PsiRtwoi-component} define a $p$-linear transformation $\big((\Psi R)^2_i\big)_{\angw;\, x_i}$ for the following reasons.
\begin{itemize}
\item The unity axiom \cref{niitransformationunity} for $\big((\Psi R)^2_i\big)_{\angw;\, x_i}$ says that, if any $y_j = \pu$ in $\D_j$, then 
\[(R^2_i)_{\angw,\angy;\, x_i} = 1_\pu \inspace \B.\]
This holds by the constraint unity axiom \cref{constraintunity} for $R^2_i$.
\item Using \cref{PplusQtwoi} for $(\Psi R)\angw \oplus (\Psi R)\angwp$ and \cref{Rtwonplusj}, the constraint compatibility axiom \cref{eq:monoidal-in-each-variable} for $\big((\Psi R)^2_i\big)_{\angw;\, x_i}$ is the following diagram in $\B$.
\[\begin{tikzpicture}[vcenter]
\def\h{3.5} \def\u{1} \def\v{1.8}
\draw[0cell=.75]
(0,0) node[align=left] (a1) {$\phantom{\oplus\,}R\big(\angw,\angy\big) \oplus R\big(\angwp,\angy\big)$\\ $\oplus\, R\big(\angw,\angyp\big) \oplus R\big(\angwp,\angyp\big)$}
(a1)+(\h,\u) node (a2) {R\big(\angwpp,\angy\big) \oplus R\big(\angwpp,\angyp\big)}
(a1)+(\h+.5,-\v/2) node (a3) {R\big(\angwpp,\angypp\big)}
(a1)+(0,-\v) node[align=left] (b1) {$\phantom{\oplus\,}R\big(\angw,\angy\big) \oplus R\big(\angw,\angyp\big)$\\ $\oplus\, R\big(\angwp,\angy\big) \oplus R\big(\angwp,\angyp\big)$}
(b1)+(\h,-\u) node (b2) {R\big(\angw,\angypp\big) \oplus R\big(\angwp,\angypp\big)}
;
\draw[1cell=.75]
(a1) edge[transform canvas={xshift=-2em}] node {(R^2_i)_{\angw,\angy;\, x_i} \oplus (R^2_i)_{\angw,\angyp;\, x_i}} (a2)
(a2) edge node[pos=.7] {(R^2_{n+j})_{\angwpp,\angy;\, z_j}} (a3)
(a1) edge node[swap] {1 \oplus \xi \oplus 1} (b1)
(b1) edge[transform canvas={xshift=-2em}] node[swap] {(R^2_{n+j})_{\angw,\angy;\, z_j} \oplus (R^2_{n+j})_{\angwp,\angy;\, z_j}} (b2)
(b2) edge node[swap,pos=.7] {(R^2_i)_{\angw,\angypp;\, x_i}} (a3)
;
\end{tikzpicture}\]
This diagram commutes by the constraint 2-by-2 axiom \cref{eq:f2-2by2} for $R$.
\end{itemize}
Thus $\big((\Psi R)^2_i\big)_{\angw;\, x_i}$ is a $p$-linear transformation.

\medskip
\emph{Step \cref{chiinv-i}\cref{chiinv-i-iv}: Multilinearity Axioms}.  The data $\big(\Psi R, \{(\Psi R)^2_i\}\big)$ satisfy the axioms \cref{nlinearunity,constraintunity,eq:ml-f2-assoc,eq:ml-f2-symm,eq:f2-2by2} for an $n$-linear functor for the following reasons.
\begin{itemize}
\item The unity axiom \cref{nlinearunity} for $\Psi R$ follows from
\begin{itemize}
\item the definitions \cref{PsiRangw,Rtwonplusj,PsiRangf-component},
\item the unity axiom \cref{nlinearunity}, and the constraint unity axiom \cref{constraintunity} for $R$.
\end{itemize} 
\item Each of the constraint unity, associativity, symmetry, and 2-by-2 axioms---\cref{constraintunity,eq:ml-f2-assoc,eq:ml-f2-symm,eq:f2-2by2}---for $\Psi R$ follows from the same axiom for $R$ and the definition \cref{PsiRtwoi-component}.
\end{itemize}
This finishes the construction of the $n$-linear functor $\big(\Psi R, \{(\Psi R)^2_i\}\big)$ in \cref{PsiR-nlinear} and completes step \cref{chiinv-i}.

\stepproof{\cref{chiinv-ii}}{$\chi$ and $\Psi$ are Mutual Inverses}

This consists of the following two steps.
\begin{romenumerate}
\item\label{chiinv-ii-i} We show the equality
\[\Psi \chi = 1 \cn \psu\left(\angC \sscs \clpsu\pangDB\right) \to \psu\left(\angC \sscs \clpsu\pangDB\right).\]
\item\label{chiinv-ii-ii} We show the equality
\[\chi \Psi = 1 \cn \psu\scmap{\angC,\angD;\B} \to \psu\scmap{\angC,\angD;\B}.\]
\end{romenumerate}

\medskip
\emph{Step \cref{chiinv-ii}\cref{chiinv-ii-i}}.  For an $n$-linear functor
\[\big(P, \{P^2_i\}_{i=1}^n\big) \cn \txprod_{i=1}^n \C_i \to \clpsu\pangDB\]
as in \cref{PPtwoi-nlinear}, we want to show that $\Psi \chi P$ is equal to $P$ as $n$-linear functors.
\begin{itemize}
\item $\Psi \chi P$ is equal to $P$ on the objects of $\txprod_{i=1}^n \C_i$ by \cref{chiP-objects,PsiRangw}.
\item $\Psi \chi P$ is equal to $P$ on each morphism $\angf$ in $\txprod_{i=1}^n \C_i$ because 
\[(\Psi \chi P)\angf_? = (\chi P)\big(\angf, 1_?\big) = P\angf_?\]
by \cref{chiP-morphisms,PsiRangf-component}. 
\item For each $i \in \{1,\ldots,n\}$, the $i$-th linearity constraints of $\Psi \chi P$ and $P$ are equal by \cref{chiPtwoi,PsiRtwoi-component}.
\end{itemize}
This shows that $\Psi \chi$ is equal to the identity function.

\medskip
\emph{Step \cref{chiinv-ii}\cref{chiinv-ii-ii}}.  For an $(n+p)$-linear functor
\[\big(R, \{R^2_r\}_{r=1}^{n+p}\big) \cn \txprod_{i=1}^n \C_i \times \txprod_{j=1}^p \D_j \to \B\]
as in \cref{Rnplusp-linear}, we want to show that $\chi \Psi R$ is equal to $R$ as $(n+p)$-linear functors.
\begin{itemize}
\item $\chi\Psi R$ is equal to $R$ on the objects of $\txprod_{i=1}^n \C_i \times \txprod_{j=1}^p \D_j$ by \cref{chiP-objects,PsiRangw}.
\item $\chi\Psi R$ is equal to $R$ on each morphism 
\[\big(\angf,\angg\big) \cn \big(\angw,\angy\big) \to \big(\angx,\angz\big) \inspace \txprod_{i=1}^n \C_i \times \txprod_{j=1}^p \D_j\]
by the following computation.
\[\begin{aligned}
&\phantom{==}(\chi\Psi R)\big(\angf,\angg\big) & \phantom{MM} &\\ 
&= (\Psi R)\angx\angg \circ (\Psi R)\angf_{\angy} && \text{by \cref{chiP-morphisms}}\\
&= R\big(1_{\angx},\angg\big) \circ R\big(\angf,1_{\angy}\big) && \text{by \cref{PsiRangw,PsiRangf-component}}\\
&= R\big(\angf,\angg\big) && \text{by functoriality of $R$}
\end{aligned}\]
\item To check that $\chi\Psi R$ and $R$ have the same $n+p$ linearity constraints, we consider the following two cases.
\begin{itemize}
\item For each $i \in \{1,\ldots,n\}$, their $i$-th linearity constraints are equal by \cref{chiPtwoi,PsiRtwoi-component}.
\item For each $j \in \{1,\ldots,p\}$, their $(n+j)$-th linearity constraints are equal by \cref{chiPtwonplusj,Rtwonplusj}.
\end{itemize} 
\end{itemize}
This shows that $\chi\Psi$ is equal to the identity function and completes step \cref{chiinv-ii}.  The proof of the Lemma is now complete.
\end{proof}

\begin{lemma}\label{clp-evbij-equiv}
The equivariance axiom for evaluation bijection \cref{eval-bij-eq} holds for
\[\big(\permcatsu \scs \clp \scs \ev \big)\]
defined in \cref{thm:permcatmulticat,clp-internal-hom,clp-sigma-action,clp-ev-angcd}.
\end{lemma}

\begin{proof}
For small permutative categories $\B$, $\angC = \ang{\C_i}_{i=1}^n$, and $\angD = \ang{\D_j}_{j=1}^p$ and permutations $(\sigma,\vsi) \in \Sigma_n \times \Sigma_p$, the desired commutative diagram is the following.
\begin{equation}\label{clp-evbij-eq-diagram}
\begin{tikzpicture}[vcenter]
\def\h{5.5} \def\v{-1.3}
\draw[0cell=.9]
(0,0) node (a) {\psu\left(\angC \sscs \clpsu\scmap{\angD;\B} \right)}
(a)+(\h,0) node (b) {\psu\scmap{\angC, \angD; \B}}
(a)+(0,\v) node (c) {\psu\left(\angC\sigma \sscs \clpsu\scmap{\angD;\B} \right)}
(c)+(0,\v) node (d) {\psu\left(\angC\sigma \sscs \clpsu\scmap{\angD\vsi;\B} \right)}
(d)+(\h,0) node (e) {\psu\scmap{\angC\sigma, \angD\vsi; \B}}
;
\draw[1cell=.85]
(a) edge node {\chi_{\angC;\, \angD;\, \B}} (b)
(b) edge node {\sigma \times \vsi} (e)
(a) edge node[swap] {\sigma} (c)
(c) edge node[swap] {\ga(\vsi \sscs -)} (d)
(d) edge node {\chi_{\angC\sigma;\, \angD\vsi;\, \B}} (e)
;
\end{tikzpicture}
\end{equation}
An object in the upper left node in \cref{clp-evbij-eq-diagram} is an $n$-linear functor 
\[\big(P, \{P^2_i\}_{i=1}^n\big) \cn \txprod_{i=1}^n \C_i \to \clpsu\pangDB\]
as in \cref{PPtwoi-nlinear}.  We need to check the following two statements.
\begin{romenumerate}
\item\label{clp-evbij-eq-i} The two composites in \cref{clp-evbij-eq-diagram} applied to $P$ are equal as functors.
\item\label{clp-evbij-eq-ii} The two functors in \cref{clp-evbij-eq-i} have the same $n+p$ linearity constraints.
\end{romenumerate}

\medskip
\emph{Statement \cref{clp-evbij-eq-i}}.  To show that the two composites in \cref{clp-evbij-eq-diagram} applied to $P$ are equal as functors, we consider the following diagram.
\[\begin{tikzpicture}
\def\g{2} \def\u{1.3}
\draw[0cell=.8]
(0,0) node (a1) {\txprod_{i=1}^n \C_{\sigma(i)} \times \txprod_{j=1}^p \D_{\vsi(j)}}
(a1)+(\g,\u) node (a2) {\txprod_{i=1}^n \C_i \times \txprod_{j=1}^p \D_j}
(a2)+(\g,-\u) node (a3) {\clpsu\scmap{\angD;\B} \times \txprod_{j=1}^p \D_j}
(a3)+(\g,\u) node (a4) {\B}
(a4)+(\g,-\u) node (a5) {\clpsu\scmap{\angD\vsi;\B} \times \txprod_{j=1}^p \D_{\vsi(j)}}
(a1)+(\g,-\u) node (b1) {\txprod_{i=1}^n \C_i \times \txprod_{j=1}^p \D_{\vsi(j)}}
(a5)+(-\g,-\u) node (b2) {\clpsu\scmap{\angD;\B} \times \txprod_{j=1}^p \D_{\vsi(j)}}
;
\draw[1cell=.8]
(a1) edge node[pos=.4] {\sigma \times \vsi} (a2)
(a2) edge node {P \times 1} (a3)
(a3) edge node {\ev} (a4)
(a1) edge node[swap,pos=.3] {\sigma \times 1} (b1)
(b1) edge node {P \times 1} (b2)
(b2) edge node[swap,pos=.7] {\vsi \times 1} (a5)
(a5) edge node[swap] {\ev} (a4)
(b1) edge node[pos=.3] {P \times \vsi} (a3)
(a3) edge node[pos=.7] {1 \times \vsi^\inv} (b2)
;
\end{tikzpicture}\]
By definitions \cref{permiicatsigmaaction,chi-of-P}, the following statements hold for the diagram above.
\begin{itemize}
\item Denoting $\chi = \chi_{\angC;\, \angD;\, \B}$, the top boundary composite is the functor
\begin{equation}\label{chiP-sigmavsi}
(\chi P)(\sigma \times \vsi).
\end{equation}
Denoting $\chi' = \chi_{\angC\sigma;\, \angD\vsi;\, \B}$, the other boundary composite is the functor
\begin{equation}\label{chiprimeP}
\chi'\big(\ga(\vsi \sscs P\sigma)\big).
\end{equation}
\item The left quadrilateral and the lower middle triangle commute by functoriality of the Cartesian product.
\item The right quadrilateral commutes by the definition \cref{permiicatsigmaaction} of the right $\vsi$-action on $\psu\scmap{\angD;\B}$, which is the same as the right $\vsi$-action on $\clpsu\scmap{\angD;\B}$ \cref{clp-sigma-action}.
\end{itemize}
This proves statement \cref{clp-evbij-eq-i}.

\medskip
\emph{Statement \cref{clp-evbij-eq-ii}}.  We consider indices $r \in \{1,\ldots,n\}$ and $t \in \{1,\ldots,p\}$ and objects
\[\big(\anga \scs \angb\big) \in \txprod_{i=1}^n \C_{\sigma(i)} \times \txprod_{j=1}^p \D_{\vsi(j)}
\andspace (a_r' \scs b_t') \in \C_{\sigma(r)} \times \D_{\vsi(t)}.\]
The following equalities follow from \cref{fsigmatwoj,chiPtwoi,chiPtwonplusj}.
\begin{align}
  \begin{split}
    \left(\big((\chi P)(\sigma \times \vsi)\big)^2_r\right)_{\anga,\angb;\, a_r'} 
    &= \left((P^2_{\sigma(r)})_{\sigma\anga;\, a_r'}\right)_{\vsi\angb}\\
    &= \left( \big(\chi' (\ga(\vsi \sscs P\sigma))\big)^2_r\right)_{\anga,\angb;\, a_r'}
  \end{split}\label{eq:evbij-lin-r}\\
  \begin{split} 
    \left(\big((\chi P)(\sigma \times \vsi)\big)^2_{n+t}\right)_{\anga,\angb;\, b_t'}
    &= \left( P(\sigma\anga)^2_{\vsi(t)} \right)_{\vsi\angb;\, b_t'}\\
    &= \left( \big(\chi' (\ga(\vsi \sscs P\sigma))\big)^2_{n+t}\right)_{\anga,\angb;\, b_t'}
  \end{split}\label{eq:evbij-lin-n+t}
\end{align} 
Thus the functors in \cref{chiP-sigmavsi,chiprimeP} have the same $n+p$ linearity constraints, proving statement \cref{clp-evbij-eq-ii}.
\end{proof}

We are now ready to show that small permutative categories form a closed multicategory (\cref{def:closed-multicat}).

\begin{theorem}\label{perm-closed-multicat}\index{closed!multicategory!- of permutative categories}\index{multicategory!closed - of permutative categories}\index{permutative category!closed multicategory of}
There is a closed multicategory
\[\big(\permcatsu \scs \clp \scs \ev \big)\]
consisting of the following data.
\begin{itemize}
\item The underlying multicategory is $\permcatsu$ in \cref{thm:permcatmulticat}.
\item The internal hom objects are the permutative categories $\clp\scmap{\angC;\D}$ in \cref{clp-angcd-permutative}.
\item The symmetric group action on internal hom is given in \cref{clp-angcd-strict}.
\item The multicategorical evaluations are the multilinear functors
\[\clp\scmap{\angC;\D} \times \txprod_{i=1}^n \C_i \fto{\ev_{\angC;\,\D}} \D\]
in \cref{clp-evaluation}.
\end{itemize}
\end{theorem}

\begin{proof}
The data are well defined by the indicated results, \cref{thm:permcatmulticat,clp-angcd-permutative,clp-angcd-strict,clp-evaluation}.  The closed multicategory axioms, \cref{right-id-action,hom-equivariance,eval-bijection,eval-bij-eq}, hold by \cref{clp-angcd-strict,clp-evbij-axiom,clp-evbij-equiv}.
\end{proof}

\section{Closed Multicategories of Lax and Strong Multilinear Functors}
\label{sec:lax-strong-multicat}

In this section we describe closed multicategory structures for $\permcat$, $\permcatsg$, and $\permcatsus$.
This material extends that of \cref{sec:multpermcat,sec:permcatsu-clmulti} to multilinear functors with generally non-trivial unit constraints.
Our treatment below describes the relevant data and axioms generalizing those of the strictly unital case.

Recall from \cref{proposition:n-lin-equiv} that multilinear functors can be identified with multifunctors out of a smash product of endomorphism multicategories.
In \cref{lemma:lax-nlinear-tensor,proposition:laxmultilin-tensor} we show that lax multilinear functors have a corresponding description using tensor products in place of smash products.
The main results of this section are stated in \cref{theorem:permcat-catmulticat,lax-perm-closed-multicat}.

\subsection*{$\Cat$-Multicategory Structure}
Recall from \cref{notation:compk} that $\ang{x \compj y}$ denotes the tuple obtained by replacing $x_j$ with $y$.
For the unit constraint below, we require the following similar notation, where the existence of $x_j$ is not assumed.
\begin{notation}\label{notation:insj}
  Suppose given $j \in \{1,\ldots,n\}$ together with symbols
  \[
    x_i \forspace 1 \leq i \leq n \withspace i \ne j,
  \]
  and another symbol $y$.
  We denote by\label{not:insj}
  \begin{equation}\label{eq:insj}
    \ang{x_i}_{i \ne j} \insj y 
    = \big(\underbracket[0.5pt]{x_1, \ldots, x_{j-1}}_{\text{empty if $j=1$}}, y, \underbracket[0.5pt]{x_{j+1}, \ldots, x_n}_{\text{empty if $j=n$}}\big)
  \end{equation}
  the $n$-tuple that has $x_i$ in each entry $i \ne j$, and has $y$ as its $j$-th entry.
  Thus, in the tuple \cref{eq:insj}, no entry of the $(n-1)$-tuple $\ang{x_i}_{i \ne j}$ is \emph{replaced} by $y$.
  Instead, $y$ is inserted as a new entry of $\ang{x_i}_{i \ne j}$, shifting the positions of the entries in positions $i > j$.

  We also use the following further variants:
  \begin{align}
    & \ang{x_i}_{i \ne j} \compk x_k' \insj \pu
    & & \text{has $x_k'$ in position $k$ and $\pu$ in position $j$,}\\
    & \ang{x_i}_{i \ne j} \compk (x_k \oplus x_k') \insj \pu
    & & \text{has $x_k \oplus x_k'$ in position $k$ and $\pu$ in position $j$, and}\\
    & \ang{x_i}_{i \ne j,k} \insk \pu \insj \pu
    & & \text{has $\pu$ in positions $k$ and $j$.}
  \end{align}
\end{notation}

The following generalizes \cref{def:nlinearfunctor}.
\begin{definition}\label{def:lax-nlinearfunctor}
  A \index{n-linear@$n$-linear!functor!lax}\index{functor!n-linear@$n$-linear!lax}\emph{lax $n$-linear functor}
\[\begin{tikzcd}[column sep=3.3cm]
    \txprod_{j=1}^n \C_j \ar{r}{\big(P,\, \{P^2_j\}_{j=1}^n,\, \{P^0_j\}_{j=1}^n\big)} & \D
\end{tikzcd}\]
consists of the following data.
\begin{itemize}
\item $P \cn \prod_{j=1}^n \C_j \to \D$ is a functor.
\item For each $j\in \{1,\ldots,n\}$, $P^2_j$ is a natural transformation, called the \index{constraint!lax linearity}\emph{$j$-th linearity constraint}, with component morphisms
\begin{equation}\label{lax-laxlinearityconstraints}
\begin{tikzcd}[column sep=large]
P\ang{x \compj x_j} \oplus P\ang{x \compj x_j'} \ar{r}{P^2_j}
& P\bang{x \compj (x_j \oplus x_j')} \in \D
\end{tikzcd}
\end{equation}
for objects $\ang{x} \in \txprod_{j=1}^n \C_j$ and $x_j' \in \C_j$.
\item For each $j \in \{1,\ldots,n\}$, $P^0_j$ is a natural transformation, called the \index{constraint!lax unity}\emph{$j$-th unit constraint}, with component morphisms
  \begin{equation}\label{eq:lax-unitconstraints}
    \begin{tikzcd}[column sep=large]
      \pu \ar{r}{P^0_j}
      & P\big( \ang{x_i}_{i \ne j} \insj \pu \big) \in \D
    \end{tikzcd}
  \end{equation}
  for objects $\ang{x_i}_{i\ne j} \in \txprod_{i \ne j} \C_i$.
  In \cref{eq:lax-unitconstraints}, $\pu$ denotes the monoidal unit of $\D$, respectively $\C_j$, in the domain, respectively codomain of $P^0_j$.
\end{itemize}
These data are required to satisfy six axioms.
Of these, three are axioms for $n$-linear functors:
constraint associativity \cref{eq:ml-f2-assoc}, constraint symmetry \cref{eq:ml-f2-symm}, and constraint 2-by-2 \cref{eq:f2-2by2}.
The remaining three axioms are as follows, using \cref{notation:insj}.
\begin{description}
\item[Lax Unity]\index{constraint!- axiom!lax unity} For each $j \in \{1,\ldots,n\}$, and for objects $\ang{x_i}_{i \ne j} \in \txprod_{i \ne j}\C_i$, the following diagrams commute in $\D$.
  \begin{equation}\label{lax-constraintunity}

  \end{equation}
\item[Constraint 0-by-2]\index{constraint!- axiom!0-by-2} For each pair $j,k \in \{1,\ldots,n\}$ with $j \ne k$, and for objects $\ang{x_i}_{i \ne j} \in \txprod_{i \ne j}\C_i$ and $x_k' \in \C_k$, the following diagram commutes in $\D$.
  \begin{equation}\label{lax-0by2}
    %
 
  \end{equation}
\item[Constraint 0-by-0]\index{constraint!- axiom!0-by-0} For each pair $j,k \in \{1,\ldots,n\}$ with $j \ne k$, and for objects $\ang{x_i}_{i \ne j,k} \in \txprod_{i \ne j,k}\C_i$, the following two morphisms are equal:
  \begin{equation}\label{lax-0by0}
    %
 
  \end{equation}
  These are the component of $P^0_j$ with $x_k = \pu$ and the component of $P^0_k$ with $x_j = \pu$.
  
\end{description}
This finishes the definition of a lax $n$-linear functor.  

Moreover, we define the following.
\begin{itemize}
\item A \index{multilinear functor!lax}\index{lax multilinear functor}\emph{lax multilinear functor} is a lax $n$-linear functor for some $n \geq 0$.
\item A lax $n$-linear functor $(P, \{P^2_j\}, \{P^0_j\})$ is
\begin{itemize}
\item \index{lax multilinear functor!strongly unital}\index{multilinear functor!lax!strongly unital}\index{strongly unital!lax multilinear functor}\emph{strongly unital} if each $P^0_j$ is a natural isomorphism,
\item \index{lax multilinear functor!strictly unital}\index{multilinear functor!lax!strictly unital}\index{strictly unital!lax multilinear functor}\emph{strictly unital} if each $P^0_j$ is an identity,
\item \index{lax multilinear functor!strongly monoidal}\index{multilinear functor!lax!strongly monoidal}\index{strongly monoidal!lax multilinear functor}\emph{strongly unital and strongly monoidal} if each $P^0_j$ and each $P^2_j$ is a natural isomorphism,
\item \emph{strictly unital and strong} if each $P^0_j$ is an identity and each $P^2_j$ is a natural isomorphism, and
\item \emph{strict} if each $P^0_j$ and each $P^2_j$ is an identity natural transformation.\dqed
\end{itemize}
\end{itemize}
\end{definition}

\begin{explanation}[Multilinearity and Strict Units]\label{explanation:lax-vs-strict-nlinear}
  In the context of \cref{def:lax-nlinearfunctor} above, note that $(P, \{P^2_j\}, \{P^0_j\})$ is strictly unital if and only if $(P, \{P^2_j\})$ is an $n$-linear functor in the sense of \cref{def:nlinearfunctor}.
  When each $P^0_j$ is an identity, then the lax unity diagrams correspond to the constraint unity axiom \cref{constraintunity} in the case $i = j$, and the constraint 0-by-2 diagrams \cref{lax-0by2} correspond to the constraint unity axiom \cref{constraintunity} in the cases $i \ne j$.
  The other strictly unital variants above can be identified with corresponding variants at the end of \cref{def:nlinearfunctor}.
\end{explanation}

\begin{remark}\label{remark:right-unity-redundant}
  In the context of \cref{def:lax-nlinearfunctor} above, commutativity of each lax unity diagram \cref{lax-constraintunity} follows from the other, by the constraint symmetry axiom \cref{eq:ml-f2-symm} for $P^2_{j}$ together with the equalities
  \[
    \xi_{\pu,?} = 1_? = \xi_{?,\pu}
  \]
  in any permutative category.
\end{remark}

\begin{lemma}\label{lemma:lax-nlinear-tensor}
  In the context of \cref{def:lax-nlinearfunctor}, there is a bijection between multifunctors
  \[
    F \cn \txotimes_{i=1}^n \End(\C_i) \to \End(\D)
  \]
  and lax $n$-linear functors
  \[
    P \cn \txprod_{i=1}^n \C_i \to \D.
  \]
\end{lemma}
\begin{proof}
  Recall from \cref{explanation:unpacking-tensor} that a multifunctor out of a tensor product is characterized by (i) multifunctoriality in each variable separately and (ii) the interchange equality \cref{eq:F-interchange-relation}.
  In the context of \cref{def:lax-nlinearfunctor}, this means that a multifunctor
  \begin{equation}\label{eq:F-otimes-multifunctor}
    F \cn \bigotimes_{i=1}^n \End(\C_i) \to \End(\D)
  \end{equation}
  is characterized by an underlying functor of categories
  \begin{equation}\label{eq:P-prod-nlin}
    P \cn \prod_{i=1}^n \C_i \to \D
  \end{equation}
  such that the following conditions hold.
  \begin{romenumerate}
  \item\label{it:F-nlin-eachvar} In each variable separately, $F$ determines a multifunctor
    \[
      \End(\C_i) \to \End(\D).
    \]
  \item\label{it:F-nlin-interchange} The interchange equality
    \begin{equation}\label{eq:F-nlin-interchage}
      F\big( \txotimes_i \phi_i \big) = F\big( \txotimes^\transp_i \phi_i \big) \cdot \xi^\otimes
    \end{equation}
    holds for multimorphisms $\phi_i \in \End(\C_i)\scmap{\ang{x_i};y_i}$, where $\xi^\otimes$ denotes the bijection \cref{eq:xitimes} induced by transposition of tensor products.
  \end{romenumerate}

  To explain how these conditions relate to the data and axioms of lax multilinearity for $P$, note that each multimorphism
  \[
    \phi_i\in \End(\C_i)\scmap{\ang{x_i};y_i} = \C_i\big( \txoplus_{\ell} x_{i,\ell} \scs y_i\big)
  \]
  with input profile
  \[
    \ang{x_i} = \ang{x_{i,\ell}}_{\ell} = (x_{i,1},\ldots,x_{i,n_i}) \in \Prof(\C_i)
  \]
  is determined uniquely as the composite of a characteristic multimorphism
  \[
    \iota_{\ang{x_i}} \in \End(\C_i)\scmap{\ang{x_i};\txoplus_{\ell} x_{i,\ell}},
  \]
  given by the identity morphism $1_{\oplus_{\ell} x_{i,\ell}}$ in $\C_i$, and a morphism
  \[
    \phi_i\cn \txoplus_{\ell} x_{i,\ell} \to y_i \inspace \C_i.
  \]
  Furthermore, by associativity of the monoidal sum in $\C_i$, each characteristic multimorphism $\iota_{\ang{x_i}}$ decomposes as an iterated composite in $\End(\C_i)$ of multimorphisms $\iota_{\ang{}}$, $\iota_{x_{i,\ell}} = 1_{x_{i,\ell}}$, and $\iota_{(x_{i,\ell}, x_{i,\ell + 1})}$.
  Of those cases, the first two are necessary when $\ang{x_i}$ has length zero or one.

  For each $j \in \{1, \ldots, n\}$, each tuple
  \[
    \ang{x_j} = \ang{x_{j,\ell}}_{\ell} = (x_{j,1},\ldots,x_{j,n_j}) \in \Prof(\C_j),
  \]
  and each $n-1$ tuple of objects
  \[
    x_i \in \C_i \withspace i \ne j,
  \]
  let
  \[
    \iota^j_{\ang{x_j}} =
    1_{x_1} \otimes \cdots \otimes 1_{x_{j-1}}
    \otimes \iota_{\ang{x_j}} \otimes
    1_{x_{j+1}} \otimes \cdots \otimes 1_{x_n}.
  \]
  The two cases of interest are where $\ang{x_j}$ has length zero or two.
  For a pair of objects $x_j,x_j' \in \C_j$, define the following:
  \begin{align*}
    P^0_j & = F\iota^j_{\ang{}} \andspace \\
    P^2_j & = F\iota^j_{(x_j,x_j')}.
  \end{align*}
  With this notation, the reasoning above shows that the assignment of $F$ on multimorphisms $\otimes_i \phi_i$ is determined by
  \begin{itemize}
  \item its underlying functor $P$,
  \item the multimorphisms $P^0_j$ and $P^2_j$, and
  \item multifunctoriality of $F$.
  \end{itemize}
  The remainder of this proof identifies how the multifunctoriality of $F$ corresponds to the six multilinearity axioms of \cref{def:lax-nlinearfunctor}.

  Condition \cref{it:F-nlin-eachvar}, multifunctoriality in each variable separately, is equivalent to the condition that the data $\{P^0_j\}$ and $\{P^2_j\}$ are natural with respect to the variables $x_i$ for $i \neq j$ and satisfy the axioms for lax unity \cref{lax-constraintunity}, constraint associativity \cref{eq:ml-f2-assoc}, and constraint symmetry \cref{eq:ml-f2-symm}.

  For condition \cref{it:F-nlin-interchange},  the interchange equality \cref{eq:F-nlin-interchage} holds in general if and only if it holds when, for each pair $j \ne k$ in $\{1,\ldots,n\}$, the multimorphisms $\phi_j$ and $\phi_k$ are characteristic multimorphisms:
  \begin{align*}
    \phi_j = \iota_{\ang{}} &  \orspace \phi_j = \iota_{(x_j,x_j')}, 
                              \\
    \phi_k = \iota_{\ang{}} &  \orspace \phi_k = \iota_{(x_k,x_k')},
  \end{align*}
  and each other $\phi_i$ is a colored unit.
  By symmetry in $j$ and $k$, the reduction above results in three distinct cases.
  The interchange equality \cref{it:F-nlin-interchange} for those three cases corresponds to the constraint axioms 2-by-2 \cref{eq:f2-2by2}, 0-by-2 \cref{lax-0by2}, and 0-by-0 \cref{lax-0by0} for $P$.

  Thus, a multifunctor $F$ in \cref{eq:F-otimes-multifunctor} determines and is uniquely determined by an underlying functor $P$ in \cref{eq:P-prod-nlin} with linearity and unit constraints
  \[
    P^2_j = F\iota^j_{(x_j,x_j')} \andspace
    P^0_j = F\iota^j_{\ang{}}
  \]
  satisfying the six multilinearity axioms of \cref{def:lax-nlinearfunctor}.
  This completes the proof.
\end{proof}

The following generalizes \cref{def:nlineartransformation}.
\begin{definition}\label{def:lax-nlineartransformation}
Suppose $P,Q$ are lax $n$-linear functors as displayed below.
\begin{equation}\label{lax-nlineartransformation}
\begin{tikzpicture}[xscale=3,yscale=2.5,baseline={(a.base)}]
\draw[0cell=.9]
(0,0) node (a) {\prod_{j=1}^n \C_j}
(a)++(1,0) node (b) {\D}
;
\draw[1cell=.9]  
(a) edge[bend left=25] node[pos=.45] {\big( P, \{P^2_j\}, \{P^0_j\}\big)} (b)
(a) edge[bend right=25] node[swap,pos=.45] {\big( Q, \{Q^2_j\}, \{Q^0_j\}\big)} (b)
;
\draw[2cell] 
node[between=a and b at .47, shift={(0,0)}, rotate=-90, 2label={above,\theta}] {\Rightarrow}
;
\end{tikzpicture}
\end{equation}
A \index{n-linear@$n$-linear!transformation!lax}\emph{lax $n$-linear transformation} $\theta \cn P \to Q$ is a natural transformation of underlying functors that satisfies the constraint compatibility condition \cref{eq:monoidal-in-each-variable} together with the following lax unity axiom.
\begin{description}
\item[Lax Unity] 
\begin{equation}\label{lax-niitransformationunity}
  \begin{tikzpicture}[x=32mm,y=15mm,vcenter]
    \tikzset{0cell/.append style={nodes={scale=.9}}}
    \tikzset{1cell/.append style={nodes={scale=.9}}}
    \draw[0cell] 
    (0,-.5) node (a) {
      \pu
    }
    (1,0) node (b) {
      P\big( \ang{x_i}_{i \ne j} \insj \pu \big)
    }
    (1,-1) node (c) {
      Q\big( \ang{x_i}_{i \ne j} \insj \pu \big)
    }
    ;
    \draw[1cell] 
    (a) edge node {P^0_j} (b)
    (a) edge['] node {Q^0_j} (c)
    (b) edge node {\theta} (c)
    ;
  \end{tikzpicture}
\end{equation}
\end{description}
This finishes the definition of a lax $n$-linear transformation.
Moreover, we define the following.
\begin{itemize}
\item A \index{lax multilinear transformation}\index{multilinear transformation!lax}\index{transformation!multilinear!lax}\emph{lax multilinear transformation} is a lax $n$-linear transformation for some $n \geq 0$. 
\item Identities and compositions of lax multilinear transformations are defined componentwise.
\end{itemize}
Note that a lax $n$-linear transformation between strictly unital lax $n$-linear functors is the same as an $n$-linear transformation in the sense of \cref{def:nlineartransformation}.
\end{definition}

\begin{definition}[Multimorphism Categories]\label{definition:permcat-homcat}
  We define the following categories of lax $n$-linear functors and transformations.
  \begin{itemize}
  \item $\permcat\scmap{\ang{\C};\D}$ is the category with
    \begin{itemize}
    \item lax $n$-linear functors $\ang{\C} \to \D$ as objects and
    \item lax $n$-linear transformations between them as morphisms.
    \end{itemize} 
  \item $\permcatsg\scmap{\ang{\C};\D}$ is the full subcategory of strongly unital and strongly monoidal $n$-linear functors.\defmark
  \end{itemize}
  Combined with the multimorphism categories of \cref{definition:permcatsus-homcat}, these fit into the following commutative diagram of full subcategory inclusions.
  \begin{equation}\label{eq:lax-multimorphismcats}

  \end{equation} 
\end{definition}

The following generalizes \cref{definition:permcat-action}.
\begin{definition}[Symmetric Group Action]\label{definition:lax-permcat-action}
  Suppose given lax $n$-linear functors $P$ and $Q$ together with a lax $n$-linear transformation $\theta$ as displayed below.
\begin{equation}\label{lax-permiicatcddata}
%
\end{equation}
For a permutation $\sigma \in \Sigma_n$, the symmetric group action 
\begin{equation}\label{lax-permiicatsymgroupaction}
\begin{tikzcd}[column sep=large]
\permcat\mmap{\D; \ang{\C}} \ar{r}{\sigma}[swap]{\iso} & \permcat\mmap{\D; \ang{\C}\sigma}
\end{tikzcd}
\end{equation}
sends the data \cref{lax-permiicatcddata} to the following composites and whiskerings, where $\si$ permutes the coordinates according to $\sigma$.
\begin{equation}\label{lax-permiicatsigmaaction}
%
\end{equation}
For objects 
\[\ang{a} = \ang{a_j}_{j=1}^n \in \prod_{j=1}^n \C_{\sigma(j)} \andspace a_j' \in \C_{\sigma(j)},\] 
the $j$-th linearity constraint of $P^\sigma = P \circ \sigma$ has component given as in \cref{fsigmatwoj}:
\[
  (P^\si)^2_j = P^2_{\si(j)}.
\]
The $j$-th unit constraint is likewise given by
\[
  (P^\si)^0_j = P^0_{\si(j)}.
\]
As with \cref{definition:permcat-action}, each of the variants in \cref{def:lax-nlinearfunctor} is preserved by this action.
\end{definition}

The following generalizes \cref{definition:permcat-comp}.
\begin{definition}[Multicategorical Composition]\label{definition:lax-permcat-comp}
Suppose given, for each $j \in \{1,\ldots,n\}$,
\begin{itemize}
\item permutative categories $\ang{\B_j} = \ang{\B_{j,i}}_{i=1}^{k_j}$,
\item lax $k_j$-linear functors $P'_j, Q'_j \cn \ang{\B_j} \to \C_j$, and
\item a lax $k_j$-linear transformation $\theta_j \cn P'_j \to Q'_j$ as follows.
\end{itemize} 
\begin{equation}\label{lax-permiicatbcdata}
    \begin{tikzpicture}[xscale=3,yscale=3,vcenter]
      \draw[0cell=.9]
      (0,0) node (a) {\prod_{i=1}^{k_j} \B_{j,i}}
      (a)++(1,0) node (b) {\C_j}
      ;
      \draw[1cell=.9]  
      (a) edge[bend left=20] node[pos=.45] {P'_j} (b)
      (a) edge[bend right=20] node[swap,pos=.45] {Q'_j} (b)
      ;
      \draw[2cell] 
      node[between=a and b at .5, shift={(0,0)}, rotate=-90, 2label={below,\theta_j}] {\Rightarrow}
      ;
    \end{tikzpicture}
  \end{equation}
With $\ang{\B} = \bang{\ang{\B_j}}_{j=1}^n$, the multicategorical composition functor
\begin{equation}\label{lax-permiicatgamma}
    \begin{tikzcd}[cells={nodes={scale=.85}}]
      \permcat\mmap{\D;\ang{\C}} \times \txprod_{j=1}^n \permcat\mmap{\C_j;\ang{\B_j}} \ar{r}{\gamma}
      & \permcat\mmap{\D;\ang{\B}}
    \end{tikzcd}
  \end{equation}
  sends the data \cref{lax-permiicatcddata,lax-permiicatbcdata} to the composites
  \begin{equation}\label{lax-compositemodification}
    \begin{tikzpicture}[xscale=5,yscale=4,vcenter]
      \draw[0cell=.9]
      (0,0) node (a) {\prod_{j=1}^n \prod_{i=1}^{k_j} \B_{j,i}}
      (a)++(1,0) node (b) {\D}
      ;
      \draw[1cell=.9]  
      (a) edge[bend left=20] node[pos=.45] {P \circ \txprod_j P'_j} (b)
      (a) edge[bend right=20] node[swap,pos=.45] {Q \circ \txprod_j Q'_j} (b)
      ;
      \draw[2cell=.9] 
      node[between=a and b at .6, shift={(0,0)}, rotate=-90, 2label={below,\theta \otimes (\textstyle\prod_j \theta_j)}] {\Longrightarrow}
      ;
    \end{tikzpicture}
  \end{equation}
  defined as follows.

  \begin{description}
  \item[Composite Lax Multilinear Functor]
    The underlying functor $\big(P \comp \prod_j P_j'\big)$ and linearity constraints $\big(P \comp \prod_j P_j'\big)^2_\ell$ are defined as in \cref{definition:permcat-comp}.
    The unit constraints $\big(P \comp \prod_j P_j'\big)^0_\ell$ are defined as follows.
    
    Suppose
    \[
      \ell = k_1 + \cdots + k_{a-1} + b,
    \]
    for some $a \in \{1,\ldots,n\}$ and $b \in \{1,\ldots,k_a\}$,
    and suppose given
    \begin{equation}\label{lax-angwj}
      \begin{alignedat}{2}
        \ang{w_j}
        &= \ang{w_{j,i}}_{i=1}^{k_j} \in \txprod_{i=1}^{k_j} \B_{j,i} 
        &\quad & \forspace j \in \{1,\ldots,n\} \withspace j \ne a \andspace \\
        w_{a,i}
        & \in \B_{a,i}
        &\quad & \forspace i \in \{ 1,\ldots,k_a \} \withspace i \ne b.
      \end{alignedat}
    \end{equation} 
    Then there is the tuple
    \begin{equation}\label{lax-compositefangw}
        \textstyle
        \bang{P_j' \ang{w_j}}_{j \ne a}
        = \big(P_1'\ang{w_1},\ldots, P_{a-1}'\ang{w_{a-1}}, P_{a+1}'\ang{w_{a+1}}, \ldots, P_n'\ang{w_n}\big) \inspace \prod_{j \ne a}\C_j.
    \end{equation}
    The unit constraint $\big(P \comp \prod_j P_j'\big)^0_\ell$ is the following composite in $\D$:
    \begin{equation}\label{eq:laxunit-composite}
      \begin{tikzpicture}[x=20ex,y=10ex]
        \draw[0cell] 
        (0,0) node (a) {\pu}
        (a)+(.8,0) node (b) {
          P\big(\bang{P'_jw_j}_{j \ne a} \insa \pu\big)
        }
        (b)+(2,0) node (c) {
          P\Big(
          \bang{P'_jw_j}_{j \ne a} \insa
          P'_a\big( \ang{w_{a,i}}_{i \ne b} \insb \pu \big)
          \Big).
        }
        ;
        \draw[1cell] 
        (a) edge node {P^0_a} (b)
        (b) edge node {P\big( \bang{1 \insa (P_a')^0_b} \big)} (c)
        ;
      \end{tikzpicture}
    \end{equation}
    If the monoidal constraints $P^0_a$ and $(P'_a)^0_b$ are isomorphisms, respectively identities, then so is the composite \cref{eq:laxunit-composite}.
    Therefore, as with \cref{definition:permcat-comp}, each of the variants in \cref{def:lax-nlinearfunctor} is preserved by composition.

  \item[Composite Lax Multinatural Transformation]
    The lax $n$-linear transformation $\theta \otimes \big(\prod_j \theta_j\big)$ in \cref{lax-compositemodification} is the horizontal composite of the natural transformations $\prod_j \theta_j$ and $\theta$, as in \cref{thetaprodthetaw}.
  \end{description}
  The finishes the definition of the multicategorical composition in $\permcat$.

  Verifying the lax multilinearity axioms for the composite $P \comp \prod_j P'_j$ is similar to that of the strictly unital case in \cite[Section 6.6]{cerberusIII}.
  In particular, the axioms \cref{lax-constraintunity,lax-0by2,lax-0by0} hold by the corresponding axioms for $P$ and $P'_j$, along with naturality and multifunctoriality of the data involved.
  For example, in the constraint 0-by-2 axiom \cref{lax-0by2} for $\big(P \comp \prod_j P_j'\big)^0_\ell$ and $\big(P \comp \prod_j P_j'\big)^2_m$ with
  \begin{align*}
    \ell & = k_1 + \cdots + k_{a-1} + b 
    & \forspace & a \in \{1,\ldots,n\} \quad \text{and} \quad b \in \{1,\ldots,k_a\}
    \\
    m & = k_1 + \cdots + k_{c-1} + d 
    & \forspace &  c \in \{1,\ldots,n\} \quad \text{and} \quad d \in \{1,\ldots,k_c\},
  \end{align*}
  there are two cases to consider.
  If $c \ne a$, then verifying this case uses \cref{lax-0by2} for $P^0_a$ along with naturality of $P^2_c$, naturality of $P^0_a$, and multifunctoriality of $P$.
  If $c = a$ and $b \ne d$, then verifying this case uses \cref{lax-constraintunity} for $P^0_a$ and \cref{lax-0by2} for $(P'_a)^0_b$ and $(P'_a)^2_d$, along with naturality of $P^2_a$ and multifunctoriality of $P$.
  
  Similarly, the lax unity axiom \cref{lax-niitransformationunity} for $\theta \otimes (\txprod_j \theta_j)$ follows from that of $\theta$ and $\theta_j$ individually.
\end{definition}

\begin{proposition}\label{proposition:laxmultilin-tensor}
  For small permutative categories $\ang{\C_j}_{j=1}^n$ and $\D$, the 2-functor 
  \[\End \cn \permcat \to \Multicat\]
  in \cref{endtwofunctor} induces an isomorphism of categories
  \begin{align*}\label{eq:laxmultilin-via-otimes}
    \permcat\scmap{\ang{\C};\D} \fto[\iso]{\End}
    & \Multicat\scmap{\ang{\End(\C)};\End(\D)}\\
    & = \Multicat\big(\txotimes_{i=1}^n \End(\C_i) \scs \End(\D)\big)
  \end{align*}
  between
  \begin{itemize}
  \item the category of lax $n$-linear functors and transformations $\ang{\C} \to \D$ and
  \item the category of multifunctors 
    \[\txotimes_{i=1}^n \End(\C_i) \to \End(\D)\]
    and multinatural transformations.
  \end{itemize} 
\end{proposition}
\begin{proof}
  This proof is an unpointed analog of \cref{expl:endst-catmulti}.
  The bijection between lax multilinear functors from $\txprod_i \C_i$ to $\D$ and multifunctors from $\txotimes_i \End(\C_i)$ to $\End(\D)$ is given by \cref{lemma:lax-nlinear-tensor}.

  The argument that $\End$ induces a bijection between
  lax $n$-linear transformations (\cref{def:lax-nlineartransformation})
  \[\theta\cn P \to Q \inspace \permcat\scmap{\ang{\C};\D}\]
  and multinatural transformations (\cref{def:enr-multicat-natural-transformation})
  \[\omega\cn \End(P) \to \End(Q) \inspace \Multicat\big(\txotimes_{i=1}^n \End(\C_i) \scs \End(\D)\big)\]
  is similar.
  Indeed, for
  \[\theta \cn P \to Q \inspace \permcat\scmap{\ang{\C};\D},\]
  the multinatural transformation
  \[\End(\theta) \cn \End(P) \to \End(Q)\]
  is uniquely determined by the component morphisms
  \begin{equation}\label{eq:End-theta-component}
    \theta_{\angc} \cn P\angc \to Q\angc \inspace \D,
  \end{equation}
  for $\angc \in \prod_{i=1}^n \C_i$.
  Using the same reduction to characteristic multimorphisms in the proof of \cref{lemma:lax-nlinear-tensor},
  the multinaturality axioms of $\omega = \End(\theta)$ determine and are uniquely determined by the two lax multilinearity axioms of $\theta$.
\end{proof}

The following generalizes \cref{thm:permcatmulticat}.
\begin{theorem}\label{theorem:permcat-catmulticat}
  There is a $\Cat$-multicategory
  \[
    \permcat
  \]
  defined by the following data.
  \begin{itemize}
  \item The objects are small permutative categories.
  \item The multimorphism categories are in \cref{definition:permcat-homcat}.
  \item The colored units are identity symmetric monoidal functors.
  \item The symmetric group action is in \cref{definition:lax-permcat-action}.
  \item The multicategorical composition is in \cref{definition:lax-permcat-comp}.
  \end{itemize}
  Moreover, the multimorphism categories \cref{eq:lax-multimorphismcats} give the following commutative diagram of sub-$\Cat$-multicategories.
  \begin{equation}\label{eq:perm-catmulticats}
    \begin{tikzpicture}[x=16ex,y=4ex,vcenter]
      \draw[0cell] 
      (0,0) node (a) {\permcatst}
      (a)+(1.2,1) node (b) {\permcatsg}
      (a)+(1.2,-1) node (sub) {\permcatsus}
      (b)+(1.1,0) node (c) {\permcat}
      (sub)+(1.1,0) node (suc) {\permcatsu}
      ;
      \draw[1cell] 
      (a) edge[right hook->] node {} (b)
      (a) edge[right hook->] node {} (sub)
      (b) edge[right hook->] node {} (c)
      (sub) edge[right hook->] node {} (suc)
      (sub) edge[right hook->] node {} (b)
      (suc) edge[right hook->] node {} (c)
      ;
    \end{tikzpicture}
  \end{equation} 
\end{theorem}

The remainder of this section describes the closed structure for $\permcat$, $\permcatsg$, and $\permcatsus$.

\subsection*{Internal Hom Objects}

The following generalizes the monoidal sum of multilinear functors.
\begin{explanation}[Monoidal Sum of Lax Multilinear Functors]\label{explanation:lax-sum}
  Suppose given small permutative categories $\D$ and $\ang{\C} = \ang{\C_i}_{i=1}^n$ for $n \ge 0$.
  The monoidal sum on $\permcatsu\scmap{\ang{\C};\D}$, from \cref{def:clp-angcd}, generalizes to a monoidal sum on $\permcat\scmap{\ang{\C};\D}$.

  To explain this, first note that 
  \[
    \oplus \cn \D \times \D \to \D
  \]
  is a strictly unital strong symmetric monoidal functor.
  Details for this appear in the proof of \cref{psucd-hom-permcat} statement \cref{psucd-hom-i}, taking $F = G = 1_{\D}$.
  The monoidal constraint of $\oplus$ is given by
  \begin{equation}\label{eq:1xi1}
    (a \oplus b) \oplus (c \oplus d) \fto{1 \oplus \xi \oplus 1}
    (a \oplus c) \oplus (b \oplus d)
  \end{equation} 
  for each quadruple $a, b, c, d \in \D$.
  The associativity and braiding axioms for $\oplus$ both follow from coherence for symmetric monoidal categories \cite[XI.1 Theorem 1]{maclane}.

  For multilinear functors $P$ and $Q$ in $\permcatsu\scmap{\ang{\C};\D}$, the underlying functor and linearity constraints of $P \oplus Q$ are given by those of the composite $\oplus \circ (P \times Q)$.
  This same definition applies more generally to lax multilinear functors $P$ and $Q$ in $\permcat\scmap{\ang{\C};\D}$.
  Since $\oplus$ is strictly unital, the unit constraints of the sum are given by
  \begin{equation}\label{eq:PQ0-P0Q0}
    (P \oplus Q)^0_j = P^0_j \oplus Q^0_j. 
  \end{equation} 
  The linearity constraints are given as in \cref{PplusQtwoi}.

  Lax multilinearity of $P \oplus Q$ is a special case of general composition of lax multilinear functors, but can also be verified directly.
  The proof of \cref{clp-angcd-permutative} verifies the axioms \cref{eq:ml-f2-assoc,eq:ml-f2-symm,eq:f2-2by2} for $P \oplus Q$.
  The remaining axioms, \cref{lax-constraintunity,lax-0by2,lax-0by0}, are similar.
  In each case, one makes use of the corresponding axioms for $P$ and $Q$, together with strictness of the unit $\pu$.
  For \cref{lax-constraintunity,lax-0by2}, one also uses naturality of $\xi$ and the identities
  \[
    \xi_{?,\pu} = 1_{?} = \xi_{\pu,?}
  \]
  that hold in any permutative category.
\end{explanation}

The monoidal sum from \cref{explanation:lax-sum} gives the following generalization of \cref{def:clp-angcd}.
\begin{definition}\label{definition:laxinternalhom}
  For small permutative categories $\D$ and $\ang{\C} = \ang{\C_i}_{i=1}^n$ for $n \ge 0$, the \index{permutative category!internal hom -}\index{internal hom!permutative category}\emph{internal hom permutative category}
  \begin{equation}\label{eq:laxinternalhom}
    \Big(\clpermcat\scmap{\angC;\D} \scs \oplus \scs \clpu \scs \clxi \Big)
  \end{equation}
  is given by the underlying categories
  \[
    \permcat\scmap{\angC;\D}.
  \]
  The monoidal sum $\oplus$ is given on objects (lax multilinear functors) by composition with the strictly unital strong symmetric monoidal functor
  \begin{equation}\label{eq:oplusD}
    \oplus \cn \D \times \D \to \D
  \end{equation} 
  as described in \cref{explanation:lax-sum}.
  The monoidal sum on morphisms (lax multilinear transformations) is given by whiskering with $\oplus$.
  The monoidal unit $\clpu$ and braiding $\clxi$ from \cref{def:clp-angcd} also provide a monoidal unit and braiding for $\permcat\scmap{\ang{\C};\D}$.
  Verification that this defines a permutative category \cref{eq:laxinternalhom} is similar to the proof of \cref{clp-angcd-permutative}, extended to the case of not-necessarily-strict unit constraints.

  This finishes the definition of the permutative structure on $\clpermcat\scmap{\angC;\D}$.
  Moreover, because $\oplus$ \cref{eq:oplusD} is a strictly unital strong symmetric monoidal functor, composition with $\oplus$ also defines permutative structures in the strong and strictly unital strong cases.
  These are denoted
  \begin{equation}\label{eq:lax-sgsus-internalhom}
    \Big(\clpermcatsg\scmap{\angC;\D} \scs \oplus \scs \clpu \scs \clxi \Big)
    \andspace
    \Big(\clpermcatsus\scmap{\angC;\D} \scs \oplus \scs \clpu \scs \clxi \Big)
  \end{equation}
  respectively.
\end{definition}

\begin{remark}\label{remark:strict-notinternalhom}
  The permutative structure described above does not specialize to strict multilinear functors
  \[
    P,Q \in \permcatst\scmap{\ang{\C};\D}.
  \]
  This is because, when $P$ and $Q$ are strict, the monoidal sum $P \oplus Q$ generally has nontrivial monoidal constraint determined by that of $\oplus\cn \D \times \D \to \D$ \cref{eq:1xi1}.
  As noted in \cref{explanation:lax-sum}, the latter is an isomorphism---generally not an identity---determined by the symmetry isomorphism of $\D$.
\end{remark}

\subsection*{Symmetric Group Action on Internal Hom}

The following generalizes \cref{def:clp-ancd-sigma}.

\begin{definition}\label{def:lax-clp-ancd-sigma}\index{symmetric group!action}
  For small permutative categories $\D$ and $\angC = \ang{\C_i}_{i=1}^n$ for $n \geq 0$ and a permutation $\sigma \in \Sigma_n$, we define the functor
  \begin{equation}\label{lax-clp-sigma-action}
    \clpermcat\scmap{\angC;\D} \fto[\iso]{\sigma} \clpermcat\scmap{\angC\sigma;\D}
  \end{equation}
  as the isomorphism of underlying categories
  \[\permcat\scmap{\angC;\D} \fto[\iso]{\sigma} \permcat\scmap{\angC\sigma;\D}\]
  that is given, as in \cref{permiicatsymgroupaction}, by precomposition and whiskering with the permutation of factors
  \[
    \txprod_{i=1}^n \C_{\si(i)} \fto{\si} \txprod_{i=1}^n \C_i.\dqed
  \]
  An argument similar to that of \cref{clp-angcd-strict} shows that \cref{lax-clp-sigma-action} is a strict symmetric monoidal isomorphism and that the equivariance axioms for internal hom objects, \cref{right-id-action,hom-equivariance}, hold.

  This describes the symmetric group action for $\clpermcat$.
  Moreover, the same action induces a symmetric group action in the strong and strictly unital strong cases \cref{eq:lax-sgsus-internalhom} above.
  The same arguments verify the equivariance axioms in these cases.
\end{definition}

\subsection*{Multicategorical Evaluation}

The following generalizes \cref{def:clp-eval}.
\begin{definition}\label{def:lax-clp-eval}\index{multilinear!evaluation functor}\index{evaluation!multilinear functor}\index{permutative category!multilinear evaluation}
  For small permutative categories $\D$ and $\angC = \ang{\C_i}_{i=1}^n$ with $n \geq 0$, we define the data of a lax $(n+1)$-linear functor
  \begin{equation}\label{lax-clp-ev-angcd}
    \clpermcat\scmap{\angC;\D} \times \txprod_{i=1}^n \C_i \fto{\ev_{\angC;\,\D}} \D
  \end{equation}
  as follows.
  \begin{description}
  \item[Underlying Functor] The underlying functor is that of  \cref{evangcd-object,evangcd-thetaangf}:
    \[
      \ev_{\angC;\,\D} \big(P,\angx\big) = P\angx
      \andspace
      \ev_{\angC;\,\D} \big(\theta,\angf\big)
      = Q(f) \circ \theta_{\ang{x}} = \theta_{\ang{y}} \circ P(f).
    \]
  \item[Linearity Constraints] The linearity constraints are given by \cref{evangcd-first-linearity,evangcd-iplusone-linearity}:
    \[
      (\ev_{\angCD})^2_1 = 1
      \andspace
      (\ev_{\angCD})^2_{i+1} = P^2_i
      \forspace i \in \{1,\ldots,n\}.
    \]
  \item[Unit Constraints] The unit constraints are given by
    \[
      (\ev_{\angCD})^0_1 = 1
      \andspace
      (\ev_{\angCD})^0_{i+1} = P^0_i
      \forspace i \in \{1,\ldots,n\}.
    \]
  \end{description}
  This finishes the definition of $\ev_{\angC;\, \D}$.
  Verification that these data define a lax $(n+1)$-linear functor is similar to the proof of \cref{clp-evaluation}.
  The axioms \cref{lax-constraintunity,lax-0by2,lax-0by0} for lax unity, constraint 0-by-2, and constraint 0-by-0, hold by those for $P$ in the cases $j > 0$ and $k > 0$.
  For the cases $j=0$ or $k=0$, these axioms hold because the constant functor $\clpu$ is strictly monoidal and because $(P \oplus Q)^0_j = P^0_j \oplus Q^0_j$ in \cref{eq:PQ0-P0Q0}.

  To verify the evaluation bijection axiom \cref{eval-bijection}, \cref{lax-clp-partners} below describes how the definitions of $\chi$ and $\Psi$ in \cref{expl:clp-partners,clp-evbij-axiom} generalize to the lax multilinear case.
  The equivariance axiom \cref{eval-bij-eq} for evaluation bijection follows from the same argument in the proof of \cref{clp-evbij-equiv}.
  Verification that the corresponding unit constraints are equal uses the same permutation of indices as in \cref{eq:evbij-lin-r,eq:evbij-lin-n+t}, with the definitions \cref{eq:chiP0i,eq:chiP0n+j} from \cref{eval-bijection} below.

  Because the unit and linearity constraints of $\ev$ depend on those of its first argument, $P$, the same definition specializes to the strong and strictly unital strong cases \cref{eq:lax-sgsus-internalhom} above.
  The same arguments verify the evaluation bijection and equivariance axioms in these cases too.
\end{definition}

The following generalizes the functions $\chi$ and $\Psi$ from \cref{expl:clp-partners,clp-evbij-axiom}, respectively.
\begin{explanation}[The Functions $\chi$ and $\Psi$]\label{lax-clp-partners}
  For small permutative categories $\B$, $\ang{\C} = \ang{\C_i}_{i=1}^n$, and $\ang{\D} = \ang{\D_j}_{j=1}^p$, there are inverse functions
  \begin{equation}\label{lax-chiPsi}

  \end{equation}
  generalizing those of \cref{expl:clp-partners,clp-evbij-axiom}.

  For
  \[
    P = \big( P , \{ P^2_i \}, \{ P^0_i \} \big) \inspace
    \permcat\left( \angC \sscs \clpermcat\pangDB\right),
  \]
  $\chi P$ has linearity constraints given by \cref{chiPtwoi,chiPtwonplusj}.
  The unit constraints of $\chi P$ are given by corresponding components of $P^0_i$ and $(P\ang{w})^0_{n+j}$ as shown in the following displays, with $\ang{w} \in \txprod_k \C_k$, $\ang{y} \in \txprod_k \D_k$, $i \in \{ 1, \ldots, n \}$, and $j \in \{ 1, \ldots, p \}$.
  \begin{equation}\label{eq:chiP0i}
    %
  \end{equation}

  Each of the first three lax multilinearity axioms for $\chi P$, \cref{eq:ml-f2-assoc,eq:ml-f2-symm,eq:f2-2by2}, follows as in \cref{expl:clp-partners}.
  Verifying the other lax multilinearity axioms,
  \cref{lax-constraintunity,lax-0by2,lax-0by0}, is similar, using the corresponding axioms for $\chi P$.
  Verification of the axioms involving indices $i \in \{1,\ldots,n\}$ and $n+j \in \{n+1,\ldots,n+p\}$ uses the lax multinaturality axioms \cref{eq:monoidal-in-each-variable,lax-niitransformationunity} for $P^0_i$ and $P^2_i$.

  For
  \[
    R = \big( R, \{ R^2_r \}, \{ R^0_r \} \big) \inspace
    \permcat\big( \angC,\angD \sscs \angB \big),
  \]
  and for each $\ang{w} \in \txprod_k \C_k$,
  the linearity constraints of $\Psi R \ang{w}$ and $\Psi R$ are given by \cref{Rtwonplusj,PsiRtwoi-component}, respectively.
  The unit constraints of $\Psi R \ang{w}$ and $\Psi R$ are determined by corresponding components of $R^0$ as shown in the following displays, with $\ang{y}$, $i$, and $j$ as above.
  \begin{equation}\label{eq:PsiRw0j}
    \begin{tikzpicture}[x=35ex,y=6ex,vcenter]
      \draw[0cell=.9] 
      (0,0) node (a) {\pu}
      (0,-1) node (a') {\pu}
      (a)+(1,0) node (b) {
        (\Psi R)\ang{w}(\ang{y_k}_{k \ne j} \insj \pu)
      }
      (a')+(1,0) node (b') {
        R\big(\ang{w},(\ang{y_k}_{k \ne j} \insj \pu)\big)
      }
      ;
      \draw[1cell=.9] 
      (b) edge[equal] node {} (b')
      (a) edge[equal] node {} (a')
      (a) edge node {\big( (\Psi R)\ang{w}^0_j\big)_{\ang{y_k}_{k \ne j}}} (b)
      (a') edge node {(R^0_{n+j})_{\ang{w},\ang{y_k}_{k \ne j}}} (b')
      ;
    \end{tikzpicture}
  \end{equation}
  \begin{equation}\label{eq:PsiR0i}
    \begin{tikzpicture}[x=35ex,y=6ex,vcenter]
      \draw[0cell=.9] 
      (0,0) node (a) {\pu}
      (0,-1) node (a') {\pu}
      (a)+(1,0) node (b) {
        (\Psi R)(\ang{w_k}_{k \ne i} \insi \pu)\ang{y}
      }
      (a')+(1,0) node (b') {
        R\big((\ang{w_k}_{k \ne i} \insi \pu),\ang{y}\big)
      }
      ;
      \draw[1cell=.9] 
      (b) edge[equal] node {} (b')
      (a) edge[equal] node {} (a')
      (a) edge node {\Big( \big( (\Psi R)^0_i\big)_{\ang{w_k}_{k \ne i}} \Big)_{\ang{y}}} (b)
      (a') edge node {(R^0_i)_{\ang{w_k}_{k \ne i},\ang{y}}} (b')
      ;
    \end{tikzpicture}
  \end{equation}

  Verification of the necessary axioms follows the same structure as in the proof of \cref{clp-evbij-axiom}.
  For $\ang{w} \in \txprod_k \C_k$, each of the lax $p$-multilinearity axioms for $\Psi R \ang{w}$, $(\Psi R \ang{w})^0_j$, and $(\Psi R \ang{w})^2_j$ with $j \in \{1, \ldots, p\}$, \cref{eq:ml-f2-assoc,eq:ml-f2-symm,eq:f2-2by2,lax-constraintunity,lax-0by2,lax-0by0}, follows from the corresponding axiom for $R$, $R^0_{n+j}$, and $R^2_{n+j}$.
  For $i \in \{1,\ldots,n\}$, the lax multinaturality axioms for $(\Psi R)^0_i$ and $(\Psi R)^2_i$ hold by the 2-by-2, 0-by-2, and 0-by-0 axioms for $R$, with indices $i$ and $n+j \in \{n+1,\ldots,n+p\}$.
  The remaining lax multilinearity axioms for $\Psi R$ then follow from the corresponding axioms for $R$.

  Lastly, it remains to verify that $\chi$ and $\Psi$ are inverse functions.
  This follows the same argument given in the proof of \cref{clp-evbij-axiom}, using \cref{eq:chiP0i,eq:chiP0n+j,eq:PsiRw0j,eq:PsiR0i} along with \cref{evangcd-first-linearity,evangcd-iplusone-linearity,chiPtwoi,chiPtwonplusj}.
\end{explanation}

With these extensions to lax multilinear functors, the following result generalizes \cref{perm-closed-multicat}.
\begin{theorem}\label{lax-perm-closed-multicat}\index{closed!multicategory!- of permutative categories}\index{multicategory!closed - of permutative categories}\index{permutative category!closed multicategory of}
  There is a closed multicategory
  \[
    \big(\permcat \scs \clpermcat \scs \ev \big)
  \]
  consisting of the following data.
  \begin{itemize}
  \item The underlying multicategory is $\permcat$ in \cref{theorem:permcat-catmulticat}.
  \item The internal hom objects are the permutative categories $\clpermcat\scmap{\angC;\D}$ in \cref{definition:laxinternalhom}.
  \item The symmetric group action on internal hom is given in \cref{def:lax-clp-ancd-sigma}.
  \item The multicategorical evaluations are the lax multilinear functors
    \[\clpermcat\scmap{\angC;\D} \times \txprod_{i=1}^n \C_i \fto{\ev_{\angC;\,\D}} \D\]
    in \cref{def:lax-clp-eval}.
  \end{itemize}
  Furthermore, the data above specialize to define closed multicategory structures
  \[
    \big(\permcatsg \scs \clpermcatsg \scs \ev \big)
    \andspace
    \big(\permcatsus \scs \clpermcatsus \scs \ev \big).
  \]
\end{theorem}

\begin{remark}\label{remark:strict-notclosed}
  Recall from \cref{remark:strict-notinternalhom} that the monoidal sum of strict multilinear functors is generally not strict, and therefore the internal hom structures in \cref{definition:laxinternalhom} do not specialize to $\permcatst$.
  The same obstruction (generally nontrivial symmetry in the target, $\D$) also prevents \cref{lax-perm-closed-multicat} from specializing to $\permcatst$.
\end{remark}

\chapter[Self-Enrichment and Standard Enrichment]{Self-Enrichment and Standard Enrichment of Closed Multicategories}
\label{ch:std_enrich}
This chapter develops the definitions and basic theory of
\begin{itemize}
\item self-enrichment for closed multicategories and
\item standard enrichment for multifunctors between closed multicategories.
\end{itemize}
For a non-symmetric multifunctor
\[
  \big(\M,\clM,\evM\big) \fto{F} \big(\N,\clN,\evN\big),
\]
\cref{def:clmulti-clcat,cl-multi-cl-cat} describe corresponding self-enriched categories for $\M$ and $\N$.
\cref{def:gspectra-thm-iii-context,gspectra-thm-iii} describe the induced $\N$-functor
\[
  \Fse \cn \M_F \to \N,
\]
where $\M_F$ is the $\N$-category obtained by applying the change-of-enrichment 2-functor (\cref{mult-change-enrichment})
\[
  \dF \cn \MCat \to \NCat
\]
to $\M$ as an $\M$-category.
The definition of $\Fse$ is given on objects by $F$ and on hom objects via $F$ together with the evaluation $\ev^\M$ from the closed structure of $\M$.

Compositionality of the standard enrichment construction is treated in \cref{sec:fun-std-enr-multi}.
For a composable pair of non-symmetric multifunctors between closed non-symmetric multicategories,
\[
  \big(\M,\clM,\evM\big) \fto{F} \big(\N,\clN,\evN\big) \fto{G} \big(\P,\clP,\evP\big),
\]
\cref{gspectra-thm-iv} shows that
the following diagram of $\P$-functors commutes.
\begin{equation}\label{intro-gspectra-iv-diagram}
  \begin{tikzpicture}[vcenter]
    \def\v{-1.2}
    \draw[0cell]
    (0,0) node (a) {\M_{GF}}
    (a)+(3,0) node (b) {\P}
    (a)+(0,\v) node (c) {(\M_F)_G}
    (b)+(0,\v) node (d) {\N_G}
    ;
    \draw[1cell=.9]
    (a) edge node {\GFse} (b)
    (a) edge[-,double equal sign distance] (c)
    (c) edge node {\Fse_G} (d)
    (d) edge node[swap] {\Gse} (b)
    ;
  \end{tikzpicture}
\end{equation}

\cref{sec:factor-Kemse} applies this to the factorization of Elmendorf-Mandell $K$-theory, $\Kem$ from \cref{Kem}.
\begin{equation}\label{intro-Kem-four}
  \begin{tikzpicture}[vcenter]
    \def\v{-1.3} \def\h{6}
    \draw[0cell=.9]
    (0,0) node (p) {\permcatsu}
    (p)+(\h,0) node (sp) {\Sp}
    (p)+(0,\v) node (m) {\MoneMod}
    (m)+(\h/2,0) node (gc) {\Gstarcat}
    (sp)+(0,\v) node (gs) {\Gstarsset}
    ;
    \draw[1cell=.9]
    (p) edge node {\Kem} (sp)
    (p) edge node[swap] {\Endm} (m)
    (m) edge node {\Jt} (gc)
    (gc) edge node {\Ner_*} (gs)
    (gs) edge node[swap] {\Kg} (sp)
    ;
  \end{tikzpicture}
\end{equation}
The introduction of \cref{sec:factor-Kemse} gives further review of the closed multicategories and multifunctors above.
The factorization above yields a corresponding factorization of the standard enrichment $\Kemse$ into four spectrally-enriched functors.
\cref{gspectra-thm-xi} gives a precise statement, and the remainder of \cref{sec:factor-Kemse} gives further details.

\subsection*{Connection with Other Chapters}

The diagram change of enrichment theory in \cref{ch:gspectra_Kem} depends on the self-enrichments and standard enrichments developed here.
That is then used to study the homotopy theory of enriched diagram categories, in \cref{ch:mackey,ch:mackey_eq}.
The standard enrichment of $\Kem$ and its factorization are used in \cref{sec:presheaf-K,sec:mult-mackey-spectra} for the development of corresponding spectral Mackey functors.

\subsection*{Background}

The development of self-enrichment and standard enrichment in \cref{sec:selfenr-clmulti,sec:std-enr-multifunctor,sec:fun-std-enr-multi} depends on the theory of multicategorical enrichment and closed multicategories, from \cref{ch:menriched,ch:change_enr,ch:gspectra}.
\cref{sec:factor-Kemse} uses the definition and factorization of Elmendorf-Mandell $K$-theory, $\Kem$, from \cref{sec:segalEMK}.

\subsection*{Chapter Summary}

\cref{sec:selfenr-clmulti} defines the self-enrichment of non-symmetric closed multicategories.
\cref{sec:std-enr-multifunctor} defines the standard enrichment of non-symmetric multifunctors between non-symmetric closed multicategories.
\cref{sec:fun-std-enr-multi} shows that the standard enrichment construction respects composition of non-symmetric multifunctors.
This is applied in \cref{sec:factor-Kemse} to factor the standard enrichment of the Elmendorf-Mandell $K$-theory functor, $\Kem$, into four spectrally enriched functors.
Here is a summary table.
\reftable{.9}{
  self-enrichment of a closed multicategory
  & \ref{def:clmulti-clcat} and \ref{cl-multi-cl-cat}
  \\ \hline
  application to permutative categories
  & \ref{ex:perm-closed-multicat}
  \\ \hline
  application to symmetric monoidal closed categories
  & \ref{ex:cl-multi-cl-cat}
  \\ \hline
  standard enrichment of a multifunctor
  & \ref{def:gspectra-thm-iii-context} and
  \ref{gspectra-thm-iii}
  \\ \hline
  examples of standard enrichment
  & \ref{ex:Fstse} and \ref{std-enr-monoidal}
  \\ \hline
  compositionality of standard enrichment
  & \ref{gspectra-thm-iv}
  \\ \hline
  application to Elmendorf-Mandell $K$-theory
  & \ref{gspectra-thm-xi},
  \ref{expl:Kemse},
  \ref{expl:Endmse},
  \ref{expl:Jtse}, and
  \ref{expl:Nerstarse}
  \\
}

We remind the reader of \cref{conv:universe} about universes and \cref{expl:leftbracketing} about left normalized bracketing for iterated products.

\section{Self-Enrichment of Closed Multicategories}
\label{sec:selfenr-clmulti}

In this section we observe that each non-symmetric closed multicategory is enriched in itself; see \cref{cl-multi-cl-cat}.  In this context of self-enrichment, we prove two consistency results.
\begin{enumerate}
\item In \cref{ex:perm-closed-multicat} we observe that, for the closed multicategory $\permcatsu$ of small permutative categories, the self-enrichment obtained from the closed multicategory structure coincides with the one in \cref{permcat-selfenr}.
\item For each symmetric monoidal closed category $\V$, the self-enrichment is the same whether $\V$ is regarded as a symmetric monoidal closed category or as a closed multicategory; see \cref{ex:cl-multi-cl-cat}.  
\end{enumerate}
The self-enrichment of a non-symmetric closed multicategory is an integral part of several key constructions later in this work, including:
\begin{itemize}
\item the standard enrichment of a multifunctor (\cref{def:gspectra-thm-iii-context}),
\item diagrams enriched in a multicategory \cref{mcat-cm}, and
\item Mackey functors enriched in a multicategory \cref{mcat-copm}.
\end{itemize} 
In short, self-enrichment is one of the key features of non-symmetric closed multicategories, and we will make full use of it.

\subsection*{Canonical Self-Enrichment}

Recall $\M$-categories for a non-symmetric multicategory $\M$ (\cref{def:menriched-cat}).  For a non-symmetric closed multicategory $\M$ (\cref{def:closed-multicat}), recall that two multimorphisms are called \emph{partners} if they correspond to each other under the evaluation bijection \cref{eval-bijection}.
\[\begin{tikzcd}[/tikz/column 1/.append style={anchor=base east},
/tikz/column 2/.append style={anchor=base west}, column sep=huge, row sep=0ex]
\M\left(\angx \sscs \clM\scmap{\angy;z} \right) \ar{r}{\chi_{\angx;\, \angy;\, z}}[swap]{\iso} 
& \M\scmap{\angx, \angy; z}\\
f \ar[maps to]{r} & \ga\big(\ev_{\angy;\, z} \sscs f, \ang{1_{y_j}}_{j=1}^p\big)
\end{tikzcd}\]  
Now we define the canonical self-enrichment.

\begin{definition}\label{def:clmulti-clcat}
Suppose $(\M,\clM,\ev)$ is a non-symmetric closed multicategory.  We define the data of an $\M$-category 
\[(\M,\circ,i),\]
which is called the \index{canonical self-enrichment}\index{self-enrichment}\index{multicategory!closed!self-enrichment}\emph{canonical self-enrichment of $\M$}, as follows.
\begin{description}
\item[Objects] The objects are those of $\M$.
\item[Hom Objects] For each pair of objects $x,y \in \M$, the morphism object is the unary internal hom object $\clM\scmap{x;y}$ in $\M$.  This is the $n=1$ case of \cref{clMangxy}.
\item[Identities] For each object $x$ in $\M$, the identity 
\begin{equation}\label{selfenrm-id}
i_x \cn \ang{} \to \clM\scmap{x;x},
\end{equation}
which is a nullary multimorphism in $\M$, is defined as the partner of the $x$-colored unit $\opu_x \in \M\scmap{x;x}$.  
\item[Composition]
For objects $x,y,z \in \M$, the composition
\begin{equation}\label{selfenrm-comp}
\left(\clM\scmap{y;z} \scs \clM\scmap{x;y}\right) \fto{\comp} \clM\scmap{x;z},
\end{equation}
which is a binary multimorphism in $\M$, is defined as the partner of the following 3-ary multimorphism.
\begin{equation}\label{selfenrm-comp-partner}
\Big(\clM\scmap{y;z}\scs \clM\scmap{x;y} \scs x\Big) 
\fto{(\opu, \ev_{x;\,y})} \big(\clM\scmap{y;z}\scs y\big)
\fto{\ev_{y;\,z}} z
\end{equation}
\end{description}
This finishes the definition of the canonical self-enrichment of $\M$.  \cref{cl-multi-cl-cat} proves that it is an $\M$-category.
\end{definition}

\begin{explanation}[Canonical Self-Enrichment]\label{expl:selfenrm}
\cref{def:clmulti-clcat} does not use anything about symmetric group action, either on $\M$ or on the internal hom objects.  Thus it makes sense for a \emph{non-symmetric} closed multicategory.

By definition \cref{eval-bijection}, the composition $\comp$ in \cref{selfenrm-comp} and the identity $i_x$ in \cref{selfenrm-id} are the \emph{unique} binary, respectively, nullary, multimorphisms in $\M$ that make the following two diagrams in $\M$ commute.
\begin{equation}\label{selfenrm-diagrams}
\begin{tikzpicture}[vcenter]
\def\h{3.7} \def\v{-1.3}
\draw[0cell=.8] 
(0,0) node (a) {\Big(\clM\scmap{y;z}\scs \clM\scmap{x;y} \scs x\Big)}
(a)+(\h,0) node (b) {\big(\clM\scmap{x;z} \scs x\big)}
(a)+(0,\v) node (c) {\big(\clM\scmap{y;z} \scs y\big)}
(b)+(0,\v) node (d) {z}
;
\draw[1cell=.85]
(a) edge node {(\comp \scs \opu_x)} (b)
(b) edge node {\ev_{x;\,z}} (d)
(a) edge node[swap] {\left(\opu \scs \ev_{x;\,y}\right)} (c)
(c) edge node {\ev_{y;\,z}} (d)
;
\begin{scope}[shift={(5.5,0)}]
\draw[0cell=.8]
(0,0) node (a) {\big(\ang{} \scs x\big)}
(a)+(2.7,0) node (b) {\big(\clM\scmap{x;x} \scs x\big)}
(b)+(0,\v) node (d) {x}
;
\draw[1cell=.85]
(a) edge node {\left(i_x \scs \opu_x\right)} (b)
(b) edge node {\ev_{x;\,x}} (d)
(a) edge node[swap] {\opu_x} (d)
;
\end{scope}
\end{tikzpicture}
\end{equation}
We call these the \emph{associativity diagram} and the \emph{unity diagram}, respectively, of the canonical self-enrichment of $\M$.  These diagrams are the analogs of those in \cref{eq:adj-comp} for a symmetric monoidal closed category.
\end{explanation}

We now show that the canonical self-enrichment is well defined.

\begin{theorem}\label{cl-multi-cl-cat}
For each non-symmetric closed multicategory $(\M,\clM,\ev)$, the canonical self-enrichment of $\M$ in \cref{def:clmulti-clcat} is an $\M$-category.
\end{theorem}

\begin{proof}
We need to prove the associativity axiom \cref{menriched-cat-assoc} and the unity axiom \cref{menriched-cat-unity} for the canonical self-enrichment of $\M$.

\medskip
\emph{Associativity \cref{menriched-cat-assoc}}.  Consider objects $w,x,y,z \in \M$.  Since taking partners is a bijection \cref{eval-bijection}, it suffices to show that the two composites in the associativity diagram \cref{menriched-cat-assoc} for $\M$ have the same partners.  These two partners are the left and right boundary composites of the following diagram in $\M$, which we want to show is commutative.  We abbreviate $\clM\scmap{x;y}$ to $\clM_{x;\, y}$.
\[\begin{tikzpicture}
\def\g{2} \def\h{4} \def\u{-1} \def\v{-1.4} \def\t{5}
\draw[0cell=.8]
(0,0) node (a) {(\clM_{y;\, z} \scs \clM_{x;\, y} \scs \clM_{w;\, x} \scs w)}
(a)+(-\h,\u/2) node (b1) {(\clM_{y;\, z} \scs \clM_{w;\, y} \scs w)}
(a)+(0,\u) node (b2) {(\clM_{y;\, z} \scs \clM_{x;\, y} \scs x)}
(a)+(\h,\u/2) node (b3) {(\clM_{x;\, z} \scs \clM_{w;\, x} \scs w)}
(b1)+(0,\v+\u) node (c1) {(\clM_{w;\, z} \scs w)}
(b2)+(-\g,\u) node (c2) {(\clM_{y;\, z} \scs y)}
(b2)+(\g,\u) node (c3) {(\clM_{x;\, z} \scs x)}
(b3)+(0,\v+\u) node (c4) {(\clM_{w;\, z} \scs w)}
(c1)+(\h,\u/2) node (d) {z}
(a)+(\h/3,\u/2) node () {\pentagram}
;
\draw[1cell=.8]
(a) edge[bend right=\t] node[swap,pos=.3] {(\opu,\comp,\opu)} (b1)
(a) edge node[swap,pos=.6] {(\opu,\opu,\ev_{w;\, x})} (b2)
(a) edge[bend left=\t] node[pos=.3] {(\comp,\opu,\opu)} (b3)
(b1) edge node[swap] {(\comp,\opu)} (c1)
(b1) edge node[pos=.5] {(\opu,\ev_{w;\, y})} (c2)
(b2) edge node[pos=.6] {(\opu,\ev_{x;\, y})} (c2)
(b2) edge node[swap,pos=.6] {(\comp,\opu)} (c3)
(b3) edge node[swap,pos=.5] {(\opu,\ev_{w;\, x})} (c3)
(b3) edge node {(\comp,\opu)} (c4)
(c1) edge[bend right=\t] node[pos=.3] {\ev_{w;\, z}} (d)
(c2) edge node {\ev_{y;\, z}} (d)
(c3) edge node[swap] {\ev_{x;\, z}} (d)
(c4) edge[bend left=\t] node[swap,pos=.3] {\ev_{w;\, z}} (d)
;
\end{tikzpicture}\]
The following statements hold for the diagram above.
\begin{itemize}
\item The sub-region labeled $\pentagram$ commutes by definition.
\item The other four sub-regions commute by the associativity diagram in \cref{selfenrm-diagrams}.
\end{itemize}
This proves the associativity axiom \cref{menriched-cat-assoc} for the canonical self-enrichment of $\M$.

\medskip
\emph{Unity \cref{menriched-cat-unity}}.  Similar to associativity, it suffices to show that the partners of the composites in the unity diagram \cref{menriched-cat-unity} are equal.  These partners are the boundary composites and the middle $\ev_{x;\, y}$ in the following diagram in $\M$, which we want to show is commutative.
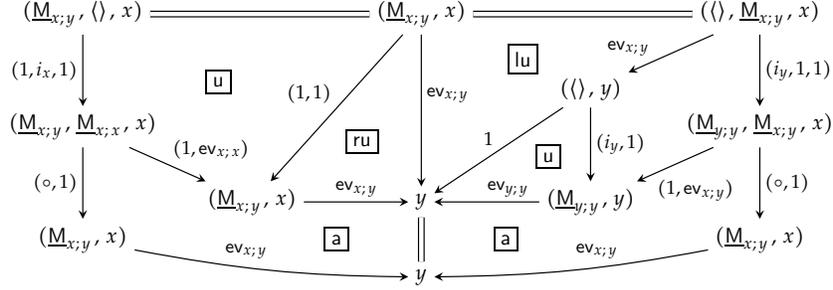
\[\begin{tikzpicture}
\def\g{2} \def\h{4.5} \def\u{-1} \def\v{-1.5} \def\t{5}
\draw[0cell=.8]
(0,0) node (a1) {(\clM_{x;\, y} \scs \ang{} \scs x)}
(a1)+(\h,0) node (a2) {(\clM_{x;\, y} \scs x)}
(a2)+(\h,0) node (a3) {(\ang{} \scs \clM_{x;\, y} \scs x)}
(a1)+(0,\v) node (b1) {(\clM_{x;\, y} \scs \clM_{x;\, x} \scs x)}
(a3)+(0,\v) node (b2) {(\clM_{y;\, y} \scs \clM_{x;\, y} \scs x)}
(b1)+(0,\v) node (c1) {(\clM_{x;\, y} \scs x)}
(b2)+(0,\v) node (c3) {(\clM_{x;\, y} \scs x)}
(c1)+(\h,-.5) node (c2) {y}
node[between=b1 and c2 at .5] (d1) {(\clM_{x;\, y} \scs x)}
(a2)+(\h/2,\u) node (d2) {(\ang{} \scs y)}
node[between=b2 and c2 at .5] (d3) {(\clM_{y;\, y} \scs y)}
node[between=d1 and d3 at .5] (c2') {y}
node[between=b1 and a2 at .4] () {\boxed{\mathsf{u}}}
node[between=d2 and d3 at .6, shift={(-.7,0)}] () {\boxed{\mathsf{u}}}
node[between=d1 and c2' at .65, shift={(0,1)}] () {\boxed{\mathsf{ru}}}
node[between=a2 and d2 at .6] () {\boxed{\mathsf{lu}}}
node[between=d1 and c2 at .5] () {\boxed{\mathsf{a}}}
node[between=d3 and c2 at .5] () {\boxed{\mathsf{a}}}
;
\draw[1cell=.8]
(a1) edge[-,double equal sign distance] (a2)
(a2) edge[-,double equal sign distance] (a3)
(c2') edge[-,double equal sign distance] (c2)
(a1) edge node[swap] {(\opu,i_x,\opu)} (b1)
(a3) edge node {(i_y,\opu,\opu)} (b2)
(b1) edge node[swap] {(\comp,\opu)} (c1)
(b2) edge node {(\comp,\opu)} (c3)
(c1) edge[bend right=\t] node[pos=.3] {\ev_{x;\, y}} (c2)
(c3) edge[bend left=\t] node[swap,pos=.3] {\ev_{x;\, y}} (c2)
(a2) edge node[pos=.4] {\ev_{x;\, y}} (c2')
(b1) edge node {(\opu, \ev_{x;\, x})} (d1)
(a2) edge node[swap] {(\opu,\opu)} (d1)
(d1) edge node {\ev_{x;\, y}} (c2')
(a3) edge node[swap,pos=.7] {\ev_{x;\, y}} (d2)
(d2) edge node[swap] {\opu} (c2')
(d2) edge node {(i_y,\opu)} (d3)
(d3) edge node[swap,pos=.3] {\ev_{y;\, y}} (c2')
(b2) edge node[pos=.8] {(\opu,\ev_{x;\, y})} (d3)
;
\end{tikzpicture}\]
The following statements hold for the diagram above.
\begin{itemize}
\item The two sub-regions labeled $\boxed{\mathsf{a}}$ commute by the associativity diagram in \cref{selfenrm-diagrams}.
\item The two sub-regions labeled $\boxed{\mathsf{u}}$ commute by the unity diagram in \cref{selfenrm-diagrams}.
\item The sub-regions labeled $\boxed{\mathsf{ru}}$ and $\boxed{\mathsf{lu}}$ commute by, respectively, the right unity \cref{enr-multicategory-right-unity} and left unity \cref{enr-multicategory-left-unity} of $\M$.
\item The remaining unlabeled quadrilateral commutes by definition.
\end{itemize}
This proves the unity axiom \cref{menriched-cat-unity} for the canonical self-enrichment of $\M$.
\end{proof}

\subsection*{Self-Enrichment of $\permcatsu$}

Consider the closed multicategory $\psu = \permcatsu$ of small permutative categories (\cref{perm-closed-multicat}).  By \cref{cl-multi-cl-cat} $\psu$ has a canonical self-enrichment.  In other words, $\psu$ has the structure of a $\psu$-category.  Moreover, we previously established a $\psu$-category structure on $\psu$ in \cref{permcat-selfenr}.  Now we observe that these $\psu$-categories are the same.

\begin{proposition}\label{ex:perm-closed-multicat}\index{canonical self-enrichment!for permutative categories}\index{self-enrichment!of permutative categories}\index{permutative category!self-enrichment of multicategory of}
For the closed multicategory $\permcatsu$,
\begin{itemize}
\item the self-enrichment in \cref{permcat-selfenr} and
\item the canonical self-enrichment in \cref{cl-multi-cl-cat}
\end{itemize} 
are equal as $\psu$-categories.
\end{proposition}

\begin{proof}
The two $\psu$-categories in question are the same for the following reasons.
\begin{romenumerate}
\item\label{perm-clmulti-i} In each of \cref{permcat-selfenr,cl-multi-cl-cat}, the $\psu$-category has small permutative categories as objects.
\item\label{perm-clmulti-ii} For small permutative categories $\C$ and $\D$, the hom object in \cref{def:clmulti-clcat} is the small permutative category $\clpsu(\C;\D)$ in \cref{def:clp-angcd}.  This is the same as the permutative category $\psu(\C,\D)$ in \cref{psucd-hom-permcat}, as we pointed out in \cref{expl:clp-zero}.
\item\label{perm-clmulti-iii} The bilinear evaluation $\ev_{\C;\,\D}$ in \cref{clp-evaluation} coincides with the bilinear evaluation $\ev_{\C,\D}$ in \cref{ev-bilinear}, as we pointed out in \cref{expl:clp-eval-zero}.
\item\label{perm-clmulti-iv} The identity of a small permutative category $\C$ in the sense of \cref{selfenrm-id} is the unique 0-linear functor
\[i_\C \cn \boldone \to \clpsu(\C;\C),\]
which means a strictly unital symmetric monoidal functor $\C \to \C$, that makes the unity diagram in \cref{selfenrm-diagrams} commutative.  By \cref{def:clp-eval} for $\ev_{\C;\C}$, the identity symmetric monoidal functor $1_\C$ makes the unity diagram in \cref{selfenrm-diagrams} commutative.  Thus uniqueness implies that $i_\C$ is given by $1_\C$.  This is the same as the identity in \cref{def:perm-selfenr}.
\item\label{perm-clmulti-v} For small permutative categories $\B$, $\C$, and $\D$, the composition $\comp$ in \cref{selfenrm-comp} is the unique bilinear functor 
\[\clpsu(\C;\D) \times \clpsu(\B;\C) \fto{\comp} \clpsu(\B;\D)\]
that makes the associativity diagram in \cref{selfenrm-diagrams} commutative.  By \cref{ev-comp} the composition bilinear functor $\mcomp_{\B,\C,\D}$ in \cref{psu-mBCD} also makes the associativity diagram in \cref{selfenrm-diagrams} commutative.  Thus uniqueness implies 
\[\comp = \mcomp_{\B,\C,\D},\] 
which is the composition in \cref{def:perm-selfenr}.
\end{romenumerate}
This finishes the proof.
\end{proof}

\subsection*{Self-Enrichment of Symmetric Monoidal Closed Categories}

For a symmetric monoidal closed category $\V$, by \cref{smclosed-closed-multicat} there is an endomorphism closed multicategory 
\[\big(\End\,\V \scs \clEndV \scs \ev\big)\] 
with internal hom objects and evaluation induced by those of $\V$.  There are two self-enrichment constructions in this context:
\begin{enumerate}
\item $\V$ has a canonical self-enrichment (\cref{theorem:v-closed-v-sm}).
\item By \cref{cl-multi-cl-cat} the endomorphism multicategory $\End\,\V$ has a canonical self-enrichment, which is an $(\End\,\V)$-category.  By \cref{EndV-enriched} enrichment in $\V$ and in $\End\,\V$ are the same thing.  Thus we may also regard the canonical self-enrichment of $\End\,\V$ as a $\V$-category (\cref{def:enriched-category}).
\end{enumerate} 
Now we observe that these two self-enrichment constructions are the same.

\begin{proposition}\label{ex:cl-multi-cl-cat}\index{canonical self-enrichment!symmetric monoidal}\index{symmetric monoidal category!closed!self-enrichment}
For a symmetric monoidal closed category $\big(\V,\otimes,\tu,[,]\big)$, 
\begin{itemize}
\item the canonical self-enrichment of $\V$ in \cref{theorem:v-closed-v-sm} and
\item the canonical self-enrichment of $\End\,\V$ in \cref{cl-multi-cl-cat}
\end{itemize} 
are equal as $\V$-categories.
\end{proposition}

\begin{proof}
We compare \cref{definition:canonical-v-enrichment} for the canonical self-enrichment of $\V$ and \cref{def:clmulti-clcat} for the canonical self-enrichment of $\EndV$.

\medskip
\emph{Objects and Morphism Objects}. The objects of $\End\,\V$ are those of $\V$.  For objects $x,y \in \End\,\V$, the morphism object in $\End\,\V$ is the internal hom object
\[\clEndV\scmap{x;y} = [x,y] \inspace \V.\]
This is the same hom object as in \cref{definition:canonical-v-enrichment}.

\medskip
\emph{Identities}. 
The identity of an object $x \in \End\,\V$ is the nullary multimorphism \cref{selfenrm-id}
\[i_x \in (\End\,\V)\big(\ang{} \sscs \clEndV(x;x)\big) = \V\big(\tu, [x,x]\big)\]
that makes the left diagram in $\End\,\V$ below commutative.
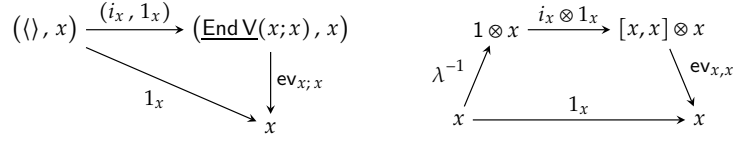
\[\begin{tikzpicture}
\def\h{.5} \def\v{-1.2}
\draw[0cell=.85]
(0,0) node (a) {\big(\ang{} \scs x\big)}
(a)+(3,0) node (b) {\big(\clEndV(x;x) \scs x\big)}
(b)+(0,-1.3) node (d) {x}
;
\draw[1cell=.9]
(a) edge node {\left(i_x \scs \opu_x\right)} (b)
(b) edge node {\ev_{x;\,x}} (d)
(a) edge node[swap] {\opu_x} (d)
;
\begin{scope}[shift={(6,0)}]
\draw[0cell=.85]
(0,0) node (a) {\tu \otimes x}
(a)+(2.2,0) node (b) {[x,x] \otimes x}
(a)+(-\h,\v) node (c) {x}
(b)+(\h,\v) node (d) {x}
;
\draw[1cell=.9]
(c) edge node {1_x} (d)
(c) edge node[pos=.3] {\la^\inv} (a)
(a) edge node {i_x \otimes 1_x} (b)
(b) edge node[pos=.6] {\ev_{x,x}} (d)
;
\end{scope}
\end{tikzpicture}\]
The left diagram in $\EndV$ above means the right commutative diagram in $\V$, with $\lambda$ denoting the left unit isomorphism.  Comparing the right diagram above with the right diagram in \cref{eq:adj-comp}, the uniqueness of adjoints implies that $i_x$ is equal to the identity of $x$ in the canonical self-enrichment of $\V$.

\medskip
\emph{Composition}. 
For objects $x,y,z \in \End\,\V$, the composition in $\End\,\V$ is the binary multimorphism \cref{selfenrm-comp}
\[\begin{split}
\comp \in &\, (\End\,\V)\Big( \clEndV(y;z) \scs \clEndV(x;y) \sscs \clEndV(x;z) \Big)\\
&= \V\big( [y,z] \otimes [x,y] \scs [x,z] \big)
\end{split}\]
that makes the left diagram in $\End\,\V$ below commutative.
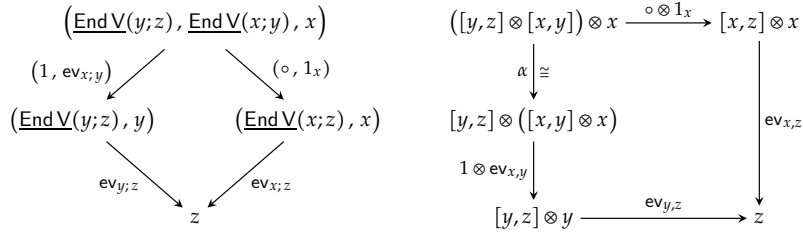
\[\begin{tikzpicture}[vcenter]
\def\h{1.5} \def\v{-1.3}
\draw[0cell=.8] 
(0,0) node (a) {\Big(\clEndV(y;z) \scs \clEndV(x;y) \scs x\Big)}
(a)+(\h,\v) node (b) {\big(\clEndV(x;z) \scs x\big)}
(a)+(-\h,\v) node (c) {\big(\clEndV(y;z) \scs y\big)}
(a)+(0,2*\v) node (d) {z}
;
\draw[1cell=.8]
(a) edge node[pos=.7] {(\comp \scs \opu_x)} (b)
(b) edge node {\ev_{x;\,z}} (d)
(a) edge node[swap,pos=.8] {\left(\opu \scs \ev_{x;\,y}\right)} (c)
(c) edge node[swap] {\ev_{y;\,z}} (d)
;
\begin{scope}[shift={(4.5,0)}]
\draw[0cell=.8] 
(0,0) node (l) {\big([y,z] \otimes [x,y]\big) \otimes x}
(l)+(0,\v) node (r) {[y,z] \otimes \big([x,y] \otimes x\big)}
(r)+(0,\v) node (yzy) {[y,z] \otimes y}
(l)+(3,0) node (xzx) {[x,z] \otimes x}
(xzx)+(0,2*\v) node (z) {z}
;
\draw[1cell=.8]
(l) edge node {\comp \otimes 1_x} (xzx)
(xzx) edge node {\ev_{x,z}} (z)
(l) edge node[swap] {\al} node {\iso} (r)
(r) edge['] node {1 \otimes \ev_{x,y}} (yzy)
(yzy) edge node {\ev_{y,z}} (z)
;
\end{scope}
\end{tikzpicture}\]
The left diagram in $\EndV$ above means the right commutative diagram in $\V$, with $\alpha$ denoting the associativity isomorphism.  Comparing the right diagram above with the left diagram in \cref{eq:adj-comp}, the uniqueness of adjoints implies that $\circ$ is equal to the composition $\mcomp$ in the canonical self-enrichment of $\V$.
\end{proof}

\section{Standard Enrichment of a Multifunctor}
\label{sec:std-enr-multifunctor}


In \cref{cl-multi-cl-cat} we showed that each non-symmetric closed multicategory $\M$ has a canonical self-enrichment.  In this section we use the canonical self-enrichment to show that each non-symmetric multifunctor $F$ between non-symmetric closed multicategories induces a multicategorically enriched functor $\Fse$, called the standard enrichment of $F$.  We discuss further functoriality properties and an application to $K$-theory in \cref{sec:fun-std-enr-multi,sec:factor-Kemse}.  In subsequent chapters, the standard enrichment is one of the two key constructions for change of enrichment for enriched diagrams \cref{diag-change-enr-assign} and enriched Mackey functors \cref{presheaf-change-enr}.  

Here is an outline of this section.
\begin{itemize}
\item The standard enrichment is constructed in \cref{def:gspectra-thm-iii-context} and verified in \cref{gspectra-thm-iii}.
\item As an illustration of the construction, \cref{ex:Fstse} describes the standard enrichment of the non-symmetric multifunctor 
\[\Fst \cn \pMulticat \to \permcatsu\]
in \cref{ptmulticat-xvii}.
\item As a consistency check, in \cref{std-enr-monoidal} we prove that, for a monoidal functor between symmetric monoidal closed categories, the two standard enrichment constructions in \cref{proposition:U-std-enr,gspectra-thm-iii} are the same.
\end{itemize}

\subsection*{Standard Enrichment}

Recall from \cref{def:closed-multicat} that, in a non-symmetric closed multicategory, two multimorphisms that correspond to each other under the evaluation bijection \cref{eval-bijection} are called \emph{partners}.  The partner of a multimorphism $f$ is denoted $\pnf$.  

\begin{definition}[Standard Enrichment]\label{def:gspectra-thm-iii-context}
For a non-symmetric multifunctor (\cref{def:enr-multicategory-functor}) between non-symmetric closed multicategories
\[F \cn \big(\M,\clM,\ev^\M\big) \to \big(\N,\clN,\ev^\N\big),\]
we define the data of an $\N$-functor (\cref{def:mfunctor})
\begin{equation}\label{std-enr-F}
\Fse \cn \M_F \to \N,
\end{equation}
which is called the \index{enrichment!standard}\index{standard enrichment!multifunctor}\index{multifunctor!standard enrichment}\emph{standard enrichment of $F$}, as follows.
\begin{description}
\item[Domain] The domain of $\Fse$ is the $\N$-category $\M_F$ obtained from the canonical self-enrichment of $\M$ (\cref{cl-multi-cl-cat}), which is an $\M$-category, by applying the change-of-enrichment 2-functor (\cref{mult-change-enrichment})
\[\dF \cn \MCat \to \NCat.\]
\item[Codomain] The codomain of $\Fse$ is the canonical self-enrichment of $\N$ (\cref{cl-multi-cl-cat}), which is an $\N$-category.
\item[Object Assignment] $\Fse$ has the same object assignment as $F$.
\item[Component Morphisms] For each pair of objects $x,y \in \M$, the $(x,y)$-component unary multimorphism
\begin{equation}\label{Fprimexy}
\Fse_{x,y} = \pn{\big(F(\ev^\M_{x;\,y}) \big)} \cn F\clM(x;y) \to \clN(Fx;Fy) \inspace \N
\end{equation}
is defined as the partner of the binary multimorphism
\begin{equation}\label{Fev}
F(\ev^\M_{x;\,y}) \cn \big( F\clM(x;y), Fx \big) \to Fy \inspace \N.
\end{equation}
This is the image under $F$ of the evaluation binary multimorphism \cref{evangxy}
\[\ev^\M_{x;\,y} \cn \big(\clM(x;y), x\big) \to y\]
at $(x;\, y)$, which is part of the closed multicategory structure of $\M$.
\end{description}
This finishes the definition of the standard enrichment $\Fse$.  We verify that $\Fse$ is an $\N$-functor in \cref{gspectra-thm-iii} below.
\end{definition}

Before we show that $\Fse$ is an $\N$-functor, let us explain some aspects of the standard enrichment construction.

\begin{explanation}[Symmetry is Not Required]\label{expl:Fse-symmetry}
In \cref{def:gspectra-thm-iii-context}, even if $\M$ and $\N$ are closed multicategories---as opposed to non-symmetric ones---$F$ is \emph{not} required to preserve the symmetric group action on (i) $\M$ and $\N$ \cref{enr-multifunctor-equivariance} and (ii) their internal hom objects \cref{sigma-clMxy}.  This is possible because the change-of-enrichment 2-functor $\dF$ (\cref{mult-change-enrichment}) and the canonical self-enrichment (\cref{cl-multi-cl-cat}) do not require symmetry.  This point is important for our applications in \cref{ch:mackey_eq} where we consider the non-symmetric multifunctors $\Fst$ and $\Fm$ in \cref{ptmulticat-xvii,Fm-multi-def}; see \cref{mackey-xiv-pmulticat,mackey-xiv-mone}.  We discuss the standard enrichment of $\Fst$ in \cref{ex:Fstse} below.
\end{explanation}

\begin{explanation}[Domain of $\Fse$]\label{expl:M-sub-F}
Interpreting the change of enrichment (\cref{def:mult-change-enr}) for the canonical self-enrichment of $\M$ (\cref{cl-multi-cl-cat}), we unpack the $\N$-category $\M_F$, which is the domain of the standard enrichment $\Fse$, as follows. 
\begin{itemize}
\item $\M_F$ has the same objects as $\M$. 
\item For each pair of objects $x,y \in \M_F$, the hom object is
\[(\M_F)(x,y) = F\clM(x;y) \inspace \N.\]
\item The identity of an object $x \in \M_F$ is the nullary multimorphism
\begin{equation}\label{MsubF-identity}
\ang{} \fto{F(i_x)} F\clM(x;x) \inspace \N
\end{equation}
given by applying $F$ to the identity \cref{selfenrm-id} 
\[\ang{} \fto{i_x} \clM(x;x) \inspace \M.\] 
\item For objects $x,y,z \in \M_F$, the composition binary multimorphism
\begin{equation}\label{MsubF-composition}
\big( F\clM(y;z) \scs F\clM(x;y) \big) \fto{F(\comp)} F\clM(x;z) \inspace \N
\end{equation}
is given by applying $F$ to the composition \cref{selfenrm-comp} 
\[\big(\clM(y;z) \scs \clM(x;y)\big) \fto{\comp} \clM(x;z) \inspace \M.\]
\end{itemize} 
Moreover, applying the non-symmetric multifunctor $F \cn \M \to \N$ to the associativity and unity diagrams in \cref{selfenrm-diagrams} yields the following commutative diagrams in $\N$.  
\begin{equation}\label{FseM-diagrams}
\begin{tikzpicture}[vcenter]
\def\h{3.7} \def\v{-1.3}
\draw[0cell=.7] 
(0,0) node (a) {\Big(F\clM\scmap{y;z} \scs F\clM\scmap{x;y} \scs Fx\Big)}
(a)+(\h,0) node (b) {\big(F\clM\scmap{x;z} \scs Fx\big)}
(a)+(0,\v) node (c) {\big(F\clM\scmap{y;z} \scs Fy\big)}
(b)+(0,\v) node (d) {Fz}
;
\draw[1cell=.7]
(a) edge node {(F(\comp) \scs \opu)} (b)
(b) edge node {F(\ev^\M_{x;\,z})} (d)
(a) edge node[swap] {\left(\opu \scs F(\ev^\M_{x;\,y})\right)} (c)
(c) edge node {F(\ev^\M_{y;\,z})} (d)
;
\begin{scope}[shift={(5.5,0)}]
\draw[0cell=.7]
(0,0) node (a) {\big(\ang{} \scs Fx\big)}
(a)+(2.7,0) node (b) {\big(F\clM\scmap{x;x} \scs Fx\big)}
(b)+(0,\v) node (d) {Fx}
;
\draw[1cell=.7]
(a) edge node {\left(F(i_x) \scs \opu\right)} (b)
(b) edge node {F(\ev^\M_{x;\,x})} (d)
(a) edge node[swap] {\opu} (d)
;
\end{scope}
\end{tikzpicture}
\end{equation}
We emphasize that the two commutative diagrams in \cref{FseM-diagrams} use the fact that $F$ preserves colored units \cref{enr-multifunctor-unit} and composition \cref{v-multifunctor-composition}, but they do \emph{not} require $F$ to preserve the symmetric group action even if $\M$ and $\N$ are multicategories.
\end{explanation}

\begin{explanation}[Component Morphisms of $\Fse$]\label{expl:Fse-xy}
By definition \cref{Fprimexy}, for objects $x,y\in\M$, the $(x,y)$-component unary multimorphism $\Fse_{x,y}$ is the partner of $F\big(\ev^\M_{x;y}\big)$.  By the definition of the evaluation bijection \cref{eval-bijection}, this means that $\Fse_{x,y}$ is the \emph{unique} unary multimorphism that makes the following diagram in $\N$ commutative.
\begin{equation}\label{Fse-xy-diag}
\begin{tikzpicture}[vcenter]
\draw[0cell=.9]
(0,0) node (a) {\big( F\clM(x;y) \scs Fx \big)}
(a)+(4.5,0) node (b) {\big( \clN(Fx;Fy) \scs Fx \big)}
(b)+(0,-1.3) node (c) {Fy}
;
\draw[1cell=.9]
(a) edge node {\left( \Fse_{x,y} \scs 1_{Fx} \right)} (b)
(b) edge node {\ev^\N_{Fx;\, Fy}} (c)
(a) edge node[pos=.4,swap] {F(\ev^\M_{x;\,y})} (c)
;
\end{tikzpicture}
\end{equation}
The diagram \cref{Fse-xy-diag} is the analog of the first diagram in \cref{expl:std-enr} in the context of the standard enrichment of a monoidal functor between symmetric monoidal closed categories.  We make this connection precise in \cref{std-enr-monoidal} below.
\end{explanation}

Now we show that the standard enrichment is a well-defined enriched functor (\cref{def:mfunctor}).

\begin{theorem}\label{gspectra-thm-iii}
For each non-symmetric multifunctor between non-symmetric closed multicategories
\[F \cn \big(\M,\clM,\ev^\M\big) \to \big(\N,\clN,\ev^\N\big),\]
the standard enrichment in \cref{def:gspectra-thm-iii-context}
\[\Fse \cn \M_F \to \N\]
is an $\N$-functor.
\end{theorem}

\begin{proof}
We need to show that $\Fse$ preserves composition and identities in the sense of \cref{mfunctor-diagrams}.  For $\Fse$ those diagrams are the diagrams in $\N$ in \cref{Fse-functor-diagrams} below for objects $x,y,z \in \M$.  We use \cref{MsubF-identity,MsubF-composition} for the identities and composition of $\M_F$, and we denote $\clM(x;y)$ by $\clM_{x;y}$.
\begin{equation}\label{Fse-functor-diagrams}

\end{equation}
Since taking partners is a bijection \cref{eval-bijection}, it suffices to show that, in each diagram in \cref{Fse-functor-diagrams}, the partners of the two composites are equal.

\medskip
\emph{Preservation of Composition}.  For the left diagram in \cref{Fse-functor-diagrams}, the partners of the two composites are the two boundary composites in the following diagram in $\N$, which we want to show is commutative.
\[%
\]
The diagram above is commutative for the following reasons.
\begin{itemize}
\item The triangle labeled $\clubsuit$ is commutative by definition.
\item The sub-regions labeled $\spadesuit$ and $\filledlozenge$ are commutative by the left diagrams in, respectively, \cref{selfenrm-diagrams,FseM-diagrams}.
\item The remaining three sub-regions are commutative by \cref{Fse-xy-diag}.
\end{itemize}
This proves that the left diagram in \cref{Fse-functor-diagrams} is commutative.

\medskip
\emph{Preservation of Identities}.  For the right diagram in \cref{Fse-functor-diagrams}, the partners of the two composites are the two boundary composites in the following diagram in $\N$.
\[\begin{tikzpicture}
\def\g{.8} \def\h{5} \def\v{-1.4} \def\s{10}
\draw[0cell=.8]
(0,0) node (a1') {}
(a1')+(\h,0) node (a2') {}
(a1')+(\g,0) node (a1) {\big(\ang{} \scs Fx\big)}
(a2')+(-\g,0) node (a2) {\big(F\clM_{x;\,x} \scs Fx\big)}
(a1')+(0,\v) node (b1) {\big(\clN_{Fx;\,Fx} \scs Fx\big)}
(a2')+(0,\v) node (b2) {\big(\clN_{Fx;\,Fx} \scs Fx\big)}
(b1)+(\h/2,-.5) node (c) {Fx}
;
\draw[1cell=.8]
(a1) edge node[swap,pos=.6] {(i_{Fx} \scs \opu)} (b1)
(b1) edge[bend right=\s] node[pos=.2] {\evN} (c)
(a1) edge node {(F(i_x) \scs \opu)} (a2)
(a2) edge node[pos=.7] {(\Fse_{x,x} \scs \opu)} (b2)
(b2) edge[bend left=\s] node[swap,pos=.2] {\evN} (c)
(a1) edge node {\opu} (c)
(a2) edge node[swap,pos=.4] {F(\evM)} (c)
;
\end{tikzpicture}\]
From left to right, the three sub-regions in the diagram above are commutative by, respectively, the unity diagram in \cref{selfenrm-diagrams}, the right diagram in \cref{FseM-diagrams}, and \cref{Fse-xy-diag}.  This proves that the right diagram in \cref{Fse-functor-diagrams} is commutative.
\end{proof}

\subsection*{Examples of Standard Enrichment}

\begin{example}[Standard Enrichment of $\Fst$]\label{ex:Fstse}\index{category!free permutative!pointed}\index{permutative category!free!pointed}\index{free!permutative category!pointed}\index{pointed!free permutative category}
As an illustration of \cref{def:gspectra-thm-iii-context}, consider the non-symmetric multifunctor (\cref{ptmulticat-xvii})
\[\Fst \cn \pMulticat \to \permcatsu\]
and its standard enrichment $\permcatsu$-functor
\begin{equation}\label{Fstse-definition}
\Fstse \cn (\pMulticat)_{\Fst} \to \permcatsu.
\end{equation}
For the context, first recall the following.
\begin{itemize}
\item $\pMulticat = \pM$ is a symmetric monoidal closed category (\cref{thm:pmulticat-smclosed}). 
\begin{itemize}
\item Its objects are small pointed multicategories (\cref{def:ptd-multicat}).
\item Its monoidal product is the smash product, $\sma$, in \cref{eq:multicat-smash-pushout}.
\item Its internal hom is the pointed internal hom, $\pHom$, in \cref{eq:multicat-pHom}.  
\end{itemize} 
Via its endomorphism multicategory, $\pM$ is a closed multicategory (\cref{smclosed-closed-multicat}).  It is enriched in itself (\cref{ex:cl-multi-cl-cat}).
\item $\permcatsu = \psu$ is a closed multicategory (\cref{perm-closed-multicat}).
\begin{itemize}
\item Its objects are small permutative categories (\cref{def:symmoncat}). 
\item Its multicategory structure is discussed in \cref{sec:multpermcat}. 
\item Its internal hom objects, their symmetric group action, and multicategorical evaluation are constructed in, respectively, \cref{def:clp-angcd,def:clp-ancd-sigma,def:clp-eval}
\item We proved its evaluation bijection axiom in \cref{clp-evbij-axiom}.  In particular, the inverse of the function $\chi$ is the function $\Psi$ in \cref{chiPsi}.
\end{itemize}  
As a closed multicategory, $\psu$ is enriched in itself (\cref{ex:perm-closed-multicat}).
\item $\Fst$ is a non-symmetric $\Cat$-multifunctor, hence also a non-symmetric multifunctor.  It is genuinely non-symmetric because its construction \cref{Sst-Fbst} involves $\Fst^n$ in \cref{eq:Frn-Fstn}, which is induced by $\Fr^n$ in \cref{def:S-multi}.  As we mentioned in \cref{expl:Fnotsymmetric} \cref{Fnotsym-i}, $\Fr^n$ is not compatible with permutations.  So $\Fst^n$ is also not compatible with permutations, leading to the non-symmetry of $\Fst$.
\item As we discussed in detail in \cref{ex:dFst}, the non-symmetric multifunctor $\Fst$ has an associated change-of-enrichment 2-functor 
\[\dFst \cn \pMulticatcat \to \permcatsucat.\]
When we apply $\dFst$ to the $\pM$-category $\pM$, as in \cref{expl:M-sub-F}, we obtain the $\psu$-category $(\pM)_{\Fst}$, which is the domain of $\Fstse$.
\end{itemize}

The object assignment of the standard enrichment $\Fstse$ is the same as that of $\Fst$ in \cref{def:Fst-permutative}.  For each pair of small pointed multicategories $\X$ and $\Y$, the $(\X,\Y)$-component is the strictly unital symmetric monoidal functor
\begin{equation}\label{Fstse-XY-def}
\Fstse_{\X,\Y} = \pn{\big(\Fst(\ev^{\pM}_{\X,\Y})\big)} \cn 
\Fst\pHom(\X,\Y) \to \psu(\Fst\X,\Fst\Y).
\end{equation}
In other words, $\Fstse_{\X,\Y}$ is the image under $\Psi$ \cref{chiPsi} of the bilinear functor
\[\Fst(\ev^{\pM}_{\X,\Y}) \cn \Fst\pHom(\X,\Y) \times \Fst\X \to \Fst\Y.\]
This bilinear functor is the image under $\Fst$ \cref{Fst-multimorphism-cat} of the evaluation pointed multifunctor
\[\ev^{\pM}_{\X,\Y} \cn \pHom(\X,\Y) \sma \X \to \Y\]
in \cref{evaluation}, which is regarded as a binary multimorphism in $\pM$.
\end{example}

Recall from \cref{EndP} that each monoidal functor induces a non-symmetric multifunctor via the endomorphism construction.  The next observation is a consistency result.  It says that for a monoidal functor between symmetric monoidal closed categories, the standard enrichment in \cref{proposition:U-std-enr,gspectra-thm-iii} are the same.

\begin{proposition}\label{std-enr-monoidal}
For each monoidal functor between symmetric monoidal closed categories
\[(U,U^2,U^0) \cn \big(\V,\otimes,[,]\big) \to \big(\W,\otimes,[,]\big),\]
the following two $\W$-functors are the same:
\begin{itemize}
\item The standard enrichment 
\[\Use \cn \V_U \to \W\]
of $U$ in \cref{proposition:U-std-enr}.
\item The standard enrichment 
\[\EndUse \cn (\End\,\V)_{\End\,U} \to \End\,\W\]
in \cref{gspectra-thm-iii} of the non-symmetric multifunctor
\[\End\,U \cn \End\,\V \to \End\,\W.\]
\end{itemize}
\end{proposition}

\begin{proof}
This assertion follows by combining the following facts.
\begin{itemize}
\item By \cref{EndV-enriched} enrichment in the symmetric monoidal category $\V$ and in the multicategory $\End\,\V$ are the same thing. 
\item By \cref{smclosed-closed-multicat} $\End\,\V$ is a closed multicategory, with internal hom objects and evaluation induced by those of $(\V,[,])$. 
\item By \cref{ex:cl-multi-cl-cat} the canonical self-enrichment of the symmetric monoidal closed category $\V$ is equal to the canonical self-enrichment of the closed multicategory $\End\,\V$.  These statements also hold for $\W$.
\item By \cref{mon-change-enrichment} the change-of-enrichment 2-functors $\dU$ and $\dEndU$ are equal.  This yields an equality of $\W$-categories
\[\V_U = (\End\,\V)_{\End\,U}.\]
So the $\W$-functors in question, $\Use$ and $\EndUse$, have the same domain and the same codomain.
\item Both $\Use$ and $\EndUse$ have the same object assignment as $U$.
\end{itemize}

It remains to show that $\Use$ and $\EndUse$ have the same component morphisms.  For each pair of objects $x,y \in \V$, the diagram \cref{Fse-xy-diag} for $F = \End\,U$ is the diagram
\begin{equation}\label{EndUsexy}
\begin{tikzpicture}[vcenter]
\def\v{-1.3} \def\s{.9}
\draw[0cell=\s]
(0,0) node (a) {U[x,y] \otimes Ux}
(a)+(4.5,0) node (b) {[Ux,Uy] \otimes Ux}
(a)+(0,\v) node (c) {U\big([x,y] \otimes x \big)}
(b)+(0,\v) node (d) {Uy}
;
\draw[1cell=\s]
(a) edge node {\EndUse_{x,y} \otimes 1} (b)
(b) edge node {\ev} (d)
(a) edge node[swap] {U^2} (c)
(c) edge node {U(\ev)} (d)
;
\end{tikzpicture}
\end{equation}
in $\W$ because, by \cref{PtwoPf},
\[(\End\,U)(\ev) = U(\ev) \circ U^2.\]
The diagram \cref{EndUsexy} coincides with the first diagram in \cref{expl:std-enr}, which has $\Use_{x,y}$ in place of $\EndUse_{x,y}$.  Since \cref{Fse-xy-diag} uniquely defines $\EndUse_{x,y}$, we conclude that it is equal to $\Use_{x,y}$.  This proves that $\Use$ and $\EndUse$ have the same component morphisms.
\end{proof}

\section{Compositionality of Standard Enrichment}
\label{sec:fun-std-enr-multi}

For each non-symmetric multifunctor
\[F \cn \M \to \N\]
between non-symmetric closed multicategories, in \cref{gspectra-thm-iii} we constructed its standard enrichment $\N$-functor
\[\Fse \cn \M_F \to \N.\]
In this section we show that the standard enrichment construction, $F \mapsto \Fse$, respects composition of non-symmetric multifunctors; see \cref{gspectra-thm-iv}.  In particular, this compositionality property holds for monoidal functors between symmetric monoidal closed categories (\cref{ex:stdenr-func-monfunctor}).  In \cref{sec:factor-Kemse} we apply \cref{gspectra-thm-iv} to factor the standard enrichment of Elmendorf-Mandell $K$-theory, $\Kemse$, into four spectrally enriched functors.

\subsection*{Context of Compositionality}

For the set up, consider
\begin{equation}\label{MNP-FG}
\M \fto{F} \N \fto{G} \P
\end{equation}
consisting of
\begin{itemize}
\item non-symmetric closed multicategories $\M$, $\N$, and $\P$ (\cref{def:closed-multicat}) and
\item non-symmetric multifunctors $F$ and $G$ (\cref{def:enr-multicategory-functor}).
\end{itemize}
By \cref{func-change-enr} the following diagram of change-of-enrichment 2-functors is commutative.
\begin{equation}\label{gspectra-iv-enr-cat}
\begin{tikzpicture}[vcenter]
\def\h{2.5} \def\w{.6}
\draw[0cell=1]
(0,0) node (a) {\MCat}
(a)+(\h,0) node (b) {\NCat}
(b)+(\h,0) node (c) {\PCat}
;
\draw[1cell=.9]
(a) edge node {\dF} (b)
(b) edge node {\dG} (c)
;
\draw[1cell=.9]
(a) [rounded corners=3pt] |- ($(b)+(-1,\w)$)
-- node {\dGF} ($(b)+(1,\w)$) -| (c)
;
\end{tikzpicture}
\end{equation}
We consider the following three $\P$-functors.
\begin{romenumerate}
\item The standard enrichment of $G$ in \eqref{std-enr-F} is the $\P$-functor
\[\Gse \cn \N_G \to \P.\]
\item The standard enrichment of the composite non-symmetric multifunctor $GF$ in \cref{MNP-FG} is the $\P$-functor
\[\GFse \cn \M_{GF} \to \P.\]
\item Using \cref{gspectra-iv-enr-cat} and applying the change-of-enrichment 2-functor $\dG$ to the standard enrichment $\Fse \cn \M_F \to \N$ of $F$ in \eqref{std-enr-F} yield the $\P$-functor 
\[\Fse_G \cn \M_{GF} = (\M_F)_G \to \N_G.\]
\end{romenumerate}

\begin{explanation}[Change of Enrichment of Standard Enrichment]\label{expl:FseG}
In the context of \cref{MNP-FG,gspectra-iv-enr-cat}, the $\P$-functor 
\[\Fse_G \cn \M_{GF} = (\M_F)_G \to \N_G\]
has the same object assignment as $F \cn \M \to \N$.  For each pair of objects $x,y \in \M$, the $(x,y)$-component of $\Fse$ is the unary multimorphism \cref{Fprimexy}
\[\Fse_{x,y} = \pn{\big(F(\ev^\M_{x;\,y}) \big)} \cn F\clM(x;y) \to \clN(Fx;Fy) \inspace \N.\]
The $(x,y)$-component of $\Fse_G$ is the unary multimorphism
\begin{equation}\label{FseGxy}
(\Fse_G)_{x,y} = G\Fse_{x,y} \cn GF\clM(x;y) \to G\clN(Fx;Fy) \inspace \P
\end{equation}
obtained by applying $G$ to $\Fse_{x,y}$.  In particular, applying the non-symmetric multifunctor $G$ to the commutative diagram \cref{Fse-xy-diag} yields the following commutative diagram in $\P$.
\begin{equation}\label{GFse-partner}
\begin{tikzpicture}[vcenter]
\draw[0cell=.9]
(0,0) node (a) {\big( GF\clM(x;y) \scs GFx \big)}
(a)+(5.5,0) node (b) {\big( G\clN(Fx;Fy) \scs GFx \big)}
(b)+(0,-1.3) node (c) {GFy}
;
\draw[1cell=.9]
(a) edge node {\left( G\Fse_{x,y} \scs 1_{GFx} \right)} (b)
(b) edge node {G(\ev^\N_{Fx;\, Fy})} (c)
(a) edge node[pos=.4,swap] {GF(\ev^\M_{x;\,y})} (c)
;
\end{tikzpicture}
\end{equation}
This uses the fact that $G$, as a non-symmetric multifunctor, preserves composition and colored units.
\end{explanation}

The main observation of this section is that the $\P$-functors $\Fse_G$, $\Gse$, and $\GFse$ are related as follows.

\begin{theorem}\label{gspectra-thm-iv}\index{enrichment!standard - compositionality}\index{multifunctor!standard enrichment compositionality}\index{compositionality!- of standard enrichment}
For composable non-symmetric multifunctors between non-symmetric closed multicategories
\[\big(\M,\clM,\evM\big) \fto{F} \big(\N,\clN,\evN\big) \fto{G} \big(\P,\clP,\evP\big),\]
the following diagram of $\P$-functors commutes.
\begin{equation}\label{gspectra-iv-diagram}
\begin{tikzpicture}[vcenter]
\def\v{-1.2}
\draw[0cell]
(0,0) node (a) {\M_{GF}}
(a)+(3,0) node (b) {\P}
(a)+(0,\v) node (c) {(\M_F)_G}
(b)+(0,\v) node (d) {\N_G}
;
\draw[1cell=.9]
(a) edge node {\GFse} (b)
(a) edge[-,double equal sign distance] (c)
(c) edge node {\Fse_G} (d)
(d) edge node[swap] {\Gse} (b)
;
\end{tikzpicture}
\end{equation}
\end{theorem}

\begin{proof}
Each of the two composites in \cref{gspectra-iv-diagram} has the same object assignment as $GF$.  It remains to show that the two composites have the same component morphisms.  

For objects $x,y \in \M$, the $(x,y)$-component of $\GFse$ is, by definition \cref{Fprimexy}, the unary multimorphism
\begin{equation}\label{GFse-in-P}
\GFse_{x,y} = \pn{\big(GF(\ev^\M_{x;\,y})\big)} \cn GF\clM(x;y) \to \clP(GFx; GFy) \inspace \P.
\end{equation}
As the partner of $GF(\ev^\M_{x;\,y})$, it is \emph{uniquely} determined by the following commutative diagram in $\P$, which is \cref{Fse-xy-diag} applied to $GF$.
\begin{equation}\label{GFse-xy-diag}
\begin{tikzpicture}[vcenter]
\draw[0cell=.9]
(0,0) node (a) {\big( GF\clM(x;y) \scs GFx \big)}
(a)+(5.5,0) node (b) {\big( \clP(GFx;GFy) \scs GFx \big)}
(b)+(0,-1.3) node (c) {GFy}
;
\draw[1cell=.9]
(a) edge node {\left( \GFse_{x,y} \scs 1_{GFx} \right)} (b)
(b) edge node {\ev^\P_{GFx;\, GFy}} (c)
(a) edge node[pos=.4,swap] {GF(\ev^\M_{x;\,y})} (c)
;
\end{tikzpicture}
\end{equation}
On the other hand, the $(x,y)$-component of $\Gse \circ \Fse_G$ is the following composite in $\P$, with $G\Fse_{x,y}$ as in \cref{FseGxy}.
\begin{equation}\label{GseFseGxy}
\begin{tikzpicture}[baseline={(a.base)}]
\def\h{3.3} \def\w{.7}
\draw[0cell=.9]
(0,0) node (a) {GF\clM(x;y)}
(a)+(\h,0) node (b) {G\clN(Fx;Fy)}
(b)+(\h,0) node (c) {\clP(GFx;GFy)}
;
\draw[1cell=.9]
(a) edge node {G\Fse_{x,y}} (b)
(b) edge node {\Gse_{Fx,Fy}} (c)
;
\end{tikzpicture}
\end{equation}
We want to show that \cref{GFse-in-P,GseFseGxy} are equal.  By uniqueness of partners, it suffices to show that the composite in \cref{GseFseGxy} also makes the diagram \cref{GFse-xy-diag} commutative.

To check this, we consider the following diagram in $\P$, whose boundary is obtained from \cref{GFse-xy-diag} by replacing $\GFse_{x,y}$ with the composite in \cref{GseFseGxy}.
\begin{equation}\label{gaPevP}
\begin{tikzpicture}[vcenter]
\def\h{2.8} \def\v{1.2} \def\t{15}
\draw[0cell=.9]
(0,0) node (a) {\big( GF\clM(x;y) \scs GFx \big)}
(a)+(\h,\v) node (b) {\big( G\clN(Fx;Fy) \scs GFx \big)}
(b)+(\h,-\v) node (c) {\big( \clP(GFx;GFy) \scs GFx \big)}
(a)+(\h,-\v) node (d) {GFy}
;
\draw[1cell=.9]
(a) edge[bend left=\t] node[pos=.3] {\left( G\Fse_{x,y} \scs 1_{GFx} \right)} (b)
(b) edge[bend left=\t] node[pos=.7] {\left(\Gse_{Fx,Fy} \scs 1_{GFx} \right)} (c)
(a) edge[bend right=\t] node[pos=.1,swap] {GF(\ev^\M_{x;\,y})} (d)
(b) edge node[pos=.2] {G(\ev^\N_{Fx,Fy})} (d)
(c) edge[bend left=\t] node[pos=.1] {\ev^\P_{GFx;\, GFy}} (d)
;
\end{tikzpicture}
\end{equation}
\begin{itemize}
\item The left sub-region in \cref{gaPevP} is the commutative diagram \cref{GFse-partner}.
\item The right sub-region in \cref{gaPevP} is the commutative diagram \cref{Fse-xy-diag} for $\Gse_{Fx,Fy}$.
\end{itemize}
Thus the composite in \cref{GseFseGxy} also makes the diagram \cref{GFse-xy-diag} commutative.  This proves that the unary multimorphisms in \cref{GFse-in-P,GseFseGxy} are equal.
\end{proof}

\begin{example}[Monoidal Functors]\label{ex:stdenr-func-monfunctor}
Consider monoidal functors between symmetric monoidal closed categories
\[\V \fto{T} \W \fto{U} \X.\]
Passing to the endomorphism non-symmetric multifunctors \cref{EndP}
\[\End\,\V \fto{\End\,T} \End\,\W \fto{\End\,U} \End\,\X,\]
\cref{gspectra-thm-iv} yields the following commutative diagram of $(\End\,\X)$-functors.
\begin{equation}\label{EndUTse-factor}
\begin{tikzpicture}[vcenter]
\def\v{-1.2}
\draw[0cell]
(0,0) node (a) {(\End\,\V)_{(\End\,U)(\End\,T)}}
(a)+(5,0) node (b) {\End\,\X}
(a)+(0,\v) node (c) {\left((\End\,\V)_{\End\,T}\right)_{\End\,U}}
(b)+(0,\v) node (d) {(\End\,\W)_{\End\,U}}
;
\draw[1cell=.9]
(a) edge node {\se{(\End\,U)(\End\,T)}} (b)
(a) edge[-,double equal sign distance] (c)
(c) edge node {\se{\End\,T}_{\End\,U}} (d)
(d) edge node[swap] {\EndUse} (b)
;
\end{tikzpicture}
\end{equation}
By \cref{mon-change-enrichment,std-enr-monoidal,endtwofunctor}, the commutative diagram \cref{EndUTse-factor} is equal to the following diagram of $\X$-functors.
\begin{equation}\label{UTse-factor}
\begin{tikzpicture}[vcenter]
\def\v{-1.2}
\draw[0cell]
(0,0) node (a) {\V_{UT}}
(a)+(3,0) node (b) {\X}
(a)+(0,\v) node (c) {(\V_T)_U}
(b)+(0,\v) node (d) {\W_U}
;
\draw[1cell=.9]
(a) edge node {\UTse} (b)
(a) edge[-,double equal sign distance] (c)
(c) edge node {\Tse_U} (d)
(d) edge node[swap] {\Use} (b)
;
\end{tikzpicture}
\end{equation}
The diagram \cref{UTse-factor} uses the canonical self-enrichment, change of enrichment, and standard enrichment in \cref{theorem:v-closed-v-sm,proposition:U-std-enr,proposition:U-VCat-WCat} in the context of monoidal functors.
\end{example}

\section{Factorization of \texorpdfstring{$K$}{K}-Theory Standard Enrichment}
\label{sec:factor-Kemse}

In this section we illustrate \cref{gspectra-thm-iv} by applying it to Elmendorf-Mandell $K$-theory $\Kem$.  The result is a factorization of the standard enrichment $\Kemse$ into four spectrally enriched functors; see \cref{gspectra-thm-xi}.  For the relationship between the standard enrichment $\Kemse$ and the work of Bohmann-Osorno \cite{bohmann_osorno-mackey}, see \cref{rk:BO-Phi}.  Each $\Sp$-functor in the factorization of $\Kemse$ in \cref{gspectra-thm-xi} is either a standard enrichment functor (\cref{gspectra-thm-iii}) or the change of enrichment of a standard enrichment functor.  We explain these $\Sp$-functors further in \cref{expl:Kemse,expl:Kgse,expl:Endmse,expl:Jtse,expl:Nerstarse}.

\subsection*{Context}

First recall from \cref{Kem} that $\Kem$ factors into four multifunctors between closed multicategories as follows.  
\begin{equation}\label{Kem-four}
\begin{tikzpicture}[vcenter]
\def\v{-1.3} \def\h{6}
\draw[0cell=.9]
(0,0) node (p) {\permcatsu}
(p)+(\h,0) node (sp) {\Sp}
(p)+(0,\v) node (m) {\MoneMod}
(m)+(\h/2,0) node (gc) {\Gstarcat}
(sp)+(0,\v) node (gs) {\Gstarsset}
;
\draw[1cell=.9]
(p) edge node {\Kem} (sp)
(p) edge node[swap] {\Endm} (m)
(m) edge node {\Jt} (gc)
(gc) edge node {\Ner_*} (gs)
(gs) edge node[swap] {\Kg} (sp)
;
\end{tikzpicture}
\end{equation}

\emph{Closed Multicategories}\
\begin{itemize}
\item $\permcatsu$ is a closed multicategory by \cref{perm-closed-multicat}.
\item Each of the other four multicategories in \cref{Kem-four} is induced by a symmetric monoidal closed structure.  See
\begin{itemize}
\item \cref{proposition:EM2-5-1} \eqref{monebicomplete} for $\MoneMod$,
\item \cref{Gstar-V} for $\Gstarcat$ and $\Gstarsset$, and
\item \cref{SymSp} for $\Sp$.
\end{itemize}   
By \cref{smclosed-closed-multicat,ex:endc}, each of their endomorphism multicategories is a closed multicategory, which we denote by the same symbol.
\end{itemize}

\emph{Multifunctors}\
\begin{itemize}
\item $\Endm$ is the $\Cat$-multifunctor in \cref{expl:endm-catmulti}.
\item $\Jt$ is the symmetric monoidal $\Cat$-functor in \cref{Jt-smcat}.
\item $\Ner_*$ is the symmetric monoidal $\sSet$-functor in \cref{NerGs} induced by the nerve functor.
\item $\Kg$ is the symmetric monoidal $\sSet$-functor in \cref{Kg-sm}.
\end{itemize}
Each of the symmetric monoidal functors $\Jt$, $\Ner_*$, and $\Kg$ induces a multifunctor by the endomorphism construction \cref{EndP}.

\begin{theorem}\label{gspectra-thm-xi}
The factorization \cref{Kem-four} of multifunctors between closed multicategories
\[\Kem = \Kg \circ \Ner_* \circ\, \Jt \circ \Endm \cn \permcatsu \to \Sp\]
induces the following factorization of the standard enrichment $\Kemse$ into four $\Sp$-functors.
\begin{equation}\label{Kemse-factor-four}
\begin{tikzpicture}[vcenter]
\def\v{-1.3} \def\h{4}
\draw[0cell=.9]
(0,0) node (p) {(\permcatsu)_{\Kem}}
(p)+(\h,0) node (sp) {\Sp}
(p)+(0,\v) node (m) {(\MoneMod)_{\Kg \Ner_* \Jt}}
(m)+(\h/2,-1.2) node (gc) {(\Gstarcat)_{\Kg \Ner_*}}
(sp)+(0,\v) node (gs) {(\Gstarsset)_{\Kg}}
;
\draw[1cell=.9]
(p) edge node {\Kemse} (sp)
(p) edge node[swap] {(\Endmse)_{\Kg \Ner_* \Jt}} (m)
(m) edge[transform canvas={xshift=-1em}] node[swap,pos=.3] {(\Jtse)_{\Kg \Ner_*}\!\!} (gc)
(gc) edge[transform canvas={xshift=1em}, shorten <=-1ex] node[swap,pos=.6] {(\Nerstarse)_{\Kg}} (gs)
(gs) edge node[swap] {\Kgse} (sp)
;
\end{tikzpicture}
\end{equation}
\end{theorem}

\begin{proof}
The desired factorization of $\Kemse$ in \cref{Kemse-factor-four} is obtained by applying \cref{gspectra-thm-iv} three times to the factorization \cref{Kem-four} of $\Kem$.  More precisely, we compute as follows, where we denote some standard enrichment $\se{?}$ by $?^{\wedge}$\label{not:se} to improve readability.
\[\begin{split}
\Kemse &= \big(\Kg \circ \Ner_* \circ\, \Jt \circ \Endm\big)^\wedge\\
&= \Kgse \circ \big(\Ner_* \circ\, \Jt \circ \Endm\big)^\wedge_{\Kg}\\
&= \Kgse \circ (\Nerstarse)_{\Kg} \circ \big(\Jt \circ \Endm\big)^\wedge_{\Kg \Ner_*}\\
&= \Kgse \circ (\Nerstarse)_{\Kg} \circ (\Jtse)_{\Kg \Ner_*} \circ (\Endmse)_{\Kg \Ner_* \Jt}
\end{split}\]
For the last two equalities above, we also use the compositionality of change-of-enrichment 2-functors in \cref{func-change-enr}.
\end{proof}

We explain the $\Sp$-functors in \cref{Kemse-factor-four} in more detail after the following remark.

\begin{remark}[Work of Bohmann-Osorno]\label{rk:BO-Phi}
In \cite[Theorem 6.2]{bohmann_osorno-mackey} there is a $\Sp$-functor $\Phi$\label{not:BOPhi} that is categorically similar to the standard enrichment $\Kemse$ in \cref{Kemse-factor-four}.  The important difference is that, while $\Kemse$ is constructed from Elmendorf-Mandell $K$-theory $\Kem$, the Bohmann-Osorno $\Sp$-functor $\Phi$ is the standard enrichment of the $K$-theory non-symmetric multifunctor $\bbK$\label{not:KGM} in the Guillou-May \cref{theorem:GM} \cite{guillou_may,gmmo}.  As far as the authors know, there is no known multiplicative comparison between $\Kem$ and $\bbK$.  Thus we also do not know how $\Kemse$ is related to $\Phi$.
\end{remark}

The rest of this section explains the $\Sp$-functors in \cref{Kemse-factor-four} in more detail.  We use the shortened notation
\[\psu = \permcatsu \andspace \clpsu = \clp.\]

\begin{explanation}[The $\Sp$-functor $\Kemse$]\label{expl:Kemse}
  Specifying to $\Kem$ \cref{Kem}, \cref{gspectra-thm-iii} says that the \index{standard enrichment!- of Elmendorf-Mandell $K$-theory}\index{Elmendorf-Mandell!K-theory@$K$-theory!standard enrichment}\index{K-theory@$K$-theory!Elmendorf-Mandell!standard enrichment}standard enrichment of $\Kem$ is the $\Sp$-functor
\begin{equation}\label{Kemse}
\Kemse \cn (\psu)_{\Kem} \to \Sp.
\end{equation}
Next we describe its object assignment and component morphisms.

\medskip
\emph{Object Assignment}.
The standard enrichment $\Kemse$ has the same object assignment as $\Kem$.  In other words, it sends each small permutative category $\C$ to the connective symmetric spectrum $\Kem \C$.

\medskip
\emph{Components}. 
For each pair of small permutative categories $\C$ and $\D$, by definition \cref{Fprimexy,Fev}, the component morphism 
\[\Kemse_{\C,\D} \cn \Kem \clpsu(\C; \D) \to \clsp(\Kem\C; \Kem\D) \inspace \Sp\]
is the adjoint of the morphism
\[\Kem(\ev_{\C;\,\D}) \cn \Kem\clpsu(\C;\D) \sma \Kem\C \to \Kem\D.\]
\begin{itemize}
\item $\clpsu(\C;\D)$ is the small permutative category in \cref{clp-angcd-permutative}.  Since the domain has length 1, $\clpsu(\C;\D)$ is equal to the hom object $\psu(\C,\D)$ in \cref{def:perm-selfenr}.
\item The evaluation is the bilinear functor \cref{clp-ev-angcd}
\[\ev_{\C;\,\D} \cn \clpsu(\C; \D) \times \C \to \D.\]
This is equal to the evaluation $\ev_{\C,\D}$ in \cref{evCD}.
\item $\sma$ is the smash product of symmetric spectra \cite[7.6.1]{cerberusIII}.
\end{itemize}
By uniqueness of adjoints, $\Kemse_{\C,\D}$ is characterized by the following commutative diagram in $\Sp$.
\begin{equation}\label{Kemse-diag}
\begin{tikzpicture}[vcenter]
\draw[0cell=.9]
(0,0) node (a) {\Kem\clpsu(\C;\D) \sma \Kem\C}
(a)+(3.8,0) node (b) {\phantom{\Kem\D}}
(b)+(1.3,0) node (b') {\clsp(\Kem\C; \Kem\D) \sma \Kem\C}
(b)+(0,-1.3) node (c) {\Kem\D}
;
\draw[1cell=.9]
(a) edge node {\Kemse_{\C,\D} \sma 1} (b)
(b) edge[shorten <=.5ex] node {\ev} (c)
(a) edge node[pos=.4,swap] {\Kem(\ev)} (c)
;
\end{tikzpicture}
\end{equation}
This is the commutative diagram \cref{Fse-xy-diag} for $F = \Kem$.
\end{explanation}

\begin{explanation}[The $\Sp$-Functor $\Kgse$]\label{expl:Kgse}
In the diagram \cref{Kemse-factor-four}, the $\Sp$-functor 
\[\Kgse \cn (\Gstarsset)_{\Kg} \to \Sp\]
is the standard enrichment of the symmetric monoidal functor in \cref{Kg-sm},
\[\Kg \cn \Gstarsset \to \Sp.\]
The standard enrichment $\Sp$-functor $\Kgse$ exists by either \cref{proposition:U-std-enr} or \cref{gspectra-thm-iii}, which yield the same $\Sp$-functor by \cref{std-enr-monoidal}.  
\end{explanation}

The other three constituent $\Sp$-functors in \cref{Kemse-factor-four} are the following.
\[\begin{split}
(\Endmse)_{\Kg \Ner_* \Jt} & \cn (\psu)_{\Kem} \to (\MoneMod)_{\Kg \Ner_* \Jt}\\
(\Jtse)_{\Kg \Ner_*} & \cn (\MoneMod)_{\Kg \Ner_* \Jt} \to (\Gstarcat)_{\Kg \Ner_*}\\
(\Nerstarse)_{\Kg} & \cn (\Gstarcat)_{\Kg \Ner_*} \to (\Gstarsset)_{\Kg}
\end{split}\]
Each of these three $\Sp$-functors is obtained from the indicated standard enrichment $\se{?}$ by applying the change of enrichment in the subscript, as in \cref{expl:FseG}.  We describe them more explicitly in \cref{expl:Endmse,expl:Jtse,expl:Nerstarse} below.

\begin{explanation}[The $\Sp$-Functor $(\Endmse)_{\Kg \Ner_* \Jt}$]\label{expl:Endmse}
We obtain an explicit description of the $\Sp$-functor 
\[(\Endmse)_{\Kg \Ner_* \Jt} \cn (\psu)_{\Kem} \to (\MoneMod)_{\Kg \Ner_* \Jt}\]
in \cref{Kemse-factor-four} by interpreting \cref{expl:FseG} with the multifunctors between closed multicategories
\[\begin{split}
F = \Endm & \cn \psu \to \MoneMod \andspace\\ 
G = \Kg \Ner_* \Jt & \cn \MoneMod \to \Sp
\end{split}\]
in \cref{expl:endm-catmulti,eq:Ksummary}, respectively.  Note that $\Kem = GF$ by definition \cref{Kem}.

\medskip
\emph{Object Assignment}.
$(\Endmse)_{\Kg \Ner_* \Jt}$ sends each small permutative category $\C$ to 
\[(\Endmse)_{\Kg \Ner_* \Jt} (\C) = \Endm\C \inspace \MoneMod.\]
This is the endomorphism left $\Mone$-module in \cref{ex:endmc}.

\medskip
\emph{Standard Enrichment of $\Endm$}.
For the component morphisms, we first consider the standard enrichment $\MoneMod$-functor of $\Endm$ (\cref{gspectra-thm-iii}),
\begin{equation}\label{Endm-stdenr}
\Endmse \cn (\psu)_{\Endm} \to \MoneMod.
\end{equation}
For small permutative categories $\C$ and $\D$, the $(\C,\D)$-component of $\Endmse$ is the following morphism in $\MoneMod$.
\begin{equation}\label{Endm-ev-partner}
\scalebox{.9}{$\pn{\big(\Endm(\ev_{\C;\, \D})\big)}$} \cn 
\scalebox{.9}{$\Endm \clpsu(\C;\D)$} \to 
\scalebox{.9}{$\pHom\big(\Endm\C; \Endm\D\big)$}
\end{equation}
This is adjoint to the following morphism in $\MoneMod$.
\[\Endm(\ev_{\C;\, \D}) \cn \Endm \clpsu(\C;\D) \sma \Endm\C \to \Endm\D\]
Here 
\[\ev_{\C;\, \D} \cn \clpsu(\C;\D) \times \C \to \D\]
is the evaluation bilinear functor in \cref{clp-ev-angcd}, which is the same as $\ev_{\C,\D}$ in \cref{evCD}.  The smash product, $\sma$, and the internal hom, $\pHom$, are part of the symmetric monoidal closed structure of $\MoneMod$ in \cref{proposition:EM2-5-1} \eqref{monebicomplete}.

\medskip
\emph{Components of $(\Endmse)_{\Kg \Ner_* \Jt}$}.
The composite symmetric monoidal functor $\Kg \Ner_* \Jt$ induces a change-of-enrichment 2-functor
\begin{equation}\label{KNJ-change-enr}
\dsub{\Kg \Ner_* \Jt} \cn \Monemodcat \to \Spcat.
\end{equation}
This change of enrichment exists by either \cref{proposition:U-VCat-WCat} or \cref{mult-change-enrichment}, which yield the same 2-functor by \cref{mon-change-enrichment}.

Applying the change of enrichment \cref{KNJ-change-enr} to the standard enrichment $\Endmse$ in \cref{Endm-stdenr} yields the $\Sp$-functor in the upper left of the diagram \cref{Kemse-factor-four}:
\begin{equation}\label{EndmseKNJ}
\begin{tikzpicture}[vcenter]
\def\v{-1.3} \def\h{-3}
\draw[0cell=.9]
(0,0) node (p) {\big((\psu)_{\Endm}\big)_{\Kg \Ner_* \Jt}}
(p)+(\h,0) node (p') {(\psu)_{\Kem}}
(p)+(0,\v) node (m) {(\MoneMod)_{\Kg \Ner_* \Jt}}
;
\draw[1cell=.9]
(p) edge node[swap] {(\Endmse)_{\Kg \Ner_* \Jt}} (m)
(p) edge[-,double equal sign distance] (p')
;
\end{tikzpicture}
\end{equation}
The equality of $\Sp$-categories at the top of \cref{EndmseKNJ} follows from \cref{func-change-enr} applied to the factorization \cref{Kem-four} of $\Kem$.  For small permutative categories $\C$ and $\D$, the $(\C,\D)$-component of $(\Endmse)_{\Kg \Ner_* \Jt}$ is the following morphism in $\Sp$.
\[\begin{tikzpicture}
\draw[0cell=.9]
(0,0) node (a) {\Kg \Ner_* \Jt \Endm \clpsu(\C;\D)}
(a)+(-3.7,0) node (a') {\Kem \clpsu(\C;\D)}
(a)+(0,-1.3) node (b) {\Kg \Ner_* \Jt \pHom\big(\Endm\C; \Endm\D\big)}
;
\draw[1cell=.9]
(a) edge node[swap] {\Kg \Ner_* \Jt \pn{\big(\Endm(\ev_{\C;\, \D})\big)}} (b)
(a) edge[-,double equal sign distance] (a')
;
\end{tikzpicture}\]
This is obtained from the morphism $\big(\Endmse\big)_{\C,\D}$ in \cref{Endm-ev-partner} by applying the functor $\Kg \Ner_* \Jt$.
\end{explanation}

\begin{explanation}[The $\Sp$-Functor $(\Jtse)_{\Kg \Ner_*}$]\label{expl:Jtse}
The symmetric monoidal functor in \cref{Jt-smcat},
\[\Jt \cn \MoneMod \to \Gstarcat,\]
has a standard enrichment $(\Gstarcat)$-functor
\[\Jtse \cn (\MoneMod)_{\Jt} \to \Gstarcat\]
by \cref{proposition:U-std-enr,gspectra-thm-iii,std-enr-monoidal}.  The symmetric monoidal functor
\[\Kg \Ner_* \cn \Gstarcat \to \Gstarsset \to \Sp\]
induces a change-of-enrichment 2-functor
\[\dsub{\Kg \Ner_*} \cn \Gstarcatcat \to \Spcat\]
by \cref{proposition:U-VCat-WCat,mult-change-enrichment,mon-change-enrichment}.  Applying this change of enrichment to the standard enrichment $\Jtse$ yields the following $\Sp$-functor in \cref{Kemse-factor-four}.
\[\begin{tikzpicture}[vcenter]
\def\v{-1.3} \def\h{-3.3}
\draw[0cell=.9]
(0,0) node (m) {\big((\MoneMod)_{\Jt}\big)_{\Kg \Ner_*}}
(m)+(\h,0) node (m') {(\MoneMod)_{\Kg \Ner_* \Jt}}
(m)+(4.5,0) node (gc) {(\Gstarcat)_{\Kg \Ner_*}}
;
\draw[1cell=.9]
(m) edge node {(\Jtse)_{\Kg \Ner_*}} (gc)
(m) edge[-,double equal sign distance] (m')
;
\end{tikzpicture}\]
On objects it sends each left $\Mone$-module $\N$ to 
\[(\Jtse)_{\Kg \Ner_*}(\N) = \Jt \N \inspace \Gstarcat.\]

For left $\Mone$-modules $\N$ and $\P$, the $(\N,\P)$-component morphism in $\Sp$
\begin{equation}\label{JtseKgN}
\begin{tikzpicture}[baseline={(m.base)}]
\def\v{-1.3} \def\h{6.2}
\draw[0cell=.9]
(0,0) node (m) {\Kg \Ner_* \Jt \pHom(\N,\P)}
(m)+(\h,0) node (gc) {\Kg \Ner_* \Homgstar\left(\Jt\N, \Jt\P\right)}
;
\draw[1cell=.9]
(m) edge node {\big((\Jtse)_{\Kg \Ner_*}\big)_{\N,\P}} (gc)
;
\end{tikzpicture}
\end{equation}
is obtained by applying the functor $\Kg\Ner_*$ to the adjoint---taken in the symmetric monoidal closed category $\big(\Gstarcat, \sma, \Homgstar\!\big)$ in \cref{Gstar-V}---of the following composite morphism.
\begin{equation}\label{JttwoJtev}
\begin{tikzpicture}[baseline={(a.base)}]
\draw[0cell=.9]
(0,0) node (a) {\Jt\pHom(\N,\P) \sma \Jt\N}
(a)+(4.3,0) node (b) {\Jt\big(\pHom(\N,\P) \sma \N\big)}
(b)+(3.3,0) node (c) {\Jt\P}
;
\draw[1cell=.9]
(a) edge node {(\Jt)^2} (b)
(b) edge node {\Jt(\ev)} (c)
;
\end{tikzpicture}
\end{equation}
\begin{itemize}
\item In \cref{JtseKgN,JttwoJtev}, $\left(\MoneMod,\sma,\pHom\right)$ is the symmetric monoidal closed structure on $\MoneMod$ in \cref{proposition:EM2-5-1} \eqref{monebicomplete}.
\item $(\Jt)^2$ is the monoidal constraint of the symmetric monoidal functor $\Jt$.  See \cite[10.3.11]{cerberusIII} for a detailed description.
\item $\ev$ is the evaluation \cref{evaluation} in $(\MoneMod,\sma,\pHom)$.
\item The first $\sma$ in \cref{JttwoJtev} and $\Homgstar$ in \cref{JtseKgN} are the pointed Day convolution and pointed hom for $\Gstar$-categories in \cref{eq:GstarV-convolution,eq:GstarV-hom}.\defmark
\end{itemize}
\end{explanation}

\begin{explanation}[The $\Sp$-Functor $(\Nerstarse)_{\Kg}$]\label{expl:Nerstarse}
The symmetric monoidal functor in \cref{NerGs},
\[\Ner_* \cn \Gstarcat \to \Gstarsset,\]
is induced levelwise by the nerve functor, $\Ner$.  Its standard enrichment $(\Gstarsset)$-functor
\[\Nerstarse \cn (\Gstarcat)_{\Ner_*} \to \Gstarsset\]
exists by \cref{proposition:U-std-enr,gspectra-thm-iii,std-enr-monoidal}.  The symmetric monoidal functor
\[\Kg \cn \Gstarsset \to \Sp\]
in \cref{Kg-sm} induces a change-of-enrichment 2-functor
\[\dsub{\Kg} \cn \Gstarssetcat \to \Spcat\]
by \cref{proposition:U-VCat-WCat,mult-change-enrichment,mon-change-enrichment}.  Applying this change of enrichment to the standard enrichment $\Nerstarse$ yields the following $\Sp$-functor in \cref{Kemse-factor-four}.
\[\begin{tikzpicture}[vcenter]
\def\v{-1.3} \def\h{-3}
\draw[0cell=.9]
(0,0) node (gc) {\big((\Gstarcat)_{\Ner_*}\big)_{\Kg}}
(gc)+(\h,0) node (gc') {(\Gstarcat)_{\Kg \Ner_*}}
(gc)+(4,0) node (gs) {(\Gstarsset)_{\Kg}}
;
\draw[1cell=.9]
(gc) edge node {(\Nerstarse)_{\Kg}} (gs)
(gc) edge[-,double equal sign distance] (gc')
;
\end{tikzpicture}\]
On objects it sends each $\Gstar$-category $A$ to 
\[(\Nerstarse)_{\Kg}(A) = \Ner_* A \inspace \Gstarsset.\]

For $\Gstar$-categories $A$ and $B$, the $(A,B)$-component morphism in $\Sp$
\begin{equation}\label{NseAB}
\begin{tikzpicture}[baseline={(a.base)}]
\def\v{-1.3} \def\h{6}
\draw[0cell=.9]
(0,0) node (a) {\Kg \Ner_* \Homgstar(A,B)}
(a)+(\h,0) node (b) {\Kg \Homgstar\left(\Ner_*A , \Ner_*B\right)}
;
\draw[1cell=.9]
(a) edge node {\big((\Nerstarse)_{\Kg}\big)_{A,B}} (b)
;
\end{tikzpicture}
\end{equation}
is obtained by applying the functor $\Kg$ to the adjoint---taken in the symmetric monoidal closed category $\big(\Gstarsset, \sma, \Homgstar\!\big)$ in \cref{Gstar-V}---of the following composite morphism.
\begin{equation}\label{NtwoNev}
\begin{tikzpicture}[baseline={(a.base)}]
\draw[0cell=.85]
(0,0) node (a) {\Ner_*\Homgstar(A,B) \sma \Ner_*A}
(a)+(4.4,0) node (b) {\Ner_*\big(\Homgstar(A,B) \sma A\big)}
(b)+(3.4,0) node (c) {\Ner_*B}
;
\draw[1cell=.8]
(a) edge node {\Ner_*^2} (b)
(b) edge node {\Ner_*(\ev)} (c)
;
\end{tikzpicture}
\end{equation}
\begin{itemize}
\item In \cref{NseAB} $\Homgstar$ in the domain is the pointed hom in $\Gstarcat$ \cref{eq:GstarV-hom}.  The $\Homgstar$ in the codomain is the one in $\Gstarsset$.
\item In \cref{NtwoNev} $\Ner_*^2$ is the monoidal constraint of the symmetric monoidal functor $\Ner_*$.  See \cite[3.8.4]{cerberusIII} for a detailed description.
\item $\ev$ is the evaluation \cref{evaluation} in $(\Gstarcat,\sma,\Homgstar\!)$.\defmark
\end{itemize}
\end{explanation}

\chapter[Enriched Mackey Functors of Closed Multicategories]{Enriched Diagrams and Mackey Functors of Closed Multicategories}
\label{ch:gspectra_Kem}
This chapter studies categories of enriched functors 
\[
  \MCat(\C,\M) \andspace \MCat(\Cop,\M),
\]
where $\C$ is an $\M$-category.
At left, $\M$ is assumed to be a non-symmetric closed multicategory.
At right, $\M$ is assumed to have the additional symmetric structure of a closed multicategory, so that $\Cop$ is defined (\cref{def:opposite-mcat}).

For a non-symmetric multifunctor between non-symmetric closed multicategories
\[
  F \cn \M \to \N,
\]
the change-of-enrichment 2-functors from \cref{mult-change-enrichment} induce diagram change of enrichment (\cref{gspectra-thm-v})
\[
  \Fdg \cn \MCat(\C,\M) \to \NCat\big( \C_F,\N \big).
\]
When $F$ is a multifunctor between closed multicategories, there is a presheaf change of enrichment (\cref{dF-opposite,gspectra-thm-v-cor})
\[
  \Fdg \cn \MCat(\Cop,\M) \to \NCat\big( (\C_F)^\op,\N \big).
\]
These results apply to the following diagram of $K$-theory multifunctors from
\cref{eq:Ksummary}; see \cref{ex:gsp-v-examples}.
\begin{equation}\label{intro-Kem-functors}
\begin{tikzpicture}[xscale=1.2,yscale=1,vcenter]
\def\t{.35} \def\v{-1.5} \def\u{.6} \def\w{-.6}
\draw[0cell=.8]
(0,0) node (a1) {\permcatsu}
(a1)+(3,0) node (a2) {\Gacat}
(a2)+(2,0) node (a3) {\Gasset}
(a3)+(2,0) node (a4) {\Sp}
(a2)+(0,\v) node (b2) {\Gstarcat}
(a3)+(0,\v) node (b3) {\Gstarsset}
(b2)+(-2.3,0) node (b1) {\MoneMod}
;
\draw[1cell=.8]
(a2) edge node {\Ner_*} (a3)
(a3) edge node {\Kf} (a4)
(a1) edge node[pos=.3] {\Jem} (b2)
(a1) edge node[pos=.85] {\Endm} (b1)
(b1) edge node[pos=.4] {\Jt} (b2)
(b2) edge node {\Ner_*} (b3)
(b3) edge node [pos=.4,swap] {\Kg} (a4)
(a2) edge[transform canvas={xshift=0ex}] node[swap,pos=.3] {\sma^*} (b2)
(a3) edge[transform canvas={xshift=0ex}] node[pos=.3] {\sma^*} (b3)
;
\draw[1cell=.85]
(a1) [rounded corners=3pt] |- ($(b2)+(0,\w)$)
-- node[pos=\t] {\Kem} ($(b3)+(0,\w)$) -| (a4)
;
\end{tikzpicture}
\end{equation}
Each arrow in this diagram is a multifunctor between closed multicategories.

The compositionality of change-of-diagram and change-of-presheaf involves that of both the change of enrichment (\cref{sec:functoriality-change-enr}) and the standard enrichment (\cref{sec:fun-std-enr-multi}).
The resulting theory is explained in \cref{sec:fun-change-enr-diag}.

The main application in this chapter is \cref{gspectra-thm-xiv}, which shows that
the factorization of $\Kem$ induces the following factorization of the diagram change-of-enrichment functor $\Kemdg$, where $\C$ is a small $\permcatsu$-category.
\begin{equation}\label{intro-Kemdg-factor-four}
\begin{tikzpicture}[vcenter]
\def\v{-1.3} \def\h{5}
\draw[0cell=.8]
(0,0) node (p) {\permcatsucat\big(\C, \permcatsu\big)}
(p)+(\h,0) node (sp) {\Spcat\big(\C_{\Kem}, \Sp\big)}
(p)+(0,\v) node (m) {\Monemodcat\Big(\C_{\Endm} \scs \MoneMod \Big)}
(m)+(\h/2,-1.2) node (gc) {\Gstarcatcat\Big(\C_{\Jt\Endm} \scs \Gstarcat \Big)}
(sp)+(0,\v) node (gs) {\Gstarssetcat\Big(\C_{\Ner_* \Jt\Endm} \scs \Gstarsset \Big)}
;
\draw[1cell=.9]
(p) edge node {\Kemdg} (sp)
(p) edge node[swap] {\Endmdg} (m)
(m) edge[transform canvas={xshift=-1.5em}, shorten >=-1ex] node[swap,pos=.2] {\Jtdg} (gc)
(gc) edge[transform canvas={xshift=1.5em}, shorten <=-1ex] node[swap,pos=.8] {\Nerstardg} (gs)
(gs) edge node[swap] {\Kgdg} (sp)
;
\end{tikzpicture}
\end{equation}
A similar factorization holds with $\C$ and each $\C_?$ above replaced by $\Cop$ and $(\C_?)^\op$, respectively. 

\subsection*{Connection with Other Chapters}

\cref{ch:mackey} develops homotopical properties of the enriched diagram and Mackey functor change-of-enrichment functors $\Fdg$.
\cref{ch:mackey_eq} then gives the corresponding applications to $K$-theory multifunctors.

\subsection*{Background}

The material in this chapter is a culmination of the multifunctorial enrichment and change of enrichment from \cref{ch:menriched,ch:change_enr,ch:gspectra,ch:std_enrich}.
In particular, the applications to $K$-theory in \cref{sec:presheaf-K,sec:mult-mackey-spectra} depend on the descriptions in \cref{sec:factor-Kemse}.

\subsection*{Chapter Summary}

\cref{sec:enr-diag-psh} gives the basic definitions of enriched diagrams and enriched Mackey functors.
\cref{sec:enr-diag-change-enr} defines change-of-enrichment functors for enriched diagrams.
\cref{sec:diag-psh-change-enr-functors} contains the proofs that these are functorial, along with key $K$-theoretic applications.
In \cref{sec:fun-change-enr-diag} we show that the diagram change-of-enrichment construction respects composition.
\cref{sec:presheaf-K} applies the general results from \cref{sec:diag-psh-change-enr-functors} to the Elmendorf-Mandell $K$-theory functor $\Kem$.
\cref{sec:mult-mackey-spectra} applies the general results from \cref{sec:fun-change-enr-diag} to factor the change of enrichment given by $\Kem$.
Here is a summary table.
\reftable{.9}{
  enriched diagrams and enriched Mackey functors
  & \ref{def:enr-diag-cat}
  \\ \hline
  diagram change of enrichment definition
  & \ref{def:enr-diag-change-enr}, \ref{expl:Fdg-object}, and \ref{expl:Fdg-morphism}
  \\ \hline
  diagram change of enrichment functoriality
  & \ref{gspectra-thm-v} and \ref{gspectra-thm-v-cor}
  \\ \hline
  diagram change of enrichment compositions
  & \ref{gspectra-thm-vii} and \ref{gspectra-thm-vii-cor}
  \\ \hline
  application to Elmendorf-Mandell $K$-theory
  &   \ref{thm:Kemdg}, \ref{expl:Kemdg-object}, and \ref{expl:Kemdg-morphism}
  \\ \hline
  factorization of Elmendorf-Mandell change of enrichment
  &  \ref{gspectra-thm-xiv} and \ref{expl:Monediagrams}
  \\
}


\section{Enriched Diagrams and Mackey Functors as Modules}
\label{sec:enr-diag-psh}

In this section we introduce categories of enriched diagrams and enriched presheaves, also called Mackey functors, with respect to a (non-symmetric) closed multicategory.
\begin{itemize}
\item Diagrams and presheaves enriched in a (non-symmetric) closed multicategory are in \cref{def:enr-diag-cat}.
\item We characterize the objects, morphisms, and composition in these categories in terms of partners in \cref{C-diagram-partner,C-diag-morphism-pn,mcat-vertical-comp}.  The upshot is that we may consider these categories as categories of modules over an enriched category. 
\item In \cref{ex:schwede-shipley,ex:guillou-may} we discuss enriched presheaves in the work of
\begin{itemize}
\item \cite{schwede-shipley_stable} about stable model categories and
\item \cite{guillou_may} about genuine equivariant spectra.
\end{itemize}   
\end{itemize}  
We discuss change of enrichment of enriched diagram and presheaf categories in the remaining sections of this chapter.

\subsection*{Defining Enriched Diagrams and Mackey Functors}

\begin{definition}\label{def:enr-diag-cat}
Suppose $(\M,\clM,\ev)$ is a non-symmetric closed multicategory (\cref{def:closed-multicat}).  We also regard $\M$ as an $\M$-category with the canonical self-enrichment $(\M,\comp,i)$ (\cref{cl-multi-cl-cat}).  Suppose $(\C,\mcomp,i)$ is an $\M$-category (\cref{def:menriched-cat}).
\begin{itemize}
\item The category
\begin{equation}\label{mcat-cm}
\MCat(\C,\M)
\end{equation}
is called the \emph{$\C$-diagram category}\index{enriched!diagram}\index{diagram!enriched} of $\M$.  An object in this category is called a \emph{$\C$-diagram} enriched in $\M$. 
\item Suppose, in addition, $\M$ is a closed multicategory, and $\Cop$ is the opposite $\M$-category (\cref{opposite-mcat}).  The $\Cop$-diagram category
\begin{equation}\label{mcat-copm}
\MCat(\Cop,\M)
\end{equation}
is also called the \emph{$\C$-presheaf category}\index{enriched!presheaf}\index{presheaf!enriched} of $\M$ and the \emph{$\C$-Mackey functor category}\index{enriched!Mackey functor}\index{Mackey functor!enriched} of $\M$.  An object in this category is also called a \emph{$\C$-presheaf} and a \emph{$\C$-Mackey functor} enriched in $\M$.
\end{itemize}  
This finishes the definition.
\end{definition}

Using \cref{opposite-mcat,dF-opposite}, the discussion below about $\C$-diagrams also applies to $\C$-Mackey functors.

\begin{explanation}[Size]\label{expl:enr-diag-cat-small}\index{Grothendieck Universe}\index{universe}\index{Axiom of Universes}\index{convention!universe}
To define the category $\MCat(\C,\M)$, technically we need $\C$ and $\M$ to be small.  We can deal with this issue in one of two ways.
\begin{enumerate}
\item As in \cref{conv:universe}, if necessary we can move to a larger universe where $\C$ and $\M$ are small.
\item We observe that all of our proofs and assertions are about $\M$-functors $\C \to \M$ and $\M$-natural transformations between them.  These notions are defined without assuming that $\C$ and $\M$ are small.  We refer to the category $\MCat(\C,\M)$ simply because it provides a convenient context to phrase functorial and naturality properties.  Thus, $\C$ and $\M$ do not need to be small.\defmark
\end{enumerate}
\end{explanation}

\subsection*{Diagrams as Modules}

Recall that in a (non-symmetric) closed multicategory, taking partner, denoted $f \mapsto \pn{f}$, is a bijection \cref{eval-bijection}.  In the rest of this section, we characterize the objects, morphisms, and composition in the category $\MCat(\C,\M)$ in terms of partners.

\begin{explanation}[Unpacking $\C$-Diagrams]\label{expl:enr-diag-cat}
In \cref{mcat-cm} $\MCat(\C,\M)$ is a hom category in the 2-category $\MCat$ (\cref{mcat-iicat}).  An object in $\MCat(\C,\M)$ is an $\M$-functor $A \cn \C \to \M$ (\cref{def:mfunctor}).  Such an $\M$-functor has an object assignment
\begin{equation}\label{A-obj-assignment}
A \cn \Ob\C \to \Ob\M.
\end{equation}
For objects $x,y \in \C$, its $(x,y)$-component is a unary multimorphism
\[A_{x,y} \cn \C(x,y) \to \clM(Ax;Ay) \inspace \M.\]
Its partner \cref{eval-bijection} is a binary multimorphism
\begin{equation}\label{Axy-pn-binary}
\pn{A_{x,y}} \cn \big(\C(x,y) \scs Ax\big) \to Ay \inspace \M.
\end{equation}
In \cref{C-diagram-partner} below, we interpret the $\M$-functor axioms \cref{mfunctor-diagrams} for $A$ in terms of these componentwise partners.
\end{explanation}

The following observation allows us to regard a $\C$-diagram in $\M$ as a left $\C$-module.  

\begin{proposition}\label{C-diagram-partner}
In the context of \cref{def:enr-diag-cat}, an $\M$-functor $A \cn \C \to \M$ is uniquely determined by 
\begin{itemize}
\item an object assignment as in \cref{A-obj-assignment} and
\item component binary multimorphisms $\big\{\pn{A_{x,y}}\big\}_{x,y\in\C}$ as in \cref{Axy-pn-binary}
\end{itemize}
such that the following two diagrams in $\M$ commute for all objects $x,y,z \in \C$, with $\C(x,y)$ abbreviated to $\C_{x,y}$.
\begin{equation}\label{c-diag-pn}

\end{equation}
\end{proposition}

\begin{proof}
Since taking partner is a bijection \cref{eval-bijection}, it suffices to show that the diagrams in \cref{c-diag-pn} are the partners of the diagrams in \cref{mfunctor-diagrams}, which are the axioms for an $\M$-functor $\C \to \M$.  We abbreviate $\clM(?;?)$ and $\ev^\M$ to $\clM_{?;?}$ and $\ev$, respectively.  For objects $x,y \in \C$, the component $A_{x,y}$ and its partner $\pn{A_{x,y}}$ determine each other via the following commutative diagram in $\M$.
\begin{equation}\label{Axy-Axypn-diagram}
%
\end{equation}

To see that the left diagram in \cref{mfunctor-diagrams} yields the left diagram in \cref{c-diag-pn} upon taking partners, we consider the following diagram in $\M$.
\begin{equation}\label{ACM-axiom-parnter}
%
\end{equation}
The sub-region labeled $\pentagram$, when composed with the lower right $\ev$, yields the partners \cref{eval-bijection} of the two composites in the left diagram in \cref{mfunctor-diagrams} for $A \cn \C \to \M$.  The boundary of \cref{ACM-axiom-parnter} is the left diagram in \cref{c-diag-pn}.  Thus it suffices to show that the other four sub-regions in \cref{ACM-axiom-parnter} commute.
\begin{itemize}
\item The middle diamond region commutes by the definition of $\comp$ in the canonical self-enrichment of $\M$ \cref{selfenrm-diagrams}.
\item The other three sub-regions commute by \cref{Axy-Axypn-diagram}.
\end{itemize}

To see that the right diagram in \cref{mfunctor-diagrams} yields the right diagram in \cref{c-diag-pn} upon taking partners, we consider the following diagram in $\M$.
\begin{equation}\label{ACM-axiom-ii-pn}
\begin{tikzpicture}[vcenter]
\def\v{-1.3} \def\g{1.5}
\draw[0cell=.85]
(0,0) node (a1) {\big(\ang{} \scs Ax\big)}
(a1)+(4,0) node (a2) {\big(\C_{x,x} \scs Ax\big)}
(a2)+(0,\v) node (b) {\big(\clM_{Ax;Ax} \scs Ax\big)}
(b)+(0,\v) node (c) {Ax}
;
\draw[1cell=.85]
(a1) edge node[swap,pos=.7] {(i_{Ax} \scs \opu)} (b)
(a1) edge node {(i_x \scs \opu)} (a2)
(a2) edge[shorten >=-.5ex] node[swap,pos=.4] {(A_{x,x} \scs \opu)} (b)
(b) edge node[swap] {\ev} (c)
;
\draw[1cell=.85]
(a1) edge[out=-90,in=180] node[swap,pos=.3] {\opu} (c)
;
\draw[1cell=.85]
(a2) [rounded corners=3pt] -| ($(b)+(\g,1)$)
-- node[swap,pos=.2] {\pn{A_{x,x}}} ($(b)+(\g,-1)$) |- (c)
;
\end{tikzpicture}
\end{equation}
The top triangle, when composed with $\ev$, yields the partners of the two composites in the right diagram in \cref{mfunctor-diagrams} for $A \cn \C \to \M$.  The boundary of \cref{ACM-axiom-ii-pn} is the right diagram in \cref{c-diag-pn}.  Thus it suffices to show that the other two sub-regions in \cref{ACM-axiom-ii-pn} commute.
\begin{itemize}
\item The left sub-region commutes by the definition of $i_{Ax}$ in the canonical self-enrichment of $\M$ \cref{selfenrm-diagrams}.
\item The right sub-region commutes by \cref{Axy-Axypn-diagram}.
\end{itemize}
This finishes the proof.
\end{proof}

\subsection*{Mackey Functors as Modules}

Recall that each category $\C$ enriched in a multicategory has an opposite $(\Cop,\mcompop,i)$ (\cref{opposite-mcat}).  The following observation is \cref{C-diagram-partner} applied to $\Cop$.

\begin{proposition}\label{C-presheaf-partner}
In the context of \cref{mcat-copm}, an $\M$-functor $A \cn \Cop \to \M$ is uniquely determined by 
\begin{itemize}
\item an object assignment $A \cn \Ob\C \to \Ob\M$ and
\item for each pair of objects $x,y \in \C$, a component binary multimorphism
\[\pn{A_{x,y}} \cn \big(\C_{y,x} \scs Ax\big) \to Ay \inspace \M\]
\end{itemize}
such that the following two diagrams in $\M$ commute for all objects $x,y,z \in \C$, with $\tau \in \Sigma_2$ the nonidentity permutation. 
\begin{equation}\label{c-presheaf-pn}
\begin{tikzpicture}[vcenter]
\def\v{-1.4} \def\h{3} \def\t{20}
\draw[0cell=.85]
(0,0) node (a) {\big(\C_{z,y} \scs \C_{y,x} \scs Ax\big)}
(a)+(\h,0) node (b) {\big(\C_{z,x} \scs Ax\big)}
(a)+(\h/2,1.2) node (e) {\big(\C_{y,x} \scs \C_{z,y} \scs Ax\big)}
(a)+(0,\v) node (c) {\big(\C_{z,y} \scs Ay\big)}
(b)+(0,\v) node (d) {Az}
;
\draw[1cell=.85]
(a) edge[bend left=\t,transform canvas={xshift=-.8ex}] node[pos=.3] {(\tau \scs \opu)} (e)
(e) edge[bend left=\t,transform canvas={xshift=.8ex}] node[pos=.7] {(\mcomp \scs \opu)} (b)
(b) edge node {\pn{A_{x,z}}} (d)
(a) edge node[swap] {(\opu \scs \pn{A_{x,y}})} (c)
(c) edge node {\pn{A_{y,z}}} (d)
;
\begin{scope}[shift={(5,0)}]
\draw[0cell=.85]
(0,0) node (a) {\big(\ang{} \scs Ax\big)}
(a)+(2.5,0) node (b) {\big(\C_{x,x} \scs Ax\big)}
(a)+(0,\v) node (c) {Ax}
(b)+(0,\v) node (d) {Ax}
;
\draw[1cell=.85]
(a) edge node {(i_x \scs \opu)} (b)
(b) edge node {\pn{A_{x,x}}} (d)
(c) edge node {\opu} (d)
;
\draw[1cell=2]
(a) edge[-,double equal sign distance] (c);
\end{scope}
\end{tikzpicture}
\end{equation}
\end{proposition}

\subsection*{Diagram Morphisms as Module Morphisms}

\begin{explanation}\label{expl:c-diag-morphism}
A morphism in $\MCat(\C,\M)$ is an $\M$-natural transformation between $\M$-functors (\cref{def:mnaturaltr}) as follows.
\[\begin{tikzpicture}
\def\h{2} \def\t{25}
\draw[0cell]
(0,0) node (a) {\C}
(a)+(\h,0) node (b) {\M}
;
\draw[1cell=.8]
(a) edge[bend left=\t] node {A} (b)
(a) edge[bend right=\t] node[swap] {B} (b)
;
\draw[2cell]
node[between=a and b at .45, rotate=-90, 2label={above,\theta}] {\Rightarrow}
;
\end{tikzpicture}\]
Such an $\M$-natural transformation consists of, for each object $x$ in $\C$, a component nullary multimorphism
\[\theta_x \cn \ang{} \to \clM(Ax;Bx) \inspace \M\]
that satisfies the naturality axiom \cref{mnaturality-diag}.
\begin{itemize}
\item Identity morphisms in $\MCat(\C,\M)$ are identity $\M$-natural transformations \cref{id-mnat}, where each component $\theta_x$ is the identity $i_{Ax}$ \cref{selfenrm-id}.
\item Composition is given by vertical composition of $\M$-natural transformations (\cref{def:mnaturaltr-vcomp}).  We characterize composition in \cref{mcat-vertical-comp} below.
\end{itemize}
The partner \cref{eval-bijection} of $\theta_x$ is a unary multimorphism
\begin{equation}\label{partner-of-theta-x}
\pn{\theta_x} \cn Ax \to Bx \inspace \M.
\end{equation}
In \cref{C-diag-morphism-pn} below, we interpret the naturality axiom \cref{mnaturality-diag} for $\theta$ in terms of these componentwise partners and those of $A$ and $B$ in \cref{Axy-pn-binary}.
\end{explanation}

\cref{C-diagram-partner} above interprets an object in $\MCat(\C,\M)$ as a left $\C$-module.  The following observation interprets a morphism in $\MCat(\C,\M)$ as a morphism of left $\C$-modules.

\begin{proposition}\label{C-diag-morphism-pn}
In the context of \cref{def:enr-diag-cat}, an $\M$-natural transformation $\theta \cn A \to B$ is uniquely determined by component unary multimorphisms as in \cref{partner-of-theta-x} 
\[\left\{\pn{\theta_x} \cn Ax \to Bx\right\}_{x\in\C}\]
such that the following diagram in $\M$ commutes for all objects $x,y \in \C$, with $\pn{A_{x,y}}$ and $\pn{B_{x,y}}$ as in \cref{Axy-pn-binary}.
\begin{equation}\label{c-diag-morphism-pn}
\begin{tikzpicture}[vcenter]
\def\v{-1.4}
\draw[0cell]
(0,0) node (a) {\big(\C(x,y) \scs Ax\big)}
(a)+(3.5,0) node (b) {Ay}
(a)+(0,\v) node (c) {\big(\C(x,y) \scs Bx\big)}
(b)+(0,\v) node (d) {By}
;
\draw[1cell=.9]
(a) edge node {\pn{A_{x,y}}} (b)
(b) edge node {\pn{\theta_y}} (d)
(a) edge node[swap] {(\opu \scs \pn{\theta_x})} (c)
(c) edge node {\pn{B_{x,y}}} (d)
;
\end{tikzpicture}
\end{equation}
\end{proposition}

\begin{proof}
Since taking partner is a bijection \cref{eval-bijection}, it suffices to show that the naturality diagram \cref{mnaturality-diag} for $\theta \cn A \to B$ yields the diagram \cref{c-diag-morphism-pn} upon taking partners.  We use the same abbreviations as in the proof of \cref{C-diagram-partner}, so
\begin{equation}\label{abbreviate-CMev}
\C(?,?) = \C_{?,?}, \quad \clM(?;?) = \clM_{?;?}, \andspace \ev = \ev^\M.
\end{equation}
For each object $x$ in $\C$, the $x$-component $\theta_x$ and its partner $\pn{\theta_x}$ determine each other via the following commutative diagram in $\M$.
\begin{equation}\label{thx-thxpn-diagram}

\end{equation}
The boundary of the following diagram in $\M$ is \cref{c-diag-morphism-pn}.
\begin{equation}\label{thxpn-naturality}
%
\end{equation}
The sub-region labeled $\pentagram$, when composed with the lower middle $\ev$, yields the partners \cref{eval-bijection} of the two composites in the naturality diagram \cref{mnaturality-diag} for $\theta \cn A \to B$.  Thus it suffices to show that the other sub-regions in \cref{thxpn-naturality} are commutative.
\begin{itemize}
\item The top left and bottom right sub-regions commute by \cref{thx-thxpn-diagram}.
\item The top right and bottom left sub-regions commute by \cref{Axy-Axypn-diagram}.
\item The two remaining sub-regions commute by the definition of $\comp$ in the canonical self-enrichment of $\M$ \cref{selfenrm-diagrams}.
\end{itemize}
This finishes the proof.
\end{proof}

\subsection*{Partner Characterization of Composition}

\cref{mcat-vertical-comp} below characterizes vertical composition of $\M$-natural transformations (\cref{def:mnaturaltr-vcomp}) in the category $\MCat(\C,\M)$ in terms of partners.  It says that, at each component, composition commutes with taking partners.

\begin{proposition}\label{mcat-vertical-comp}\index{partner!characterization of composition}
In the context of \cref{def:enr-diag-cat}, suppose $\theta$ and $\psi$ are vertically composable $\M$-natural transformations as follows.
\[
\]
Then for each object $x$ in $\C$, there is an equality of unary multimorphisms
\begin{equation}\label{pn-psi-theta-x}
\pn{(\psi\theta)_x} = \ga^\M\scmap{\pn{\psi_x}; \pn{\theta_x}} \cn Fx \to Hx \inspace \M.
\end{equation}
\end{proposition}

\begin{proof}
We use the abbreviations in \cref{abbreviate-CMev}.  By definition \cref{mnat-vcomp-component} the $x$-component of $\psi\theta$ is the following composite nullary multimorphism in $\M$.
\begin{equation}\label{psi-theta-x-component}
%
\end{equation}
Its partner \cref{thx-thxpn-diagram} is the left-bottom composite in the diagram \cref{vcomp-partner} in $\M$ below.  The right-hand side of the desired equality \cref{pn-psi-theta-x} is the top-right composite in \cref{vcomp-partner}.
\begin{equation}\label{vcomp-partner}
%
\end{equation}
The three sub-regions in \cref{vcomp-partner} are commutative for the following reasons.
\begin{itemize}
\item The bottom left sub-region commutes by the definition of $\comp$ in \cref{selfenrm-diagrams}.
\item The top left sub-region commutes by the definition of $\pn{\theta_x}$ \cref{thx-thxpn-diagram}.
\item The right sub-region commutes by the definition of $\pn{\psi_x}$ \cref{thx-thxpn-diagram}.
\end{itemize}
The commutative diagram \cref{vcomp-partner} proves the desired equality \cref{pn-psi-theta-x}.
\end{proof}

\subsection*{Examples of Mackey Functor Categories}

\begin{example}[Stable Model Categories]\label{ex:schwede-shipley}
  Suppose $\M$ is a simplicial, cofibrantly generated, proper, and stable model category (\cref{definition:model-variants}).
  The Schwede-Shipley Characterization \cref{theorem:schwede-shipley} shows that, if $P$ is a set of compact generators of $\M$, then there is a chain of simplicial Quillen equivalences
  \[\M \hty_Q \Spcat\big(\EP^\op, \Sp\big)\]
  between $\M$ and the $\EP$-presheaf category of $\Sp$ in the sense of \cref{mcat-copm}. 
  On the right-hand side, $\EP$ is the spectral endomorphism category (\cref{definition:EP}) and $\EP^\op$ is its opposite $\Sp$-category as in \cref{opposite-mcat}.
\end{example}

\begin{example}[Genuine Equivariant Spectra]\label{ex:guillou-may}
  Recall from \cref{definition:Burnside2} the \index{Burnside!2-category}\index{2-category!Burnside}Burnside 2-category $\GE$ for a finite group $G$.
  The Guillou-May \cref{theorem:GM} gives a zigzag of Quillen equivalences
  \begin{equation}\label{guillou-may-eq}
    \GSp \hty_Q \Spcat\big((\GE_\bbK)^\op \scs \Sp\big)
  \end{equation}
  between the category of genuine equivariant $G$-spectra, $\GSp$, and the
  category of spectral Mackey functors for $\bbK$ (\cref{definition:Sp-G-Mackey}).
  As in \cref{rk:BO-Phi}, $\bbK$ denotes the $K$-theory non-symmetric multifunctor in \cite{guillou_may,gmmo}.

  We emphasize that in \cref{guillou-may-eq} the opposite in $(\GE_\bbK)^\op$ is taken in $\Spcat$ (\cref{opposite-mcat}) \emph{after} the change of enrichment along $\bbK$.  
  In \cref{rk:BO7.5} we further discuss the relationship between
  \begin{itemize}
  \item the Guillou-May Quillen equivalence \cref{guillou-may-eq},
  \item the work of Bohmann-Osorno \cite{bohmann_osorno-mackey}, and
  \item our diagram change-of-enrichment functor in \cref{gspectra-thm-v}.
  \end{itemize} 
  There, we note and discuss the nontrivial distinction between the $\Sp$-categories $(\GE_\bbK)^\op$ and $(\GE^\op)_\bbK$.
\end{example}

\section{Change of Enrichment of Enriched Diagrams and Mackey Functors}
\label{sec:enr-diag-change-enr}

In this section we construct change-of-enrichment functors on enriched diagram and Mackey functor categories associated to (non-symmetric) multifunctors between (non-symmetric) closed multicategories (\cref{def:enr-multicategory-functor,def:closed-multicat,def:enr-diag-cat}). 
\begin{itemize}
\item The change-of-enrichment construction is in \cref{def:enr-diag-change-enr}.
\item \cref{expl:Fdg-object,expl:Fdg-morphism} unpack the change-of-enrichment construction on objects and morphisms.
\end{itemize}
We defer the proof that change of enrichment is a functor to \cref{sec:diag-psh-change-enr-functors}; see \cref{gspectra-thm-v,gspectra-thm-v-cor}.  In \cref{sec:fun-change-enr-diag} we show that these change-of-enrichment functors are compatible with composition of (non-symmetric) multifunctors.

\subsection*{Defining $\Fdg$}

\begin{definition}\label{def:enr-diag-change-enr}
Suppose given a non-symmetric multifunctor between non-symmetric closed multicategories
\[F \cn \M \to \N\]
and a small $\M$-category $\C$ (\cref{def:menriched-cat}).  We define the data of a functor
\begin{equation}\label{diag-change-enr}
\Fdg \cn \MCat(\C,\M) \to \NCat(\C_F,\N),
\end{equation}
called the \index{enrichment!diagram change of -}\index{change of enrichment!diagram}\index{diagram!change of enrichment}\emph{diagram change of enrichment} of $F$ at $\C$, as follows.
\begin{description}
\item[Domain] The domain of $\Fdg$ is the $\C$-diagram category of $\M$ in \cref{mcat-cm}.
\item[Codomain] $\C_F$ is the $\N$-category obtained from $\C$ by applying the change-of-enrichment 2-functor in \cref{mult-change-enrichment}
\[\dF \cn \MCat \to \NCat.\]
The codomain of $\Fdg$ is the $\C_F$-diagram category of $\N$ in \cref{mcat-cm}.
\item[Object and Morphism Assignments] Suppose given
\begin{itemize}
\item $\M$-functors $A$ and $B$ and
\item an $\M$-natural transformation $\psi \cn A \to B$ in $\MCat(\C,\M)$
\end{itemize}
as in the left diagram below.
\begin{equation}\label{diag-change-enr-assign}
\begin{tikzpicture}[baseline={(a.base)}]
\def\t{25} \def\d{.8}
\draw[0cell]
(0,0) node (a) {\C}
(a)+(1.7,0) node (b) {\M}
(b)+(\d,0) node (x) {}
(x)+(1,0) node (y) {}
(y)+(\d,0) node (a') {\C_F}
(a')+(2,0) node (b') {\M_F}
(b')+(1.5,0) node (c) {\N}
;
\draw[1cell]
(a) edge[bend left=30] node {A} (b)
(a) edge[bend right=30] node[swap] {B} (b)
(x) edge[|->] node {\Fdg} (y)
;
\draw[2cell]
node[between=a and b at .42, rotate=-90, 2label={above,\psi}] {\Rightarrow}
;
\draw[1cell]
(a') edge[bend left=\t] node {A_F} (b')
(a') edge[bend right=\t] node[swap] {B_F} (b')
(b') edge node {\Fse} (c)
;
\draw[2cell]
node[between=a' and b' at .4, rotate=-90, 2label={above,\psi_F}] {\Rightarrow}
;
\end{tikzpicture}
\end{equation}
Then $\Fdg$ sends $A$, $B$, and $\psi$ to the composites and whiskering as in the right diagram in \cref{diag-change-enr-assign}, with $\Fse$ the standard enrichment of $F$ in \cref{gspectra-thm-iii}. 
\end{description}
This finishes the definition of $\Fdg$.  We also call $\Fdg$ a \index{functor!diagram change-of-enrichment}\index{diagram!change-of-enrichment functor}\emph{diagram change-of-enrichment functor}.

Moreover, we define the following.
\begin{itemize}
\item If we want to emphasize $\C$, then we write $\Fdg^{\C}$ instead of $\Fdg$.
\item Suppose $F$ is a multifunctor between closed multicategories and $\Cop$ is the opposite $\M$-category of $\C$ (\cref{opposite-mcat}).  Using \cref{dF-opposite} to identify the $\N$-categories $(\Cop)_F$ and $(\C_F)^\op$, we call
\begin{equation}\label{presheaf-change-enr}
\Fdg \cn \MCat(\Cop,\M) \to \NCat\big((\C_F)^\op,\N\big)
\end{equation}
the \index{functor!presheaf change-of-enrichment}\index{presheaf!change of enrichment}\emph{presheaf change of enrichment} of $F$.  We also call this $\Fdg$ a \emph{presheaf change-of-enrichment functor}.
\end{itemize} 
This finishes the definition.  \cref{gspectra-thm-v,gspectra-thm-v-cor} prove that $\Fdg$ in \cref{diag-change-enr,presheaf-change-enr} are functors.
\end{definition}

The assignment $\Fdg$ in \cref{diag-change-enr-assign} is
\begin{itemize}
\item the change of enrichment $\dF$ (\cref{mult-change-enrichment}) followed by
\item the standard enrichment $\Fse$ (\cref{gspectra-thm-iii}).
\end{itemize}
We describe $\Fdg$ in more detail in \cref{expl:Fdg-object,expl:Fdg-morphism} below.

\subsection*{Unpacking $\Fdg$}

\begin{explanation}[$\Fdg$ on Objects]\label{expl:Fdg-object}
For an $\M$-functor $A \cn \C \to \M$, in \cref{diag-change-enr-assign} the $\N$-functor
\begin{equation}\label{FdgA}
\Fdg A \cn \C_F \fto{A_F} \M_F \fto{\Fse} \N
\end{equation}
is the composite of
\begin{itemize}
\item $A_F$, which is the image of $A$ under the change of enrichment $\dF$, and
\item the standard enrichment $\Fse$ of $F$ in \cref{gspectra-thm-iii}.
\end{itemize}

\medskip
\emph{Object Assignment}.
More explicitly, the $\N$-functor $A_F$ has the same object assignment as $A \cn \C \to \M$.  The standard enrichment $\Fse$ has the same object assignment as $F$.  Thus the object assignment of the $\N$-functor $\Fdg A$ is given by
\begin{equation}\label{FAx}
(\Fdg A)(x) = F\big(A(x)\big) \in \N \forspace x \in \C.
\end{equation}

\medskip
\emph{Components}.
To describe $\Fdg A$ on hom objects, suppose given a pair of objects $x,y \in \C$.  The component $(\Fdg A)_{x,y}$ is the following composite in $\N$ of unary multimorphisms.
\begin{equation}\label{FAxy}
\begin{tikzpicture}[baseline={(a.base)}]
\def\h{3.5} \def\w{.65}
\draw[0cell=1]
(0,0) node (a) {F \C(x,y)}
(a)+(\h-.3,0) node (b) {F \clM\scmap{Ax;Ay}}
(b)+(\h+.3,0) node (c) {\clN\scmap{FAx;FAy}}
;
\draw[1cell=.9]
(a) edge node {F A_{x,y}} (b)
(b) edge node {\Fse_{Ax, Ay}} (c)
;
\draw[1cell=.9]
(a) [rounded corners=3pt] |- ($(b)+(-1,\w)$)
-- node {(\Fdg A)_{x,y}} ($(b)+(1,\w)$) -| (c)
;
\end{tikzpicture}
\end{equation}
The two constituent arrows in \cref{FAxy} are as follows.
\begin{itemize}
\item The unary multimorphism 
\begin{equation}\label{A-sub-xy}
A_{x,y} \cn \C(x,y) \to \clM\scmap{Ax;Ay} \inspace \M
\end{equation}
is the $(x,y)$-component of the $\M$-functor $A$ (\cref{def:mfunctor}).  The left arrow in \cref{FAxy} is the image of $A_{x,y}$ under $F$, which is a unary multimorphism in $\N$.
\item The right arrow in \cref{FAxy} is the $(Ax,Ay)$-component unary multimorphism of the standard enrichment $\Fse$ in \cref{Fprimexy}.  Its construction uses the closed structure on both $\M$ and $\N$.
\end{itemize}

\medskip
\emph{Partner Characterization}.
We can characterize $(\Fdg A)_{x,y}$ in terms of its partner \cref{eval-bijection}, using the following commutative diagram in $\N$.
\begin{equation}\label{FdgAxy-partner}
\begin{tikzpicture}[vcenter]
\def\v{-1.3} \def\w{-2.4} \def\t{18}
\draw[0cell=.9]
(0,0) node (a) {\big( F\C(x,y) \scs FAx \big)}
(a)+(0,\v) node (b) {\big( F\clM\scmap{Ax;Ay} \scs FAx \big)}
(b)+(0,\v) node (c) {\big( \clN\scmap{FAx;FAy} \scs FAx \big)}
(b)+(4,0) node (d) {FAy}
;
\draw[1cell=.9]
(a) edge node[swap] {\left( FA_{x,y} \scs 1_{FAx} \right)} (b)
(b) edge node[swap] {\left(\Fse_{Ax,Ay} \scs 1_{FAx}\right)} (c)
(a) edge[bend left=\t] node[pos=.6] {F(\pn{A_{x,y}})} (d)
(b) edge node[pos=.4] {F(\ev^\M_{Ax;\, Ay})} (d)
(c) edge[bend right=\t] node[swap,pos=.7] {\ev^\N_{FAx;\, FAy}} (d)
;
\draw[1cell=.9]
(a) [rounded corners=3pt] -| ($(b)+(\w,1)$)
-- node[swap,pos=0] {\left( (\Fdg A)_{x,y} \scs 1_{FAx} \right)} ($(b)+(\w,-1)$) |- (c)
;
\end{tikzpicture}
\end{equation}
The diagram \cref{FdgAxy-partner} commutes for the following reasons.
\begin{itemize}
\item $\pn{A_{x,y}}$ is the partner \cref{eval-bijection} of $A_{x,y}$ in \cref{A-sub-xy}.  By definition it is the following composite in $\M$.
\begin{equation}\label{Axy-partner}
\begin{tikzpicture}[vcenter]
\def\v{-1.3} \def\t{18}
\draw[0cell=.9]
(0,0) node (a) {\big( \C(x,y) \scs Ax \big)}
(a)+(0,\v) node (b) {\big( \clM\scmap{Ax;Ay} \scs Ax \big)}
(b)+(3.5,0) node (d) {Ay}
;
\draw[1cell=.9]
(a) edge node[swap] {\left( A_{x,y} \scs 1_{Ax} \right)} (b)
(a) edge[bend left=\t] node[pos=.6] {\pn{A_{x,y}}} (d)
(b) edge node[pos=.4] {\ev^\M_{Ax;\, Ay}} (d)
;
\end{tikzpicture}
\end{equation}
The upper right region in \cref{FdgAxy-partner} commutes because it is the image under $F$ of \cref{Axy-partner}.  This uses the fact that $F$, as a non-symmetric multifunctor, preserves colored units and composition.
\item The bottom right region in \cref{FdgAxy-partner} commutes by the definition of $\Fse_{Ax,Ay}$; see \cref{Fse-xy-diag}.
\item The left region in \cref{FdgAxy-partner} commutes by \cref{FAxy}.
\end{itemize}
In \cref{FdgAxy-partner} the left-bottom composite is, by definition, the partner of $(\Fdg A)_{x,y}$.  Thus \cref{FdgAxy-partner} yields the equality of binary multimorphisms
\begin{equation}\label{FdgAxy-pn-equality}
\pn{(\Fdg A)_{x,y}} = F(\pn{A_{x,y}}) \cn \big( F\C(x,y) \scs FAx \big) \to FAy \inspace \N.
\end{equation}
So the partner of $(\Fdg A)_{x,y}$ is the image under $F$ of the partner of $A_{x,y}$.
\end{explanation}

\begin{explanation}[$\Fdg$ on Morphisms]\label{expl:Fdg-morphism}
Suppose $\psi$ is an $\M$-natural transformation as follows.
\[\begin{tikzpicture}[baseline={(a.base)}]
\def\t{25} \def\d{.8}
\draw[0cell]
(0,0) node (a) {\C}
(a)+(1.7,0) node (b) {\M}
;
\draw[1cell=.8]
(a) edge[bend left=30] node {A} (b)
(a) edge[bend right=30] node[swap] {B} (b)
;
\draw[2cell]
node[between=a and b at .42, rotate=-90, 2label={above,\psi}] {\Rightarrow}
;
\end{tikzpicture}\]
In the definition \cref{diag-change-enr-assign},
\begin{equation}\label{Fdgpsi}
\Fdg \psi = 1_{\Fse} * \psi_F \cn \Fdg A = \Fse \circ A_F \to \Fdg B = \Fse \circ B_F
\end{equation}
is the horizontal composite (\cref{def:mnaturaltr-hcomp}) below.
\[\begin{tikzpicture}[baseline={(a.base)}]
\def\h{2} \def\t{28} \def\d{.8}
\draw[0cell]
(0,0) node (a) {\C_F}
(a)+(\h,0) node (b) {\M_F}
(b)+(\h,0) node (c) {\N}
;
\draw[1cell=.8]
(a) edge[bend left=\t] node {A_F} (b)
(a) edge[bend right=\t] node[swap] {B_F} (b)
(b) edge[bend left=\t] node {\Fse} (c)
(b) edge[bend right=\t] node[swap] {\Fse} (c)
;
\draw[2cell]
node[between=a and b at .42, rotate=-90, 2label={above,\psi_F}] {\Rightarrow}
node[between=b and c at .42, rotate=-90, 2label={above,1_{\Fse}}] {\Rightarrow}
;
\end{tikzpicture}\]
\begin{itemize}
\item The $\N$-natural transformation $\psi_F$ is the image of $\psi$ under the change of enrichment $\dF$.
\item $1_{\Fse}$ is the identity $\N$-natural transformation of the standard enrichment $\Fse$ \cref{id-mnat}.
\end{itemize}

\medskip
\emph{Components}.
More explicitly, for each object $x \in \C$, the $x$-component of $\psi \cn A \to B$ is a nullary multimorphism
\begin{equation}\label{psi-xcomponent}
\psi_x \cn \ang{} \to \clM\scmap{Ax;Bx} \inspace \M,
\end{equation}
with $\ang{}$ denoting the empty sequence.  After the change of enrichment $\dF$, the $x$-component of $\psi_F \cn A_F \to B_F$ is the nullary multimorphism
\[(\psi_F)_x = F \psi_x \cn \ang{} \to F\clM\scmap{Ax;Bx} \inspace \N.\]
The $x$-component of $\Fdg\psi$ is the nullary multimorphism given by the following composite in $\N$.
\begin{equation}\label{Fdgpsi-x}

\end{equation}
The right arrow in \cref{Fdgpsi-x} is the $(Ax,Bx)$-component unary multimorphism of the standard enrichment $\Fse$ in \cref{Fprimexy}.

\medskip
\emph{Partner Characterization}.
We can describe the component $(\Fdg\psi)_x$ in terms of its partner \cref{eval-bijection}, using the following commutative diagram in $\N$.
\begin{equation}\label{Fdgpsix-partner}
%
\end{equation}
The diagram \cref{Fdgpsix-partner} commutes for the following reasons.
\begin{itemize}
\item $\pn{\psi_x}$ is the partner \cref{eval-bijection} of $\psi_x$ in \cref{psi-xcomponent}.  By definition it is the following composite in $\M$.
\begin{equation}\label{psix-partner}
%
\end{equation}
The upper right region in \cref{Fdgpsix-partner} commutes because it is the image under $F$ of \cref{psix-partner}.  This uses the fact that $F$, as a non-symmetric multifunctor, preserves colored units and composition.
\item The bottom right region in \cref{Fdgpsix-partner} commutes by the definition of $\Fse_{Ax,Bx}$; see \cref{Fse-xy-diag}.
\item The left region in \cref{Fdgpsix-partner} commutes by \cref{Fdgpsi-x}.
\end{itemize}
In \cref{Fdgpsix-partner} the left-bottom composite is, by definition, the partner of $(\Fdg \psi)_x$.  Thus \cref{Fdgpsix-partner} yields the equality of unary multimorphisms
\begin{equation}\label{Fdg-psix-pn-equality}
\pn{(\Fdg \psi)_x} = F(\pn{\psi_x}) \cn FAx \to FBx \inspace \N.
\end{equation}
So the partner of $(\Fdg \psi)_x$ is the image under $F$ of the partner of $\psi_x$.
\end{explanation}

\section{Diagram and Mackey Functor Change-of-Enrichment Functors}
\label{sec:diag-psh-change-enr-functors}

This section has two purposes.
\begin{enumerate}
\item We show that the diagram and Mackey functor change of enrichment in \cref{diag-change-enr,presheaf-change-enr} are functors (\cref{gspectra-thm-v,gspectra-thm-v-cor}).
\item We illustrate them with $K$-theoretic functors in \cref{ex:P-diag-change-enr,ex:gsp-v-examples}.
\end{enumerate}

\subsection*{Diagram Change of Enrichment is a Functor}

\begin{theorem}\label{gspectra-thm-v}
For each non-symmetric multifunctor between non-symmetric closed multicategories
\[F \cn \M \to \N\]
and each small $\M$-category $\C$, the diagram change of enrichment
\[\Fdg \cn \MCat(\C,\M) \to \NCat(\C_F,\N)\]
in \cref{diag-change-enr} is a functor.
\end{theorem}

\begin{proof}
The assignments of $\Fdg$
\begin{itemize}
\item on objects
\[A \mapsto \Fdg A = \Fse \circ A_F\]
in \cref{FdgA} and
\item on morphisms
\[\psi \mapsto \Fdg \psi = 1_{\Fse} * \psi_F\]
in \cref{Fdgpsi}
\end{itemize}  
are well defined because they are given by composition of $\N$-functors and horizontal composition of $\N$-natural transformations, respectively.

\medskip
\emph{Preservation of Identity Morphisms}.
Suppose $\psi = 1_A$ is the identity $\M$-natural transformation of an $\M$-functor $A \cn \C \to \M$.  Then the 2-functoriality of the change of enrichment $\dF$ (\cref{mult-change-enrichment}) implies the following equalities.
\[\Fdg 1_A = 1_{\Fse} * (1_A)_F = 1_{\Fse} * 1_{(A_F)} = 1_{\Fdg A}\]

\medskip
\emph{Preservation of Composition}.
Suppose given $\M$-functors $A,B,D \cn \C \to \M$ and vertically composable $\M$-natural transformations 
\[A \fto{\psi} B \fto{\phi} D.\]
The following computation, using the 2-functoriality of $\dF$ in the second equality, shows that $\Fdg$ preserves composition.
\[\begin{split}
\Fdg(\phi \psi) &= 1_{\Fse} * (\phi\psi)_F\\
&= 1_{\Fse} * (\phi_F \psi_F)\\
&= \big(1_{\Fse} * \phi_F \big) \big(1_{\Fse} * \psi_F\big)\\
&= (\Fdg\phi) (\Fdg\psi)
\end{split}\]
This finishes the proof.
\end{proof}

\begin{example}[Inverse $K$-Theory and Free Permutative Categories]\label{ex:P-diag-change-enr}
\cref{gspectra-thm-v} is applicable to the non-symmetric multifunctors in the following diagram.
\begin{equation}\label{PFFF}
\begin{tikzpicture}[baseline={(pm.base)}]
\def\h{3} \def\f{1} \def\v{1}
\draw[0cell]
(0,0) node (pm) {\pMulticat}
(pm)+(\h+.3,0) node (p') {\phantom{\permcatsu}}
(p')+(0,.04) node (p) {\permcatsu}
(p')+(\h,0) node (g') {\phantom{\Gacat}}
(g')+(0,.03) node (g) {\Gacat}
(pm)+(\f,\v) node (m) {\Multicat}
(pm)+(\f,-\v) node (mone) {\MoneMod}
;
\draw[1cell]
(m) edge node[pos=.4] {\Fr} (p)
(pm) edge node[pos=.4] {\Fst} (p')
(mone) edge node[swap,pos=.4] {\Fm} (p)
(g') edge node[swap] {\cP} (p')
;
\end{tikzpicture}
\end{equation}
\begin{itemize}
\item $\permcatsu$ is a closed multicategory by \cref{perm-closed-multicat}.
\item $\Multicat$, $\pMulticat$, $\MoneMod$, and $\Gacat$ are symmetric monoidal closed categories by \cref{theorem:multicat-sm-closed,thm:pmulticat-smclosed,proposition:EM2-5-1,GammaV}.  Thus they are closed multicategories by \cref{smclosed-closed-multicat}.
\item $\cP$ is a non-symmetric multifunctor by \cite[1.3]{johnson-yau-invK}.
\item $\Fr$, $\Fst$, and $\Fm$ are non-symmetric multifunctors by \cref{theorem:F-multi,ptmulticat-xvii,Fm-multi-def}. 
\end{itemize}

For example, \cref{gspectra-thm-v} applied to inverse $K$-theory $\cP$ says that, for each small $(\Gacat)$-category $\C$ (\cref{def:menriched-cat,def:enriched-category}), there is a diagram change-of-enrichment functor\label{not:cPdg}
\[\cPdg \cn \Gacatcat\big(\C,\Gacat\big) \to \permcatsucat\big(\C_{\cP},\permcatsu\big)\]
defined as in \cref{diag-change-enr-assign}.  More explicitly, the functor $\cPdg$ sends each $(\Gacat)$-functor (\cref{def:mfunctor,def:enriched-functor})
\[A \cn \C \to \Gacat\]
to the composite $\permcatsu$-functor\label{not:APPse}
\[\cPdg A \cn \C_{\cP} \fto{A_{\cP}} (\Gacat)_{\cP} \fto{\cPse} \permcatsu.\]
\begin{itemize}
\item $A_{\cP}$ is the image of $A$ under the change-of-enrichment 2-functor (\cref{mult-change-enrichment}) 
\[\dcP \cn \Gacatcat \to \permcatsucat\]
along $\cP$.
\item The $\permcatsu$-functor
\[\cPse \cn (\Gacat)_{\cP} \to \permcatsu\]
is the standard enrichment of $\cP$ (\cref{gspectra-thm-iii}).
\end{itemize} 
The morphism assignment of $\cPdg$ sends a $(\Gacat)$-natural transformation (\cref{def:mnaturaltr,def:enriched-natural-transformation}) $\psi$ to the whiskering
\[\cPdg \psi = 1_{\cPse} * \psi_{\cP}\]
as in \cref{Fdgpsi}.  The diagram change-of-enrichment functors $\Frdg$, $\Fstdg$, and $\Fmdg$ admit analogous description.
\end{example}

We provide further examples and applications of \cref{gspectra-thm-v} in \cref{ex:gsp-v-examples,ch:mackey,ch:mackey_eq}.  In particular, in \cref{ch:mackey_eq} we use $\Fstdg$ and $\Fmdg$ to construct equivalences of homotopy theories from diagrams and presheaves enriched in $\pMulticat$ and $\MoneMod$ to those enriched in $\permcatsu$.

\subsection*{Mackey Functor Change of Enrichment is a Functor}

\begin{theorem}\label{gspectra-thm-v-cor}\index{enrichment!Mackey functor change of -}\index{change of enrichment!Mackey functor}\index{Mackey functor!change of enrichment}
For each multifunctor between closed multicategories 
\[F\cn \M \to \N\]
and each small $\M$-category $\C$, the presheaf change of enrichment
\[\Fdg \cn \MCat(\Cop,\M) \to \NCat\big( (\C_F)^\op,\N \big)\]
in \cref{presheaf-change-enr} is a functor.
\end{theorem}

\begin{proof}
This is \cref{gspectra-thm-v} applied to the opposite $\M$-category $\Cop$.  We use \cref{dF-opposite} to obtain the equality 
\[(\Cop)_F = (\C_F)^\op\]
of $\N$-categories.
\end{proof}

We stress that \cref{gspectra-thm-v-cor} does \emph{not} apply to non-symmetric multifunctors.  For the equality in its proof to hold, $F$ needs to preserve the symmetric group action as in the second equality in \cref{CopF-composition}.

\begin{example}[$K$-Theory Multifunctors]\label{ex:gsp-v-examples}\index{enrichment!Mackey functor change of -!Elmendorf-Mandell $K$-theory}\index{change of enrichment!Mackey functor!Elmendorf-Mandell $K$-theory}\index{Mackey functor!change of enrichment!Elmendorf-Mandell $K$-theory}\index{Elmendorf-Mandell!K-theory@$K$-theory!Mackey functor change of enrichment}\index{K-theory@$K$-theory!Elmendorf-Mandell!Mackey functor change of enrichment}
\cref{gspectra-thm-v,gspectra-thm-v-cor} are applicable to the following multifunctors in \cref{eq:Ksummary}.
\begin{equation}\label{Kem-functors}
\begin{tikzpicture}[xscale=1.2,yscale=1,vcenter]
\def\t{.35} \def\v{-1.5} \def\u{.6} \def\w{-.6}
\draw[0cell=.8]
(0,0) node (a1) {\permcatsu}
(a1)+(3,0) node (a2) {\Gacat}
(a2)+(2,0) node (a3) {\Gasset}
(a3)+(2,0) node (a4) {\Sp}
(a2)+(0,\v) node (b2) {\Gstarcat}
(a3)+(0,\v) node (b3) {\Gstarsset}
(b2)+(-2.3,0) node (b1) {\MoneMod}
;
\draw[1cell=.8]
(a2) edge node {\Ner_*} (a3)
(a3) edge node {\Kf} (a4)
(a1) edge node[pos=.3] {\Jem} (b2)
(a1) edge node[pos=.85] {\Endm} (b1)
(b1) edge node[pos=.4] {\Jt} (b2)
(b2) edge node {\Ner_*} (b3)
(b3) edge node [pos=.4,swap] {\Kg} (a4)
(a2) edge[transform canvas={xshift=0ex}] node[swap,pos=.3] {\sma^*} (b2)
(a3) edge[transform canvas={xshift=0ex}] node[pos=.3] {\sma^*} (b3)
;
\draw[1cell=.85]
(a1) [rounded corners=3pt] |- ($(b2)+(0,\w)$)
-- node[pos=\t] {\Kem} ($(b3)+(0,\w)$) -| (a4)
;
\end{tikzpicture}
\end{equation}
Each arrow in \cref{Kem-functors} is a multifunctor between closed multicategories.  
\begin{itemize}
\item $\permcatsu$, $\Gacat$, and $\Gasset$ are closed multicategories by \cref{GammaV,smclosed-closed-multicat,perm-closed-multicat}, as discussed in \cref{ex:P-diag-change-enr}.
\item $\MoneMod$, $\Gstarcat$, $\Gstarsset$, and $\Sp$ are symmetric monoidal closed categories (\cref{proposition:EM2-5-1,Gstar-V,SymSp}), hence also closed multicategories (\cref{smclosed-closed-multicat}).
\item $\Jem$, $\Endm$, and $\Kem$ are multifunctors.
\item The other arrows in \cref{Kem-functors} are symmetric monoidal functors, hence also multifunctors via the endomorphism construction \cref{EndP}.
\end{itemize}
We discuss the case for $\Kem$ in more detail in \cref{sec:presheaf-K,sec:mult-mackey-spectra}.  We emphasize that \cref{gspectra-thm-v-cor} does \emph{not} apply to the arrows in \cref{PFFF}---namely, $\cP$, $\Fr$, $\Fst$, and $\Fm$---because those are \emph{non-symmetric} multifunctors.
\end{example}

\section{Composition of Diagram Change-of-Enrichment Functors}
\label{sec:fun-change-enr-diag}

In \cref{gspectra-thm-v} we observe that there is a diagram change-of-enrichment functor
\[\Fdg = \Fdg^\C \cn \MCat(\C,\M) \to \NCat(\C_F,\N)\]
for each
\begin{itemize}
\item non-symmetric multifunctor $F \cn \M \to \N$ between non-symmetric closed multicategories and 
\item small $\M$-category $\C$.
\end{itemize}   
In this section we show that the construction, $F \mapsto \Fdg$, respects composition of non-symmetric multifunctors; see \cref{gspectra-thm-vii}.  The version for enriched Mackey functors is \cref{gspectra-thm-vii-cor}.  We discuss applications of \cref{gspectra-thm-vii,gspectra-thm-vii-cor} to Elmendorf-Mandell $K$-theory in \cref{sec:mult-mackey-spectra}.

\begin{theorem}\label{gspectra-thm-vii}\index{compositionality!- of diagram change of enrichment}\index{functor!diagram change-of-enrichment!compositionality of}\index{diagram!change-of-enrichment functor!compositionality}
Suppose given non-symmetric multifunctors between non-symmetric closed multicategories
\[\M \fto{F} \N \fto{G} \P\]
and a small $\M$-category $\C$.  Then the following diagram of functors commutes.
\begin{equation}\label{gspectra-vii-diag}
\begin{tikzpicture}[vcenter]
\def\h{4.3}
\draw[0cell]
(0,0) node (a) {\MCat(\C,\M)}
(a)+(\h/2,-1.2) node (b) {\NCat(\C_F,\N)}
(a)+(\h,0) node (c) {\PCat(\C_{GF},\P)}
;
\draw[1cell]
(a) edge node {\GFdg^\C} (c)
(a) edge node[swap,pos=.2] {\Fdg^\C} (b)
(b) edge node[swap,pos=.8] {\Gdg^{\C_F}} (c)
;
\end{tikzpicture}
\end{equation}
\end{theorem}

\begin{proof}
By \cref{func-change-enr} the following diagram of change-of-enrichment 2-functors commutes.
\begin{equation}\label{MNPCat}
\begin{tikzpicture}[baseline={(a.base)}]
\def\h{2.5} \def\w{.6}
\draw[0cell=1]
(0,0) node (a) {\MCat}
(a)+(\h,0) node (b) {\NCat}
(b)+(\h,0) node (c) {\PCat}
;
\draw[1cell=.9]
(a) edge node {\dF} (b)
(b) edge node {\dG} (c)
;
\draw[1cell=.9]
(a) [rounded corners=3pt] |- ($(b)+(-1,\w)$)
-- node {\dGF} ($(b)+(1,\w)$) -| (c)
;
\end{tikzpicture}
\end{equation}
This gives an equality of $\P$-categories
\[(\C_F)_G = \C_{GF},\]
so the arrow $\Gdg^{\C_F}$ in \cref{gspectra-vii-diag} is well defined.

To prove that \cref{gspectra-vii-diag} is commutative, suppose $A \cn \C \to \M$ is an $\M$-functor.  By definition \cref{FdgA} $\Fdg A$ is the composite $\N$-functor
\[\Fdg A \cn \C_F \fto{A_F} \M_F \fto{\Fse} \N\]
with
\begin{itemize}
\item $\dF$ the change of enrichment in \cref{MNPCat} and
\item $\Fse$ the standard enrichment of $F$ (\cref{gspectra-thm-iii}).
\end{itemize}   
Applying the change-of-enrichment 2-functor $\dG$ to the above composite and composing with the standard enrichment $\Gse$, we obtain the composite $\P$-functor along the top of the following diagram.
\begin{equation}\label{diagram-change-func}
\begin{tikzpicture}[baseline={(a.base)}]
\def\g{2.7} \def\h{1.7} \def\w{.6}
\draw[0cell=1]
(0,0) node (a) {\C_{GF}}
(a)+(\g,0) node (b) {\M_{GF} = (\M_F)_G}
(b)+(\g,0) node (c) {\N_G}
(c)+(\h,0) node (d) {\P}
;
\draw[1cell=.9]
(a) edge node {A_{GF}} (b)
(b) edge node {\Fse_G} (c)
(c) edge node {\Gse} (d)
;
\draw[1cell=.9]
(a) [rounded corners=3pt] |- ($(b)+(-1,\w)$)
-- node {(\Fdg A)_G} ($(b)+(1,\w)$) -| (c)
;
\draw[1cell=.9]
(b) [rounded corners=3pt] |- ($(c)+(-1,-\w)$)
-- node[pos=.2] {\GFse} ($(c)+(1,-\w)$) -| (d)
;
\end{tikzpicture}
\end{equation}
\cref{gspectra-thm-iv} gives the equality of $\P$-functors
\[\GFse = \Gse \circ \Fse_G \cn \M_{GF} \to \P.\]
In the commutative diagram \cref{diagram-change-func},
\begin{itemize}
\item the composite along the top is $\Gdg\Fdg(A)$, and
\item the composite along the bottom is $\GFdg(A)$. 
\end{itemize} 
This proves that the diagram \cref{gspectra-vii-diag} is commutative on objects.

Replacing $A$ by an $\M$-natural transformation in $\MCat(\C,\M)$, the previous paragraph also proves that the diagram \cref{gspectra-vii-diag} is commutative on morphisms.
\end{proof}

The following result shows that presheaf change-of-enrichment functors are closed under composition.  Recall from \cref{def:enr-multicategory-functor} that multifunctors are required to preserve the colored units, composition, and symmetric group action.

\begin{theorem}\label{gspectra-thm-vii-cor}
Suppose given multifunctors between closed multicategories
\[\M \fto{F} \N \fto{G} \P\]
and a small $\M$-category $\C$.  Then the following diagram of presheaf change-of-enrichment functors commutes.
\begin{equation}\label{gspectra-vii-cor-diag}
\begin{tikzpicture}[vcenter]
\def\h{4.7}
\draw[0cell]
(0,0) node (a) {\MCat(\Cop,\M)}
(a)+(\h/2,-1.2) node (b) {\NCat\big((\C_F)^\op,\N\big)}
(a)+(\h,0) node (c) {\PCat\big((\C_{GF})^\op,\P\big)}
;
\draw[1cell]
(a) edge node {\GFdg^{\Cop}} (c)
(a) edge[transform canvas={xshift=-1em}] node[swap,pos=.2] {\Fdg^{\Cop}} (b)
(b) edge[transform canvas={xshift=1em}] node[swap,pos=.8] {\Gdg^{(\C_F)^\op}} (c)
;
\end{tikzpicture}
\end{equation}
\end{theorem}

\begin{proof}
This is \cref{gspectra-thm-vii} applied to the opposite $\M$-category $\Cop$ (\cref{opposite-mcat}).  The equality of $\N$-categories 
\[(\Cop)_F = (\C_F)^\op\]
and the equalities of $\P$-categories 
\[\big((\C_F)^\op\big)_G = \big((\C_F)_G\big)^\op = (\C_{GF})^\op = (\Cop)_{GF}\]
are from \cref{func-change-enr,dF-opposite}.  
\end{proof}

We emphasize that \cref{gspectra-thm-vii-cor} does \emph{not} apply to non-symmetric multifunctors because the equalities in the proof above require that $F$ and $G$ preserve the symmetric group action.

\section{Spectral Mackey Functors from \texorpdfstring{$K$}{K}-Theory}
\label{sec:presheaf-K}

Recall from \cref{Kem} that Elmendorf-Mandell $K$-theory
\[\Kem \cn \permcatsu \to \Sp\]
is a multifunctor in the $\sSet$-enriched sense.  In this section we apply our general results about diagram and presheaf change of enrichment, \cref{gspectra-thm-v,gspectra-thm-v-cor}, to obtain spectrally enriched diagrams and Mackey functors via Elmendorf-Mandell $K$-theory.  

\begin{theorem}\label{thm:Kemdg}
Suppose $\C$ is a small $\permcatsu$-category.  Then $\Kem$ induces a diagram change-of-enrichment functor
\begin{equation}\label{Kemdg-functor}
\permcatsucat\big(\C, \permcatsu\big) \fto{\Kemdg} \Spcat\big(\C_{\Kem}, \Sp\big)
\end{equation}
and a presheaf change-of-enrichment functor 
\begin{equation}\label{Kemdg-presheaf}
\permcatsucat\big(\Cop, \permcatsu\big) \fto{\Kemdg} \Spcat\big((\C_{\Kem})^\op, \Sp\big).
\end{equation}
\end{theorem}

\begin{proof}
The existence of the functor $\Kemdg$ in \cref{Kemdg-functor} is \cref{gspectra-thm-v} applied to the multifunctor $F = \Kem$ between closed multicategories.
\begin{itemize}
\item $\permcatsu$ is a closed multicategory by \cref{perm-closed-multicat}.
\item The symmetric monoidal closed category $\Sp$ of symmetric spectra is a closed multicategory by \cref{smclosed-closed-multicat}.  
\end{itemize}  
The existence of the functor in \cref{Kemdg-presheaf} is \cref{gspectra-thm-v-cor} applied to $\Kem$.
\end{proof}

We describe the functor $\Kemdg$ in more detail in \cref{expl:Kemdg-object,expl:Kemdg-morphism} after \cref{rk:Kemdg-opposite,rk:BO7.5}.  In \cref{sec:mult-mackey-spectra} we factor the functors in  \cref{Kemdg-functor,Kemdg-presheaf} through enriched diagram and Mackey functor categories defined on left $\Mone$-modules, $\Gstar$-categories, and $\Gstar$-simplicial sets.

\subsection*{Symmetry and Opposite}

\begin{remark}[Symmetry of $\Kem$]\label{rk:Kemdg-opposite}
For the presheaf change-of-enrichment functor $\Kemdg$ in \cref{Kemdg-presheaf}, it is crucial that $\Kem$ is a multifunctor in the \emph{symmetric} sense in order to identify the $\Sp$-categories 
\[(\Cop)_{\Kem} \andspace (\C_{\Kem})^\op\] 
using \cref{dF-opposite}.  Without this symmetry property of $\Kem$, we would have to use $(\Cop)_{\Kem}$ in the codomain of $\Kemdg$, which is \cref{Kemdg-functor} for $\Cop$.  Thus, the fact that the codomain of $\Kemdg$ in \cref{Kemdg-presheaf} is a category of spectrally enriched \emph{presheaves}---as opposed to enriched diagrams---depends on the symmetry of $\Kem$.
\end{remark}

\begin{remark}[Symmetric and Non-Symmetric $K$-Theory Multifunctors]\label{rk:BO7.5}
  The diagram change-of-enrichment functor $\Kemdg$ in \cref{Kemdg-functor} is the $\Kem$ variant of the main result in \cite[Theorem 7.5]{bohmann_osorno-mackey}, which states the following using our notation:
  \begin{quote}
    Let $G$ be a finite group.  Then there is a functor
    \begin{equation}\label{bKdg}
      \bKdg \cn \permcatsucat\big(\GEop, \permcatsu\big) \to \Spcat\big((\GEop)_{\bbK} \scs \Sp\big).
    \end{equation}
  \end{quote}
  In \cite{bohmann_osorno-mackey} $\bKdg$, $\permcatsu$, and $\Sp$ are denoted $K_G$, $Perm$, and $Spec$, respectively.
  As in \cref{theorem:GM,rk:BO-Phi,ex:guillou-may}, $\bbK$ is the $K$-theory \emph{non-symmetric} multifunctor in \cite{guillou_may,gmmo}, and $\GE$ is the permutative Burnside category (\cref{definition:Burnside2}).
  The diagram change-of-enrichment functor $\bKdg$ in \cref{bKdg} exists by \cref{gspectra-thm-v} applied to
\begin{itemize}
\item the non-symmetric multifunctor $\bbK$ and
\item the $\permcatsu$-category $\GEop$, which is the opposite $\permcatsu$-category of $\GE$ (\cref{opposite-mcat}). 
\end{itemize} 

We emphasize the following regarding the diagram change-of-enrichment functor $\bKdg$ in \cref{bKdg}.
\begin{romenumerate}
\item Unlike Elmendorf-Mandell $K$-theory, $\bbK$ does not preserve symmetry \cite[Theorem 8.12]{gmmo}. 
  Thus we \emph{cannot} use \cref{dF-opposite} to identify the $\Sp$-categories
  \begin{itemize}
  \item $(\GEop)_{\bbK}$, with opposite taken in $\permcatsucat$, and
  \item $(\GE_{\bbK})^\op$, with opposite taken in $\Spcat$.
  \end{itemize} 
\item Recall from \cref{remark:GE-non-self-dual} that the $\permcatsu$-category $\GE$ and its opposite $\GEop$ are \emph{not} known to be equivalent as $\permcatsu$-categories.
  The assignment that sends a span of $G$-sets $(f,g)$ to $(g,f)$ as in \cref{eq:GB-self-dual} does not define a $\permcatsu$-functor 
  \[\GE \to \GEop\] 
  because it does not preserve composition as in \cref{mfunctor-diagrams}.

\item The $\Sp$-category $\GE_{\bbK}$ and its opposite $(\GE_{\bbK})^\op$ are also \emph{not} known to be equivalent.
\end{romenumerate}
Thus the codomain of $\bKdg$ has to be $(\GEop)_{\bbK}$ and \emph{not} $(\GE_{\bbK})^\op$ as stated in \cite[7.5]{bohmann_osorno-mackey}, where $(\GE_{\bbK})^\op$ is denoted $\GBop$.  As far as the authors know, the codomain of $\bKdg$ is a category of spectrally enriched diagrams but \emph{not} spectrally enriched presheaves as in the Guillou-May Quillen equivalence \cref{guillou-may-eq}.
\end{remark}

\subsection*{Description of $\Kemdg$}\index{enrichment!diagram change of -!Elmendorf-Mandell $K$-theory}\index{change of enrichment!diagram!Elmendorf-Mandell $K$-theory}\index{diagram!change of enrichment!Elmendorf-Mandell $K$-theory}\index{Elmendorf-Mandell!K-theory@$K$-theory!diagram change of enrichment}\index{K-theory@$K$-theory!Elmendorf-Mandell!diagram change of enrichment}

The rest of this section describes the functor $\Kemdg$ in \cref{Kemdg-functor} in detail.  The same discussion also applies to $\Cop$ and $(\Cop)_{\Kem} = (\C_{\Kem})^\op$; see \cref{rk:Kemdg-opposite}.  As defined in \cref{diag-change-enr-assign}, the diagram change-of-enrichment functor $\Kemdg$ 
\begin{itemize}
\item first applies the change-of-enrichment 2-functor (\cref{mult-change-enrichment})
\begin{equation}\label{Kem-change-enr}
\dKem \cn \permcatsucat \to \Spcat.
\end{equation}
along $\Kem$ and then
\item composes or whiskers with the standard enrichment $\Sp$-functor $\Kemse$ (\cref{expl:Kemse}).
\end{itemize}

\begin{explanation}[$\Kemdg$ on Objects]\label{expl:Kemdg-object}
Consider a $\permcatsu$-functor
\[A \cn \C \to \permcatsu.\]
Applying \cref{expl:Fdg-object} to the context of \cref{thm:Kemdg}, the $\Sp$-functor $\Kemdg A$ is the following composite.
\begin{equation}\label{KemdgX}
\begin{tikzpicture}[baseline={(a.base)}]
\def\g{3} \def\w{.6}
\draw[0cell=1]
(0,0) node (a) {\C_{\Kem}}
(a)+(\g,0) node (b) {(\permcatsu)_{\Kem}}
(b)+(\g,0) node (c) {\Sp}
;
\draw[1cell=.9]
(a) edge node {A_{\Kem}} (b)
(b) edge node {\Kemse} (c)
;
\draw[1cell=.9]
(a) [rounded corners=3pt] |- ($(b)+(-1,\w)$)
-- node {\Kemdg A} ($(b)+(1,\w)$) -| (c)
;
\end{tikzpicture}
\end{equation}
\begin{itemize}
\item $A_{\Kem}$ is the image of $A$ under the change-of-enrichment 2-functor $\dKem$ in \cref{Kem-change-enr}.
\item $\Kemse$ is the standard enrichment of $\Kem$ in \cref{expl:Kemse}.
\end{itemize}
Next we describe its object assignment and component morphisms.

\medskip
\emph{Object Assignment}. 
For an object $x \in \C$, the object assignment is
\[(\Kemdg A)(x) = \Kem(Ax) \inspace \Sp.\]
This is the Elmendorf-Mandell $K$-theory of the small permutative category $Ax$.

\medskip
\emph{Components}. 
For objects $x,y \in \C$, the component $(\Kemdg A)_{x,y}$ is the following composite morphism in $\Sp$.
\[\begin{tikzpicture}[vcenter]
\draw[0cell=.9]
(0,0) node (a) {\Kem \C(x,y)}
(a)+(0,-1.3) node (b) {\Kem\clp(Ax; Ay)}
(a)+(4.5,0) node (c) {\clsp\big(\Kem (Ax); \Kem (Ay)\big)}
;
\draw[1cell=.9]
(a) edge node {(\Kemdg A)_{x,y}} (c)
(a) edge node[swap] {\Kem (A_{x,y})} (b)
(b) edge[bend right=15] node[swap,pos=.7] {\Kemse_{Ax,Ay}} (c)
;
\end{tikzpicture}\]
The adjoint of $(\Kemdg A)_{x,y}$ in $\Sp$ is the following composite morphism.
\begin{equation}\label{KemdgAxy-adj}
\begin{tikzpicture}[vcenter]
\draw[0cell=.9]
(0,0) node (a) {\Kem \C(x,y) \sma \Kem (Ax)}
(a)+(0,-1.3) node (b) {\Kem\clp(Ax; Ay) \sma \Kem Ax}
(a)+(4,0) node (c) {\Kem (Ay)}
;
\draw[1cell=.9]
(a) edge node {\Kem(\pn{A_{x,y}})} (c)
(a) edge node[swap] {\Kem (A_{x,y}) \sma 1} (b)
(b) edge[bend right=15] node[swap,pos=.7] {\Kem(\ev_{Ax;\, Ay})} (c)
;
\end{tikzpicture}
\end{equation}
Here $\pn{A_{x,y}}$ is the partner of $A_{x,y}$.  It is defined in \cref{Axy-partner}.
\end{explanation}

\begin{explanation}[$\Kemdg$ on Morphisms]\label{expl:Kemdg-morphism}
Consider a $\permcatsu$-natural transformation (\cref{expl:perm-nattr}) $\psi$ as follows.
\[\begin{tikzpicture}[baseline={(a.base)}]
\def\t{25}
\draw[0cell]
(0,0) node (a') {\C}
(a')+(0,0) node (a) {\phantom{X}}
(a)+(1.8,0) node (b) {\phantom{X}}
(b)+(.65,0.02) node (b') {\permcatsu}
;
\draw[1cell=.8]
(a) edge[bend left=\t] node {A} (b)
(a) edge[bend right=\t] node[swap] {B} (b)
;
\draw[2cell]
node[between=a and b at .44, rotate=-90, 2label={above,\psi}] {\Rightarrow}
;
\end{tikzpicture}\]
For each object $x \in \C$, the $x$-component of $\psi$ is a nullary multimorphism
\[\psi_x \cn \ang{} \to \clp\scmap{Ax;Bx} \inspace \permcatsu\]
with $\ang{}$ denoting the empty sequence.  This means a choice of an object in $\clp\scmap{Ax;Bx}$.  In other words, $\psi_x$ is a strictly unital symmetric monoidal functor 
\begin{equation}\label{psix-sharp}
\psi_x \cn Ax \to Bx \inspace \permcatsu.
\end{equation}

Applying \cref{expl:Fdg-morphism} to the context of \cref{thm:Kemdg}, the image of $\psi$ under the diagram change-of-enrichment functor $\Kemdg$, 
\[\Kemdg \psi \cn \Kemdg A = \Kemse \circ A_{\Kem} \to \Kemdg B = \Kemse \circ B_{\Kem},\]
is the following whiskering of the $\Sp$-natural transformation $\psi_{\Kem}$ with the $\Sp$-functor $\Kemse$ in \cref{expl:Kemse}.
\[\begin{tikzpicture}[baseline={(a.base)}]
\def\t{25} \def\d{.8}
\draw[0cell]
(0,0) node (a') {\C_{\Kem}}
(a')+(.25,0) node (a) {\phantom{X}}
(a)+(2,0) node (b) {\phantom{X}}
(b)+(1.08,0) node (b') {(\permcatsu)_{\Kem}}
(b')+(2.7,0) node (c) {\Sp}
;
\draw[1cell=.8]
(a) edge[bend left=30] node {A_{\Kem}} (b)
(a) edge[bend right=30] node[swap] {B_{\Kem}} (b)
(b') edge node {\Kemse} (c)
;
\draw[2cell]
node[between=a and b at .36, rotate=-90, 2label={above,\psi_{\Kem}}] {\Rightarrow}
;
\end{tikzpicture}\]
For each object $x \in \C$, the $x$-component of $\Kemdg \psi$ is a nullary multimorphism
\[(\Kemdg \psi)_x \cn \ang{} \to \clsp\big(\Kem(Ax); \Kem(Bx)\big) \inspace \Sp.\]
Since the multicategory structure on $\Sp$ is induced by a symmetric monoidal structure, this nullary multimorphism is a morphism
\[(\Kemdg \psi)_x \cn S \to \clsp\big(\Kem(Ax); \Kem(Bx)\big) \inspace \Sp,\]
where $S$ denotes the symmetric sphere spectrum \cite[7.4.1]{cerberusIII}.  The adjoint of this morphism in $\Sp$ is the composite
\begin{equation}\label{Kemdgpsix-adj}
\begin{tikzpicture}[baseline={(a.base)}]
\def\w{.6}
\draw[0cell=1]
(0,0) node (a) {S \sma \Kem(Ax)}
(a)+(2.7,0) node (b) {\Kem(Ax)}
(b)+(3.3,0) node (c) {\Kem(Bx)}
;
\draw[1cell=.9]
(a) edge node {\iso} (b)
(b) edge node {\Kem(\psi_x)} (c)
;
\draw[1cell=.9]
(a) [rounded corners=3pt] |- ($(b)+(-.7,\w)$)
-- node {\pn{(\Kemdg\psi)_x}} ($(b)+(1.3,\w)$) -| (c)
;
\end{tikzpicture}
\end{equation}
with $\psi_x$ in \cref{psix-sharp}.
\end{explanation}

\section{Spectral Mackey Functors from Multicategorical Mackey Functors}
\label{sec:mult-mackey-spectra}

In this section we apply \cref{gspectra-thm-vii,gspectra-thm-vii-cor} along with the factorization of Elmendorf-Mandell $K$-theory to factor the diagram and presheaf change-of-enrichment functors $\Kemdg$ in \cref{thm:Kemdg}.  For the context, recall from \cref{Kem} that Elmendorf-Mandell $K$-theory $\Kem$ is the following composite.
\begin{equation}\label{Kem-factor-four}
\begin{tikzpicture}[vcenter]
\def\v{-1.3} \def\h{6}
\draw[0cell=.9]
(0,0) node (p) {\permcatsu}
(p)+(\h,0) node (sp) {\Sp}
(p)+(0,\v) node (m) {\MoneMod}
(m)+(\h/2,0) node (gc) {\Gstarcat}
(sp)+(0,\v) node (gs) {\Gstarsset}
;
\draw[1cell=.9]
(p) edge node {\Kem} (sp)
(p) edge node[swap] {\Endm} (m)
(m) edge node {\Jt} (gc)
(gc) edge node {\Ner_*} (gs)
(gs) edge node[swap] {\Kg} (sp)
;
\end{tikzpicture}
\end{equation}
As we explain under \cref{Kem-four}, the commutative diagram \cref{Kem-factor-four} consists of multifunctors between closed multicategories.

\begin{theorem}\label{gspectra-thm-xiv}
Suppose $\C$ is a small $\permcatsu$-category.  Then the factorization \cref{Kem-factor-four} of $\Kem$ induces the following factorization of the diagram change-of-enrichment functor $\Kemdg$ in \cref{Kemdg-functor}.
\begin{equation}\label{Kemdg-factor-four}
\begin{tikzpicture}[vcenter]
\def\v{-1.3} \def\h{5}
\draw[0cell=.8]
(0,0) node (p) {\permcatsucat\big(\C, \permcatsu\big)}
(p)+(\h,0) node (sp) {\Spcat\big(\C_{\Kem}, \Sp\big)}
(p)+(0,\v) node (m) {\Monemodcat\Big(\C_{\Endm} \scs \MoneMod \Big)}
(m)+(\h/2,-1.2) node (gc) {\Gstarcatcat\Big(\C_{\Jt\Endm} \scs \Gstarcat \Big)}
(sp)+(0,\v) node (gs) {\Gstarssetcat\Big(\C_{\Ner_* \Jt\Endm} \scs \Gstarsset \Big)}
;
\draw[1cell=.9]
(p) edge node {\Kemdg} (sp)
(p) edge node[swap] {\Endmdg} (m)
(m) edge[transform canvas={xshift=-1.5em}, shorten >=-1ex] node[swap,pos=.2] {\Jtdg} (gc)
(gc) edge[transform canvas={xshift=1.5em}, shorten <=-1ex] node[swap,pos=.8] {\Nerstardg} (gs)
(gs) edge node[swap] {\Kgdg} (sp)
;
\end{tikzpicture}
\end{equation}
Moreover, there is an analogous factorization of the presheaf change-of-enrichment functor $\Kemdg$ in \cref{Kemdg-presheaf} with $\C$ and each $\C_?$ in \cref{Kemdg-factor-four} replaced by $\Cop$ and $(\C_?)^\op$, respectively. 
\end{theorem}

\begin{proof}
The factorization \cref{Kemdg-factor-four} is the result of applying \cref{gspectra-thm-vii} three times to the commutative diagram \cref{Kem-factor-four}.  The second assertion about $\Cop$ follows similarly from \cref{gspectra-thm-vii-cor}.
\end{proof}

\begin{explanation}[Diagram Change of Enrichment]\label{expl:Kemdg-change-enr}\index{functor!diagram change-of-enrichment}\index{diagram!change-of-enrichment functor}\index{compositionality!- of diagram change of enrichment}\index{functor!diagram change-of-enrichment!compositionality of}\index{diagram!change-of-enrichment functor!compositionality}
By definition \cref{diag-change-enr-assign} each diagram change-of-enrichment functor $\Fdg$ is given by
\begin{itemize}
\item applying the change-of-enrichment 2-functor $\dF$ (\cref{mult-change-enrichment}) and then
\item composing or whiskering the result with the standard enrichment $\Fse$ in \cref{gspectra-thm-iii}.
\end{itemize}
In \cref{expl:Kemdg-object,expl:Kemdg-morphism} we describe the diagram change-of-enrichment functor $\Kemdg$.  The other four diagram change-of-enrichment functors in \cref{Kemdg-factor-four} admit analogous description.  Just like $\Sp$, each of the closed multicategories 
\[\MoneMod, \quad \Gstarcat, \andspace \Gstarsset\]
is a symmetric monoidal closed category by \cref{proposition:EM2-5-1,Gstar-V}.  In each case we can use the adjunction between the monoidal product and the closed structure to obtain the analogs of the adjoint description in \cref{KemdgAxy-adj,Kemdgpsix-adj}.
\end{explanation}

\begin{explanation}[Multicategorical Mackey Functors]\label{expl:Monediagrams}
The second assertion of \cref{gspectra-thm-xiv} gives, for each small $\permcatsu$-category $\C$, the following factorization of the presheaf change-of-enrichment functor $\Kemdg$ in \cref{Kemdg-presheaf}.
\begin{equation}\label{Kemdg-factor-four-op}
\begin{tikzpicture}[vcenter]
\def\v{-1.3} \def\h{5.2}
\draw[0cell=.8]
(0,0) node (p) {\permcatsucat\big(\Cop, \permcatsu\big)}
(p)+(\h,0) node (sp) {\Spcat\big((\C_{\Kem})^\op, \Sp\big)}
(p)+(0,\v) node (m) {\Monemodcat\Big(\big(\C_{\Endm})^\op \scs \MoneMod \Big)}
(m)+(\h/2,-1.2) node (gc) {\Gstarcatcat\Big(\big(\C_{\Jt\Endm})^\op \scs \Gstarcat \Big)}
(sp)+(0,\v) node (gs) {\Gstarssetcat\Big(\big(\C_{\Ner_* \Jt\Endm})^\op \scs \Gstarsset \Big)}
;
\draw[1cell=.9]
(p) edge node {\Kemdg} (sp)
(p) edge node[swap] {\Endmdg} (m)
(m) edge[transform canvas={xshift=-1.5em}, shorten >=-1ex] node[swap,pos=.2] {\Jtdg} (gc)
(gc) edge[transform canvas={xshift=1.5em}, shorten <=-1ex] node[swap,pos=.8] {\Nerstardg} (gs)
(gs) edge node[swap] {\Kgdg} (sp)
;
\end{tikzpicture}
\end{equation}
Via the composite $\Kgdg \circ \Nerstardg \circ \Jtdg$, we obtain $\Sp$-enriched Mackey functors from $\MoneMod$-enriched Mackey functors.  This is what the title of this section refers to.  
\end{explanation}

\begin{explanation}[Equivalences of Homotopy Theories]\label{expl:Endmdg-heq}
In \cref{mackey-xiv-mone} we observe that $\Endmdg$ in each of \cref{Kemdg-factor-four,Kemdg-factor-four-op} is an equivalence of homotopy theories.  This implies that the homotopy theories of modules in $\permcatsu$ and in $\MoneMod$ are equivalent via $\Endmdg$.  An analogous equivalence of homotopy theories also holds with $\MoneMod$ replaced by $\pMulticat$; see \cref{mackey-xiv-pmulticat}.  These equivalences of homotopy theories are all instances of the much more general \cref{mackey-gen-xiv,mackey-xiv-cor}, which hold at the level of (non-symmetric) multifunctors between (non-symmetric) closed multicategories.  We do \emph{not} know whether $\Kemdg$, $\Jtdg$, $\Nerstardg$, $\Kgdg$, or any of their composites are equivalences of homotopy theories or not.  See \cref{qu:Kemdg}.
\end{explanation}

\part{Homotopy Theory of Enriched Diagrams and Mackey Functors}
\label{part:homotopy-mackey}

\chapter[Homotopy Equivalences between Enriched Diagram Categories]{Homotopy Equivalences between Enriched Diagram and Mackey Functor Categories}
\label{ch:mackey}
Throughout this chapter we suppose given a pair of non-symmetric multifunctors between non-symmetric closed multicategories
\[F \cn \M \lradj \N \cn E,\]
together with multinatural transformations as follows.
\[\begin{tikzpicture}[baseline={(a.base)}]
    \def\t{25} \def\s{27}
    \draw[0cell]
    (0,0) node (a) {\M}
    (a)+(1.7,0) node (b) {\M}
    ;
    \draw[1cell=.9]
    (a) edge[bend left=\t] node {1_\M} (b)
    (a) edge[bend right=\t] node[swap] {EF} (b)
    ;
    \draw[2cell]
    node[between=a and b at .44, rotate=-90, 2label={above,\uni}] {\Rightarrow}
    ;
    \begin{scope}[shift={(3.5,0)}]
      \draw[0cell]
      (0,0) node (a) {\N}
      (a)+(1.7,0) node (b) {\N}
      ;
      \draw[1cell=.9]
      (a) edge[bend left=\s] node {1_\N} (b)
      (a) edge[bend right=\s] node[swap] {FE} (b)
      ;
      \draw[2cell]
      node[between=a and b at .44, rotate=-90, 2label={above,\cou}] {\Rightarrow}
      ;
    \end{scope}
  \end{tikzpicture}\]
The main results of this chapter identify conditions under which these data will induce, for each small $\N$-category $\C$, inverse equivalences of homotopy theories between enriched diagram categories (\cref{mackey-gen-xiv})
\[\begin{tikzpicture}[baseline={(a.base)}]
    \draw[0cell]
    (0,0) node (a) {\MCat(\C_E,\M)}
    (a)+(3.6,0) node (b) {\NCat(\C,\N)}
    ;
    \draw[1cell]
    (a) edge[transform canvas={yshift=.6ex}] node {\Fdgr} (b)
    (b) edge[transform canvas={yshift=-.4ex}] node {\Edg} (a) 
    ;
  \end{tikzpicture}\]
and enriched Mackey functor categories (\cref{mackey-xiv-cor})
\[\begin{tikzpicture}[baseline={(a.base)}]
    \draw[0cell]
    (0,0) node (a) {\MCat((\C_E)^\op,\M)}
    (a)+(3.6,0) node (b) {\NCat(\C^\op,\N),}
    ;
    \draw[1cell]
    (a) edge[transform canvas={yshift=.6ex}] node {\Fdgr} (b)
    (b) edge[transform canvas={yshift=-.4ex}] node {\Edg} (a) 
    ;
  \end{tikzpicture}\]
with the latter requiring that $E$ be a multifunctor, not merely a non-symmetric multifunctor.
The functors $\Edg$ are the diagram change of enrichment for $E$ (\cref{gspectra-thm-v}).
The functors $\Fdgr$ are similar, and are described in \cref{sec:comparing-enr-diag-cat}.

It is important to note that, while $E$ is required to satisfy the additional symmetric group action axiom of a multifunctor in \cref{mackey-xiv-cor}, $F$ is not.
In the applications, \cref{mackey-xiv-pmulticat,mackey-xiv-mone} below, $E$ is an endomorphism multifunctor and $F$ is a corresponding free non-symmetric multifunctor.

\subsection*{Connection with Other Chapters}
The results in this chapter provide a general approach to three main results in \cref{ch:mackey_eq}, as follows.
\begin{itemize}
\item The application $\brb{F,E} = \brb{\Fst,\Est}$ and $(\uni\scs\cou) = (\etast\scs\vrhost)$ is described in \cref{sec:mackey-pmulticat,sec:Endstdg,sec:Fstdgr}.
\item The application $\brb{F,E} = \brb{\Fm,\Em}$ and $(\uni\scs\cou) = (\etam\scs\vrhom)$ is described in \cref{sec:mackey-mone,sec:Endmdg-Fmdgr}.
\item The application $\brb{F,E} = \brb{\Mone \sma -, \Um}$ and $(\uni\scs\cou) = (\etahat\scs \epzhatinv)$ is described in \cref{sec:mackey-pmult-mone,sec:Monesmadgr-Umdg}.
\end{itemize}
For each application, we recall further background and context in the respective sections of \cref{ch:mackey_eq}.

\subsection*{Background}

\Cref{sec:enr-diag-psh} describes enriched diagram categories and enriched Mackey functor categories.
\Cref{sec:enr-diag-change-enr,sec:diag-psh-change-enr-functors,sec:fun-change-enr-diag}
give the definitions and basic properties of the functors $\Edg$ and $\Fdg$ for diagram change of enrichment.
That material, and many of the details later in this chapter, depends on the basic theory of (self-)enrichment in a non-symmetric multicategory, from \cref{ch:menriched,ch:change_enr,ch:std_enrich}.

\subsection*{Chapter Summary}

The precise context and assumptions for this chapter are given in \cref{def:mackey-gen-context}.
The remainder of \cref{sec:comparing-enr-diag-cat} gives the definition and further explanations of the functor $\Fdgr$.
\Cref{sec:unidg,sec:coudg} construct two natural transformations, denoted $\unidg$ and $\coudg$, that compare the composites $\Edg\Fdgr$ and $\Fdgr\Edg$ to the respective identity functors.
The main result, that $\Fdgr$ and $\Edg$ are inverse equivalences of homotopy theories in the sense of \cref{def:inverse-heq}, is in \cref{sec:heq-enrdiag-cat}; see \cref{mackey-gen-xiv}.
With the further assumption that $E$ is a multifunctor (not merely a non-symmetric multifunctor), \cref{mackey-xiv-cor} yields inverse equivalences of homotopy theories between enriched Mackey functor categories.
Here is a summary table.
\reftable{.9}{
  definition and explanations of $\Fdgr$
  & \ref{def:mackey-gen-context}, \ref{expl:Fdgr}, and \ref{expl:Fdgr-morphisms}
  \\ \hline
  $\unidg \cn 1 \to \Edg\Fdgr$
  & \ref{def:unidg}, \ref{unidgA-natural}, and \ref{unidg-natural}
  \\ \hline
  $\coudg \cn 1 \to \Fdgr\Edg$
  & \ref{def:coudg}, \ref{coudgP-natural}, and \ref{coudg-natural}
  \\ \hline
  componentwise stable equivalences for enriched diagrams
  & \ref{def:enr-diagcat-relative}, \ref{def:mackey-heq-context}, \ref{cWtri}, and \ref{MFinvcX}
  \\ \hline
  $\brb{\Fdgr,\Edg}$ inverse equivalence of homotopy theories
  & \ref{mackey-gen-xiv} and \ref{mackey-xiv-cor}
  \\ \hline
  $\brb{F,E}$ inverse equivalence of homotopy theories
  & \ref{FE-inv-heq}
  \\
}

\section[Comparing Enriched Diagram Categories]{Comparing Enriched Diagram and Mackey Functor Categories}
\label{sec:comparing-enr-diag-cat}

In this section we define a pair of functors
\[\begin{tikzpicture}[baseline={(a.base)}]
\draw[0cell]
(0,0) node (a) {\MCat(\C_E,\M)}
(a)+(3.6,0) node (b) {\NCat(\C,\N)}
;
\draw[1cell]
(a) edge[transform canvas={yshift=.6ex}] node {\Fdgr} (b)
(b) edge[transform canvas={yshift=-.4ex}] node {\Edg} (a) 
;
\end{tikzpicture}\]
that compare two categories of enriched diagrams in the sense of \cref{def:enr-diag-cat}.  The functor $\Edg$ is the diagram change of enrichment of $E \cn \N \to \M$ (\cref{gspectra-thm-v}).  The functor $\Fdgr$ and the context for this chapter are discussed in \cref{def:mackey-gen-context}.  After that definition we unravel the functor $\Fdgr$.
\begin{itemize}
\item \cref{expl:Fdgr} describes $\Fdgr$ on objects.
\item \cref{expl:Fdgr-morphisms} describes $\Fdgr$ on morphisms.
\end{itemize}

\subsection*{Defining the Functor $\Fdgr$}

\begin{definition}\label{def:mackey-gen-context}
Suppose given the data \cref{mackey-context-i,mackey-context-ii,mackey-context-iii,mackey-context-iv,mackey-context-v} below.
\begin{romenumerate}
\item\label{mackey-context-i} $\big(\M,\clM,\ev^\M\big)$ and $\big(\N,\clN,\ev^\N\big)$ are non-symmetric closed multicategories (\cref{def:closed-multicat}).
\item\label{mackey-context-ii} $(\C,\mcomp,i)$ is a small $\N$-category (\cref{def:menriched-cat}).
\item\label{mackey-context-iii}  $F$ and $E$ are non-symmetric multifunctors (\cref{def:enr-multicategory-functor}) as follows.
\[F \cn \M \lradj \N \cn E\] 
\item\label{mackey-context-iv} $\uni$ and $\cou$ are multinatural transformations (\cref{def:enr-multicat-natural-transformation}) as follows.
\[
\]
\item\label{mackey-context-v} For each pair of objects $x,y\in \C$, the two unary multimorphisms in $\M$ 
\begin{equation}\label{counit-compatibility}
%
\end{equation}
are equal, where $\C(x,y) \in \Ob\N$ is a morphism object of $\C$.
\end{romenumerate}
We define the functor $\Fdgr$ as the composite
\begin{equation}\label{Fdgr-def}
%
\end{equation}
involving the following data.
\begin{itemize}
\item Each of $\M$ and $\N$ is equipped with the canonical self-enrichment (\cref{cl-multi-cl-cat}).
\item The change-of-enrichment 2-functors (\cref{mult-change-enrichment,func-change-enr})
\begin{equation}\label{Fdg-Edg}
%
\end{equation}
are induced by $F$, $E$, and $FE$.  The $\M$-category $\C_E$ and the $\N$-category $\C_{FE}$ are obtained from $\C$ by applying, respectively, $\dE$ and $\dFE$.
\item $\Fdg$ in \cref{Fdgr-def} is the diagram change-of-enrichment functor of $F$ at the $\M$-category $\C_E$ (\cref{gspectra-thm-v}).
\item The $\N$-functor 
\begin{equation}\label{C-sub-cou}
\C_\cou \cn \C \to \C_{FE}
\end{equation}
is the $\C$-component \cref{C-sub-theta} of the 2-natural transformation
\[%
\]
induced by $\cou \cn 1_\N \to FE$ (\cref{dtheta-twonat}).
\item The functor $\C_\cou^*$ in \cref{Fdgr-def} is given by pre-composing and whiskering with the $\N$-functor $\C_\cou$ in \cref{C-sub-cou}.
\end{itemize}
This finishes the definition of the functor $\Fdgr$.
\end{definition}

\begin{remark}\label{rk:eta-not-used}
The functor $\Fdgr$ in \cref{Fdgr-def} does \emph{not} use the multinatural transformation $\uni$ in \cref{def:mackey-gen-context} \cref{mackey-context-iv} and condition \cref{mackey-context-v}.  For the discussion below, it is more convenient to state $\uni$ and $\cou$ together in one place.  We use $\uni$ in \cref{def:unidg} below.  We use condition \cref{mackey-context-v} in \cref{unidgA-natural-partners}; see also \cref{rk:diff-unidgA}.  Instances of condition \cref{mackey-context-v} include
\begin{itemize}
\item \cref{etaEEvrho}, which is used in the proof of \cref{mackey-xiv-pmulticat};
\item \cref{etavrhoMone}, which is used in the proof of \cref{mackey-xiv-mone}; and
\item the right triangle identity of the 2-adjunction $(\Monesma,\Um)$ in \cref{MonesmaUm-iiadj}.\defmark
\end{itemize}
\end{remark}

In \cref{expl:Fdgr,expl:Fdgr-morphisms}, we explicitly describe the object and morphism assignments of the functor
\[\Fdgr \cn \MCat(\C_E,\M) \fto{\Fdg} \NCat(\C_{FE},\N) \fto{\C_\cou^*} \NCat(\C,\N)\]
in \cref{Fdgr-def}.

\subsection*{Unraveling the Functor $\Fdgr$}

\begin{explanation}[$\Fdgr$ on Objects]\label{expl:Fdgr}
For an $\M$-functor $A \cn \C_E \to \M$ (\cref{def:mfunctor}), the $\N$-functor
\begin{equation}\label{Fdg-of-A}

\end{equation}
is, by \cref{def:enr-diag-change-enr}, the composite of two $\N$-functors.
\begin{itemize}
\item $A_F$ is the image of $A$ under the change-of-enrichment 2-functor $\dF$ (\cref{mult-change-enrichment}).
\item $\Fse$ is the standard enrichment of $F \cn \M \to \N$ (\cref{gspectra-thm-iii}).
\end{itemize}
Its object assignment is given by
\[(\Fdg A)x = FAx \forspace x \in \C.\]
For objects $x,y \in \C$, the $(x,y)$-component of $\Fdg A$ is the composite unary multimorphism in $\N$
\begin{equation}\label{FdgAxy-diagram}
%
\end{equation}
By definition \cref{Fprimexy},
\begin{equation}\label{Fse-AxAy}
\Fse_{Ax,Ay} = \pn{\big(F (\ev^\M_{Ax;\, Ay})\big)}
\end{equation}
is the partner of the binary multimorphism
\[F\big(\ev^\M_{Ax;\, Ay}\big) \cn \big( F\clM\scmap{Ax; Ay} \scs FAx \big) \to FAy \inspace \N.\]

By definition \cref{Fdgr-def} the $\N$-functor
\begin{equation}\label{FdgrA-diag}
%
\end{equation}
is the composite of $\C_\cou$ in \cref{C-sub-cou} and $\Fdg A$ in \cref{Fdg-of-A}.  Since $\C_\cou$ is the identity on objects, the object assignment of $\Fdgr A$ is given by
\begin{equation}\label{FdgrAx}
(\Fdgr A)x = FAx \forspace x \in \C.
\end{equation}
For objects $x,y \in \C$, the $(x,y)$-component of $\Fdgr A$ is the composite unary multimorphism in $\N$
\begin{equation}\label{FdgrA-component}
%
\end{equation}
of the $\C(x,y)$-component of $\cou$ and $(\Fdg A)_{x,y}$ in \cref{FdgAxy-diagram}.
\end{explanation}

\begin{explanation}[$\Fdgr$ on Morphisms]\label{expl:Fdgr-morphisms}
Consider an $\M$-natural transformation $\psi$ (\cref{def:mnaturaltr}) as in the left diagram below.
\begin{equation}\label{Fdgr-psi}
%
\end{equation}
Then the $\N$-natural transformation $\Fdgr\psi$, as in the right diagram in \cref{Fdgr-psi}, is given by the following whiskering, where $\psi_F$ is the image of $\psi$ under the change of enrichment $\dF$ (\cref{mult-change-enrichment}).
\begin{equation}\label{Fdgrpsi-def}
%
\end{equation}
This diagram is obtained from \cref{FdgrA-diag} by replacing $A$ by $\psi$.  

For each object $x \in \C$, since $\C_\cou$ is the identity on objects, the $x$-component of $\Fdgr\psi$ in \cref{Fdgrpsi-def} is the nullary multimorphism in $\N$ given by the composite
\begin{equation}\label{Fdgrpsi-component}
%
\end{equation}
precisely as in \cref{Fdgpsi-x}.
\begin{itemize}
\item $F\psi_x$ is the image under $F$ of the $x$-component of $\psi$, which is a nullary multimorphism
\[\ang{} \fto{\psi_x} \clM(Ax;Bx) \inspace \M.\]
\item The unary multimorphism
\[\Fse_{Ax,Bx} = \pn{\big(F(\ev^\M_{Ax;\, Bx})\big)}\]
is the partner of the binary multimorphism
\[F\big(\ev^\M_{Ax;\, Bx}\big) \cn \big(F\clM\scmap{Ax; Bx} \scs FAx\big) \to FBx \inspace \N.\]
\end{itemize}
As in \cref{Fdg-psix-pn-equality}, the partner of $(\Fdgr\psi)_x$ in \cref{Fdgrpsi-component} is the unary multimorphism
\begin{equation}\label{Fdgrpsix-partner}
\pn{(\Fdgr \psi)_x} = F(\pn{\psi_x}) \cn FAx \to FBx \inspace \N,
\end{equation}
where $\pn{\psi_x} \cn Ax \to Bx$ is the partner of $\psi_x$ \cref{psix-partner}.
\end{explanation}

In the context of \cref{def:mackey-gen-context}, there are two functors
\begin{equation}\label{Fdgr-Edg}
\begin{tikzpicture}[baseline={(a.base)}]
\draw[0cell]
(0,0) node (a) {\MCat(\C_E,\M)}
(a)+(3.6,0) node (b) {\NCat(\C,\N)}
;
\draw[1cell]
(a) edge[transform canvas={yshift=.6ex}] node {\Fdgr} (b)
(b) edge[transform canvas={yshift=-.4ex}] node {\Edg} (a) 
;
\end{tikzpicture}
\end{equation}
as follows.
\begin{itemize}
\item $\Fdgr = \C_\cou^* \Fdg$ is the functor in \cref{Fdgr-def}.
\item $\Edg$ is the diagram change-of-enrichment functor of $E \cn \N \to \M$ at the $\N$-category $\C$ (\cref{gspectra-thm-v}).
\end{itemize} 
In \cref{sec:unidg,sec:coudg} we relate the composite functors $\Edg\Fdgr$ and $\Fdgr\Edg$ to the respective identity functors using the multinatural transformations $\uni$ and $\cou$.

\section{Comparing \texorpdfstring{$\Edg\Fdgr$}{EF} and the Identity}
\label{sec:unidg}

Throughout this section we assume the same context as in \cref{def:mackey-gen-context}.  In this section we extend the multinatural transformation in \cref{def:mackey-gen-context} \cref{mackey-context-iv}
\[\begin{tikzpicture}[baseline={(a.base)}]
\def\t{25} \def\s{27}
\draw[0cell]
(0,0) node (a) {\M}
(a)+(1.7,0) node (b) {\M}
;
\draw[1cell=.9]
(a) edge[bend left=\t] node {1_\M} (b)
(a) edge[bend right=\t] node[swap] {EF} (b)
;
\draw[2cell]
node[between=a and b at .44, rotate=-90, 2label={above,\uni}] {\Rightarrow}
;
\end{tikzpicture}\]
to a natural transformation $\unidg$ comparing $\Edg\Fdgr$ in \cref{Fdgr-Edg} and the identity functor.  This section is organized as follows.
\begin{itemize}
\item $\unidg$ is in \cref{def:unidg}.
\item To show that $\unidg$ has the desired naturality properties, in \cref{expl:unidg,expl:EdgFdgr-objects,EdgFdgrA-commutes} we discuss the codomain $\Edg\Fdgr$ of $\unidg$ and its object assignment.
\item We describe $\Edg\Fdgr$ on morphisms in \cref{expl:EdgFdgr-morphisms,EdgFdgr-morphism-partner}.
\item We show that $\unidg$ is a natural transformation in \cref{unidgA-natural,unidg-natural}.
\end{itemize}

\subsection*{The Natural Transformation $\unidg$}

\begin{definition}\label{def:unidg}
In the context of \cref{def:mackey-gen-context,Fdgr-Edg}, we define the data of a natural transformation
\begin{equation}\label{unidg-def}
\begin{tikzpicture}[baseline={(a.base)}]
\def\s{25}
\draw[0cell]
(0,0) node (a) {\phantom{\M}}
(a)+(2,0) node (b) {\phantom{\N}}
(a)+(-.9,0) node (a') {\MCat(\C_E,\M)}
(b)+(.9,0) node (b') {\MCat(\C_E,\M)}
;
\draw[1cell=.9]
(a) edge[bend left=\s] node {1} (b)
(a) edge[bend right=\s] node[swap] {\Edg\Fdgr} (b)
;
\draw[2cell]
node[between=a and b at .45, rotate=-90, 2label={above,\unidg}] {\Rightarrow}
;
\end{tikzpicture}
\end{equation}
as follows.  For an $\M$-functor $A \cn \C_E \to \M$ (\cref{def:mfunctor}), the $A$-component of $\unidg$ is the $\M$-natural transformation (\cref{def:mnaturaltr})
\begin{equation}\label{unidg-A}
\begin{tikzpicture}[baseline={(a.base)}]
\def\s{25}
\draw[0cell]
(0,0) node (a) {\C_E}
(a)+(2.3,0) node (b) {\M}
;
\draw[1cell=.9]
(a) edge[bend left=\s] node {A} (b)
(a) edge[bend right=\s] node[swap] {\Edg\Fdgr A} (b)
;
\draw[2cell]
node[between=a and b at .42, rotate=-90, 2label={above,\unidg_A}] {\Rightarrow}
;
\end{tikzpicture}
\end{equation}
with, for each object $x \in \C$, $x$-component given by the nullary multimorphism
\begin{equation}\label{unidgAx}
(\unidg_A)_x = \pn{\uni_{Ax}} \cn \ang{} \to \clM\scmap{Ax;EFAx} \inspace \M.
\end{equation}
This is the partner \cref{eval-bij-zero} of the $(Ax)$-component 
\[\uni_{Ax} \cn Ax \to EFAx\]
of the multinatural transformation $\uni \cn 1_\M \to EF$, which is a unary multimorphism in $\M$.  This finishes the definition of $\unidg$.  We check that
\begin{itemize}
\item $\unidg_A$ is an $\M$-natural transformation in \cref{unidgA-natural} and 
\item $\unidg$ is a natural transformation in \cref{unidg-natural}.\defmark
\end{itemize} 
\end{definition}

Before we prove the $\M$-naturality of $\unidg_A$ and the naturality of $\unidg$, we first discuss the codomain $\Edg\Fdgr$ of $\unidg$ in detail in \cref{expl:unidg,expl:EdgFdgr-objects,expl:EdgFdgr-morphisms}.

\begin{explanation}[Codomain of $\unidg$]\label{expl:unidg}
In the context of \cref{def:mackey-gen-context,def:unidg}, by \cref{diag-change-enr-assign,FdgrA-diag} the codomain of $\unidg$ in \cref{unidg-def} is the composite of the functors
\[\begin{split}
\Fdgr &= \Fse \circ \dF \circ \C_\cou \andspace\\
\Edg &= \Ese \circ \dE
\end{split}\]
involving the following.
\begin{itemize}
\item The change-of-enrichment 2-functors $\dE$ and $\dF$ are as in \cref{Fdg-Edg}. 
\item $\C_\cou \cn \C \to \C_{FE}$ is the $\N$-functor in \cref{C-sub-cou}.
\item The standard enrichment
\[\Fse \cn \M_F \to \N \andspace \Ese \cn \N_E \to \M\]
are from in \cref{gspectra-thm-iii}.\defmark
\end{itemize} 
\end{explanation}

\begin{explanation}[$\Edg\Fdgr$ on Objects]\label{expl:EdgFdgr-objects} 
For an $\M$-functor $A \cn \C_E \to \M$ (\cref{def:mfunctor}), the functoriality of $\dE$ implies that the codomain of $\unidg_A$ in \cref{unidg-A} is the following $\M$-functor.
\[\begin{split}
\Edg \Fdgr A &= \Ese \circ (\Fdgr A)_E\\
&= \Ese \circ \big(\Fse \circ A_F \circ \C_\cou \big)_E\\
&= \Ese \circ \Fse_E \circ (A_F)_E \circ (\C_\cou)_E \cn \C_E \to \M
\end{split}\]
By \cref{func-change-enr}, $\Edg \Fdgr A$ is the following composite $\M$-functor.
\begin{equation}\label{EdgFdgr-A-diag}
\begin{tikzpicture}[vcenter]
\def\h{3} \def\v{1.3}
\draw[0cell=1]
(0,0) node (a) {\C_E}
(a)+(0,-\v) node (b) {\C_{EFE}}
(b)+(\h,0) node (c) {\M_{EF}}
(c)+(\h,0) node (d) {\N_E}
(d)+(0,\v) node (e) {\M}
;
\draw[1cell=.9]
(a) edge node[swap] {(\C_\cou)_E} (b)
(b) edge node {A_{EF}} (c)
(c) edge node {\Fse_E} (d)
(d) edge node[swap] {\Ese} (e)
(a) edge node {\Edg\Fdgr A} (e)
;
\end{tikzpicture}
\end{equation}
We explain its object assignment and components in \cref{EdgFdgrA-object,EdgFdgrA-component}, respectively.

\emph{Object Assignment}. 
By \cref{FAx} applied to $\Edg$ and \cref{FdgrAx}, the object assignment of $\Edg\Fdgr A$ in \cref{EdgFdgr-A-diag} is given by
\begin{equation}\label{EdgFdgrA-object}
(\Edg\Fdgr A)x = EFAx \forspace x \in \C.
\end{equation}
This implies that the $x$-component $(\unidg_A)_x$ in \cref{unidgAx} is well defined.

\emph{Components}. 
By \cref{EdgFdgr-A-diag}, for objects $x,y \in \C$, the $(x,y)$-component of $\Edg\Fdgr A$ is the composite unary multimorphism along the boundary of the following diagram in $\M$.
\begin{equation}\label{EdgFdgrA-component}
\begin{tikzpicture}[vcenter]
\def\g{-.5}\def\h{3} \def\v{1.5} \def\u{1}
\draw[0cell=.9]
(0,0) node (a) {E\C(x,y)}
(a)+(\g,-\v) node (b) {EFE\C(x,y)}
(b)+(\h,-\u) node (c) {EF\clM\scmap{Ax;Ay}}
(c)+(\h,\u) node (d) {E\clN\scmap{FAx;FAy}}
(d)+(\g,\v) node (e) {\clM\scmap{EFAx;EFAy}}
;
\draw[1cell=.9]
(a) edge node {(\Edg\Fdgr A)_{x,y}} (e)
(a) edge[transform canvas={xshift=-1em}] node[swap,pos=.8] {E\cou_{\C(x,y)}} (b)
(b) edge[transform canvas={xshift=-2em}, shorten >=-1.7em] node[swap,pos=.6] {EFA_{x,y}} (c)
(c) edge[transform canvas={xshift=2em}, shorten <=-1.7em] node[swap,pos=.4] {E\Fse_{Ax,Ay}} (d)
(d) edge[transform canvas={xshift=1em}] node[swap,pos=.1] {\Ese_{FAx,FAy}} (e)
;
\draw[1cell=.8]
(b) edge node[pos=.3] {\pn{\big(EF(\pn{A_{x,y}})\big)}} (e)
(c) edge node[pos=.15] {\pn{\big(EF(\ev^\M_{Ax;\, Ay})\big)}\!\!} (e)
;
\end{tikzpicture}
\end{equation}
By definition \cref{Fprimexy}, in the lower right arrow, $\Fse_{Ax,Ay}$ is as in \cref{Fse-AxAy}.  The upper right arrow
\[\Ese_{FAx,FAy} = \pn{\big(E(\ev^\N_{FAx;\, FAy}) \big)}\]
is the partner of the binary multimorphism
\[E\big(\ev^\N_{FAx;\, FAy} \big) \cn \big( E\clN\scmap{FAx;FAy} \scs EFAx \big) \to EFAy\]
in $\M$.  In \cref{EdgFdgrA-commutes} we discuss the two interior arrows in \cref{EdgFdgrA-component}.
\end{explanation}

\begin{lemma}\label{EdgFdgrA-commutes}
The diagram \cref{EdgFdgrA-component} is commutative.
\end{lemma}

\begin{proof}
Since the boundary of the diagram \cref{EdgFdgrA-component} commutes by \cref{EdgFdgr-A-diag}, it suffices to prove the following two equalities regarding its two interior arrows.
\begin{equation}\label{EFevM-pn}
\begin{aligned}
\ga^\M\scmap{\Ese_{FAx,FAy}; E\Fse_{Ax,Ay}} 
&= \pn{\big(EF(\ev^\M_{Ax;\, Ay})\big)}\\
\ga^\M\left(\pn{\big(EF(\ev^\M_{Ax;\, Ay})\big)} \sscs EFA_{x,y} \right) 
&= \pn{\big(EF(\pn{A_{x,y}})\big)}
\end{aligned}
\end{equation}
To prove these equalities, we consider the following diagram in $\M$.
\begin{equation}\label{EFevM-pn-diag}
\begin{tikzpicture}[vcenter]
\def\v{-1.3} \def\t{15}
\draw[0cell=.9]
(0,0) node (a) {\big(EFE\C(x,y) \scs EFAx\big)}
(a)+(0,\v) node (b) {\big(EF\clM\scmap{Ax;Ay} \scs EFAx\big)}
(b)+(0,\v) node (c) {\big(E\clN\scmap{FAx;FAy} \scs EFAx\big)}
(c)+(0,\v) node (d) {\big(\clM\scmap{EFAx;EFAy} \scs EFAx\big)}
(a)+(5,0) node (e1) {EFAy}
(e1)+(0,\v) node (e2) {EFAy}
(e2)+(0,\v) node (e3) {EFAy}
(e3)+(0,\v) node (e4) {EFAy}
;
\draw[1cell=.9]
(a) edge node[swap] {(EFA_{x,y} \scs 1)} (b)
(b) edge node[swap] {(E\Fse_{Ax,Ay} \scs 1)} (c)
(c) edge node[swap] {(\Ese_{FAx,FAy} \scs 1)} (d)
(a) edge node {EF(\pn{A_{x,y}})} (e1)
(b) edge node {EF(\ev^\M_{Ax;\, Ay})} (e2)
(c) edge node {E(\ev^\N_{FAx;\, FAy})} (e3)
(d) edge node {\ev^\M_{EFAx;\, EFAy}} (e4)
(e1) edge[-,double equal sign distance] (e2)
(e2) edge[-,double equal sign distance] (e3)
(e3) edge[-,double equal sign distance] (e4)
;
\end{tikzpicture}
\end{equation}
The three sub-regions in \cref{EFevM-pn-diag} are commutative for the following reasons.
\begin{itemize}
\item In the top sub-region, $\pn{A_{x,y}}$ is the partner \cref{eval-bijection} of $A_{x,y}$.  By definition the following diagram in $\M$ commutes.  
\begin{equation}\label{Axy-partner-diagram}
\begin{tikzpicture}[vcenter]
\draw[0cell=.9]
(0,0) node (a) {\big(E\C(x,y) \scs Ax\big)}
(a)+(0,-1.3) node (b) {\big(\clM\scmap{Ax;Ay} \scs Ax\big)}
(b)+(4,0) node (c) {Ay}
;
\draw[1cell=.9]
(a) edge node[swap] {(A_{x,y} \scs \opu)} (b)
(b) edge node[pos=.4] {\ev^\M_{Ax;\, Ay}} (c)
(a) edge[bend left=15] node {\pn{A_{x,y}}} (c)
;
\end{tikzpicture}
\end{equation}
Applying the non-symmetric multifunctor $EF$ to this commutative diagram yields the top sub-region in \cref{EFevM-pn-diag}.
\item The middle sub-region in \cref{EFevM-pn-diag} is obtained from the commutative diagram \cref{Fse-xy-diag} defining $\Fse_{Ax,Ay}$ by applying the non-symmetric multifunctor $E$.
\item The bottom sub-region in \cref{EFevM-pn-diag} is the commutative diagram \cref{Fse-xy-diag} that defines $\Ese_{FAx,FAy}$.
\end{itemize}

Since taking partner is a bijection \cref{eval-bijection}, the commutativity of the bottom two sub-regions in \cref{EFevM-pn-diag} proves the first desired equality in \cref{EFevM-pn}.  The second equality in \cref{EFevM-pn} follows from the first equality and the boundary of the commutative diagram \cref{EFevM-pn-diag}.
\end{proof}

\begin{explanation}[$\Edg\Fdgr$ on Morphisms]\label{expl:EdgFdgr-morphisms}
Suppose $\psi \cn A \to B$ is an $\M$-natural transformation (\cref{def:mnaturaltr}) as in the left diagram below.
\[
\]
By \cref{expl:unidg} the $\M$-natural transformation $\Edg\Fdgr\psi$, as in the right diagram above, is the following whiskering.
\begin{equation}\label{EdgFdgr-psi-diag}
%
\end{equation}
Since $(\C_\cou)_E$ has the identity object assignment, by \cref{mnat-change-enrichment,EdgFdgr-psi-diag}, for each object $x \in \C$ the $x$-component of $\Edg\Fdgr\psi$ is the following composite nullary multimorphism in $\M$.
\begin{equation}\label{EdgFdgr-psi-component}
%
\end{equation}
In \cref{EdgFdgr-psi-component},
\begin{itemize}
\item $\psi_x \cn \ang{} \to \clM(Ax;Bx)$ is the $x$-component of $\psi$ \cref{mnatural-component}, and
\item the component morphisms $\Fse_{Ax,Bx}$ and $\Ese_{FAx,FBx}$ are as in \cref{Fprimexy}.
\end{itemize}
Taking the partner \cref{eval-bijection} of $(\Edg\Fdgr\psi)_x$ yields the top unary multimorphism in $\M$ below.
\begin{equation}\label{EdgFdgr-psi-partner}
%
\end{equation}
The bottom arrow is the image under $EF$ of the unary multimorphism
\[\pn{\psi_x} \cn Ax \to Bx \inspace \M,\]
which is the partner of $\psi_x$.  \cref{EdgFdgr-morphism-partner} proves that they are the same.
\end{explanation}

\begin{lemma}\label{EdgFdgr-morphism-partner}
The two unary multimorphisms in \cref{EdgFdgr-psi-partner} are equal.
\end{lemma}

\begin{proof}
Since $(\Edg\Fdgr\psi)_x$ is the composite in \cref{EdgFdgr-psi-component}, its partner is, by definition \cref{eval-bijection}, the left-bottom composite in the following diagram in $\M$.
\[\begin{tikzpicture}
\def\h{5} \def\v{-1.3}
\draw[0cell=.9]
(0,0) node (a1) {\big(\ang{} \scs EFAx\big)}
(a1)+(0,\v) node (b1) {\big(EF\clM\scmap{Ax;Bx} \scs EFAx\big)}
(b1)+(0,\v) node (c1) {\big(E\clN\scmap{FAx;FBx} \scs EFAx\big)}
(c1)+(0,\v) node (d1) {\big(\clM\scmap{EFAx;EFBx} \scs EFAx\big)}
(a1)+(\h,0) node (a2) {EFBx}
(a2)+(0,\v) node (b2) {EFBx}
(b2)+(0,\v) node (c2) {EFBx}
(c2)+(0,\v) node (d2) {EFBx}
;
\draw[1cell=.9]
(a1) edge node {EF(\pn{\psi_x})} (a2)
(b1) edge node {EF(\ev^\M_{Ax;\, Bx})} (b2)
(c1) edge node {E(\ev^\N_{FAx;\, FBx})} (c2)
(d1) edge node {\ev^\M_{EFAx;\, EFBx}} (d2)
(a1) edge node[swap] {(EF(\psi_x) \scs \opu)} (b1)
(b1) edge node[swap] {(E\Fse_{Ax,Bx} \scs \opu)} (c1)
(c1) edge node[swap] {(\Ese_{FAx,FBx} \scs \opu)} (d1)
(a2) edge[-,double equal sign distance] (b2)
(b2) edge[-,double equal sign distance] (c2)
(c2) edge[-,double equal sign distance] (d2)
;
\end{tikzpicture}\]
The three sub-regions in the above diagram are commutative for the following reasons.
\begin{itemize}
\item The top sub-region is $EF$ applied to the commutative diagram in $\M$
\[\begin{tikzpicture}[vcenter]
\draw[0cell=.9]
(0,0) node (a) {\big(\ang{} \scs Ax\big)}
(a)+(0,-1.3) node (b) {\big(\clM\scmap{Ax;Bx} \scs Ax\big)}
(b)+(4,0) node (c) {Bx}
;
\draw[1cell=.9]
(a) edge node[swap] {(\psi_x \scs \opu)} (b)
(b) edge node[pos=.4] {\ev^\M_{Ax;\, Bx}} (c)
(a) edge[bend left=15] node {\pn{\psi_x}} (c)
;
\end{tikzpicture}\]
that defines the partner $\pn{\psi_x}$.
\item The middle sub-region is $E$ applied to the commutative diagram \cref{Fse-xy-diag} that defines $\Fse_{Ax,Bx}$.
\item The bottom sub-region is the commutative diagram \cref{Fse-xy-diag} that defines $\Ese_{FAx,FBx}$.
\end{itemize}
This proves that $\pn{(\Edg\Fdgr\psi)_x}$ is equal to $EF(\pn{\psi_x})$.
\end{proof}

\subsection*{Naturality of $\unidg$}

To prove that $\unidg$ is a natural transformation, we first check that its components are well defined.

\begin{lemma}\label{unidgA-natural}
In the context of \cref{def:mackey-gen-context,def:unidg}, for each $\M$-functor $A \cn \C_E \to \M$,
\[
\]
in \cref{unidg-A} is an $\M$-natural transformation.
\end{lemma}

\begin{proof}
We must prove that the naturality diagram \cref{mnaturality-diag} for $\unidg_A$, which is the diagram in $\M$ below for objects $x,y \in \C$, is commutative.
\begin{equation}\label{unidgA-nat-diag}
%
\end{equation}
Since taking partner in $\M$ is a bijection \cref{eval-bijection}, it suffices to show that the two composites in \cref{unidgA-nat-diag} have the same partner.  We compute these partners in \cref{unidgA-nat-top-right,unidgA-nat-left-bot}.  Then we observe that they are equal to finish the proof.

\emph{Top-Right Composite}. 
The partner of the top-right composite in \cref{unidgA-nat-diag} is, by definition, the left-bottom composite in the following diagram in $\M$.
\begin{equation}\label{unidgA-nat-top-right}
%
\end{equation}
The four sub-regions in \cref{unidgA-nat-top-right} are commutative for the following reasons.
\begin{itemize}
\item By \cref{def:mult-change-enr}, there is an equality of objects in $\M$ 
\[\C_E(x,y) = E\C(x,y).\]
The top trapezoid commutes by the naturality \cref{enr-multinat} of the multinatural transformation $\uni \cn 1_\M \to EF$.
\item The top left triangle commutes by the definition of $\pn{A_{x,y}}$ in \cref{Axy-partner-diagram} and the definition of the $y$-component of $\unidg_A$ in \cref{unidgAx},
\[(\unidg_A)_y = \pn{\uni_{Ay}} \cn \ang{} \to \clM\scmap{Ay; EFAy},\]
as the partner of the $(Ay)$-component of $\uni$. 
\item The bottom rectangle commutes by the definition of the composition $\comp$ in the canonical self-enrichment of $\M$ in \cref{selfenrm-diagrams}.
\item The right sub-region commutes by the definition of $\pn{\uni_{Ay}}$ as the partner of $\uni_{Ay}$ \cref{eval-bijection}.
\end{itemize}
This proves that the diagram \cref{unidgA-nat-top-right} is commutative.

\emph{Left-Bottom Composite}.
The partner of the left-bottom composite in \cref{unidgA-nat-diag} is the left-bottom composite in the following diagram in $\M$.
\begin{equation}\label{unidgA-nat-left-bot}
\begin{tikzpicture}[baseline={(b2.base)}]
\def\h{3.5} \def\g{3} \def\v{-1.4} \def\w{.6}
\draw[0cell=.8]
(0,0) node (a1) {\big( E\C(x,y) \scs \ang{} \scs Ax\big)}
(a1)+(\h,0) node (a2) {\big( E\C(x,y) \scs EFAx\big)}
(a1)+(0,\v) node (b1) {\scalebox{.9}{$\big(\clM(EFAx;EFAy) \scs \clM(Ax;EFAx) \scs Ax\big)$}}
(b1)+(\h,\v) node (b2) {\big(\clM(EFAx;EFAy) \scs EFAx\big)}
(b1)+(\h+\g,0) node (b3) {\big( EFE\C(x,y) \scs EFAx\big)}
(b1)+(0,2*\v) node (c1) {\big(\clM(Ax;EFAy) \scs Ax\big)}
(c1)+(\h+\g,0) node (c2) {EFAy}
(b1)+(.8,.6) node (bo) {\boxbox}
(b2)+(.8,.6) node (di) {\diamonddiamond}
(b3)+(-.4,-.7) node (star) {\pentagram}
;
\draw[1cell=.8]
(a1) edge node[swap] {\big((\Edg\Fdgr A)_{x,y} \scs (\unidg_A)_x \scs \opu\big)} (b1)
(b1) edge node[swap] {(\comp \scs \opu)} (c1)
(c1) edge node[pos=.4] {\ev^\M} (c2)
(b1) edge node[pos=.5] {(\opu \scs \ev^\M)} (b2)
(b2) edge node[pos=.5] {\ev^\M} (c2)
(a1) edge node {(\opu \scs \uni_{Ax})} (a2)
(a2) edge node[swap,pos=.2] {\big((\Edg\Fdgr A)_{x,y} \scs \opu\big)} (b2)
(a2) edge node[sloped] {\pn{(\Edg\Fdgr A)_{x,y}}} (c2)
(a2) edge node[pos=.3] {(E\cou_{\C(x,y)} \scs \opu)} (b3)
(b3) edge node {EF(\pn{A_{x,y}})} (c2)
;
\draw[1cell=.8]
(a1) [rounded corners=3pt] |- ($(a2)+(-1,\w)$)
-- node {(E\cou_{\C(x,y)} \scs \uni_{Ax})} ($(a2)+(1,\w)$) -| (b3)
;
\end{tikzpicture}
\end{equation}
The five sub-regions in \cref{unidgA-nat-left-bot} are commutative for the following reasons.
\begin{itemize}
\item The top sub-region commutes by the unity properties, \cref{enr-multicategory-right-unity,enr-multicategory-left-unity}, in $\M$.
\item The bottom triangle commutes by the definition of $\comp$ in \cref{selfenrm-diagrams}.
\item The sub-region labeled $\boxbox$ commutes by the definition of $(\unidg_A)_x$ as the partner of $\uni_{Ax}$ \cref{unidgAx}.
\item The triangle labeled $\diamonddiamond$ commutes by the definition of $\pn{(\Edg\Fdgr A)_{x,y}}$ as the partner of $(\Edg\Fdgr A)_{x,y}$ \cref{eval-bijection}.
\item The triangle labeled $\pentagram$ is the boundary of the following diagram. 
\[\begin{tikzpicture}
\def\h{6} \def\v{-1.4} \def\t{15}
\draw[0cell=.85]
(0,0) node (a1) {\big(E\C(x,y) \scs EFAx\big)}
(a1)+(\h,0) node (a2) {EFAy}
(a1)+(\h/2,\v) node (b) {\big(\clM(EFAx;EFAy) \scs EFAx\big)}
(b)+(0,\v) node (c) {\big(EFE\C(x,y) \scs EFAx\big)}
;
\draw[1cell=.85]
(a1) edge node {\pn{(\Edg\Fdgr A)_{x,y}}} (a2)
(a1) edge node {((\Edg\Fdgr A)_{x,y} \scs \opu)} (b)
(b) edge node {\ev^\M} (a2)
(a1) edge[out=-90,in=180] node[swap] {\big(E\cou_{\C(x,y)} \scs \opu\big)} (c)
(c) edge[out=0,in=-90] node[swap] {EF(\pn{A_{x,y}})} (a2)
(c) edge node[swap] {\big(\pn{(EF(\pn{A_{x,y}}))} \scs \opu \big)} (b)
;
\end{tikzpicture}\]
This diagram is commutative for the following reasons.
\begin{itemize}
\item The top triangle commutes by the definition of $\pn{(\Edg\Fdgr A)_{x,y}}$ as the partner of $(\Edg\Fdgr A)_{x,y}$ \cref{eval-bijection}.
\item The right sub-region commutes by the definition of $\pn{(EF(\pn{A_{x,y}}))}$ as the partner of $EF(\pn{A_{x,y}})$ \cref{eval-bijection}.
\item The left sub-region commutes by the top commutative triangle in \cref{EdgFdgrA-component}, which is proved in \cref{EdgFdgrA-commutes}.
\end{itemize}
\end{itemize}
This proves that the diagram \cref{unidgA-nat-left-bot} is commutative.

\emph{Comparing Partners}.
By the commutative diagrams \cref{unidgA-nat-top-right,unidgA-nat-left-bot}, the partners of the top-right composite and of the left-bottom composite in the desired diagram \cref{unidgA-nat-diag} are the following binary multimorphisms in $\M$.
\begin{equation}\label{unidgA-natural-partners}
\begin{aligned}
\ga^\M\scmap{EF(\pn{A_{x,y}}); \uni_{E\C(x,y)} \scs \uni_{Ax}}\\
\ga^\M\scmap{EF(\pn{A_{x,y}}); E\cou_{\C(x,y)} \scs \uni_{Ax}}
\end{aligned}
\end{equation}
The two binary multimorphisms in \cref{unidgA-natural-partners} are equal by \cref{def:mackey-gen-context} \cref{mackey-context-v}, which assumes the equality
\[\uni_{E\C(x,y)} = E\cou_{\C(x,y)} \cn E\C(x,y) \to EFE\C(x,y)\]
for all objects $x,y \in \C$.  Since taking partner is a bijection \cref{eval-bijection}, we conclude that the diagram \cref{unidgA-nat-diag} is commutative.
\end{proof}

\begin{lemma}\label{unidg-natural}
In the context of \cref{def:mackey-gen-context,def:unidg}, 
\[
\]
in \cref{unidg-def} is a natural transformation.
\end{lemma}

\begin{proof}
\cref{unidgA-natural} proves that each component of $\unidg$ is a well-defined $\M$-natural transformation.  Naturality for $\unidg$ means that, for each $\M$-natural transformation $\psi$ as in
\[%
\]
the following diagram of $\M$-natural transformations in $\MCat(\C_E,\M)$ commutes.
\begin{equation}\label{unidg-natural-psi}
%
\end{equation}
Since each $\M$-natural transformation is determined by its components, it suffices to show that, for each object $x \in \C$, the two vertical composites in \cref{unidg-natural-psi} have the same $x$-components.  By \cref{def:mnaturaltr-vcomp} these two $x$-components are the following two composite nullary multimorphisms in $\M$.
\begin{equation}\label{unidg-natural-psix}
%
\end{equation}
Since taking partner is a bijection \cref{eval-bijection}, it suffices to show that the two composites in \cref{unidg-natural-psix} have the same partner.  We compute these partners in \cref{unidg-nat-topright,unidg-nat-leftbot}.  Then we observe that they are equal to finish the proof.

\emph{Top-Right Composite}. 
The partner of the top-right composite in \cref{unidg-natural-psix} is the left-bottom composite unary multimorphism in the following diagram in $\M$.
\begin{equation}\label{unidg-nat-topright}
%
\end{equation}
The three sub-regions in \cref{unidg-nat-topright} are commutative for the following reasons.
\begin{itemize}
\item The bottom left sub-region commutes by the definition of $\comp$ in \cref{selfenrm-diagrams}.
\item The top left sub-region commutes by 
\begin{itemize}
\item the definition of $(\unidg_A)_x$ as the partner of $\uni_{Ax}$ in \cref{unidgAx} and
\item \cref{EdgFdgr-morphism-partner}, which implies that $(\Edg\Fdgr\psi)_x$ is the partner of $EF(\pn{\psi_x})$.
\end{itemize}
\item The right sub-region commutes by the definition of $\pn{(EF(\pn{\psi_x}))}$ as the partner of $EF(\pn{\psi_x})$.
\end{itemize}
This proves that the diagram \cref{unidg-nat-topright} is commutative.

\emph{Left-Bottom Composite}.
The partner of the left-bottom composite in \cref{unidg-natural-psix} is the left-bottom composite unary multimorphism in the following diagram in $\M$.
\begin{equation}\label{unidg-nat-leftbot}
\begin{tikzpicture}[vcenter]
\def\h{4} \def\v{-1.5} \def\u{-1.2} \def\t{2} \def\w{1}
\draw[0cell=.8]
(0,0) node (a1) {\big(\ang{} \scs \ang{} \scs Ax\big) = Ax}
(a1)+(\h,0) node (a2) {\big(\ang{} \scs Bx\big)}
(a1)+(0,\v) node (b1) {\big(\clM\scmap{Bx;EFBx} \scs \clM\scmap{Ax;Bx} \scs Ax \big)}
(b1)+(\h,\u) node (b2) {\big(\clM\scmap{Bx;EFBx} \scs Bx \big)}
(b1)+(0,2*\u) node (c1) {\big(\clM\scmap{Ax;EFBx} \scs Ax\big)}
(c1)+(\h,0) node (c2) {EFBx}
;
\draw[1cell=.8]
(a1) edge node {\pn{\psi_x}} (a2)
(b1) edge node[pos=.7] {(\opu \scs \ev^\M)} (b2)
(c1) edge node {\ev^\M} (c2)
(a1) edge node[swap] {\big( (\unidg_B)_x \scs \psi_x \scs \opu \big)} (b1)
(b1) edge node[swap] {(\comp \scs \opu)} (c1)
(a2) edge node[swap,pos=.2] {\big( \pn{\uni_{Bx}} \scs \opu\big)} (b2)
(b2) edge node[swap] {\ev^\M} (c2)
;
\draw[1cell=.8]
(a2) [rounded corners=3pt] -| ($(b2)+(\t,\w)$)
-- node[swap,pos=0] {\uni_{Bx}} ($(b2)+(\t,-\w)$) |- (c2)
;
\end{tikzpicture}
\end{equation}
The three sub-regions in \cref{unidg-nat-leftbot} are commutative for the following reasons.
\begin{itemize}
\item The bottom left sub-region commutes by the definition of $\comp$ in \cref{selfenrm-diagrams}.
\item The top left sub-region commutes by 
\begin{itemize}
\item the definition of $(\unidg_B)_x$ as the partner of $\uni_{Bx}$ in \cref{unidgAx} and
\item the definition of $\pn{\psi_x}$ as the partner of $\psi_x$.
\end{itemize}
\item The right sub-region commutes by the definition of $\pn{\uni_{Bx}}$ as the partner of $\uni_{Bx}$.
\end{itemize}
This proves that the diagram \cref{unidg-nat-leftbot} is commutative.

\emph{Comparing Partners}.
By the commutative diagrams \cref{unidg-nat-topright,unidg-nat-leftbot}, the partners of the two composites in the desired diagram \cref{unidg-natural-psix} are the composites in $\M$ as follows.
\[\begin{tikzpicture}
\def\v{-1.3}
\draw[0cell=.9]
(0,0) node (a) {Ax}
(a)+(3,0) node (b) {EFAx}
(a)+(0,\v) node (c) {Bx}
(b)+(0,\v) node (d) {EFBx}
;
\draw[1cell=.9]
(a) edge node {\uni_{Ax}} (b)
(b) edge node {EF(\pn{\psi_x})} (d)
(a) edge node[swap] {\pn{\psi_x}} (c)
(c) edge node {\uni_{Bx}} (d)
;
\end{tikzpicture}\]
This diagram commutes by the naturality condition \cref{enr-multinat} of the multinatural transformation $\uni \cn 1_\M \to EF$.
\end{proof}

\section{Comparing \texorpdfstring{$\Fdgr\Edg$}{FE} and the Identity}
\label{sec:coudg}

Throughout this section we assume the same context as in \cref{def:mackey-gen-context} and consider the functors in \cref{Fdgr-Edg}:
\[\begin{tikzpicture}[baseline={(a.base)}]
\draw[0cell]
(0,0) node (a) {\MCat(\C_E,\M)}
(a)+(3.6,0) node (b) {\NCat(\C,\N).}
;
\draw[1cell]
(a) edge[transform canvas={yshift=.6ex}] node {\Fdgr} (b)
(b) edge[transform canvas={yshift=-.4ex}] node {\Edg} (a) 
;
\end{tikzpicture}\]
In this section we extend the multinatural transformation in \cref{def:mackey-gen-context} \cref{mackey-context-iv} 
\[\begin{tikzpicture}[baseline={(a.base)}]
\def\t{25} \def\s{27}
\draw[0cell]
(0,0) node (a) {\N}
(a)+(1.7,0) node (b) {\N}
;
\draw[1cell=.9]
(a) edge[bend left=\s] node {1_\N} (b)
(a) edge[bend right=\s] node[swap] {FE} (b)
;
\draw[2cell]
node[between=a and b at .44, rotate=-90, 2label={above,\cou}] {\Rightarrow}
;
\end{tikzpicture}\]
to a natural transformation $\coudg$ comparing $\Fdgr\Edg$ and the identity functor.  This section is organized as follows.
\begin{itemize}
\item $\coudg$ is in \cref{def:coudg}.
\item To show that $\coudg$ has the desired naturality properties, in \cref{expl:FdgrEdg-objects,FdgrEdgP-commutes} we describe the codomain $\Fdgr\Edg$ of $\coudg$ on objects.
\item We describe $\Fdgr\Edg$ on morphisms in \cref{expl:FdgrEdg-morphisms,FdgrEdg-morphism-partner}.
\item We show that $\coudg$ is a natural transformation in \cref{coudgP-natural,coudg-natural}.
\end{itemize}

\subsection*{The Natural Transformation $\coudg$}

\begin{definition}\label{def:coudg}
In the context of \cref{def:mackey-gen-context,Fdgr-Edg}, we define the data of a natural transformation
\begin{equation}\label{coudg-def}
\begin{tikzpicture}[baseline={(a.base)}]
\def\s{25}
\draw[0cell]
(0,0) node (a) {\phantom{\M}}
(a)+(2,0) node (b) {\phantom{\N}}
(a)+(-.75,0) node (a') {\NCat(\C,\N)}
(b)+(.8,0) node (b') {\NCat(\C,\N)}
;
\draw[1cell=.9]
(a) edge[bend left=\s] node {1} (b)
(a) edge[bend right=\s] node[swap] {\Fdgr\Edg} (b)
;
\draw[2cell]
node[between=a and b at .45, rotate=-90, 2label={above,\coudg}] {\Rightarrow}
;
\end{tikzpicture}
\end{equation}
as follows.  For an $\N$-functor $P \cn \C \to \N$, the $P$-component of $\coudg$ is the $\N$-natural transformation
\begin{equation}\label{coudg-P}
\begin{tikzpicture}[baseline={(a.base)}]
\def\s{25}
\draw[0cell]
(0,0) node (a) {\C}
(a)+(2.2,0) node (b) {\N}
;
\draw[1cell=.9]
(a) edge[bend left=\s] node {P} (b)
(a) edge[bend right=\s] node[swap] {\Fdgr\Edg P} (b)
;
\draw[2cell]
node[between=a and b at .42, rotate=-90, 2label={above,\coudg_P}] {\Rightarrow}
;
\end{tikzpicture}
\end{equation}
with, for each object $x \in \C$, $x$-component given by the nullary multimorphism
\begin{equation}\label{coudgPx}
(\coudg_P)_x = \pn{\cou_{Px}} \cn \ang{} \to \clN\scmap{Px;FEPx} \inspace \N.
\end{equation}
This is the partner \cref{eval-bij-zero} of the $(Px)$-component 
\[\cou_{Px} \cn Px \to FEPx\]
of the multinatural transformation $\cou \cn 1_\N \to FE$, which is a unary multimorphism in $\N$.  This finishes the definition of $\coudg$.  We check that
\begin{itemize}
\item $\coudg_P$ is an $\N$-natural transformation in \cref{coudgP-natural} and 
\item $\coudg$ is a natural transformation in \cref{coudg-natural}.\defmark
\end{itemize} 
\end{definition}

Before we prove the $\N$-naturality of $\coudg_P$ and the naturality of $\coudg$, we first discuss the codomain $\Fdgr\Edg$ of $\coudg$ in detail in \cref{expl:FdgrEdg-objects,expl:FdgrEdg-morphisms}.

\begin{explanation}[$\Fdgr\Edg$ on Objects]\label{expl:FdgrEdg-objects}
The codomain of $\coudg$ in \cref{coudg-def} is the composite of the functors
\[\begin{split}
\Edg &= \Ese \circ \dE\andspace\\
\Fdgr &= \Fse \circ \dF \circ \C_\cou
\end{split}\]
as in \cref{expl:unidg}.  For an $\N$-functor $P \cn \C \to \N$ (\cref{def:mfunctor}), the codomain of $\coudg_P$ in \cref{coudg-P} is the following composite $\N$-functor.
\begin{equation}\label{FdgrEdg-P-diag}
\begin{tikzpicture}[vcenter]
\def\h{3} \def\v{1.3}
\draw[0cell=1]
(0,0) node (a) {\C}
(a)+(0,-\v) node (b) {\C_{FE}}
(b)+(\h,0) node (c) {\N_{FE} = (\N_E)_F}
(c)+(\h,0) node (d) {\M_F}
(d)+(0,\v) node (e) {\N}
;
\draw[1cell=.9]
(a) edge node[swap] {\C_\cou} (b)
(b) edge node {P_{FE}} (c)
(c) edge node {\Ese_F} (d)
(d) edge node[swap] {\Fse} (e)
(a) edge node {\Fdgr\Edg P} (e)
;
\end{tikzpicture}
\end{equation}

\emph{Object Assignment}. 
Since $\C_\cou$ \cref{C-sub-cou} has the identity object assignment, by \cref{mfunctor-change-enrichment,def:gspectra-thm-iii-context} the object assignment of $\Fdgr\Edg P$ in \cref{FdgrEdg-P-diag} is given by
\begin{equation}\label{FdgrEdgP-object}
(\Fdgr\Edg P)x = FEPx \forspace x \in \C.
\end{equation}
This implies that the $x$-component $(\coudg_P)_x$ in \cref{coudgPx} is well defined.

\emph{Components}. 
By \cref{FdgrEdg-P-diag}, for objects $x,y \in \C$, the $(x,y)$-component of $\Fdgr\Edg P$ is the composite unary multimorphism along the boundary of the following diagram in $\N$.
\begin{equation}\label{FdgrEdgP-component}
\begin{tikzpicture}[vcenter]
\def\g{-.5}\def\h{3} \def\v{1.5} \def\u{1}
\draw[0cell=.9]
(0,0) node (a) {\C(x,y)}
(a)+(\g,-\v) node (b) {FE\C(x,y)}
(b)+(\h,-\u) node (c) {FE\clN\scmap{Px;Py}}
(c)+(\h,\u) node (d) {F\clM\scmap{EPx;EPy}}
(d)+(\g,\v) node (e) {\clN\scmap{FEPx;FEPy}}
;
\draw[1cell=.9]
(a) edge node {(\Fdgr\Edg P)_{x,y}} (e)
(a) edge[transform canvas={xshift=-1em}] node[swap,pos=.8] {\cou_{\C(x,y)}} (b)
(b) edge[transform canvas={xshift=-2em}, shorten >=-1.7em] node[swap,pos=.6] {FEP_{x,y}} (c)
(c) edge[transform canvas={xshift=2em}, shorten <=-1.7em] node[swap,pos=.4] {F\Ese_{Px,Py}} (d)
(d) edge[transform canvas={xshift=1em}] node[swap,pos=.1] {\Fse_{EPx,EPy}} (e)
;
\draw[1cell=.8]
(b) edge node[pos=.3] {\pn{\big(FE(\pn{P_{x,y}})\big)}} (e)
(c) edge node[pos=.15] {\pn{\big(FE(\ev^\N_{Px;\, Py})\big)}\!\!} (e)
;
\end{tikzpicture}
\end{equation}
In \cref{FdgrEdgP-commutes} we discuss the two interior arrows in \cref{FdgrEdgP-component}.
\end{explanation}

\begin{lemma}\label{FdgrEdgP-commutes}
The diagram \cref{FdgrEdgP-component} is commutative.
\end{lemma}

\begin{proof}
Since the boundary of the diagram \cref{FdgrEdgP-component} commutes by \cref{FdgrEdg-P-diag}, it suffices to prove the following two equalities regarding its two interior arrows.
\begin{equation}\label{FEevN-pn}
\begin{aligned}
\ga^\N\scmap{\Fse_{EPx,EPy}; F\Ese_{Px,Py}} 
&= \pn{\big(FE(\ev^\N_{Px;\, Py})\big)}\\
\ga^\N\left(\pn{\big(FE(\ev^\N_{Px;\, Py})\big)} \sscs FEP_{x,y} \right) 
&= \pn{\big(FE(\pn{P_{x,y}})\big)}
\end{aligned}
\end{equation}
To prove these equalities, we consider the following diagram in $\N$.
\begin{equation}\label{FEevN-pn-diag}
\begin{tikzpicture}[vcenter]
\def\v{-1.3} \def\t{15}
\draw[0cell=.9]
(0,0) node (a) {\big(FE\C(x,y) \scs FEPx\big)}
(a)+(0,\v) node (b) {\big(FE\clN\scmap{Px;Py} \scs FEPx\big)}
(b)+(0,\v) node (c) {\big(F\clM\scmap{EPx;EPy} \scs FEPx\big)}
(c)+(0,\v) node (d) {\big(\clN\scmap{FEPx;FEPy} \scs FEPx\big)}
(a)+(5,0) node (e1) {FEPy}
(e1)+(0,\v) node (e2) {FEPy}
(e2)+(0,\v) node (e3) {FEPy}
(e3)+(0,\v) node (e4) {FEPy}
;
\draw[1cell=.9]
(a) edge node[swap] {(FEP_{x,y} \scs 1)} (b)
(b) edge node[swap] {(F\Ese_{Px,Py} \scs 1)} (c)
(c) edge node[swap] {(\Fse_{EPx,EPy} \scs 1)} (d)
(a) edge node {FE(\pn{P_{x,y}})} (e1)
(b) edge node {FE(\ev^\N_{Px;\, Py})} (e2)
(c) edge node {F(\ev^\M_{EPx;\, EPy})} (e3)
(d) edge node {\ev^\N_{FEPx;\, FEPy}} (e4)
(e1) edge[-,double equal sign distance] (e2)
(e2) edge[-,double equal sign distance] (e3)
(e3) edge[-,double equal sign distance] (e4)
;
\end{tikzpicture}
\end{equation}
The three sub-regions in \cref{FEevN-pn-diag} are commutative for the following reasons.
\begin{itemize}
\item In the top sub-region, $\pn{P_{x,y}}$ is the partner \cref{eval-bijection} of $P_{x,y}$.  By definition the following diagram in $\N$ commutes.  
\begin{equation}\label{Pxy-partner-diagram}
\begin{tikzpicture}[vcenter]
\draw[0cell=.9]
(0,0) node (a) {\big(\C(x,y) \scs Px\big)}
(a)+(0,-1.3) node (b) {\big(\clN\scmap{Px;Py} \scs Px\big)}
(b)+(4,0) node (c) {Py}
;
\draw[1cell=.9]
(a) edge node[swap] {(P_{x,y} \scs \opu)} (b)
(b) edge node[pos=.4] {\ev^\N_{Px;\, Py}} (c)
(a) edge[bend left=15] node {\pn{P_{x,y}}} (c)
;
\end{tikzpicture}
\end{equation}
Applying the non-symmetric multifunctor $FE$ to this commutative diagram yields the top sub-region in \cref{FEevN-pn-diag}.
\item The middle sub-region in \cref{FEevN-pn-diag} is obtained from the commutative diagram \cref{Fse-xy-diag} defining $\Ese_{Px,Py}$ by applying the non-symmetric multifunctor $F$.
\item The bottom sub-region in \cref{FEevN-pn-diag} is the commutative diagram \cref{Fse-xy-diag} that defines $\Fse_{EPx,EPy}$.
\end{itemize}

Since taking partner is a bijection \cref{eval-bijection}, the commutativity of the bottom two sub-regions in \cref{FEevN-pn-diag} proves the first desired equality in \cref{FEevN-pn}.  The second equality in \cref{FEevN-pn} follows from the first equality and the boundary of the commutative diagram \cref{FEevN-pn-diag}.
\end{proof}

\begin{explanation}[$\Fdgr\Edg$ on Morphisms]\label{expl:FdgrEdg-morphisms}
Suppose $\theta \cn P \to Q$ is an $\N$-natural transformation (\cref{def:mnaturaltr}) as in the left diagram below.
\[
\]
By \cref{expl:FdgrEdg-objects} the $\N$-natural transformation $\Fdgr\Edg\theta$, as in the right diagram above, is the following whiskering.
\begin{equation}\label{FdgrEdg-theta-diag}
%
\end{equation}
Since $\C_\cou$ has the identity object assignment, by \cref{mnat-change-enrichment,FdgrEdg-theta-diag}, for each object $x \in \C$ the $x$-component of $\Fdgr\Edg\theta$ is the following composite nullary multimorphism in $\N$.
\begin{equation}\label{FdgrEdg-theta-component}
%
\end{equation}
In \cref{FdgrEdg-theta-component},
\begin{itemize}
\item $\theta_x \cn \ang{} \to \clN(Px;Qx)$ is the $x$-component of $\theta$ \cref{mnatural-component}, and
\item the component morphisms $\Ese_{Px,Qx}$ and $\Fse_{EPx,EQx}$ are as in \cref{Fprimexy}.
\end{itemize}
Taking the partner \cref{eval-bijection} of $(\Fdgr\Edg\theta)_x$ yields the top unary multimorphism in $\N$ below.
\begin{equation}\label{FdgrEdg-theta-partner}
%
\end{equation}
The bottom arrow is the image under $FE$ of the unary multimorphism
\[\pn{\theta_x} \cn Px \to Qx \inspace \N,\]
which is the partner of $\theta_x$.  \cref{FdgrEdg-morphism-partner} proves that they are the same.
\end{explanation}

\begin{lemma}\label{FdgrEdg-morphism-partner}
The two unary multimorphisms in \cref{FdgrEdg-theta-partner} are equal.
\end{lemma}

\begin{proof}
Since $(\Fdgr\Edg\theta)_x$ is the composite in \cref{FdgrEdg-theta-component}, its partner is, by definition \cref{eval-bijection}, the left-bottom composite in the following diagram in $\N$.
\[\begin{tikzpicture}
\def\h{5} \def\v{-1.3}
\draw[0cell=.9]
(0,0) node (a1) {\big(\ang{} \scs FEPx\big)}
(a1)+(0,\v) node (b1) {\big(FE\clN\scmap{Px;Qx} \scs FEPx\big)}
(b1)+(0,\v) node (c1) {\big(F\clM\scmap{EPx;EQx} \scs FEPx\big)}
(c1)+(0,\v) node (d1) {\big(\clN\scmap{FEPx;FEQx} \scs FEPx\big)}
(a1)+(\h,0) node (a2) {FEQx}
(a2)+(0,\v) node (b2) {FEQx}
(b2)+(0,\v) node (c2) {FEQx}
(c2)+(0,\v) node (d2) {FEQx}
;
\draw[1cell=.9]
(a1) edge node {FE(\pn{\theta_x})} (a2)
(b1) edge node {FE(\ev^\N_{Px;\, Qx})} (b2)
(c1) edge node {F(\ev^\M_{EPx;\, EQx})} (c2)
(d1) edge node {\ev^\N_{FEPx;\, FEQx}} (d2)
(a1) edge node[swap] {(FE(\theta_x) \scs \opu)} (b1)
(b1) edge node[swap] {(F\Ese_{Px,Qx} \scs \opu)} (c1)
(c1) edge node[swap] {(\Fse_{EPx,EQx} \scs \opu)} (d1)
(a2) edge[-,double equal sign distance] (b2)
(b2) edge[-,double equal sign distance] (c2)
(c2) edge[-,double equal sign distance] (d2)
;
\end{tikzpicture}\]
The three sub-regions in the above diagram are commutative for the following reasons.
\begin{itemize}
\item The top sub-region is $FE$ applied to the commutative diagram in $\N$
\[\begin{tikzpicture}[vcenter]
\draw[0cell=.9]
(0,0) node (a) {\big(\ang{} \scs Px\big)}
(a)+(0,-1.3) node (b) {\big(\clN\scmap{Px;Qx} \scs Px\big)}
(b)+(4,0) node (c) {Qx}
;
\draw[1cell=.9]
(a) edge node[swap] {(\theta_x \scs \opu)} (b)
(b) edge node[pos=.4] {\ev^\N_{Px;\, Qx}} (c)
(a) edge[bend left=15] node {\pn{\theta_x}} (c)
;
\end{tikzpicture}\]
that defines the partner $\pn{\theta_x}$.
\item The middle sub-region is $F$ applied to the commutative diagram \cref{Fse-xy-diag} that defines $\Ese_{Px,Qx}$.
\item The bottom sub-region is the commutative diagram \cref{Fse-xy-diag} that defines $\Fse_{EPx,EQx}$.
\end{itemize}
This proves that $\pn{(\Fdgr\Edg\theta)_x}$ is equal to $FE(\pn{\theta_x})$.
\end{proof}

\subsection*{Naturality of $\coudg$}

To prove that $\coudg$ is a natural transformation, we first check that its components are well defined.

\begin{lemma}\label{coudgP-natural}
In the context of \cref{def:mackey-gen-context,def:coudg}, for each $\N$-functor $P \cn \C \to \N$, 
\[
\]
in \cref{coudg-P} is an $\N$-natural transformation.
\end{lemma}

\begin{proof}
We must prove that the naturality diagram \cref{mnaturality-diag} for $\coudg_P$, which is the diagram in $\N$ below for objects $x,y \in \C$, is commutative.
\begin{equation}\label{coudgP-nat-diag}
%
\end{equation}
Since taking partner in $\N$ is a bijection \cref{eval-bijection}, it suffices to show that the two composites in \cref{coudgP-nat-diag} have the same partner.  We compute these partners in \cref{coudgP-nat-top-right,coudgP-nat-left-bot}.  Then we observe that they are equal to finish the proof.

\emph{Top-Right Composite}. 
The partner of the top-right composite in \cref{coudgP-nat-diag} is, by definition, the left-bottom composite in the following diagram in $\N$.
\begin{equation}\label{coudgP-nat-top-right}
%
\end{equation}
The four sub-regions in \cref{coudgP-nat-top-right} are commutative for the following reasons.
\begin{itemize}
\item The top trapezoid commutes by the naturality \cref{enr-multinat} of the multinatural transformation $\cou \cn 1_\N \to FE$.
\item The top left triangle commutes by the definition of $\pn{P_{x,y}}$ in \cref{Pxy-partner-diagram} and the definition of the $y$-component of $\coudg_P$ in \cref{coudgPx},
\[(\coudg_P)_y = \pn{\cou_{Py}} \cn \ang{} \to \clN\scmap{Py; FEPy},\]
as the partner of the $(Py)$-component of $\cou$. 
\item The bottom rectangle commutes by the definition of the composition $\comp$ in the canonical self-enrichment of $\N$ in \cref{selfenrm-diagrams}.
\item The right sub-region commutes by the definition of $\pn{\cou_{Py}}$ as the partner of $\cou_{Py}$ \cref{eval-bijection}.
\end{itemize}
This proves that the diagram \cref{coudgP-nat-top-right} is commutative.

\emph{Left-Bottom Composite}.
The partner of the left-bottom composite in \cref{coudgP-nat-diag} is the left-bottom composite in the following diagram in $\N$.
\begin{equation}\label{coudgP-nat-left-bot}
\begin{tikzpicture}[baseline={(b2.base)}]
\def\h{3.5} \def\g{3} \def\v{-1.4} \def\w{.6}
\draw[0cell=.8]
(0,0) node (a1) {\big( \C(x,y) \scs \ang{} \scs Px\big)}
(a1)+(\h,0) node (a2) {\big( \C(x,y) \scs FEPx\big)}
(a1)+(0,\v) node (b1) {\scalebox{.9}{$\big(\clN(FEPx;FEPy) \scs \clN(Px;FEPx) \scs Px\big)$}}
(b1)+(\h,\v) node (b2) {\big(\clN(FEPx;FEPy) \scs FEPx\big)}
(b1)+(\h+\g,0) node (b3) {\big( FE\C(x,y) \scs FEPx\big)}
(b1)+(0,2*\v) node (c1) {\big(\clN(Px;FEPy) \scs Px\big)}
(c1)+(\h+\g,0) node (c2) {FEPy}
(b1)+(.8,.6) node (bo) {\boxbox}
(b2)+(.8,.6) node (di) {\diamonddiamond}
(b3)+(-.4,-.7) node (star) {\pentagram}
;
\draw[1cell=.8]
(a1) edge node[swap] {\big((\Fdgr\Edg P)_{x,y} \scs (\coudg_P)_x \scs \opu\big)} (b1)
(b1) edge node[swap] {(\comp \scs \opu)} (c1)
(c1) edge node[pos=.4] {\ev^\N} (c2)
(b1) edge node[pos=.5] {(\opu \scs \ev^\N)} (b2)
(b2) edge node[pos=.5] {\ev^\N} (c2)
(a1) edge node {(\opu \scs \cou_{Px})} (a2)
(a2) edge node[swap,pos=.2] {\big((\Fdgr\Edg P)_{x,y} \scs \opu\big)} (b2)
(a2) edge node[sloped] {\pn{(\Fdgr\Edg P)_{x,y}}} (c2)
(a2) edge node[pos=.3] {(\cou_{\C(x,y)} \scs \opu)} (b3)
(b3) edge node {FE(\pn{P_{x,y}})} (c2)
;
\draw[1cell=.8]
(a1) [rounded corners=3pt] |- ($(a2)+(-1,\w)$)
-- node {(\cou_{\C(x,y)} \scs \cou_{Px})} ($(a2)+(1,\w)$) -| (b3)
;
\end{tikzpicture}
\end{equation}
The five sub-regions in \cref{coudgP-nat-left-bot} are commutative for the following reasons.
\begin{itemize}
\item The top sub-region commutes by the unity properties, \cref{enr-multicategory-right-unity,enr-multicategory-left-unity}, in $\M$.
\item The bottom triangle commutes by the definition of $\comp$ in \cref{selfenrm-diagrams}.
\item The sub-region labeled $\boxbox$ commutes by the definition of $(\coudg_P)_x$ as the partner of $\cou_{Px}$ \cref{coudgPx}.
\item The triangle labeled $\diamonddiamond$ commutes by the definition of $\pn{(\Fdgr\Edg P)_{x,y}}$ as the partner of $(\Fdgr\Edg P)_{x,y}$ \cref{eval-bijection}.
\item The triangle labeled $\pentagram$ is the boundary of the following diagram. 
\[\begin{tikzpicture}
\def\h{6} \def\v{-1.4} \def\t{15}
\draw[0cell=.85]
(0,0) node (a1) {\big(\C(x,y) \scs FEPx\big)}
(a1)+(\h,0) node (a2) {FEPy}
(a1)+(\h/2,\v) node (b) {\big(\clN(FEPx;FEPy) \scs FEPx\big)}
(b)+(0,\v) node (c) {\big(FE\C(x,y) \scs FEPx\big)}
;
\draw[1cell=.85]
(a1) edge node {\pn{(\Fdgr\Edg P)_{x,y}}} (a2)
(a1) edge node {((\Fdgr\Edg P)_{x,y} \scs \opu)} (b)
(b) edge node {\ev^\N} (a2)
(a1) edge[out=-90,in=180] node[swap] {\big(\cou_{\C(x,y)} \scs \opu\big)} (c)
(c) edge[out=0,in=-90] node[swap] {FE(\pn{P_{x,y}})} (a2)
(c) edge node[swap] {\big(\pn{(FE(\pn{P_{x,y}}))} \scs \opu \big)} (b)
;
\end{tikzpicture}\]
This diagram is commutative for the following reasons.
\begin{itemize}
\item The top triangle commutes by the definition of $\pn{(\Fdgr\Edg P)_{x,y}}$ as the partner of $(\Fdgr\Edg P)_{x,y}$ \cref{eval-bijection}.
\item The right sub-region commutes by the definition of $\pn{(FE(\pn{P_{x,y}}))}$ as the partner of $FE(\pn{P_{x,y}})$ \cref{eval-bijection}.
\item The left sub-region commutes by the top commutative triangle in \cref{FdgrEdgP-component}, which is proved in \cref{FdgrEdgP-commutes}.
\end{itemize}
\end{itemize}
This proves that the diagram \cref{coudgP-nat-left-bot} is commutative.

\emph{Comparing Partners}.
By the commutative diagrams \cref{coudgP-nat-top-right,coudgP-nat-left-bot}, the partner of each composite in the desired diagram \cref{coudgP-nat-diag} is the following binary multimorphism in $\N$.
\[\ga^\N\scmap{FE(\pn{P_{x,y}}); \cou_{\C(x,y)} \scs \cou_{Px}}\]
Since taking partner is a bijection \cref{eval-bijection}, we conclude that the diagram \cref{coudgP-nat-diag} is commutative.
\end{proof}

\begin{remark}[Difference with $\unidg_A$]\label{rk:diff-unidgA}
We structure the proofs of \cref{unidgA-natural,coudgP-natural} in a way that highlights their conceptual similarity.  There is, however, one nontrivial difference between these two proofs of enriched naturality.  In the last paragraph of the proof of \cref{unidgA-natural}, to conclude that the two binary multimorphisms in \cref{unidgA-natural-partners} are the same, we need to use the assumption in \cref{def:mackey-gen-context} \cref{mackey-context-v}.  On the other hand, in the last paragraph of the proof of \cref{coudgP-natural}, no such assumption is needed.
\end{remark}

\begin{lemma}\label{coudg-natural}
In the context of \cref{def:mackey-gen-context,def:coudg}, 
\[
\]
in \cref{coudg-def} is a natural transformation.
\end{lemma}

\begin{proof}
\cref{coudgP-natural} proves that each component of $\coudg$ is a well-defined $\N$-natural transformation.  Naturality for $\coudg$ means that, for each $\N$-natural transformation $\theta$ as in
\[%
\]
the following diagram of $\N$-natural transformations in $\NCat(\C,\N)$ commutes.
\begin{equation}\label{coudg-natural-theta}
%
\end{equation}
Since each $\N$-natural transformation is determined by its components, it suffices to show that, for each object $x \in \C$, the two vertical composites in \cref{coudg-natural-theta} have the same $x$-components.  By \cref{def:mnaturaltr-vcomp} these two $x$-components are the following two composite nullary multimorphisms in $\N$.
\begin{equation}\label{coudg-natural-thetax}
%
\end{equation}
Since taking partner is a bijection \cref{eval-bijection}, it suffices to show that the two composites in \cref{coudg-natural-thetax} have the same partner.  We compute these partners in \cref{coudg-nat-topright,coudg-nat-leftbot}.  Then we observe that they are equal to finish the proof.

\emph{Top-Right Composite}. 
The partner of the top-right composite in \cref{coudg-natural-thetax} is the left-bottom composite unary multimorphism in the following diagram in $\N$.
\begin{equation}\label{coudg-nat-topright}
%
\end{equation}
The three sub-regions in \cref{coudg-nat-topright} are commutative for the following reasons.
\begin{itemize}
\item The bottom left sub-region commutes by the definition of $\comp$ in \cref{selfenrm-diagrams}.
\item The top left sub-region commutes by 
\begin{itemize}
\item the definition of $(\coudg_P)_x$ as the partner of $\cou_{Px}$ in \cref{coudgPx} and
\item \cref{FdgrEdg-morphism-partner}, which implies that $(\Fdgr\Edg\theta)_x$ is the partner of $FE(\pn{\theta_x})$.
\end{itemize}
\item The right sub-region commutes by the definition of $\pn{(FE(\pn{\theta_x}))}$ as the partner of $FE(\pn{\theta_x})$.
\end{itemize}
This proves that the diagram \cref{coudg-nat-topright} is commutative.

\emph{Left-Bottom Composite}.
The partner of the left-bottom composite in \cref{coudg-natural-thetax} is the left-bottom composite unary multimorphism in the following diagram in $\N$.
\begin{equation}\label{coudg-nat-leftbot}
\begin{tikzpicture}[vcenter]
\def\h{4} \def\v{-1.5} \def\u{-1.2} \def\t{2} \def\w{1}
\draw[0cell=.8]
(0,0) node (a1) {\big(\ang{} \scs \ang{} \scs Px\big) = Px}
(a1)+(\h,0) node (a2) {\big(\ang{} \scs Qx\big)}
(a1)+(0,\v) node (b1) {\big(\clN\scmap{Qx;FEQx} \scs \clN\scmap{Px;Qx} \scs Px \big)}
(b1)+(\h,\u) node (b2) {\big(\clN\scmap{Qx;FEQx} \scs Qx \big)}
(b1)+(0,2*\u) node (c1) {\big(\clN\scmap{Px;FEQx} \scs Px\big)}
(c1)+(\h,0) node (c2) {FEQx}
;
\draw[1cell=.8]
(a1) edge node {\pn{\theta_x}} (a2)
(b1) edge node[pos=.7] {(\opu \scs \ev^\N)} (b2)
(c1) edge node {\ev^\N} (c2)
(a1) edge node[swap] {\big( (\coudg_Q)_x \scs \theta_x \scs \opu \big)} (b1)
(b1) edge node[swap] {(\comp \scs \opu)} (c1)
(a2) edge node[swap,pos=.2] {\big( \pn{\cou_{Qx}} \scs \opu\big)} (b2)
(b2) edge node[swap] {\ev^\N} (c2)
;
\draw[1cell=.8]
(a2) [rounded corners=3pt] -| ($(b2)+(\t,\w)$)
-- node[swap,pos=0] {\cou_{Qx}} ($(b2)+(\t,-\w)$) |- (c2)
;
\end{tikzpicture}
\end{equation}
The three sub-regions in \cref{coudg-nat-leftbot} are commutative for the following reasons.
\begin{itemize}
\item The bottom left sub-region commutes by the definition of $\comp$ in \cref{selfenrm-diagrams}.
\item The top left sub-region commutes by 
\begin{itemize}
\item the definition of $(\coudg_Q)_x$ as the partner of $\cou_{Qx}$ in \cref{coudgPx} and
\item the definition of $\pn{\theta_x}$ as the partner of $\theta_x$.
\end{itemize}
\item The right sub-region commutes by the definition of $\pn{\cou_{Qx}}$ as the partner of $\cou_{Qx}$.
\end{itemize}
This proves that the diagram \cref{coudg-nat-leftbot} is commutative.

\emph{Comparing Partners}.
By the commutative diagrams \cref{coudg-nat-topright,coudg-nat-leftbot}, the partners of the two composites in the desired diagram \cref{coudg-natural-thetax} are the composites in $\N$ as follows.
\[\begin{tikzpicture}
\def\v{-1.3}
\draw[0cell=.9]
(0,0) node (a) {Px}
(a)+(3,0) node (b) {FEPx}
(a)+(0,\v) node (c) {Qx}
(b)+(0,\v) node (d) {FEQx}
;
\draw[1cell=.9]
(a) edge node {\cou_{Px}} (b)
(b) edge node {FE(\pn{\theta_x})} (d)
(a) edge node[swap] {\pn{\theta_x}} (c)
(c) edge node {\cou_{Qx}} (d)
;
\end{tikzpicture}\]
This diagram commutes by the naturality condition \cref{enr-multinat} of the multinatural transformation $\cou \cn 1_\N \to FE$.
\end{proof}

\section[Homotopy Equivalent Mackey Functor Categories]{Homotopy Equivalent Enriched Diagram and Mackey Functor Categories}
\label{sec:heq-enrdiag-cat}

In this section we apply the results in previous sections to prove that, if the data $(F,E,\uni,\cou)$ in \cref{def:mackey-gen-context} are inverse equivalences of homotopy theories (\cref{def:inverse-heq}), then they lift to inverse equivalences of homotopy theories between the enriched diagram categories in \cref{Fdgr-Edg}
\[\begin{tikzpicture}[baseline={(a.base)}]
\draw[0cell]
(0,0) node (a) {\MCat(\C_E,\M)}
(a)+(3.6,0) node (b) {\NCat(\C,\N).}
;
\draw[1cell]
(a) edge[transform canvas={yshift=.6ex}] node {\Fdgr} (b)
(b) edge[transform canvas={yshift=-.4ex}] node {\Edg} (a) 
;
\end{tikzpicture}\]
This section is organized as follows.
\begin{itemize}
\item Stable equivalences in categories of enriched diagrams are defined componentwise.  We make this precise in \cref{def:enr-diagcat-relative,def:mackey-heq-context}.  \cref{cWtri,MFinvcX} contain some basic properties of these componentwise stable equivalences.
\item \cref{mackey-gen-xiv} is the main result of this chapter.  It states that, under suitable conditions, $\Fdgr$ and $\Edg$ are inverse equivalences of homotopy theories.
\item \cref{mackey-xiv-cor} is the variant that involves an opposite $\N$-category $\Cop$. 
\item \cref{FE-inv-heq} shows that, under the assumptions of \cref{mackey-gen-xiv}, $F$ and $E$ are inverse equivalences of homotopy theories between the underlying categories.  Therefore, the inverse equivalences of homotopy theories $\Fdgr$ and $\Edg$ are, in fact, lifted from the underlying categories.
\end{itemize}

\subsection*{Componentwise Relative Structure on Enriched Diagram Categories}

To consider the homotopy theory of an enriched diagram category, we first define its relative category structure (\cref{definition:rel-cat}).  Recall that a multicategory has an \emph{underlying category} (\cref{ex:unarycategory}), which we denote by the same symbol.  Recall that a subcategory is \emph{wide} if it contains all the objects of the larger category.  The next definition is an adaptation of \cref{def:nonsymalgebra} to the current setting of enriched diagram categories.

\begin{definition}\label{def:enr-diagcat-relative}
Suppose $\P$ is a non-symmetric closed multicategory, and $\D$ is a $\P$-category (\cref{def:menriched-cat,def:closed-multicat}).  Suppose the underlying category of $\P$ is equipped with the structure of a relative category $(\P,\cW)$.  For the category $\PCat(\D,\P)$ in \cref{mcat-cm}, we define a wide subcategory
\begin{equation}\label{cWtri-subcat}
\cWtri \bigsubset \PCat(\D,\P)
\end{equation}
as follows.  A $\P$-natural transformation $\theta$ (\cref{def:mnaturaltr}) as in
\[\begin{tikzpicture}
\def\h{2} \def\t{25} \def\s{23}
\draw[0cell]
(0,0) node (a) {\D}
(a)+(\h,0) node (b) {\P}
;
\draw[1cell=.8]
(a) edge[bend left=\t] node {A} (b)
(a) edge[bend right=\t] node[swap] {B} (b)
;
\draw[2cell]
node[between=a and b at .45, rotate=-90, 2label={above,\theta}] {\Rightarrow}
;
\end{tikzpicture}\]
is in $\cWtri$ if, for each object $x$ in $\D$, the unary multimorphism
\begin{equation}\label{cWtri-def}
\pn{\theta_x} \cn Ax \to Bx \quad\text{is in $\cW$}.
\end{equation}
Here $\pn{\theta_x}$ is the partner \cref{eval-bij-zero} of the $x$-component of $\theta$, 
\[\theta_x \cn \ang{} \to \clP(Ax;Bx),\]
which is a nullary multimorphism in $\P$.  We regard the pair
\[\big(\PCat(\D,\P) \scs \cWtri\big)\]
as a relative category.
\end{definition}
In other words, $\theta$ is in $\cWtri$ if each component of $\theta$ has its partner in $\cW$.

Recall that a \emph{category with weak equivalences} (\cref{definition:rel-cat} \eqref{relcat-vi}) is a relative category $(\P,\cW)$ such that $\cW$ contains all the isomorphisms in $\P$ and has the 2-out-of-3 property.

\begin{lemma}\label{cWtri}
In the context of \cref{def:enr-diagcat-relative}, statements \cref{cWtri-i,cWtri-ii,cWtri-iii,cWtri-iv} below hold.
\begin{romenumerate}
\item\label{cWtri-i} The subcategory $\cWtri$ in \cref{cWtri-subcat} is well defined.
\item\label{cWtri-ii} If $\cW$ has the 2-out-of-3 property, then so does $\cWtri$.
\item\label{cWtri-iii} If $\cW$ contains all the isomorphisms, then so does $\cWtri$.
\item\label{cWtri-iv} If $(\P,\cW)$ is a category with weak equivalences, then so is $\big(\PCat(\D,\P) \scs \cWtri\big)$.
\end{romenumerate}
\end{lemma}

\begin{proof}
Suppose $\theta$ and $\psi$ are $\P$-natural transformations as in the left diagram below.
\[\begin{tikzpicture}
\def\h{2.5} \def\t{40} \def\s{30}
\draw[0cell]
(0,0) node (a) {\D}
(a)+(\h,0) node (b) {\P}
;
\draw[1cell=.8]
(a) edge[bend left=\t] node {A} (b)
(a) edge node[pos=.2] {B} (b)
(a) edge[bend right=\t] node[swap] {C} (b)
;
\draw[2cell]
node[between=a and b at .45, shift={(0,.3)}, rotate=-90, 2label={above,\theta}] {\Rightarrow}
node[between=a and b at .45, shift={(0,-.3)}, rotate=-90, 2label={above,\psi}] {\Rightarrow}
;
\begin{scope}[shift={(4,0)}]
\draw[0cell]
(0,0) node (a) {\D}
(a)+(2,0) node (b) {\P}
;
\draw[1cell=.9]
(a) edge[bend left=\s] node {A} (b)
(a) edge[bend right=\s] node[swap] {C} (b)
;
\draw[2cell]
node[between=a and b at .4, rotate=-90, 2label={above,\psi\theta}] {\Rightarrow}
;
\end{scope}
\end{tikzpicture}\]
For each object $x$ in $\D$, the $x$-component of the vertical composite $\psi\theta$, as in the right diagram above, has partner \cref{eval-bijection} given by 
\begin{equation}\label{psithetax-pn-comp}
\pn{(\psi\theta)_x} = \ga^\P\scmap{\pn{\psi_x}; \pn{\theta_x}} \cn Ax \to Cx
\end{equation}
by \cref{mcat-vertical-comp}.  Thus if $\theta$ and $\psi$ are in $\cWtri$, then so is $\psi\theta$, proving statement \cref{cWtri-i}.  The equality \cref{psithetax-pn-comp} also shows that $\cWtri$ has the 2-out-of-3 property whenever $\cW$ does, proving statement \cref{cWtri-ii}.  

To prove statement \cref{cWtri-iii}, suppose $\theta$ and $\psi$ are inverses of each other.  By \cref{selfenrm-id} this means that, for each object $x$ in $\D$, there are equalities
\begin{equation}\label{psitheta-inverses}
\pn{(\psi\theta)_x} = 1_{Ax} \andspace \pn{(\theta\psi)_x} = 1_{Bx}.
\end{equation}
Together with \cref{psithetax-pn-comp} and the variant for $\pn{(\theta\psi)_x}$, the equalities in \cref{psitheta-inverses} imply that $\pn{\psi_x}$ and $\pn{\theta_x}$ are mutually inverse isomorphisms in the underlying category of $\P$.  This proves statement \cref{cWtri-iii}.

Statement \cref{cWtri-iv} follows from statements \cref{cWtri-i,cWtri-ii,cWtri-iii}.
\end{proof}

Next we apply \cref{def:enr-diagcat-relative} to the context of the previous sections.  We consider the underlying functor of the non-symmetric multifunctor $F \cn \M \to \N$ (\cref{ex:un-v-functor}), which we denote by the same symbol.

\begin{definition}\label{def:mackey-heq-context}
In the context of \cref{def:mackey-gen-context}, suppose, in addition, the underlying category of $\N$ is equipped with the structure of a category with weak equivalences $(\N,\cX)$ (\cref{definition:rel-cat} \eqref{relcat-vi}).  We define the following.
\begin{itemize}
\item We define the wide subcategory
\begin{equation}\label{FinvcX}
F^\inv\cX \bigsubset \M
\end{equation}
as the preimage of $\cX$ under the underlying functor of $F \cn \M \to \N$.  We refer to morphisms in $F^\inv\cX$ as \emph{$F$-stable equivalences}.  We regard the pair
\begin{equation}\label{MFinvX-relative}
\big(\M, F^\inv\cX\big)
\end{equation}
as a relative category.
\item Applying \cref{def:enr-diagcat-relative} with $(\D,\P,\cW)$ given by
\begin{itemize}
\item $\big(\C_E, \M, F^\inv\cX\big)$ and
\item $(\C,\N,\cX)$,
\end{itemize}
we obtain the following two relative categories.
\begin{equation}\label{diag-enr-rel-cat}
\big(\MCat(\C_E,\M) \scs \FinvcXtri\big) \qquad \big(\NCat(\C,\N) \scs \cXtri\big) 
\end{equation}
\end{itemize}  
This finishes the definition.
\end{definition}

\begin{explanation}[Unpacking $\FinvcXtri$]\label{expl:FinvXtri}
An $\M$-natural transformation
\[\begin{tikzpicture}
\def\h{2} \def\t{25} \def\s{23}
\draw[0cell]
(0,0) node (a) {\C_E}
(a)+(\h,0) node (b) {\M}
;
\draw[1cell=.8]
(a) edge[bend left=\t] node {A} (b)
(a) edge[bend right=\t] node[swap] {B} (b)
;
\draw[2cell]
node[between=a and b at .45, rotate=-90, 2label={above,\psi}] {\Rightarrow}
;
\end{tikzpicture}\]
is in the subcategory $\FinvcXtri$ if, for each object $x$ in $\C$, the unary multimorphism 
\[\pn{\psi_x} \cn Ax \to Bx \quad \text{is in $F^\inv \cX \bigsubset \M$}.\]  
This means that the unary multimorphism 
\begin{equation}\label{Fpsix-pn-X}
F(\pn{\psi_x}) \cn FAx \to FBx \quad \text{is in $\cX \bigsubset \N$}.\defmark
\end{equation}
\end{explanation}

\begin{lemma}\label{MFinvcX}
The relative categories in \cref{MFinvX-relative,diag-enr-rel-cat} are categories with weak equivalences.
\end{lemma}

\begin{proof}
The wide subcategory $F^\inv\cX \bigsubset \M$ in \cref{FinvcX} contains all the isomorphisms and has the 2-out-of-3 property for the following two reasons. 
\begin{itemize}
\item By assumption $\cX \bigsubset \N$ has these two properties.
\item The underlying functor of $F \cn \M \to \N$ preserves identity morphisms and composition of morphisms.
\end{itemize}
The relative categories in \cref{diag-enr-rel-cat} are categories with weak equivalences by \cref{cWtri} \cref{cWtri-iv} applied to $\big(\C_E, \M, F^\inv\cX\big)$ and $(\C,\N,\cX)$.
\end{proof}

\subsection*{Equivalences of Homotopy Theories}

Recall that \index{equivalence!inverse - of homotopy theories}\index{homotopy theory!inverse equivalence of}\index{inverse equivalence!of homotopy theories}\emph{inverse equivalences of homotopy theories} are two relative functors in opposite direction such that each composite is connected to the identity functor by a zigzag of relative natural transformations (\cref{definition:rel-cat,def:inverse-heq}).  Such functors are equivalences of homotopy theories by \cref{gjo29}.  \cref{mackey-gen-xiv} below is the main result of this chapter.  It provides checkable criteria that guarantee that inverse equivalences of homotopy theories (\cref{FE-inv-heq})
\[\begin{tikzpicture}
\draw[0cell]
(0,0) node (a) {\big(\M, F^\inv\cX)}
(a)+(3,0) node (b) {(\N,\cX)}
;
\draw[1cell]
(a) edge[transform canvas={yshift=.6ex}] node {F} node[swap] {\sim} (b)
(b) edge[transform canvas={yshift=-.6ex}] node {E} (a)
;
\end{tikzpicture}\]
lift to the categories of enriched diagrams in \cref{diag-enr-rel-cat}.  \cref{expl:mackey-xiv-assumptions} summarizes all the hypotheses in \cref{mackey-gen-xiv}.  The variant for the opposite $\N$-category $\Cop$ is \cref{mackey-xiv-cor}.  In \cref{ch:mackey_eq} we apply \cref{mackey-gen-xiv} to the (non-symmetric) multifunctors connecting $\permcatsu$, $\pMulticat$, and $\MoneMod$; see \cref{mackey-xiv-pmulticat,mackey-xiv-mone,mackey-pmulti-mone,qu:morita}.

\begin{theorem}\label{mackey-gen-xiv}
In the context of \cref{def:mackey-gen-context,def:mackey-heq-context}, suppose, in addition, the components of the multinatural transformations
\[\begin{tikzpicture}[baseline={(a.base)}]
\def\t{25} \def\s{27}
\draw[0cell]
(0,0) node (a) {\M}
(a)+(1.7,0) node (b) {\M}
(b)+(1.5,0) node (and) {\text{and}}
;
\draw[1cell=.9]
(a) edge[bend left=\t] node {1_\M} (b)
(a) edge[bend right=\t] node[swap] {EF} (b)
;
\draw[2cell]
node[between=a and b at .44, rotate=-90, 2label={above,\uni}] {\Rightarrow}
;
\begin{scope}[shift={(4.7,0)}]
\draw[0cell]
(0,0) node (a) {\N}
(a)+(1.7,0) node (b) {\N}
;
\draw[1cell=.9]
(a) edge[bend left=\s] node {1_\N} (b)
(a) edge[bend right=\s] node[swap] {FE} (b)
;
\draw[2cell]
node[between=a and b at .44, rotate=-90, 2label={above,\cou}] {\Rightarrow}
;
\end{scope}
\end{tikzpicture}\]
are in $F^\inv\cX$ and $\cX$, respectively.  Then the functors in \cref{Fdgr-Edg}
\begin{equation}\label{mackey-xiv-functors}
\begin{tikzpicture}[baseline={(a.base)}]
\draw[0cell]
(0,0) node (a) {\big(\MCat(\C_E,\M) \scs \FinvcXtri\big)}
(a)+(5,0) node (b) {\big(\NCat(\C,\N) \scs \cXtri\big)} 
;
\draw[1cell]
(a) edge[transform canvas={yshift=.6ex}] node {\Fdgr} node[swap] {\sim} (b)
(b) edge[transform canvas={yshift=-.6ex}] node {\Edg} (a)
;
\end{tikzpicture}
\end{equation}
are inverse equivalences of homotopy theories.
\end{theorem}

\begin{proof}
The two relative categories in \cref{mackey-xiv-functors} are defined in \cref{diag-enr-rel-cat}.  They are categories with weak equivalences by \cref{MFinvcX}.

We prove statements \cref{mackey-xiv-i,mackey-xiv-ii,mackey-xiv-iii} below.  By \cref{def:inverse-heq} these statements imply that $\Fdgr$ and $\Edg$ are inverse equivalences of homotopy theories.
\begin{romenumerate}
\item\label{mackey-xiv-i} $\Fdgr$ and $\Edg$ in \cref{mackey-xiv-functors} are relative functors.
\item\label{mackey-xiv-ii} The natural transformation in \cref{unidg-def}
\[\begin{tikzpicture}[baseline={(a.base)}]
\def\s{25}
\draw[0cell]
(0,0) node (a) {\phantom{\M}}
(a)+(2,0) node (b) {\phantom{\N}}
(a)+(-.9,0) node (a') {\MCat(\C_E,\M)}
(b)+(.9,0) node (b') {\MCat(\C_E,\M)}
;
\draw[1cell=.9]
(a) edge[bend left=\s] node {1} (b)
(a) edge[bend right=\s] node[swap] {\Edg\Fdgr} (b)
;
\draw[2cell]
node[between=a and b at .45, rotate=-90, 2label={above,\unidg}] {\Rightarrow}
;
\end{tikzpicture}\]
has each component in $\FinvcXtri$.
\item\label{mackey-xiv-iii} The natural transformation in \cref{coudg-def}
\[\begin{tikzpicture}[baseline={(a.base)}]
\def\s{25}
\draw[0cell]
(0,0) node (a) {\phantom{\M}}
(a)+(2,0) node (b) {\phantom{\N}}
(a)+(-.8,0) node (a') {\NCat(\C,\N)}
(b)+(.8,0) node (b') {\NCat(\C,\N)}
;
\draw[1cell=.9]
(a) edge[bend left=\s] node {1} (b)
(a) edge[bend right=\s] node[swap] {\Fdgr\Edg} (b)
;
\draw[2cell]
node[between=a and b at .45, rotate=-90, 2label={above,\coudg}] {\Rightarrow}
;
\end{tikzpicture}\]
has each component in $\cXtri$.
\end{romenumerate}

\medskip
\emph{Statement \cref{mackey-xiv-i}: $\Fdgr$ is a Relative Functor}. 
This means that, for each $\M$-natural transformation $\psi \in \FinvcXtri$ as in \cref{psi-Fdgrpsi}, the $\N$-natural transformation $\Fdgr\psi$ is in $\cXtri$.  
\begin{equation}\label{psi-Fdgrpsi}
\begin{tikzpicture}[baseline={(a.base)}]
\def\h{2} \def\t{25} \def\s{25}
\draw[0cell]
(0,0) node (a) {\C_E}
(a)+(\h,0) node (b) {\M}
;
\draw[1cell=.8]
(a) edge[bend left=\t] node {A} (b)
(a) edge[bend right=\t] node[swap] {B} (b)
;
\draw[2cell]
node[between=a and b at .45, rotate=-90, 2label={above,\psi}] {\Rightarrow}
;
\begin{scope}[shift={(4,0)}]
\draw[0cell]
(0,0) node (a) {\C}
(a)+(2.2,0) node (b) {\N}
;
\draw[1cell=.8]
(a) edge[bend left=\s] node {\Fdgr A} (b)
(a) edge[bend right=\s] node[swap] {\Fdgr B} (b)
;
\draw[2cell]
node[between=a and b at .37, rotate=-90, 2label={above,\Fdgr\psi}] {\Rightarrow}
;
\end{scope}
\end{tikzpicture}
\end{equation}
By \cref{Fdgrpsix-partner,cWtri-def}, the desired condition $\Fdgr\psi \in \cXtri$ means that, for each object $x$ in $\C$, the unary multimorphism
\begin{equation}\label{partner-Fdgrpsix}
\pn{(\Fdgr\psi)_x} = F(\pn{\psi_x}) \quad \text{is in $\cX$}.
\end{equation}
This is true by \cref{Fpsix-pn-X}.  Thus $\Fdgr$ is a relative functor.

\medskip
\emph{Statement \cref{mackey-xiv-i}: $\Edg$ is a Relative Functor}. 
This means that, for each $\N$-natural transformation $\theta \in \cXtri$ as in \cref{theta-Edgtheta}, the $\M$-natural transformation $\Edg\theta$ is in $\FinvcXtri$.  
\begin{equation}\label{theta-Edgtheta}
\begin{tikzpicture}[baseline={(a.base)}]
\def\h{2} \def\t{25} \def\s{23}
\draw[0cell]
(0,0) node (a) {\C}
(a)+(\h,0) node (b) {\N}
;
\draw[1cell=.8]
(a) edge[bend left=\t] node {P} (b)
(a) edge[bend right=\t] node[swap] {Q} (b)
;
\draw[2cell]
node[between=a and b at .45, rotate=-90, 2label={above,\theta}] {\Rightarrow}
;
\begin{scope}[shift={(4,0)}]
\draw[0cell]
(0,0) node (a) {\C_E}
(a)+(2.2,0) node (b) {\M}
;
\draw[1cell=.8]
(a) edge[bend left=\s] node {\Edg P} (b)
(a) edge[bend right=\s] node[swap] {\Edg Q} (b)
;
\draw[2cell]
node[between=a and b at .4, rotate=-90, 2label={above,\Edg\theta}] {\Rightarrow}
;
\end{scope}
\end{tikzpicture}
\end{equation}
By \cref{Fpsix-pn-X} and \cref{Fdg-psix-pn-equality} applied to $\Edg\theta$, the desired condition $\Edg\theta \in \FinvcXtri$ means that, for each object $x$ in $\C$, the unary multimorphism
\begin{equation}\label{FEthetaxpn-in-X}
F\big((\pn{\Edg\theta)_x}\big) = FE(\pn{\theta_x}) \cn FEPx \to FEQx \quad \text{is in $\cX$}.
\end{equation}
To prove \cref{FEthetaxpn-in-X}, we use the naturality of $\cou \cn 1_\N \to FE$ \cref{enr-multinat} to obtain the following commutative diagram in $\N$.
\begin{equation}\label{cou-natural-thetaxpn}
\begin{tikzpicture}[vcenter]
\def\v{-1.4}
\draw[0cell]
(0,0) node (a) {Px}
(a)+(3,0) node (b) {FEPx}
(a)+(0,\v) node (c) {Qx}
(b)+(0,\v) node (d) {FEQx}
;
\draw[1cell=.9]
(a) edge node {\cou_{Px}} (b)
(c) edge node {\cou_{Qx}} (d)
(a) edge node[swap] {\pn{\theta_x}} (c)
(b) edge node {FE(\pn{\theta_x})} (d)
;
\end{tikzpicture}
\end{equation}
\begin{itemize}
\item The assumption $\theta \in \cXtri$ means that each $\pn{\theta_x}$ is in $\cX$. 
\item The components $\cou_{Px}$ and $\cou_{Qx}$ are in $\cX$ by the assumption on $\cou$. 
\end{itemize}
The commutative diagram \cref{cou-natural-thetaxpn} and the 2-out-of-3 property of $\cX$ imply that $FE(\pn{\theta_x})$ is in $\cX$, proving the desired condition \cref{FEthetaxpn-in-X}.  Thus $\Edg$ is a relative functor.  This finishes the proof of statement \cref{mackey-xiv-i}.

\medskip
\emph{Statement \cref{mackey-xiv-ii}}. 
Suppose $A \cn \C_E \to \M$ is an $\M$-functor.  We want to show that the $\M$-natural transformation in \cref{unidg-A}
\[\begin{tikzpicture}[baseline={(a.base)}]
\def\s{25}
\draw[0cell]
(0,0) node (a) {\C_E}
(a)+(2.3,0) node (b) {\M}
;
\draw[1cell=.9]
(a) edge[bend left=\s] node {A} (b)
(a) edge[bend right=\s] node[swap] {\Edg\Fdgr A} (b)
;
\draw[2cell]
node[between=a and b at .42, rotate=-90, 2label={above,\unidg_A}] {\Rightarrow}
;
\end{tikzpicture}\]
is in $\FinvcXtri$.  By \cref{unidgAx,Fpsix-pn-X}, the desired condition $\unidg_A \in \FinvcXtri$ means that, for each object $x$ in $\C$, the unary multimorphism
\begin{equation}\label{FunidgAx-FuniAx}
F\big(\pn{(\unidg_A)_x}\big) = F(\uni_{Ax}) \cn FAx \to FEFAx \quad \text{is in $\cX$}.
\end{equation}
This is true because each component of $\uni \cn 1_\M \to EF$ is in $F^\inv\cX$ by assumption.  This proves statement \cref{mackey-xiv-ii}.

\medskip
\emph{Statement \cref{mackey-xiv-iii}}.
Suppose $P \cn \C \to \N$ is an $\N$-functor.  We want to show that the $\N$-natural transformation in \cref{coudg-P}
\[\begin{tikzpicture}[baseline={(a.base)}]
\def\s{25}
\draw[0cell]
(0,0) node (a) {\C}
(a)+(2.2,0) node (b) {\N}
;
\draw[1cell=.9]
(a) edge[bend left=\s] node {P} (b)
(a) edge[bend right=\s] node[swap] {\Fdgr\Edg P} (b)
;
\draw[2cell]
node[between=a and b at .42, rotate=-90, 2label={above,\coudg_P}] {\Rightarrow}
;
\end{tikzpicture}\]
is in $\cXtri$.  By \cref{coudgPx,cWtri-def}, the desired condition $\coudg_P \in \cXtri$ means that, for each object $x$ in $\C$, the unary multimorphism
\begin{equation}\label{coudgPxpn-couPx}
\pn{(\coudg_P)_x} = \cou_{Px} \cn Px \to FEPx \quad\text{is in $\cX$}.
\end{equation}
This is true by the assumption that each component of $\cou$ is in $\cX$.  This proves statement \cref{mackey-xiv-iii}.
\end{proof}

\begin{explanation}[Summary]\label{expl:mackey-xiv-assumptions}
We summarize the assumptions for \cref{mackey-gen-xiv} in \eqref{mackey-xiv-assum-i} through \eqref{mackey-xiv-assum-iii} below.
\begin{enumerate}
\item\label{mackey-xiv-assum-i} We assume \cref{mackey-context-i,mackey-context-ii,mackey-context-iii,mackey-context-iv,mackey-context-v} in \cref{def:mackey-gen-context}.  These assumptions are multicategorical in nature.  They do not involve relative category structures.  Using these assumptions we construct
\begin{itemize}
\item the functors 
\[\begin{tikzpicture}[baseline={(a.base)}]
\draw[0cell]
(0,0) node (a) {\MCat(\C_E,\M)}
(a)+(3.5,0) node (b) {\NCat(\C,\N)} 
;
\draw[1cell]
(a) edge[transform canvas={yshift=.6ex}] node {\Fdgr} (b)
(b) edge[transform canvas={yshift=-.4ex}] node {\Edg} (a)
;
\end{tikzpicture}\]
in \cref{gspectra-thm-v} (applied to $\Edg$) and \cref{Fdgr-def} and
\item the natural transformations 
\[\begin{tikzpicture}[baseline={(a.base)}]
\def\s{25}
\draw[0cell]
(0,0) node (a) {\phantom{\M}}
(a)+(2,0) node (b) {\phantom{\N}}
(a)+(-.9,0) node (a') {\MCat(\C_E,\M)}
(b)+(.9,0) node (b') {\MCat(\C_E,\M)}
;
\draw[1cell=.9]
(a) edge[bend left=\s] node {1} (b)
(a) edge[bend right=\s] node[swap] {\Edg\Fdgr} (b)
;
\draw[2cell]
node[between=a and b at .45, rotate=-90, 2label={above,\unidg}] {\Rightarrow}
;
\end{tikzpicture}\]
\[\begin{tikzpicture}[baseline={(a.base)}]
\def\s{25}
\draw[0cell]
(0,0) node (a) {\phantom{\M}}
(a)+(2,0) node (b) {\phantom{\N}}
(a)+(-.8,0) node (a') {\NCat(\C,\N)}
(b)+(.8,0) node (b') {\NCat(\C,\N)}
;
\draw[1cell=.9]
(a) edge[bend left=\s] node {1} (b)
(a) edge[bend right=\s] node[swap] {\Fdgr\Edg} (b)
;
\draw[2cell]
node[between=a and b at .45, rotate=-90, 2label={above,\coudg}] {\Rightarrow}
;
\end{tikzpicture}\]
in \cref{unidg-def,coudg-def}.
\end{itemize} 
\item\label{mackey-xiv-assum-ii} We assume that $(\N,\cX)$ is a category with weak equivalences (\cref{def:mackey-heq-context}).  Using this assumption we define the wide subcategories
\begin{itemize}
\item $\cXtri \bigsubset \NCat(\C,\N)$ in \cref{cWtri-subcat}, 
\item $F^\inv\cX \bigsubset \M$ in \cref{FinvcX}, and
\item $\FinvcXtri \bigsubset \MCat(\C_E,\M)$ in \cref{diag-enr-rel-cat}.
\end{itemize} 
The relative categories in \cref{mackey-xiv-functors} are defined using $\FinvcXtri$ and $\cXtri$.
\item\label{mackey-xiv-assum-iii} We assume that 
\begin{itemize}
\item each component of $\uni \cn 1_\M \to EF$ is in $F^\inv\cX$ and
\item each component of $\cou \cn 1_\N \to FE$ is in $\cX$.
\end{itemize}
In the proof of \cref{mackey-gen-xiv}, we use the assumption about $\uni$ to prove that $\unidg$ is a relative natural transformation \cref{FunidgAx-FuniAx}.  We use the assumption about $\cou$ to prove that
\begin{itemize}
\item $\Edg$ is a relative functor \cref{cou-natural-thetaxpn} and
\item $\coudg$ is a relative natural transformation \cref{coudgPxpn-couPx}.
\end{itemize} 
\end{enumerate}
As we explain in \cref{partner-Fdgrpsix}, the relative functoriality of $\Fdgr$ is a consequence of the definitions of $\FinvcXtri$ and $\Fdgr\psi$. 
\end{explanation}

Recall that each category $\C$ enriched in a multicategory $\N$ has an opposite $\N$-category $\Cop$ (\cref{opposite-mcat}), whose composition involves the symmetric group action on $\N$.  In the next result, we assume that $E$ is a multifunctor (\cref{def:enr-multicategory-functor}), so it strictly preserves the symmetric group action.  This result gives an equivalence of homotopy theories between enriched Mackey functor categories.

\begin{theorem}\label{mackey-xiv-cor}
In the context of \cref{mackey-gen-xiv}, suppose, furthermore, that $E \cn \N \to \M$ is a multifunctor between multicategories.  Then the functors in \cref{Fdgr-Edg} applied to $\Cop$,
\[\begin{tikzpicture}[baseline={(a.base)}]
\draw[0cell]
(0,0) node (a) {\Big(\MCat\big((\C_E)^\op,\M\big) \scs \FinvcXtri\Big)}
(a)+(5.5,0) node (b) {\big(\NCat(\Cop,\N) \scs \cXtri\big),} 
;
\draw[1cell]
(a) edge[transform canvas={yshift=.6ex}] node {\Fdgr} node[swap] {\sim} (b)
(b) edge[transform canvas={yshift=-.6ex}] node {\Edg} (a)
;
\end{tikzpicture}\]
are inverse equivalences of homotopy theories.
\end{theorem}

\begin{proof}
This is \cref{mackey-gen-xiv} applied to the opposite $\N$-category $\Cop$.  By \cref{dF-opposite} the multifunctoriality of $E$ yields the equality
\[(\Cop)_E = (\C_E)^\op\]
of $\M$-categories.
\end{proof}

For completeness we end this section with the following observation that says that the functors $F$ and $E$ are inverse equivalences of homotopy theories.

\begin{proposition}\label{FE-inv-heq}
Under the assumptions of \cref{mackey-gen-xiv}, the functors
\[\begin{tikzpicture}
\draw[0cell]
(0,0) node (a) {\big(\M, F^\inv\cX)}
(a)+(3,0) node (b) {(\N,\cX)}
;
\draw[1cell]
(a) edge[transform canvas={yshift=.6ex}] node {F} node[swap] {\sim} (b)
(b) edge[transform canvas={yshift=-.6ex}] node {E} (a)
;
\end{tikzpicture}\]
are inverse equivalences of homotopy theories.
\end{proposition}

\begin{proof}
This is a much simpler variant of the proof of \cref{mackey-gen-xiv}.
\begin{itemize}
\item $F$ is a relative functor by the definition of $F^\inv\cX$ in \cref{FinvcX}.
\item To see that $E$ is a relative functor, suppose $f \cn a \to b$ is a morphism in $\cX$.  We want to show that $Ef$ is in $F^\inv\cX$, which means $FEf \in \cX$.  The naturality of $\xi$ yields the following commutative diagram in $\N$.
\[\begin{tikzpicture}[vcenter]
\def\v{-1.3}
\draw[0cell]
(0,0) node (a) {a}
(a)+(3,0) node (b) {FEa}
(a)+(0,\v) node (c) {b}
(b)+(0,\v) node (d) {FEb}
;
\draw[1cell=.9]
(a) edge node {\cou_{a}} (b)
(c) edge node {\cou_{b}} (d)
(a) edge node[swap] {f} (c)
(b) edge node {FEf} (d)
;
\end{tikzpicture}\]
Similar to \cref{cou-natural-thetaxpn},
\begin{itemize}
\item the assumption $f \in \cX$,
\item the assumption that each component of $\cou$ is in $\cX$, and
\item the 2-out-of-3 property of $\cX$
\end{itemize}  
imply that $FEf$ is in $\cX$.
\item The components of $\uni$ and $\cou$ are in $F^\inv\cX \bigsubset \M$ and $\cX \bigsubset \N$, respectively, by assumption.  So $\uni$ and $\cou$ are relative natural transformations.
\end{itemize}
Thus, by \cref{def:inverse-heq}, the functors $F$ and $E$ are inverse equivalences of homotopy theories.
\end{proof}

\chapter{Applications to Multicategories and Permutative Categories}
\label{ch:mackey_eq}
This chapter develops three main applications from \cref{ch:mackey}.
These are summarized below, where $\C$ is a small $\permcatsu$-category, $\D$ is a small $\MoneMod$-category, and each wide subcategory of stable equivalences is created as in \cref{cWtri-subcat}.
In each case, we state the equivalence of homotopy theories between enriched diagram categories that follows from \cref{mackey-gen-xiv}.
There are also variants for enriched Mackey functors that follow from \cref{mackey-xiv-cor} because, in each of the three applications, the reverse functor $E$ is a multifunctor in the symmetric sense (\cref{def:enr-multicategory-functor}).
The more precise statements in the body of this chapter give further details.

\subsection*{Pointed Multicategories and Permutative Categories}
The first application concerns the following data.
\begin{equation}\label{eq:pMult-permsu-application}
  \begin{gathered}

  \end{gathered} 
\end{equation}
\Cref{mackey-xiv-pmulticat} shows that these induce an inverse equivalence of homotopy theories between enriched diagram categories
\[
  %
\]

\subsection*{$\Mone$-Modules and Permutative Categories}
The second application concerns the following data.
\begin{equation}\label{eq:MMod-permsu-application}
  \begin{gathered}
  %
\end{gathered}
\end{equation}
\Cref{mackey-xiv-mone} shows that these induce an inverse equivalence of homotopy theories between enriched diagram categories
\[
  %
\]

\subsection*{Pointed Multicategories and $\Mone$-Modules}
The third application concerns the following data.
\begin{equation}\label{eq:pMult-MMod-application}
  \begin{gathered}
  %
\end{gathered}
\end{equation}
\Cref{mackey-pmulti-mone} shows that these induce an inverse equivalence of homotopy theories between enriched diagram categories
\[
  %
\]

\subsection*{Connection with Other Chapters}

The results in this chapter relate to broader work in homotopy theory of diagram spectra and spectral Mackey functors via the constructions in \Cref{sec:presheaf-K,sec:mult-mackey-spectra}.
Those sections develop spectral Mackey functors from enriched Mackey functors in $\permcatsu$ and $\MoneMod$.

\subsection*{Background}

In addition to the results of \cref{ch:mackey}, the main applications in this chapter depend on the following context.
The underlying inverse equivalences of homotopy theories, 
\[
  \brb{\Fst,\Endst}, \quad \brb{\Fm,\Endm}, \andspace \brb{\Mone\sma-,\Um}
\]
are developed in \cref{ch:ptmulticat-sp,ch:ptmulticat-alg}.
The closed multicategory structure for $\permcatsu$ is developed in \cref{ch:gspectra}.
The corresponding structures for $\pMulticat$ and $\MoneMod$ follow from the general discussion in \cref{sec:closed-multicat} about closed multicategorical structure on symmetric monoidal closed categories.

\subsection*{Chapter Summary}

\Cref{sec:mackey-pmulticat} describes the context and application for the data in \cref{eq:pMult-permsu-application}.
\Cref{sec:Endstdg,sec:Fstdgr} further unpack and explain the details of the functors involved.

\cref{sec:mackey-mone} describes the context and application for the data in \cref{eq:MMod-permsu-application}.
\Cref{sec:Endmdg-Fmdgr} further unpacks and explains the details of the functors involved.

\cref{sec:mackey-pmult-mone} describes the context and application for the data in \cref{eq:pMult-MMod-application}.
\Cref{sec:Monesmadgr-Umdg} further unpacks and explains the details of the functors involved.

Here is a summary table.
\reftablestretch{.98}{1.6}{
  $\brb{\Fstdgr,\Endstdg}$ inverse equivalence of homotopy theories
  & \ref{mackey-xiv-pmulticat}
  \\ \hline
  explanations of $\Endstdg$ 
  & \ref{expl:Endstdg} and \ref{expl:Endstdg-ii}
  \\ \hline
  explanations of $\Fstdgr$
  & \ref{expl:Fstdgr}, \ref{expl:Fstdgr-object}, and \ref{expl:Fstdgr-morphism} 
  \\ \hline
  $\brb{\Fmdgr,\Endmdg}$ inverse equivalence of homotopy theories
  & \ref{mackey-xiv-mone}
  \\ \hline
  explanations of $\Endmdg$
  & \ref{expl:endmdg}, \ref{expl:Endmdg-objects}, and \ref{expl:endmdg-morphisms}
  \\ \hline
  explanations of $\Fmdgr$
  & \ref{expl:Fmdgr}, \ref{expl:Fmdgr-objects}, and \ref{expl:Fmdgr-morphisms}
  \\ \hline
  $\brb{\Monesmadgr,\Umdg}$ inverse equivalence of homotopy theories
  & \ref{mackey-pmulti-mone}
  \\ \hline
  explanation of $\Umdg$
  & \ref{expl:Umdg}
  \\ \hline
  explanations of $\Monesmadgr$
  & \ref{expl:Monesmadgr} and \ref{expl:Monesmadgr-morphism}
  \\
}
 
\section[Multicategorical and Permutative Enriched Diagrams]{Homotopy Equivalent Multicategorical and Permutative Enriched Diagrams}
\label{sec:mackey-pmulticat}

In this section we apply \cref{mackey-gen-xiv,mackey-xiv-cor} to show that the categories of enriched diagrams and Mackey functors in pointed multicategories and permutative categories are connected by inverse equivalences of homotopy theories.  See \cref{mackey-xiv-pmulticat}.  As a result, left modules in $\permcatsu$ and left modules in $\pMulticat$ have equivalent homotopy theories; see \cref{expl:mackey-pmult-module}.  In \cref{sec:Endstdg,sec:Fstdgr} we explain in detail the functors that constitute this pair of inverse equivalences of homotopy theories.

\subsection*{Context}

For the context first recall the diagram 
\begin{equation}\label{mackey-xiv-setting}
\begin{tikzpicture}[baseline={(a.base)}]
\draw[0cell]
(0,0) node (a') {\pMulticat}
(a')+(.7,0) node (a) {}
(a)+(2,0) node (b) {}
(b)+(.8,.05) node (b') {\permcatsu}
;
\draw[1cell=.9]
(a) edge[transform canvas={yshift=.6ex}] node {\Fst} (b)
(b) edge[transform canvas={yshift=-.4ex}] node {\Endst} (a)
;
\end{tikzpicture}
\end{equation}
in \cref{FstEst-context} consisting of
\begin{itemize}
\item the $\Cat$-multicategory $\pMulticat$ in \cref{expl:ptmulticatcatmulticat},
\item the $\Cat$-multicategory $\permcatsu$ in \cref{thm:permcatmulticat},
\item the $\Cat$-multifunctor $\Endst$ in \cref{expl:endst-catmulti}, and
\item the non-symmetric $\Cat$-multifunctor $\Fst$ in \cref{ptmulticat-xvii}.
\end{itemize}
The two composites in \cref{mackey-xiv-setting} are connected to the respective identity functors via the following $\Cat$-multinatural transformations from \cref{ptmulticat-xx,ptmulticat-xxii}.
\begin{equation}\label{mackey-xiv-setting-ii}
\begin{tikzpicture}[baseline={(a.base)}]
\def\t{30} \def\s{27}
\draw[0cell]
(0,0) node (a) {\phantom{X}}
(a)+(-.53,0) node (a') {\pMulticat}
(a)+(1.7,0) node (b) {\phantom{X}}
(b)+(.55,.0) node (b') {\pMulticat}
;
\draw[1cell=.9]
(a) edge[bend left=\t] node {1} (b)
(a) edge[bend right=\t] node[swap] {\Endst\Fst} (b)
;
\draw[2cell]
node[between=a and b at .42, rotate=-90, 2label={above,\etast}] {\Rightarrow}
;
\begin{scope}[shift={(5.2,.05)}]
\draw[0cell]
(0,0) node (a) {\phantom{X}}
(a)+(-.65,0) node (a') {\permcatsu}
(a)+(1.7,0) node (b) {\phantom{X}}
(b)+(.65,.0) node (b') {\permcatsu}
;
\draw[1cell=.9]
(a) edge[bend left=\t] node {1} (b)
(a) edge[bend right=\t] node[swap] {\Fst\Endst} (b)
;
\draw[2cell]
node[between=a and b at .42, rotate=-90, 2label={above,\vrhost}] {\Rightarrow}
;
\end{scope}
\end{tikzpicture}
\end{equation}
The underlying categories of $\pMulticat$ and $\permcatsu$ are equipped with the relative category structures
\begin{equation}\label{pMultPerm-steq}
\brb{\pMulticat,\cSst} \andspace 
\brb{\permcatsu,\cS}
\end{equation}
in \cref{SM-Sst-def} and \cref{perm-steq}.
\begin{itemize}
\item The wide subcategory of stable equivalences 
\[\cS \bigsubset \permcatsu\]
is created by Segal $K$-theory $\Kse$ \cref{Kse}.  So a strictly unital symmetric monoidal functor $P$ is in $\cS$ if and only if $\Kse P$ is a stable equivalence of symmetric spectra.  For a small $\permcatsu$-category $\C$, the wide subcategory in \cref{mackey-xiv-pmult-functors} below
\begin{equation}\label{sStri-definition}
\cStri \bigsubset \permcatsucat\big(\C,\permcatsu\big)
\end{equation}
is defined as in \cref{cWtri-subcat} using $\cS$. 
\item The wide subcategory of $\Fst$-stable equivalences
\[\cSst = \Fst^\inv(\cS) \bigsubset \pMulticat\] 
is created by $\Fst$.  The wide subcategory in \cref{mackey-xiv-pmult-functors} below
\begin{equation}\label{sSsttri-definition}
\cSsttri \bigsubset \pMulticatcat\big(\C_{\Endst},\pMulticat\big)
\end{equation}
is defined as in \cref{cWtri-subcat} using $\cSst$. 
\end{itemize}

\subsection*{Equivalences of Homotopy Theories}

In \cref{ptmulticat-xxiii} we observe that the pair $(\Fst,\Endst)$ induces inverse equivalences of homotopy theories between the respective categories of non-symmetric $\Q$-algebras for each small non-symmetric $\Cat$-multicategories $\Q$.  The following observation extends the inverse equivalences of homotopy theories $(\Fst,\Endst)$ to categories of enriched diagrams and Mackey functors.  In \cref{sec:Endstdg,sec:Fstdgr} we further explain the functors $\Endstdg$ and $\Fstdgr$.

\begin{theorem}\label{mackey-xiv-pmulticat}
Suppose $\C$ is a small $\permcatsu$-category.  Then the functors
\begin{equation}\label{mackey-xiv-pmult-functors}
\begin{tikzpicture}[baseline={(a.base)}]
\draw[0cell=.8]
(0,0) node (a) {\Big(\pMulticatcat\big(\C_{\Endst},\pMulticat\big) \scs \cSsttri\Big)}
(a)+(6,0) node (b) {\Big(\permcatsucat\big(\C,\permcatsu\big) \scs \cStri\Big),} 
;
\draw[1cell=.8]
(a) edge[transform canvas={yshift=.6ex}] node (f) {\Fstdgr} node[swap] {\sim} (b)
(b) edge[transform canvas={yshift=-.5ex}] node (e) {\Endstdg} (a)
;
\end{tikzpicture}
\end{equation}
defined by the data in \cref{mackey-xiv-setting,mackey-xiv-setting-ii,pMultPerm-steq,sStri-definition,sSsttri-definition}, are inverse equivalences of homotopy theories.  

Moreover, the variant with $(\C_{\Endst})^\op$ and $\Cop$ replacing, respectively, $\C_{\Endst}$ and $\C$ is also true.
\end{theorem}

\begin{proof}
The first assertion is an instance of \cref{mackey-gen-xiv}, which is applicable in the current setting as we now explain.  Following the summary in \cref{expl:mackey-xiv-assumptions}, first we verify that \cref{def:mackey-gen-context} \cref{mackey-context-i,mackey-context-ii,mackey-context-iii,mackey-context-iv,mackey-context-v} are satisfied in the current context.
\begin{romenumerate}
\item\label{mackey-xiv-pm-i} $\M = \pMulticat$ is a closed multicategory by
\begin{itemize}
\item \cref{smclosed-closed-multicat} and
\item the fact that it is a symmetric monoidal closed category (\cref{thm:pmulticat-smclosed}).
\end{itemize}  
By \cref{perm-closed-multicat}, $\N = \permcatsu$ is a closed multicategory.
\item\label{mackey-xiv-pm-ii} $\C$ is, by assumption, a small $\permcatsu$-category.
\item\label{mackey-xiv-pm-iii} $F = \Fst$ in \cref{mackey-xiv-setting} is a non-symmetric multifunctor by \cref{ptmulticat-xvii}, and $E = \Endst$ is a multifunctor by \cref{proposition:n-lin-equiv}.
\item\label{mackey-xiv-pm-iv} $\uni = \etast$ and $\cou = \vrhost$ in \cref{mackey-xiv-setting-ii} are multinatural transformations by \cref{ptmulticat-xx,ptmulticat-xxii}, respectively.
\item\label{mackey-xiv-pm-v} In the current setting, the condition \cref{counit-compatibility} is the equality of the following two pointed multifunctors for each pair of objects $x,y \in \C$. 
\[\begin{tikzpicture}
\draw[0cell=1]
(0,0) node (a) {\Endst \C(x,y)}
(a)+(4.5,0) node (b) {\Endst \Fst \Endst \C(x,y)}
;
\draw[1cell=.9]
(a) edge[transform canvas={yshift=.6ex}] node {\etast_{\Endst \C(x,y)}} (b)
(a) edge[transform canvas={yshift=-.4ex}] node[swap] {\Endst \vrhost_{\C(x,y)}} (b)
;
\end{tikzpicture}\]
This equality holds by \cref{etaEEvrho} because each hom object $\C(x,y)$ is a small permutative category.
\end{romenumerate}
Thus \cref{def:mackey-gen-context} \cref{mackey-context-i,mackey-context-ii,mackey-context-iii,mackey-context-iv,mackey-context-v} hold in the context of \cref{mackey-xiv-setting,mackey-xiv-setting-ii,pMultPerm-steq,sStri-definition,sSsttri-definition}.

Next, the only assumption in \cref{def:mackey-heq-context} is that the relative category
\[(\N, \cX) = \brb{\permcatsu,\cS}\]
is a category with weak equivalences (\cref{definition:rel-cat} \eqref{relcat-vi}).  This is true for the following two reasons.
\begin{itemize}
\item The wide subcategory $\cS \bigsubset \permcatsu$ is created by a functor, namely, Segal $K$-theory \cref{perm-steq}
\[\Kse \cn \permcatsu \to \Spc.\]
\item The class of stable equivalences in $\Spc$ contains all the isomorphisms and has the 2-out-of-3 property.
\end{itemize} 
In the current setting there are equalities of wide subcategories
\[F^\inv\cX = \Fst^\inv(\cS) = \cSst \bigsubset \pMulticat.\]
The data in \cref{mackey-xiv-pmult-functors} are those in \cref{mackey-xiv-functors} in the current context.

Finally, each component of $\vrhost$ is a stable equivalence in $\permcatsu$ by \cref{remark:steq} \cref{it:steq-2} because it admits a left adjoint by \cref{ptmulticat-prop-viii}.  Moreover, in the proof of \cref{ptmulticat-thm-x} we explain that each component of $\etast$ is an $\Fst$-stable equivalence in $\pMulticat$.  Thus \cref{mackey-gen-xiv} is applicable in the current setting, proving the first assertion.

The second assertion about $(\C_{\Endst})^\op$ and $\Cop$ is an instance of \cref{mackey-xiv-cor}.  It is applicable because $\Endst$ is a multifunctor (\cref{proposition:n-lin-equiv}).
\end{proof}

\begin{explanation}[Homotopy Equivalent Categories of Modules]\label{expl:mackey-pmult-module}
By \cref{C-diagram-partner,C-diag-morphism-pn}, for each small $\permcatsu$-category $\C$, the functors in \cref{mackey-xiv-pmult-functors} are inverse equivalences of homotopy theories between
\begin{itemize}
\item left $\C$-modules in $\permcatsu$ and
\item left $\C_{\Endst}$-modules in $\pMulticat$.
\end{itemize}
We explain the functors $\Endstdg$ and $\Fstdgr$ in detail in \cref{sec:Endstdg,sec:Fstdgr}.
\end{explanation}

\begin{remark}[Non-Existence of Unpointed Version]\label{rk:unptd-mackey}
We do not know of any analog of \cref{mackey-xiv-pmulticat} for the symmetric monoidal closed category $\Multicat$ (\cref{theorem:multicat-sm-closed}), whose objects are small multicategories.  The reason is that the multinatural transformation $\vrhost$ in \cref{mackey-xiv-setting-ii} is necessary to define the functor $\Fstdgr$.  As we discuss in \cref{rk:rhonotunital}, for each permutative category $\C$, the symmetric monoidal functor $\vrho_\C$ is \emph{not} strictly unital.  Thus we \emph{cannot} use $\vrho$ to define a multinatural transformation $1_{\permcatsu} \to \Fr\,\End$.
\end{remark}

\section[Multicategorical Enriched Diagrams]{Permutative to Multicategorical Enriched Diagrams}
\label{sec:Endstdg}

In this section we explain in detail the equivalence of homotopy theories in \cref{mackey-xiv-pmult-functors}
\[\begin{tikzpicture}[baseline={(a.base)}]
\draw[0cell=.9]
(0,0) node (a) {\pMulticatcat\big(\C_{\Endst} \scs \pMulticat\big)}
(a)+(6,0) node (b) {\permcatsucat\big(\C,\permcatsu\big)} 
;
\draw[1cell=.9]
(b) edge node {\sim} node[swap] {\Endstdg} (a)
;
\end{tikzpicture}\]
that produces pointed multicategorical enriched diagrams from permutative enriched diagrams.  To simplify the notation, we use the following abbreviations throughout this section.
\begin{equation}\label{pMpsu}
\pM = \pMulticat \qquad\qquad \psu = \permcatsu
\end{equation}
This section is organized as follows.
\begin{itemize}
\item \cref{expl:Endstdg} describes $\Endstdg$ in terms of $\dEndst$ and $\Endstse$.
\item \cref{expl:dEndst} describes the 2-functor $\dEndst$.
\item \cref{expl:Endstse} describes the standard enrichment $\Endstse$.
\item \cref{expl:Endstdg-ii} summarizes \cref{expl:Endstdg,expl:dEndst,expl:Endstse}.
\end{itemize}

\begin{explanation}[Unpacking $\Endstdg$]\label{expl:Endstdg}\index{endomorphism!multicategory!diagram change of enrichment}\index{multicategory!endomorphism - diagram change of enrichment}\index{permutative category!endomorphism multicategory!diagram change of enrichment}\index{diagram!change of enrichment!endomorphism multicategory -}
The diagram change-of-enrichment functor in \cref{mackey-xiv-pmult-functors}
\[\Endstdg \cn \psucat(\C,\psu) \to \pMcat(\C_{\Endst},\pM)\]
is defined in \cref{diag-change-enr-assign} and verified in \cref{gspectra-thm-v}.  To understand its assignments on objects and morphisms, consider
\begin{itemize}
\item $\psu$-functors $A,B \cn \C \to \psu$ (\cref{expl:perm-enr-functors}) and
\item a $\psu$-natural transformation $\psi \cn A \to B$  (\cref{expl:perm-nattr})
\end{itemize}
as in the left diagram below.
\begin{equation}\label{Endstdg-def}
\begin{tikzpicture}[baseline={(a.base)}]
\def\t{25} \def\d{.7}
\draw[0cell=.9]
(0,0) node (a) {\C}
(a)+(1.5,0) node (b1) {\phantom{C}}
(b1)+(.1,0) node (b) {\psu}
(b)+(\d,0) node (x) {}
(x)+(1.4,0) node (y) {}
(y)+(\d,0) node (a') {\C_{\Endst}}
(a')+(.25,0) node (a'') {\phantom{C}}
(a'')+(2,0) node (b2) {\phantom{C}}
(b2)+(.55,0) node (b') {(\psu)_{\Endst}}
(b')+(2,0) node (c) {\pM}
;
\draw[1cell=.85]
(a) edge[bend left=30] node {A} (b1)
(a) edge[bend right=30] node[swap] {B} (b1)
(x) edge[|->] node {\Endstdg} (y)
;
\draw[2cell=.9]
node[between=a and b1 at .42, rotate=-90, 2label={above,\psi}] {\Rightarrow}
;
\draw[1cell=.85]
(a'') edge[bend left=\t] node {A_{\Endst}} (b2)
(a'') edge[bend right=\t] node[swap] {B_{\Endst}} (b2)
(b') edge node {\Endstse} (c)
;
\draw[2cell=.9]
node[between=a'' and b2 at .35, rotate=-90, 2label={above,\psi_{\Endst}}] {\Rightarrow}
;
\end{tikzpicture}
\end{equation}
Then $\Endstdg$ sends $A$, $B$, and $\psi$ to the composites and whiskering as in the right diagram in \cref{Endstdg-def}.  In other words, the functor $\Endstdg$ 
\begin{itemize}
\item first applies the change of enrichment $\dEndst$ and then
\item composes or whiskers with the standard enrichment $\Endstse$.
\end{itemize}
We describe $\dEndst$ and $\Endstse$ further in \cref{expl:dEndst,expl:Endstse} below.  Then we summarize the discussion in \cref{expl:Endstdg-ii}.
\end{explanation}

\begin{explanation}[Unpacking $\dEndst$]\label{expl:dEndst}
We describe the change-of-enrichment 2-functor in \cref{Endstdg-def}
\[\dEndst \cn \psucat \to \pMcat\]
by interpreting \cref{def:mult-change-enr} for the multifunctor (\cref{expl:endst-catmulti})
\[\Endst \cn \psu \to \pM.\]
The existence of $\dEndst$ is an instance of \cref{mult-change-enrichment}.

\medskip
\emph{Objects}.
First we consider a $\psu$-category $(\D,\mcomp^\D)$ (\cref{expl:perm-enr-cat}).
\begin{itemize}
\item The $\pM$-category $\D_{\Endst}$ has the same objects as $\D$.  So $(\psu)_{\Endst}$ has small permutative categories as objects.
\item For each pair of objects $a,b \in \D$, its hom object is
\[(\D_{\Endst})(a,b) = \Endst\D(a,b) \inspace \pM.\]
So for small permutative categories $\X$ and $\Y$, by \cref{permcat-selfenr} there is a hom object
\[(\psu)_{\Endst}(\X,\Y) = \Endst\psu(\X,\Y) \inspace \pM.\]
\begin{itemize}
\item In $\psu(\X,\Y)$ the objects are strictly unital symmetric monoidal functors $\X \to \Y$.
\item The morphisms are monoidal natural transformations.
\item The monoidal structure is defined pointwise in the codomain $\Y$.
\end{itemize}
$\Endst\psu(\X,\Y)$ is the pointed multicategory associated to the permutative category $\psu(\X,\Y)$.
\item For objects $a,b,c \in \D$, the composition binary multimorphism in $\pM$
\[\Endst\big(\mcomp^\D_{a,b,c}\big) \cn \big(\Endst \D(b,c) \scs \Endst \D(a,b)\big) \to \Endst\D(a,c)\]
is the image under $\Endst$ (\cref{proposition:n-lin-equiv}) of the composition in $\D$, 
\[\mcomp^\D_{a,b,c} \cn \D(b,c) \times \D(a,b) \to \D(a,c).\]
The latter is a bilinear functor of permutative categories (\cref{def:nlinearfunctor}).  The composition in $\psu$ as a $\psu$-category (\cref{def:perm-selfenr}) extends composition of strictly unital symmetric monoidal functors and is a bilinear functor.
\end{itemize} 

\medskip
\emph{1-Cells}.
For a $\psu$-functor $A \cn \C \to \psu$ as in \cref{Endstdg-def}, the $\pM$-functor
\[A_{\Endst} \cn \C_{\Endst} \to (\psu)_{\Endst}\]
has the same object assignment as $A$.  In other words, for each object $x$ in $\C$,
\begin{equation}\label{A-Endst-x}
(A_{\Endst})x = Ax \inspace \psu.
\end{equation}
For objects $x,y \in \C$, the $(x,y)$-component pointed multifunctor
\begin{equation}\label{A-Endst-xy}
(A_{\Endst})_{x,y} = \Endst(A_{x,y}) \cn \Endst\C(x,y) \to \Endst\psu(Ax,Ay)
\end{equation}
is the image under $\Endst$ of the $(x,y)$-component strictly unital symmetric monoidal functor
\[\big(A_{x,y} \scs A_{x,y}^2\big) \cn \C(x,y) \to \psu(Ax,Ay),\]
as defined in \cref{ex:endstc} \cref{EndstP}.  The same explanation also applies to the $\psu$-functor $B \cn \C \to \psu$.

\medskip
\emph{2-Cells}.
A $\psu$-natural transformation $\psi \cn A \to B$ as in \cref{Endstdg-def} is determined by its components.  For each object $x \in \C$, the $x$-component is a nullary multimorphism
\[\psi_x \cn \ang{} \to \psu(Ax,Bx) \inspace \psu.\]
By \cref{def:nlinearfunctor,definition:permcatsus-homcat}, such a nullary multimorphism $\psi_x$ is a 0-linear functor to $\psu(Ax,Bx)$.  This, in turn, means a choice of an object in the permutative category $\psu(Ax,Bx)$.  In other words, each component 
\begin{equation}\label{psix-Ax-to-Bx}
\psi_x \cn Ax \to Bx
\end{equation}
is a strictly unital symmetric monoidal functor.  

Under the change of enrichment $\dEndst$, the $\pM$-natural transformation in \cref{Endstdg-def}
\[\psi_{\Endst} \cn A_{\Endst} \to B_{\Endst}\]
has, for each object $x$ in $\C$, $x$-component nullary multimorphism
\begin{equation}\label{psi-Endst-x}
(\psi_{\Endst})_x \cn \ang{} \to \Endst\psu(Ax,Bx) \inspace \pM.
\end{equation}
Since the multicategory structure on $\pM$ is induced by its symmetric monoidal structure, such a nullary multimorphism is a pointed multifunctor
\[(\psi_{\Endst})_x \cn \Mtup = \Mtu \bincoprod \Mterm \to \Endst\psu(Ax,Bx)\]
from the smash unit $\Mtup$ in \cref{eq:smashunit}.  Preservation of basepoints and $\Mtu$ being the initial operad imply that $(\psi_{\Endst})_x $ is a choice of an object in $\Endst\psu(Ax,Bx)$, which means an object in $\psu(Ax,Bx)$.  This, in turn, means a strictly unital symmetric monoidal functor $Ax \to Bx$, which is given by $\psi_x$ in \cref{psix-Ax-to-Bx}.

In summary, the components of the $\pM$-natural transformation $\psi_{\Endst}$ are the components of the $\psu$-natural transformation $\psi$.
\end{explanation}

\begin{explanation}[Unpacking $\Endstse$]\label{expl:Endstse}
We describe the last arrow in \cref{Endstdg-def}
\[\Endstse \cn (\psu)_{\Endst} \to \pM.\]
This is the standard enrichment of $\Endst$, which is an $\pM$-functor (\cref{gspectra-thm-iii}).  Its object assignment is the same as that of $\Endst$ (\cref{ex:endstc}).  In other words, 
\begin{equation}\label{Endstse-X}
\Endstse\,\X = \Endst\X \inspace \pM
\end{equation}
for each small permutative category $\X$.

For small permutative categories $\X$ and $\Y$, the $(\X,\Y)$-component of $\Endstse$ is a pointed multifunctor
\begin{equation}\label{Endstse-XY}
(\Endstse)_{\X,\Y} \cn \Endst\psu(\X,\Y) \to \pHom(\Endst\X, \Endst\Y).
\end{equation}
In the codomain, $\pHom$ is the pointed internal hom \cref{eq:multicat-pHom} in $\pM$.
\begin{itemize}
\item The objects of $\pHom(?,?)$ are pointed multifunctors.
\item Its multimorphisms are pointed transformations (\cref{ex:pointedhommulticat}).
\end{itemize}
Next we describe the component pointed multifunctor $(\Endstse)_{\X,\Y}$ in two equivalent ways.

First, $(\Endstse)_{\X,\Y}$ is uniquely determined by its adjoint, which is the arrow $\Endst(\ev_{\X,\Y})$ in the commutative diagram \cref{EndstseXY-diag} in $\pM$.
\begin{equation}\label{EndstseXY-diag}
\begin{tikzpicture}[vcenter]
\def\h{3} \def\v{1.3} \def\t{20}
\draw[0cell=.9]
(0,0) node (a) {\Endst\psu(\X,\Y) \sma \Endst\X}
(a)+(\h,\v) node (b) {\pHom(\Endst\X, \Endst\Y) \sma \Endst\X}
(b)+(\h,-\v) node (c) {\Endst\Y} 
;
\draw[1cell=.9]
(a) edge node {\Endst(\ev_{\X,\Y})} (c)
(a) edge[bend left=\t] node[pos=.3] {(\Endstse)_{\X,\Y} \sma \opu} (b)
(b) edge[bend left=\t] node {\ev^{\pM}} (c)
;
\end{tikzpicture}
\end{equation}
The diagram \cref{EndstseXY-diag} is the diagram \cref{Fse-xy-diag} for $\Endstse$.
\begin{itemize}
\item In \cref{EndstseXY-diag} $\sma$ is the smash product \cref{eq:multicat-smash-pushout} in $\pM$.
\item The evaluation for permutative categories
\[\ev_{\X,\Y} \cn \psu(\X,\Y) \times \X \to \Y\] 
is the bilinear functor in \cref{evCD}.
\item $\ev^{\pM}$ is the evaluation \cref{evaluation} in the symmetric monoidal closed category $\pM$ (\cref{thm:pmulticat-smclosed}).
\end{itemize}

Alternatively, we obtain from \cref{EndstseXY-diag} a direct description of the pointed multifunctor $(\Endstse)_{\X,\Y}$ as follows.  An object in $\Endst\psu(\X,\Y)$ is an object in $\psu(\X,\Y)$, which is a strictly unital symmetric monoidal functor
\[(Q,Q^2) \cn \X \to \Y.\]
Its image under $(\Endstse)_{\X,\Y}$ is the pointed multifunctor
\begin{equation}\label{Endstse-Q}
(\Endstse)_{\X,\Y}(Q,Q^2) = \Endst(Q,Q^2) \cn \Endst\X \to \Endst\Y
\end{equation}
obtained from $(Q,Q^2)$ by applying $\Endst$ (\cref{ex:endstc} \cref{EndstP}).

For $n \geq 0$ an $n$-ary multimorphism
\[\begin{split}
\theta \in &~ \big(\Endst\psu(\X,\Y)\big) \Big( \bang{(Q_i, Q_i^2)}_{i=1}^n \sscs (Q,Q^2) \Big)\\
&= \psu(\X,\Y)\big( \txoplus_{i=1}^n (Q_i, Q_i^2) \scs (Q,Q^2) \big)
\end{split}\]
is a monoidal natural transformation as follows.
\[\begin{tikzpicture}
\def\t{25}
\draw[0cell]
(0,0) node (a) {\X}
(a)+(1.8,0) node (b) {\Y}
;
\draw[1cell=.8]
(a) edge[bend left=\t] node {\txoplus_{i=1}^n Q_i} (b)
(a) edge[bend right=\t] node[swap] {Q} (b) 
;
\draw[2cell]
node[between=a and b at .45, rotate=-90, 2label={above,\theta}] {\Rightarrow}
;
\end{tikzpicture}\]
The domain of $\theta$ is the strictly unital symmetric monoidal functor 
\[\big( \txoplus_{i=1}^n Q_i \scs (\txoplus_{i=1}^n Q_i)^2 \big) \cn \X \to \Y.\]
\begin{itemize}
\item As a functor the sum is taken pointwise in $\Y$, so
\[\big( \txoplus_{i=1}^n Q_i \big)x = \txoplus_{i=1}^n (Q_i x) \forspace x \in \X.\]
\item $\txoplus_{i=1}^n Q_i$ is strictly unital because each $Q_i$ is so.
\item For objects $x,x' \in \X$, the $(x,x')$-component of its monoidal constraint is the following composite in $\Y$, with the isomorphism permuting the objects using the braiding in $\Y$.
\[\begin{tikzpicture}
\def\h{2.75} \def\v{1.1} \def\t{15}
\draw[0cell=.9]
(0,0) node (a) {\big(\txoplus_{i=1}^n Q_ix\big) \oplus \big(\txoplus_{i=1}^n Q_ix'\big)}
(a)+(\h,-\v) node (b) {\txoplus_{i=1}^n \big( Q_ix \oplus Q_i x'\big)}
(b)+(\h,\v) node (c) {\txoplus_{i=1}^n Q_i (x \oplus x')}
;
\draw[1cell=.9]
(a) edge node {(\txoplus_{i=1}^n Q_i)^2_{x,x'}} (c)
(a) edge[bend right=\t] node[swap,pos=.4] {\iso} (b)
(b) edge[bend right=\t] node[swap,pos=.7] {\txoplus_{i=1}^n (Q_i^2)_{x,x'}} (c)
;
\end{tikzpicture}\]
\end{itemize}
For each object $x \in \X$, the $x$-component of $\theta$ is a morphism
\begin{equation}\label{theta-x-component}
\theta_x \cn \txoplus_{i=1}^n Q_ix \to Qx \inspace \Y.
\end{equation}

Applying $(\Endstse)_{\X,\Y}$ to $\theta$ yields the $n$-ary pointed transformation (\cref{ex:pointedhommulticat})
\[(\Endstse)_{\X,\Y}(\theta) \in \pHom(\Endst\X,\Endst\Y) \big(\ang{\Endst Q_i}_{i=1}^n \sscs \Endst Q \big).\]
For each object $x \in \Endst\X$, meaning $x \in \X$, the $x$-component of $(\Endstse)_{\X,\Y}(\theta)$ is the $n$-ary multimorphism
\begin{equation}\label{Endstse-XY-thetax}
\begin{aligned}
(\Endstse)_{\X,\Y}(\theta)_x 
\in &~ (\Endst \Y)\Big( \bang{(\Endst Q_i)x}_{i=1}^n \sscs (\Endst Q)x \Big)\\
& = \Y\big(\txoplus_{i=1}^n Q_ix \scs Qx\big)
\end{aligned}
\end{equation}
given by the $x$-component $\theta_x$ in \cref{theta-x-component}.  In summary, the components of $(\Endstse)_{\X,\Y}(\theta)$ are the components of $\theta$.
\end{explanation}

\begin{explanation}[Back to $\Endstdg$]\label{expl:Endstdg-ii}
We summarize \cref{expl:Endstdg,expl:dEndst,expl:Endstse}.  For a $\psu$-functor $A \cn \C \to \psu$ as in \cref{Endstdg-def}, the composite $\pM$-functor
\[
\]
has, for each object $x \in \C$, object assignment
\[\big(\Endstdg A\big)x = \Endst(Ax) \inspace \pM\]
by \cref{A-Endst-x,Endstse-X}.  For objects $x,y \in \C$, the $(x,y)$-component pointed multifunctor of $\Endstdg A$ is the composite
\begin{equation}\label{EndstdgA-xy-comp}
%
\end{equation}
by \cref{A-Endst-xy,Endstse-XY}.

For a $\psu$-natural transformation $\psi \cn A \to B$ as in \cref{Endstdg-def}, the $\pM$-natural transformation $\Endstdg\psi$ is the whiskering below.
\[%
\]
For each object $x \in \C$, its $x$-component nullary multimorphism is the following composite in $\pM$.
\[%
\]
By \cref{psi-Endst-x,Endstse-Q}, this is given by the pointed multifunctor
\begin{equation}\label{Endst-psi-x}
\big(\Endstdg \psi\big)_x = \Endst(\psi_x) \cn \Endst(Ax) \to \Endst(Bx)
\end{equation}
obtained from the strictly unital symmetric monoidal functor $\psi_x$ in \cref{psix-Ax-to-Bx} by applying $\Endst$.
\end{explanation}

\section[Permutative Enriched Diagrams]{Multicategorical to Permutative Enriched Diagrams}
\label{sec:Fstdgr}

We continue to use the abbreviations in \cref{pMpsu}, so
\[\pM = \pMulticat \andspace \psu = \permcatsu.\]
In this section we explain in detail the equivalence of homotopy theories in \cref{mackey-xiv-pmult-functors}
\[\begin{tikzpicture}[baseline={(a.base)}]
\draw[0cell=1]
(0,0) node (a) {\pMcat\big(\C_{\Endst} \scs \pM\big)}
(a)+(4.5,0) node (b) {\psucat\big(\C,\psu\big)} 
;
\draw[1cell=.9]
(a) edge node[swap] {\sim} node {\Fstdgr} (b)
;
\end{tikzpicture}\]
that produces permutative enriched diagrams from pointed multicategorical enriched diagrams.  This section is organized as follows.
\begin{itemize}
\item \cref{expl:Fstdgr} describes $\Fstdgr$ in terms of $\Fstdg$ and $\C^*_{\vrhost}$.
\item \cref{expl:Fstdgr-object} describes $\Fstdgr$ on objects.
\item \cref{expl:Fstdgr-morphism} describes $\Fstdgr$ on morphisms.
\end{itemize}

\begin{explanation}[Unpacking $\Fstdgr$]\label{expl:Fstdgr}
The functor $\Fstdgr$ is an instance of the functor $\Fdgr$ \cref{Fdgr-def} defined with the non-symmetric multifunctor (\cref{ptmulticat-xvii})
\[F = \Fst \cn \pM \to \psu\]
and the multinatural transformation $\cou = \vrhost$ given by 
\[\begin{tikzpicture}[baseline={(a.base)}]
\def\t{23}
\draw[0cell]
(0,0) node (a) {\psu}
(a)+(2,0) node (b) {\psu}
;
\draw[1cell=.85]
(a) edge[bend left=\t] node {1} (b)
(a) edge[bend right=\t] node[swap] {\Fst\Endst} (b)
;
\draw[2cell]
node[between=a and b at .42, rotate=-90, 2label={above,\vrhost}] {\Rightarrow}
;
\end{tikzpicture}\]
in \cref{ptmulticat-xxii}.  
\begin{itemize}
\item $\Fst$ is defined on objects and multimorphisms in \cref{def:Fst-permutative,def:Fst-multi}, respectively.
\item The components of $\vrhost$ are defined in \cref{vrhostC-def}.
\end{itemize}
By definition \cref{Fdgr-def} the functor $\Fstdgr$ is the following composite.
\begin{equation}\label{Fstdgr-definition}
\begin{tikzpicture}[vcenter]
\def\h{2.5} \def\v{1.1} \def\t{15} \def\w{.6}
\draw[0cell=.9]
(0,0) node (a) {\pMcat(\C_{\Endst} \scs \pM)}
(a)+(\h,-\v) node (b) {\psucat(\C_{\Fst\Endst} \scs \psu)}
(b)+(\h,\v) node (c) {\psucat(\C,\psu)}
;
\draw[1cell=.9]
(a) edge node {\Fstdgr} (c)
(a) edge[bend right=\t] node[swap,pos=.2] {\Fstdg} (b)
(b) edge[bend right=\t] node[swap,pos=.8] {\C_{\vrhost}^*} (c)
;
\end{tikzpicture}
\end{equation}
The two constituent functors in \cref{Fstdgr-definition} are as follows.
\begin{itemize}
\item $\Fstdg$ is the diagram change-of-enrichment functor (\cref{gspectra-thm-v}) of $\Fst$ at the $\pM$-category $\C_{\Endst}$.  The latter is the image of the $\psu$-category $\C$ under the change of enrichment $\dEndst$ (\cref{expl:dEndst}).
\item $\C_{\vrhost}^*$ is defined by pre-composition and whiskering with the $\psu$-functor
\begin{equation}\label{C-vrhost}
\C_{\vrhost} \cn \C \to \C_{\Fst\Endst} = (\C_{\Endst})_{\Fst}.
\end{equation}
This $\psu$-functor is the $\C$-component of the 2-natural transformation (\cref{dtheta-twonat})
\[\begin{tikzpicture}[baseline={(a.base)}]
\def\s{25}
\draw[0cell]
(0,0) node (a) {\phantom{\M}}
(a)+(2.5,0) node (b) {\phantom{\N}}
(a)+(-.4,0) node (a') {\psucat}
(b)+(.4,0) node (b') {\psucat}
;
\draw[1cell=.85]
(a) edge[bend left=\s] node {1} (b)
(a) edge[bend right=\s] node[swap] {\dFstEndst = \dFst \,\dEndst} (b)
;
\draw[2cell]
node[between=a and b at .39, rotate=-90, 2label={above,\dvrhost}] {\Rightarrow}
;
\end{tikzpicture}\]
induced by the multinatural transformation $\vrhost \cn 1_{\psu} \to \Fst\Endst$.  Here $\dFst$ is the change-of-enrichment 2-functor (\cref{ex:dFst})
\begin{equation}\label{dFst-pmcat-psucat}
\dFst \cn \pMcat \to \psucat
\end{equation}
induced by the non-symmetric multifunctor $\Fst \cn \pM \to \psu$.
\end{itemize}
We describe $\Fstdgr$ on objects and morphisms in \cref{expl:Fstdgr-object,expl:Fstdgr-morphism}, respectively.
\end{explanation}

\begin{explanation}[$\Fstdgr$ on Objects]\label{expl:Fstdgr-object}
Consider an $\pM$-functor (\cref{def:enriched-functor})
\[A \cn \C_{\Endst} \to \pM.\]
By \cref{FdgrA-diag} the $\psu$-functor $\Fstdgr A$ is the following composite.
\begin{equation}\label{FstdgrA-diag}
\begin{tikzpicture}[baseline={(a.base)}]
\def\h{2} \def\w{.6}
\draw[0cell=1]
(0,0) node (a) {\C} 
(a)+(\h,0) node (b) {\C_{\Fst\Endst}}
(b)+(1.8,0) node (b') {(\C_{\Endst})_{\Fst}}
(b')+(\h+.5,0) node (c) {(\pM)_{\Fst}}
(c)+(\h,0) node (d) {\psu}
;
\draw[1cell=2]
(b) edge[-,double equal sign distance] (b')
;
\draw[1cell=.9]
(a) edge node {\C_{\vrhost}} (b)
(b') edge node {A_{\Fst}} (c)
(c) edge node {\Fstse} (d)
;
\draw[1cell=.9]
(a) [rounded corners=3pt] |- ($(b')+(0,\w)$)
-- node[pos=.25] {\Fstdgr A} ($(c)+(0,\w)$) -| (d)
;
\end{tikzpicture}
\end{equation}
The constituent $\psu$-functors in \cref{FstdgrA-diag} are as follows.
\begin{romenumerate}
\item $\C_{\vrhost}$ is the $\psu$-functor in \cref{C-vrhost}.  
\begin{itemize}
\item Its object assignment is the identity function.
\item For objects $x,y \in \C$, its $(x,y)$-component is the strictly unital symmetric monoidal functor
\begin{equation}\label{vrhost-Cxy}
\vrhost_{\C(x,y)} \cn \C(x,y) \to \Fst\Endst\C(x,y)
\end{equation}
given by the $\C(x,y)$-component of $\vrhost$ \cref{vrhostC-def}.
\end{itemize}
\item $A_{\Fst}$ is the change of enrichment of $A$ under $\dFst$ (\cref{ex:dFst}).
\begin{itemize}
\item Its object assignment is the same as that of $A$.
\item For objects $x,y \in \C$, its $(x,y)$-component is the strict symmetric monoidal functor (\cref{def:Fst-onecells})
\begin{equation}\label{Fst-Axy}
\Fst A_{x,y} \cn \Fst\Endst\C(x,y) \to \Fst\pHom(Ax,Ay)
\end{equation}
given by applying $\Fst$ to the $(x,y)$-component pointed multifunctor of $A$.  Here $\pHom$ is the pointed internal hom for small pointed multicategories \cref{eq:multicat-pHom}.
\end{itemize}
\item $\Fstse$ is the standard enrichment $\psu$-functor of $\Fst \cn \pM \to \psu$ in \cref{Fstse-definition}.
\begin{itemize}
\item Its object assignment is the same as that of $\Fst$.
\item For small pointed multicategories $\X$ and $\Y$, its $(\X,\Y)$-component is the strictly unital symmetric monoidal functor
\begin{equation}\label{Fstse-XY}
\Fstse_{\X,\Y} = \pn{\big(\Fst(\ev^{\pM}_{\X,\Y})\big)} \cn 
\Fst\pHom(\X,\Y) \to \psu(\Fst\X,\Fst\Y).
\end{equation}
\end{itemize}
\end{romenumerate}

Combining \cref{FstdgrA-diag,vrhost-Cxy,Fst-Axy,Fstse-XY}, the object assignment of $\Fstdgr A$ is given by, for each object $x \in \C$,
\[\big(\Fstdgr A\big)x = \Fst(Ax) \inspace \psu.\]
For objects $x,y \in \C$, its $(x,y)$-component strictly unital symmetric monoidal functor is the following composite.
\begin{equation}\label{FstdgrA-component}

\end{equation}
The diagram \cref{FstdgrA-component} is the diagram \cref{FdgrA-component} in the current context.
\end{explanation}

\begin{explanation}[$\Fstdgr$ on Morphisms]\label{expl:Fstdgr-morphism}
Consider an $\pM$-natural transformation $\psi$ (\cref{def:enriched-natural-transformation}) as in the left diagram below.
\[%
\]
The $\psu$-natural transformation $\Fstdgr\psi$, as in the right diagram above, is the following whiskering.
\begin{equation}\label{Fstdgrpsi-def}
%
\end{equation}
The whiskering \cref{Fstdgrpsi-def} is the one in \cref{Fdgrpsi-def} in the current context.  
\begin{itemize}
\item The $\psu$-functors $\C_{\vrhost}$ and $\Fstse$ are as in \cref{FstdgrA-diag}.
\item $\psi_{\Fst}$ is the change of enrichment of $\psi$ under $\dFst$ in \cref{ex:dFst}.
\end{itemize}

Next we describe the components of $\Fstdgr\psi$ explicitly.   For each object $x \in \C$, the $x$-component of $\psi$ is a pointed multifunctor
\[\psi_x \cn \Mtup = \Mtu \bincoprod \Mterm \to \pHom(Ax,Bx)\]
from the smash unit $\Mtup$ in \cref{eq:smashunit}.  Preservation of basepoint and $\Mtu$ being the initial operad imply that $\psi_x$ is equivalent to a choice of an object in $\pHom(Ax,Bx)$.  This means a pointed multifunctor
\begin{equation}\label{psi-xcomponent-AxBx}
\psi_x \cn Ax \to Bx.
\end{equation}
The $x$-component of $\Fstdgr\psi$ is a nullary multimorphism
\[\big(\Fstdgr\psi\big)_x \cn \ang{} \to \psu\big(\Fst(Ax),\Fst(Bx)\big) \inspace \psu.\]
This means a choice of an object in $\psu\big(\Fst(Ax),\Fst(Bx)\big)$.  This, in turn, means a strictly unital symmetric monoidal functor $\Fst(Ax) \to \Fst(Bx)$.  Using \cref{Fdgrpsix-partner} in the current context, we obtain the $x$-component
\[\big(\Fstdgr\psi\big)_x = \Fst(\psi_x) \cn \Fst(Ax) \to \Fst(Bx)\]
by applying $\Fst$ to $\psi_x$ in \cref{psi-xcomponent-AxBx}.  By \cref{def:Fst-onecells} $\Fst(\psi_x)$ is a \emph{strict} symmetric monoidal functor.
\end{explanation}

\section[$\Mone$-Modules and Permutative Enriched Diagrams]{Homotopy Equivalent \texorpdfstring{$\Mone$}{M1}-Modules and Permutative Enriched Diagrams}
\label{sec:mackey-mone}

In this section we apply \cref{mackey-gen-xiv,mackey-xiv-cor} to show that the categories of enriched diagrams and Mackey functors in left $\Mone$-modules and permutative categories are connected by inverse equivalences of homotopy theories.  See \cref{mackey-xiv-mone}.  \cref{expl:mackey-mone} summarizes \cref{mackey-xiv-pmulticat,mackey-xiv-mone} in terms of left modules in $\MoneMod$, $\pMulticat$, and $\permcatsu$.  We explain these equivalences of homotopy theories further in \cref{sec:Endmdg-Fmdgr}.

\subsection*{Context}

For the context first recall the diagram 
\begin{equation}\label{mackey-mone-setting}
\begin{tikzpicture}[baseline={(a.base)}]
\draw[0cell]
(0,0) node (a') {\MoneMod}
(a')+(.5,-.06) node (a) {}
(a)+(2,0) node (b) {}
(b)+(.8,.05) node (b') {\permcatsu}
;
\draw[1cell=.9]
(a) edge[transform canvas={yshift=.6ex}] node {\Fm} (b)
(b) edge[transform canvas={yshift=-.4ex}] node {\Endm} (a)
;
\end{tikzpicture}
\end{equation}
in \cref{FmEmcontext} consisting of
\begin{itemize}
\item the $\Cat$-multicategory $\MoneMod$ in \cref{expl:monemodcatmulticat},
\item the $\Cat$-multicategory $\permcatsu$ in \cref{thm:permcatmulticat},
\item the $\Cat$-multifunctor $\Endm$ in \cref{expl:endm-catmulti}, and
\item the non-symmetric $\Cat$-multifunctor $\Fm = \Fst\Um$ in \cref{Fm-multi-def}.
\end{itemize}
The two composites in \cref{mackey-mone-setting} are connected to the respective identity functors via the following $\Cat$-multinatural transformations from \cref{def:etam-multi,def:vrhom-multi}.
\begin{equation}\label{mackey-mone-setting-ii}
\begin{tikzpicture}[baseline={(a.base)}]
\def\t{30} \def\s{27}
\draw[0cell]
(0,0) node (a) {\phantom{X}}
(a)+(-.45,0) node (a') {\MoneMod}
(a)+(2,0) node (b) {\phantom{X}}
(b)+(.45,.0) node (b') {\MoneMod}
;
\draw[1cell=.9]
(a) edge[bend left=\t] node {1} (b)
(a) edge[bend right=\t] node[swap] {\Endm\Fm} (b)
;
\draw[2cell]
node[between=a and b at .38, rotate=-90, 2label={above,\etam}] {\Rightarrow}
;
\begin{scope}[shift={(5,0)}]
\draw[0cell]
(0,0) node (a) {\phantom{X}}
(a)+(-.65,0) node (a') {\permcatsu}
(a)+(2,0) node (b) {\phantom{X}}
(b)+(.65,.0) node (b') {\permcatsu}
;
\draw[1cell=.9]
(a) edge[bend left=\t] node {1} (b)
(a) edge[bend right=\t] node[swap] {\Fm\Endm} (b)
;
\draw[2cell]
node[between=a and b at .37, rotate=-90, 2label={above,\vrhom}] {\Rightarrow}
;
\end{scope}
\end{tikzpicture}
\end{equation}
The underlying categories of $\MoneMod$ and $\permcatsu$ are equipped with the relative category structures
\begin{equation}\label{MonemodPerm-steq}
\brb{\MoneMod,\cSM} \andspace 
\brb{\permcatsu,\cS}
\end{equation}
in \cref{SM-Sst-def} and \cref{perm-steq}.
\begin{itemize}
\item The wide subcategory of stable equivalences 
\[\cS \bigsubset \permcatsu\]
is created by Segal $K$-theory $\Kse$ \cref{Kse}.  For a small $\permcatsu$-category $\C$, the wide subcategory in \cref{mackey-xiv-mone-functors} below
\begin{equation}\label{sStri-def}
\cStri \bigsubset \permcatsucat\big(\C,\permcatsu\big)
\end{equation}
is defined as in \cref{cWtri-subcat} using $\cS$. 
\item The wide subcategory of $\Fm$-stable equivalences
\[\cSM = \Fm^\inv(\cS) \bigsubset \MoneMod\] 
is created by $\Fm$.  The wide subcategory in \cref{mackey-xiv-mone-functors} below
\begin{equation}\label{sSmtri-def}
\cSmtri \bigsubset \Monemodcat\big(\C_{\Endm},\MoneMod\big)
\end{equation}
is defined as in \cref{cWtri-subcat} using $\cSM$. 
\end{itemize}

\subsection*{Equivalences of Homotopy Theories}

In \cref{ptmulticat-xxv} we observe that the pair $(\Fm,\Endm)$ induces inverse equivalences of homotopy theories between the respective categories of non-symmetric $\Q$-algebras for each small non-symmetric $\Cat$-multicategories $\Q$.  The following observation extends the inverse equivalences of homotopy theories $(\Fm,\Endm)$ to categories of enriched diagrams and Mackey functors.  It is the $\MoneMod$ analog of \cref{mackey-xiv-pmulticat}.

\begin{theorem}\label{mackey-xiv-mone}
Suppose $\C$ is a small $\permcatsu$-category.  Then the functors
\begin{equation}\label{mackey-xiv-mone-functors}
\begin{tikzpicture}[baseline={(a.base)}]
\draw[0cell=.8]
(0,0) node (a) {\Big(\Monemodcat\big(\C_{\Endm},\MoneMod\big) \scs \cSmtri\Big)}
(a)+(6,0) node (b) {\Big(\permcatsucat\big(\C,\permcatsu\big) \scs \cStri\Big),} 
;
\draw[1cell=.8]
(a) edge[transform canvas={yshift=.6ex}] node (f) {\Fmdgr} node[swap] {\sim} (b)
(b) edge[transform canvas={yshift=-.5ex}] node (e) {\Endmdg} (a)
;
\end{tikzpicture}
\end{equation}
defined by the data in \cref{mackey-mone-setting,mackey-mone-setting-ii,MonemodPerm-steq,sStri-def,sSmtri-def}, are inverse equivalences of homotopy theories.  

Moreover, the variant with $(\C_{\Endm})^\op$ and $\Cop$ replacing, respectively, $\C_{\Endm}$ and $\C$ is also true.
\end{theorem}

\begin{proof}
The first assertion is an instance of \cref{mackey-gen-xiv}, which is applicable in the current setting as we now explain.  Following the summary in \cref{expl:mackey-xiv-assumptions}, first we verify that \cref{def:mackey-gen-context} \cref{mackey-context-i,mackey-context-ii,mackey-context-iii,mackey-context-iv,mackey-context-v} are satisfied in the current context.
\begin{romenumerate}
\item\label{mackey-xiv-mone-i} $\M = \MoneMod$ is a closed multicategory by
\begin{itemize}
\item \cref{smclosed-closed-multicat} and
\item the fact that it is a symmetric monoidal closed category (\cref{proposition:EM2-5-1} \cref{monebicomplete}).
\end{itemize}  
By \cref{perm-closed-multicat}, $\N = \permcatsu$ is a closed multicategory.
\item\label{mackey-xiv-mone-ii} $\C$ is, by assumption, a small $\permcatsu$-category.
\item\label{mackey-xiv-mone-iii} $F = \Fm$ in \cref{mackey-mone-setting} is a non-symmetric multifunctor by definition \cref{Fm-multi-def}, and $E = \Endm$ is a multifunctor by \cref{endfactorization}.
\item\label{mackey-xiv-mone-iv} $\uni = \etam$ and $\cou = \vrhom$ in \cref{mackey-mone-setting-ii} are multinatural transformations by
\begin{itemize}
\item \cref{expl:etam-multi} for $\etam$ and
\item \cref{ptmulticat-xxii,def:vrhom-multi} for $\vrhom$.
\end{itemize} 
\item\label{mackey-xiv-mone-v} In the current setting, the condition \cref{counit-compatibility} is the equality of the following two left $\Mone$-module morphisms for each pair of objects $x,y \in \C$. 
\[\begin{tikzpicture}
\draw[0cell=1]
(0,0) node (a) {\Endm \C(x,y)}
(a)+(5.5,0) node (b) {\Endm \Fm \Endm \C(x,y)}
;
\draw[1cell=.9]
(a) edge[transform canvas={yshift=.6ex}] node {\etam_{\Endm \C(x,y)}} (b)
(a) edge[transform canvas={yshift=-.4ex}] node[swap] {\Endm \vrhom_{\C(x,y)}} (b)
;
\end{tikzpicture}\]
This equality holds by \cref{etavrhoMone} because each hom object $\C(x,y)$ is a small permutative category.
\end{romenumerate}
Thus \cref{def:mackey-gen-context} \cref{mackey-context-i,mackey-context-ii,mackey-context-iii,mackey-context-iv,mackey-context-v} hold in the context of \cref{mackey-mone-setting,mackey-mone-setting-ii,MonemodPerm-steq,sStri-def,sSmtri-def}.

Next, the only assumption in \cref{def:mackey-heq-context} is that the relative category
\[(\N, \cX) = \brb{\permcatsu,\cS}\]
is a category with weak equivalences (\cref{definition:rel-cat} \eqref{relcat-vi}).  This is true for the following two reasons.
\begin{itemize}
\item The wide subcategory $\cS \bigsubset \permcatsu$ is created by a functor, namely, Segal $K$-theory \cref{perm-steq}
\[\Kse \cn \permcatsu \to \Spc.\]
\item The class of stable equivalences in $\Spc$ contains all the isomorphisms and has the 2-out-of-3 property.
\end{itemize} 
In the current setting there are equalities of wide subcategories
\[F^\inv\cX = \Fm^\inv(\cS) = \cSM \bigsubset \MoneMod.\]
The data in \cref{mackey-xiv-mone-functors} are those in \cref{mackey-xiv-functors} in the current context.

Finally, the components of $\etam$ are $\Fm$-stable equivalences, as explained in the proof of \cref{ptmulticat-thm-xi}.  Moreover, by definition \cref{vrhom-multi-def} each component of $\vrhom$ is a component of $\vrhost$.  The latter is a stable equivalence in $\permcatsu$ by \cref{remark:steq} \cref{it:steq-2} because it has a left adjoint by \cref{ptmulticat-prop-viii}.  Thus \cref{mackey-gen-xiv} is applicable in the current setting, proving the first assertion.

The second assertion about $(\C_{\Endm})^\op$ and $\Cop$ is an instance of \cref{mackey-xiv-cor}.  It is applicable because $\Endm$ is a multifunctor (\cref{endfactorization}).
\end{proof}

\begin{explanation}[Homotopy Equivalent Categories of Modules]\label{expl:mackey-mone}
The categories of enriched diagrams in \cref{mackey-xiv-pmult-functors,mackey-xiv-mone-functors} are categories of left modules by \cref{C-diagram-partner,C-diag-morphism-pn}.  Thus, \cref{mackey-xiv-pmulticat,mackey-xiv-mone} together assert that, for each small $\permcatsu$-category $\C$, the functors in the diagram
\begin{equation}\label{left-C-modules}
\begin{tikzpicture}[vcenter]
\def\h{0} \def\v{-1.3}
\draw[0cell=.9]
(0,0) node (a) {\Big(\Monemodcat\big(\C_{\Endm} \scs \MoneMod\big) \scs \cSmtri\Big)}
(a)+(\h,\v) node (b) {\Big(\permcatsucat\big(\C \scs \permcatsu\big) \scs \cStri\Big)}
(a)+(0,2*\v) node (c) {\Big(\pMulticatcat\big(\C_{\Endst} \scs \pMulticat\big) \scs \cSsttri\Big)}
;
\draw[1cell=.8]
(a) edge[transform canvas={xshift=1ex}] node[pos=.5] {\Fmdgr} (b)
(b) edge[transform canvas={xshift=-1ex}] node[pos=.5] {\Endmdg} (a)
(c) edge[transform canvas={xshift=1ex}] node[swap,pos=.5] {\Fstdgr} (b)
(b) edge[transform canvas={xshift=-1ex}] node[swap,pos=.6] {\Endstdg} (c)
;
\draw[2cell]
node[between=a and b at .5, rotate=-90] {\sim}
node[between=b and c at .5, rotate=-90] {\sim}
;
\end{tikzpicture}
\end{equation}
are equivalences of homotopy theories between
\begin{itemize}
\item left $\C$-modules in $\permcatsu$,
\item left $\C_{\Endm}$-modules in $\MoneMod$, and
\item left $\C_{\Endst}$-modules in $\pMulticat$.
\end{itemize}
Furthermore, these equivalences of homotopy theories still hold if $\C$, $\C_{\Endm}$, and $\C_{\Endst}$ are replaced by $\Cop$, $(\C_{\Endm})^\op$, and $(\C_{\Endst})^\op$, respectively.
\end{explanation}

\section{Explanation of the Equivalences of Homotopy Theories}
\label{sec:Endmdg-Fmdgr}

In this section we
\begin{enumerate}
\item explain in detail the inverse equivalences of homotopy theories $\Endmdg$ and $\Fmdgr$ in \cref{left-C-modules} and
\item compare them with the inverse equivalences of homotopy theories $\Endstdg$ and $\Fstdgr$.
\end{enumerate} 
This section is organized as follows.
\begin{itemize}
\item \cref{expl:endmdg} describes how $\Endstdg$ factors through $\Endmdg$.
\item \cref{expl:Endmdg-objects,expl:endmdg-morphisms} describe $\Endmdg$ on objects and morphisms, respectively.
\item \cref{expl:Fmdgr} describes how $\Fmdgr$ factors through $\Fstdgr$.
\item \cref{expl:Fmdgr-objects,expl:Fmdgr-morphisms} describe $\Fmdgr$ on objects and morphisms, respectively.
\end{itemize}
Throughout this section we use the shortened notation in \cref{pMpsu}, so
\[\pM = \pMulticat \andspace \psu = \permcatsu.\]

\subsection*{The Functor $\Endmdg$}

\begin{explanation}[Factoring $\Endstdg$ Through $\Endmdg$]\label{expl:endmdg}
By \cref{endufactor} there is a factorization of multifunctors
\begin{equation}\label{Endst-factorization}
\Endst \cn \psu \fto{\Endm} \MoneMod \fto{\Um} \pM.
\end{equation}
By \cref{gspectra-thm-vii} this factorization yields the following commutative diagram.
\begin{equation}\label{EndstUmEndmdg}
\begin{tikzpicture}[vcenter]
\def\h{3} \def\v{-1.3}
\draw[0cell=.9]
(0,0) node (a) {\Monemodcat\big(\C_{\Endm} \scs \MoneMod\big)}
(a)+(\h,\v) node (b) {\psucat\big(\C \scs \psu\big)}
(a)+(0,2*\v) node (c) {\pMcat\big(\C_{\Endst} \scs \pM\big)}
;
\draw[1cell=.9]
(b) edge[transform canvas={xshift=-1em}] node[swap,pos=.2] {\Endmdg} (a)
(a) edge[transform canvas={xshift=0em}] node[swap,pos=.5] {\Umdg} (c)
(b) edge[transform canvas={xshift=-1em}] node[pos=.2] {\Endstdg} (c)
;
\end{tikzpicture}
\end{equation}
In \cref{EndstUmEndmdg} each arrow is a diagram change-of-enrichment functor (\cref{gspectra-thm-v}).  We describe $\Endmdg$ on objects and morphisms in \cref{expl:Endmdg-objects,expl:endmdg-morphisms} below.
\end{explanation}

\begin{explanation}[$\Endmdg$ on Objects]\label{expl:Endmdg-objects}
The functor $\Endmdg$ in \cref{EndstUmEndmdg} has an analogous description as $\Endstdg$ in \cref{expl:Endstdg-ii}.  More explicitly, the functor $\Endmdg$ sends a $\psu$-functor $A \cn \C \to \psu$ to the following composite $\MoneMod$-functor.
\begin{equation}\label{Endmdg-A-composite}
\begin{tikzpicture}[baseline={(a.base)}]
\def\h{3} \def\w{.7}
\draw[0cell=1]
(0,0) node (a) {\C_{\Endm}}
(a)+(\h,0) node (b) {(\psu)_{\Endm}}
(b)+(\h,0) node (c) {\phantom{\MoneMod}}
(c)+(0,.1) node (c') {\MoneMod}
;
\draw[1cell=.8]
(a) edge node {A_{\Endm}} (b)
(b) edge node {\Endmse} (c)
;
\draw[1cell=.9]
(a) [rounded corners=3pt,shorten <=-.3ex] |- ($(b)+(-1,\w)$)
-- node {\Endmdg A} ($(b)+(1,\w)$) -| (c')
;
\end{tikzpicture}
\end{equation}
\begin{itemize}
\item $A_{\Endm}$ is the image of $A$ under the change-of-enrichment 2-functor
\[\dEndm \cn \psucat \to \Monemodcat\]
along $\Endm$ (\cref{expl:endm-catmulti,mult-change-enrichment}).
\item $\Endmse$ is the standard enrichment of $\Endm$ (\cref{gspectra-thm-iii}).
\end{itemize}
The $\MoneMod$-functor $\Endmdg A$ in \cref{Endmdg-A-composite} sends each object $x \in \C$ to 
\[\big(\Endmdg A\big)x = \Endm(Ax) \inspace \MoneMod.\]

For objects $x,y \in \C$, the $(x,y)$-component left $\Mone$-module morphism of $\Endmdg A$ is the following composite, where $\pHom$ is the internal hom in $\MoneMod$ (\cref{proposition:EM2-5-1} \cref{monebicomplete}).
\[\begin{tikzpicture}
\def\h{2.8} \def\v{1.1} \def\t{15}
\draw[0cell=.9]
(0,0) node (a) {\Endm\C(x,y)}
(a)+(\h,-\v) node (b) {\Endm\psu(Ax,Ay)}
(b)+(\h+.5,\v) node (c) {\pHom\big(\Endm(Ax) \scs \Endm(Ay)\big)}
(b)+(\h,\v) node (c') {\phantom{M}}
;
\draw[1cell=.9]
(a) edge node {\big(\Endmdg A\big)_{x,y}} (c)
(a) edge[bend right=\t] node[swap,pos=.3] {\Endm(A_{x,y})} (b)
(b) edge[bend right=\t,shorten >=1ex] node[swap,pos=.7] {(\Endmse)_{Ax,Ay}} (c')
;
\end{tikzpicture}\]
The underlying pointed multifunctor of this left $\Mone$-module morphism is equal to the one in \cref{EndstdgA-xy-comp} by \cref{proposition:EM2-5-1} \cref{monebicomplete} and the factorization \cref{Endst-factorization}.
\end{explanation}

\begin{explanation}[$\Endmdg$ on Morphisms]\label{expl:endmdg-morphisms}
For a $\psu$-natural transformation $\psi \cn A \to B$ as in \cref{Endstdg-def}, the $\MoneMod$-natural transformation $\Endmdg\psi$ is the whiskering below.
\[\begin{tikzpicture}[baseline={(a.base)}]
\def\t{25} \def\d{.8}
\draw[0cell=1]
(0,0) node (a') {\C_{\Endm}}
(a')+(.4,0) node (a'') {\phantom{C}}
(a'')+(2.5,0) node (b2) {\phantom{C}}
(b2)+(.7,0) node (b') {(\psu)_{\Endm}}
(b')+(3,0) node (c) {\phantom{\MoneMod}}
(c)+(0,.1) node (c') {\MoneMod}
;
\draw[1cell=.8]
(a'') edge[bend left=\t] node {A_{\Endm}} (b2)
(a'') edge[bend right=\t] node[swap] {B_{\Endm}} (b2)
(b') edge[transform canvas={yshift=.5ex}] node {\Endmse} (c)
;
\draw[2cell=.9]
node[between=a'' and b2 at .32, rotate=-90, 2label={above,\psi_{\Endm}}] {\Rightarrow}
;
\end{tikzpicture}\]
For each object $x \in \C$, its $x$-component is given by the left $\Mone$-module morphism
\begin{equation}\label{Endm-psi-x}
\big(\Endmdg \psi\big)_x = \Endm(\psi_x) \cn \Endm(Ax) \to \Endm(Bx).
\end{equation}
This is obtained from the strictly unital symmetric monoidal functor in \cref{psix-Ax-to-Bx}
\[\psi_x \cn Ax \to Bx\]
by applying $\Endm$.  The underlying pointed multifunctor of $\Endm(\psi_x)$ in \cref{Endm-psi-x} is equal to $\Endst(\psi_x)$ in \cref{Endst-psi-x} by the factorization \cref{Endst-factorization}.
\end{explanation}

\subsection*{The Functor $\Fmdgr$}

\begin{explanation}[Factoring $\Fmdgr$ Through $\Fstdgr$]\label{expl:Fmdgr}
By definition \cref{Fdgr-def} the functor $\Fmdgr$ is the composite along the top of the following diagram.
\begin{equation}\label{Fmdgr-factorization}
\begin{tikzpicture}[baseline={(a.base)}]
\def\h{3} \def\v{1.1} \def\w{.5}
\draw[0cell=.85]
(0,0) node (a) {\Monemodcat\big(\C_{\Endm} \scs \MoneMod\big)}
(a)+(-.5,0) node (a') {\phantom{\MoneMod}}
(a)+(\h,\v) node (b) {\psucat\big(\C_{\Fm\Endm} \scs \psu\big)}
(b)+(\h,-\v) node (c) {\psucat(\C,\psu)}
(c)+(.4,0) node (c') {\phantom{\psu}}
(a)+(\h,-\v) node (d) {\pMcat\big(\C_{\Endst} \scs \pM\big)}
;
\draw[1cell=.8]
(a) edge[shorten <=-1ex] node[pos=.3] {\Fmdg} (b)
(b) edge node[pos=.7] {\C^*_{\vrhom}} (c)
(a) edge[shorten <=-1ex] node[swap,pos=.3] {\Umdg} (d)
(d) edge[transform canvas={xshift=0em}] node[swap] {\Fstdg} (b)
(d) edge node[swap,pos=.7] {\Fstdgr} (c)
;
\draw[1cell=.85]
(a') [rounded corners=3pt,shorten <=1ex,shorten >=1ex] |- ($(b)+(-1,\w)$)
-- node {\Fmdgr} ($(b)+(1,\w)$) -| (c')
;
\end{tikzpicture}
\end{equation}
The two triangles in \cref{Fmdgr-factorization} are commutative for the following reasons.
\begin{enumerate}
\item The diagram change of enrichment $\Umdg$ is well defined by the factorization \cref{Endst-factorization}
\[\Endst = \Um \Endm\]
and the functoriality of change of enrichment with respect to composition (\cref{func-change-enr}).  The diagram change of enrichment $\Fstdg$ is well defined by the equalities
\[\Fm\Endm = \Fst\Um\Endm = \Fst\Endst,\]
with the first equality given by the definition \cref{Fm-multi-def} of $\Fm$ as the composite
\[\Fm \cn \MoneMod \fto{\Um} \pM \fto{\Fst} \psu.\]
The left triangle in \cref{Fmdgr-factorization} commutes by the functoriality of diagram change of enrichment with respect to composition (\cref{gspectra-thm-vii}).
\item The $\psu$-functor
\begin{equation}\label{C-vrhom}
\C_{\vrhom} \cn \C \to \C_{\Fm\Endm}
\end{equation}
is the $\C$-component of the 2-natural transformation (\cref{dtheta-twonat})
\[\begin{tikzpicture}[baseline={(a.base)}]
\def\s{22}
\draw[0cell]
(0,0) node (a) {\phantom{\M}}
(a)+(3,0) node (b) {\phantom{\N}}
(a)+(-.4,0) node (a') {\psucat}
(b)+(.4,0) node (b') {\psucat}
;
\draw[1cell=.85]
(a) edge[bend left=\s] node {1} (b)
(a) edge[bend right=\s] node[swap] {\dFmEndm = \dFm \,\dEndm} (b)
;
\draw[2cell]
node[between=a and b at .35, rotate=-90, 2label={above,\dvrhom}] {\Rightarrow}
;
\end{tikzpicture}\]
induced by the multinatural transformation \cref{vrhom-multi-def}
\[\vrhom \cn 1_{\psu} \to \Fm\Endm = \Fst\Endst.\]
There is an equality of $\psu$-functors
\begin{equation}\label{Cvrhom-Cvrhost}
\C_{\vrhom} = \C_{\vrhost} \cn \C \to \C_{\Fm\Endm} = \C_{\Fst\Endst},
\end{equation}
with $\C_{\vrhost}$ the $\psu$-functor in \cref{C-vrhost}, for the following two reasons.
\begin{romenumerate}
\item Both $\C_{\vrhom}$ and $\C_{\vrhost}$ are the identity function on objects.
\item For objects $x,y \in \C$, the $(x,y)$-component strictly unital symmetric monoidal functors
\[\begin{tikzpicture}
\draw[0cell]
(0,0) node (a) {\C(x,y)}
(a)+(5,0) node (b) {\Fm\Endm\C(x,y) = \Fst\Endst\C(x,y)}
;
\draw[1cell=.9]
(a) edge[transform canvas={yshift=.6ex}] node {\vrhom_{\C(x,y)}} (b)
(a) edge[transform canvas={yshift=-.4ex}] node[swap] {\vrhost_{\C(x,y)}} (b)
;
\end{tikzpicture}\]
are equal by the definition \cref{vrhom-multi-def} of $\vrhom$.
\end{romenumerate}
The equality \cref{Cvrhom-Cvrhost} implies the equality of functors
\[\C_{\vrhom}^* = \C_{\vrhost}^* \cn \psucat\big(\C_{\Fm\Endm} \scs \psu\big) \to \psucat(\C,\psu).\]
Thus the right triangle in \cref{Fmdgr-factorization} commutes by the definition \cref{Fstdgr-definition} of $\Fstdgr$.
\end{enumerate}
In summary, the diagram \cref{Fmdgr-factorization} factors $\Fmdgr$ through $\Fstdgr$.
\end{explanation}

\begin{explanation}[$\Fmdgr$ on Objects]\label{expl:Fmdgr-objects}
The functor $\Fmdgr$ in \cref{Fmdgr-factorization} has an analogous description as $\Fstdgr$ in \cref{expl:Fstdgr-object,expl:Fstdgr-morphism} by replacing $\Endst$, $\Fst$, $\vrhost$, and $\pM$ with $\Endm$, $\Fm$, $\vrhom$, and $\MoneMod$, respectively.  More explicitly, $\Fmdgr$ sends a $\MoneMod$-functor
\begin{equation}\label{monemod-functor-A}
A \cn \C_{\Endm} \to \MoneMod
\end{equation}
to the following composite $\psu$-functor.
\begin{equation}\label{Fmdgr-A-composite}
\begin{tikzpicture}[vcenter]
\def\g{.3} \def\v{1.3} \def\w{.6}
\draw[0cell=.9]
(0,0) node (a) {\C} 
(a)+(\g,-\v) node (b) {\C_{\Fm\Endm}}
(b)+(2.2,0) node (b') {(\C_{\Endm})_{\Fm}}
(b')+(3.5,0) node (c) {(\MoneMod)_{\Fm}}
(c)+(\g,\v) node (d) {\psu}
;
\draw[1cell=2]
(b) edge[-,double equal sign distance] (b')
;
\draw[1cell=.9]
(a) edge node {\Fmdgr A} (d)
(a) edge node[pos=.2,swap] {\C_{\vrhom}} (b)
(b') edge node {A_{\Fm}} (c)
(c) edge node[pos=.8,swap] {\Fmse} (d)
;
\end{tikzpicture}
\end{equation}
\begin{itemize}
\item $\C_{\vrhom}$ is the $\psu$-functor in \cref{C-vrhom}.
\item $A_{\Fm}$ is the image of $A$ under the change-of-enrichment 2-functor 
\begin{equation}\label{dFm}
\dFm \cn \Monemodcat \to \psucat
\end{equation}
along $\Fm$ (\cref{ptmulticat-thm-ii,mult-change-enrichment}).
\item $\Fmse$ is the standard enrichment $\psu$-functor of $\Fm$ (\cref{gspectra-thm-iii}).
\end{itemize}
The $\psu$-functor in \cref{Fmdgr-A-composite} sends each object $x \in \C$ to
\begin{equation}\label{Fm-Ax}
\big(\Fmdgr A\big)x = \Fm(Ax) = \Fst\Um(Ax) \inspace \psu,
\end{equation}
with the second equality from the definition \cref{Fm-multi-def} of $\Fm$.

For objects $x,y \in \C$, the $(x,y)$-component of the $\psu$-functor in \cref{Fmdgr-A-composite} is the following composite strictly unital symmetric monoidal functor.
\begin{equation}\label{Fmdgr-Axy}

\end{equation}
By the factorization \cref{Fmdgr-factorization}, there is an equality of strictly unital symmetric monoidal functors
\[\big(\Fmdgr A\big)_{x,y} = \big(\Fstdgr (\Umdg A)\big)_{x,y}\]
with the left-hand side from \cref{Fmdgr-Axy} and the right-hand side from \cref{FstdgrA-component}.
\end{explanation}

\begin{explanation}[$\Fmdgr$ on Morphisms]\label{expl:Fmdgr-morphisms}
Consider a $\MoneMod$-natural transformation $\psi$ (\cref{def:enriched-natural-transformation}) as in the left diagram below.
\[%
\]
The $\psu$-natural transformation $\Fmdgr\psi$, as in the right diagram above, is the following whiskering.
\begin{equation}\label{Fmdgrpsi-def}
%
\end{equation}
The whiskering \cref{Fmdgrpsi-def} is the one in \cref{Fdgrpsi-def} in the current context and is an analog of \cref{Fstdgrpsi-def}.
\begin{itemize}
\item The $\psu$-functors $\C_{\vrhom}$ and $\Fmse$ are as in \cref{Fmdgr-A-composite}.
\item $\psi_{\Fm}$ is the change of enrichment of $\psi$ under $\dFm$ in \cref{dFm}.
\end{itemize}

Next we describe the components of $\Fmdgr\psi$ explicitly.  For each object $x \in \C$, the $x$-component of $\psi$ is a left $\Mone$-module morphism
\[\psi_x \cn \Mone \to \pHom(Ax,Bx).\]
By the $\sma$-$\pHom$ adjunction in $\MoneMod$ (\cref{proposition:EM2-5-1} \cref{monebicomplete}),
$\psi_x$ is uniquely determined by its adjoint, which is a left $\Mone$-module morphism that we also denote by
\begin{equation}\label{psi-AxBx-mone}
\psi_x \cn Ax \to Bx.
\end{equation}
On the other hand, by \cref{Fm-Ax} the $x$-component of $\Fmdgr\psi$ is a nullary multimorphism
\[\big(\Fmdgr\psi\big)_x \cn \ang{} \to \psu\big(\Fm(Ax),\Fm(Bx)\big) \inspace \psu.\]
This means a choice of an object in $\psu\big(\Fm(Ax),\Fm(Bx)\big)$.  This, in turn, means a strictly unital symmetric monoidal functor $\Fm(Ax) \to \Fm(Bx)$.  Using \cref{Fdgrpsix-partner} in the current context, we obtain the $x$-component
\[\begin{split}
\big(\Fmdgr\psi\big)_x &= \Fm(\psi_x)\\
&= \Fst\Um(\psi_x) \cn \Fm(Ax) \to \Fm(Bx)
\end{split}\]
by applying $\Fm = \Fst\Um$ to $\psi_x$ in \cref{psi-AxBx-mone}.  By \cref{def:Fst-onecells} $\Fst\Um(\psi_x)$ is a \emph{strict} symmetric monoidal functor.
\end{explanation}

\section[Multicategorical and $\Mone$-Modules Enriched Diagrams]{Homotopy Equivalent Multicategorical and \texorpdfstring{$\Mone$}{M1}-Modules Enriched Diagrams}
\label{sec:mackey-pmult-mone}

As we discuss in \cref{left-C-modules}, enriched diagrams in $\MoneMod$ and $\pMulticat$ are connected by two zigzags of equivalences of homotopy theories:
\[\big(\Fmdgr \scs \Fstdgr\big) \andspace \big(\Endmdg \scs \Endstdg\big).\]
Each of these two zigzags goes through enriched diagrams in $\permcatsu$.  In this section we apply \cref{mackey-gen-xiv,mackey-xiv-cor} to show that the categories of enriched diagrams and Mackey functors in $\MoneMod$ and $\pMulticat$ are directly connected by inverse equivalences of homotopy theories.  See \cref{mackey-pmulti-mone} and the summary in \cref{expl:left-mod-heq}.  We explain these functors further in \cref{sec:Monesmadgr-Umdg}.

\subsection*{Context}

For the context first recall the diagram 
\begin{equation}\label{mackey-pmult-mone-setting}
\begin{tikzpicture}[baseline={(a.base)}]
\draw[0cell]
(0,0) node (a') {\pMulticat}
(a')+(.7,0) node (a) {}
(a)+(2,0) node (b) {}
(b)+(.55,.06) node (b') {\MoneMod}
;
\draw[1cell=.9]
(a) edge[transform canvas={yshift=.6ex}] node {\Monesma} (b)
(b) edge[transform canvas={yshift=-.4ex}] node {\Um} (a)
;
\end{tikzpicture}
\end{equation}
consisting of
\begin{itemize}
\item the $\Cat$-multicategory $\pMulticat$ in \cref{expl:ptmulticatcatmulticat},
\item the $\Cat$-multicategory $\MoneMod$ in \cref{expl:monemodcatmulticat},
\item the $\Cat$-multifunctor $\Monesma$ in \cref{Monesma-CatSM}.
\item the $\Cat$-multifunctor $\Um$ in \cref{expl:Um-catmulti}, and
\end{itemize}
The two composite functors in \cref{mackey-pmult-mone-setting} are connected to the respective identity functors via the following $\Cat$-multinatural transformations from \cref{etahat-epzhat-monCatnat}, where $\epzhatinv$ denotes the inverse of $\epzhat$.
\begin{equation}\label{mackey-pmult-mone-setting-ii}
\begin{tikzpicture}[baseline={(a.base)}]
\def\t{25} \def\s{25}
\draw[0cell]
(0,0) node (a) {\phantom{X}}
(a)+(2,0) node (b) {\phantom{X}}
(a)+(-.53,0) node (a') {\pMulticat}
(b)+(.55,.0) node (b') {\pMulticat}
;
\draw[1cell=.9]
(a) edge[bend left=\t] node {1} (b);
\draw[1cell=.8]
(a) edge[bend right=\t] node[swap] {\Um(\Monesma)} (b)
;
\draw[2cell=1]
node[between=a and b at .42, rotate=-90, 2label={above,\etahat}] {\Rightarrow}
;
\begin{scope}[shift={(4.8,0.05)}]
\draw[0cell]
(0,0) node (a) {\phantom{X}}
(a)+(2,0) node (b) {\phantom{X}}
(a)+(-.43,0) node (a') {\MoneMod}
(b)+(.43,.0) node (b') {\MoneMod}
;
\draw[1cell=.9]
(a) edge[bend left=\t] node {1} (b);
\draw[1cell=.8]
(a) edge[bend right=\t] node[swap] {(\Monesma)\Um} (b)
;
\draw[2cell=1]
node[between=a and b at .43, rotate=-90, 2label={above,\epzhatinv}] {\Rightarrow}
;
\end{scope}
\end{tikzpicture}
\end{equation}
The underlying categories of $\pMulticat$ and $\MoneMod$ are equipped with the relative category structures in \cref{SM-Sst-def}:
\begin{equation}\label{pMultMone-steq}
\brb{\pMulticat,\cSst} \andspace 
\brb{\MoneMod,\cSM}.
\end{equation}
By \cref{ptmulticat-thm-vi} the functors in \cref{mackey-pmult-mone-setting} are inverse equivalences of homotopy theories.  The algebra version is \cref{Monesma-Um-algebra}. 
\begin{itemize}
\item The wide subcategory of $\Fst$-stable equivalences
\[\cSst = \Fst^\inv(\cS) \bigsubset \pMulticat\] 
is created by the functor
\[\Fst \cn \pMulticat \to \permcatsu.\]
The wide subcategory $\cS \bigsubset \permcatsu$ is created by Segal $K$-theory $\Kse$; see \cref{perm-steq}.  For a small $\MoneMod$-category $\D$, the wide subcategory in \cref{mackey-pmult-mone-functors} below
\begin{equation}\label{sSsttri-defin}
\cSsttri \bigsubset \pMulticatcat\big(\D_{\Um},\pMulticat\big)
\end{equation}
is defined as in \cref{cWtri-subcat} using $\cSst$. 
\item The wide subcategory of $\Fm$-stable equivalences
\[\cSM = \Fm^\inv(\cS) \bigsubset \MoneMod\] 
is created by the functor
\[\Fm \cn \MoneMod \to \permcatsu.\]  
The wide subcategory in \cref{mackey-pmult-mone-functors} below
\begin{equation}\label{sSmtri-defin}
\cSmtri \bigsubset \Monemodcat\big(\D,\MoneMod\big)
\end{equation}
is defined as in \cref{cWtri-subcat} using $\cSM$. 
\end{itemize}

\subsection*{Equivalences of Homotopy Theories}

In \cref{ptmulticat-thm-vi} we observe that the pair $(\Monesma,\Um)$ is an adjoint equivalence of homotopy theories.  The following observation extends these equivalences of homotopy theories to categories of enriched diagrams and Mackey functors.

\begin{theorem}\label{mackey-pmulti-mone}
Suppose $\D$ is a small $\MoneMod$-category.  Then the functors
\begin{equation}\label{mackey-pmult-mone-functors}
\begin{tikzpicture}[baseline={(a.base)}]
\draw[0cell=.8]
(0,0) node (a) {\Big(\pMulticatcat\big(\D_{\Um},\pMulticat\big) \scs \cSsttri\Big)}
(a)+(6,0) node (b) {\Big(\Monemodcat\big(\D,\MoneMod\big) \scs \cSmtri\Big),} 
;
\draw[1cell=.7]
(a) edge[transform canvas={yshift=.6ex}] node {\Monesmadgr} node[swap] {\sim} (b)
;
\draw[1cell=.8]
(b) edge[transform canvas={yshift=-.5ex}] node {\Umdg} (a)
;
\end{tikzpicture}
\end{equation}
defined by the data in \cref{mackey-pmult-mone-setting,mackey-pmult-mone-setting-ii,pMultMone-steq,sSsttri-defin,sSmtri-defin}, are inverse equivalences of homotopy theories.  

Moreover, the variant with $(\D_{\Um})^\op$ and $\D^\op$ replacing, respectively, $\D_{\Um}$ and $\D$ is also true.
\end{theorem}

\begin{proof}
The first assertion is an instance of \cref{mackey-gen-xiv}, which is applicable in the current setting as we now explain.  Following the summary in \cref{expl:mackey-xiv-assumptions}, first we verify that \cref{def:mackey-gen-context} \cref{mackey-context-i,mackey-context-ii,mackey-context-iii,mackey-context-iv,mackey-context-v} are satisfied in the current context.
\begin{romenumerate}
\item $\M = \pMulticat$ and $\N = \MoneMod$ are closed multicategories by
\begin{itemize}
\item \cref{smclosed-closed-multicat} and
\item the fact that they are symmetric monoidal closed categories (\cref{thm:pmulticat-smclosed,proposition:EM2-5-1}).
\end{itemize}  
\item $\D$ is, by assumption, a small $\MoneMod$-category.
\item $F = \Monesma$ in \cref{mackey-pmult-mone-setting} is a multifunctor by \cref{Monesma-CatSM}, and $E = \Um$ is a multifunctor by \cref{expl:Um-catmulti}.
\item $\uni = \etahat$ and $\cou = \epzhatinv$ in \cref{mackey-pmult-mone-setting-ii} are multinatural transformations by \cref{etahat-epzhat-monCatnat}.
\item In the current setting, the condition \cref{counit-compatibility} is the equality of the following two pointed multifunctors for each pair of objects $x,y \in \D$. 
\[\begin{tikzpicture}
\draw[0cell]
(0,0) node (a) {\Um\D(x,y)}
(a)+(5,0) node (b) {\Um(\Monesma)\Um \D(x,y)}
;
\draw[1cell=.9]
(a) edge[transform canvas={yshift=.6ex}] node {\etahat_{\Um\D(x,y)}} (b)
(a) edge[transform canvas={yshift=-.4ex}] node[swap] {\Um \epzhat^\inv_{\D(x,y)}} (b)
;
\end{tikzpicture}\]
Since $\D(x,y)$ is a left $\Mone$-module, the equality of these two arrows follows from the right triangle identity \cref{triangleidentities} of the 2-adjunction 
\begin{equation}\label{MonesmaUm-iiadj}
\big((\Monesma), \Um, \etahat, \epzhat\big) \cn \pMulticat \to \MoneMod
\end{equation}
in \cref{MonesmaUmadj}.
\end{romenumerate}
Thus \cref{def:mackey-gen-context} \cref{mackey-context-i,mackey-context-ii,mackey-context-iii,mackey-context-iv,mackey-context-v} hold in the context of \cref{mackey-pmult-mone-setting,mackey-pmult-mone-setting-ii,pMultMone-steq,sSsttri-defin,sSmtri-defin}.

Next, the only assumption in \cref{def:mackey-heq-context} is that the relative category
\[(\N, \cX) = \brb{\MoneMod,\cSM}\]
is a category with weak equivalences (\cref{definition:rel-cat} \eqref{relcat-vi}).  This is true for the following two reasons.
\begin{itemize}
\item By \cref{def:ptmulti-stableeq} the wide subcategory $\cSM \bigsubset \MoneMod$ is created by the functor
\[\Fm \cn \MoneMod \to \permcatsu.\]
\item The class of stable equivalences $\cS \bigsubset \permcatsu$ contains all the isomorphisms and has the 2-out-of-3 property.
\end{itemize} 
By \cref{FstFmMonesma,SM-Sst-def} there are equalities of wide subcategories as follows.
\[\begin{split}
F^\inv\cX &= (\Monesma)^\inv(\cSM)\\
&= (\Monesma)^\inv \Fm^\inv(\cS)\\
&= \Fst^\inv(\cS)\\
&= \cSst \bigsubset \pMulticat
\end{split}\]
The data in \cref{mackey-pmult-mone-functors} are those in \cref{mackey-xiv-functors} in the current context.

Finally, we observe that the components of $\etahat$ and $\epzhatinv$ are stable equivalences.
\begin{itemize}
\item Each component of $\epzhatinv$ \cref{epzm-def} is an isomorphism, hence also an $\Fm$-stable equivalence in $\MoneMod$.
\item The left triangle identity \cref{triangleidentities} for the 2-adjunction \cref{MonesmaUm-iiadj} and the 2-out-of-3 property imply that each component of $\etahat$ is an $\Fst$-stable equivalence in $\pMulticat$.
\end{itemize}
Thus \cref{mackey-gen-xiv} is applicable in the current setting, proving the first assertion.

The second assertion about $(\D_{\Um})^\op$ and $\D^\op$ is an instance of \cref{mackey-xiv-cor}.  It is applicable because $\Um$ is a multifunctor (\cref{expl:Um-catmulti}).
\end{proof}

\begin{explanation}[Homotopy Equivalent Categories of Modules]\label{expl:left-mod-heq}
The categories of enriched diagrams in \cref{mackey-pmult-mone-functors} are categories of left modules by \cref{C-diagram-partner,C-diag-morphism-pn}.  Thus, \cref{mackey-pmulti-mone} asserts that, for each small $\MoneMod$-category $\D$, the functors in \cref{mackey-pmult-mone-functors} are equivalences of homotopy theories between
\begin{itemize}
\item left $\D$-modules in $\MoneMod$ and
\item left $\D_{\Um}$-modules in $\pMulticat$.
\end{itemize}
Furthermore, these equivalences of homotopy theories still hold if $\D$ and $\D_{\Um}$ are replaced by $\D^\op$ and $(\D_{\Um})^\op$, respectively.  We explain the functors $\Umdg$ and $\Monesmadgr$ in \cref{mackey-pmult-mone-functors} in more detail in \cref{sec:Monesmadgr-Umdg} below.

We use the abbreviations in \cref{pMpsu}, so
\[\pM = \pMulticat \andspace \psu = \permcatsu.\]
\cref{mackey-xiv-pmulticat,mackey-xiv-mone,mackey-pmulti-mone} together yield the following diagram of inverse equivalences of homotopy theories for each small $\permcatsu$-category $\C$.
\begin{equation}\label{left-mod-heq}
\begin{tikzpicture}[vcenter]
\def\h{3} \def\v{-1.3}
\draw[0cell=.9]
(0,0) node (a) {\Big(\Monemodcat\big(\C_{\Endm} \scs \MoneMod\big) \scs \cSmtri\Big)}
(a)+(\h,\v) node (b) {\Big(\psucat\big(\C \scs \psu\big) \scs \cStri\Big)}
(a)+(0,2*\v) node (c) {\Big(\pMcat\big(\C_{\Endst} \scs \pM\big) \scs \cSsttri\Big)}
;
\draw[1cell=.8]
(a) edge[transform canvas={xshift=-2em}] node {\Umdg} (c)
(c) edge[transform canvas={xshift=-2.5em}] node {\Monesmadgr} (a)
(a) edge[transform canvas={xshift=3em}] node[pos=.8] {\Fmdgr} (b)
(b) edge[transform canvas={xshift=2em}] node[pos=.9] {\Endmdg} (a)
(c) edge[transform canvas={xshift=3em}] node[swap,pos=.8] {\Fstdgr} (b)
(b) edge[transform canvas={xshift=2em}] node[swap,pos=.9] {\Endstdg} (c)
;
\end{tikzpicture}
\end{equation}
Moreover, by \cref{EndstUmEndmdg,Fmdgr-factorization}, the factorizations 
\[\begin{split}
\Endstdg &= \Umdg \Endmdg \andspace\\
\Fmdgr &= \Fstdgr \Umdg
\end{split}\]
hold in \cref{left-mod-heq}.
\end{explanation}

\section{Explanation of the Equivalences of Homotopy Theories}
\label{sec:Monesmadgr-Umdg}

In this section we explain in detail the inverse equivalences of homotopy theories $\Umdg$ and $\Monesmadgr$ in \cref{mackey-pmult-mone-functors}.  We use the abbreviation
\[\pM = \pMulticat\]
and denote by $\D$ a small $\MoneMod$-category (\cref{def:enriched-category}).  This section is organized as follows.
\begin{itemize}
\item \cref{expl:Umdg} describes the functor $\Umdg$.
\item \cref{expl:Monesmadgr} describes $\Monesmadgr$ on objects.
\item \cref{expl:Monesmadgr-morphism} describes $\Monesmadgr$ on morphisms.
\end{itemize}

\begin{explanation}[The Functor $\Umdg$]\label{expl:Umdg}
The diagram change-of-enrichment functor
\[\begin{tikzpicture}
\draw[0cell]
(0,0) node (a) {\Monemodcat\big(\D \scs \MoneMod\big)}
(a)+(5,0) node (b) {\phantom{\pMcat\big(\D_{\Um} \scs \pM\big)}}
(b)+(0,-.03) node (b') {\pMcat\big(\D_{\Um} \scs \pM\big)}
;
\draw[1cell=.9]
(a) edge node {\Umdg} (b)
;
\end{tikzpicture}\]
is defined in \cref{diag-change-enr-assign} and verified in \cref{gspectra-thm-v}.  To understand its assignments on objects and morphisms, consider
\begin{itemize}
\item $\MoneMod$-functors $A,B \cn \D \to \MoneMod$ (\cref{def:enriched-functor}) and
\item a $\MoneMod$-natural transformation $\psi \cn A \to B$ (\cref{def:enriched-natural-transformation})
\end{itemize}
as in the left diagram below.
\begin{equation}\label{Umdg-def}
\begin{tikzpicture}[baseline={(a.base)}]
\def\t{25} \def\d{.8}
\draw[0cell=.85]
(0,0) node (a) {\D}
(a)+(1.5,0) node (b1) {\phantom{C}}
(b1)+(.35,0) node (b) {\phantom{\MoneMod}}
(b)+(0,.02) node (mone) {\MoneMod}
(b)+(\d,0) node (x) {}
(x)+(1.1,0) node (y) {}
(y)+(\d-.2,0) node (a') {\phantom{\D_{\Um}}}
(a')+(0,-.05) node (a1) {\D_{\Um}}
(a')+(.25,0) node (a'') {\phantom{C}}
(a'')+(1.8,0) node (b2) {\phantom{C}}
(b2)+(.75,0) node (b') {\phantom{(\MoneMod)_{\Um}}}
(b')+(0,-.03) node (b3) {(\MoneMod)_{\Um}}
(b')+(2.2,0) node (c) {\phantom{\pM}}
(c)+(0,.01) node (c') {\pM}
;
\draw[1cell=.8]
(a) edge[bend left=30] node {A} (b1)
(a) edge[bend right=30] node[swap] {B} (b1)
(x) edge[|->] node {\Umdg} (y)
;
\draw[2cell=.9]
node[between=a and b1 at .42, rotate=-90, 2label={above,\psi}] {\Rightarrow}
;
\draw[1cell=.8]
(a'') edge[bend left=\t] node {A_{\Um}} (b2)
(a'') edge[bend right=\t] node[swap] {B_{\Um}} (b2)
(b') edge node {\Umse} (c)
;
\draw[2cell=.8]
node[between=a'' and b2 at .35, rotate=-90, 2label={above,\psi_{\Um}}] {\Rightarrow}
;
\end{tikzpicture}
\end{equation}
Then $\Umdg$ sends $A$, $B$, and $\psi$ to the composites and whiskering as in the right diagram in \cref{Umdg-def}.  In other words, the functor $\Umdg$ 
\begin{itemize}
\item first applies the change-of-enrichment 2-functor along the multifunctor $\Um$ (\cref{Um-multifunctor,mult-change-enrichment})
\[\dUm \cn \Monemodcat \to \pMcat\]
and then
\item composes or whiskers with the standard enrichment $\Umse$ (\cref{gspectra-thm-iii}).
\end{itemize}

Next we describe the object assignment and components of $\Umdg A$.  For each object $x \in \D$, the object assignment of $\Umdg A$ is
\[\big(\Umdg A\big)x = \Um(Ax) \inspace \pM.\]
For objects $x,y \in \D$, the $(x,y)$-component of $A$ is a left $\Mone$-module morphism
\[A_{x,y} \cn \D(x,y) \to \pHom(Ax,Ay).\]
Its adjoint in $\MoneMod$ is a morphism
\[\pn{A_{x,y}} \cn \D(x,y) \sma Ax \to Ay.\]
The $(x,y)$-component of $\Umdg A$ is the following composite in $\MoneMod$.
\[\begin{tikzpicture}
\def\h{2.5} \def\v{1.1} \def\t{15}
\draw[0cell=.9]
(0,0) node (a) {\Um\D(x,y)}
(a)+(\h,-\v) node (b) {\Um\pHom(Ax,Ay)}
(b)+(\h+.2,\v) node (c) {\pHom\big(\Um(Ax) \scs \Um(Ay)\big)}
(b)+(\h,\v) node (c') {\phantom{M}}
;
\draw[1cell=.9]
(a) edge node {\big(\Umdg A\big)_{x,y}} (c)
(a) edge[bend right=\t] node[swap,pos=.3] {\Um(A_{x,y})} (b)
(b) edge[bend right=\t,shorten >=1ex] node[swap,pos=.7] {(\Umse)_{Ax,Ay}} (c')
;
\end{tikzpicture}\]
Its adjoint in $\MoneMod$ is the following composite, with $\Um^2$ the monoidal constraint of $\Um$, which is the identity by \cref{expl:moneinclusion}.
\[\begin{tikzpicture}
\def\h{2.5} \def\v{1.1} \def\t{15}
\draw[0cell=.9]
(0,0) node (a) {\Um\D(x,y) \sma \Um(Ax)}
(a)+(\h,-\v) node (b) {\Um\big(\D(x,y) \sma Ax\big)}
(b)+(\h,\v) node (c) {\Um(Ay)}
;
\draw[1cell=.9]
(a) edge node {\pn{\big(\Umdg A\big)_{x,y}}} (c)
(a) edge[bend right=\t] node[swap,pos=.3] {\Um^2} (b)
(b) edge[bend right=\t] node[swap,pos=.7] {\Um(\pn{A_{x,y}})} (c)
;
\end{tikzpicture}\]

The $x$-component of $\Umdg \psi$ is the morphism
\[\big(\Umdg \psi\big)_x = \Um(\psi_x) \cn \Um(Ax) \to \Um(Bx) \inspace \pM.\]
This is obtained from the $x$-component of $\psi$ by applying $\Um$.
\end{explanation}

\begin{explanation}[$\Monesmadgr$ on Objects]\label{expl:Monesmadgr}
We abbreviate the top functor in \cref{mackey-pmult-mone-functors} to
\begin{equation}\label{sfH-Monesmadgr}
\sfH = \Monesmadgr \cn \pMcat\big(\D_{\Um} \scs \pM\big) \to \Monemodcat(\D,\MoneMod).
\end{equation}
The functor $\sfH$ sends an $\pM$-functor
\[P \cn \D_{\Um} \to \pM\]
to the following composite $\MoneMod$-functor.
\begin{equation}\label{HP-composite}
\begin{tikzpicture}[vcenter]
\def\g{.3} \def\v{1.3} \def\w{.6}
\draw[0cell=.9]
(0,0) node (a) {\D} 
(a)+(\g,-\v) node (b) {\D_{(\Monesma)\Um}}
(b)+(2.2,0) node (b') {(\D_{\Um})_{\Monesma}}
(b')+(3.5,0) node (c) {(\pM)_{\Monesma}}
(c)+(\g,\v) node (d) {\MoneMod}
;
\draw[1cell=2]
(b) edge[-,double equal sign distance] (b')
;
\draw[1cell=.9]
(a) edge node {\sfH P} (d)
(a) edge node[pos=.2,swap] {\D_{\epzhatinv}} (b)
(b') edge node {P_{\Monesma}} (c)
(c) edge node[pos=.8,swap] {\Monesmase} (d)
;
\end{tikzpicture}
\end{equation}
\begin{itemize}
\item The $\MoneMod$-functor $\D_{\epzhatinv}$ is the $\D$-component of the 2-natural transformation (\cref{dtheta-twonat})
\[\begin{tikzpicture}[baseline={(a.base)}]
\def\s{22}
\draw[0cell]
(0,0) node (a) {\phantom{\M}}
(a)+(3,0) node (b) {\phantom{\N}}
(a)+(-.7,0) node (a') {\Monemodcat}
(b)+(.7,0) node (b') {\Monemodcat}
;
\draw[1cell=.85]
(a) edge[bend left=\s] node {1} (b)
(a) edge[bend right=\s] node[swap] {\dMonesmaUm = \dMonesma \,\dUm} (b)
;
\draw[2cell]
node[between=a and b at .35, rotate=-90, 2label={above,\depzhatinv}] {\Rightarrow}
;
\end{tikzpicture}\]
induced by the multinatural transformation (\cref{etahat-epzhat-monCatnat})
\[\epzhatinv \cn 1_{\MoneMod} \to (\Monesma)\Um.\]
\item $P_{\Monesma}$ is the image of $P$ under the change-of-enrichment 2-functor 
\begin{equation}\label{dMonesmas}
\dMonesma \cn \pMcat \to \Monemodcat
\end{equation}
along the multifunctor $\Monesma$ (\cref{Monesma-CatSM,mult-change-enrichment}).
\item $\Monesmase$ is the standard enrichment $\MoneMod$-functor of $\Monesma$ (\cref{gspectra-thm-iii}).
\end{itemize}

The $\MoneMod$-functor $\sfH P$ in \cref{HP-composite} sends each object $x \in \D$ to
\[(\sfH P)x = \Mone \sma Px \inspace \MoneMod.\]
For objects $x,y \in \D$, the $(x,y)$-component of $P$ is a pointed multifunctor
\[P_{x,y} \cn \Um \D(x,y) \to \pHom(Px,Py).\]
Its adjoint in $\pM$ is a pointed multifunctor
\[\pn{P_{x,y}} \cn \Um\D(x,y) \sma Px \to Py.\]
The $(x,y)$-component of $\sfH P$ is the following composite morphism in $\MoneMod$.
\[
\]
The adjoint of $(\sfH P)_{x,y}$ in $\MoneMod$ is the following composite, with $\binprod_{\ord{1},\ord{1}}$ the partition product in \cref{eq:part-prod-1}.
\[%
\]
The bottom horizontal isomorphism permutes the middle two factors.
\end{explanation}

\begin{explanation}[$\Monesmadgr$ on Morphisms]\label{expl:Monesmadgr-morphism}
Consider an $\pM$-natural transformation $\psi$ (\cref{def:enriched-natural-transformation}) as in the left diagram below.
\[%
\]
With the notation in \cref{sfH-Monesmadgr}, the $\MoneMod$-natural transformation $\sfH\psi$, as in the right diagram above, is the following whiskering.
\begin{equation}\label{sfHpsi-def}
%
\end{equation}
The whiskering \cref{sfHpsi-def} is the one in \cref{Fdgrpsi-def} in the current context.  It is obtained from \cref{HP-composite} by replacing $P$ with $\psi$.
\begin{itemize}
\item The $\MoneMod$-functors $\D_{\epzhatinv}$ and $\Monesmase$ are as in \cref{HP-composite}.
\item $\psi_{\Monesma}$ is the change of enrichment of $\psi$ under $\dMonesma$ in \cref{dMonesmas}.
\end{itemize}
For each object $x \in \D$, the $x$-component of $\psi$ is a pointed multifunctor
\[\psi_x \cn Px \to Qx.\]
The $x$-component of $\sfH\psi$ is the following morphism in $\MoneMod$.
\[(\sfH\psi)_x = \Mone \sma \psi_x \cn \Mone \sma Px \to \Mone \sma Qx\]
This is obtained from $\psi_x$ by applying $\Monesma$ (\cref{MonesmaUmadj}).
\end{explanation}

\part*{Appendices}
\label{part:appendices}
\appendix

\chapter{Categories}
\label{ch:prelim}
In this appendix we review monoidal categories and 2-categories.  The following table summarizes the main content in this appendix.
\smallskip
\begin{center}
\resizebox{.8\width}{!}{
{\renewcommand{\arraystretch}{1.4}%
{\setlength{\tabcolsep}{1ex}
\begin{tabular}{|c|c|}\hline
\multicolumn{2}{|c|}{\appnumname{sec:monoidalcat}} \\ \hline
Grothendieck universes & \ref{def:universe} and \ref{conv:universe} \\ \hline
monoidal categories (braided, symmetric, closed) & \ref{def:monoidalcategory} (\ref{def:braidedmoncat}, \ref{def:symmoncat}, \ref{def:closedcat}) \\ \hline
monoids, modules, and commutative monoids & \ref{def:monoid}, \ref{def:modules}, and \ref{def:commonoid} \\ \hline
diagram categories & \ref{def:diagramcat} \\ \hline
monoidal functors and natural transformations & \ref{def:monoidalfunctor} and \ref{def:monoidalnattr} \\ \hline\hline
\multicolumn{2}{|c|}{\appnumname{sec:twocategories}} \\ \hline
2-categories, 2-functors, and 2-natural transformations & \ref{def:twocategory}, \ref{def:twofunctor}, and \ref{def:twonaturaltr} \\ \hline
2-category of small (permutative, 2-) categories & \ref{ex:catastwocategory} (\ref{def:permcat}, \ref{ex:iicat}) \\ \hline
2-adjunctions & \ref{def:twoadjunction} \\ \hline
\end{tabular}}}}
\end{center}
\medskip
References for \cref{sec:monoidalcat,sec:twocategories} are \cite{joyal-street,maclane,cerberusI,cerberusII} and \cite{johnson-yau}, respectively.

\section{Monoidal Categories}
\label{sec:monoidalcat}

In this section we review Grothendieck universes, monoidal categories, monoidal functors, and monoidal natural transformations.

\begin{definition}\label{def:universe}
A \index{universe}\emph{universe} is a set\label{not:universe} $\calu$ that satisfies \cref{univ1,univ2,univ3,univ4} below:
\begin{romenumerate}
\item\label{univ1} If $a \in \calu$ and $b \in a$, then $b\in\calu$.
\item\label{univ2} If $a \in \calu$, then $\Pset(a)\in\calu$, where $\Pset(a)$ is the set of subsets of $a$.
\item\label{univ3} If $a \in \calu$ and $x_j \in \calu$ for each $j \in a$, then the union $\bigcup_{j\in a} x_j \in \calu$.
\item\label{univ4} $\bbN \in \calu$, where $\bbN$ is the set of finite ordinals.\defmark
\end{romenumerate}
\end{definition}

\begin{convention}[Universe]\label{conv:universe}\index{Grothendieck Universe}\index{universe}\index{Axiom of Universes}\index{convention!universe}
We assume Grothendieck's \emph{Axiom of Universes}: 
\begin{quote}
Every set belongs to some universe.  
\end{quote}
We fix a universe $\cU$.  An element in $\cU$ is called a \emph{set}.  A subset of $\cU$ is called a \emph{class}.  A categorical structure is called \emph{small} if it has a set of objects.  We automatically replace $\cU$ by a larger universe $\cV$ in which $\cU$ is a set whenever necessary.  For more discussion of universes, see \cite[Section 1.1]{johnson-yau} and \cite{agv,maclane-foundation}.
\end{convention}

\begin{definition}\label{def:monoidalcategory}\index{monoidal category}\index{category!monoidal}
A \emph{monoidal category} is a sextuple\label{not:monoidalcat} 
\[(\C,\otimes,\tu,\alpha,\lambda,\rho)\]
consisting of the following data.
\begin{itemize}
\item $\C$ is a category.
\item \label{notation:monoidal-product}$\otimes \cn \C \times \C \to \C$ is a functor, which is called the \index{monoidal category!monoidal product}\emph{monoidal product}.
\item \label{not:monoidalunit}$\tu \in \C$ is an object, which is called the \index{monoidal category!monoidal unit}\emph{monoidal unit}.
\item \label{not:associativityiso}\label{not:unitisos}$\alpha$, $\lambda$, and $\rho$ are natural isomorphisms with the following components for objects $x,y,z \in \C$.
\[\begin{tikzcd}[column sep=huge,row sep=tiny]
(x \otimes y) \otimes z \ar{r}{\alpha_{x,y,z}}[swap]{\iso} & x \otimes (y \otimes z)
\end{tikzcd}\]\vspace{-1em}
\[\begin{tikzcd}[column sep=large]
\tu \otimes x \ar{r}{\lambda_x}[swap]{\iso} & x & x \otimes \tu \ar{l}{\iso}[swap]{\rho_x}
\end{tikzcd}\]
They are called the \index{associativity!isomorphism}\index{monoidal category!associativity isomorphism}\emph{associativity isomorphism}, the \index{left unit isomorphism}\index{monoidal category!left unit isomorphism}\emph{left unit isomorphism}, and the \index{right unit isomorphism}\index{monoidal category!right unit isomorphism}\emph{right unit isomorphism}, respectively.
\end{itemize}
The data above are required to make the following \emph{middle unity}\index{monoidal category!middle unity axiom} and \emph{pentagon}\index{monoidal category!pentagon axiom}\index{pentagon axiom} diagrams commute for objects $w,x,y,z \in \C$, where $\otimes$ is omitted to save space.
\begin{equation}\label{monoidalcataxioms}
\begin{tikzpicture}[xscale=1,yscale=1,vcenter]
\tikzset{0cell/.append style={nodes={scale=.85}}}
\tikzset{1cell/.append style={nodes={scale=.85}}}
\def\h{1.6} \def\g{1.3} \def\v{.9} \def\u{2}
\draw[0cell]
(0,0) node (a) {(x \tu) y}
(a)++(2,0) node (b) {x (\tu y)}
(a)++(1,-1.3) node (c) {xy}
;
\draw[1cell]  
(a) edge node[swap,pos=.4] {\rho_x 1_y} (c)
(a) edge node {\alpha_{x,\tu,y}} (b)
(b) edge node[pos=.4] {1_x \lambda_y} (c)
;
\begin{scope}[shift={(6,.5)}]
\draw[0cell]
(0,0) node (x0) {(wx)(yz)}
(x0)++(-\h,-\v) node (x11) {((wx)y)z}
(x0)++(\h,-\v) node (x12) {w(x(yz))}
(x0)++(-\g,-\u) node (x21) {(w(xy))z}
(x0)++(\g,-\u) node (x22) {w((xy)z)}
;
\draw[1cell]
(x11) edge node[pos=.2] {\al_{wx,y,z}} (x0)
(x0) edge node[pos=.8] {\al_{w,x,yz}} (x12)
(x11) edge node[swap,pos=.2] {\al_{w,x,y} 1_z} (x21)
(x21) edge node {\al_{w,xy,z}} (x22)
(x22) edge node[swap,pos=.8] {1_w \al_{x,y,z}} (x12)
;
\end{scope}
\end{tikzpicture}
\end{equation}
This finishes the definition of a monoidal category.  Moreover, we define the following.
\begin{itemize}
\item We call a monoidal category \emph{strict}\index{monoidal category!strict}\index{strict!monoidal category} if $\alpha$, $\lambda$, and $\rho$ are identity natural transformations.
\item We also call a (monoidal) category a \index{monoidal 1-category}\index{1-category}\emph{(monoidal) 1-category}.\defmark
\end{itemize}
\end{definition}

\begin{remark}[Unity Properties]\label{rk:monoidalidentities}
In each monoidal category, the unit isomorphisms agree at the monoidal unit:\index{monoidal category!unity properties} 
\begin{equation}\label{lambda=rho}
\lambda_{\tensorunit} = \rho_{\tensorunit} \cn \tensorunit \otimes \tensorunit \to \tensorunit.
\end{equation}
Moreover, the following unity diagrams commute.  
\begin{equation}\label{moncat-other-unit-axioms}
\begin{tikzcd}[column sep=0ex]
(\tensorunit \otimes x) \otimes y \arrow{dr}[swap]{\lambda_x \otimes 1_y} \arrow{rr}{\alpha_{\tensorunit,x,y}}
&& \tensorunit \otimes (x \otimes y) \arrow{dl}{\lambda_{x \otimes y}}\\ 
& x \otimes y & \end{tikzcd}\qquad
\begin{tikzcd}[column sep=0ex]
(x \otimes y) \otimes \tensorunit \arrow{dr}[swap]{\rho_{x \otimes y}} \arrow{rr}{\alpha_{x,y,\tensorunit}}
&& x \otimes (y \otimes \tensorunit) \arrow{dl}{1_x \otimes \rho_y}\\ 
& x \otimes y & \end{tikzcd}
\end{equation} 
See \cite[Section 2.2]{johnson-yau} for the proofs.
\end{remark}

\begin{definition}\label{def:monoid}
A \emph{monoid}\index{monoid} in a monoidal category $\C$ is a triple\label{notation:monoid} $(x,\mcomp,i)$ consisting of the following data.
\begin{itemize}
\item $x \in \C$ is an object.
\item $\mcomp \cn x \otimes x \to x$ is a morphism, which is called the \index{multiplication}\emph{multiplication}.
\item $i \cn \tensorunit \to x$ is a morphism, which is called the \index{unit}\emph{unit}.
\end{itemize}
The data above are required to make the following associativity and unity diagrams commute.
\[\begin{tikzcd}[column sep=large]
(x \otimes x) \otimes x \arrow{dd}[swap]{\mcomp \otimes 1_x} \rar{\alpha} & x \otimes (x \otimes x) \dar{1_x \otimes \mcomp}\\ 
& x \otimes x \dar{\mcomp}\\  
x \otimes x \arrow{r}{\mcomp} & x
\end{tikzcd}\qquad
\begin{tikzcd}[column sep=large]
\tensorunit \otimes x \ar{d}[swap]{i \otimes 1_x} \ar{dr}{\lambda_x} & \\ 
x \otimes x \ar{r}[pos=.3]{\mcomp} & x \\
x \otimes \tensorunit \ar{u}{1_x \otimes i} \ar{ur}[swap]{\rho_x} &
\end{tikzcd}\]
A \emph{morphism} of monoids
\[\begin{tikzcd}[column sep=large]
(x,\mcomp^x,i^x) \ar{r}{f} & (y,\mcomp^y,i^y)
\end{tikzcd}\] 
is a morphism $f \cn x \to y$ in $\C$ such that the diagrams
\[\begin{tikzcd}[column sep=large]
x \otimes x \dar[swap]{\mcomp^x} \rar{f\otimes f} 
& y \otimes y \dar{\mcomp^y}\\ 
x \rar{f} & y\end{tikzcd} \qquad
\begin{tikzcd}[column sep=large]
\tensorunit \arrow{dr}[swap]{i^y} \rar{i^x} & x \dar{f}\\ 
& y\end{tikzcd}\]
commute.
\end{definition}

\begin{definition}\label{def:modules}
Suppose $(x,\mcomp,i)$ is a monoid in a monoidal category $(\C,\otimes,\tu)$.
\begin{enumerate}
\item A \emph{left $x$-module}\index{module} is a pair $(a,\mu)$\label{not:leftmodule} consisting of the following data.
\begin{itemize}
\item $a$ is an object in $\C$.
\item $\mu \cn x \otimes a \to a$ is a morphism, called the \emph{structure morphism}.
\end{itemize}
The data above are required to make the following associativity and unity diagrams commute.
\[\begin{tikzcd}[column sep=large]
(x \otimes x) \otimes a \arrow{dd}[swap]{\mcomp \otimes 1_a} \rar{\alpha} & x \otimes (x \otimes a) \dar{1_x \otimes \mu}\\ 
& x \otimes a \dar{\mu}\\  
x \otimes a \arrow{r}{\mu} & a
\end{tikzcd}\qquad
\begin{tikzcd}[column sep=large]
\tensorunit \otimes a \ar{d}[swap]{i \otimes 1_a} \ar{dr}{\lambda_a} & \\ 
x \otimes a \ar{r}[pos=.3]{\mu} & a
\end{tikzcd}\]
\item A \emph{morphism} of left $x$-modules
\[\begin{tikzcd}[column sep=large]
(a,\mu^a) \ar{r}{f} & (b,\mu^b)
\end{tikzcd}\] 
is a morphism $f \cn a \to b$ in $\C$ such that the following diagram commutes.
\[\begin{tikzcd}[column sep=large]
x \otimes a \dar[swap]{\mu^a} \rar{1_x \otimes f} 
& x \otimes b \dar{\mu^b}\\ 
a \rar{f} & b
\end{tikzcd}\]
\item \emph{Right $x$-modules}, with structure morphisms $a \otimes x \to a$, and their morphisms are defined similarly.\defmark
\end{enumerate} 
\end{definition}

\begin{definition}\label{def:braidedmoncat}\index{braided monoidal category}\index{monoidal category!braided}\index{category!braided monoidal}
A \emph{braided monoidal category} is a pair $(\C,\xi)$ consisting of the following data.
\begin{itemize}
\item $\C$ is a monoidal category (\cref{def:monoidalcategory}).
\item $\xi$ is a natural isomorphism, which is called the \index{braiding}\emph{braiding}, with components\label{notation:symmetry-iso} 
\[x \otimes y \fto[\iso]{\xi_{x,y}} y \otimes x \forspace x,y \in \C.\]
\end{itemize}
The hexagon diagrams\index{hexagon diagram}
\begin{equation}\label{hexagon-braided}
\def\sb{\scalebox{.85}}
\begin{tikzpicture}[scale=.75,commutative diagrams/every diagram]
\node (P0) at (0:2cm) {\sb{$y \otimes (z \otimes x)$}};
\node (P1) at (60:2cm) {\makebox[3ex][l]{\sb{$y \otimes (x \otimes z)$}}};
\node (P2) at (120:2cm) {\makebox[3ex][r]{\sb{$(y \otimes x) \otimes z$}}};
\node (P3) at (180:2cm) {\sb{$(x \otimes y) \otimes z$}};
\node (P4) at (240:2cm) {\makebox[3ex][r]{\sb{$x \otimes (y \otimes z)$}}};
\node (P5) at (300:2cm) {\makebox[3ex][l]{\sb{$(y \otimes z) \otimes x$}}};
\path[commutative diagrams/.cd, every arrow, every label]
(P3) edge node[pos=.25] {\sb{$\xitimes_{x,y}\otimes 1_z$}} (P2)
(P2) edge node {\sb{$\alpha$}} (P1)
(P1) edge node[pos=.75] {\sb{$1_y \otimes \xitimes_{x,z}$}} (P0)
(P3) edge node[swap,pos=.4] {\sb{$\alpha$}} (P4)
(P4) edge node {\sb{$\xitimes_{x,y \otimes z}$}} (P5)
(P5) edge node[swap,pos=.6] {\sb{$\alpha$}} (P0);
\end{tikzpicture}
\qquad
\begin{tikzpicture}[scale=.75, commutative diagrams/every diagram]
\node (P0) at (0:2cm) {\sb{$(z \otimes x) \otimes y$}};
\node (P1) at (60:2cm) {\makebox[3ex][l]{\sb{$(x \otimes z) \otimes y$}}};
\node (P2) at (120:2cm) {\makebox[3ex][r]{\sb{$x \otimes (z \otimes y)$}}};
\node (P3) at (180:2cm) {\sb{$x \otimes (y \otimes z)$}};
\node (P4) at (240:2cm) {\makebox[3ex][r]{\sb{$(x \otimes y) \otimes z$}}};
\node (P5) at (300:2cm) {\makebox[3ex][l]{\sb{$z \otimes (x \otimes y)$}}};
\path[commutative diagrams/.cd, every arrow, every label]
(P3) edge node[pos=.25] {\sb{$1_x \otimes \xitimes_{y,z}$}} (P2)
(P2) edge node {\sb{$\alpha^\inv$}} (P1)
(P1) edge node[pos=.75] {\sb{$\xitimes_{x,z}\otimes 1_y$}} (P0)
(P3) edge node[swap,pos=.4] {\sb{$\alpha^\inv$}} (P4)
(P4) edge node {\sb{$\xi_{x \otimes y, z}$}} (P5)
(P5) edge node[swap,pos=.6] {\sb{$\alpha^\inv$}} (P0);
\end{tikzpicture}
\end{equation}
are required to commute for objects $x,y,z \in \C$. 
\end{definition}

\begin{remark}[Unity Properties]\label{rk:braidedunity}
In each braided monoidal category, the following unity diagrams\index{braided monoidal category!unity properties} commute for each object $x$.
\begin{equation}\label{braidedunity}
\begin{tikzcd}[column sep=tiny]
x \otimes \tensorunit \ar{dr}[swap]{\rho_x} \ar{rr}{\xi_{x,\tensorunit}} 
&& \tensorunit \otimes x \ar{dl}{\lambda_x}\\ & x &
\end{tikzcd}\qquad
\begin{tikzcd}[column sep=tiny] 
\tu \otimes x \ar{dr}[swap]{\lambda_x} \ar{rr}{\xi_{\tu,x}} 
&& x \otimes \tu \ar{dl}{\rho_x}\\ & x &
\end{tikzcd}
\end{equation}
See \cite[1.3.21]{cerberusII} for the proof. 
\end{remark}

\begin{definition}\label{def:symmoncat}\index{symmetric monoidal category}\index{monoidal category!symmetric}\index{category!symmetric monoidal}
A \emph{symmetric monoidal category} is a pair $(\C,\xi)$ consisting of the following data.
\begin{itemize}
\item $\C$ is a monoidal category (\cref{def:monoidalcategory}).
\item $\xi$ is a natural isomorphism, which is called the \emph{symmetry isomorphism} or the \emph{braiding}, with components
\[x \otimes y \fto[\iso]{\xi_{x,y}} y \otimes x \forspace x,y \in \C.\]
\end{itemize}
The data above are required to make the following symmetry\index{symmetry axiom} and \index{hexagon diagram}hexagon diagrams commute for objects $x,y,z \in \C$. 
\begin{equation}\label{symmoncatsymhexagon}
\begin{tikzpicture}[xscale=3,yscale=1.2,vcenter]
\def\h{.1}
\draw[0cell=.8] 
(0,0) node (a) {x \otimes y}
(a)++(.7,0) node (c) {x \otimes y}
(a)++(.35,-1) node (b) {y \otimes x}
;
\draw[1cell=.8] 
(a) edge node {1_{x \otimes y}} (c)
(a) edge node [swap,pos=.3] {\xi_{x,y}} (b)
(b) edge node [swap,pos=.7] {\xi_{y,x}} (c)
;
\begin{scope}[shift={(1.5,.5)}]
\draw[0cell=.8] 
(0,0) node (x11) {(y \otimes x) \otimes z}
(x11)++(.9,0) node (x12) {y \otimes (x \otimes z)}
(x11)++(-\h,-1) node (x21) {(x \otimes y) \otimes z}
(x12)++(\h,-1) node (x22) {y \otimes (z \otimes x)}
(x11)++(0,-2) node (x31) {x \otimes (y \otimes z)}
(x12)++(0,-2) node (x32) {(y \otimes z) \otimes x}
;
\draw[1cell=.8]
(x21) edge node[pos=.25] {\xi_{x,y} \otimes 1_z} (x11)
(x11) edge node {\alpha} (x12)
(x12) edge node[pos=.75] {1_y \otimes \xi_{x,z}} (x22)
(x21) edge node[swap,pos=.25] {\alpha} (x31)
(x31) edge node {\xi_{x,y \otimes z}} (x32)
(x32) edge node[swap,pos=.75] {\alpha} (x22)
;
\end{scope}
\end{tikzpicture}
\end{equation}
A \index{permutative category}\index{category!permutative}\emph{permutative category} is a strict symmetric monoidal category.  For a generic permutative category, we often write its monoidal product and monoidal unit as \label{not:opluse}$\oplus$ and $\pu$, respectively.
\end{definition}

\begin{remark}[Symmetry Implies Braided]\label{rk:smcat}
The symmetry axiom, $\xi_{y,x} \xi_{x,y} = 1$, implies that the hexagons \cref{hexagon-braided} are equivalent.  Thus, a symmetric monoidal category is precisely a braided monoidal category that satisfies the symmetry axiom.  As a result, the unity diagrams \cref{braidedunity} commute in each symmetric monoidal category.
\end{remark}

The following permutative category plays an important role in both Segal and Elmendorf-Mandell $K$-theory (\cref{sec:Gamma-cat,sec:Gstar-cat}).
\begin{definition}[Pointed Finite Sets]\label{def:ordn}
We define the permutative category\index{finite sets!pointed}\index{pointed!finite sets}\label{not:Fskel}
\[\big(\Fskel, \sma, \ord{1}, \xi\big)\]
as follows.
\begin{description}
\item[Objects] The objects in $\Fskel$ are the pointed finite sets\label{not:ordn}
\[\ord{n} = \{0,\ldots,n\} \forspace n \geq 0\]
with basepoint 0.
\item[Morphisms] For $m, n \geq 0$, the set of morphisms $\Fskel(\ord{m}, \ord{n})$ is the set of pointed functions $\ord{m} \to \ord{n}$, that is, functions that preserve the basepoint 0.
\item[Monoidal Product] It is given on objects by the smash product of pointed finite sets and the lexicographic ordering,
\[\ord{m} \sma \ord{n} = \ord{mn}.\]
It is given explicitly by the identification
\[\ord{m} \sma \ord{n} \ni (i,j) \mapsto 
\begin{cases}
0 \in \ord{mn} & \text{if either $i$ or $j$ is 0, and}\\
(i-1)n + j \in \ord{mn} & \text{if $i,j>0$}.
\end{cases}\]
This extends to pointed functions by functoriality of the smash product of pointed finite sets.  The monoidal product is strictly associative.
\item[Monoidal Unit] The strict monoidal unit is $\ord{1} = \{0,1\}$.
\item[Braiding] Its component at $\ord{m}, \ord{n}$ is the pointed bijection
\[\begin{tikzcd}[column sep=large]
\ord{m} \sma \ord{n} \ar{r}{\xi_{\ord{m},\ord{n}}}[swap]{\iso} & \ord{n} \sma \ord{m}
\end{tikzcd}\]
given by
\[\ord{m} \sma \ord{n} \ni (i,j) \mapsto 
\begin{cases}
0 & \text{if either $i$ or $j$ is 0, and}\\
(j,i) & \text{if $i,j>0$}.
\end{cases}\]
\end{description}
This finishes the definition of $\Fskel$.
\end{definition}

\begin{definition}\label{def:commonoid}
In a symmetric monoidal category $(\C,\xi)$, a \emph{commutative monoid}\index{commutative monoid}\index{monoid!commutative} is a monoid $(x,\mcomp,i)$ as in \cref{def:monoid} such that the diagram
\[\begin{tikzcd}[column sep=1ex]
x \otimes x \arrow{dr}[swap]{\mcomp} \arrow{rr}{\xi_{x,x}}
&& x \otimes x \arrow{dl}{\mcomp}\\ 
& x & \end{tikzcd}\]
commutes.
\end{definition}

\begin{definition}\label{def:closedcat}\index{closed!category}\index{category!closed}\index{internal hom}\index{hom!internal}
A symmetric monoidal category $(\C,\otimes)$ is \emph{closed} if, for each object $x \in \C$, the functor 
\[- \otimes x \cn \C \to \C\] 
admits a right adjoint, which is called an \emph{internal hom}.  A right adjoint of $- \otimes x$ is denoted by $\Hom(x,-)$ or \label{notation:internal-hom}$[x,-]$.
\end{definition}

\begin{definition}[Diagrams]\label{def:diagramcat}\index{diagram!category}\index{category!diagram}
For a small category $\cB$ and a category $\C$, the \emph{diagram category} \label{not:DC}$\BC$ is defined by the following data. 
\begin{itemize}
\item Its objects are functors $\cB \to \C$.
\item Its morphisms are natural transformations between such functors.
\item Identity morphisms are identity natural transformations.
\item Composition is vertical composition of natural transformations.  
\end{itemize}
Moreover, we define the following.
\begin{itemize}
\item A \emph{$\cB$-diagram in $\C$} is a functor $\cB \to \C$.
\item For the category $\Cat$ of small categories and functors, a functor $\cB \to \Cat$ is called a \index{indexed category}\index{category!indexed}\emph{$\cB$-indexed category}.\defmark
\end{itemize}
\end{definition}

\begin{example}[Small Categories]\label{ex:cat}\index{category!of small categories}
\label{not:cat}$(\Cat, \times, \boldone, [,])$ is a symmetric monoidal closed category, where $\Cat$ is the category of small categories and functors.
\begin{itemize}
\item The monoidal product is the Cartesian product, denoted $\times$.  
\item The monoidal unit is the terminal category $\boldone$ with only one object $*$ and its identity morphism $1_*$.  
\item The closed structure $[,]$ is given by diagram categories (\cref{def:diagramcat}). 
\end{itemize} 
The category $\Cat$ is both complete and cocomplete.  For an elementary proof of its cocompleteness, see \cite[Section 1.4]{yau-involutive}.
\end{example}

\subsection*{Monoidal Functors and Natural Transformations}

\begin{definition}\label{def:monoidalfunctor}\index{monoidal functor}\index{functor!monoidal}
Suppose $\C$ and $\D$ are monoidal categories.  A \emph{monoidal functor} 
\[(F, F^2, F^0) \cn \C \to \D\]
is a triple consisting of the following data.
\begin{itemize}
\item $F \cn \C \to \D$ is a functor.
\item $F^0 \cn \tu \to F\tu$ is a morphism in $\D$, which is called the \index{unit constraint}\index{constraint!unit}\emph{unit constraint}.
\item $F^2$ is a natural transformation, which is called the \index{monoidal constraint}\index{constraint!monoidal}\emph{monoidal constraint}, with components
\[\begin{tikzcd}[column sep=large]
Fx \otimes Fy \ar{r}{F^2_{x,y}} & F(x \otimes y)
\end{tikzcd}
\forspace x,y \in \C.\]
\end{itemize}
The data above are required to make the following unity and associativity diagrams commute for objects $x,y,z \in \C$.  
\begin{equation}\label{monoidalfunctorunity}
\begin{tikzcd}[cells={nodes={scale=.9}}]
\tensorunit \otimes Fx \dar[swap]{F^0 \otimes 1_{Fx}} \rar{\lambda_{Fx}} & Fx \\ 
F\tensorunit \otimes Fx \rar{F^2} & F(\tensorunit \otimes x)
\uar[swap]{F\lambda_x}
\end{tikzcd}
\qquad
\begin{tikzcd}[cells={nodes={scale=.9}}]
Fx \otimes \tensorunit \dar[swap]{1_{Fx} \otimes F^0} \rar{\rho_{Fx}} & Fx \\ 
Fx \otimes F\tensorunit \rar{F^2} & F(x \otimes \tensorunit)
\uar[swap]{F\rho_x}
\end{tikzcd}
\end{equation}
\begin{equation}\label{monoidalfunctorassoc}
\begin{tikzcd}[column sep=large,cells={nodes={scale=.9}}]
\bigl(Fx \otimes Fy\bigr) \otimes Fz \rar{\alpha} \dar[swap]{F^2 \otimes 1_{Fz}} 
& Fx \otimes \bigl(Fy \otimes Fz\bigr) \dar{1_{Fx} \otimes F^2}\\
F(x \otimes y) \otimes Fz \dar[swap]{F^2} & Fx \otimes F(y \otimes z) \dar[d]{F^2}\\
F\bigl((x \otimes y) \otimes z\bigr) \rar{F\alpha} &
F\bigl(x \otimes (y \otimes z)\bigr)
\end{tikzcd}
\end{equation}
A monoidal functor $(F,F^2,F^0)$ is
\begin{itemize}
\item \index{monoidal functor!strictly unital}\index{strictly unital!monoidal functor}\emph{strictly unital} if $F^0$ is the identity morphism;
\item \index{monoidal functor!strong}\index{strong monoidal!functor}\emph{strong} if $F^0$ and $F^2$ are isomorphisms; and
\item \index{monoidal functor!strict}\index{strict monoidal!functor}\emph{strict} if $F^0$ and $F^2$ are identities.
\end{itemize}
An \emph{identity monoidal functor} has $F$, $F^0$, and $F^2$ all given by identities.

Moreover, a monoidal functor $(F,F^2,F^0)$ between braided monoidal categories $\C$ and $\D$ (\cref{def:braidedmoncat}) is a \index{braided monoidal functor}\index{monoidal functor!braided}\emph{braided monoidal functor} if the following diagram commutes for objects $x,y \in \C$.  
\begin{equation}\label{monoidalfunctorbraiding}
\begin{tikzcd}[column sep=large]
Fx \otimes Fy \dar[swap]{F^2} \rar{\xi_{Fx,Fy}}[swap]{\cong} & Fy \otimes Fx \dar{F^2} \\ 
F(x \otimes y) \rar{F\xi_{x,y}}[swap]{\cong} & F(y \otimes z)
\end{tikzcd}
\end{equation}
A \index{symmetric monoidal functor}\index{monoidal functor!symmetric}\emph{symmetric monoidal functor} is a braided monoidal functor between symmetric monoidal categories (\cref{def:symmoncat}).
\end{definition}

\begin{definition}\label{def:mfunctor-comp}
Suppose given monoidal functors
\[\C \fto{(F,F^2,F^0)} \D \fto{(G,G^2,G^0)} \E.\]
The \emph{composite monoidal functor}
\[\C \fto{\big(GF, (GF)^2, (GF)^0\big)} \E\]
has unit constraint given by the composite
\[\tu \fto{G^0} G\tu \fto{G(F^0)} GF\tu\]
and monoidal constraint given by the composite
\[GFx \otimes GFy \fto{G^2_{Fx,Fy}} G(Fx \otimes Fy) \fto{G(F^2_{x,y})} GF(x \otimes y)\]
for objects $x,y \in \C$.
\end{definition}

\begin{definition}\label{def:monoidalnattr}
Suppose
\[(F, F^2, F^0) \andspace (G,G^2,G^0) \cn \C \to \D\]
are monoidal functors between monoidal categories $\C$ and $\D$.  A \index{monoidal natural transformation}\index{natural transformation!monoidal}\emph{monoidal natural transformation} $\theta \cn F \to G$ is a natural transformation of the underlying functors such that the following unit constraint and monoidal constraint diagrams in $\D$ commute for objects $x,y \in \C$.
\begin{equation}\label{monnattr}
\begin{tikzpicture}[xscale=2,yscale=1.3,vcenter]
\draw[0cell=.9]
(0,0) node (x11) {\tu}
(x11)++(1,0) node (x12) {F\tu}
(x12)++(0,-1) node (x2) {G\tu}
(x12)++(1.2,0) node (y11) {Fx \otimes Fy}
(y11)++(1.5,0) node (y12) {Gx \otimes Gy}
(y11)++(0,-1) node (y21) {F(x \otimes y)}
(y12)++(0,-1) node (y22) {G(x \otimes y)}
;
\draw[1cell=.9]  
(x11) edge node {F^0} (x12)
(x11) edge node[swap] {G^0} (x2)
(x12) edge node {\theta_{\tu}} (x2)
(y11) edge node {\theta_x \otimes \theta_y} (y12)
(y12) edge node {G^2} (y22)
(y11) edge node[swap] {F^2} (y21)
(y21) edge node {\theta_{x \otimes y}} (y22)
;
\end{tikzpicture}
\end{equation}
A (monoidal) natural transformation is also denoted by
\begin{equation}\label{twocellnotation}
\begin{tikzpicture}[xscale=2,yscale=1.7,baseline={(x1.base)}]
\draw[0cell=.9]
(0,0) node (x1) {\C}
(x1)++(1,0) node (x2) {\D}
;
\draw[1cell=.9]  
(x1) edge[bend left] node {F} (x2)
(x1) edge[bend right] node[swap] {G} (x2)
;
\draw[2cell]
node[between=x1 and x2 at .47, rotate=-90, 2label={above,\theta}] {\Rightarrow}
;
\end{tikzpicture}
\end{equation}
and is called a \index{2-cell!notation}\emph{2-cell}.
\end{definition}

Detailed discussion of pasting diagrams involving 2-cells is in \cite[Ch.\! 3]{johnson-yau}.

\begin{convention}[Left Normalized Bracketing]\label{expl:leftbracketing}\index{bracketing!left normalized}\index{convention!left normalized bracketing}
An iterated monoidal product is \index{left normalized}\emph{left normalized} with the left half of each pair of parentheses at the far left, unless a different bracketing is specified.  An empty monoidal product is the monoidal unit.
\end{convention}

\begin{explanation}\label{expl:lnbracketing}
As an example of \cref{expl:leftbracketing}, we denote
\[w \otimes x \otimes y \otimes z = \big((w \otimes x) \otimes y\big) \otimes z.\]
We omit parentheses for iterated monoidal products and tacitly insert the necessary associativity and unit isomorphisms.  This is justified by Mac Lane's Strictification Theorem \cite[XI.3.1]{maclane}: Each monoidal category $\C$ admits a canonical strong monoidal adjoint equivalence $\C \rightleftarrows \C_\st$ with $\C_\st$ a strict monoidal category.  The analogous braided and symmetric strictification theorems are in \cite[21.3.1 and 21.6.1]{yau-inf-operad}.  In each case, since an equivalence is full and faithful, the strict diagrams commute if and only if their preimages in $\C$ commute.
\end{explanation}

\section{2-Categories}
\label{sec:twocategories}

In this section we review
\begin{itemize}
\item 2-categories, 2-functors, 2-natural transformations, 2-adjunctions, and
\item the 2-categories of small categories, small permutative categories, and small 2-categories.
\end{itemize} 

\begin{definition}\label{def:twocategory}\index{2-category}\index{category!2-}
A \emph{2-category} $\A$ consists of the following data.
\begin{description}
\item[Objects] It is equipped with a class $\A_0$ of \emph{objects}.
\item[1-Cells] For each pair of objects $a,b\in\A_0$, it is equipped with a class 
\[\A_1(a,b)\]
of \index{1-cell}\emph{1-cells} from $a$ to $b$.  Such a 1-cell is denoted $a \to b$.
\item[2-Cells] For 1-cells $f,f' \in \A_1(a,b)$, it is equipped with a set 
\[\A_2(f,f')\] 
of \index{2-cell}\emph{2-cells} from $f$ to $f'$.  Such a 2-cell is denoted $f \to f'$.
\item[Identities]  $\A$ is equipped with
\begin{itemize}
\item an \emph{identity 1-cell} 
\[1_a \in \A_1(a,a)\]
for each object $a$ and
\item an \emph{identity 2-cell} 
\[1_f \in \A_2(f,f)\]
for each 1-cell $f \in \A_1(a,b)$.
\end{itemize}
\item[Compositions]
In the following compositions, $a$, $b$, and $c$ denote objects in $\A$.
\begin{itemize}
\item For 1-cells $f,f',f'' \in \A_1(a,b)$, it is equipped with an assignment\label{not:vcompiicell}
\[\begin{tikzcd}
\A_2(f',f'') \times \A_2(f,f') \rar{v} & \A_2(f,f'')
\end{tikzcd},\qquad v(\theta',\theta) = \theta'\theta\] 
called the \index{vertical composition!2-category}\emph{vertical composition} of 2-cells.
\item It is equipped with an assignment\label{not:hcompicell}
\[\begin{tikzcd}
\A_1(b,c) \times \A_1(a,b) \rar{h_1} & \A_1(a,c)
\end{tikzcd},\qquad h_1(g,f) = gf\]
called the \index{horizontal composition!2-category}\emph{horizontal composition} of 1-cells.
\item For 1-cells $f,f' \in \A_1(a,b)$ and $g,g' \in \A_1(b,c)$, it is equipped with an assignment\label{not:hcompiicell}
\[\begin{tikzcd}
\A_2(g,g') \times \A_2(f,f') \rar{h_2} & \A_2(gf,g'f')
\end{tikzcd},\quad h_2(\phi,\theta) = \phi * \theta\]
called the \emph{horizontal composition} of 2-cells.
\end{itemize}
\end{description}
The data above are required to satisfy \cref{twocat-i,twocat-ii,twocat-iii,twocat-iv} below:
\begin{romenumerate}
\item\label{twocat-i} Vertical composition is associative and unital for identity 2-cells.
\item\label{twocat-ii} Horizontal composition preserves vertical composition and identity 2-cells.
\item\label{twocat-iii} Horizontal composition of 1-cells is associative and unital for identity 1-cells.
\item\label{twocat-iv} Horizontal composition of 2-cells is associative and unital for identity 2-cells of identity 1-cells.
\end{romenumerate}
This finishes the definition of a 2-category.

Moreover, we define the following.
\begin{itemize}
\item For objects $a$ and $b$ in a 2-category $\A$, the \index{hom!- category in a 2-category}\emph{hom category} \label{not:homcat}$\A(a,b)$ is the category defined by the following data.
\begin{itemize}
\item Objects are 1-cells from $a$ to $b$.
\item Morphisms are 2-cells between 1-cells $a \to b$.
\item Composition is vertical composition of 2-cells.
\item Identities are identity 2-cells.
\end{itemize} 
\item A 2-category is \index{locally small}\emph{locally small} if each hom category is small.
\item A 2-category is \emph{small} if it has a set of objects and is locally small.
\item The \index{underlying 1-category}\index{1-category!underlying}\emph{underlying 1-category} of a 2-category is defined as follows.
\begin{itemize}
\item It has the same class of objects.
\item Morphisms are 1-cells.
\item Composition is horizontal composition of 1-cells.
\item Identity morphisms are identity 1-cells.
\end{itemize} 
\end{itemize} 
We also use the 2-cell notation \cref{twocellnotation} for 2-cells in a 2-category.
\end{definition}

\begin{example}[Small Categories]\label{ex:catastwocategory}\index{2-category!of small categories}
The category $\Cat$ in \cref{ex:cat} is the underlying 1-category of a 2-category, in which the 2-cells are natural transformations.  Horizontal, respectively vertical, composition of 2-cells in $\Cat$ is the same as that of natural transformations. 
\end{example}

\begin{definition}[Small Permutative Categories]\label{def:permcat}
Denote by 
\[\permcat\]
the 2-category with
\begin{itemize}
\item small permutative categories as objects,
\item symmetric monoidal functors as 1-cells, and
\item monoidal natural transformations as 2-cells.
\end{itemize}
Moreover, we define the following locally-full sub-2-categories of $\permcat$ with the same objects but restricting the 1-cells.
\begin{itemize}
\item $\permcatsu$ has 1-cells given by strictly unital symmetric monoidal functors.
\item $\permcatsg$ has 1-cells given by strong symmetric monoidal functors.
\item $\permcatsus$ has 1-cells given by strictly unital strong symmetric monoidal functors.
\item $\permcatst$ has 1-cells given by strict symmetric monoidal functors.
\end{itemize}
In each case the 2-cells are given by monoidal natural transformations.
\end{definition}

The extensions of $\permcatsu$, $\permcatsus$, and $\permcatst$ to $\Cat$-multicategories are in \cref{thm:permcatmulticat}.

\begin{definition}\label{def:twofunctor}
For 2-categories $\A$ and $\B$, a \index{2-functor}\index{functor!2-}\emph{2-functor} $F \cn \A \to \B$ consists of the following data.
\begin{description}
\item[Object Assignment] It is equipped with a function
\[\begin{tikzcd}[column sep=large]
\A_0 \ar{r}{F_0} & \B_0.
\end{tikzcd}\]
\item[1-Cell Assignment] For each pair of objects $a,b$ in $\A$, it is equipped with a function
\[\begin{tikzcd}[column sep=large]
\A_1(a,b) \ar{r}{F_1} & \B_1(F_0 a,F_0 b).
\end{tikzcd}\] 
\item[2-Cell Assignment] For each pair of objects $a,b$ in $\A$ and 1-cells $f,f' \in \A_1(a,b)$, it is equipped with a function
\[\begin{tikzcd}[column sep=large]
\A_2(f,f') \ar{r}{F_2} & \B_2(F_1 f,F_1 f').
\end{tikzcd}\] 
\end{description}
The data above are required to satisfy \cref{twofunctor-i,twofunctor-ii,twofunctor-iii} below, with each of $F_0$, $F_1$, and $F_2$ abbreviated to $F$:
\begin{romenumerate}
\item\label{twofunctor-i} The object and 1-cell assignments of $F$ form a functor from the underlying 1-category of $\A$ to the underlying 1-category of $\B$.
\item\label{twofunctor-ii} For each pair of objects $a,b$ in $\A$, the 1-cell and 2-cell assignments of $F$ form a functor between hom categories:
\[\begin{tikzcd}[column sep=large]
\A(a,b) \ar{r}{F} & \B(Fa,Fb).
\end{tikzcd}\] 
\item\label{twofunctor-iii} $F$ preserves horizontal composition of 2-cells.
\end{romenumerate}
This finishes the definition of a 2-functor.

Moreover, we define the following.
\begin{itemize}
\item The \emph{identity 2-functor} $1_{\A} \cn \A \to \A$ is defined by the identity assignments on objects, 1-cells, and 2-cells.
\item Given a 2-functor $G \cn \B \to \C$, the \emph{composite} $GF$ is the 2-functor 
\[\begin{tikzcd}[column sep=large]
\A \ar{r}{GF} & \C
\end{tikzcd}\]
defined by separately composing the assignments on objects, 1-cells, and 2-cells.\defmark
\end{itemize}
\end{definition}

\begin{example}\label{ex:permcattwofunctors}
In the context of \cref{def:permcat}, there are inclusion 2-functors as follows.
\begin{equation}\label{permcatinclusion}
\begin{tikzpicture}[xscale=1,yscale=1,baseline={(a.base)}]
\def\h{2.8} \def\g{1.5} \def\v{1} \def\u{2}
\draw[0cell=.9]
(0,0) node (a) {\permcatst}
(a)++(\h,0) node (b) {\permcatsus}
(b)++(\g,\v) node (c) {\permcatsu}
(b)++(\g,-\v) node (d) {\permcatsg}
(b)++(2*\g,0) node (e) {\permcat}
;
\draw[1cell=.9]  
(a) edge node {} (b)
(b) edge node {} (c)
(b) edge node {} (d)
(c) edge node {} (e)
(d) edge node {} (e)
;
\end{tikzpicture}
\end{equation}
Each 2-functor in \cref{permcatinclusion} is the identity on
\begin{itemize}
\item objects, which are small permutative categories, and
\item 2-cells between each pair of 1-cells, which are monoidal natural transformations.\defmark
\end{itemize} 
\end{example}

\begin{definition}\label{def:twonaturaltr}\index{2-natural!transformation}\index{natural transformation!2-}
For 2-functors $F,G \cn \A\to\B$ between 2-categories $\A$ and $\B$, a \emph{2-natural transformation} $\vphi \cn F \to G$ consists of, for each object $a$ in $\A$, a \emph{component 1-cell} 
\[\begin{tikzcd}[column sep=large]
Fa \ar{r}{\vphi_a} & Ga \inspace \B
\end{tikzcd}\]
such that the axioms \cref{onecellnaturality,twocellnaturality} below hold:
\begin{description}
\item[1-Cell Naturality] For each 1-cell $f \cn a \to b$ in $\A$, the following two composite 1-cells in $\B(Fa,Gb)$ are equal.
\begin{equation}\label{onecellnaturality}
\begin{tikzcd}[column sep=large]
Fa \ar{d}[swap]{Ff} \ar{r}{\vphi_a} & Ga \ar{d}{Gf}\\
Fb \ar{r}{\vphi_b} & Gb
\end{tikzcd}
\end{equation}
\item[2-Cell Naturality] 
For each 2-cell $\theta \cn f \to g$ in $\A(a,b)$, the two whiskered 2-cells in the following diagram are equal.
\begin{equation}\label{twocellnaturality}
\begin{tikzpicture}[xscale=3,yscale=1.4,vcenter]
\def\a{30} \def\s{.85}
\draw[0cell=\s]
(0,0) node (x11) {Fa}
(x11)++(1,0) node (x12) {Ga}
(x11)++(0,-1) node (x21) {Fb}
(x12)++(0,-1) node (x22) {Gb}
;
\draw[1cell=\s]  
(x11) edge node {\vphi_a} (x12)
(x21) edge node {\vphi_b} (x22)
(x11) edge[bend right] node[swap] {Ff} (x21)
(x11) edge[bend left] node {Fg} (x21)
(x12) edge[bend right] node[swap] {Gf} (x22)
(x12) edge[bend left] node {Gg} (x22)
;
\draw[2cell=.9]
node[between=x11 and x21 at .6, 2label={above,F\theta}] {\Rightarrow}
node[between=x12 and x22 at .6, 2label={above,G\theta}] {\Rightarrow}
;
\end{tikzpicture}
\end{equation}
The axiom \cref{twocellnaturality} means the following 2-cell equality in $\B(Fa,Gb)$:
\[G\theta * 1_{\vphi_a} = 1_{\vphi_b} * F\theta.\]
\end{description}
This finishes the definition of a 2-natural transformation. 

Moreover, we define the following.
\begin{itemize}
\item We extend the 2-cell notation \cref{twocellnotation} to 2-natural transformations. 
\item A \index{2-natural!isomorphism}\index{natural isomorphism!2-}\emph{2-natural isomorphism} is a 2-natural transformation such that each component 1-cell is an isomorphism in the underlying 1-category.
\item Identity 2-natural transformations, horizontal composition of 2-natural transformations, and vertical composition of 2-natural transformations are defined componentwise.\defmark
\end{itemize}
\end{definition}

\begin{example}[Small 2-Categories]\label{ex:iicat}
Analogous to \cref{ex:catastwocategory}, there is a 2-category $\iicat$ with
\begin{itemize}
\item small 2-categories as objects,
\item 2-functors as 1-cells, and
\item 2-natural transformations as 2-cells.
\end{itemize}
Horizontal, respectively vertical, composition of 2-cells in $\iicat$ is given by that of 2-natural transformations.  We also use the notation $\iicat$ for its underlying 1-category.
\end{example}

\begin{definition}\label{def:twoadjunction}
Suppose $\A$ and $\B$ are 2-categories.  A \index{2-adjunction}\index{adjunction!2-}\emph{2-adjunction} from $\A$ to $\B$ is a quadruple
\[(F,G,\eta,\epz) \cn \A \to \B\]
consisting of the following data.
\begin{itemize}
\item $F \cn \A \to \B$ is a 2-functor, which is called the \emph{left adjoint}.
\item $G \cn \B \to \A$ is a 2-functor, which is called the \emph{right adjoint}.
\item $\eta \cn 1_{\A} \to GF$ is a 2-natural transformation, which is called the \emph{unit}.
\item $\epz \cn FG \to 1_{\B}$ is a 2-natural transformation, which is called the \emph{counit}.
\end{itemize}
The data above are required to make the following two diagrams commute, where $*$ denotes horizontal composition of 2-natural transformations.
\begin{equation}\label{triangleidentities}
\begin{tikzpicture}[xscale=1,yscale=1.3,vcenter]
\def\s{.9} \def\h{1.8}
\draw[0cell=\s]
(0,0) node (x11) {F}
(x11)++(\h,0) node (x12) {FGF}
(x12)++(0,-1) node (x22) {F}
(x12)++(\h,0) node (y11) {G}
(y11)++(\h,0) node (y12) {GFG}
(y12)++(0,-1) node (y22) {G}
;
\draw[1cell=\s]  
(x11) edge node {1_F * \eta} (x12)
(x12) edge node {\epz * 1_F} (x22)
(x11) edge node[swap] {1_F} (x22)
(y11) edge node {\eta * 1_G} (y12)
(y12) edge node {1_G * \epz} (y22)
(y11) edge node[swap] {1_G} (y22)
;
\end{tikzpicture}
\end{equation}
The two diagrams in \cref{triangleidentities} are called, respectively, the \index{triangle identities}\emph{left triangle identity} and the \emph{right triangle identity}.  Such a 2-adjunction is also denoted $F \dashv G$.
\end{definition}

\chapter{Enriched Category Theory}
\label{ch:prelim_enriched}
In this appendix we review elements of enriched category theory.  The following table summarizes the main content in this appendix.
\smallskip
\begin{center}
\resizebox{.8\width}{!}{
{\renewcommand{\arraystretch}{1.4}%
{\setlength{\tabcolsep}{1ex}
\begin{tabular}{|c|c|}\hline
\multicolumn{2}{|c|}{\appnumname{sec:enrichedcat}} \\ \hline
enriched categories, functors, and natural transformations & \ref{def:enriched-category}, \ref{def:enriched-functor}, and \ref{def:enriched-natural-transformation} \\ \hline
unit and opposite $\V$-categories & \ref{definition:unit-vcat} and \ref{definition:vcat-op} \\ \hline
$\Cat$-categories as locally small 2-categories & \ref{locallysmalltwocat} \\ \hline
2-category $\VCat$ & \ref{ex:vcatastwocategory} \\ \hline\hline
\multicolumn{2}{|c|}{\appnumname{sec:enrmonoidalcat}} \\ \hline
tensor products and monoidal category of $\V$-categories & \ref{definition:vtensor-0} and \ref{theorem:vcat-mon} \\ \hline
monoidal $\V$-categories (braided, symmetric) & \ref{definition:monoidal-vcat} (\ref{definition:braided-monoidal-vcat}, \ref{definition:symm-monoidal-vcat}) \\
\hline
\end{tabular}}}}
\end{center}
\smallskip
The material in this chapter is adapted from \cite{cerberusIII}; see also \cite{kelly-enriched}.  We remind the reader of \cref{conv:universe,expl:leftbracketing}.

\section{Enriched Categories}
\label{sec:enrichedcat}

In this section we review enriched variants of categories, functors, natural transformations, and opposite categories.  Suppose 
\[(\V,\otimes,\tu,\alpha,\lambda,\rho)\]
is a monoidal category (\cref{def:monoidalcategory}).  A braiding on $\V$ is not needed until \cref{definition:vcat-op} of the opposite $\V$-category.  The material in this section is adapted from \cite[Sections 1.1 and 1.2]{cerberusIII}.

\begin{definition}\label{def:enriched-category}
A \emph{$\V$-category} $\C$, which is also called a \index{category!enriched}\index{enriched category}\emph{category enriched in $\V$}, consists of the following data.
\begin{description}
\item[Objects] It is equipped with a class \index{object!enriched category}$\Ob(\C)$, whose elements are called \emph{objects}.
\item[Hom Objects] Each pair of objects $x,y$ in $\C$ is equipped with a \index{hom object}\index{object!hom}\emph{hom object} 
\[\C(x,y) \in \V.\]
\item[Composition] For objects $x,y,z$ in $\C$, it is equipped with a morphism in $\V$\label{not:enrcomposition}
\begin{equation}\label{enr-composition}
\begin{tikzcd}[column sep=huge] 
\C(y,z) \otimes \C(x,y) \rar{\mcomp_{x,y,z}} & \C(x,z)
\end{tikzcd}
\end{equation}
called the \index{composition!enriched category}\emph{composition}.
\item[Identities] Each object $x$ in $\C$ is equipped with a morphism in $\V$\label{not:enridentity}
\begin{equation}\label{enr-identity}
\begin{tikzcd}[column sep=large] 
\tensorunit \rar{i_x} & \C(x,x)
\end{tikzcd}
\end{equation}
called the \emph{identity} of $x$.
\end{description}
The data above are required to make the following \index{associativity!enriched category}\emph{associativity diagram} and \index{unity!enriched category}\emph{unity diagram} commute for objects $w,x,y,z$ in $\C$.
\begin{equation}\label{enriched-cat-associativity}
\begin{tikzcd}[cells={nodes={scale=.9}}]
  \big(\C(y,z) \otimes \C(x,y)\big) \otimes \C(w,x) \arrow{dd}[swap]{\mcomp \otimes 1}
  \rar{\alpha}
  & \C(y,z) \otimes \big(\C(x,y) \otimes \C(w,x)\big)  \dar{1 \otimes \mcomp} \\
& \C(y,z) \otimes \C(w,y) \dar{\mcomp}\\
\C(x,z) \otimes \C(w,x) \rar{\mcomp} & \C(w,z)  
\end{tikzcd}
\end{equation}
\begin{equation}\label{enriched-cat-unity}
\begin{tikzcd}[cells={nodes={scale=.9}}]
\tensorunit \otimes \C(x,y) \dar[swap]{i_y \otimes 1} \arrow[bend left=20]{dr}{\lambda} & & \C(x,y) \otimes \tensorunit \arrow[bend right=20]{dl}[swap]{\rho} \dar{1 \otimes i_x}\\
\C(y,y) \otimes \C(x,y) \rar{\mcomp} & \C(x,y) & \C(x,y) \otimes \C(x,x) \lar[swap]{\mcomp}
\end{tikzcd}
\end{equation}
This finishes the definition of a $\V$-category.  Moreover, a $\V$-category $\C$ is \emph{small}\index{small!enriched category}\index{enriched category!small} if $\Ob(\C)$ is a set.
\end{definition}

\begin{definition}\label{definition:unit-vcat}
The \emph{unit $\V$-category},\index{unit enriched category}\index{enriched category!unit}\index{tensor product!enriched category!unit}\index{enriched!tensor product!unit} $\vtensorunit$, is the one-object $\V$-category whose unique hom object is the monoidal unit, $\tensorunit$, of $\V$.  The composition and identity structure morphisms are given, respectively, by the left unit isomorphism $\lambda_{\tu}$ and the identity morphism $1_{\tensorunit}$.
\end{definition}

\begin{proposition}\label{locallysmalltwocat}
Regarding $(\Cat, \times, \boldone)$ as a monoidal category, a locally small 2-category is precisely a $\Cat$-category.
\end{proposition}

\begin{definition}\label{def:enriched-functor}
A \index{functor!enriched}\index{enriched!functor}\emph{$\V$-functor} 
\[F \cn \C \to \D\]
between $\V$-categories $\C$ and $\D$ consists of the following data.
\begin{description}
\item[Object Assignment] It is equipped with a function 
\[F \cn \Ob(\C) \to \Ob(\D).\]
\item[Component Morphisms] It is equipped with, for each pair of objects $x,y$ in $\C$, a morphism
\[\begin{tikzcd}[column sep=large]
\C(x,y) \rar{F_{x,y}} & \D\bigl(Fx,Fy\bigr)\end{tikzcd} \inspace \V.\] 
\end{description}
The data above is required to make the following two diagrams commute for objects $x,y,z$ in $\C$.\index{composition!enriched functor}
\begin{equation}\label{eq:enriched-composition}
\begin{tikzcd}[cells={nodes={scale=.9}}]
\C(y,z) \otimes \C(x,y) \rar{\mcomp} \dar[swap]{F \otimes F} & \C(x,z) \dar{F}\\
\D(Fy,Fz) \otimes \D(Fx,Fy) \rar{\mcomp} & \D(Fx,Fz)\end{tikzcd}
\qquad
\begin{tikzcd}[cells={nodes={scale=.9}}]
\tensorunit \arrow{dr}[swap]{i_{Fx}} \rar{i_x} & \C(x,x) \dar{F}\\
& \D(Fx,Fx)
\end{tikzcd}
\end{equation}
Moreover, we define the following.
\begin{itemize}
\item An \index{identity!enriched functor}identity $\V$-functor is given by the identity object assignment and identity component morphisms.
\item Composition of $\V$-functors is defined by separately composing the object assignments and the component morphisms.\defmark
\end{itemize}
\end{definition}

\begin{definition}\label{def:enriched-natural-transformation}
For $\V$-functors $F,G \cn \C\to\D$ between $\V$-categories $\C$ and $\D$, a \index{enriched!natural transformation}\index{natural transformation!enriched}\emph{$\V$-natural transformation} $\theta \cn F\to G$ consists of, for each object $x$ in $\C$, a morphism
\[\begin{tikzcd}[column sep=large]
\tensorunit \ar{r}{\theta_x} & \D(Fx,Gx)
\end{tikzcd} \inspace \V,\]
called the \index{component}\emph{$x$-component} of $\theta$.  The following \index{naturality diagram}\emph{naturality diagram} in $\V$ is required to commute for objects $x,y$ in $\C$.
\begin{equation}\label{enr-naturality}
\begin{tikzpicture}[x=40mm,y=12mm,vcenter]
\def\s{.85}
  \draw[0cell=\s] 
  (0,0) node (a) {
    \C(x,y)
  }
  (.25,1) node (b1) {
    \tensorunit \otimes \C(x,y)    
  }
  (b1)++(1,0) node (c1) {
    \D(Fy,Gy) \otimes \D(Fx,Fy)
  }
  (.25,-1) node (b2) {
    \C(x,y) \otimes \tensorunit
  }
  (b2)++(1,0) node (c2) {
    \D(Gx,Gy) \otimes \D(Fx,Gx)
  }
  (c1)++(.25,-1) node (d) {
    \D(Fx,Gy)
  }
  ;
  \draw[1cell=\s] 
  (a) edge node[pos=.3] {\la^\inv} node['] {\iso}(b1)
  (a) edge node[swap,pos=.2] {\rho^\inv} node {\iso} (b2)
  (b1) edge node {\theta_y \otimes F} (c1)
  (b2) edge node {G \otimes \theta_x} (c2)
  (c1) edge node[pos=.6] {\mcomp} (d)
  (c2) edge node[swap,pos=.7] {\mcomp} (d)
  ;
\end{tikzpicture}
\end{equation}
Moreover, each of the following notions is defined componentwise:
\begin{itemize}
\item \index{identity!enriched natural transformation}identity $\V$-natural transformations,
\item \index{horizontal composition!enriched natural transformation}horizontal composition of $\V$-natural transformations, and
\item \index{vertical composition!enriched natural transformation}vertical composition of $\V$-natural transformations.
\end{itemize}
We use the 2-cell notation \cref{twocellnotation} for $\V$-natural transformations.
\end{definition}

\begin{example}[Small Enriched Categories]\label{ex:vcatastwocategory}\index{2-category!of small enriched categories}
Each monoidal category $\V$ has an associated 2-category $\VCat$ defined by the following data.
\begin{itemize}
\item Objects are small $\V$-categories.
\item 1-cells are $\V$-functors.
\item 2-cells are $\V$-natural transformations.
\end{itemize}
The 2-category $\Cat$ in \cref{ex:catastwocategory} is the special case for $\V = (\Set, \times, *)$.
\end{example}

For the following definition, we assume that $\V$ is a braided monoidal category (\cref{def:braidedmoncat}).  The next definition is \cite[1.2.16]{cerberusIII}.

\begin{definition}\label{definition:vcat-op}
  Suppose $\C$ is a $\V$-category with $(\V,\otimes,\tu,\xi)$ a braided monoidal category.  The \index{opposite!enriched category}\index{enriched category!opposite}\emph{opposite $\V$-category}, $\Cop$, is defined as follows.
\begin{description}
\item[Objects] $\Cop$ has the same class of objects as $\C$.
\item[Hom Objects] Each pair of objects $x,y$ in $\Cop$ is equipped with the hom object
\[\Cop(x,y) = \C(y,x).\]
\item[Composition] The composition in $\Cop$ is defined for each triple of objects $x,y,z$ in $\Cop$ as the following composite in $\V$ using the braiding $\xi$ of $\V$ and the composition $\mcomp$ of $\C$:
\begin{equation}\label{opposite-comp}
  \begin{tikzpicture}[x=28mm,y=8mm,vcenter]
    \draw[0cell] 
    (0,0) node (a) {
      \Cop(y,z) \otimes \Cop(x,y)
    }
    (0,-1) node (b) {
      \C(z,y) \otimes \C(y,x)
    }
    (1.3,-1) node (c) {
      \C(y,x) \otimes \C(z,y)
    }
    (2.3,-1) node (d) {
      \C(z,x)
    }
    (2.3,0) node (e) {
      \Cop(x,z)
    }
    ;
    \draw[1cell] 
    (a) edge[equal] node {} (b)
    (b) edge node {\xi} (c)
    (c) edge node {\mcomp} (d)
    (d) edge[equal] node {} (e)
    ;
  \end{tikzpicture}
\end{equation}
\item[Identities] The identity of each object $x$ in $\Cop$ is the same as the identity of $x$ in $\C$:
  \[\tensorunit \fto{i_x} \C(x,x) = \Cop(x,x).\]
\end{description} 
This finishes the definition of the opposite $\V$-category.  

Moreover, we extend the opposite construction to $\V$-functors and $\V$-natural transformations as follows.
\begin{itemize}
\item For a $\V$-functor $F \cn \C \to \D$, the \index{opposite!enriched functor}\emph{opposite $\V$-functor}
\[\Fop \cn \Cop \to \Dop\]
has
\begin{itemize}
\item the same object assignment as $F$ and
\item the following component morphisms for objects $x,y$ in $\Cop$:
\[\Fop_{x,y} = F_{y,x} \cn \Cop(x,y) \to \Dop(Fx,Fy).\]
\end{itemize} 
\item For a $\V$-natural transformation $\theta \cn F \to G$ with $F,G \cn \C \to \D$ both $\V$-functors, the \index{opposite!enriched natural transformation}\emph{opposite $\V$-natural transformation} 
\[\thetaop \cn \Gop \to \Fop\]
is defined by the component morphisms
\[\thetaop_x = \theta_x \cn \tu \to \Dop(Gx,Fx)\]
for objects $x$ in $\Cop$.\defmark
\end{itemize}
\end{definition}

\section{Enriched Monoidal Categories}
\label{sec:enrmonoidalcat}

In this section we review monoidal categories, functors, and natural transformations enriched in a braided monoidal category $(\V,\otimes,\tu,\alpha,\lambda,\rho,\xi)$ (\cref{def:braidedmoncat}).  Whenever we need $\V$ to be symmetric monoidal, we state so explicitly.  The material in this section is adapted from \cite[Chapter 1]{cerberusIII}.

\subsection*{Tensor Product of Enriched Categories}

\begin{definition}\label{definition:vtensor-0}
Suppose $\C$ and $\D$ are $\V$-categories.  The \emph{tensor product}\index{tensor product!enriched category}\index{enriched category!tensor product}, $\C \otimes \D$, is the $\V$-category defined by the following data.
\begin{description}
\item[Objects] Its class of objects is 
\[\Ob(\C \otimes \D) = \Ob(\C) \times \Ob(\D).\]  
Objects in $\C \otimes \D$ are denoted $x \otimes y$ for $x \in \C$ and $y \in \D$.
\item[Hom Objects] For objects $x \otimes y$ and $x' \otimes y'$, the hom object is the monoidal product
\[(\C \otimes \D)(x \otimes y, x' \otimes y') = \C(x,x') \otimes \D(y,y') \inspace \V.\]
\item[Composition] For objects $x \otimes y$, $x'\otimes y'$, and $x''\otimes y''$, the composition is the following composite in $\V$, where $\ximid$ interchanges the middle two factors using the associativity isomorphism and braiding in $\V$.\label{not:ximid}
\[\begin{tikzpicture}[x=30mm,y=15mm]
\draw[0cell=.8]
(0,0) node (a) {\big(\C(x',x'')\otimes\D(y',y'')\big) \otimes \big(\C(x,x')\otimes\D(y,y')\big)}
(1,-1) node (b) {\big(\C(x',x'')\otimes\C(x,x')\big) \otimes \big(\D(y',y'')\otimes\D(y,y')\big)}
(2,0) node (c) {\C(x,x'') \otimes \D(y,y'')}
;
\draw[1cell=.9]
(a) edge[transform canvas={xshift=-1em}] node[swap,pos=.2] {\ximid} (b)
(b) edge[transform canvas={xshift=1em}] node[swap,pos=.7] {\mcomp \otimes \mcomp} (c);
\end{tikzpicture}\]
\item[Identities] The identity of an object $x \otimes y$ is the following composite in $\V$.
    \[\tensorunit \fto[\iso]{\la^\inv} \tensorunit \otimes \tensorunit \fto{i_x \otimes i_y} \C(x,x)
    \otimes \D(y,y) = (\C \otimes \D)(x \otimes y, x \otimes y)\]
\end{description}
This finishes the definition of the $\V$-category $\C \otimes \D$.  The tensor product $\otimes$ extends to $\V$-functors and $\V$-natural transformations componentwise.
\end{definition}

Recall from \cref{ex:vcatastwocategory} the 2-category $\VCat$ of small $\V$-categories, $\V$-functors, and $\V$-natural transformations.

\begin{proposition}\label{vtensoriifunctor}
The tensor product is a 2-functor
\[\begin{tikzcd}[column sep=large]
\VCat \times \VCat \ar{r}{\otimes} & \VCat.
\end{tikzcd}\]
\end{proposition}

Recall the unit $\V$-category $\vtu$ in \cref{definition:unit-vcat}.  The tensor product on $\VCat$ is part of a monoidal structure, with the following unit and associativity isomorphisms.  

\begin{definition}\label{definition:vtensor-unitors}
We define the \emph{left unitor}\index{left unitor!enriched tensor product}\index{tensor product!enriched category!left unitor}\index{enriched category!tensor product!left unitor} $\ell^\otimes$ and the \emph{right unitor}\index{right unitor!enriched tensor product}\index{tensor product!enriched category!right unitor}\index{enriched category!tensor product!right unitor} $r^\otimes$ as the 2-natural isomorphisms
  \[
  \begin{tikzpicture}[x=25mm,y=13mm,vcenter]
    \draw[0cell=.85]
    (0,1) node (a) {\VCat^2}
    (-.5,0) node (b) {\VCat}
    (.5,0) node (b') {\VCat}
    ;
    \draw[1cell=.9] 
    (b) edge node {\vtensorunit \times 1} (a)
    (a) edge node {\otimes} (b')
    (b) edge['] node (z) {1} (b')
    ;
    \draw[2cell] 
    node[between=a and z at .55, shift={(-.05,0)}, rotate=-90, 2label={above,\ell^\otimes}] {\Rightarrow}
    ;
  \end{tikzpicture}
  \qquad \text{and} \qquad
  \begin{tikzpicture}[x=25mm,y=13mm,vcenter]
    \draw[0cell=.85] 
    (0,1) node (a) {\VCat^2}
    (-.5,0) node (b) {\VCat}
    (.5,0) node (b') {\VCat}
    ;
    \draw[1cell=.9] 
    (b) edge node {1 \times \vtensorunit} (a)
    (a) edge node {\otimes} (b')
    (b) edge['] node (z) {1} (b')
    ;
    \draw[2cell] 
    node[between=a and z at .55, shift={(-.05,0)}, rotate=-90, 2label={above,r^\otimes}] {\Rightarrow}
    ;
  \end{tikzpicture}
  \]
as follows.  The unitor
  components at a $\V$-category $\C$ are the $\V$-functors
  \[
  \vtensorunit \otimes \C \fto{\ell^\otimes_{\C}} \C \xleftarrow{r^\otimes_{\C}} \C
  \otimes \vtensorunit
  \]
given
\begin{itemize}
\item on objects by the unitors for the Cartesian product and
\item on hom objects by the unit isomorphisms
  \[
  \tensorunit \otimes \C(x,x') \fto{\la} \C(x,x') \xleftarrow{\rho}
  \C(x,x') \otimes \tensorunit \inspace \V
  \]
for objects $x,x' \in \C$.\defmark
\end{itemize} 
\end{definition}

\begin{definition}\label{definition:vtensor-assoc}
We define the \emph{associator}\index{associator!enriched tensor product}\index{tensor product!enriched category!associator}\index{enriched!tensor product!associator} $a^\otimes$ as the 2-natural isomorphism
  \[
  \begin{tikzpicture}[x=25mm,y=15mm]
    \draw[0cell=.85] 
    (0,0) node (a) {\VCat^3}
    (1,0) node (b) {\VCat^2}
    (0,-1) node (c) {\VCat^2}
    (1,-1) node (d) {\VCat}
    ;
    \draw[1cell=.9] 
    (a) edge node {\otimes \times 1} (b)
    (c) edge['] node {\otimes} (d)
    (a) edge['] node {1 \times \otimes} (c)
    (b) edge node {\otimes} (d)
    ;
    \draw[2cell] 
    node[between=a and d at .5, rotate=225, 2label={below,a^\otimes}] {\Rightarrow}
    ;
  \end{tikzpicture}
  \]
as follows.  For small $\V$-categories $\C$, $\D$, and $\E$, the associator component is the $\V$-functor
  \[\begin{tikzcd}[column sep=huge]
  (\C \otimes \D) \otimes \E \ar{r}{a^\otimes_{\C,\D,\E}} & \C \otimes (\D \otimes \E)
  \end{tikzcd}\]
given
\begin{itemize}
\item on objects by the associativity isomorphism of the Cartesian product and
\item on hom objects by the associativity isomorphism
  \[\begin{tikzcd}
  \big(\C(x,x') \otimes \D(y,y')\big) \otimes \E(z,z') \ar{r}{\al} & \C(x,x') \otimes \big(\D(y,y') \otimes \E(z,z')\big)
  \end{tikzcd}\]
in $\V$ for objects $x,x' \in \C$, $y,y' \in \D$, and $z,z'\in \E$.\defmark
\end{itemize} 
\end{definition}

The Cartesian product on objects and the monoidal structure of $\V$ both satisfy the unity and pentagon axioms in \cref{monoidalcataxioms,lambda=rho,moncat-other-unit-axioms}.  Thus the data $(\otimes, \vtu, a^\otimes, \ell^\otimes, r^\otimes)$ in \cref{definition:unit-vcat,definition:vtensor-0,definition:vtensor-unitors,definition:vtensor-assoc} also satisfy these axioms.

\begin{definition}\label{definition:vtensor-beta}
Suppose that $(\V,\xi)$ is a symmetric monoidal category.  We define the
  \emph{braiding}\index{braiding!enriched tensor product}\index{tensor product!enriched category!braiding}\index{enriched!tensor product!braiding} $\beta^\otimes$ as the 2-natural isomorphism
  \[
  \begin{tikzpicture}[xscale=3, yscale=1.1]
    \def\v{-1} \def\h{1} \def\m{1} \def\q{15}
    \draw[0cell=.85] 
    (0,0) node (x11) {\VCat^2}
    ($(x11)+(\h,0)$) node (x12) {\VCat}
    ($(x11)+(\h/2,\v)$) node (x2) {\VCat^2}
    ;
    \draw[1cell] 
    (x11) edge node (i) {\otimes} (x12)
    (x11) edge[bend right=\q] node[swap,pos=.5] (a) {\tau} (x2)
    (x2) edge[bend right=\q] node[swap,pos=.5] {\otimes} (x12)
    ;
    \draw[2cell] 
    node[between=i and x2 at .5, shift={(-.1,0)}, rotate=-90, 2label={above,\beta^\otimes}] {\Rightarrow}
    ;
  \end{tikzpicture}
  \]
as follows, with $\tau$ permuting the two arguments.  For small $\V$-categories $\C$ and $\D$, the braiding component is the $\V$-functor
\[\begin{tikzcd}[column sep=large]
\C \otimes \D \ar{r}{\beta^\otimes_{\C,\D}} & \D \otimes \C
\end{tikzcd}\]
given
\begin{itemize}
\item on objects by the braiding of the Cartesian product and
\item on hom objects by the braiding
\[\begin{tikzcd}[column sep=large]
\C(x,x') \otimes \D(y,y') \ar{r}{\xi} & \D(y,y') \otimes \C(x,x')
\end{tikzcd}\]
in $\V$ for objects $x,x' \in \C$ and $y,y' \in \D$.\defmark
\end{itemize} 
\end{definition}

The Cartesian product on objects and the symmetric monoidal structure of $\V$ both satisfy the hexagon, unity, and symmetry axioms in \cref{hexagon-braided,braidedunity,symmoncatsymhexagon}.   Thus the data $(\otimes, \vtu, a^\otimes, \ell^\otimes, r^\otimes,\beta^\otimes)$ also satisfy these axioms.

\begin{theorem}\label{theorem:vcat-mon}
Suppose $\V = (\V,\otimes,\xi)$ is a braided monoidal category.  Then\index{enriched category!2-category}\index{2-category!of small enriched categories} 
\[\big( \VCat,\otimes,\vtensorunit,a^\otimes,\ell^\otimes,r^\otimes \big)\]
is a monoidal category.\index{monoidal category!of small enriched categories}\index{enriched category!monoidal}\index{tensor product!enriched category!monoidal}  If $\V$ is symmetric monoidal, then so is $(\VCat,\beta^\otimes)$.
\end{theorem}

\subsection*{The Monoidal $\Cat$-Category of Enriched Categories}

\cref{theorem:vcat-mon} is improved to (symmetric) monoidal $\Cat$-categories in  \cref{theorem:vcat-cat-mon} below.  To make the necessary definitions, we
\begin{itemize}
\item abbreviate the tensor product of $\V$-categories to juxtaposition and
\item use superscript with a minus sign, \label{not:minusotimes}$?^{-\otimes}$, to denote the inverse $(?^\otimes)^{-1}$.
\end{itemize} 

\begin{definition}\label{definition:monoidal-vcat}
  A \emph{monoidal $\V$-category}\index{enriched category!monoidal}\index{enriched!monoidal - category}\index{category!enriched!monoidal}\index{monoidal category!enriched} is a sextuple
\[\big( \K,\vmtimes,\vmunitbox,a^\vmtimes,\ell^\vmtimes,r^\vmtimes \big)\]
  consisting of the following data.
  \begin{description}
  \item[Base $\V$-category] It has a $\V$-category $\K$, called the
    \emph{base $\V$-category}.\index{base enriched category!monoidal enriched category}
  \item[Monoidal Composition] It has a $\V$-functor
    \[
    \K \otimes \K \fto{\vmtimes} \K
    \]
    called the \emph{monoidal composition}.\index{monoidal composition!monoidal enriched category}\index{composition!monoidal enriched category}
  \item[Monoidal Identity] It has a $\V$-functor
    \[
    \vtensorunit \fto{\vmunitbox} \K
    \]
    called the \emph{monoidal identity}.\index{monoidal identity!monoidal enriched category}  The image of the unique
    object in $\vtensorunit$ is also denoted $\vmunitbox$ and called the
    \emph{identity object}.\index{identity object!monoidal enriched category}
  \item[Monoidal Unitors] It has $\V$-natural isomorphisms
    \begin{equation}\label{eq:monoidal-unitors}
    \begin{tikzpicture}[x=40mm,y=13mm,vcenter]
    \def\s{.9}
      \draw[0cell=\s] 
      (0,0) node (a) {\K}
      (.25,1) node (b) {\vtensorunit\K}
      (.75,1) node (c) {\K^2}
      (1,0) node (d) {\K}
      ;
      \draw[1cell=\s] 
      (a) edge node {\ell^{-\otimes}} (b)
      (b) edge node (z) {\vmunitbox 1_{\K}} (c)
      (c) edge node {\vmtimes} (d)
      (a) edge['] node (w) {1_{\K}} (d)
      ;
      \draw[2cell] 
      node[between=z and w at .5, rotate=-90, 2label={above,\ell^\vmtimes}] {\Rightarrow}
      ;
    \end{tikzpicture}
    \qquad
    \begin{tikzpicture}[x=40mm,y=13mm,vcenter]
    \def\s{.9}
      \draw[0cell=\s] 
      (0,0) node (a) {\K}
      (.25,1) node (b) {\K\vtensorunit}
      (.75,1) node (c) {\K^2}
      (1,0) node (d) {\K}
      ;
      \draw[1cell=\s] 
      (a) edge node {r^{-\otimes}} (b)
      (b) edge node (z) {1_{\K} \vmunitbox} (c)
      (c) edge node {\vmtimes} (d)
      (a) edge['] node (w) {1_{\K}} (d)
      ;
      \draw[2cell] 
      node[between=z and w at .5, rotate=-90, 2label={above,\,r^\vmtimes}] {\Rightarrow}
      ;
    \end{tikzpicture}
    \end{equation}
called the \emph{left monoidal unitor}\index{left monoidal unitor!monoidal enriched category}\index{unitor!monoidal enriched category} and the \emph{right monoidal unitor},\index{right monoidal unitor!monoidal enriched category} respectively.  Their components at an object $x \in \K$ are, respectively,
    \[
    \tensorunit \fto{\ell^\gmtimes_x} \K(\vmunitbox \gmtimes x, x) \andspace
    \tensorunit \fto{r^\gmtimes_x} \K(x \gmtimes \vmunitbox, x).
    \]
  \item[Monoidal Associator] It has a $\V$-natural isomorphism
    \begin{equation}\label{eq:monoidal-assoc}
    \begin{tikzpicture}[x=20mm,y=20mm,vcenter]
    \def\s{.9}
      \draw[0cell=\s] 
      (0,0) node (a) {(\K^{2})\K}
      (225:.707) node (a') {\K(\K^{2})}
      (a)++(1,0) node (b) {\K^{2}}
      (a')++(0,-.6) node (c) {\K^{2}}
      (b)++(0,-1.1) node (d) {\K}
      ;
      \draw[1cell=\s] 
      (a) edge node {\vmtimes 1_\K} (b)
      (c) edge['] node {\vmtimes} (d)
      (a') edge['] node {1_\K \vmtimes} (c)
      (b) edge node {\vmtimes} (d)
      (a) edge['] node {a^\otimes} (a')
      ;
      \draw[2cell] 
      node[between=b and c at .5, rotate=225, 2label={below,a^\vmtimes}] {\Rightarrow}
      ;
    \end{tikzpicture}
    \end{equation}
called the \emph{monoidal
      associator}.\index{associator!monoidal enriched category}\index{monoidal associator!monoidal enriched category}  Its component at a triple of objects $x,y,z \in
    \K$ is a morphism in $\V$
    \[
    \tensorunit \fto{a^\vmtimes_{x,y,z}}
    \K\big((x \vmtimes y) \vmtimes z \scs x \vmtimes (y \vmtimes z)\big).
    \]
  \end{description}
  These data are required to satisfy the following two axioms, with $1$ denoting the identity $\V$-functor.
  \begin{description}
  \item[Unity Axiom]\index{unity!monoidal enriched category}
    The composites of the following two \index{unity!monoidal enriched category}\emph{middle
      unity pasting diagrams}\index{middle unity diagram!monoidal enriched category} are equal.
    \begin{equation}\label{vmonoidal-unit}

    \end{equation}
    In the first diagram in \cref{vmonoidal-unit}, the unlabeled rectangle commutes by
    naturality of $a^\otimes$.  The
    region labeled $\medstar$ commutes by
    the middle unity for $\ell^\otimes$ and $r^\otimes$.
  \item[Pentagon Axiom]\index{pentagon axiom!monoidal enriched category}
    The composites of the following two \emph{pentagon pasting
      diagrams} are equal.
    \begin{equation}\label{vmonoidal-pentagon-axiom}
      %
    \end{equation}
    The central square in the first diagram in \cref{vmonoidal-pentagon-axiom} commutes by
    2-functoriality of $\otimes$ in each variable.  The other unmarked
    quadrilaterals in the two diagrams in \cref{vmonoidal-pentagon-axiom} commute by 2-naturality
    of $a^\otimes$.  The pentagon labeled
    $\medstar$ commutes by the pentagon axiom for $a^\otimes$.
  \end{description}
  This finishes the definition of a monoidal $\V$-category.  
\end{definition}

For the definition of a braided monoidal $\V$-category, we use a
mate of $a^\vmtimes$ similar to the mates of a pentagonator\index{pentagonator}\index{pentagonator!mate}\index{mate!pentagonator} in \cite[12.1.4]{johnson-yau}.

\begin{definition}\label{eq:a-vmtimes-inv-mate}\index{mate!associator}\index{associator!mate}
For a monoidal $\V$-category $\K$, we denote by $a^\vmtimes_1$ the mate of $a^\vmtimes$ given by the inverse of $a^\otimes$, as shown below.
  
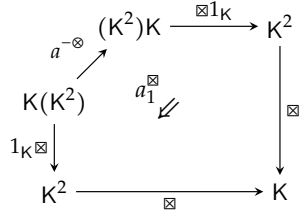
\[
  \begin{tikzpicture}[x=20mm,y=20mm]
  \def\s{.9}
    \draw[0cell=\s] 
    (0,0) node (a) {(\K^{2})\K}
    (225:.707) node (a') {\K(\K^{2})}
    (a)++(1,0) node (b) {\K^{2}}
    (a')++(0,-.6) node (c) {\K^{2}}
    (b)++(0,-1.1) node (d) {\K}
    ;
    \draw[1cell=\s] 
    (a) edge node {\vmtimes 1_\K} (b)
    (c) edge['] node {\vmtimes} (d)
      (a') edge['] node {1_\K \vmtimes} (c)
    (b) edge node {\vmtimes} (d)
    (a') edge node {a^{-\otimes}} (a)
    ;
    \draw[2cell] 
    node[between=b and c at .5, rotate=225, 2label={below,a^\vmtimes_1}] {\Rightarrow}
    ;
  \end{tikzpicture}
  \]
  We denote by $a^{-\vmtimes}_1$ the inverse of $a^\vmtimes_1$.
\end{definition}

\begin{definition}\label{definition:braided-monoidal-vcat}
For a symmetric monoidal category $\V$, a \emph{braided
    monoidal $\V$-category}\index{braided monoidal category!enriched}\index{enriched!braided monoidal category}\index{enriched category!braided monoidal} is a pair $(\K,\beta^\vmtimes)$ consisting of the following data.
  \begin{itemize}
  \item $\K$ is a monoidal $\V$-category (\cref{definition:monoidal-vcat}).
  \item $\beta^\vmtimes$ is a $\V$-natural isomorphism
    \[

    \]
called the \emph{braiding}\index{braiding!braided monoidal enriched category} of $\K$.
  \end{itemize}
  These data are required to satisfy the following two axioms.
\begin{description}
\item[Left Hexagon Axiom] The composites of the following two
  \emph{left hexagon pasting diagrams} are equal.
  \begin{equation}\label{hexagon-bvmL}
    %
  \end{equation}
  In \cref{hexagon-bvmL}, the unlabeled quadrilateral commutes by
  2-naturality of $\beta^\otimes$.  The
  hexagon labeled $\medstar$ commutes by the left hexagon axiom for
  $\beta^\otimes$.
\item[Right Hexagon Axiom] The composites of the following two
  \emph{right hexagon pasting diagrams}\index{hexagon diagram!braided monoidal enriched category} are equal.
  \begin{equation}\label{hexagon-bvmR}
    %
  \end{equation}
  In \cref{hexagon-bvmR}, the unlabeled quadrilateral commutes by
  2-naturality of $\beta^\otimes$.  The
  hexagon labeled $\medstar$ commutes by the right hexagon axiom for
  $\beta^\otimes$.  The 2-cell isomorphism
  $a^{-\vmtimes}_1$ is the inverse of $a^{\vmtimes}_1$ (\cref{eq:a-vmtimes-inv-mate}).
\end{description}
This finishes the definition of a braided monoidal $\V$-category.
\end{definition}

\begin{definition}\label{definition:symm-monoidal-vcat}
For a symmetric monoidal category $\V$, a \emph{symmetric
    monoidal $\V$-category}\index{symmetric monoidal category!enriched}\index{enriched!symmetric monoidal category}\index{enriched category!symmetric monoidal} is a braided monoidal $\V$-category $(\K,\beta^\vmtimes)$ that satisfies the following axiom.
  \begin{description}
  \item[Symmetry Axiom]\index{symmetry axiom!symmetric monoidal enriched category} The composites of the following two
    \emph{symmetry pasting diagrams} are equal.
    
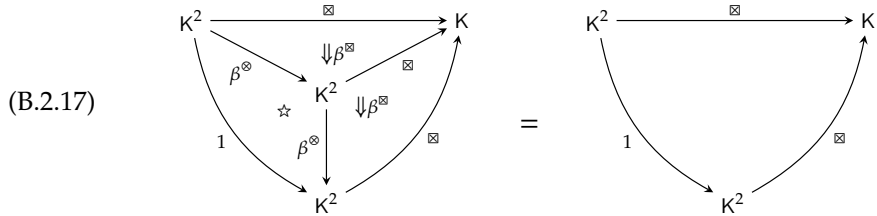
\begin{equation}\label{vmonoidal-symmetry}
    \begin{tikzpicture}[x=18mm,y=14mm,vcenter]
      \def\w{2} 
      \def\h{2} 
      \def\m{.5} 
      \def\s{.8}
      \draw[font=\Large] (\w+\m,0) node (eq) {=}; 
      \newcommand{\boundary}{
        \draw[0cell=\s] 
        (0,0) node (a) {\K^2}
        (a)++(-60:2) node (c) {\K^2}
        (a)++(0:2) node (d) {\K}
        ; 
        \draw[1cell=\s] 
        (a) edge[',bend right=25] node {1} (c)
        (a) edge node {\vmtimes} (d)
        (c) edge[',bend right=25] node {\vmtimes} (d)
        ;
      }
      \begin{scope}[shift={(0,\h/2)}]
        \boundary
        \draw[0cell=\s]
        (c)++(90:2*.57735-.1) node (b) {\K^2}
        ;
        \draw[1cell=\s] 
        (a) edge['] node {\beta^\otimes} (b)
        (b) edge['] node {\beta^\otimes} (c)
        (b) edge['] node {\vmtimes} (d)
        ;
        \draw[2cell=\s] 
        (b)++(90:.57735/2+.07) node[rotate=-90, 2label={above,\beta^\vmtimes}] {\Rightarrow}
        (b)++(-30:.57735/2) node[rotate=-90, 2label={above,\beta^\vmtimes}] {\Rightarrow}
        (b)++(210:.57735/2+.07) node {\medstar}
        ;
      \end{scope}
      \begin{scope}[shift={(\w+\m+\m,\h/2)}]
        \boundary
      \end{scope}
    \end{tikzpicture}
    \end{equation}
  \end{description}
  The right hand diagram in \cref{vmonoidal-symmetry} is the identity $\V$-natural
  transformation.  In the left hand diagram, the region labeled
  $\medstar$ commutes by the symmetry axiom for $\beta^\otimes$.  When this axiom holds, $\beta^\vmtimes$ is also called the \emph{symmetry}\index{symmetry!symmetric monoidal enriched category} of $\K$.
\end{definition}

The following result uses the symmetric monoidal category $(\Cat, \times, \boldone)$ in \cref{ex:cat}.

\begin{theorem}\label{theorem:vcat-cat-mon}
Suppose $\V = (\V,\otimes,\xi)$ is a braided monoidal category.
\begin{enumerate}
\item There is a monoidal $\Cat$-category (\cref{definition:monoidal-vcat})
\[\big( \VCat,\otimes,\vtensorunit, a^\otimes, \ell^\otimes, r^\otimes \big).\]
\item If $\V$ is a symmetric monoidal category, then $(\VCat,\beta^\otimes)$ is a symmetric monoidal $\Cat$-category (\cref{definition:symm-monoidal-vcat}).
\end{enumerate}
\end{theorem}

\subsection*{Enriched Monoidal Functors and Natural Transformations}

\begin{definition}\label{eq:lr-vmtimes-mates}
Suppose $\K$ is a monoidal $\V$-category with $(\V,\xi)$ a braided monoidal category.  We denote by $\ell^\vmtimes_1$\index{left unitor!mate}\index{mate!left unitor} and $r^\vmtimes_1$\index{right unitor!mate}\index{mate!right unitor} the mates of $\ell^\vmtimes$ and $r^\vmtimes$ given,
  respectively, by using $\ell^\otimes$ and $r^\otimes$ in place of
  their inverses, as shown below.
    \[
\]\dqed
\end{definition}

\begin{definition}\label{definition:monoidal-V-fun}
  Suppose $\K$ and $\sfL$ are monoidal $\V$-categories with $\V$ a braided monoidal category.  A
  \emph{monoidal $\V$-functor}\index{monoidal functor!enriched}\index{enriched!monoidal functor}\index{functor!monoidal enriched} 
  \[(F,F^2,F^0)\cn \K \to \sfL\]
is a triple consisting of the following data.
  \begin{itemize}
  \item $F\cn \K \to \sfL$ is a $\V$-functor.
  \item $F^2$ and $F^0$ are $\V$-natural transformations as follows.
  \[
  %
  \]
They are called the \emph{monoidal constraint}\index{monoidal constraint!enriched} and the \index{unit constraint!enriched}\emph{unit constraint}, respectively.
  \end{itemize}
These data are required to satisfy the following associativity and unity axioms.
  \begin{description}
  \item[Associativity]\index{associativity!enriched monoidal functor} The composites of the following two
    \emph{associativity pasting diagrams} are equal.
    \begin{equation}\label{enrmonfunctor-ass}
    %
    \end{equation}
    In the right hand diagram in \cref{enrmonfunctor-ass}, the unlabeled parallelogram 
    commutes by naturality of $a^\otimes$. 
  \item[Left Unity]\index{left unity!enriched monoidal functor} The composites of the following two \emph{left unity pasting
    diagrams} are equal.
    \begin{equation}\label{enrmonfunctor-lunity}
    %
    \end{equation}
    In the right hand diagram in \cref{enrmonfunctor-lunity}, the lower unlabeled quadrilateral
    commutes by naturality of $\ell^\otimes$.  The upper unlabeled
    region commutes by 2-functorality of $\otimes$.  The 2-cell
    isomorphisms labeled $\ell^\vmtimes_1$ are each the mate of
    $\ell^\vmtimes$ (\cref{eq:lr-vmtimes-mates}).
  \item[Right Unity]\index{right unity!enriched monoidal functor} The composites of the following two \emph{right unity pasting
    diagrams} are equal.
    \begin{equation}\label{enrmonfunctor-runity}
    %
    \end{equation}
    In the right hand diagram in \cref{enrmonfunctor-runity}, the lower unlabeled quadrilateral
    commutes by naturality of $r^\otimes$.  The upper unlabeled
    region commutes by 2-functorality of $\otimes$.  The 2-cell
    isomorphisms labeled $r^\vmtimes_1$ are each the mate of
    $r^\vmtimes$ (\cref{eq:lr-vmtimes-mates}).
  \end{description}
  This finishes the definition of a monoidal $\V$-functor.  

Moreover, we define the following variants.
  \begin{itemize}
  \item A \emph{unital monoidal $\V$-functor}\index{unital monoidal!enriched functor} is one for which $F^0$ is invertible.
  \item A \emph{strictly unital monoidal $\V$-functor}\index{strictly unital!monoidal enriched functor} is one for which $F^0$ is an identity.
  \item A \emph{strong monoidal $\V$-functor}\index{strong monoidal!enriched functor} is one for which both $F^0$ and $F^2$ are invertible.
  \item A \emph{strict monoidal $\V$-functor}\index{strict monoidal!enriched functor} is one for which both $F^0$ and $F^2$ are identities.
  \end{itemize}
For each variant, composition of composable monoidal $\V$-functors is defined by composing the $\V$-functors, pasting the monoidal constraints, and pasting the unit constraint.
\end{definition}

\begin{definition}\label{definition:braided-monoidal-vfunctor}
  Suppose $\K$ and $\sfL$ are braided monoidal $\V$-categories with
  $\V$ a symmetric monoidal category.  A \emph{braided monoidal}\index{braided monoidal functor!enriched}\index{enriched!braided monoidal functor}\index{functor!braided monoidal enriched}\index{monoidal functor!braided enriched} $\V$-functor
  \[
  (F,F^2,F^0)\cn \K \to \sfL
  \]
  is a monoidal $\V$-functor that satisfies the following axiom.
  \begin{description}
  \item[Braid Axiom]\index{braid axiom!braided monoidal enriched functor} The composites of the following two \emph{braiding pasting diagrams}\index{braiding diagram!braided monoidal enriched functor} are equal.
    \begin{equation}\label{enrmonfunctor-braided}
    \begin{tikzpicture}[x=17mm,y=10mm,baseline={(eq.base)}]
      \def\wl{1} 
      \def\wr{1} 
      \def\h{1} 
      \def\m{.5} 
      \def\s{.8}
      \draw[font=\Large] (0,0) node (eq) {=}; 
      \newcommand{\boundary}{
        \draw[0cell=\s] 
        (0,0) node (a) {\K^2}
        (-1,-1) node (b) {\K^2}
        (0,-2) node (c) {\K}
        (1,0) node (a') {\sfL^2}
        (1,-2) node (c') {\sfL}
        ; 
        \draw[1cell=\s] 
        (a) edge['] node {\beta^\otimes} (b)
        (b) edge['] node {\vmtimes} (c)
        (a) edge node {FF} (a')
        (c) edge node {F} (c')
        (a') edge node {\vmtimes} (c')
        ;
      }
      \begin{scope}[shift={(-\wl-\m,\h)}]
        \boundary
        \draw[0cell=\s]
        (0,-1) node (b') {\sfL^2}
        ;
        \draw[1cell=\s] 
        (a') edge['] node {\beta^\otimes} (b')
        (b') edge['] node {\vmtimes} (c')
        (b) edge node {FF} (b')
        ;
        \draw[2cell=\s] 
        (b')++(.6,0) node[rotate=180, 2label={below,\beta^\vmtimes}] {\Rightarrow}
        node[between=b' and c at .5, rotate=240, 2label={below,F^2}] {\Rightarrow}
        ;
      \end{scope}
      \begin{scope}[shift={(\wr+\m,\h)}]
        \boundary
        \draw[1cell=\s] 
        (a) edge node {\vmtimes} (c)
        ;
        \draw[2cell=\s] 
        (b)++(.6,0) node[rotate=180, 2label={below,\beta^\vmtimes}] {\Rightarrow}        
        node[between=c and a' at .55, rotate=225, 2label={below,F^2}] {\Rightarrow}
        ;
      \end{scope}
    \end{tikzpicture}
    \end{equation}
  \end{description}
  In the left hand diagram in \cref{enrmonfunctor-braided}, the unlabeled quadrilateral commutes
  by naturality of $\beta^\otimes$.  This finishes the definition of a
  braided monoidal $\V$-functor.

  If $\K$ and $\sfL$ are symmetric monoidal $\V$-categories, then we
  say that $F$ is a \emph{symmetric monoidal $\V$-functor}.\index{symmetric monoidal functor!enriched}\index{enriched!symmetric monoidal - functor}\index{functor!symmetric monoidal enriched}\index{monoidal functor!symmetric enriched}
\end{definition}

\begin{definition}\label{definition:monoidal-V-nt}
Suppose $F,G\cn \K \to \sfL$ are monoidal $\V$-functors between monoidal $\V$-categories with $\V$ a braided monoidal category.  A \emph{monoidal $\V$-natural transformation}\index{monoidal natural transformation!enriched}\index{enriched!monoidal natural transformation}\index{natural transformation!monoidal enriched}
  \[
  \theta\cn F \to G
  \]
  is a $\V$-natural transformation of underlying $\V$-functors (\cref{def:enriched-natural-transformation}) that
  satisfies the following two additional axioms.
  \begin{description}
  \item[Monoidal Naturality]\index{monoidal naturality!enriched} The composites of the following two
    \emph{monoidal naturality pasting diagrams} are equal.
    \[

    \]

  \item[Unit Naturality]\index{unit naturality!enriched} The composites of the following two
    \emph{unit naturality pasting diagrams} are equal.
    \[
    %
    \]
  \end{description}
  This finishes the definition of a monoidal $\V$-natural
  transformation.  Identity and composites of monoidal $\V$-natural transformations
  are defined via underlying $\V$-natural transformations.
\end{definition}

\begin{theorem}\label{thm:vsmcat}
Suppose $\V = (\V,\otimes,\xi)$ is a braided monoidal category.  For items \eqref{it:def-VBMCat} and \eqref{it:def-VSMCat} below, suppose that $\V$ is a symmetric monoidal category.  
  \begin{enumerate}
  \item\label{it:def-VMCat}\index{monoidal category!enriched!2-category}\index{enriched category!monoidal!2-category}\index{2-category!of small monoidal enriched categories}\index{enriched!monoidal category!2-category} There exists a 2-category $\VMCat$ with small monoidal $\V$-categories as objects, monoidal $\V$-functors as 1-cells, and monoidal $\V$-natural transformations as 2-cells.
  \item\label{it:def-VBMCat}\index{braided monoidal category!enriched!2-category}\index{enriched category!braided monoidal!2-category}\index{2-category!of small braided monoidal enriched categories}\index{enriched!braided monoidal category!2-category} There exists a 2-category $\VBMCat$ with small braided monoidal $\V$-categories as objects, braided monoidal $\V$-functors as 1-cells, and monoidal $\V$-natural transformations as 2-cells.
  \item\label{it:def-VSMCat}\index{symmetric monoidal category!enriched!2-category}\index{enriched category!symmetric monoidal!2-category}\index{2-category!of small symmetric monoidal enriched categories}\index{enriched!symmetric monoidal category!2-category} There exists a 2-category $\VSMCat$ with small symmetric monoidal $\V$-categories as objects, symmetric monoidal $\V$-functors as 1-cells, and monoidal $\V$-natural transformations as 2-cells.
  \end{enumerate}
Moreover, there exist forgetful 2-functors
\[\VSMCat \to \VBMCat \to \VMCat \to \VCat.\dqed\]
\end{theorem}

\section{Self-Enriched Symmetric Monoidal Categories}
\label{sec:selfenr-smc}

In this section we review the self-enriched symmetric monoidal structure of a symmetric monoidal closed category (\cref{def:closedcat}).  The material in this section is adapted from \cite[Chapter 3]{cerberusIII}.  Throughout this section, we assume that
\[\big(\V,\otimes,\tu,\alpha,\lambda,\rho,\xi,[,]\big)\]
is a symmetric monoidal closed category.
  
\begin{definition}[Evaluation and Coevaluation]\label{definition:eval}
Suppose $x$ is an object in $\V$.
\begin{itemize}
\item The \index{evaluation!at $x$}\emph{evaluation at $x$} is the counit
\begin{equation}\label{evaluation}
[x,-] \otimes x \fto{\ev_{x,-}} 1_\V.
\end{equation}
\item The \index{coevaluation!at $x$}\index{evaluation!at $x$!co-}\emph{coevaluation at $x$} is the unit
\begin{equation}\label{coevaluation}
1_\V \fto{\coev_{x,-}} [x, - \otimes x].
\end{equation}
\end{itemize} 
These natural transformations refer to the adjunction 
\[- \otimes x \cn \V \lradj \V \cn [x,-]\]
that is part of the closed structure of $\V$.
\end{definition}

Recall the notion of a $\V$-category in \cref{def:enriched-category}.

\begin{definition}[Canonical Self-Enrichment]\label{definition:canonical-v-enrichment}
We define the data of a $\V$-category $\Vse$, called the \index{canonical self-enrichment}\index{self-enrichment}\index{symmetric monoidal category!closed!self-enrichment}\index{category!symmetric monoidal!closed!self-enrichment}\emph{canonical self-enrichment} of $\V$, as follows.
\begin{description}
\item[Objects] $\Ob(\Vse) = \Ob(\V)$.
\item[Hom Objects] Each pair of objects $x,y \in \Vse$ is equipped with the hom object 
\[\Vse(x,y) = [x,y] \in \V.\]
\item[Composition] For objects $x,y,z \in \Vse$, the composition morphism
\[[y,z] \otimes [x,y] \fto{\mcomp_{x,y,z}} [x,z]\]
is the adjoint of the following composite morphism in $\V$.
\begin{equation}\label{eq:m-adj}
\begin{tikzpicture}[x=45mm,y=14mm,vcenter]
\draw[0cell]
(0,0) node (a) {([y,z] \otimes [x,y]) \otimes x}
(0,-1) node (b) {[y,z] \otimes ([x,y] \otimes x)}
(1,-1) node (c) {[y,z] \otimes y}
(1,0) node (d) {z}
;
\draw[1cell]
(a) edge node[swap] {\al} node {\iso} (b)
(b) edge node {1 \otimes \ev} (c)
(c) edge['] node {\ev} (d)
;
\end{tikzpicture}
\end{equation}
\item[Identities] The identity
\[\tu \fto{i_x} [x,x] \forspace x \in \Vse\]
is adjoint to the left unit isomorphism 
\begin{equation}\label{i-adjoint}
\tu \otimes x \fto[\iso]{\la} x \inspace \V.
\end{equation}
\end{description}
This finishes the definition of $\Vse$.  If there is no danger of confusion, we abbreviate $\Vse$ to $\V$.
\end{definition}

Recall from \cref{definition:monoidal-vcat,definition:braided-monoidal-vcat,definition:symm-monoidal-vcat} the notion of a symmetric monoidal $\V$-category.  The following result combines \cite[3.1.11 and 3.3.2]{cerberusIII}.

\begin{theorem}\label{theorem:v-closed-v-sm}\index{canonical self-enrichment!symmetric monoidal}\index{symmetric monoidal category!closed!self-enrichment}
Suppose $\V$ is a symmetric monoidal closed category.  Then the following statements hold.
\begin{enumerate}
\item\label{vse-i}
$\Vse$ in \cref{definition:canonical-v-enrichment} is a $\V$-category.
\item\label{vse-ii}
The symmetric monoidal structure on $\V$ extends to $\Vse$ such that $\Vse$ is a symmetric monoidal $\V$-category.
\end{enumerate}
\end{theorem}

\begin{explanation}[Canonical Self-Enrichment]\label{explanation:recall-adjoints}
In the $\V$-category $\Vse$, the uniqueness of adjoints implies that the composition $\mcomp$ and the identity $i$ are uniquely characterized by the following two diagrams in $\V$.
\begin{equation}\label{eq:adj-comp}
\begin{tikzpicture}[x=35mm,y=13mm,vcenter]
\draw[0cell=.9] 
(0,0) node (l) {\big([y,z] \otimes [x,y]\big) \otimes x}
(0,-1) node (r) {[y,z] \otimes \big([x,y] \otimes x\big)}
(0,-2) node (yzy) {[y,z] \otimes y}
(1,0) node (xzx) {[x,z] \otimes x}
(1,-2) node (z) {z}
;
\draw[1cell=.9]
(l) edge node {\mcomp \otimes 1} (xzx)
(xzx) edge node {\ev} (z)
(l) edge node[swap] {\al} node {\iso} (r)
(r) edge['] node {1 \otimes \ev} (yzy)
(yzy) edge node {\ev} (z)
;
\begin{scope}[shift={(1.5,0)}]
\draw[0cell=.9]
(0,0) node (a) {\tu \otimes x}
(.7,0) node (b) {[x,x] \otimes x}
(b)+(0,-1) node (c) {x}
;
\draw[1cell=.9]
(a) edge node {i \otimes 1} (b)
(b) edge node {\ev} (c)
(a) edge node[swap,pos=.4] {\la} node {\iso} (c)
;
\end{scope}
\end{tikzpicture}
\end{equation}
We use these diagrams in \cref{ex:cl-multi-cl-cat}.
\end{explanation}

\begin{explanation}[Symmetric $\V$-Monoidal Structure]\label{expl:otimes-enriched}
In \cref{theorem:v-closed-v-sm} \eqref{vse-ii}, the symmetric monoidal $\V$-category structure on the $\V$-category $\Vse$ is given as follows.
\begin{description}
\item[Monoidal Composition]
The $\V$-functor
\[\vmtimes\cn \Vse \otimes \Vse \to \Vse\]
has object assignment
\[x \vmtimes y = x \otimes y \forspace x,y \in \V.\]
We use the notation $\vmtimes$ to avoid confusion with the monoidal product $\otimes$ of $\V$ and the tensor product of $\V$-categories.  

For a pair of objects $(x, x'), (y, y') \in \Vse \otimes \Vse$, the morphism
\[\vmtimes_{(x,x'), (y,y')}\cn [x,y] \otimes [x',y'] \to [x \otimes x', y \otimes y'] \inspace \V\]
is adjoint to each of the following two equal composites, with $\ximid$ interchanging the middle two factors.
\begin{equation}\label{vmtimes-adj}
\begin{tikzpicture}[vcenter]
\def\u{-1.4}
\draw[0cell=.9]
(0,0) node (a) {([x,y] \otimes [x',y']) \otimes (x \otimes x')}
(a)+(0,\u) node (b) {([x,y] \otimes x) \otimes ([x',y'] \otimes x')}
(a)+(5.2,0) node (c) {[x \otimes x', y \otimes y'] \otimes (x \otimes x')}
(c)+(0,\u) node (d) {y \otimes y'}
;
\draw[1cell=.9]
(a) edge node[swap,pos=.5] {\ximid} (b)
(b) edge node {\ev \otimes \ev} (d)
(a) edge node {\vmtimes \otimes 1} (c)
(c) edge node {\ev} (d)
;
\end{tikzpicture}
\end{equation}
\item[Monoidal Identity]
The $\V$-functor
\[\vmunit \cn \vtensorunit \to \Vse\]
has object assignment
\[\vmunit(*) = \tensorunit.\]
The morphism between hom objects
\[\vmunit \cn \tensorunit \to [\tensorunit,\tensorunit] \inspace \V\]
is adjoint to each of the following two equal composites.
\begin{equation}\label{vmunit-adj}
\begin{tikzpicture}[vcenter]
\draw[0cell]
(0,0) node (a) {\tu \otimes \tu}
(a)+(3,0) node (b) {[\tu,\tu] \otimes \tu}
(b)+(0,-1.3) node (c) {\tu}
;
\draw[1cell]
(a) edge node[swap,pos=.4] {\la} node {\iso} (c)
(a) edge node {\vmunit \otimes 1} (b)
(b) edge node {\ev} (c)
;
\end{tikzpicture}
\end{equation}
\item[Other Structure]
The monoidal associator, monoidal unitors, and braiding have component morphisms
\begin{equation}\label{Vse-other}
\left\{\begin{aligned}
\al^\adj_{x,y,z} & \cn \tu \to \big[(x \vmtimes y) \vmtimes z, x \vmtimes (y \vmtimes z)\big] \\
\la^\adj_{x} & \cn \tu \to [\tu \vmtimes x, x] \\
\rho^\adj_{x} & \cn \tu \to [x \vmtimes \tu, x] \\
\xi^\adj_{x,y} & \cn \tu \to [x \vmtimes y, y \vmtimes x]
\end{aligned}\right.
\end{equation}
adjoint to the composites of $\la$ with the corresponding components of $\al$, $\la$, $\rho$, and $\xi$, respectively, in $\V$.  In each case, the adjoint component is equal to the composite $\ev \comp (?^\adj \otimes 1)$, similar to the top-right composites in \cref{vmtimes-adj,vmunit-adj}.
\end{description}
This finishes the description of the symmetric monoidal $\V$-category $\Vse$.
\end{explanation}

\section{Change of Enrichment}
\label{sec:change-enrichment}

Recall from \cref{ex:vcatastwocategory} that each monoidal category $\V$ has an associated 2-category $\VCat$ of small $\V$-categories, $\V$-functors, and $\V$-natural transformations.  In this section we review properties of changing the enriching monoidal category $\V$.  Whenever we need the enriching monoidal category to be braided or symmetric (\cref{def:braidedmoncat,def:symmoncat}), we state so explicitly.  The material in this section is adapted from \cite[Chapters 2 and 3]{cerberusIII}.  

\begin{definition}\label{definition:U-VCat-WCat}
Suppose given a monoidal functor between monoidal categories
\[(U,U^2,U^0) \cn (\V,\otimes,\tu) \to (\W,\otimes,\tu).\]
We define the data of a 2-functor\label{not:dU}
\[(-)_U \cn \VCat \to \WCat,\] 
called the \index{change of enrichment}\index{enriched category!change of enrichment}\index{monoidal functor!change of enrichment}\index{symmetric monoidal functor!change of enrichment}\index{2-functor!change of enrichment}\emph{change of enrichment}, as follows.
\begin{description}
\item[Object Assignment] For a $\V$-category $(\C,\mcomp,i)$, the $\W$-category 
\[\big(\C_U, \mcomp_U, i_U\big)\]
has objects
\[\Ob(\C_U) = \Ob(\C)\]
and hom objects
\[\C_U(x,y) = U\C(x,y) \in \W \forspace x,y \in \C_U.\]
The composition in $\C_U$ is the following composite in $\W$ for $x,y,z \in \C_U$.
\begin{equation}\label{eq:comp-U}

\end{equation}
The identity of an object $x \in \C_U$ is the following composite.
\begin{equation}\label{eq:id-U}
%
\end{equation}
\item[1-Cell Assignment] For a $\V$-functor $F\cn \C \to \D$, the $\W$-functor 
\[F_U \cn \C_U \to \D_U\]
has the same object assignment as $F$.  On hom objects it is the morphism
\begin{equation}\label{FUxy}
(F_U)_{x,y} = U(F_{x,y}) \cn U\C(x,y) \to U\D(Fx,Fy) \inspace \W
\end{equation}
for $x,y \in \C_U$.
\item[2-Cell Assignment] For a $\V$-natural transformation $\theta$ as in the left diagram below
\[%
\]
the $\W$-natural transformation $\tha_U$, as in the right diagram above, has component morphism at $x \in \C_U$ given by the following composite.
\begin{equation}\label{thetaUx}
%
\end{equation}
\end{description}
This finishes the definition of $(-)_U$.
\end{definition}

The following is \cite[2.1.2]{cerberusIII}.

\begin{proposition}\label{proposition:U-VCat-WCat}
In the context of \cref{definition:U-VCat-WCat}, 
\[(-)_U \cn \VCat \to \WCat\] 
is a 2-functor.
\end{proposition}

Change of enrichment is compatible with composition of monoidal functors (\cref{def:mfunctor-comp}), as in the following result from \cite[2.2.4]{cerberusIII}.

\begin{proposition}\label{proposition:change-enr-horiz-comp}
Given monoidal functors between monoidal categories
\[\V_1 \fto{U_1} \V_2 \fto{U_2} \V_3,\]
the following diagram of change-of-enrichment 2-functors commutes.
\[\begin{tikzpicture}[baseline={(a.base)}]
\def\u{.7}
\draw[0cell]
(0,0) node (a) {\V_1\mh\Cat}
(a)+(3,0) node (b) {\V_2\mh\Cat}
(b)+(3,0) node (c) {\V_3\mh\Cat}
;
\draw[1cell=.9]
(a) edge node {(-)_{U_1}} (b)
(b) edge node {(-)_{U_2}} (c)
;
\draw[1cell=.9]
(a) [rounded corners=3pt] |- ($(b)+(-1,\u)$)
-- node {(-)_{U_2U_1}} ($(b)+(1,\u)$) -| (c)
;
\end{tikzpicture}\]
\end{proposition}

\subsection*{Compatibility with Enriched Tensor Product}

For a braided monoidal category $\V$, $(\VCat,\otimes)$ is a monoidal $\Cat$-category, which is, furthermore, symmetric if $\V$ is symmetric (\cref{theorem:vcat-cat-mon}).  Change of enrichment is compatible with the tensor product of enriched categories (\cref{definition:vtensor-0}), using the following definitions.

\begin{definition}\label{definition:U2-VCat-WCat}
Suppose given a braided monoidal functor between braided monoidal categories
\[(U, U^2, U^0) \cn \V \to \W.\]
We define monoidal constraint $(-)_U^2$ and unit constraint $(-)_U^0$ for the change-of-enrichment 2-functor $(-)_U$ as follows.
\begin{description}
\item[Monoidal Constraint] Its component $\W$-functor at small $\V$-categories $\C$ and $\D$\index{change of enrichment!monoidal constraint}\index{monoidal constraint!change of enrichment},
\[(-)^2_U\cn \C_U \otimes \D_U \to (\C \otimes \D)_U,\]
has the identity object assignment.  On hom objects it is given by the following morphism in $\W$ for $x,x'\in\C$ and $y,y' \in \D$.
\[\begin{tikzpicture}[x=54mm,y=9mm]
\draw[0cell]
(0,0) node (a) {(\C_U \otimes \D_U)(x \otimes y, x' \otimes y')}
(0,-1) node (a') {U\C(x,x')\otimes U\D(y,y')}
(1,0) node (b) {(\C \otimes \D)_U(x \otimes y, x' \otimes y')}
(1,-1) node (b') {U\big(\C(x,x') \otimes \D(y,y')\big)}
;
\draw[1cell=1.5] 
node[between=a and a' at .5, rotate=90] {=}
node[between=b and b' at .5, rotate=90] {=}
;
\draw[1cell] 
(a') edge node {U^2} (b')
;
\end{tikzpicture}\]
\item[Unit Constraint] It is the $\W$-functor\index{change of enrichment!unit constraint}\index{unit constraint!change of enrichment}
\[(-)_U^0\cn \vtensorunit \to \vtensorunit_U\]
given by the identity on the unique object and the morphism
\[\tensorunit \fto{U^0} U\tensorunit \inspace \W\]
on the unique hom object.
\end{description}
This finishes the definition of $(-)_U^2$ and $(-)_U^0$ .
\end{definition}

The following is \cite[2.3.7]{cerberusIII}, which uses \cref{theorem:vcat-cat-mon,definition:monoidal-V-fun,definition:braided-monoidal-vfunctor}.

\begin{theorem}\label{theorem:U-braided-mon}
For each braided monoidal functor between braided monoidal categories
\[(U, U^2, U^0) \cn \V \to \W,\]
the triple in \cref{definition:U-VCat-WCat,definition:U2-VCat-WCat}
\[\big((-)_U, (-)_U^2, (-)_U^0\big) \cn (\VCat, \otimes) \to (\WCat, \otimes)\]
is a monoidal $\Cat$-functor.  Moreover, if $U$ is a symmetric monoidal functor between symmetric monoidal categories, then $(-)_U$ is a symmetric monoidal $\Cat$-functor.
\end{theorem}

\subsection*{Compatibility with Enriched Monoidal Structure}

Change of enrichment preserves enriched monoidal structure, as in the following result from \cite[2.4.10]{cerberusIII}.

\begin{theorem}\label{theorem:KU-monoidal}
Suppose given a braided monoidal functor between braided monoidal categories
\[U \cn \V \to \W.\]
For \eqref{it:KU-monoidal} and \eqref{it:FU-monoidal} below, the braided and symmetric monoidal cases assume that $U$, $\V$, and $\W$ are symmetric monoidal.
\begin{enumerate}
\item\label{it:KU-monoidal} If $\K$ is a (braided, respectively symmetric) monoidal $\V$-category, then $\K_U$ is a (braided, respectively symmetric) monoidal $\W$-category.
\item\label{it:FU-monoidal} If $F\cn\K \to \sfL$ is a (braided, respectively symmetric) monoidal $\V$-functor between (braided, respectively symmetric) monoidal $\V$-categories, then
\[F_U\cn \K_U \to \sfL_U\]
is a (braided, respectively symmetric) monoidal $\W$-functor.
\item\label{it:thU-monoidal} If $\theta\cn F \to G$ is a monoidal $\V$-natural transformation between monoidal $\V$-functors $F$ and $G$, then
\[\theta_U\cn F_U \to G_U\]
is a monoidal $\W$-natural transformation.
\end{enumerate}
\end{theorem}

\cref{theorem:v-closed-v-sm} \eqref{vse-ii} and \cref{theorem:KU-monoidal} \eqref{it:KU-monoidal} yield the following.

\begin{corollary}\label{VseU}
Suppose given a symmetric monoidal functor between symmetric monoidal categories
\[(U, U^2, U^0) \cn \V \to \W\]
with $(\V,[,])$ closed.  Then $\Vse_U$ is a symmetric monoidal $\W$-category.
\end{corollary}

\begin{explanation}\label{expl:VseU}
In \cref{VseU} the symmetric monoidal $\W$-category $\Vse_U$ is given explicitly as follows.
\begin{description}
\item[Underlying $\W$-Category] The $\W$-category 
\[\big(\Vse_U, \mcomp_U, i_U\big)\] 
is obtained from the $\V$-category $(\Vse,\mcomp,i)$ in \cref{definition:canonical-v-enrichment} by applying the change-of-enrichment $(-)_U$ in \cref{definition:U-VCat-WCat}.  In other words, it has objects
\[\Ob(\Vse_U) = \Ob(\Vse) = \Ob(\V)\]
and hom objects
\[\Vse_U(x,y) = U\Vse(x,y) = U[x,y] \in \W \forspace x,y \in \V.\]
The composition is the following composite, with $\mcomp$ the composition in $\Vse$ in \cref{eq:adj-comp}.
\[
\]
The identity of an object $x \in \Vse_U$ is the following composite, with $i$ the identity in $\Vse$ in \cref{eq:adj-comp}.
\[%
\]
\item[Monoidal Composition] The $\W$-functor
\[\Vse_U \otimes \Vse_U \fto{\vmtimes_U} \Vse_U\]
is given on objects by
\[x \vmtimes_U y = x \otimes y \forspace x,y \in \V.\]
On hom objects it is the composite
\[%
\]
for $x \otimes x', y \otimes y' \in \Vse_U \otimes \Vse_U$, with $\vmtimes$ as in \cref{vmtimes-adj}.
\item[Monoidal Unit] The $\W$-functor
\[\vmunit_U \cn \vtensorunit \to \Vse_U\]
has object assignment
\[\vmunit_U(*) = \tu.\]
On hom objects it is the following composite, with $\vmunit$ as in \cref{vmunit-adj}.
\[%
\]
\item[Other Structure] The monoidal associator, monoidal unitors, and braiding have component morphisms as follows, with $\al^\adj$, $\la^\adj$, $\rho^\adj$, and $\xi^\adj$ as in \cref{Vse-other}.
\[\left\{
\begin{aligned}
\tu \fto{U^0} U\tu & \fto{U\al^\adj_{x,y,z}} U\big[(x \vmtimes y) \vmtimes z, x \vmtimes (y \vmtimes z)\big] \\
\tu \fto{U^0} U\tu & \fto{U\la^\adj_{x}} U[\tu \vmtimes x, x] \\
\tu \fto{U^0} U\tu & \fto{U\rho^\adj_{x}} U[x \vmtimes \tu, x] \\
\tu \fto{U^0} U\tu & \fto{U\xi^\adj_{x,y}} U[x \vmtimes y, y \vmtimes x]
\end{aligned}\right.\]
\end{description}
This finishes the description of the symmetric monoidal $\W$-category $\Vse_U$.
\end{explanation}

\subsection*{Standard Enrichment}

The following definition uses the symmetric monoidal $\W$-category $\Vse_U$ in \cref{expl:VseU}.

\begin{definition}\label{definition:U-std-enr}
Suppose given a monoidal functor 
\[(U,U^2,U^0) \cn (\V,\otimes,\tu,[,]) \to (\W,\otimes,\tu,[,])\] 
between symmetric monoidal closed categories.  We define the data of a monoidal $\W$-functor\label{not:Use}
\[\big(\Use, \Use^2, \Use^0\big) \cn \Vse_U \to \Wse,\]
called the \index{enrichment!standard}\index{standard enrichment!symmetric monoidal functor}\index{symmetric monoidal functor!standard enrichment}\emph{standard enrichment} of $U$, as follows.
\begin{description}
\item[Object Assignment] The object assignment of $\Use$ is the same as that of $U$.
\item[Component Morphisms] For objects $x,y \in \V$, the component morphism
\[\Use_{x,y} \cn U[x,y] \to [Ux,Uy] \inspace \W\]
is adjoint to the composite
\begin{equation}\label{UtwoUev}
U[x,y] \otimes Ux \fto{U^2} U\big([x,y] \otimes x\big) \fto{U(\ev)} Uy.
\end{equation}
\item[Unit Constraint] The morphism
\[\Use^0 \cn \tu \to [\tu,U\tu] \inspace \W\]
is adjoint to the composite
\begin{equation}\label{laUzero}
\tu \otimes \tu \fto[\iso]{\la} \tu \fto{U^0} U\tu.
\end{equation}
\item[Monoidal Constraint] Its component morphism
\[\Use^2_{x \otimes x'} \cn \tu \to \big[Ux \otimes Ux', U(x\otimes x')\big] \forspace x \otimes x' \in \V \otimes \V\]
is adjoint to the composite 
\begin{equation}\label{laUtwo}
\tu \otimes (Ux \otimes Ux') \fto[\iso]{\la} Ux \otimes Ux' \fto{U^2_{x,x'}} U(x \otimes x').
\end{equation}
\end{description}
This finishes the definition of the standard enrichment of $U$.
\end{definition}

The following is \cite[3.3.4]{cerberusIII}.

\begin{proposition}\label{proposition:U-std-enr}
In the context of \cref{definition:U-std-enr}, the standard enrichment
\[\big(\Use, \Use^2, \Use^0\big) \cn \Vse_U \to \Wse\]
is a monoidal $\W$-functor.  Moreover, if $U$ is a symmetric monoidal functor, then $\Use$ is a symmetric monoidal $\W$-functor.
\end{proposition}

\begin{explanation}\label{expl:std-enr}
Each of $\Use_{x,y}$, $\Use^0$, and $\Use^2_{x \otimes x'}$ is defined by its adjoint in, respectively,  \cref{UtwoUev,laUzero,laUtwo}.  Thus the uniqueness of adjoints implies that these structure morphisms are uniquely characterized by the following commutative diagrams.
\[
\]
We use this adjoint characterization of $\Use_{x,y}$ in \cref{std-enr-monoidal}.
\end{explanation}

\chapter{Multicategories}
\label{ch:prelim_multicat}
In this appendix we review basic elements of multicategory theory.  The following table summarizes the main content in this appendix.
\smallskip
\begin{center}
\resizebox{.95\textwidth}{!}{
{\renewcommand{\arraystretch}{1.4}%
{\setlength{\tabcolsep}{1ex}
\begin{tabular}{|c|c|}\hline
\multicolumn{2}{|c|}{\appnumname{sec:enrmulticat} to \appnumname{sec:mult-from-smc}} \\ \hline
$\V$-multicategories (multifunctors, multinatural transformations) & \ref{def:enr-multicategory} (\ref{def:enr-multicategory-functor}, \ref{def:enr-multicat-natural-transformation}) \\ \hline
underlying $\V$-category & \ref{ex:unarycategory} \\ \hline
2-category $\VMulticat$ of $\V$-multicategories & \ref{v-multicat-2cat} \\ \hline
$\Cat$-multinatural transformations & \ref{expl:catmultitransformation} \\ \hline
endomorphism multicategories $\End$ (enriched) & \ref{ex:endc} (\ref{definition:EndK}) \\ \hline
\hline
\multicolumn{2}{|c|}{\appnumname{sec:ptmulticatiicat}} \\ \hline
pointed multicategories, multifunctors, and multinatural transformations & \ref{def:ptd-multicat} \\ \hline
pointed endomorphism multicategories $\Endst$ & \ref{ex:endstc} \\ \hline
2-category $\pMulticat$ of pointed multicategories & \ref{thm:pmulticat} \\ \hline
free-forgetful 2-adjunction $\dplus \cn \Multicat \lradj \pMulticat \cn \Ust$ & \ref{dplustwofunctor} \\ \hline
\end{tabular}}}}
\end{center}
\medskip
The main reference for this chapter is \cite{cerberusIII}; see also \cite{yau-operad}.  \cref{conv:universe,expl:leftbracketing} are still in effect.

\section{Enriched Multicategories}
\label{sec:enrmulticat}

In this section we review enriched multicategories, multifunctors, and multinatural transformations.  The material in this section is adapted from \cite[Section 6.1]{cerberusIII}.  Throughout this section we assume that 
\[\big(\V,\otimes,\tu,\alpha,\lambda,\rho,\xi\big)\]
is a symmetric monoidal category (\cref{def:symmoncat}).  We use the following notation for finite tuples of objects.

\begin{definition}\label{def:profile}
Suppose $S$\label{notation:s-class} is a class.  
\begin{description}
\item[Profiles] The class of finite tuples in $S$ is denoted by\index{profile}\label{notation:profs}
\[\Prof(S) = \coprodover{k \geq 0}\ S^{k}.\] 
\begin{itemize}
\item An element in $\Prof(S)$ is called an \emph{$S$-profile}.  
\item An $S$-profile of length\index{length!of a profile} $n=\len\angx$ is denoted by\label{notation:us}  
\[\angx = (x_1, \ldots, x_n) = \ang{x_j}_{j=1}^n \in S^{n}.\]
The empty $S$-profile\index{empty profile} is denoted by $\ang{}$.
\item An element in $\Prof(S)\times S$ is denoted by\label{notation:duc} $\smscmap{\angx; y}$ with $\angx\in\Prof(S)$ and $y\in S$.
\end{itemize}
\item[Concatenation] For two $S$-profiles $\angx = \ang{x_i}_{i=1}^m$ and $\angy = \ang{y_j}_{j=1}^n$, their \emph{concatenation}\index{concatenation} is the $S$-profile\label{not:concat}
\begin{equation}\label{concatprof}
\angx \oplus \angy = \big(x_1, \ldots , x_m, y_1, \ldots , y_n\big).
\end{equation}
Concatenation is associative with the empty tuple $\ang{}$ as the strict unit.\defmark
\end{description}
\end{definition}

The \index{symmetric group}symmetric group on $n$ letters is denoted $\Sigma_n$.

\begin{definition}\label{def:enr-multicategory}
A \emph{$\V$-multicategory}\index{enriched multicategory}\index{category!enriched multi-}\index{multicategory!enriched} is a triple\label{notation:enr-multicategory}  
\[(\M, \gamma, \operadunit)\]
consisting of the following data.
\begin{description}
\item[Objects] It is equipped with a class $\ObM$, whose elements are called\index{object!enriched multicategory} \emph{objects}.  We abbreviate $\Prof(\ObM)$ to $\ProfM$.
\item[Multimorphisms] For $\smscmap{\angx;x'} \in \ProfMM$ with $\angx = \ang{x_j}_{j=1}^n$, it is equipped with an object in $\V$\label{notation:enr-cduc}
\[\M\scmap{\angx;x'} = \M\mmap{x'; x_1,\ldots,x_n}.\]
It is called the \emph{$n$-ary operation object}\index{n-ary@$n$-ary!operation object} or \index{multimorphism}\emph{$n$-ary multimorphism object} with \emph{input profile}\index{input profile} $\angx$ and \emph{output}\index{output} $x'$.
\begin{itemize}
\item We also say \emph{nullary} for 0-ary, \emph{unary} for 1-ary, and \emph{binary} for 2-ary.
\item If objects in $\V$ have underlying sets (for example, if $\V$ is $\Set$ or $\Cat$), then an element in $\M\scmap{\angx;x'}$ is called an \emph{$n$-ary multimorphism} or \emph{$n$-ary operation} and denoted
\[\angx = (x_1,\ldots,x_n) \to x'.\]
\end{itemize}   
\item[Symmetric Group Action]
For $\smscmap{\angx;x'} \in \ProfMM$ and a permutation $\sigma \in
\Sigma_n$, $\M$ is equipped with an isomorphism in $\V$
\begin{equation}\label{rightsigmaaction}
\begin{tikzcd}[column sep=large]
\M\scmap{\angx;x'} \rar{\sigma}[swap]{\cong} & \M\scmap{\angx\sigma; x'},\end{tikzcd}
\end{equation}
called the \emph{right $\sigma$-action}\index{right action} or the \index{symmetric group!action}\emph{symmetric group action}, where\label{enr-notation:c-sigma}
\[\angx\sigma  
= \left(x_{\sigma(1)}, \ldots, x_{\sigma(n)}\right) 
= \ang{x_{\sigma(j)}}_{j=1}^n \in \ProfM\]
is the right permutation\index{right permutation} of $\angx$ by $\sigma$.
\item[Units] Each object $x$ in $\M$ is equipped with a morphism\label{notation:enr-unit-c}
\begin{equation}\label{ccoloredunit}
\begin{tikzcd}[column sep=large]
\tu \ar{r}{\operadunit_x} & \M\scmap{x;x},
\end{tikzcd}
\end{equation}
called the \index{colored unit}\emph{$x$-colored unit}.
\item[Composition] Suppose given
\begin{itemize}
\item $\scmap{\ang{x'};x''} \in \ProfMM$ with $\ang{x'} = \ang{x'_j}_{j=1}^n \in \ProfM$ and 
\item $\ang{x_j} = \ang{x_{j,i}}_{i=1}^{k_j} \in \ProfM$ for each $j\in\{1,\ldots,n\}$ with $\angx = \oplus_{j=1}^n \ang{x_j}$.
\end{itemize}
Then $\M$ is equipped with a morphism in $\V$\label{notation:enr-multicategory-composition}
\begin{equation}\label{eq:enr-defn-gamma}
\begin{tikzcd}[column sep=large]
\M\mmap{x'';\ang{x'}} \otimes
\txotimes_{j=1}^n \M\mmap{x_j';\ang{x_j}} \rar{\gamma} &
\M\mmap{x'';\ang{x}},
\end{tikzcd}
\end{equation}
called the \index{multicategory!composition}\index{composition!multicategory}\emph{composition} or \emph{multicategorical composition}.  If objects in $\V$ have underlying sets, then we also denote composition diagrammatically by
\[\big(\ang{x_1}, \ldots, \ang{x_n}\big) \fto{(f_1,\ldots,f_n)} \angxp \fto{f} x''\]
for multimorphisms
\[f \in \M\mmap{x'';\ang{x'}} \andspace f_j \in \M\mmap{x_j';\ang{x_j}}.\]
\end{description}
The data above are required to satisfy the following axioms.
\begin{description}
\item[Symmetric Group Action]
The identity in $\Sigma_n$ acts as the identity morphism on $\M\scmap{\angx;x'}$ with $n=\len\ang{x}$.  Moreover, for $\sigma,\tau\in\Sigma_n$, the following diagram in $\V$ commutes. 
\begin{equation}\label{enr-multicategory-symmetry}
\begin{tikzcd}
\M\scmap{\angx;x'} \arrow{rd}[swap]{\sigma\tau} \rar{\sigma} 
& \M\scmap{\angx\sigma;x'} \dar{\tau}\\
& \M\mmap{x';\angx\sigma\tau}
\end{tikzcd}
\end{equation}
\item[Associativity]
Suppose given
\begin{itemize}
\item $\scmap{\ang{x''};x'''} \in \ProfMM$ with $\ang{x''} = \ang{x''_j}_{j=1}^n \in \ProfM$,
\item $\ang{x_j'} = \ang{x'_{j,i}}_{i=1}^{k_{j}} \in \ProfM$ for each $j \in \{1,\ldots,n\}$ with $\ang{x'} = \oplus_{j=1}^n \ang{x_j'}$ and $k_j > 0$ for at least one $j$, and
\item $\ang{x_{j,i}} = \ang{x_{j,i,p}}_{p=1}^{\ell_{j,i}} \in \ProfM$ for each $j\in\{1,\ldots,n\}$ and each $i \in \{1,\ldots,k_j\}$ with $\ang{x_j} = \oplus_{i=1}^{k_j}{\ang{x_{j,i}}}$ and $\ang{x} =
\oplus_{j=1}^n{\ang{x_{j}}}$.
\end{itemize}
Then the \index{associativity!enriched multicategory}\emph{associativity diagram} below commutes.
\begin{equation}\label{enr-multicategory-associativity}
\begin{tikzpicture}[x=40mm,y=15mm,vcenter]
  \draw[0cell=.85] 
  (0,0) node (a) {\textstyle
    \M\mmap{x''';\ang{x''}}
    \otimes
    \biggl[\bigotimes\limits_{j=1}^n \M\mmap{x''_j;\ang{x'_{j}}}\biggr]
    \otimes
    \bigotimes\limits_{j=1}^n \biggl[\bigotimes\limits_{i=1}^{k_j} \M\mmap{x'_{j,i};\ang{x_{j,i}}}\biggr] 
  }
  (1,.8) node (b) {\textstyle
    \M\mmap{x''';\ang{x'}}
    \otimes
    \bigotimes\limits_{j=1}^{n} \biggl[\bigotimes\limits_{i=1}^{k_j} \M\mmap{x'_{j,i};\ang{x_{j,i}}}\biggr]
  }
  (0,-1.2) node (a') {\textstyle
    \M\mmap{x''';\ang{x''}} \otimes
    \bigotimes\limits_{j=1}^n \biggl[\M\mmap{x_j'';\ang{x_j'}} \otimes \bigotimes\limits_{i=1}^{k_j} \M\mmap{x'_{j,i};\ang{x_{j,i}}}\biggr]
  }
  (1,-2) node (b') {\textstyle
    \M\mmap{x''';\ang{x''}} \otimes \bigotimes\limits_{j=1}^n \M\mmap{x_j'';\ang{x_{j}}}
  }
  (1.2,-.6) node (c) {\textstyle
    \M\mmap{x''';\ang{x}}
  }
  ;
  \draw[1cell=.85]
  (a) edge[shorten <=-1ex,shorten >=-1ex] node {\iso} node['] {\mathrm{permute}} (a')
  (a) edge[shorten >=-4ex,shorten <=-1ex, transform canvas={xshift=-2.5em}] node[pos=.6] {(\ga,1)} (b)
  (b) edge node {\ga} (c)
  (a') edge[shorten >=-3ex, shorten <=-1ex,transform canvas={xshift=-2.5em}] node[swap,pos=.6] {(1,\textstyle\bigotimes_j \ga)} (b')
  (b') edge['] node {\ga} (c)
  ;
\end{tikzpicture}
\end{equation}

\item[Unity]
Suppose $\scmap{\angx;x'} \in \ProfMM$ with $\angx = \ang{x_j}_{j=1}^n \in \ProfM$.\index{unity!enriched multicategory}
\begin{enumerate}
\item If $n \geq 1$, then the following \emph{right unity diagram}\index{right unity!enriched multicategory} commutes.
\begin{equation}\label{enr-multicategory-right-unity}
\begin{tikzcd} 
\M\scmap{\angx;x'} \otimes \tensorunit^{\otimes n} \dar[swap]{1 \otimes (\otimes_j \operadunit_{x_j})} \ar[bend left=15]{dr}{\rho} &\\
\M\scmap{\angx;x'} \otimes \txotimes_{j=1}^n \M\scmap{x_j;x_j} \rar{\gamma} & \M\scmap{\angx;x'}
\end{tikzcd}
\end{equation}

\item
The \index{unity!enriched multicategory}\index{left unity!enriched multicategory}\emph{left unity diagram} below commutes.
\begin{equation}\label{enr-multicategory-left-unity}
\begin{tikzcd}
\tensorunit \otimes \M\scmap{\angx;x'} \dar[swap]{\operadunit_{x'} \otimes 1} \ar[bend left=15]{dr}{\lambda} &\\
\M\mmap{x';x'} \otimes \M\scmap{\angx;x'} \rar{\gamma} & \M\scmap{\angx;x'}
\end{tikzcd}
\end{equation}
\end{enumerate}

\item[Equivariance]
Suppose $\len\ang{x_j} = k_j \geq 0$ in the definition of $\gamma$ \eqref{eq:enr-defn-gamma}.\index{equivariance!enriched multicategory}
\begin{enumerate}
\item For each $\sigma \in \Sigma_n$, the following \index{top equivariance!enriched multicategory}\emph{top equivariance diagram} commutes.
\begin{equation}\label{enr-operadic-eq-1}
\begin{tikzcd}[column sep=large,cells={nodes={scale=.8}},
every label/.append style={scale=.8}]
\M\mmap{x'';\ang{x'}} \otimes \txotimes_{j=1}^n \M\mmap{x'_j;\ang{x_j}} 
\dar[swap]{\gamma} \rar{(\sigma, \sigma^{-1})}
& \M\mmap{x'';\ang{x'}\sigma} \otimes \txotimes_{j=1}^n \M\mmap{x'_{\sigma(j)};\ang{x_{\sigma(j)}}} \dar{\gamma}\\
\M\mmap{x'';\ang{x_1},\ldots,\ang{x_n}} \rar{\sigma\langle k_{\sigma(1)}, \ldots , k_{\sigma(n)}\rangle}
& \M\mmap{x'';\ang{x_{\sigma(1)}},\ldots,\ang{x_{\sigma(n)}}}
\end{tikzcd}
\end{equation}
In \cref{enr-operadic-eq-1} the block permutation\index{block!permutation}\label{notation:enr-block-permutation}  
\begin{equation}\label{blockpermutation}
\sigma\langle k_{\sigma(1)}, \ldots , k_{\sigma(n)} \rangle \in \Sigma_{k_1+\cdots+k_n}
\end{equation}
permutes $n$ consecutive intervals of lengths $k_{\sigma(1)}$, $\ldots$, $k_{\sigma(n)}$, respectively, as $\sigma$ permutes $\{1,\ldots,n\}$, without changing the order within each interval.
\item
Given permutations $\tau_j \in \Sigma_{k_j}$ for $1 \leq j \leq n$,
the following \index{bottom equivariance!enriched multicategory}\emph{bottom equivariance diagram} commutes.
\begin{equation}\label{enr-operadic-eq-2}
\begin{tikzcd}[cells={nodes={scale=.85}},
every label/.append style={scale=.85}]
\M\mmap{x'';\ang{x'}} \otimes \bigotimes_{j=1}^n \M\mmap{x'_j;\ang{x_j}}
\dar[swap]{\gamma} \rar{(1, \otimes_j \tau_j)} & 
\M\mmap{x'';\ang{x'}} \otimes \bigotimes_{j=1}^n \M\mmap{x'_j;\ang{x_j}\tau_j}\dar{\gamma} \\
\M\mmap{x'';\ang{x_1},\ldots,\ang{x_n}} \rar{\tau_1 \times \cdots \times \tau_n}
& \M\mmap{x'';\ang{x_1}\tau_1,\ldots,\ang{x_n}\tau_n}
\end{tikzcd}
\end{equation}
In \cref{enr-operadic-eq-2} the block sum\index{block!sum}\label{notation:enr-block-sum} 
\begin{equation}\label{blocksum}
\tau_1 \times\cdots \times\tau_n \in \Sigma_{k_1+\cdots+k_n}
\end{equation} 
is the image of $(\tau_1, \ldots, \tau_n)$ under the canonical inclusion 
\[\Sigma_{k_1} \times \cdots \times \Sigma_{k_n} \to \Sigma_{k_1 + \cdots + k_n}.\]
\end{enumerate}
\end{description}
This finishes the definition of a $\V$-multicategory.  A $\V$-multicategory is \emph{small}\index{multicategory!enriched!small}\index{small!enriched multicategory} if it has a set of objects.

Moreover, we define the following variants.
\begin{itemize}
\item A \emph{non-symmetric $\V$-multicategory}\index{non-symmetric!multicategory}\index{multicategory!non-symmetric} is defined in the same way as a $\V$-multicategory by omitting the symmetric group action and the axioms \cref{enr-multicategory-symmetry,enr-operadic-eq-1,enr-operadic-eq-2} involving the symmetric group action.
\item A \index{multicategory}\emph{(non-symmetric) multicategory} is a (non-symmetric) $\Set$-multicategory, where $(\Set,\times,*)$ is the symmetric monoidal category of sets and functions with the Cartesian product as the monoidal product.
\item A \index{enriched operad}\index{enriched!operad}\index{operad!enriched}\emph{$\V$-operad} is a $\V$-multicategory with one object.  If $\M$ is a $\V$-operad, then its $n$-ary multimorphism object is denoted by \label{not:nthobject}$\M_n \in \V$.
\item An \emph{operad}\index{multicategory!one object} is a $\Set$-operad, that is, a multicategory with one object.\defmark
\end{itemize}
\end{definition}

\begin{remark}[Related Concepts]\label{rk:nsmulticat}\
\begin{enumerate}
\item In the literature, including \cite{lambek} where this concept originated, a \emph{multicategory} sometimes means a non-symmetric multicategory in the sense of \cref{def:enr-multicategory}.  Our convention is to include the symmetric group action by default.  We always include the word \emph{non-symmetric} if we are referring to the variant without the symmetric group action.
\item There are more conceptual ways to define enriched multicategories as (i) monoids in a certain monoidal category and (ii) algebras over some monad.  See \cite[Ch.\! 4]{yau-hqft}.
\item There are variants of enriched multicategories whose equivariant structure is parametrized by groups different from the symmetric groups, such as the braid groups.  See \cite{yau-inf-operad}.
\item There are many different but related generalizations of enriched multicategories whose operation objects have input and output profiles of arbitrary lengths.  See \cite[Section 2.5]{johnson-yau} for one such variant called \emph{polycategories} and \cite{bluemonster} for other variants.\defmark
\end{enumerate}
\end{remark}

\begin{example}[Underlying $\V$-Categories]\label{ex:unarycategory}\index{enriched multicategory!underlying enriched category}\index{enriched category!underlying - of an enriched multicategory}
Each non-symmetric $\V$-multicategory $(\M,\ga,\operadunit)$ has an underlying $\V$-category (\cref{def:enriched-category}) defined as follows.
\begin{itemize}
\item It has the same class of objects as $\M$.
\item For objects $x,y \in \M$, the hom object is $\M\smscmap{x;y}$.
\item The identities are the colored units in $\M$.
\item The composition is given by 
\[\begin{tikzcd}[column sep=large]
\M\scmap{y;z} \otimes \M\scmap{x;y} \ar{r}{\gamma} & \M\scmap{x;z}
\end{tikzcd}\]
for objects $x, y, z \in \M$.
\end{itemize}   
The associativity and unity diagrams, \cref{enriched-cat-associativity,enriched-cat-unity}, of a $\V$-category are the 1-ary restrictions of, respectively, the associativity and unity diagrams, \cref{enr-multicategory-associativity,enr-multicategory-right-unity,enr-multicategory-left-unity}, of a $\V$-multicategory.

Suppose $\V$ is the symmetric monoidal category $(\Cat, \times, \boldone)$ of small categories and functors with the Cartesian product.  Then each non-symmetric $\Cat$-multicategory has an underlying 2-category by \cref{locallysmalltwocat}.
\end{example}

\begin{example}\label{definition:terminal-operad-comm}  
With $\V = (\Set, \times, *)$, the \index{terminal!multicategory}\index{multicategory!terminal}\emph{terminal multicategory} $\Mterm$ consists of a single object $*$ and a single $n$-ary operation $\iota_n$ for each $n \ge 0$.
\end{example}

\begin{example}[Endomorphism Operad]\label{example:enr-End}
For each $\V$-multicategory $\M$ and $x \in \ObM$, the \emph{endomorphism $\V$-operad} \index{endomorphism!operad}\index{operad!endomorphism}$\End(x)$ consists of the single object $x$ and $n$-ary multimorphism object
\[\End(x)_n = \M\mmap{x;\ang{x}} \in \V,\]
with $\ang{x}$ the $n$-tuple of copies of $x$.  Its multicategory structure is inherited from $\M$.
\end{example}

\subsection*{The 2-Category of Enriched Multicategories}

\begin{definition}\label{def:enr-multicategory-functor}
Suppose $\M$ and $\N$ are $\V$-multicategories.  A \emph{$\V$-multifunctor}\index{multifunctor!enriched}\index{enriched!multifunctor}\index{functor!multi-} 
\[F \cn \M \to \N\]
consists of the following data.
\begin{description}
\item[Object Assignment] It is equipped with a function 
\[F \cn \ObM \to \ObN.\]
\item[Component Morphisms] 
For $\smscmap{\angx;y} \in \ProfMM$ with $\angx= \ang{x_j}_{j=1}^n$, it is equipped with a morphism 
\[\begin{tikzcd}[column sep=large]
\M\mmap{y;\ang{x}} \ar{r}{F} & \N\mmap{Fy;F\ang{x}}
\end{tikzcd} \inspace \V\]
where $F\angx = \ang{Fx_j}_{j=1}^n$.
\end{description}
The data above are required to satisfy the following axioms.
\begin{description}
\item[Symmetric Group Action] 
For $\smscmap{\angx;y} \in \ProfMM$ and $\sigma \in \Sigma_n$, the following diagram\index{equivariance!enriched multifunctor} in $\V$ commutes.
\begin{equation}\label{enr-multifunctor-equivariance}
\begin{tikzcd}[column sep=large]
\M\mmap{y;\ang{x}} \ar{d}{\cong}[swap]{\sigma} \ar{r}{F} 
& \N\mmap{Fy;F\ang{x}} \ar{d}{\cong}[swap]{\sigma}\\
\M\mmap{y;\ang{x}\sigma} \ar{r}{F} & \N\mmap{y;F\ang{x}\sigma}\end{tikzcd}
\end{equation}
\item[Units] 
For each $x \in \ObM$, the following diagram in $\V$ commutes.
\begin{equation}\label{enr-multifunctor-unit}
\begin{tikzpicture}[x=25mm,y=15mm,vcenter]
  \draw[0cell] 
  (0,0) node (a) {\tu}
  (1,.5) node (b) {\M\mmap{x;x}}
  (1,-.5) node (b') {\N\mmap{Fx;Fx}}
  ;
  \draw[1cell] 
  (a) edge node {\operadunit_x} (b)
  (a) edge node[swap,pos=.6] {\operadunit_{Fx}} (b')
  (b) edge node {F} (b')
  ;
\end{tikzpicture}
\end{equation} 
\item[Composition] 
For $x''$, $\ang{x'}$, and $\ang{x} = \oplus_{j=1}^n\ang{x_j}$ as in \eqref{eq:enr-defn-gamma}, the following diagram in $\V$ commutes.
\begin{equation}\label{v-multifunctor-composition}
\begin{tikzcd}[column sep=large,cells={nodes={scale=.8}},
every label/.append style={scale=.9}]
\M\mmap{x'';\ang{x'}} \otimes \bigotimes_{j=1}^n \M\mmap{x'_j;\ang{x_j}} \dar[swap]{\gamma} \ar{r}{(F,\otimes_j F)} & \N\mmap{Fx'';F\ang{x'}} \otimes \bigotimes_{j=1}^n \N\mmap{Fx'_j;F\ang{x_j}} \dar{\gamma}\\  
\M\mmap{x'';\ang{x}} \ar{r}{F} & \N\mmap{Fx'';F\ang{x}}
\end{tikzcd}
\end{equation}
\end{description}
This finishes the definition of a $\V$-multifunctor.  

Moreover, we define the following.
\begin{itemize}
\item For a $\V$-multifunctor $G \cn \N\to\P$, the \index{composition!enriched multifunctor}\emph{composition} 
\[GF \cn \M\to\P\]
is the $\V$-multifunctor with object assignment given by the composite function 
\[\begin{tikzcd} \ObM \ar{r}{F} & \ObN \ar{r}{G} & \ObP
\end{tikzcd}\]
and component morphisms given by the composites
\begin{equation}\label{composite-multifunctor}
\begin{tikzcd}
\M\mmap{y;\ang{x}} \ar{r}{F} & \N\mmap{Fy;F\ang{x}} \ar{r}{G} & \P\mmap{GFy;GF\ang{x}}.
\end{tikzcd}
\end{equation}
\item The \index{identity!enriched multifunctor}\emph{identity $\V$-multifunctor} $1_{\M} \cn \M\to\M$ has the identity object assignment and identity component morphisms.
\item A \emph{$\V$-operad morphism}\index{operad!morphism}\index{enriched operad!morphism}\index{morphism!enriched operad} is a $\V$-multifunctor between $\V$-multicategories with one object.
\item A \emph{non-symmetric $\V$-multifunctor}\index{non-symmetric!multifunctor}\index{multifunctor!non-symmetric} $F \cn \M \to \N$ between non-symmetric $\V$-multicategories is defined in the same way as a $\V$-multifunctor but without the symmetric group action axiom \cref{enr-multifunctor-equivariance}.  Composition and identities are defined as above.
\item A \emph{(non-symmetric) multifunctor} is a (non-symmetric) $\Set$-multifunctor.\defmark
\end{itemize}
\end{definition}

\begin{example}[Underlying $\V$-Functors]\label{ex:un-v-functor}
Continuing \cref{ex:unarycategory}, each (non-symmetric) $\V$-multifunctor (\cref{def:enr-multicategory-functor}) restricts to a $\V$-functor between the underlying $\V$-categories.  The compatibility diagrams \cref{eq:enriched-composition} of a $\V$-functor are the unit diagram \cref{enr-multifunctor-unit} and the unary restriction of the composition diagram \cref{v-multifunctor-composition}.
\end{example}

\begin{definition}\label{def:enr-multicat-natural-transformation}
Suppose $F,G \cn \M\to\N$ are $\V$-multifunctors.  A \emph{$\V$-multinatural transformation}\index{enriched!multinatural transformation}\index{natural transformation!enriched multi-}\index{multinatural transformation!enriched} $\theta \cn F \to G$ consists of, for each $x \in \ObM$, a component morphism 
\[\begin{tikzcd}[column sep=large]
\tu \ar{r}{\theta_x} & \N\mmap{Gx;Fx}
\end{tikzcd} \inspace \V\]
such that the following \emph{$\V$-naturality diagram} commutes for $\smscmap{\ang{x};y} \in \ProfMM$ with $\angx= \ang{x_j}_{j=1}^n$.
\begin{equation}\label{enr-multinat}
\begin{tikzpicture}[x=22mm,y=12mm,baseline={(a.base)}]
  \draw[0cell=.8]
  (.5,0) node (a) {\M\mmap{y;\ang{x}}}
  (.7,1) node (b) {\tu \otimes \M\mmap{y;\ang{x}}}
  (3.3,1) node (c) {\N\mmap{Gy;Fy} \otimes \N\mmap{Fy;F\ang{x}}}
  (3.5,0) node (d) {\N\mmap{Gy;F\ang{x}}}
  (.7,-1) node (b') {\M\mmap{y;\ang{x}}\otimes \txotimes_{j=1}^n \tu}
  (3.3,-1) node (c') {\N\mmap{Gy;G\ang{x}} \otimes \txotimes_{j=1}^n \N\mmap{Gx_j;Fx_j}}
  ;
  \draw[1cell=.8] 
  (a) edge node[pos=.3] {\la^\inv} (b)
  (a) edge node[swap,pos=.1] {\rho^\inv} (b')
  (b) edge node {\theta_{y} \otimes F} (c)
  (c) edge node[pos=.7] {\ga} (d)
  (b') edge node {G \otimes \txotimes_{j=1}^n \theta_{x_j}} (c')
  (c') edge['] node[pos=.7] {\ga} (d)
  ;
\end{tikzpicture}
\end{equation}
This finishes the definition of a $\V$-multinatural transformation.  

Moreover, we define the following.
\begin{itemize}
\item The \emph{identity $\V$-multinatural transformation}\index{multinatural transformation!identity} $1_F \cn F\to F$ has each component given by a colored unit: 
\[(1_F)_x = \operadunit_{Fx} \forspace x\in\ObM.\]
\item A \emph{multinatural transformation} is a $\Set$-multinatural transformation.
\item A \emph{$\V$-multinatural transformation}\index{multinatural transformation!non-symmetric}\index{non-symmetric!multinatural transformation} $\theta \cn F \to G$ between non-symmetric $\V$-multifunctors $F, G \cn \M \to \N$ is defined as above.  In this case, we also call $\theta$ a \emph{non-symmetric $\V$-multinatural transformation} if we want to emphasize that its domain and codomain are non-symmetric $\V$-multifunctors.
\end{itemize} 
We use the 2-cell notation \cref{twocellnotation} for (non-symmetric) $\V$-multinatural transformation.
\end{definition}

\begin{example}[Underlying $\V$-Natural Transformations]\label{ex:un-v-nat}
Continuing \cref{ex:unarycategory,ex:un-v-functor}, each $\V$-multinatural transformation (\cref{def:enr-multicat-natural-transformation}) is also a $\V$-natural transformation (\cref{def:enriched-natural-transformation}) between the underlying $\V$-functors between the underlying $\V$-categories.  The naturality diagram \cref{enr-naturality} of a $\V$-natural transformation is the unary restriction of the naturality diagram \cref{enr-multinat} of a $\V$-multinatural transformation.
\end{example}

\begin{definition}\label{def:enr-multinatural-composition}
Suppose $\M$, $\N$, and $\P$ are $\V$-multicategories.
\begin{enumerate}
\item Suppose $\theta$ and $\psi$ are $\V$-multinatural transformations as in the left diagram below. 
\begin{equation}\label{multinatvcomp}

\end{equation} 
The \emph{vertical composition}\index{vertical composition!enriched multinatural transformation}\label{notation:enr-operad-vcomp} $\psi\theta$, as in the right diagram above, is the $\V$-multinatural transformation with component at each $x \in \ObM$ given by the following composite in $\V$.
\begin{equation}\label{multinat-vcomp-def}
%
\end{equation}
\item Suppose $\theta$ and $\theta'$ are $\V$-multinatural transformations as in the left diagram below.  
\begin{equation}\label{multinathcomp}
%
\end{equation} 
The \emph{horizontal composition}\index{horizontal composition!enriched multinatural transformation}\label{notation:enr-operad-hcomp} $\theta' \ast \theta$, as in the right diagram above, is the $\V$-multinatural transformation with component at each $x \in \ObM$ given by the following composite in $\V$.
\begin{equation}\label{multinat-hcomp-def}
%
\end{equation}
\end{enumerate}
Vertical and horizontal compositions of non-symmetric $\V$-multinatural transformations are defined as above.
\end{definition}

\begin{theorem}\label{v-multicat-2cat}
Suppose $\V$ is a symmetric monoidal category.
\begin{enumerate}
\item\label{vmulti-i} There is a 2-category\index{2-category!of small enriched multicategories}\index{enriched multicategory!2-category}\index{multicategory!enriched!2-category} 
\[\VMulticat\]
consisting of the following data.
\begin{itemize}
\item Its objects are small $\V$-multicategories.
\item For small $\V$-multicategories $\M$ and $\N$, the hom category 
\[\VMulticat(\M,\N)\]
is defined as follows.
\begin{itemize}
\item Its objects are $\V$-multifunctors $\M \to \N$.
\item Its morphisms are $\V$-multinatural transformations.
\item Identity morphisms are identity $\V$-multinatural transformations.
\item Composition is vertical composition of $\V$-multinatural transformations.
\end{itemize}
\item The identity 1-cell $1_{\M}$ is the identity $\V$-multifunctor $1_{\M}$.
\item Horizontal composition of 1-cells is the composition of $\V$-multifunctors.
\item Horizontal composition of 2-cells is that of $\V$-multinatural transformations.
\end{itemize}
\item\label{vmulti-ii} There is an analogous 2-category 
\[\VMulticatns\]
with
\begin{itemize}
\item non-symmetric small $\V$-multicategories as objects,
\item non-symmetric $\V$-multifunctors as 1-cells, and
\item non-symmetric $\V$-multinatural transformations as 2-cells.
\end{itemize} 
\item\label{vmulti-iii} Suppose, furthermore, $\V$ is a complete and cocomplete symmetric monoidal closed category.  Then the underlying 1-categories of $\VMulticat$ and $\VMulticatns$ are complete and cocomplete.  
\end{enumerate} 
\end{theorem}

\begin{proof}
Assertions \eqref{vmulti-i} and \eqref{vmulti-iii} for $\VMulticat$ are proved using essentially the same proofs as \cite[2.4.26]{johnson-yau} and \cite[5.5.14]{cerberusIII}, respectively, which deal with the case $\V = \Set$.  The analogous statements for the non-symmetric case use the same proofs by ignoring the symmetric group action.
\end{proof}

We define
\begin{equation}\label{multicat-def}
\begin{aligned}
\Multicat &= \Set\mh\Multicat \andspace\\
\Multicatns &= \Set\mh\Multicatns
\end{aligned}
\end{equation} 
which are, respectively, $\VMulticat$ and $\VMulticatns$ with $(\V,\otimes,\tu) = (\Set,\times,*)$.  In \cref{sec:multicatclosed} we extend the 2-category $\Multicat$ to 
\begin{itemize}
\item a symmetric monoidal $\Cat$-category (\cref{theorem:multicat-symmon}) and
\item a $\Cat$-multicategory (\cref{expl:multicatcatmulticat}).
\end{itemize} 

\begin{example}[Initial and Terminal Objects]\label{ex:vmulticatinitialterminal}\index{initial!enriched multicategory}\index{enriched multicategory!initial}\index{terminal!enriched multicategory}\index{enriched multicategory!terminal}\
\begin{romenumerate}
\item\label{ex:initialoperad} With $(\V,\otimes) = (\Set,\times)$, the \emph{initial operad} $\Mtu$ has a single object $*$ and a single unit operation $\opu_* \in \Mtu_1$.
\item\label{ex:initialvmulticat} The \emph{initial $\V$-multicategory} has an empty set of objects.
\item\label{ex:terminalvmulticat} If $\V$ has a terminal object \label{not:bt}$\bt$, then a \emph{terminal $\V$-multicategory} $\Mterm$ has a single object $*$ and $n$-ary multimorphism object
\[\Mterm_n = \Mterm(\overbracket[.5pt]{*, \ldots, *}^{\text{$n$ terms}} ; *) = \bt\]
for each $n \geq 0$.\defmark
\end{romenumerate}
\end{example}

\section{Categorically-Enriched Multicategories}
\label{sec:cat-multicat}

In this section we review multicategories enriched in $\Cat$.  Recall from \cref{ex:cat} the symmetric monoidal closed category $(\Cat,\times,\boldone)$ of small categories and functors with the monoidal product given by the Cartesian product.

\begin{definition}\label{def:catmulticat}
Suppose $\angx = \ang{x_j}_{j=1}^n, y, z$ are objects in a $\Cat$-multicategory $(\M,\ga,\opu)$.  The category 
\[\M\scmap{\angx;y}\]
is called a \index{category!multimorphism}\index{multimorphism!category}\emph{multimorphism category}.
\begin{itemize}
\item An object in $\M\scmap{\angx;y}$ is called an \index{1-cell!n-ary@$n$-ary}\index{n-ary@$n$-ary!1-cell}\emph{$n$-ary 1-cell} and is denoted \label{not:naryonecell}$\angx \to y$.
\item A morphism 
\[\theta \cn f \to g \inspace \M\scmap{\angx;y}\]
is called an \index{2-cell!n-ary@$n$-ary}\index{n-ary@$n$-ary!2-cell}\emph{$n$-ary 2-cell}.  We extend the 2-cell notation \cref{twocellnotation} to $n$-ary 2-cells.
\end{itemize}
The same terminology applies to non-symmetric $\Cat$-multicategories.
\end{definition}

\begin{explanation}[$\Cat$-Multinatural Transformations]\label{expl:catmultitransformation}
Suppose $\M$ and $\N$ are $\Cat$-multicategories, and 
\[F,G \cn \M \to \N\]
are $\Cat$-multifunctors.  A \index{multinatural transformation!Cat-@$\Cat$-}\index{Cat-multinatural@$\Cat$-multinatural!transformation}\emph{$\Cat$-multinatural transformation} $\theta \cn F \to G$ consists of, for each $x \in \Ob\M$, a component 1-ary 1-cell
\begin{equation}\label{thetaccomponent}
\begin{tikzcd}[column sep=large]
Fx \ar{r}{\theta_x} & Gx
\end{tikzcd}
\inspace \N
\end{equation}
such that the following two \index{Cat-naturality conditions@$\Cat$-naturality conditions}\emph{$\Cat$-naturality conditions} hold:
\begin{description}
\item[Objects] For each $k$-ary 1-cell 
\[p \cn \angx \to x' \inspace \M\]
with $\angx = \ang{x_j}_{j=1}^k$, denote by
\begin{equation}\label{Fangcthetaangc}
\left\{\begin{aligned}
F\ang{x} &= \ang{F x_j}_{j=1}^k \in (\Ob\N)^k\\
\theta_{\ang{x}} &= \ang{\theta_{x_j}}_{j=1}^k \in \txprod_{j=1}^k \N\mmap{G x_j; F x_j}.
\end{aligned}\right.
\end{equation}
Then the following equality of $k$-ary 1-cells holds, where the composition is taken in $\N$:
\begin{equation}\label{catmultinaturality}
\gamma\scmap{Gp; \theta_{\angx}} =
\gamma\scmap{\theta_{x'}; Fp} \inspace \N\mmap{Gx';F\angx}.
\end{equation}
\item[Morphisms] For each $k$-ary 2-cell 
\[f \cn p \to q \inspace \M\mmap{x';\ang{x}},\] 
the following equality of $k$-ary 2-cells holds, with $1_{\theta_{\angx}} = \ang{1_{\theta_{x_j}}}_{j=1}^k$:
\begin{equation}\label{catmultinaturalityiicell}
\gamma\scmap{Gf; 1_{\theta_{\angx}}} =
\gamma\scmap{1_{\theta_{x'}}; Ff} \inspace \N\mmap{Gx';F\angx}.
\end{equation}
\end{description}
The conditions \cref{catmultinaturality,catmultinaturalityiicell} together comprise the $\V$-naturality condition \cref{enr-multinat} with $\V = \Cat$.

The two sides of the object $\Cat$-naturality condition \cref{catmultinaturality} use the following two compositions in $\N$ on objects.
\begin{equation}\label{catmultinaturalitycomp}
\begin{tikzpicture}[xscale=2.5,yscale=1.2,vcenter]
\draw[0cell=.85]
(0,0) node (a) {\N\mmap{Gx';G\angx} \times \txprod_{j=1}^k \N\mmap{G x_j;F x_j}}
(a)++(2,0) node (b) {\N\mmap{Gx';Fx'} \times \N\mmap{Fx'; F\angx}}
(a)++(1,-1) node (c) {\N\mmap{Gx';F\angx}}
;
\draw[1cell=.9]  
(a) edge node[swap,pos=.5] {\gamma} (c)
(b) edge node[pos=.5] {\gamma} (c)
;
\end{tikzpicture}
\end{equation}
The object $\Cat$-naturality condition \cref{catmultinaturality} is the commutative diagram
\begin{equation}\label{catmultinaturalitydiagram}
\begin{tikzcd}[column sep=large]
F\angx \ar{d}[swap]{Fp} \ar{r}{\theta_{\angx}} & G \angx \ar{d}{Gp}\\
Fx' \ar{r}{\theta_{x'}} & Gx'
\end{tikzcd}
\end{equation}
involving the compositions \cref{catmultinaturalitycomp} in $\N$.

The two sides of the morphism $\Cat$-naturality condition \cref{catmultinaturalityiicell} use the compositions in \cref{catmultinaturalitycomp} on morphisms.  The condition \cref{catmultinaturalityiicell} is the equality of \index{pasting diagram!multicategorical}multicategorical pasting diagrams
\begin{equation}\label{catmultinatiicellpasting}
\begin{tikzpicture}[xscale=2.5,yscale=1.5,vcenter]
\def\a{35} \def\s{.8} \def\h{.8}
\def\boundary{
\draw[0cell=\s]
(0,0) node (x11) {F\ang{x}}
(x11)++(\h,0) node (x12) {G\ang{x}}
(x11)++(0,-1) node (x21) {Fx'}
(x12)++(0,-1) node (x22) {Gx'}
;
\draw[1cell=\s]  
(x11) edge node {\theta_{\ang{x}}} (x12)
(x12) edge[bend left=\a] node[pos=.75] {Gq} (x22)
(x11) edge[bend right=\a] node[swap,pos=.25] {Fp} (x21)
(x21) edge node {\theta_{x'}} (x22)
;}
\draw (0,0) node[font=\LARGE] (equal) {=};
\begin{scope}[shift={(.6,.5)}]
\boundary
\draw[1cell=\s] 
(x11) edge[bend left=\a] node[pos=.4] {Fq} (x21)
;
\draw[2cell=\s]
node[between=x11 and x21 at .55, rotate=0, 2label={above,Ff}] {\Longrightarrow}
;
\end{scope}
\begin{scope}[shift={(-1.4,.5)}]
\boundary
\draw[1cell=\s] 
(x12) edge[bend right=\a] node[swap,pos=.4] {Gp} (x22)
;
\draw[2cell=\s]
node[between=x12 and x22 at .55, rotate=0, 2label={above,Gf}] {\Longrightarrow}
;
\end{scope}
\end{tikzpicture}
\end{equation}
involving the composition in $\N$.

By the object $\Cat$-naturality condition \cref{catmultinaturality} for $\theta' \cn F' \to G'$, there are two  ways to express the $x$-component of the horizontal composite \cref{multinathcomp} for $x \in \ObM$:
\begin{equation}\label{hcompexpressions}
\begin{aligned}
(\theta' * \theta)_x &= \ga\scmap{\theta'_{Gx}; F' \theta_x}\\
&= \ga\scmap{G' \theta_x; \theta'_{Fx}}.
\end{aligned}
\end{equation}
A non-symmetric $\Cat$-multinatural transformation admits the same description as above.
\end{explanation}

\section{Endomorphism Multicategories}
\label{sec:mult-from-smc}

In this section we review the endomorphism construction that goes from symmetric monoidal categories to multicategories, which mean $\Set$-multicategories.  This construction defines a 2-functor; see \cref{endtwofunctor}.  The enriched variant is in \cref{definition:EndK}.  The material in this section is adapted from \cite[Section 6.3]{cerberusIII}.

\begin{example}[Endomorphism Multicategory]\label{ex:endc}
Suppose $(\C,\oplus,\pu,\xi)$ is a permutative category (\cref{def:symmoncat}).  Then it has an associated \index{endomorphism!multicategory}\index{multicategory!endomorphism}\index{permutative category!endomorphism multicategory}\emph{endomorphism multicategory} $\End(\C)$ defined as follows.
\begin{itemize}
\item The object class is $\ObC$.
\item The $n$-ary multimorphism set is
\[\End(\C)\mmap{y;\ang{x}} = \C\big(\!\txoplus_{j=1}^n x_j \scs y\big)\]
for $\smscmap{\angx;y} \in \Profcc$ with $\ang{x} = \ang{x_j}_{j=1}^n$.  By definition, an empty $\oplus$ in $\C$ is the monoidal unit $\pu$.
\item The symmetric group action is induced by the braiding $\xi$.
\item For each object $x$ in $\C$, the $x$-colored unit is the identity morphism $1_x$.
\item The multicategorical composition in $\End(\C)$ is induced by $\oplus$ and composition in $\C$.
\end{itemize}  

There are two variants of the above endomorphism construction.
\begin{enumerate}
\item If $(\C,\otimes,\tu,\xi)$ is a symmetric monoidal category that is not necessarily strict, then $\End(\C)$ is still a multicategory.
\begin{itemize}
\item The $n$-ary multimorphism set is
\begin{equation}\label{endc-angxy}
\End(\C)\mmap{y;\ang{x}} = \C\big(\!\txotimes_{j=1}^n x_j \scs y\big)
\end{equation}
with $\txotimes_{j=1}^n x_j$ using \cref{expl:leftbracketing} of left normalized bracketing.  By definition, an empty $\otimes$ is the monoidal unit $\tu$.
\item The symmetric group action is induced by the braiding $\xi$ and the associativity isomorphism $\alpha$ in $\C$. 
\item The composition in $\End(\C)$ is induced by $\alpha$, $\otimes$, and composition in $\C$.  If nullary multimorphisms are involved, then we also use the left and right unit isomorphisms in $\C$.
\end{itemize}  
\item If $(\C,\otimes,\tu)$ is a monoidal category, then the above definitions, without the symmetric group action, yield a non-symmetric multicategory that we also denote by $\End(\C)$.
\end{enumerate}

Furthermore, the endomorphism multicategory extends to symmetric monoidal functors and monoidal natural transformations as follow.  Each symmetric monoidal functor between symmetric monoidal categories
\[(P,P^2,P^0) \cn \C \to \D\]
induces a multifunctor
\begin{equation}\label{EndP}
\End(P) \cn \End(\C) \to \End(\D)
\end{equation}
with the same object assignment as $P$.  For an $n$-ary multimorphism 
\[f \in \End(\C)\mmap{y;\ang{x}} = \C\big(\!\txotimes_{j=1}^n x_j \scs y\big)\]
with $\smscmap{\angx;y}$ as above, its image in $\End(\D)\scmap{\ang{Px};Py}$ is the composite 
\begin{equation}\label{PtwoPf}
\txotimes_{j=1}^n P x_j \fto{P^2} P\big(\!\txotimes_{j=1}^n x_j\big) \fto{Pf} Py \inspace \D.
\end{equation}
The first morphism $P^2$ in \cref{PtwoPf} means
\begin{itemize}
\item the unit constraint $P^0 \cn \tu \to P\tu$ if $n=0$,
\item the identity if $n=1$, and
\item a repeated application of the monoidal constraint $P^2$ if $n > 1$.  
\end{itemize} 

Suppose $\theta$ is a monoidal natural transformation between symmetric monoidal functors between symmetric monoidal categories, as in the left diagram below.
\[\begin{tikzpicture}
\def\h{2} \def\t{25}
\draw[0cell]
(0,0) node (a) {\C}
(a)+(\h,0) node (b) {\D}
;
\draw[1cell=.8]
(a) edge[bend left=\t] node {(P,P^2,P^0)} (b)
(a) edge[bend right=\t] node[swap] {(Q,Q^2,Q^0)} (b)
;
\draw[2cell]
node[between=a and b at .45, rotate=-90, 2label={above,\theta}] {\Rightarrow}
;
\begin{scope}[shift={(4.5,-.04)}]
\draw[0cell]
(0,0) node (a) {\phantom{X}}
(a)+(-.4,0) node (a') {\End(\C)}
(a)+(2.5,0) node (b) {\phantom{X}}
(b)+(.4,0) node (b') {\End(\D)}
;
\draw[1cell=.8]
(a) edge[bend left=\t] node {\End(P)} (b)
(a) edge[bend right=\t] node[swap] {\End(Q)} (b)
;
\draw[2cell]
node[between=a and b at .33, rotate=-90, 2label={above,\End(\theta)}] {\Rightarrow}
;
\end{scope}
\end{tikzpicture}\]
Then $\theta$ induces a multinatural transformation $\End(\theta)$, as in the right diagram above, with component morphisms
\begin{equation}\label{Endthetax}
\End(\theta)_x = \theta_x \forspace x \in \C.
\end{equation}

The non-symmetric variants of the above statements also hold.
\begin{enumerate}
\item If $P$ is a monoidal functor between monoidal categories, then $\End(P)$ is a non-symmetric multifunctor. 
\item If $\theta$ is a monoidal natural transformation between monoidal functors between monoidal categories, then $\End(\theta)$ is a non-symmetric multinatural transformation.
\end{enumerate}  
If there is no danger of confusion, we abbreviate $\End(\C)$ to $\C$ and similarly for $\End(P)$ and $\End(\theta)$.  We extend this example to the pointed context in \cref{ex:endstc} below.  
\end{example}

Recall the 2-categories
\begin{itemize}
\item $\Multicat$ of small multicategories, multifunctors, and multinatural transformations (\cref{v-multicat-2cat} with $\V = (\Set, \times, *)$) and
\item $\permcat$ of small permutative categories, symmetric monoidal functors, and monoidal natural transformations (\cref{def:permcat}).
\end{itemize}

\begin{proposition}\label{endtwofunctor}\index{endomorphism!2-functor}
The endomorphism multicategory in \cref{ex:endc} defines a 2-functor
\[\End \cn \permcat \to \Multicat.\]
\end{proposition}

We also denote by $\End$ the restriction of the 2-functor in \cref{endtwofunctor} to any one of the locally-full sub-2-categories of $\permcat$ in \cref{def:permcat}, including $\permcatsu$ and $\permcatst$.  In \cref{endufactor} we discuss an extension of $\End$, with domain $\permcatsu$, to a $\Cat$-multifunctor.

\subsection*{Enriched Endomorphism Multicategories}

For the rest of this section, we assume that $(\V,\otimes,\tu,\xi)$ is a symmetric monoidal category.  Next we review the $\V$-multicategory associated to a symmetric monoidal $\V$-category.  The endomorphism multicategory in \cref{ex:endc} is the special case $(\V,\otimes) = (\Set, \times)$.

\begin{convention}[Left Normalized Product]\label{convention:norm}\index{convention!left normalized product}
Suppose $(\K,\vmtimes)$ is a monoidal $\V$-category (\cref{definition:monoidal-vcat}), and $\ang{x} = \ang{x_j}_{j=1}^n$ is an $n$-tuple of objects of $\K$ for some $n \geq 0$.  We define the \index{left normalized!product}\index{normalized!left - product}\emph{left normalized product} as the object
\[\bigvmtimes \ang{x} = \bigvmtimes_{j=1}^n x_j = \big(\cdots \, \big((x_1 \vmtimes x_2) \vmtimes x_3 \big) \, \cdots \big) \vmtimes x_n,\]
which is the identity object of $\K$ if $n=0$.
\end{convention}

\begin{definition}\label{definition:EndK}
Suppose $(\K,\beta^\vmtimes)$ is a symmetric monoidal $\V$-category (\cref{definition:symm-monoidal-vcat}) with $\V$ a symmetric monoidal category.  The \index{endomorphism!enriched multicategory}\index{enriched multicategory!endomorphism}\emph{endomorphism $\V$-multicategory} of $\K$, denoted $\End(\K)$, consists of the following data.
\begin{description}
\item[Objects] $\Ob(\End(\K)) = \Ob\K$.
\item[Multimorphism Objects] 
For an object $x' \in \K$ and a tuple $\ang{x} \in \Prof(\K)$, we define the object
\[\End(\K)\mmap{x';\ang{x}} = \K\big( \bigvmtimes \ang{x} \scs x' \big) \inspace \V\]
with $\bigvmtimes \ang{x}$ denoting the left normalized product in \cref{convention:norm}.
\item[Symmetric Group Action]
For objects $\ang{x}$ and $x'$ as above and a permutation $\sigma \in \Si_n$, the right $\sigma$-action is defined as the following composite.
\[
      \ \qquad
      \begin{tikzpicture}[x=50mm,y=15mm,vcenter]
      \def\s{.9}
        \draw[0cell=\s] 
        (0,0) node (a) {
          \K \big( \bigvmtimes \ang{x}, x' \big)
        }
        (0,-1) node (b) {
          \K \big( \bigvmtimes \ang{x}, x' \big)
          \otimes
          \tensorunit
        }
        (1,-1) node (c) {
          \K \big( \bigvmtimes \ang{x}, x' \big)
          \otimes
          \K \big( \bigvmtimes \ang{x}\si, \bigvmtimes \ang{x} \big)
        }
        (1,0) node (d) {
          \K \big( \bigvmtimes \ang{x}\si, x' \big)
        }
        ;
        \draw[1cell=\s] 
        (a) edge['] node {\rho^\inv} (b)
        (b) edge node {1 \otimes \beta^\vmtimes_\si} (c)
        (c) edge['] node {\mcomp} (d)
        (a) edge node {\si} (d)
        ;
      \end{tikzpicture}
    \]
    In the above diagram, $\beta^\vmtimes_\si$ denotes the $\V$-natural
    isomorphism that permutes coordinates according to the
    permutation $\si$.
\item[Units] For each object $x \in \K$, the $x$-colored unit of $\End(\K)$ is the identity of $x$ in $\K$,
    \[\begin{tikzcd}[column sep=large]
    \tensorunit \ar{r}{1_x} & \K(x,x) = \End(\K)\mmap{x;x}.
    \end{tikzcd}\]    
\item[Composition] The composition $\ga$ in $\End(\K)$ is defined as the following composite for tuples of objects
\[x'', \quad \ang{x'} = \ang{x'_j}_{j=1}^n \scs \andspace
\ang{x_j} = \ang{x_{j,i}}_{i=1}^{k_j}\] 
with $j\in\{1,\ldots,n\}$ and $\ang{x} = \ang{\ang{x_j}}_{j=1}^n$. 
\[\begin{tikzpicture}[x=50mm,y=15mm,vcenter]
\def\s{.85}
        \draw[0cell=\s] 
        (0,0) node (a) {
          \K \big( \bigvmtimes\ang{x'},x'' \big) \otimes
          {\textstyle\bigotimes\limits_{j=1}^n} \K \Big( \bigvmtimes_{i=1}^{k_j} x_{j,i} \, , \, x_j' \Big)
        }
        (0,-1) node (b) {
          \K \big( \bigvmtimes\ang{x'},x'' \big) \otimes
          \K \Big( \bigvmtimes_{j=1}^n\bigvmtimes_{i=1}^{k_j} x_{j,i} \, , \,
          \bigvmtimes \ang{x'} \Big)
        }
        (1,-1) node (c) {
          \K \Big( \bigvmtimes_{j=1}^n\bigvmtimes_{i=1}^{k_j} x_{j,i} \, , \, x'' \Big)
        }
        (1,0) node (d) {
          \K \big( \bigvmtimes \ang{x},x'' \big)
        }
        ;
        \draw[1cell=\s] 
        (a) edge[transform canvas={xshift={-2ex}}, shorten <=-1ex, shorten >=-1ex] node[swap] {1 \otimes \vmtimes_{j=1}^{n-1}} (b)
        (b) edge node {\mcomp} (c)
        (c) edge node['] {\iso} (d)
        (a) edge node {\ga} (d)
        ;
\end{tikzpicture}\]
\end{description}
This finishes the definition of $\End(\K)$.  To simplify the notation, we also denote $\End(\K)$ by $\K$.
\end{definition}

The following result is \cite[6.3.6]{cerberusIII}.

\begin{proposition}\label{proposition:monoidal-v-cat-v-multicat}
In the context of \cref{definition:EndK}, $\End(\K)$ is a $\V$-multicategory.
\end{proposition}

If there is no danger of confusion, we abbreviate $\End(\K)$ to $\K$.

\section{Pointed Multicategories}
\label{sec:ptmulticatiicat}

Recall that a \emph{multicategory} means a $\Set$-multicategory (\cref{def:enr-multicategory}), where $(\Set, \times, *)$ is the symmetric monoidal category of sets with the Cartesian product.  In this section we review \emph{pointed} multicategories.  
\begin{itemize}
\item The 2-category $\pMulticat$ of small pointed multicategories is in \cref{thm:pmulticat}.
\item The free-forgetful 2-adjunction between $\Multicat$ and $\pMulticat$ is in \cref{dplustwofunctor}.
\end{itemize}
The material in this section is adapted from \cite[Section 5.3]{cerberusIII}.  Recall the terminal multicategory $\Mterm$ in \cref{definition:terminal-operad-comm}.  It has a single object $*$ and one $n$-ary operation $\iota_n$ for each $n \geq 0$.

\begin{definition}\label{def:ptd-multicat}
We define the following.
\begin{enumerate}
\item A \index{pointed!multicategory}\index{multicategory!pointed}\emph{pointed multicategory} $(\M,i)$ is a pair consisting of the following data.
\begin{itemize}
\item $\M$ is a multicategory (\cref{def:enr-multicategory}).
\item $i \cn \Mterm \to \M$ is a multifunctor (\cref{def:enr-multicategory-functor}), which is called the \emph{pointed structure}.
\end{itemize} 
We denote
\begin{itemize}
\item $i(*) \in \ObM$ by $*$, which is called the \emph{basepoint object}, and
\item $i(\iota_n) \in \M\scmap{\ang{*}_{j=1}^n;*}$ by $\iota_n$ or $\iota^n$, which is called the \emph{$n$-ary basepoint operation}, for each $n \geq 0$.
\end{itemize} 
\item For pointed multicategories $(\M,i^\M)$ and $(\N,i^\N)$, a \index{pointed!multifunctor}\index{multifunctor!pointed}\emph{pointed multifunctor}
\[F \cn (\M,i^\M) \to (\N,i^\N)\]
is a multifunctor $F \cn \M \to \N$ (\cref{def:enr-multicategory-functor}) such that the following diagram of multifunctors commutes.
\begin{equation}\label{ptmultifunctor-diagram}
\begin{tikzpicture}[xscale=1,yscale=1,vcenter]
\draw[0cell]
(0,0) node (a) {\Mterm}
(a)++(2,.7) node (m) {\M}
(a)++(2,-.7) node (n) {\N}
;
\draw[1cell]
(a) edge node[pos=.6] {i^\M} (m)
(a) edge node[swap,pos=.6] {i^\N} (n)
(m) edge node {F} (n)
;
\end{tikzpicture}
\end{equation} 
\item A \index{pointed!multinatural transformation}\index{multinatural transformation!pointed}\emph{pointed multinatural transformation} 
\[\theta \cn F \to G \cn (\M,i^\M) \to (\N,i^\N)\]
between pointed multifunctors $F$ and $G$ is a multinatural transformation (\cref{def:enr-multicat-natural-transformation}) such that the basepoint component 
\begin{equation}\label{ptmntr}
\theta_{*} = \opu_* \in \N\scmap{F(*); G(*)} = \N\scmap{*;*},
\end{equation}
which is the colored unit of the basepoint object $*$ in $\N$.
\end{enumerate}
With these definitions, composites of pointed multifunctors and multinatural transformations are again pointed.
\end{definition}

\begin{explanation}[Pointed Structure]\label{expl:ptmulticat}
For a multicategory $(\M,\ga,\opu)$, a multifunctor $i \cn \Mterm \to \M$ (\cref{def:enr-multicategory-functor}) is uniquely determined by
\begin{itemize}
\item a \emph{basepoint object} $* \in \ObM$ and
\item an \emph{$n$-ary basepoint operation} 
\[\iota_n \in \M\scmap{\ang{*}_{j=1}^n;*} \forspace \geq 0\]
\end{itemize}
such that the following three conditions hold:
\begin{description}
\item[Symmetry] For each $n \geq 0$ and permutation $\sigma \in \Sigma_n$, there is an equality
\begin{equation}\label{ptstructure-sym}
\iota_n \cdot \sigma = \iota_n.
\end{equation}
\item[Unity] There is an equality 
\begin{equation}\label{ptstructure-unity}
\iota_1 = \opu_* \in \M\smscmap{*;*},
\end{equation}
which is the colored unit of $*$ in $\M$.
\item[Composition] For $n \geq 1$ and $k_j \geq 0$ for $j \in \{1,\ldots,n\}$, there is an equality
\begin{equation}\label{ptstructure-comp}
\ga\scmap{\iota_n; \ang{\iota_{k_j}}_{j=1}^n} = \iota_{k_1 + \cdots + k_n}.
\end{equation}
\end{description}

Moreover, for a pointed multifunctor 
\[F \cn (\M,i^\M) \to (\N,i^\N),\] 
the commutative diagram \cref{ptmultifunctor-diagram} means
\begin{itemize}
\item $F(*) = *$ in $\ObN$ and
\item $F(\iota_n) = \iota_n$ for each $n \geq 0$.
\end{itemize}
In other words, a pointed multifunctor is a multifunctor that preserves the basepoint object and the basepoint operations.
\end{explanation}

\begin{example}[Pointed Endomorphism Multicategory]\label{ex:endstc}\index{endomorphism!pointed multicategory}
Each permutative category $(\C,\oplus,\pu,\xi)$ has an associated pointed multicategory
\[\Endst(\C) = \big(\End(\C), i\big)\]
defined as follows.
\begin{itemize}
\item $\End(\C)$ is the endomorphism multicategory in \cref{ex:endc}.
\item Using \cref{expl:ptmulticat}, the pointed structure is given by the multifunctor
\[i \cn \Mterm \to \End(\C)\]
determined by
\begin{itemize}
\item the basepoint object $\pu \in \C$ and
\item $n$-ary basepoint operations
\[\iota^n = 1_{\pu} \cn \txoplus_{j=1}^n \pu = \pu \to \pu \forspace n \geq 0.\]
\end{itemize}
\end{itemize}

As in \cref{ex:endc}, the above definitions still yield a pointed multicategory $\Endst(\C)$ if $(\C,\otimes,\tu)$ is a symmetric monoidal category that is not necessarily strict.  In this case, the $n$-ary basepoint operation
\[\iota^n \cn \txotimes_{j=1}^n \tu \fto{\iso} \tu\]
is
\begin{itemize}
\item the identity $1_{\tu}$ if $n=0$ or if $n=1$ and
\item an iterate of the right unit isomorphism $\rho$ in $\C$ if $n > 1$.
\end{itemize} 
The unity condition \cref{ptstructure-unity} holds by definition.  The symmetry and composition conditions, \cref{ptstructure-sym,ptstructure-comp}, hold by the Coherence Theorem for symmetric monoidal categories \cite[XI.1 Theorem 1]{maclane}.

Moreover, the following statements hold:
\begin{romenumerate}
\item\label{EndstP} Each \emph{strictly unital} symmetric monoidal functor between symmetric monoidal categories
\[(P,P^2,P^0=1) \cn \C \to \D\]
induces a pointed multifunctor
\[\Endst(P) \cn \Endst(\C) \to \Endst(\D)\]
given by the multifunctor $\End(P)$ in \cref{EndP}.  Strict unity of $P$ ensures that $\End(P)$ is pointed as in \cref{expl:ptmulticat}.
\item\label{Endst-theta} Each monoidal natural transformation between strictly unital symmetric monoidal functors between symmetric monoidal categories
\[\theta \cn (P,P^2,P^0=1) \to (Q,Q^2,Q^0=1) \cn \C \to \D\]
induces a pointed multinatural transformation
\[\Endst(\theta) \cn \Endst(P) \to \Endst(Q) \cn \Endst(\C) \to \Endst(\D)\]
with components as in \cref{Endthetax}:
\[\Endst(\theta)_x = \theta_x \forspace x \in \C.\]
The pointed condition
\[\Endst(\theta)_\tu = \theta_\tu = 1_\tu \cn \tu \to \tu \inspace \D\]
follows from
\begin{itemize}
\item the left diagram in \cref{monnattr} and
\item the assumption that both unit constraints $P^0$ and $Q^0$ are the identities.
\end{itemize} 
\end{romenumerate}
If there is no danger of confusion, we denote $\Endst(\C)$ by $\C$.  We extend this example to left $\Mone$-modules in \cref{ex:endmc}.
\end{example}

Pointed multifunctors are multifunctors satisfying the extra property \cref{ptmultifunctor-diagram}.  Likewise, pointed multinatural transformations are multinatural transformations with the extra property \cref{ptmntr}.  These extra properties are closed under the compositions of the 2-category $\Multicat$ (\cref{v-multicat-2cat}).  Therefore, the proof for the existence of the 2-category $\Multicat$ also yields the following.

\begin{theorem}\label{thm:pmulticat}
In the context of \cref{def:ptd-multicat}, there is a 2-category 
\[\pMulticat\]
defined by the following data.
\begin{itemize}
\item The objects are small pointed multicategories.
\item The 1-cells are pointed multifunctors.
\item The 2-cells are pointed multinatural transformations.
\item The horizontal and vertical compositions and identity 1-cells and 2-cells are defined as in the 2-category $\Multicat$.
\end{itemize}
\end{theorem}

In \cref{sec:ptmulticatclosed} we discuss extensions of the 2-category $\pMulticat$ to 
\begin{itemize}
\item a symmetric monoidal $\Cat$-category (\cref{thm:pmulticat-smclosed}) and
\item a $\Cat$-multicategory (\cref{expl:ptmulticatcatmulticat}).
\end{itemize} 

\begin{proposition}\label{endsttwofunctor}\index{endomorphism!2-functor}
The pointed endomorphism multicategory in \cref{ex:endstc} defines a 2-functor
\[\Endst \cn \permcatsu \to \pMulticat\]
with $\permcatsu$ the 2-category in \cref{def:permcat}.
\end{proposition}

We also denote by $\Endst$ the restriction of the 2-functor in \cref{endsttwofunctor} to the locally-full sub-2-category $\permcatst$ in \cref{def:permcat}.  In \cref{endufactor} we discuss an extension of $\Endst$ to a $\Cat$-multifunctor.

\cref{Usttwofunctor,dplustwofunctor} follow directly from the definitions.

\begin{proposition}\label{Usttwofunctor}\index{pointed!multicategory!forgetful functor}
There is a forgetful 2-functor
\[\Ust \cn \pMulticat \to \Multicat\]
that sends
\begin{itemize}
\item a small pointed multicategory $(\M,i)$ to the multicategory $\M$,
\item a pointed multifunctor $F$ to the multifunctor $F$, and
\item a pointed multinatural transformation $\theta$ to the multinatural transformation $\theta$. 
\end{itemize}
\end{proposition}

\begin{explanation}\label{expl:EndstUstEnd}
The 2-functors in \cref{endtwofunctor,endsttwofunctor,Usttwofunctor} yield a commutative diagram
\[\begin{tikzpicture}[baseline={(a.base)}]
\def\u{.6}
\draw[0cell]
(0,0) node (a') {\permcatsu}
(a')+(.65,-.05) node (a) {\phantom{X}}
(a)+(2.5,0) node (b) {\pMulticat}
(b)+(2.8,0) node (c) {\Multicat}
;
\draw[1cell=.9]
(a) edge node {\Endst} (b)
(b) edge node {\Ust} (c)
;
\draw[1cell=.9]
(a') [rounded corners=3pt] |- ($(b)+(-1,\u)$)
-- node {\End} ($(b)+(1,\u)$) -| (c)
;
\end{tikzpicture}\]
with $\End$ restricted to $\permcatsu$.
\end{explanation}

The forgetful 2-functor $\Ust$ admits a left 2-adjoint defined as follows.

\begin{definition}[Adjoining a Basepoint]\label{def:disjoint-basept}\index{multicategory!adjoined basepoint}
We define a 2-functor
\begin{equation}\label{dplus-def}
\dplus \cn \Multicat \to \pMulticat
\end{equation}
together with 2-natural transformations 
\begin{equation}\label{etast-epzst-def}
\begin{aligned}
\etaplus & \cn 1_{\Multicat} \to \Ust \circ~ \dplus \andspace\\
\epzplus & \cn \dplus \circ~ \Ust \to 1_{\pMulticat}
\end{aligned}
\end{equation}
as follows.
\begin{description}
\item[The 2-Functor $\dplus$]
This is defined by the following assignments.
\begin{itemize}
\item For a small multicategory $\M$, we define the pointed multicategory
\[\M_+ = \M \bincoprod \Mterm\]
with pointed structure given by the $\Mterm$ summand, where the coproduct is taken in $\Multicat$.
\item For a multifunctor $F \cn \M \to \N$ between small multicategories, we define the pointed multifunctor
\[F_+ = F \bincoprod 1_{\Mterm} \cn \M_+ \to \N_+.\]
\item Suppose $\theta \cn F \to G$ is a multinatural transformation for multifunctors $F,G \cn \M \to \N$ between small multicategories.  The pointed multinatural transformation
\[\theta_+ \cn F_+ \to G_+\]
is given by
\begin{itemize}
\item the same components as $\theta$ for objects in $\M$ and
\item $\opu_*$ in $\N_+$ for the basepoint object $* \in \M_+$.
\end{itemize}
\end{itemize}
\item[Unit] For a small multicategory $\M$, the unit $\etaplus$ has component multifunctor
\[\etaplus_{\M} \cn \M \to \Ust(\M_+) = \M \bincoprod \Mterm\]
given by the inclusion of the $\M$ summand.
\item[Counit] For a small pointed multicategory $(\N,i)$, the counit $\epzplus$ has component pointed multifunctor
\[\epzplus_{\N} = \big(1_{\N} \scs i\big) \cn \big(\Ust(\N,i)\big)_+ = \N \bincoprod \Mterm \to \N.\]
\end{description}
This finishes the definition.
\end{definition}

Recall the notion of a 2-adjunction in \cref{def:twoadjunction}.

\begin{proposition}\label{dplustwofunctor}
In the context of \cref{def:disjoint-basept}, there is a 2-adjunction
\[\big(\dplus, \Ust, \etaplus, \epzplus\big) \cn \Multicat \to \pMulticat.\]
\end{proposition}

\chapter{Open Questions}
\label{ch:questions}


In this chapter, we discuss open questions related to the topics of this work.

\begin{question}[Diagrams and Presheaves on $\GE$]\label{question:GE-duality}
  Considering the \index{Burnside!2-category}\index{2-category!Burnside}Burnside 2-category $\GE$ (\cref{definition:Burnside2}),
  \cref{remark:GE-non-self-dual} notes that the assignment
  \[
    (f,g) \mapsto (g,f),
  \]
  sending a span to its reverse, is not functorial on 1-cells.
  However, it may be a pseudofunctor.
  That possibility raises the following question.

  Is there an equivalence of homotopy theories
  \[
    \permcatsucat\big(\GE, \permcatsu\big) \to 
    \permcatsucat\big(\GE^\op, \permcatsu\big)
  \]
  with respect to the stable equivalences \cref{perm-steq}?
  Such an equivalence, combined with a change of enrichment as in \cref{thm:Kemdg}, would give another approach to $G$-spectra and would further inform the discussion in \cref{rk:BO7.5}. 

  The approach to equivalences of homotopy theories that we use throughout this work requires an underlying (1-)functor, as in \cref{def:inverse-heq}.
  Therefore, answering the above question in the affirmative appears to require new and likely interesting extensions of that basic approach.
  Note, moreover, that this question is a special case of \cref{qu:morita} below.
\end{question} 

\begin{question}[The Forgetful $\Ust$]\label{question:Ust-rel}\index{pointed!multicategory!forgetful functor}
  Considering the stable equivalences
  \[
    \cSst \subset \pMulticat
    \andspace
    \cSF \subset \Multicat
  \]
  in \cref{ptmulticat-page-v-summary}, is the forgetful functor
  \[
    \Ust \cn \pMulticat \to \Multicat
  \]
  a relative functor?
  If so, then the equality $\End = \Ust \circ \Endst$, together with the other results in \cref{ch:ptmulticat-sp}, implies that $\Ust$ is also an equivalence of homotopy theories.
  In definition of $\etast_\M$, \cref{etast-def} suppresses the forgetful $\Ust$.
  Including it, \cref{etast-def} defines $\etast_\M$ as the pointed multifunctor such that
  \[
    \Ust\etast_\M = \End(\pst_\M) \circ \eta_{\Ust\M}.
  \]
  Observe, by \cref{thm:F-heq,ptmulticat-thm-x}, that
  \begin{itemize}
  \item $\eta_{\Ust\M}$ and $\etast_\M$ are stable equivalences, and
  \item $\End(\pst_\M)$ is an $\Fr$-stable equivalence if and only if $\pst_\M$ is a stable equivalence.
  \end{itemize}
  Therefore, it follows that $\Ust$ is a relative functor if and only if $\pst_{\M}$ is a stable equivalence for each pointed multicategory $\M$.
  
  For
  \[
    H \cn \M \to \N \in \pMulticat,
  \]
  there is a commutative diagram induced by naturality of $\pst$ (\cref{pst-iinatural}).
  \[
    \begin{tikzpicture}[x=30mm,y=15mm]
      \draw[0cell] 
      (0,0) node (a) {\Fr\T}
      (a)+(1,0) node (b) {\Fr\M}
      (b)+(1,0) node (c) {\Fst\M}
      (a)+(0,-1) node (a') {\Fr\T}
      (b)+(0,-1) node (b') {\Fr\N}
      (c)+(0,-1) node (c') {\Fst\N}
      ;
      \draw[1cell] 
      (a) edge['] node {1} (a')
      (b) edge node {\Fr H} (b')
      (c) edge node {\Fst H} (c')
      (a) edge node {\Fr\iota^\M} (b)
      (b) edge node {\pst_\M} (c)
      (a') edge node {\Fr\iota^\N} (b')
      (b') edge node {\pst_\N} (c')
      ;
    \end{tikzpicture}
  \]
  Since the empty tuple is an initial object in $\Fr\T$, the nerve $\Ner(\Fr\T)$ is contractible.
  A quasifibration argument, such as Quillen's Theorem A \cite{quillenKI}, might be one way to approach this question.
  
  Another approach might try to show directly that $\pst_\M$ is a stable equivalence for each $\M$.
  For this second approach, show that there is an equivalence of categories
  \[
    \Fr(\M_+) \hty \Fr\M
  \]
  given by deleting the (disjoint) basepoint objects and operations from each tuple of objects and morphisms in $\Fr(\M_+)$.
  This, combined with the isomorphism
  \[
    \Fst(\M_+) \iso \Fr\M
  \]
  from \cref{FFstdplus}, implies that
  \[
    \pst_{(\M_+)}\cn \Fr(\M_+) \fto{\hty} \Fst(\M_+)
  \]
  is an equivalence of categories.
  Further work will be needed to determine whether or not there is a stable equivalence between 
  \[
    \Fst(\M) \andspace \Fst(\M_+)
  \]
  for all small pointed multicategories $\M$.
  The authors are not aware of either a proof or a counterexample to such an equivalence.
\end{question}

\begin{question}[$K$-Theoretic Equivalences of Homotopy Theories]\label{qu:Kemdg}
Are the diagram and presheaf change-of-enrichment functors $\Jtdg$, $\Nerstardg$, $\Kgdg$, and $\Kemdg$ in \cref{Kemdg-factor-four,Kemdg-factor-four-op} equivalences of homotopy theories (\cref{definition:rel-cat-pow})?
\end{question}

\begin{question}[Morita Theory for Closed Multicategories]\label{qu:morita}
Develop \index{Morita theory}Morita theory for diagram categories enriched in a non-symmetric closed multicategory $\M$ (\cref{def:enr-diag-cat}).  In other words, for $\M$-categories $\C$ and $\D$, give criteria that guarantee that the categories
\[\MCat(\C,\M) \andspace \MCat(\D,\M)\]
are
\begin{romenumerate}
\item equivalent or
\item connected by equivalences of homotopy theories.
\end{romenumerate}   
Moreover, for a closed multicategory $\M$, we ask the same questions for the enriched presheaf categories
\[\MCat(\Cop,\M) \andspace \MCat(\Dop,\M)\]
with $\Cop$ and $\Dop$ the opposite $\M$-categories of $\C$ and $\D$ (\cref{opposite-mcat}), respectively.  There is a huge literature on Morita theory in many different contexts.  See, for example,
\begin{itemize}
\item \cite[Sections 4.4 and 4.5]{cohen_algebra} for Morita theory of modules over rings, 
\item \cite{Lindner,palmquist} for Morita theory of categories enriched in a closed category, and
\item \cite[Section 4]{schwede-shipley_stable} and \cref{ex:schwede-shipley} for Morita theory of modules over symmetric ring spectra.
\end{itemize} 
Our \cref{mackey-gen-xiv,mackey-xiv-cor,mackey-xiv-pmulticat,mackey-xiv-mone,mackey-pmulti-mone} are a kind of Morita theory that involves a change of closed multicategories.\index{closed!multicategory!Morita theory}
\end{question}

\part*{Back Matter}
\backmatter

\bibliographystyle{sty/amsalpha3}
\bibliography{references}

\providecommand{\bysame}{\leavevmode\hbox to3em{\hrulefill}\thinspace}
\providecommand{\MR}{\relax\ifhmode\unskip\space\fi MR }
\providecommand{\MRhref}[2]{%
  \href{http://www.ams.org/mathscinet-getitem?mr=#1}{#2}
}
\providecommand{\nopubyear}{$\infty$}
\providecommand{\doi}[1]{%
  doi:\href{https://dx.doi.org/#1}{\nolinkurl{#1}}}
\providecommand{\arxiv}[1]{%
  arXiv:\href{https://arxiv.org/abs/#1}{#1}}
\begin{thebibliography}{Yau$\infty$bAA}

\bibitem[AGV72]{agv}
M.~Artin, A.~Grothendieck, and J.~L. Verdier, \emph{{S\'eminaire de
  g\'eom\'etrie alg\'ebrique du Bois-Marie 1963--1964. Th\'eorie des topos et
  cohomologie \'etale des sch\'emas. (SGA 4). Un s\'eminaire dirig\'e par M.
  Artin, A. Grothendieck, J. L. Verdier. Avec la collaboration de N. Bourbaki,
  P. Deligne, B. Saint-Donat. Tome 1: Th\'eorie des topos. Expos\'es I \`a IV.
  2e \'ed.}}, {Lecture Notes in Mathematics}, vol. 269, Springer, Cham, 1972
  (French). \doi{10.1007/BFb0081551}

\bibitem[BR20]{barnes_roitzheim}
D.~Barnes and C.~Roitzheim, \emph{Foundations of stable homotopy theory},
  Cambridge Studies in Advanced Mathematics, vol. 185, Cambridge University
  Press, 2020.

\bibitem[BK12]{barwick-kan}
C.~Barwick and D.~M. Kan, \emph{A characterization of simplicial localization
  functors and a discussion of {DK} equivalences}, Indag. Math. (N.S.)
  \textbf{23} (2012), no.~1-2, 69--79. \doi{10.1016/j.indag.2011.10.001}

\bibitem[Bar17]{barwick}
C.~Barwick, \emph{Spectral {Mackey} functors and equivariant algebraic
  {{$K$}}-theory. {I}.}, Adv. Math. \textbf{304} (2017), 646--727.
  \doi{10.1016/j.aim.2016.08.043}

\bibitem[BGS20]{barwick-glasman-shah}
C.~Barwick, S.~Glasman, and J.~Shah, \emph{Spectral {Mackey} functors and
  equivariant algebraic {{$K$}}-theory. {II}.}, Tunis. J. Math. \textbf{2}
  (2020), no.~1, 97--146. \doi{10.2140/tunis.2020.2.97}

\bibitem[BKP89]{bkp}
R.~Blackwell, G.~M. Kelly, and A.~J. Power, \emph{Two-dimensional monad
  theory}, J. Pure Appl. Algebra \textbf{59} (1989), no.~1, 1--41.
  \doi{10.1016/0022-4049(89)90160-6}

\bibitem[BV73]{boardman-vogt}
J.~M. Boardman and R.~M. Vogt, \emph{Homotopy invariant algebraic structures on
  topological spaces}, Lecture Notes in Mathematics, Vol. 347, Springer-Verlag,
  Berlin-New York, 1973.

\bibitem[BO15]{bohmann_osorno-mackey}
A.~M. Bohmann and A.~Osorno, \emph{Constructing equivariant spectra via
  categorical {M}ackey functors}, Algebr. Geom. Topol. \textbf{15} (2015),
  no.~1, 537--563. \doi{10.2140/agt.2015.15.537}

\bibitem[BF78]{bousfield_friedlander}
A.~K. Bousfield and E.~M. Friedlander, \emph{Homotopy theory of
  {$\Gamma$}-spaces, spectra, and bisimplicial sets}, Geometric applications of
  homotopy theory ({P}roc. {C}onf., {E}vanston, {I}ll., 1977), {II}, Lecture
  Notes in Math., vol. 658, Springer, Berlin, 1978, pp.~80--130.

\bibitem[Coh03]{cohen_algebra}
P.~M. Cohn, \emph{Further algebra and applications}, Springer-Verlag London,
  Ltd., London, 2003. \doi{10.1007/978-1-4471-0039-3}

\bibitem[DK80]{dwyer-kan}
W.~G. Dwyer and D.~M. Kan, \emph{Calculating simplicial localizations}, J. Pure
  Appl. Algebra \textbf{18} (1980), no.~1, 17--35.
  \doi{10.1016/0022-4049(80)90113-9}

\bibitem[Elm83]{elmendorf-systems}
A.~D. Elmendorf, \emph{Systems of fixed point sets}, Trans. Amer. Math. Soc.
  \textbf{277} (1983), no.~1, 275--284. \doi{10.2307/1999356}

\bibitem[EM06]{elmendorf-mandell}
A.~D. Elmendorf and M.~A. Mandell, \emph{Rings, modules, and algebras in
  infinite loop space theory}, Adv. Math. \textbf{205} (2006), no.~1, 163--228.
  \doi{10.1016/j.aim.2005.07.007}

\bibitem[EM09]{elmendorf-mandell-perm}
\bysame, \emph{Permutative categories, multicategories and algebraic
  {$K$}-theory}, Algebr. Geom. Topol. \textbf{9} (2009), no.~4, 2391--2441.
  \doi{10.2140/agt.2009.9.2391}

\bibitem[FPP75]{palmquist}
J.~Fisher-Palmquist and P.~H. Palmquist, \emph{Morita contexts of enriched
  categories}, Proc. Amer. Math. Soc. \textbf{50} (1975), 55--60.
  \doi{10.2307/2040513}

\bibitem[GMR19]{GMRenriched}
B.~Guillou, J.~P. May, and J.~Rubin, \emph{Enriched model categories in
  equivariant contexts}, Homology, Homotopy and Applications \textbf{21}
  (2019), no.~1. \doi{doi:10.4310/HHA.2019.v21.n1.a10}

\bibitem[Gui10]{guillou}
B.~J. Guillou, \emph{Strictification of categories weakly enriched in symmetric
  monoidal categories}, Theory Appl. Categ. \textbf{24} (2010), No. 20,
  564--579.

\bibitem[GM22]{guillou_may}
B.~J. Guillou and J.~P. May, \emph{Models of {$G$}-spectra as presheaves of
  spectra}, Algebraic \& Geometric Topology (2022), to appear.
  \arxiv{1110.3571v5}

\bibitem[GMMO23]{gmmo}
B.~J. Guillou, J.~P. May, M.~Merling, and A.~M. Osorno, \emph{Multiplicative
  equivariant {$K$}-theory and the {B}arratt-{P}riddy-{Q}uillen theorem},
  Advances in Mathematics \textbf{414} (2023). \doi{10.1016/j.aim.2023.108865}

\bibitem[GJO17a]{gjo-extending}
N.~Gurski, N.~Johnson, and A.~M. Osorno, \emph{Extending homotopy theories
  across adjunctions}, Homology Homotopy Appl. \textbf{19} (2017), no.~2,
  89--110. \doi{10.4310/HHA.2017.v19.n2.a6}

\bibitem[GJO17b]{gjo1}
\bysame, \emph{{$K$}-theory for 2-categories}, Adv. Math. \textbf{322} (2017),
  378--472. \doi{10.1016/j.aim.2017.10.011}

\bibitem[HHR16]{hhr-kervaire}
M.~A. Hill, M.~J. Hopkins, and D.~C. Ravenel, \emph{On the nonexistence of
  elements of {K}ervaire invariant one}, Ann. of Math. (2) \textbf{184} (2016),
  no.~1, 1--262. \doi{10.4007/annals.2016.184.1.1}

\bibitem[HHR21]{hhrbook}
M.~A. Hill, M.~J. Hopkins, and D.~C. Ravenel, \emph{Equivariant stable homotopy
  theory and the {K}ervaire invariant problem}, New Mathematical Monographs,
  vol.~40, Cambridge University Press, Cambridge, 2021.
  \doi{10.1017/9781108917278}

\bibitem[Hir03]{hirschhorn}
P.~S. Hirschhorn, \emph{Model categories and their localizations}, Mathematical
  Surveys and Monographs, vol.~99, American Mathematical Society, Providence,
  RI, 2003.

\bibitem[Hov99]{hovey}
M.~Hovey, \emph{Model categories}, Mathematical Surveys and Monographs,
  vol.~63, American Mathematical Society, Providence, RI, 1999.

\bibitem[HSS00]{hss}
M.~Hovey, B.~Shipley, and J.~Smith, \emph{Symmetric spectra}, J. Amer. Math.
  Soc. \textbf{13} (2000), no.~1, 149--208. \doi{10.1090/S0894-0347-99-00320-3}

\bibitem[JY$\infty$]{cerberusIII}
N.~Johnson and D.~Yau, \emph{{B}imonoidal {C}ategories, {$E_n$}-{M}onoidal
  {C}ategories, and {A}lgebraic {$K$}-{T}heory. {V}olume {III}: {F}rom
  {C}ategories to {S}tructured {R}ing {S}pectra}, available at
  \url{https://nilesjohnson.net/}.

\bibitem[JY21]{johnson-yau}
\bysame, \emph{{2}-{D}imensional {C}ategories}, Oxford University Press, New
  York, 2021. \doi{10.1093/oso/9780198871378.001.0001}

\bibitem[JY22a]{johnson-yau-Fmulti}
\bysame, \emph{Homotopy equivalent algebraic structures in multicategories and
  permutative categories}, Theory Appl. Categ. \textbf{38} (2022), no.~30,
  1156--1208.

\bibitem[JY22b]{johnson-yau-invK}
\bysame, \emph{{M}ultifunctorial inverse {$K$}-theory}, Annals of $K$-Theory
  \textbf{7} (2022), no.~3, 507--548. \doi{10.2140/akt.2022.7.507}

\bibitem[JY22c]{johnson-yau-multiK}
\bysame, \emph{Multifunctorial {$K$}-theory is an equivalence of homotopy
  theories}, J. Homotopy Relat. Struct. \textbf{17} (2022), no.~4, 569--592.
  \doi{10.1007/s40062-022-00317-8}

\bibitem[JY23]{johnson-yau-permmult}
\bysame, \emph{Multicategories model all connective spectra}, Homology,
  Homotopy and Applications \textbf{25} (2023), no.~1, 147--172.
  \doi{10.4310/HHA.2023.v25.n1.a8}

\bibitem[JS93]{joyal-street}
A.~Joyal and R.~Street, \emph{Braided tensor categories}, Adv. Math.
  \textbf{102} (1993), no.~1, 20--78. \doi{10.1006/aima.1993.1055}

\bibitem[Kel05]{kelly-enriched}
G.~M. Kelly, \emph{Basic concepts of enriched category theory}, Repr. Theory
  Appl. Categ. (2005), no.~10, vi+137, Reprint of the 1982 original (Cambridge
  Univ. Press, Cambridge).

\bibitem[Lam69]{lambek}
J.~Lambek, \emph{Deductive systems and categories {II}. {S}tandard
  constructions and closed categories}, Category {T}heory, {H}omology {T}heory
  and their {A}pplications, {I} ({B}attelle {I}nstitute {C}onference,
  {S}eattle, {W}ash., 1968, {V}ol. {O}ne), Springer, Berlin, 1969, pp.~76--122.

\bibitem[LMS86]{lewis-may-steinberger}
L.~G. Lewis, Jr., J.~P. May, and M.~Steinberger, \emph{Equivariant stable
  homotopy theory}, Lecture Notes in Mathematics, vol. 1213, Springer-Verlag,
  Berlin, 1986, With contributions by J. E. McClure. \doi{10.1007/BFb0075778}

\bibitem[Lin74]{Lindner}
H.~Lindner, \emph{Morita equivalences of enriched categories}, Cahiers
  Topologie G\'{e}om. Diff\'{e}rentielle \textbf{15} (1974), no.~4, 377--397,
  449--450.

\bibitem[ML69]{maclane-foundation}
S.~Mac~Lane, \emph{One universe as a foundation for category theory}, Reports
  of the {M}idwest {C}ategory {S}eminar. {III}, Lecture Notes in Mathematics,
  Vol. 106, Springer, Berlin, 1969, pp.~192--200.

\bibitem[ML98]{maclane}
\bysame, \emph{Categories for the working mathematician}, second ed., Graduate
  Texts in Mathematics, vol.~5, Springer-Verlag, New York, 1998.

\bibitem[MM19]{merling-malkiewich-Athy}
C.~Malkiewich and M.~Merling, \emph{Equivariant {$A$}-theory}, Doc. Math.
  \textbf{24} (2019), 815--855. \doi{doi:10.25537/dm.2019v24.815-855}

\bibitem[MM22]{merling-malkiewich}
\bysame, \emph{The equivariant parametrized {{\(h\)}}-cobordism theorem, the
  non-manifold part}, Adv. Math. \textbf{399} (2022), 42 (English), Id/No
  108242. \doi{10.1016/j.aim.2022.108242}

\bibitem[{Man}10]{mandell_inverseK}
M.~A. {Mandell}, \emph{{An inverse {$K$}-theory functor}}, {Doc. Math.}
  \textbf{15} (2010), 765--791.

\bibitem[{Man}12]{manzyuk}
O.~{Manzyuk}, \emph{{Closed categories vs. closed multicategories}}, {Theory
  Appl. Categ.} \textbf{26} (2012), 132--175.

\bibitem[May78]{may-permutative}
J.~P. May, \emph{The spectra associated to permutative categories}, Topology
  \textbf{17} (1978), no.~3, 225--228. \doi{10.1016/0040-9383(78)90027-7}

\bibitem[May96]{alaska-notes}
\bysame, \emph{Equivariant homotopy and cohomology theory}, CBMS Regional
  Conference Series in Mathematics, vol.~91, Published for the Conference Board
  of the Mathematical Sciences, Washington, DC; by the American Mathematical
  Society, Providence, RI, 1996, With contributions by M. Cole, G.
  Comeza\~{n}a, S. Costenoble, A. D. Elmendorf, J. P. C. Greenlees, L. G.
  Lewis, Jr., R. J. Piacenza, G. Triantafillou, and S. Waner.
  \doi{10.1090/cbms/091}

\bibitem[May99]{mayconcise}
\bysame, \emph{A concise course in algebraic topology}, Chicago Lectures in
  Mathematics, University of Chicago press, Chicago, IL, 1999.

\bibitem[Qui73]{quillenKI}
D.~G. Quillen, \emph{Higher algebraic {$K$}-theory {I}: {H}igher
  {$K$}-theories}, {P}roc. {C}onf., {B}attelle {M}emorial {I}nst., {S}eattle,
  {W}ash., (1972), Lecture Notes in Math., Vol. 341, 1973, pp.~85--147.

\bibitem[Rez01]{rezk-homotopy-theory}
C.~Rezk, \emph{A model for the homotopy theory of homotopy theory}, Trans.
  Amer. Math. Soc. \textbf{353} (2001), no.~3, 973--1007 (electronic).
  \doi{10.1090/S0002-9947-00-02653-2}

\bibitem[SS03]{schwede-shipley_stable}
S.~Schwede and B.~Shipley, \emph{Stable model categories are categories of
  modules}, Topology \textbf{42} (2003), no.~1, 103--153.
  \doi{10.1016/S0040-9383(02)00006-X}

\bibitem[Seg74]{segal}
G.~Segal, \emph{Categories and cohomology theories}, Topology \textbf{13}
  (1974), 293--312. \doi{10.1016/0040-9383(74)90022-6}

\bibitem[Tho95]{thomason}
R.~W. Thomason, \emph{Symmetric monoidal categories model all connective
  spectra}, Theory Appl. Categ. \textbf{1} (1995), No. 5, 78--118.

\bibitem[To{\"e}05]{toen-axiomatisation}
B.~To{\"e}n, \emph{Vers une axiomatisation de la th\'eorie des cat\'egories
  sup\'erieures}, $K$-Theory \textbf{34} (2005), no.~3, 233--263.
  \doi{10.1007/s10977-005-4556-6}

\bibitem[tD79]{tomdieck}
T.~tom Dieck, \emph{Transformation groups and representation theory}, Lecture
  Notes in Mathematics, vol. 766, Springer, Berlin, 1979.
  \doi{10.1007/BFb0085965}

\bibitem[Wal85]{waldhausen}
F.~Waldhausen, \emph{Algebraic {$K$}-theory of spaces}, Algebraic and geometric
  topology ({N}ew {B}runswick, {N}.{J}., 1983), Lecture Notes in Math., vol.
  1126, Springer, Berlin, 1985, pp.~318--419. \doi{10.1007/BFb0074449}

\bibitem[Web00]{webb-mackey}
P.~Webb, \emph{A guide to {M}ackey functors}, Handbook of algebra, {V}ol. 2,
  Handb. Algebr., vol.~2, Elsevier/North-Holland, Amsterdam, 2000,
  pp.~805--836. \doi{10.1016/S1570-7954(00)80044-3}

\bibitem[Yau$\infty$a]{cerberusI}
D.~Yau, \emph{{B}imonoidal {C}ategories, {$E_n$}-{M}onoidal {C}ategories, and
  {A}lgebraic {$K$}-{T}heory. {V}olume {I}: {S}ymmetric {B}imonoidal
  {C}ategories and {M}onoidal {B}icategories}, available at
  \url{https://u.osu.edu/yau.22/main/}.

\bibitem[Yau$\infty$b]{cerberusII}
\bysame, \emph{{B}imonoidal {C}ategories, {$E_n$}-{M}onoidal {C}ategories, and
  {A}lgebraic {$K$}-{T}heory. {V}olume {II}: {B}raided {B}imonoidal
  {C}ategories with {A}pplications}, available at
  \url{https://u.osu.edu/yau.22/main/}.

\bibitem[Yau16]{yau-operad}
\bysame, \emph{Colored operads}, Graduate Studies in Mathematics, vol. 170,
  American Mathematical Society, Providence, RI, 2016.

\bibitem[Yau20a]{yau-hqft}
\bysame, \emph{Homotopical quantum field theory}, Hackensack, NJ: World
  Scientific, 2020. \doi{10.1142/11626}

\bibitem[Yau20b]{yau-involutive}
\bysame, \emph{Involutive category theory}, Lecture Notes in Mathematics, vol.
  2279, Springer, Cham, 2020. \doi{10.1007/978-3-030-61203-0}

\bibitem[Yau22]{yau-inf-operad}
\bysame, \emph{Infinity operads and monoidal categories with group
  equivariance}, Singapore: World Scientific, 2022. \doi{10.1142/12687}

\bibitem[Yau24]{yau-multigro}
\bysame, \emph{{G}rothendieck {C}onstruction of {B}ipermutative-{I}ndexed
  {C}ategories}, CRC Press, 2024.

\bibitem[YJ15]{bluemonster}
D.~Yau and M.~W. Johnson, \emph{A foundation for {PROP}s, algebras, and
  modules}, Mathematical Surveys and Monographs, vol. 203, American
  Mathematical Society, Providence, RI, 2015. \doi{10.1090/surv/203}

\bibitem[Zak18]{zakharevich}
I.~Zakharevich, \emph{The category of {W}aldhausen categories is a closed
  multicategory}, New directions in homotopy theory, Contemp. Math., vol. 707,
  Amer. Math. Soc., Providence, RI, 2018, pp.~175--194.
  \doi{10.1090/conm/707/14259}

\end{thebibliography}

\newcommand{\fact}[2]{\noindent({#2}) {#1}}

\newcommand{\chapNumName}[1]{\medskip\begin{center}\textbf{\Cref{#1}}.  \nameref{#1}\end{center}}
\newcommand{\thm}[1]{\textbf{#1}.}



\chapter*{List of Main Facts}

\chapNumName{ch:motivations}
\fact{\thm{Elmendorf's Theorem} There is a Quillen equivalence between $G$-spaces and topological presheaves on $\OrbG$.}{\ref{theorem:elmendorf}}

\fact{The spans in a small category $\C$ with chosen pullbacks form 1-cells of a bicategory $\Span(\C)$.}{\ref{definition:Span}}

\fact{The Burnside category $\GB$ is self-dual via the functor that reverses spans.}{\ref{lemma:GB-self-dual}}

\fact{$\Span(\FinG)$ has a choice of pullbacks that makes the composition of 1-cells strictly associative and strictly unital on one side.}{\ref{explanation:SpanFinG}}

\fact{The self-duality of the Burnside category does \emph{not} extend to a 2-functor on the Burnside 2-category.}{\ref{remark:GE-non-self-dual}}

\fact{\thm{Guillou-May Theorem} There is a Quillen equivalence between $G$-spectra and spectral Mackey functors.}{\ref{theorem:GM}}

\fact{\thm{Schwede-Shipley Characterization Theorem} Stable model categories with certain additional hypotheses are characterized by spectral Mackey functors on spectral endomorphism categories.}{\ref{theorem:schwede-shipley}}

\chapNumName{part:background}

\medskip
\chapNumName{ch:ptmulticat}

\fact{The category of small multicategories is strictly monadic over the category of multigraphs.}{\ref{thm:multigraphcat-monadic}}

\fact{$\Multicat$ is a symmetric monoidal $\Cat$-category with the Boardman-Vogt tensor product.}{\ref{theorem:multicat-symmon}}

\fact{$\Multicat$ is a $\Cat$-multicategory.}{\ref{expl:multicatcatmulticat}}

\fact{The symmetric monoidal category $\Multicat$ is closed.}{\ref{theorem:multicat-sm-closed}}

\fact{$\pMulticat$ is a complete and cocomplete symmetric monoidal closed category.}{\ref{thm:pmulticat-smclosed}}

\fact{$\pMulticat$ is a $\Cat$-multicategory.}{\ref{expl:ptmulticatcatmulticat}}

\fact{The forgetful 2-functor $\Ust \cn \pMulticat \to \Multicat$ is a symmetric monoidal $\Cat$-functor.}{\ref{expl:Ust}}

\fact{$\Ust \cn \pMulticat \to \Multicat$ is a $\Cat$-multifunctor.}{\ref{expl:Ust-catmulti}}

\fact{The partition products $\binprod_{\ord{1},-}$ and $\binprod_{-,\ord{1}}$ are isomorphisms.}{\ref{lemma:part-prod-1}}

\fact{The partition multicategory $\cM \cn \Fskel^\op \to \pMulticat$ is a symmetric monoidal functor.}{\ref{proposition:cM-symm-mon}}

\fact{$\Mone$ is a commutative monoid in $\pMulticat$.}{\ref{def:Monecommonoid}}

\fact{There is a 2-category $\MoneMod$ of left $\Mone$-modules.}{\ref{definition:MoneMod-prelim}}

\fact{Each symmetric monoidal category has an endomorphism left $\Mone$-module.}{\ref{ex:endmc}}

\fact{$\Endm \cn \permcatsu \to \MoneMod$ is a 2-functor.}{\ref{endmtwofunctor}}

\fact{Each small pointed multicategory has at most one left $\Mone$-module structure, and the structure morphism is an isomorphism.  $\MoneMod$ is a full sub-2-category of $\pMulticat$.  $\MoneMod$ is a complete and cocomplete symmetric monoidal closed category.}{\ref{proposition:EM2-5-1}}

\fact{There is a free-forgetful adjunction $\Monesma \cn \pMulticat \lradj \MoneMod \cn \Um$.}{\ref{MonesmaUmadj}}

\fact{The counit of $\big(\Monesma, \Um\big)$ is componentwise an isomorphism.}{\ref{epzm-def}}

\fact{$\MoneMod$ is a symmetric monoidal $\Cat$-category.}{\ref{definition:MoneMod}}

\fact{$\MoneMod$ is a $\Cat$-multicategory.}{\ref{expl:monemodcatmulticat}}

\fact{$\Monesma$ is a strong symmetric $\Cat$-monoidal functor, hence also a $\Cat$-multifunctor.}{\ref{Monesma-CatSM}}

\fact{$\Um$ is a symmetric monoidal $\Cat$-functor.}{\ref{expl:moneinclusion}}

\fact{$\Um$ is a $\Cat$-multifunctor.}{\ref{expl:Um-catmulti}}

\fact{The unit and counit of $\big(\Monesma, \Um\big)$ are monoidal $\Cat$-natural transformations, hence also $\Cat$-multinatural transformations.}{\ref{etahat-epzhat-monCatnat}}

\fact{A 1-linear functor is precisely a strictly unital symmetric monoidal functor.}{\ref{ex:onelinearfunctor}}

\fact{A 1-linear transformation is precisely a monoidal natural transformation.}{\ref{ex:onelineartr}}

\fact{There are $\Cat$-multicategories $\permcatsu$ and $\permcatst$.}{\ref{thm:permcatmulticat}}

\fact{$\Endst$ induces an isomorphism between multimorphism categories.}{\ref{proposition:n-lin-equiv}}

\fact{$\Endst$ is a $\Cat$-multifunctor.}{\ref{expl:endst-catmulti}}

\fact{$\End = \Ust \Endst$ and $\Endst = \Um \Endm$.}{\ref{endfactorization}}

\fact{$\End$ is a $\Cat$-multifunctor.}{\ref{expl:end-catmulti}}

\fact{$\Endm$ is a $\Cat$-multifunctor.}{\ref{expl:endm-catmulti}}

\medskip
\chapNumName{ch:Kspectra}

\fact{The complete Segal space model structure on bisimplicial sets is a simplicial model structure whose fibrant objects are precisely the complete Segal spaces.}{\ref{theorem:css-fibrant}}

\fact{Inverse equivalences of homotopy theories are equivalences of homotopy theories.}{\ref{gjo29}}

\fact{An adjoint equivalence of homotopy theories induces equivalences of homotopy theories.}{\ref{def:heq}}

\fact{For a complete and cocomplete symmetric monoidal closed category $\C$, the category $\pC$ of pointed objects is also a complete and cocomplete symmetric monoidal closed category.}{\ref{theorem:pC-sm-closed}}

\fact{$\DstarV$ is a complete and cocomplete symmetric monoidal closed category with the pointed Day convolution.}{\ref{thm:Dgm-pv-convolution-hom}}

\fact{$\DstarV$ is enriched and (co)tensored over $\pV$.  $\DstarV$ is a $\V$-multicategory.}{\ref{theorem:diagram-omnibus}}

\fact{$\GaV$ is a complete and cocomplete symmetric monoidal closed category.}{\ref{GammaV}}

\fact{$\GstarV$ is a complete and cocomplete symmetric monoidal closed category.}{\ref{Gstar-V}}

\fact{Length-one inclusion defines a pointed functor $i \cn \Fskel \to \Gskel$.}{\ref{i-incl}}

\fact{Smash product $\sma \cn \Gskel \to \Fskel$ is a strict symmetric monoidal pointed functor.}{\ref{smagf}}

\fact{Each functor in \cref{eq:Ksummary}, except $\Jt$ and $\Kg$, is an equivalence of homotopy theories.}{\ref{eq:Ksummary}}

\fact{Segal $K$-theory is the composite functor $\Kse = \Kf \Ner_* \Jse$.}{\ref{Kse}}

\fact{Segal $J$-theory is not a multifunctor, so neither is $\Kse$.}{\ref{Jse}}

\fact{Elmendorf-Mandell $K$-theory is the multifunctor $\Kem = \Kg \Ner_* \Jt \Endm$.}{\ref{Kem}}

\fact{Segal $K$-theory is an equivalence of homotopy theories.}{\ref{Kse-heq}}

\fact{Elmendorf-Mandell $K$-theory is an equivalence of homotopy theories.}{\ref{Kem-heq}}

\medskip
\chapNumName{ch:multperm}

\fact{Each multicategory $\M$ has an associated free permutative category $\Fr\M$.}{\ref{proposition:free-perm}}

\fact{$\Fr\Mtu$ is isomorphic to the permutation category.}{\ref{example:free-Mtu}}

\fact{$\Fr\Mterm$ is isomorphic to the category of natural numbers and morphisms of finite sets.}{\ref{example:free-Mterm}}

\fact{$\Fr \cn \Multicat \to \permcatst$ is a 2-functor.}{\ref{proposition:free-perm-functor}}

\fact{$\Fr$ is a left 2-adjoint of $\End$.}{\ref{theorem:FE-adj}}

\fact{The counit $\epz$ of $(\Fr,\End)$ admits a componentwise right adjoint $\vrho$.}{\ref{proposition:epz-rho-adj}}

\fact{$\vrho_\C \cn \C \to \Fr\End(\C)$ is a symmetric monoidal functor.}{\ref{rhoc-monoidal}}

\fact{$\Fr^n$ is a strong $n$-linear functor that is 2-natural with respect to multifunctors and multinatural transformations.}{\ref{S-multifunctorial}}

\fact{$\Fr \cn \Multicat \to \permcatsu$ is a non-symmetric $\Cat$-multifunctor.}{\ref{theorem:F-multi}}

\fact{The unit $\eta \cn 1 \to \End\Fr$ is a non-symmetric $\Cat$-multinatural transformation.}{\ref{lemma:eta-mnat}}

\fact{$\Fr \cn \Multicat \lradj \permcatst \cn \End$ is an adjoint equivalence of homotopy theories.}{\ref{thm:F-heq}}

\fact{For each small non-symmetric $\Cat$-multicategory $\Q$, $(\Fr^\Q, \End^\Q)$ are inverse equivalences of homotopy theories between $\Multicat^\Q$ and $(\permcatsu)^\Q$.}{\ref{thm:alg-hty-equiv}}

\fact{$\Fr \cn \Multicat \lradj \permcatsu \cn \End$ are inverse equivalences of homotopy theories.}{\ref{thm:Fsu-heq}}

\fact{Inclusion $I \cn \permcatst \to \permcatsu$ is an equivalence of homotopy theories.}{\ref{cor:I-heq}}

\medskip
\chapNumName{part:multicat}

\medskip
\chapNumName{ch:ptmulticat-sp}

\fact{Each pointed multicategory $\M$ has an associated permutative category $\Fst\M$.}{\ref{def:Fst-permutative}}

\fact{$\Fst \cn \pMulticat \to \permcatst$ is a 2-functor.}{\ref{ptmulticat-thm-i}}

\fact{$\pst \cn \Fr \to \Fst$ is a 2-natural transformation with each component a strict symmetric monoidal functor.}{\ref{pst-iinatural}}

\fact{$\pst_\M$ is a 2-pushout of $\Fr\Mterm \to \boldone$ in $\permcatst$.}{\ref{proposition:FstM-pushout}}

\fact{For each small pointed multicategory $\M$, $\etast_\M \cn \M \to \Endst\Fst\M$ is a pointed multifunctor that is 2-natural in $\M$.}{\ref{etast-pointed}}

\fact{For each small permutative category $\C$, $\epzst_\C \cn \Fst\Endst\C \to \C$ is a strict symmetric monoidal functor.}{\ref{def:epzst}}

\fact{$\epzst_\C$ is 2-natural in $\C$.}{\ref{epzst-strict}}

\fact{There is a 2-adjunction $\Fst \cn \pMulticat \lradj \permcatst \cn \Endst$.}{\ref{ptmulticat-thm-v}}

\fact{There is a 2-natural isomorphism $\Fr \iso \Fst \circ \dplus$.}{\ref{FFstdplus}}

\fact{There is a 2-adjunction $\Fm \cn \MoneMod \lradj \permcatst \cn \Endm$ with $\Fm = \Fst\Um$.}{\ref{ptmulticat-thm-ii}}

\fact{There is a 2-natural isomorphism $\Fst \iso \Fm \circ (\Monesma)$.}{\ref{FstFmMonesma}}

\fact{$\Fst\Mterm \iso \boldone$.}{\ref{example:Fst-Mterm}}

\fact{A left $\Mone$-module structure on a small pointed multicategory $\M$ determines and is uniquely determined by binary operations $\pi^2_1(x)$ for objects $x \in \M$ that satisfy basepoint, unit, and interchange conditions.}{\ref{M1-mod-str}}

\fact{For each left $\Mone$-module $\M$, each morphism in $\Fm\M$ is represented by a length-one sequence.}{\ref{corollary:FmM-reps}}

\fact{For each small permutative category $\C$, $\vrhost_\C \cn \C \to \Fst\Endst\C$ is a strictly unital symmetric monoidal functor.}{\ref{ptmulticat-xxi}}

\fact{The adjunction $\epz_\C \dashv \vrho_\C$ extends to an adjunction $\epzst_\C \dashv \vrhost_\C$ in $\permcatsu$.}{\ref{ptmulticat-prop-viii}}

\fact{$\etast_{\Endst\C} = \Endst\vrhost_\C$.}{\ref{etaEEvrho}}

\fact{$\Fst \cn \pMulticat \lradj \permcatst \cn \Endst$ is an adjoint equivalence of homotopy theories.}{\ref{ptmulticat-thm-x}}

\fact{$\dplus \cn \Multicat \to \pMulticat$ is an equivalence of homotopy theories.}{\ref{ptmulticat-cor-page-vi}}

\fact{$\Monesma \cn \pMulticat \lradj \MoneMod \cn \Um$ is an adjoint equivalence of homotopy theories.}{\ref{ptmulticat-thm-vi}}

\fact{$\Fm \cn \MoneMod \lradj \permcatst \cn \Endm$ is an adjoint equivalence of homotopy theories.}{\ref{ptmulticat-thm-xi}}

\medskip
\chapNumName{ch:ptmulticat-alg}

\fact{$\Fst^n$ is a strong $n$-linear functor.}{\ref{ptmulticat-xv}}

\fact{$\Fst^n$ is 2-natural with respect to pointed multifunctors and pointed multinatural transformations.}{\ref{ptmulticat-xvi}}

\fact{$\Fst \cn \pMulticat \to \permcatsu$ is a non-symmetric $\Cat$-multifunctor.}{\ref{ptmulticat-xvii}}

\fact{$\pst \cn \Fr\Ust \to \Fst$ is a non-symmetric $\Cat$-multinatural transformation.}{\ref{pst-mnat}}

\fact{$\etast \cn 1 \to \Endst\Fst$ is a non-symmetric $\Cat$-multinatural transformation.}{\ref{ptmulticat-xx}}

\fact{$\vrhost \cn 1 \to \Fst\Endst$ is a non-symmetric $\Cat$-multinatural transformation.}{\ref{ptmulticat-xxii}} 

\fact{For each small non-symmetric $\Cat$-multicategory $\Q$, $(\Fst^\Q,\Endst^\Q)$ are inverse equivalences of homotopy theories between $\pMulticat^\Q$ and $(\permcatsu)^\Q$.}{\ref{ptmulticat-xxiii}}

\fact{$\Fm \cn \MoneMod \to \permcatsu$ is a non-symmetric $\Cat$-multifunctor.}{\ref{Fm-multi-def}}

\fact{$\etam \cn 1 \to \Endm\Fm$ is a non-symmetric $\Cat$-multinatural transformation.}{\ref{def:etam-multi}}

\fact{$\vrhom \cn 1 \to \Fm\Endm$ is a non-symmetric $\Cat$-multinatural transformation.}{\ref{def:vrhom-multi}}

\fact{$\etam_{\Endm\C} = \Endm\vrhom_{\C}$.}{\ref{etavrhoMone}}

\fact{For each small non-symmetric $\Cat$-multicategory $\Q$, $(\Fm^\Q,\Endm^\Q)$ are inverse equivalences of homotopy theories between $(\MoneMod)^\Q$ and $(\permcatsu)^\Q$.}{\ref{ptmulticat-xxv}}

\fact{For each small (non-)symmetric $\Cat$-multicategory $\Q$, $\big((\Monesma)^\Q,\Um^\Q\big)$ are inverse equivalences of homotopy theories between $\pMulticat^\Q$ and $(\MoneMod)^\Q$.}{\ref{Monesma-Um-algebra}}

\medskip
\chapNumName{part:enrpresheave}

\medskip
\chapNumName{ch:menriched}

\fact{For each non-symmetric multicategory $\M$, there is a 2-category $\MCat$ of small $\M$-categories, $\M$-functors, and $\M$-natural transformations.}{\ref{mcat-iicat}}

\fact{For a monoidal category $\V$, $\VCat$ and $\EndVCat$ are the same 2-categories.}{\ref{EndV-enriched}}

\fact{For small permutative categories $\C$ and $\D$, $\psu(\C,\D)$ is a small permutative category.}{\ref{psucd-hom-permcat}}

\fact{Composition $\mcomp_{\B,\C,\D}$ is a bilinear functor.}{\ref{psu-mBCD}}

\fact{$\psu$ is a $\psu$-category.}{\ref{permcat-selfenr}}

\fact{For small permutative categories $\C$ and $\D$, $\ev_{\C,\D}$ is a bilinear functor.}{\ref{ev-bilinear}}

\fact{$\mcomp_{\B,\C,\D}$ is compatible with evaluation.}{\ref{ev-comp}}

\fact{For a multicategory $\M$ and an $\M$-category $\C$, $\Cop$ is an $\M$-category.}{\ref{opposite-mcat}}

\fact{For a symmetric monoidal category $\V$, an opposite $\V$-category is the same as an opposite $(\EndV)$-category.}{\ref{v-opposite-mcat}}

\medskip
\chapNumName{ch:change_enr}

\fact{Each non-symmetric multifunctor $F$ induces a change-of-enrichment 2-functor.}{\ref{mult-change-enrichment}}

\fact{The change-of-enrichment 2-functor of a multifunctor preserves opposite enriched categories.}{\ref{dF-opposite}}

\fact{For a monoidal functor $U$, the change-of-enrichment 2-functors along $U$ and $\EndU$ are the same.}{\ref{mon-change-enrichment}}

\fact{Change-of-enrichment 2-functors of non-symmetric multifunctors are closed under composition.}{\ref{func-change-enr}}

\fact{There is a 2-functor $\Enr \cn \Multicatns \to \iicat$ that sends a small non-symmetric multicategory $\M$ to $\MCat$.}{\ref{change-enr-twofunctor}}

\medskip
\chapNumName{ch:gspectra}

\fact{A closed multicategory is a multicategory equipped with $n$-ary internal hom objects, symmetric group action on internal hom objects, and multicategorical evaluation that satisfy equivariance and evaluation bijection axioms.}{\ref{def:closed-multicat}}

\fact{For each symmetric monoidal closed category, the endomorphism multicategory is closed.}{\ref{smclosed-closed-multicat}}

\fact{Each $\clpsu\scmap{\angC;\D}$ is a permutative category.}{\ref{clp-angcd-permutative}}

\fact{The permutative categories $\clpsu\scmap{\angC;\D}$ admit symmetric group action that satisfies the equivariance axioms for internal hom objects.}{\ref{clp-angcd-strict}}

\fact{Each $\ev_{\angC;\, \D}$ is a multilinear functor.}{\ref{clp-evaluation}}

\fact{$\psu$ satisfies the evaluation bijection axiom.}{\ref{clp-evbij-axiom}}

\fact{For $\psu$, the inverse of $\chi$ is $\Psi$.}{\ref{chiPsi}}

\fact{$\psu$ satisfies the equivariance axioms for evaluation bijection.}{\ref{clp-evbij-equiv}}

\fact{$\psu$ is a closed multicategory.}{\ref{perm-closed-multicat}}

\medskip
\chapNumName{ch:std_enrich}

\fact{Each non-symmetric closed multicategory admits a canonical self-enrichment.}{\ref{cl-multi-cl-cat}}

\fact{For $\psu$, the self-enrichment coincides with the canonical self-enrichment.}{\ref{ex:perm-closed-multicat}}

\fact{For a symmetric monoidal closed category $\V$, the canonical self-enrichment of $\V$ coincides with the canonical self-enrichment of $\EndV$.}{\ref{ex:cl-multi-cl-cat}}

\fact{Each non-symmetric multifunctor $F$ admits a standard enrichment $\Fse$.}{\ref{gspectra-thm-iii}}

\fact{For a monoidal functor $U$ between symmetric monoidal closed categories, the standard enrichment of $U$ coincides with the standard enrichment of $\EndU$.}{\ref{std-enr-monoidal}}

\fact{Standard enrichment functors are closed under composition in an appropriate sense.}{\ref{gspectra-thm-iv}}

\fact{The standard enrichment of $\Kem$ factors into four spectral functors.}{\ref{gspectra-thm-xi}}

\medskip
\chapNumName{ch:gspectra_Kem}

\fact{For a non-symmetric closed multicategory $\M$ and an $\M$-category $\C$, a $\C$-diagram in $\M$ is precisely a left $\C$-module.}{\ref{C-diagram-partner}}

\fact{For a closed multicategory $\M$ and an $\M$-category $\C$, a $\C$-Mackey functor in $\M$ is precisely a left $\Cop$-module.}{\ref{C-presheaf-partner}}

\fact{An $\M$-natural transformation between $\C$-diagrams in $\M$ is precisely a left $\C$-module morphism.}{\ref{C-diag-morphism-pn}}

\fact{For vertically composable $\M$-natural transformations between $\C$-diagrams in $\M$, composition commutes with taking partners componentwise.}{\ref{mcat-vertical-comp}}

\fact{Each simplicial, cofibrantly generated, proper, and stable model category is Quillen equivalent to a category of spectral Mackey functors.}{\ref{ex:schwede-shipley}}

\fact{For each finite group $G$, the category of genuine $G$-equivariant spectra is Quillen equivalent to a spectral Mackey functor category associated to the permutative Burnside category.}{\ref{ex:guillou-may}}

\fact{For each non-symmetric multifunctor $F \cn \M \to \N$ between non-symmetric closed multicategories and a small $\M$-category $\C$, there is an induced diagram change-of-enrichment functor from $\MCat(\C,\M)$ to $\NCat(\C_F,\N)$.}{\ref{gspectra-thm-v}}

\fact{For each multifunctor $F \cn \M \to \N$ between closed multicategories and a small $\M$-category $\C$, there is an induced presheaf change-of-enrichment functor from $\MCat(\Cop,\M)$ to $\NCat((\C_F)^\op,\N)$.}{\ref{gspectra-thm-v-cor}}

\fact{Diagram change-of-enrichment functors are closed under composition.}{\ref{gspectra-thm-vii}}

\fact{Presheaf change-of-enrichment functors are closed under composition.}{\ref{gspectra-thm-vii-cor}}

\fact{$\Kem$ induces diagram and presheaf change-of-enrichment functors.}{\ref{thm:Kemdg}}

\fact{$\Kemdg$ factors into four change-of-enrichment functors.}{\ref{gspectra-thm-xiv}}

\medskip
\chapNumName{part:homotopy-mackey}

\medskip
\chapNumName{ch:mackey}

\fact{For non-symmetric multifunctors $F \cn \M \lradj \N \cn E$ between non-symmetric closed multicategories, a small $\N$-category $\C$, and a multinatural transformation $\cou \cn 1_\N \to FE$, there is an induced functor $\Fdgr$ from $\MCat(\C_E,\M)$ to $\NCat(\C,\N)$.}{\ref{def:mackey-gen-context}}

\fact{A multinatural transformation $\uni \cn 1_\M \to EF$ induces a natural transformation $\unidg \cn 1 \to \Edg\Fdgr$ on $\MCat(\C_E,\M)$.}{\ref{def:unidg}}

\fact{A multinatural transformation $\cou \cn 1_\N \to FE$ induces a natural transformation $\coudg \cn 1 \to \Fdgr\Edg$ on $\NCat(\C,\N)$.}{\ref{def:coudg}}

\fact{For a non-symmetric closed multicategory $\P$ equipped with a relative category structure $\cW$ and a $\P$-category $\D$, there is an induced relative category structure on $\PCat(\D,\P)$.}{\ref{def:enr-diagcat-relative}}

\fact{Under appropriate assumptions, inverse equivalences of homotopy theories $(F,E)$ lift to inverse equivalences of homotopy theories $(\Fdgr,\Edg)$ between $\MCat(\C_E,\M)$ and $\NCat(\C,\N)$.}{\ref{mackey-gen-xiv}}

\fact{If, furthermore, $E$ is a multifunctor, then $(\Fdgr,\Edg)$ are inverse equivalences of homotopy theories between $\MCat((\C_E)^\op,\M)$ and $\NCat(\Cop,\N)$.}{\ref{mackey-xiv-cor}}

\medskip
\chapNumName{ch:mackey_eq}

\fact{$\Fst$ and $\Endst$ induce inverse equivalences of homotopy theories between (i) $\C_{\Endst}$-diagrams in $\pMulticat$ and $\C$-diagrams in $\permcatsu$ and (ii) $\C_{\Endst}$-Mackey functors in $\pMulticat$ and $\C$-Mackey functors in $\permcatsu$.}{\ref{mackey-xiv-pmulticat}}

\fact{$\Fm$ and $\Endm$ induce inverse equivalences of homotopy theories between (i) $\C_{\Endm}$-diagrams in $\MoneMod$ and $\C$-diagrams in $\permcatsu$ and (ii) $\C_{\Endm}$-Mackey functors in $\MoneMod$ and $\C$-Mackey functors in $\permcatsu$.}{\ref{mackey-xiv-mone}}

\fact{$\Monesma$ and $\Um$ induce inverse equivalences of homotopy theories between (i) $\D_{\Um}$-diagrams in $\pMulticat$ and $\D$-diagrams in $\MoneMod$ and (ii) $\D_{\Um}$-Mackey functors in $\pMulticat$ and $\D$-Mackey functors in $\MoneMod$.}{\ref{mackey-pmulti-mone}}

\medskip
\chapNumName{ch:prelim}

\fact{\thm{Grothendieck's Axiom of Universes} Every set belongs to some universe.}{\ref{conv:universe}}

\fact{Each monoidal category satisfies $\lambda_{\tu} = \rho_{\tu}$.}{\ref{lambda=rho}}

\fact{Each monoidal category satisfies the left and right unity properties.}{\ref{moncat-other-unit-axioms}}

\fact{Each braided monoidal category satisfies $\rho = \lambda \xi_{-,\tu}$ and $\lambda = \rho \xi_{\tu,-}$}{\ref{braidedunity}}

\fact{A symmetric monoidal category is precisely a braided monoidal category that satisfies the symmetry axiom.}{\ref{rk:smcat}}

\fact{$\Cat$ is a symmetric monoidal closed category.}{\ref{ex:cat}}

\fact{Iterated monoidal products are left normalized.}{\ref{expl:leftbracketing}}

\fact{$\Cat$ is a 2-category.}{\ref{ex:catastwocategory}}

\fact{$\permcat$, $\permcatst$, and $\permcatsu$ are 2-categories.}{\ref{def:permcat}}

\medskip
\chapNumName{ch:prelim_enriched}

\fact{A locally small 2-category is precisely a $\Cat$-category.}{\ref{locallysmalltwocat}}

\fact{$\VCat$ is a 2-category.}{\ref{ex:vcatastwocategory}}

\fact{The tensor product is a 2-functor on $\VCat$.}{\ref{vtensoriifunctor}}

\fact{$\VCat$ is a monoidal category if $\V$ is braided monoidal.  It is symmetric monoidal if $\V$ is.}{\ref{theorem:vcat-mon}}

\fact{$\VCat$ is a monoidal $\Cat$-category if $\V$ is braided monoidal.  It is a symmetric monoidal $\Cat$-category if $\V$ is symmetric monoidal.}{\ref{theorem:vcat-cat-mon}}

\fact{If $\V$ is a braided monoidal category, then there is a 2-category of small monoidal $\V$-categories.  If $\V$ is a symmetric monoidal category, then there are a 2-category of small braided monoidal $\V$-categories and a 2-category of small symmetric monoidal $\V$-categories.}{\ref{thm:vsmcat}}

\fact{Evaluation is the counit of an adjunction.}{\ref{evaluation}}

\fact{For a symmetric monoidal closed category $\V$, the canonical self-enrichment is a symmetric monoidal $\V$-category.}{\ref{theorem:v-closed-v-sm}}

\fact{Each monoidal functor induces a change-of-enrichment 2-functor.}{\ref{proposition:U-VCat-WCat}}

\fact{Change-of-enrichment 2-functors are closed under composition.}{\ref{proposition:change-enr-horiz-comp}}

\fact{For a braided monoidal functor $U$, change of enrichment is a monoidal $\Cat$-functor, which is symmetric if $U$ is.}{\ref{theorem:U-braided-mon}}

\fact{Change of enrichment preserves enriched monoidal structure.}{\ref{theorem:KU-monoidal}}

\fact{For a symmetric monoidal functor $U \cn \V \to \W$ with $\V$ symmetric monoidal closed, $\Vse_U$ is a symmetric monoidal $\W$-category.}{\ref{VseU}}

\fact{For a monoidal functor $U \cn \V \to \W$ between symmetric monoidal closed categories, the standard enrichment is a monoidal $\W$-functor, which is symmetric if $U$ is.}{\ref{proposition:U-std-enr}}

\medskip
\chapNumName{ch:prelim_multicat}

\fact{Each non-symmetric $\V$-multicategory has an underlying $\V$-category.}{\ref{ex:unarycategory}}

\fact{The terminal multicategory $\Mterm$ consists of a single object and a single $n$-ary operation for each $n$.}{\ref{definition:terminal-operad-comm}}

\fact{Each object in a $\V$-multicategory generates an endomorphism $\V$-operad.}{\ref{example:enr-End}}

\fact{Each $\V$-multifunctor restricts to a $\V$-functor.}{\ref{ex:un-v-functor}}

\fact{Each $\V$-multinatural transformation restricts to a $\V$-natural transformation.}{\ref{ex:un-v-nat}}

\fact{There is a 2-category with (non-symmetric) small $\V$-multicategories as objects.}{\ref{v-multicat-2cat}}

\fact{The initial $\V$-multicategory has an empty set of objects.  The terminal $\V$-multicategory has one object and each multimorphism object given by the terminal object in $\V$.}{\ref{ex:vmulticatinitialterminal}}

\fact{A $\Cat$-multinatural transformation consists of component 1-ary 1-cells that satisfy two $\Cat$-naturality conditions for objects and morphisms.}{\ref{expl:catmultitransformation}}

\fact{Each symmetric monoidal category has an endomorphism multicategory.  Each symmetric monoidal functor induces a multifunctor.  Each monoidal natural transformation between symmetric monoidal functors induces a multinatural transformation.}{\ref{ex:endc}}

\fact{The endomorphism multicategory defines a 2-functor.}{\ref{endtwofunctor}}

\fact{Each symmetric monoidal $\V$-category induces a $\V$-multicategory.}{\ref{proposition:monoidal-v-cat-v-multicat}}

\fact{A pointed structure on a multicategory consists of a basepoint object and $n$-ary basepoint operations that satisfy symmetry, unity, and composition axioms.}{\ref{expl:ptmulticat}}

\fact{Each symmetric monoidal category induces a pointed endomorphism multicategory with the basepoint object given by the monoidal unit.  Each strictly unital symmetric monoidal functor induces a pointed multifunctor.}{\ref{ex:endstc}}

\fact{There is a 2-category with small pointed multicategories as objects.}{\ref{thm:pmulticat}}

\fact{The pointed endomorphism multicategory defines a 2-functor.}{\ref{endsttwofunctor}}

\fact{There is a forgetful 2-functor $\Ust \cn \pMulticat \to \Multicat$.}{\ref{Usttwofunctor}}

\fact{Adjoining a basepoint is a left 2-adjoint of $\Ust$.}{\ref{dplustwofunctor}}

\chapter*{List of Notations}

\newcommand{\entry}[3]{\scalebox{.75}{#1} \> \> \scalebox{.75}{\pageref{#2}} \> 
  \> \scalebox{.75}{#3}\\[-.2ex]}
\newcommand{\entryNoPage}[3]{\scalebox{.75}{#1} \> \> \> \> \scalebox{.75}{#3}\\[-.2ex]} 
\newcommand{\blob}{\> \> \> \> }
\newcommand{\stuff}[1]{\blob \scalebox{.75}{#1}\\[-.2ex]}
\newcommand{\Header}{\>\> \textbf{Page} \>\> \textbf{Description}\\}
\newcommand{\newchHeader}[1]{\blob\\ \textbf{\Cref{#1}}\Header}
\newcommand{\newpart}[1]{\blob\\ \textbf{\Cref{#1}}\blob\\}

\begin{tabbing}
\phantom{\textbf{Notation}} \= \hspace{1.5cm}\= \phantom{\textbf{Page}}\= \hspace{.5cm}\= \phantom{\textbf{Description}} \\

\blob\\
\textbf{Standard Notations}\blob\textbf{Description}\\
\entryNoPage{$\Ob(\C)$, $\Ob\C$}{not:objects}{objects in a category $\C$}
\entryNoPage{$\C(X,Y)$, $\C(X;Y)$}{not:morphisms}{set of morphisms $X \to Y$}
\entryNoPage{$1$, $1_X$}{not:idmorphism}{identity morphism}
\entryNoPage{$\dom(f)$, $\codom(f)$}{not:domain}{domain and codomain of a morphism}
\entryNoPage{$g \circ f$, $gf$}{notation:morphism-composition}{composition of morphisms}
\entryNoPage{$\iso$, $\fto{\iso}$}{not:iso}{an isomorphism}
\entryNoPage{$\sim$, $\fto{\sim}$}{not:sim}{an equivalence}
\entryNoPage{$F \cn \C \to \D$}{def:functors}{a functor}
\entryNoPage{$\Id_{\C}$, $1_{\C}$}{not:idc}{identity functor}
\entryNoPage{$\boldone$}{ex:terminal-category}{terminal category}
\entryNoPage{$(\Set,\times,*)$}{notation:set}{category of sets and functions}
\entryNoPage{$\theta_X$}{thetax}{a component of a natural transformation $\theta$}
\entryNoPage{$1_F$}{not:idf}{identity natural transformation}
\entryNoPage{$\phi\theta$}{not:vcomp}{vertical composition of natural transformations}
\entryNoPage{$\theta' * \theta$}{not:hcomp}{horizontal composition of natural transformations}
\entryNoPage{$(L,R)$, $L \dashv R$}{notation:adjunction}{an adjunction}
\entryNoPage{$\eta$, $\epz$}{adjunction-unit}{unit and counit of an adjunction}
\entryNoPage{$\varnothing$, $\varnothing^{\C}$}{not:initialobj}{an initial object}
\entryNoPage{$\coprod$, $\amalg$}{not:coprod}{a coproduct}
\entryNoPage{$\prod$, $\smallprod$}{not:coprod}{a product}
\entryNoPage{$\Sigma_n$}{not:Sigman}{symmetric group on $n$ letters}

\newchHeader{ch:motivations}

\entry{$G$}{not:G}{finite group}
\entry{$\OrbG$}{not:OrbG}{orbit category of $G$}
\entry{$\TopG$}{not:TopG}{category of $G$-spaces and equivariant morphisms}
\entry{$X^H$}{not:XH}{$H$-fixed point space of $X$}
\entry{$\Phi X$}{eq:PhiX}{fixed point functor of a $G$-space $X$}
\entry{$\Ab$}{not:Ab}{category of Abelian groups}
\entry{$\hty_Q$}{not:htyQ}{a chain of Quillen equivalences}
\entry{$\FinG$}{not:FinG}{skeleton of the category of finite $G$-sets}
\entry{$\Span(\C)$}{not:SpanC}{bicategory of spans in $\C$}
\entry{$\GB$}{not:GB}{Burnside category of $G$}
\entry{$\GA$}{not:GA}{Burnside ring of $G$}
\entry{$\Si^\infty G/H_+$}{not:suspGH}{equivariant suspension spectrum}
\entry{$(M_*,M^*)$}{not:Mstar}{covariant and contravariant functors of an Abelian Mackey functor}
\entry{$\GE$}{not:GE}{Burnside 2-category}
\entry{$\bbK$}{not:bbK}{non-symmetric $K$-theory multifunctor in \cite{guillou_may,gmmo}}
\entry{$\GSp$}{not:GSp}{category of $G$-spectra}
\entry{$\Spm$}{not:Spm}{category of symmetric spectra over $\M$}
\entry{$\EP$}{not:EP}{spectral endomorphism category}

\newpart{part:background}
\newchHeader{ch:ptmulticat}

\entry{$\angc \otimes \angd$, $\ang{c} \otimes^\transp \ang{d}$}{not:angtensor}{$\ang{\ang{(c_i,d_j)}_{i=1}^m}_{j=1}^n$ and $\ang{\ang{(c_i,d_j)}_{j=1}^n}_{i=1}^m$}
\entry{$\xitimes$, $\xitimes_{m,n}$}{xitimesmn}{transpose permutation}
\entry{$\Vt X$}{definition:multigraph}{class of vertices of a multigraph $X$}
\entry{$\MGraph$}{not:mgraph}{category of small multigraphs}
\entry{$X \amtimes Y$}{not:amtimes}{an auxiliary product of multigraphs}
\entry{$\M\shtimes\N$}{definition:sharp-prod}{sharp product of small multicategories}
\entry{$\M \otimes \N$}{bvtensor}{tensor product of small multicategories}
\entry{$(\Multicat, \otimes, \Mtu, \beta)$}{theorem:multicat-symmon}{symmetric monoidal category of small multicategories}
\entry{$\ang{F}c$}{not:angFc}{$\ang{F_ic}_{i=1}^m$ for $\angF = \ang{F_i}_{i=1}^m$}
\entry{$\ang{Fc}$, $\ang{Fc}^{\transp}$}{not:angofFc}{$\ang{\ang{F_i c_j}_{i=1}^m}_{j=1}^n$ and $\ang{\ang{F_i c_j}_{j=1}^n}_{i=1}^m$}
\entry{$\Hom(\M,\N)$}{definition:multicat-hom}{internal hom multicategory}
\entry{$\M \wed \N$}{multicat-wedge}{wedge product of small pointed multicategories}
\entry{$\vpi_{\M,\N}$}{eq:multicat-smash-pushout}{multifunctor $\M \otimes \N \to \M \sma \N$}
\entry{$\Mtup$}{eq:smashunit}{smash unit $\Mtu \bincoprod \Mterm$}
\entry{$\pHom(\M,\N)$}{eq:multicat-pHom}{pointed internal hom multicategory}
\entry{$(\pMulticat, \sma, \Mtup, \pHom)$}{thm:pmulticat-smclosed}{symmetric monoidal category of small pointed multicategories}
\entry{$a^\punc$}{aflat}{$a \setminus \{*\}$ for a pointed finite set $a$}
\entry{$\cM a$}{aflat}{partition multicategory of $a$}
\entry{$\Mone$}{ex:mofone}{partition multicategory of $\ord{1} = \{0,1\}$}
\entry{$\iota^n$, $\pi^n_j$}{not:moneoperations}{operations in $\Mone$}
\entry{$\txprod_{a,b}$}{eq:part-prod}{partition product $\cM a \sma \cM b \to \cM(a \sma b)$}
\entry{$\cM$}{proposition:cM-symm-mon}{symmetric monoidal functor $\Fskel^\op \to \pMulticat$}
\entry{$\cM^0$}{eq:cM0}{unit constraint $\Mtup \to \Mone$}
\entry{$\cM^2_{\ord{m},\ord{n}}$}{eq:cM2}{monoidal constraint $\cM\ord{m} \sma \cM\ord{n} \to \cM(\ord{mn})$}
\entry{$\MoneMod$}{definition:MoneMod-prelim}{2-category of left $\Mone$-modules}
\entry{$\Endm(-)$}{ex:endmc}{endomorphism left $\Mone$-module}
\entry{$(\MoneMod,\sma,\Mone,\pHom)$}{monebicomplete}{symmetric monoidal category of left $\Mone$-modules}
\entry{$\Monesma$}{MonesmaUmadj}{left 2-adjoint $\pMulticat \to \MoneMod$ of $\Um$}
\entry{$\etahat$}{etam-def}{unit of $(\Monesma) \dashv \Um$}
\entry{$\epzhat$}{epzm-def}{counit of $(\Monesma) \dashv \Um$}
\entry{$\ang{x \compk y}$, $\angx \compk y$}{not:compk}{replacing the $k$-th entry of $\angx$ by $y$}
\entry{$(P, \{P^2_j\}_{j=1}^n)$}{def:nlinearfunctor}{an $n$-linear functor with $j$-th linearity constraint $P^2_j$}
\entry{$\permcatsu\scmap{\angC;\D}$}{definition:permcatsus-homcat}{category of $n$-linear functors and transformations}
\entry{$\permcatst\scmap{\angC;\D}$}{definition:permcatsus-homcat}{category of strict $n$-linear functors and transformations}
\entry{$P^\sigma$, $P \circ \sigma$}{fsigmatwoj}{right $\sigma$-action on $P$}
\entry{$\permcatsu$, $\permcatst$}{thm:permcatmulticat}{$\Cat$-multicategories of small permutative categories}

\newchHeader{ch:Kspectra}

\entry{$\Ner$}{nerve}{nerve functor $\Cat \to \sSet$}
\entry{$\ord{n}$}{ordn}{pointed finite set $\{0,1,\ldots,n\}$}
\entry{$\boldtwo$}{def:bisimplicial}{nerve of the category with two isomorphic objects}
\entry{$(\C,\cW)$}{relcat-i}{a relative category with stable equivalences $\cW$}
\entry{$(\C,\cW)^\Dgm$}{def:relcat-i}{a relative diagram category}
\entry{$\Nde(\C,\cW)$}{def:relcat-ii}{classification diagram of $(\C,\cW)$}
\entry{$\pC$}{def:pointed-objects}{category of pointed objects in $\C$ with terminal object $\term$}
\entry{$(a,i^a)$}{not:aia}{a pointed object with pointed structure $i^a$}
\entry{$a \wed b$}{not:awedgeb}{wedge of pointed objects}
\entry{$a \sma b$}{eq:smash}{smash product of pointed objects}
\entry{$\stu$}{eq:smash}{smash unit $\tu \bincoprod \term$}
\entry{$\pHom(a,b)$}{eq:pHom-pullback}{pointed Hom}
\entry{$(\pC,\sma,\stu,\pHom)$}{theorem:pC-sm-closed}{symmetric monoidal category of pointed objects}
\entry{$\zob$}{definition:zero}{a zero object}
\entry{$\C^\punc(x,y)$}{Cpuncxy}{set of nonzero morphisms $\C(x,y) \setminus \{\zob\}$}
\entry{$\Dhat$}{eq:Dhat-ptd-unitary-enr}{pointed unitary enrichment}
\entry{$(\C,*)$}{def:pointed-category}{a pointed category with chosen object $*$}
\entry{$\du$}{eq:Dgm-pV-unit}{monoidal unit diagram}
\entry{$A \sma B$}{eq:Dgm-pV-convolution}{pointed Day convolution}
\entry{$\Homdstar(A,B)$}{eq:Dgm-pV-hom}{pointed hom diagram}
\entry{$\Mapdstar(A,B)$}{eq:Dgm-pV-map}{pointed mapping object}
\entry{$(L_\pu, \ev_\pu)$}{eq:pV-DstarV-adj}{an adjunction between $\pV$ and $\DstarV$}
\entry{$\GaV$}{not:GaV}{category $\FstarV$ of $\Ga$-objects in $\V$}
\entry{$\big(\GaV, \sma, \du, \Homfstar\!\big)$}{GammaV}{symmetric monoidal category of $\Ga$-objects}
\entry{$\Homfstar$, $\Mapfstar$}{eq:GaV-hom}{pointed hom diagram and mapping object of $\GaV$}
\entry{$\Inj$}{definition:ufs}{category of unpointed finite sets and injections}
\entry{$\ufs{n}$}{ufsn}{unpointed finite set $\{1,\ldots,n\}$}
\entry{$\ordtu{n}$}{ordtunpsi}{$q$-tuple $\ang{\ord{n}_k}_{k=1}^q$ of pointed finite sets}
\entry{$f_* \ordtu{n}$}{not:fstarordtun}{$\ang{\ord{n}_{f^\inv(j)}}_{j=1}^p$}
\entry{$\Fskel^{(q)}$}{definition:Fskel-roundsma}{$q$-th smash power of $\Fskel$}
\entry{$\gst$}{not:ObFzero}{basepoint object of $\Fskel^{(0)}$}
\entry{$(\Gskel,\gst)$}{definition:Gstar}{category of tuples of pointed finite sets}
\entry{$(\Gskel,\oplus,\ang{},\xi)$}{definition:concatenation-product}{permutative structure on $\Gskel$ with concatenation product $\oplus$}
\entry{$\tau_{q,q'}$}{tauqqprime}{block permutation swapping $q$ and $q'$ elements}
\entry{$\GstarV$}{not:GstarV}{category of $\Gstar$-objects in $\V$}
\entry{$\Homgstar$, $\Mapgstar$}{eq:GstarV-hom}{pointed hom diagram and mapping object of $\GstarV$}
\entry{$i$}{i-incl}{length-one inclusion $(\Fskel,\ord{0}) \to (\Gskel,\gst)$}
\entry{$\sma$}{smagf}{strict symmetric monoidal pointed functor $\Gskel \to \Fskel$}
\entry{$\Sp$, $\Spc$}{SymSp}{category of (connective) symmetric spectra}
\entry{$\Kse$}{Kse}{Segal $K$-theory}
\entry{$\Jse$}{Jse}{Segal $J$-theory}
\entry{$\Ner_*$}{NerGa}{levelwise nerve $\Gacat \to \Gasset$}
\entry{$\Kf$}{Kf-sm}{functor $\Gasset \to \Spc$}
\entry{$\Kem$}{Kem}{Elmendorf-Mandell $K$-theory}
\entry{$\Jt$}{Jt-smcat}{functor $\MoneMod \to \Gstarcat$}
\entry{$\Jem$}{Jem}{Elmendorf-Mandell $J$-theory}
\entry{$\Ner_*$}{NerGs}{levelwise nerve $\Gstarcat \to \Gstarsset$}
\entry{$\Kg$}{Kg-sm}{functor $\Gstarsset \to \Sp$}
\entry{$\Pi^*$}{not:Pistar}{a natural transformation $\sma^* \circ \Jse \to \Jem$}
\entry{$L$}{leftadjL}{left adjoint of $i^*$}
\entry{$\cS$}{perm-steq}{stable equivalences in $\permcatsu$, $\Gacat$, $\Gasset$, and $\Spc$}
\entry{$\csi$}{not:csi}{$i^*$-stable equivalences in $\Gstarcat$ and $\Gstarsset$}
\entry{$\cP$}{PinvK}{inverse $K$-theory $\Gacat \to \permcatsu$}
\entry{$S_*$}{Sstar-mandell}{homotopy inverse $\Gasset \to \Gacat$ of $\Ner_*$}
\entry{$\bA$}{A-segal}{homotopy inverse $\Spc \to \Gasset$ of $\Kf$}

\newchHeader{ch:multperm}

\entry{$\ang{x}_{f^\inv(j)}$}{eq:x-finv}{sub-tuple $\ang{x_i}_{i \in f^\inv(j)}$}
\entry{$\sigma^k_{g,f}$}{eq:sigma-kgf}{right permutation from $\bigoplus_{j \in g^\inv(k)} \angx_{f^\inv(j)}$ to $\angx_{(gf)^\inv(k)}$}
\entry{($\Fr\M,\oplus,\ang{},\xi)$}{definition:free-perm}{free permutative category of a multicategory $\M$}
\entry{$(f,\ang{\phi})$}{FMmorphism}{a morphism in $\Fr\M$}
\entry{$\Fr(\Mtu)$}{example:free-Mtu}{free permutative category of the initial operad $\Mtu$}
\entry{$\Fr\Mterm$}{example:free-Mterm}{free permutative category of the terminal multicategory $\Mterm$}
\entry{$\Fr$}{proposition:free-perm-functor}{free permutative category 2-functor $\Multicat \to \permcatst$}
\entry{$\eta_\M$}{definition:eta}{unit multifunctor $\M \to \End\Fr\M$}
\entry{$\epz_\C$}{definition:epz}{counit strict symmetric monoidal functor $\Fr\End(\C) \to \C$}
\entry{$\Fr \dashv \End$}{theorem:FE-adj}{2-adjunction $\Multicat \lradj \permcatst$}
\entry{$\vrho_\C$}{rhoc-assignments}{right adjoint $\C \to \Fr\End(\C)$ of $\epz_\C$}
\entry{$\vrho^0_\C$, $\vrho^2_\C$}{def:rhoc-monoidal}{unit and monoidal constraints of $\vrho_\C$}
\entry{$x^{\onen}_{\jonejn}$}{xonenjonen-def}{$n$-tuple $\bang{x^i_{j_i}}_{i=1}^n$}
\entry{$\ang{x^{\onen}}$}{angxonen-def}{tensor product $\txotimes_{i=1}^n \ang{x^i}$ of tuples}
\entry{$r_{\onen}$}{ronen}{$\txprod_{i=1}^n r_i$}
\entry{$f^{\onen}$}{eq:fonen}{index map induced by $\txprod_{i=1}^n f^i$}
\entry{$\ang{x^{\onen}}_{f;\,\konekn}$}{angxonen-fkonen}{$\txotimes_{i=1}^n \ang{x^i_j}_{j \in (f^i)^\inv(k_i)}$}
\entry{$\phi^{\onen}_{\konekn}$}{eq:phionenkonekn}{$\txotimes_{i=1}^n \phi^i_{k_i}$}
\entry{$\ang{\phi^{\onen}}$}{eq:phionen}{$\txotimes_{i=1}^n \ang{\phi^i}$}
\entry{$\brb{\Frn, \bang{(\Frn)^2_p}_{p=1}^n}$}{def:S-multi}{strong $n$-linear functor $\txprod_{i=1}^n \Fr\M_i \to \Fr\big(\txotimes_{i=1}^n \M_i \big)$}
\entry{$\Fr^0$}{not:Frzero}{0-linear functor $\boldone \to \Fr(\Mtu)$ determined by length-one tuple $(*)$}
\entry{$\ang{\hat{x}^\onen}$, $\ang{\tilde{x}^\onen}$}{eq:xhatonen-xtilonen}{analogs of $\ang{x^{\onen}}$ defined using $\ang{\hat{x}^p}$ and $\ang{\tilde{x}^p}$}
\entry{$\rho_{r_p,\hat{r}_p}$}{rhorbrhatb}{unique permutation determined by the (co)domain of $(\Frn)^2_p$}
\entry{$\Frbar$}{convention:Fun}{$\Fr$ on $\Multicat(\M,\N)$}
\entry{$\Fr$}{eq:thm-F-multi}{non-symmetric $\Cat$-multifunctor $\Multicat \to \permcatsu$}
\entry{$\cSI$}{FI}{stable equivalences in $\permcatst$}
\entry{$\cSF$}{FI}{$\Fr$-stable equivalences in $\Multicat$}
\entry{$\N^\P$}{not:NP}{category of $\P$-algebras in $\N$}
\entry{$\N^\Q$}{not:NQ}{category of non-symmetric $\Q$-algebras in $\N$}
\entry{$\cW^\P$}{not:cWP}{wide subcategory of morphisms with each component in $\cW$}
\entry{$\Fr^\Q$, $\End^\Q$}{thm:alg-hty-equiv}{inverse equivalences of homotopy theories induced by $\Fr$ and $\End$}

\newpart{part:multicat}
\newchHeader{ch:ptmulticat-sp}

\entry{$\ang{x}^\wedge$}{definition:removing-basepts}{sub-tuple of non-basepoint objects}
\entry{$(f', \ang{\phi}')$}{not:fphiprime}{$(f,\ang{\phi})$ with basepoint operations removed}
\entry{$\obsim$}{definition:up-to-basept-equiv}{up-to-basepoint equivalence relation on objects of $\Fr\M$}
\entry{$[\ang{x}]$}{not:obsimclass}{$\obsim$-equivalence class of $\angx$}
\entry{$\wt{\Mor}(\Fr\M)$}{not:wtmor}{$\obsim$-composable tuples of morphisms}
\entry{$\isim$}{not:isim}{relation for composition}
\entry{$\iisim$}{not:iisim}{relation for removing basepoint operations}
\entry{$[\uf]$}{not:ufclass}{equivalence class of $\uf$}
\entry{$\Fst\M$}{def:Fst-object}{pointed free permutative category of a pointed multicategory $\M$}
\entry{$(\oplus, [\ang{}], \xi)$}{def:Fst-permutative}{permutative structure on $\Fst\M$}
\entry{$\Fst$}{ptmulticat-thm-i}{2-functor $\pMulticat \to \permcatst$}
\entry{$\pst$}{eq:pst-iinatural}{2-natural transformation $\Fr \to \Fst$}
\entry{$\etast$}{etast-def}{unit $1 \to \Endst\Fst$}
\entry{$\epzst$}{epzst-def}{counit $\Fst\Endst \to 1$}
\entry{$\Fst$, $\Endst$}{ptmulticat-thm-v}{2-adjunction $\pMulticat \lradj \permcatst$}
\entry{$\Fm$, $\Endm$}{ptmulticat-thm-ii}{2-adjunction $\MoneMod \lradj \permcatst$}
\entry{$\etam$}{etam-define}{unit $1 \to \Endm\Fm$}
\entry{$\epzm$}{epzm-define}{counit $\Fm\Endm \to 1$}
\entry{$w^{\ang{x}}$}{eq:w-angx}{a morphism $\ang{x}^\wedge \to \ang{x}$ in $\Fr\M$}
\entry{$\ang{x}^{\sim}$}{eq:angxtil}{$(*)$, $\ang{}$, or $\ang{x}^\wedge$}
\entry{$c^{\ang{x}}$}{lemma:c-angx}{a morphism $\ang{x} \to \ang{x}^\sim$ in $\Fr\M$}
\entry{$\vrhost_\C$}{vrhostC-def}{strictly unital right adjoint $\C \to \Fst\Endst\C$ of $\epzst_\C$}
\entry{$\upsst$}{not:upsst}{unit $1 \to \vrhost_\C\epzst_\C$}
\entry{$\cSM$}{SM-Sst-def}{$\Fm$-stable equivalences in $\MoneMod$}
\entry{$\cSst$}{SM-Sst-def}{$\Fst$-stable equivalences in $\pMulticat$}

\newchHeader{ch:ptmulticat-alg}

\entry{$\Fstn$}{def:Sst}{functor $\txprod_{i=1}^n \Fst\M_i \to \Fst\big(\txsma_{i=1}^n \M_i\big)$}
\entry{$\Fst^0$}{not:Fstzero}{functor $\boldone \to \Fst(\Mtup)$ determined by the object $[(*)]$}
\entry{$(\Fstn)^2_p$}{eq:Fstn2p}{$p$-th linearity constraint of $\Fstn$}
\entry{$\Fbst$}{conv:Fbst}{$\Fst$ on $\pMulticat(\M,\N)$}
\entry{$\Fst$}{eq:ptmulticat-xvii}{non-symmetric $\Cat$-multifunctor $\pMulticat \to \permcatsu$}
\entry{$\Fst^\Q$, $\Endst^\Q$}{ptmulticat-xxiii}{inverse equivalences of homotopy theories induced by $\Fst$ and $\Endst$}
\entry{$\Fm$}{Fm-multi-def}{non-symmetric $\Cat$-multifunctor $\MoneMod \to \permcatsu$}
\entry{$\etam$}{etam-multi-def}{non-symmetric $\Cat$-multinatural transformation $1 \to \Endm\Fm$}
\entry{$\vrhom$}{vrhom-multi-def}{non-symmetric $\Cat$-multinatural transformation $1 \to \Fm\Endm$}
\entry{$\Fm^\Q$, $\Endm^\Q$}{ptmulticat-xxv}{inverse equivalences of homotopy theories induced by $\Fm$ and $\Endm$}
\entry{$(\Monesma)^\Q$, $\Um^\Q$}{eq:MonesmaQ-UMQ-hteq}{inverse equivalences of homotopy theories induced by $\Monesma$ and $\Um$}

\newpart{part:enrpresheave}
\newchHeader{ch:menriched}

\entry{$(\C,\mcomp,i)$}{def:menriched-cat}{a category enriched in a non-symmetric multicategory}
\entry{$\C(x,y)$, $\C_{x,y}$}{not:Cxyhom}{hom object with domain $x$ and codomain $y$}
\entry{$\mcomp_{x,y,z}$}{menriched-cat-comp}{composition for objects $x$, $y$, and $z$}
\entry{$i_x$}{menriched-cat-identity}{identity of an object $x$}
\entry{$F_{x,y}$}{Fxy-component}{$(x,y)$-component of an enriched functor $F$}
\entry{$1_\C$}{id-mfunctor}{identity enriched functor}
\entry{$GF$}{mfunctor-composition}{composite enriched functor of $F$ and $G$}
\entry{$\theta_x$}{mnatural-component}{$x$-component of an enriched natural transformation $\theta$}
\entry{$1_F$}{id-mnat}{identity enriched natural transformation}
\entry{$\psi\theta$}{mnat-vcomp-diag}{vertical composition of enriched natural transformations}
\entry{$\theta' * \theta$}{def:mnaturaltr-hcomp}{horizontal composition of enriched natural transformations}
\entry{$\MCat$}{mcat-iicat}{2-category of small $\M$-categories}
\entry{$\psu$}{permcatsu-psu}{$\permcatsu$}
\entry{$\mcomp_{x,y,z}^1$, $\mcomp_{x,y,z}^2$}{not:mcompi}{linearity constraints of $\mcomp_{x,y,z}$}
\entry{$\big(\psu(\C,\D),\oplus,\pu,\xi\big)$}{psucd-permutative}{a hom permutative category}
\entry{$(F \oplus G)^2$}{FoplusG-mon-constraint}{monoidal constraint of $F \oplus G$}
\entry{$\pu \cn \C \to \D$}{psucd-unit}{constant functor at the monoidal unit of $\D$}
\entry{$\mcomp_{\B,\C,\D}$}{mBCD}{composition bilinear functor for $\B$, $\C$, and $\D$}
\entry{$\mcomp^2_1$, $\mcomp^2_2$}{HFIF}{linearity constraints of $\mcomp_{\B,\C,\D}$}
\entry{$\ev_{\C,\D}$}{evCD}{evaluation bilinear functor for $\C$ and $\D$}
\entry{$(\ev_{\C,\D})^2_1$, $(\ev_{\C,\D})^2_2$}{evCD}{linearity constraints of $\ev_{\C,\D}$}
\entry{$(\Cop,\mcompop,i)$}{def:opposite-mcat}{opposite enriched category}

\newchHeader{ch:change_enr}

\entry{$\dF$}{not:dF}{change of enrichment along $F$}
\entry{$(\C_F,\mcomp_F,i_F)$}{mcat-change-enrichment}{change of enrichment of $(\C,\mcomp,i)$ along $F$}
\entry{$H_F$}{mfunctor-change-enrichment}{change of enrichment of $H$ along $F$}
\entry{$\theta_F$}{mnat-change-enrichment}{change of enrichment of $\theta$ along $F$}
\entry{$\dFst$}{dFst-explicit}{change of enrichment along $\Fst$}
\entry{$\pM$}{not:pM}{$\pMulticat$}
\entry{$(\Cop)_F$}{dF-opposite}{change of enrichment of $\Cop$}
\entry{$(\C_F)^\op$}{dF-opposite}{opposite of $\C_F$}
\entry{$\dU$, $\dEndU$}{mon-change-enrichment}{change of enrichment along $U$ and $\EndU$}
\entry{$\dtheta$}{dtheta-def}{2-natural transformation induced by $\theta$}
\entry{$\Enr$}{change-enr-twofunctor}{2-functor $\Multicatns \to \iicat$ sending $\M$ to $\MCat$}

\newchHeader{ch:gspectra}

\entry{$\big(\M \scs \clM \scs \ev\big)$}{def:closed-multicat}{a closed multicategory}
\entry{$\clM\scmap{\angx;y}$, $\clM_{\angx;\,y}$}{clMangxy}{an $n$-ary internal hom object}
\entry{$\ev_{\angx;\,y}$}{evangxy}{multicategorical evaluation}
\entry{$\chi_{\angx;\, \angy;\, z}$}{eval-bijection}{evaluation bijection}
\entry{$\pn{f}$}{not:partner}{partner of $f$}
\entry{$\Wald$}{ex:waldhausen}{closed multicategory of small Waldhausen categories}
\entry{$\big(\EndV \scs \clEndV \scs \ev\big)$}{smclosed-closed-multicat}{endomorphism closed multicategory}
\entry{$\clp\scmap{\angC;\D}$}{clp-internal-hom}{internal hom permutative category}
\entry{$P \oplus Q$}{PplusQfunctor}{monoidal product of $n$-linear functors $P$ and $Q$}
\entry{$\clpsu$}{clp-internal-hom}{$\clp$}
\entry{$(P \oplus Q)^2_i$}{PplusQtwoi}{$i$-th linearity constraint of $P \oplus Q$}
\entry{$\clpu$}{clpu}{constant functor $\txprod_{i=1}^n \C_i \to \D$ at the monoidal unit in $\D$}
\entry{$\clxi_{P,Q}$}{not:clxi}{$(P,Q)$-component of the braiding $\clxi$}
\entry{$\ev_{\angC;\,\D}$}{clp-ev-angcd}{multicategorical evaluation for $\permcatsu$}
\entry{$(\ev_{\angC;\,\D})^2_j$}{evangcd-iplusone-linearity}{linearity constraints of $\ev_{\angC;\,\D}$}
\entry{$\chi_{\angC;\,\angD;\,\B}$}{chi-angCDB}{evaluation bijection for $\permcatsu$}
\entry{$(\chi P)^2_r$}{chi-of-P}{linearity constraints of $\chi P$}
\entry{$\Psi$}{chiPsi}{inverse of $\chi$}
\entry{$\chi'$}{chiprimeP}{$\chi_{\angC\sigma;\, \angD\vsi;\, \B}$}

\newchHeader{ch:std_enrich}

\entry{$(\M,\comp,i)$}{def:clmulti-clcat}{canonical self-enrichment of $\M$}
\entry{$\Fse$}{std-enr-F}{standard enrichment of $F$}
\entry{$\Fse_{x,y}$}{Fprimexy}{$(x,y)$-component of $\Fse$}
\entry{$\Fstse$}{Fstse-definition}{standard enrichment of $\Fst \cn \pMulticat \to \permcatsu$}
\entry{$\EndUse$}{std-enr-monoidal}{standard enrichment of $\EndU$}
\entry{$\Fse_G$}{expl:FseG}{change of enrichment of $\Fse$ along $G$}
\entry{$\Kemse$}{Kemse-factor-four}{standard enrichment of $\Kem \cn \psu \to \Sp$}
\entry{$?^{\wedge}$}{not:se}{standard enrichment $\se{?}$}
\entry{$\Phi$}{not:BOPhi}{Bohmann-Osorno spectral functor}
\entry{$\bbK$}{not:KGM}{non-symmetric $K$-theory multifunctor}
\entry{$\Kgse$}{expl:Kgse}{standard enrichment of $\Kg \cn \Gstarsset \to \Sp$}
\entry{$\Endmse$}{Endm-stdenr}{standard enrichment of $\Endm \cn \psu \to \MoneMod$}
\entry{$\Jtse$}{expl:Jtse}{standard enrichment of $\Jt \cn \MoneMod \to \Gstarcat$}
\entry{$\Nerstarse$}{expl:Nerstarse}{standard enrichment of $\Ner_* \cn \Gstarcat \to \Gstarsset$}

\newchHeader{ch:gspectra_Kem}

\entry{$\MCat(\C,\M)$}{mcat-cm}{$\C$-diagram category}
\entry{$\MCat(\Cop,\M)$}{mcat-copm}{$\C$-presheaf category, $\C$-Mackey functor category}
\entry{$\Fdg$, $\Fdg^\C$}{diag-change-enr}{diagram change of enrichment of $F$}
\entry{$\cPdg$}{not:cPdg}{diagram change of enrichment of inverse $K$-theory $\cP$}
\entry{$\dsub{\cP}$}{not:APPse}{change of enrichment along $\cP$}
\entry{$\cPse$}{not:APPse}{standard enrichment of $\cP$}
\entry{$\Kemdg$}{Kemdg-functor}{diagram and presheaf change of enrichment of $\Kem$}
\entry{$\bKdg$}{bKdg}{diagram change of enrichment of $\bbK$}

\newpart{part:homotopy-mackey}
\newchHeader{ch:mackey}

\entry{$\uni$, $\cou$}{mackey-context-iv}{multinatural transformations $1_\M \to EF$ and $1_\N \to FE$}
\entry{$\Fdgr$}{Fdgr-def}{a functor $\MCat(\C_E,\M) \to \NCat(\C,\N)$}
\entry{$\C_\cou$}{C-sub-cou}{$\C$-component of $\dcou \cn 1_{\NCat} \to \dFE$}
\entry{$\unidg$}{unidg-def}{natural transformation $1 \to \Edg\Fdgr$}
\entry{$\coudg$}{coudg-def}{natural transformation $1 \to \Fdgr\Edg$}
\entry{$\cWtri$}{cWtri-subcat}{a wide subcategory of $\PCat(\D,\P)$ induced by $\cW \bigsubset \cP$}
\entry{$F^\inv\cX$}{FinvcX}{$F$-stable equivalences created by $F$}
\entry{$\FinvcXtri$}{diag-enr-rel-cat}{a wide subcategory of $\MCat(\C_E,\M)$ induced by $F^\inv\cX$}

\newchHeader{ch:mackey_eq}

\entry{$\cStri$}{sStri-definition}{a wide subcategory of $\permcatsucat(\C,\permcatsu)$}
\entry{$\cSsttri$}{sSsttri-definition}{a wide subcategory of $\pMulticatcat\big(\C_{\Endst},\pMulticat\big)$}
\entry{$\Fstdgr$, $\Endstdg$}{mackey-xiv-pmult-functors}{inverse equivalences of homotopy theories induced by $\Fst$ and $\Endst$}
\entry{$\dEndst$}{expl:dEndst}{change of enrichment along $\Endst \cn \psu \to \pM$}
\entry{$\Endstse$}{expl:Endstse}{standard enrichment of $\Endst$}
\entry{$\cSmtri$}{sSmtri-def}{a wide subcategory of $\Monemodcat\big(\C_{\Endm},\MoneMod\big)$}
\entry{$\Fmdgr$, $\Endmdg$}{mackey-xiv-mone-functors}{inverse equivalences of homotopy theories induced by $\Fm$ and $\Endm$}
\entry{$\dEndm$}{Endmdg-A-composite}{change of enrichment along $\Endm \cn \psu \to \MoneMod$}
\entry{$\C_{\vrhom}$}{C-vrhom}{$\C$-component of $\dvrhom \cn 1 \to \dFm \, \dEndm$}
\entry{$\dFm$}{dFm}{change of enrichment along $\Fm \cn \MoneMod \to \psu$}
\entry{$\cSsttri$}{sSsttri-defin}{a wide subcategory of $\pMulticatcat\big(\D_{\Um},\pMulticat\big)$}
\entry{$\cSmtri$}{sSmtri-defin}{a wide subcategory of $\Monemodcat\big(\D,\MoneMod\big)$}
\entry{$\Monesmadgr$, $\Umdg$}{mackey-pmult-mone-functors}{inverse equivalences of homotopy theories induced by $\Monesma$ and $\Um$}
\entry{$\dUm$}{Umdg-def}{change of enrichment along $\Um \cn \MoneMod \to \pM$}
\entry{$\Umse$}{Umdg-def}{standard enrichment of $\Um$}
\entry{$\D_{\epzhatinv}$}{HP-composite}{$\D$-component of $\depzhatinv \cn 1 \to \dMonesma\, \dUm$}
\entry{$\dMonesma$}{dMonesmas}{change of enrichment along $\Monesma \cn \pM \to \MoneMod$}

\newchHeader{ch:prelim}

\entry{$\calu$}{not:universe}{a Grothendieck universe}
\entry{$(\C,\otimes,\tu,\alpha,\lambda,\rho)$}{not:monoidalcat}{a monoidal category}
\entry{$(x,\mcomp,i)$}{notation:monoid}{a monoid}
\entry{$(a,\mu)$}{not:leftmodule}{a left module}
\entry{$\xi$}{notation:symmetry-iso}{symmetry isomorphism or braiding}
\entry{$\oplus$, $\pu$}{not:opluse}{monoidal product and monoidal unit in a permutative category}
\entry{$\big(\Fskel, \sma, \ord{1}, \xi\big)$}{not:Fskel}{permutative category of pointed finite sets}
\entry{$\ord{n}$}{not:ordn}{pointed finite set $\{0,\ldots,n\}$}
\entry{$[-,-]$, $\Hom(-,-)$}{notation:internal-hom}{internal hom}
\entry{$\BC$}{not:DC}{diagram category of functors $\cB \to \C$}
\entry{$(\Cat, \times, \boldone,[,])$}{not:cat}{category of small categories}
\entry{$(F,F^2,F^0)$}{def:monoidalfunctor}{a monoidal functor}
\entry{$\Rightarrow$}{twocellnotation}{a (multi)natural transformation, ($n$-ary) 2-cell}
\entry{$\A_0$, $\A_1$, $\A_2$}{def:twocategory}{objects, 1-cells, and 2-cells in a 2-category $\A$}
\entry{$\theta'\theta$}{not:vcompiicell}{vertical composition of 2-cells}
\entry{$gf$}{not:hcompicell}{horizontal composition of 1-cells}
\entry{$\phi * \theta$}{not:hcompiicell}{horizontal composition of 2-cells}
\entry{$\A(a,b)$}{not:homcat}{a hom category}
\entry{$\permcat$, $\permcat^?$}{def:permcat}{2-categories of small permutative categories}
\entry{$\iicat$}{ex:iicat}{2-category of small 2-categories}
\entry{$(F,G,\eta,\epz$)}{def:twoadjunction}{a 2-adjunction}

\newchHeader{ch:prelim_enriched}

\entry{$m$}{not:enrcomposition}{composition in a $\V$-category}
\entry{$i_X$}{not:enridentity}{identity in a $\V$-category}
\entry{$\vtensorunit$}{definition:unit-vcat}{unit $\V$-category}
\entry{$\VCat$}{ex:vcatastwocategory}{2-category of small $\V$-categories}
\entry{$\Cop$}{definition:vcat-op}{opposite $\V$-category}
\entry{$\C \otimes \D$}{definition:vtensor-0}{tensor product of $\V$-categories}
\entry{$\ximid$}{not:ximid}{interchanging the middle two factors}
\entry{$\ell^\otimes$, $r^\otimes$}{definition:vtensor-unitors}{left and right unitors for $\otimes$}
\entry{$a^\otimes$}{definition:vtensor-assoc}{associator for $\otimes$}
\entry{$\beta^\otimes$}{definition:vtensor-beta}{braiding for $\otimes$}
\entry{$?^{-\otimes}$}{not:minusotimes}{inverse $(?^\otimes)^\inv$}
\entry{$a^\vmtimes_1$}{eq:a-vmtimes-inv-mate}{mate of $a^\vmtimes$}
\entry{$\ell^\vmtimes_1$, $r^\vmtimes_1$}{eq:lr-vmtimes-mates}{mates of $\ell^\vmtimes$ and $r^\vmtimes$}
\entry{$\ev$}{evaluation}{evaluation}
\entry{$\coev$}{coevaluation}{coevaluation}
\entry{$\Vse$}{definition:canonical-v-enrichment}{canonical self-enrichment}
\entry{$(-)_U$}{not:dU}{change of enrichment along $U$}
\entry{$\big(\Use, \Use^2, \Use^0\big)$}{not:Use}{standard enrichment of $U$}

\newchHeader{ch:prelim_multicat}

\entry{$\Prof(S)$}{notation:profs}{class of $S$-profiles}
\entry{$\angx$, $\ang{x_j}_{j=1}^n$}{notation:us}{a length-$n$ profile $(x_1,\ldots,x_n)$}
\entry{$\ang{}$}{notation:us}{empty profile}
\entry{$\smscmap{\angx;y}$}{notation:duc}{an element in $\Prof(S) \times S$}
\entry{$\angx \oplus \angy$}{not:concat}{concatenation of profiles}
\entry{$(\M,\ga,\opu)$}{notation:enr-multicategory}{a $\V$-multicategory}
\entry{$\M\smscmap{\angx;x'}$}{notation:enr-cduc}{an $n$-ary operation object}
\entry{$\angx\sigma$}{enr-notation:c-sigma}{$\ang{x_{\sigma(j)}}_{j=1}^n$}
\entry{$\operadunit_x$}{notation:enr-unit-c}{$x$-colored unit}
\entry{$\gamma$}{notation:enr-multicategory-composition}{composition in a $\V$-multicategory}
\entry{$\sigma\ang{k_{\sigma(1)},\ldots,k_{\sigma(n)}}$}{blockpermutation}{block permutation induced by $\sigma$}
\entry{$\tau_1 \times \cdots \times \tau_n$}{blocksum}{block sum}
\entry{$\M_n$}{not:nthobject}{$n$-ary operation object of a $\V$-operad}
\entry{$\Mterm$}{definition:terminal-operad-comm}{terminal multicategory}
\entry{$\End(x)$}{example:enr-End}{endomorphism $\V$-operad of an object $x$}
\entry{$\VMulticat$}{vmulti-i}{2-category of small $\V$-multicategories}
\entry{$\VMulticatns$}{vmulti-ii}{2-category of non-symmetric small $\V$-multicategories}
\entry{$\Multicat$, $\Multicatns$}{multicat-def}{2-category of (non-symmetric) small multicategories}
\entry{$\Mtu$}{ex:initialoperad}{initial operad}
\entry{$\bt$}{not:bt}{a terminal object in $\V$}
\entry{$\angx \to y$}{not:naryonecell}{an $n$-ary 1-cell}
\entry{$F\angx$, $\theta_{\angx}$}{Fangcthetaangc}{$\ang{Fx_j}_{j=1}^k$ and $\ang{\theta_{x_j}}_{j=1}^k$}
\entry{$1_{\theta_{\angx}}$}{catmultinaturalityiicell}{$\ang{1_{\theta_{x_j}}}_{j=1}^k$}
\entry{$\End(-)$}{ex:endc}{endomorphism multicategory}
\entry{$\End(\K)$}{definition:EndK}{endomorphism $\V$-multicategory of $\K$}
\entry{$(\M,i)$}{def:ptd-multicat}{a pointed multicategory with pointed structure $i$}
\entry{$\Endst(-)$}{ex:endstc}{pointed endomorphism multicategory}
\entry{$\pMulticat$}{thm:pmulticat}{2-category of small pointed multicategories}
\entry{$\Ust$}{Usttwofunctor}{forgetful 2-functor $\pMulticat \to \Multicat$}
\entry{$\dplus$}{dplus-def}{adjoining a basepoint $\Multicat \to \pMulticat$}
\entry{$\etaplus$, $\epzplus$}{etast-epzst-def}{unit and counit for $\dplus \dashv \Ust$}

\end{tabbing}

\printindex
\end{document}